\documentclass[12pt]{caltech_thesis}
\usepackage[dvipsnames]{xcolor}
\usepackage[hyphens]{url}
\usepackage{lipsum}
\usepackage{graphicx}

\usepackage{mathtools}
\usepackage{todonotes}

\usepackage{hyperref}
\hypersetup{
    colorlinks=true,      
    linkcolor=Green,       
    citecolor=NavyBlue,      
    filecolor=Plum,    
    urlcolor=Plum         
}
\usepackage[utf8]{inputenc}
\usepackage[T1]{fontenc}
\usepackage{newtxtext,newtxmath}
\usepackage[capitalize,nameinlink]{cleveref}

\usepackage[
    backend=biber,natbib,
    backref=true,
    style=alphabetic
]{biblatex}

\AtEveryBibitem{
    \clearfield{urlyear}
    \clearfield{urldate}
    \clearfield{issn}
    \clearfield{primaryclass}
    \clearfield{publisher}
    \clearfield{address}
    \clearfield{isbn}
}

\usepackage{enumerate,enumitem}
\usepackage{subcaption}
\usepackage{soul}
\usepackage{algorithm, algpseudocode, algorithmicx}
\usepackage{comment}
\usepackage{xspace}
\usepackage{mleftright}

\usepackage[most]{tcolorbox}

\newcommand{\real}{\mathbb{R}}
\newcommand{\complex}{\mathbb{C}}
\newcommand{\field}{\mathbb{K}}
\newcommand{\onevec}{\mathbf{1}}

\DeclareMathOperator{\tr}{tr}
\DeclareMathOperator{\orth}{\mathbf{Orth}}
\DeclareMathOperator{\Orth}{\mathbf{Orth}}
\newcommand{\mat}[1]{\boldsymbol{#1}}
\renewcommand{\vec}[1]{\boldsymbol{#1}}
\newcommand{\lowrank}[1]{\mleft\llbracket #1 \mright\rrbracket}
\newcommand{\norm}[1]{\mleft\| #1 \mright\|}

\DeclareMathOperator{\range}{range}
\newcommand{\QR}{\textsf{QR}\xspace}
\newcommand{\LU}{\textsf{LU}\xspace}
\newcommand{\LAPACK}{\textsf{LAPACK}\index{LAPACK@\textsf{LAPACK}}\xspace}
\newcommand{\RandLAPACK}{\textsf{RandLAPACK}\index{RandLAPACK@\textsf{RandLAPACK}}\xspace}

\newcommand{\BLASthree}{\textsf{BLAS3}\index{BLAS3@\textsf{BLAS}}\xspace}
\newcommand{\RandBLAS}{\textsf{RandBLAS}\index{RandBLAS@\textsf{RandBLAS}}\xspace}
\DeclareMathOperator{\rank}{rank}
\newcommand{\uinorm}[1]{{\left\vert\kern-0.25ex\left\vert\kern-0.25ex\left\vert #1 
		\right\vert\kern-0.25ex\right\vert\kern-0.25ex\right\vert}}
\DeclareMathOperator{\nnz}{nnz}

\newcommand{\expmat}[1]{\begin{bmatrix} #1 \end{bmatrix}}
\newcommand{\twobytwo}[4]{\expmat{#1 & #2 \\ #3 & #4}}
\newcommand{\twobyone}[2]{\expmat{#1 \\ #2}}
\newcommand{\onebytwo}[2]{\expmat{#1 & #2}}
\newcommand{\flatonebytwo}[2]{[\begin{matrix} #1 & #2 \end{matrix}]}

\newcommand{\Id}{\mathbf{I}}
\newcommand{\evec}{\mathbf{e}}

\DeclareMathOperator{\Var}{Var}
\DeclareMathOperator{\Cov}{Cov}
\DeclareMathOperator{\Cor}{Cor}
\DeclareMathOperator{\expect}{\mathbb{E}}
\DeclareMathOperator{\prob}{\mathbb{P}}
\DeclareMathOperator{\Unif}{\textnormal{\textsc{Unif}}}
\DeclareMathOperator{\Normal}{\textnormal{\textsc{Normal}}}
\DeclareMathOperator{\GP}{GP}

\newcommand{\order}{\mathcal{O}}

\DeclareMathOperator{\poly}{poly}
\DeclareMathOperator{\polylog}{polylog}

\newcommand{\ENE}[1]{{\color{blue} [\textbf{ENE}: #1]}\errmessage{This is an ENE comment!}}
\DeclareMathOperator*{\argmax}{argmax}
\DeclareMathOperator*{\argmin}{argmin}
\newcommand{\set}[1]{\mathsf{#1}}
\newcommand{\e}{\mathrm{e}}
\newcommand{\iu}{\mathrm{i}}
\renewcommand{\Re}{\mathrm{Re}}
\renewcommand{\Im}{\mathrm{Im}}
\renewcommand{\hat}[1]{\widehat{#1}}
\renewcommand{\tilde}[1]{\widetilde{#1}}
\newcommand{\actionbox}[1]{\begin{tcolorbox}[colback=white,colframe=black,width=\columnwidth,boxsep=5pt,arc=4pt]
    #1
\end{tcolorbox}}
\newcommand{\importantresult}[1]{%
  \begin{tcolorbox}[
    colback=white,
    colframe=black,
    width=\columnwidth,
    arc=0pt,
    left=5pt,    
    right=5pt,
    top=2pt,     
    bottom=10pt,  
    boxrule=0.5pt,
  ]
    #1
  \end{tcolorbox}%
}
\renewcommand{\d}{\mathrm{d}}

\tcbuselibrary{skins, breakable}

\newtcolorbox{infobox}{
  enhanced,
  breakable,
  colback=gray!5,
  colframe=black,
  fonttitle=\bfseries,
  title={},
  sharp corners=south, 
  boxrule=0.4pt,
  left=6pt,
  right=6pt,
  top=6pt,
  bottom=6pt,
  borderline west={2pt}{0pt}{black},
  before skip=10pt,
  after skip=10pt
}

\newcommand{\infoblurb}[2]{
  \begin{infobox}
    \textbf{\emph{Information collected:}} #1

    \vspace{4pt}
    \textbf{\emph{Computational goal:}} #2
  \end{infobox}
}

\usepackage{listings}
\definecolor{bg}{rgb}{0.97,0.97,0.97}
\definecolor{mlgreen}{rgb}{0.20,0.65,0.28}
\lstdefinelanguage{Matlab}{
    keywords={break, case, catch, continue, else, elseif, end, for, function,
              global, if, otherwise, persistent, return, switch, try, while},
    keywordstyle=\color{blue},
    identifierstyle=\color{black},
    commentstyle=\color{mlgreen},
    stringstyle=\color{magenta},
    backgroundcolor=\color{bg},
    morecomment=[l][\color{mlgreen}]{\%},
    morestring=[b]"
}
\lstdefinelanguage{C}{
    keywords={auto, break, case, char, const, continue, default, do, double, else, enum, extern, float, for, goto, if, inline, int, long, register, restrict, return, short, signed, sizeof, static, struct, switch, typedef, union, unsigned, void, volatile, while},
    keywordstyle=\color{blue},
    identifierstyle=\color{black},
    commentstyle=\color{mlgreen},    
    stringstyle=\color{magenta},
    backgroundcolor=\color{bg},      
    morecomment=[l][\color{mlgreen}]{//},  
    morecomment=[s][\color{mlgreen}]{/*}{*/},  
    morestring=[b]"
}
\crefname{program}{Program}{Programs}
\crefname{conjecture}{Conjecture}{Conjectures}

\newfloat{program}{htbp}{lop}[chapter]
\floatname{program}{Program}

\newcommand{\listofprograms}{\listof{program}{List of Programs}\addcontentsline{toc}{chapter}{List of Programs}}

\ExplSyntaxOn
\NewDocumentCommand{\mytexttt}{m}
 {
  \texttt
   {
    \tl_set:Nn \l_tmpa_tl { #1 }
    \tl_replace_all:Nen \l_tmpa_tl { \char_generate:nn { `_ } { 8 } } { \_ }
    \tl_use:N \l_tmpa_tl
   }
 }
\cs_generate_variant:Nn \tl_replace_all:Nnn { Ne }
\ExplSyntaxOff

\newcommand{\myprogram}[4][t]{%
  \begin{program}[#1]
  \caption[\mytexttt{#4.m}. #2]{\mytexttt{#4.m}. #2 #3}
  \lstinputlisting[
    frame=single,
    label={prog:#4}
  ]{code/#4.m}
  \end{program}
}

\usepackage{amsthm}
\usepackage{thmtools}
\declaretheorem[name=Theorem,numberwithin=chapter]{theorem}
\declaretheorem[name=Conjecture,numberwithin=chapter,sibling=theorem]{conjecture}
\declaretheorem[name=Proposition,sibling=theorem]{proposition}
\declaretheorem[name=Fact,sibling=theorem]{fact}
\declaretheorem[name=Lemma,sibling=theorem]{lemma}
\declaretheorem[name=Corollary,sibling=theorem]{corollary}
\declaretheorem[name=Definition,numberwithin=chapter,style=definition,sibling=theorem]{definition}
\declaretheorem[name=Example,numberwithin=chapter,style=definition,sibling=theorem,qed={$\diamond$}]{example}
\declaretheorem[name=Remark,numberwithin=chapter,style=remark,sibling=theorem,qed={$\diamond$}]{remark}

\crefname{fact}{Fact}{Facts}

\numberwithin{theorem}{chapter}

\usepackage{epigraph}
    \usepackage{etoolbox}
    \makeatletter
    \newlength\epitextskip
    \pretocmd{\@epitext}{\em}{}{}
    \apptocmd{\@epitext}{\em}{}{}
    \patchcmd{\epigraph}{\@epitext{#1}\\}{\@epitext{#1}\\[\epitextskip]}{}{}
    \makeatother

\newcommand\partepigraph[2][60pt]{
\renewcommand{\afterpartskip}{%
\vskip#1
\begin{center}
    \textit{#2}
\end{center}
\vfil}
}

\DefineBibliographyStrings{english}{
  backrefpage={Cited on p.},
  backrefpages={Cited on pp.}
}

\newcommand{\hatbold}[1]{\skew{4}\widehat{\smash{\boldsymbol{#1}}\mathstrut}}
\newcommand{\tildebold}[1]{\skew{4}\widetilde{\smash{\boldsymbol{#1}}\mathstrut}}
\newcommand{\tildevector}[1]{\skew{4}\widetilde{\smash{\boldsymbol{#1}}}}
\newcommand{\hatvector}[1]{\skew{1}\widehat{\smash{\boldsymbol{#1}}}\mskip1mu}
\newcommand{\Bhat}{\smash{\hatbold{B}}}
\newcommand{\Ahat}{\smash{\hatbold{A}}}
\crefname{equation}{}{}
\DeclareMathOperator{\diag}{\mathbf{diag}}
\DeclareMathOperator{\Diag}{\mathbf{Diag}}
\DeclareMathOperator{\srn}{\mathbf{srn}}
\DeclareMathOperator{\scn}{\mathbf{scn}}
\newcommand{\sphere}{\mathbb{S}}
\newcommand{\unitcircle}{\mathbb{T}}
\newcommand{\myparagraph}[1]{\textbf{\textit{#1.}}}
\newcommand{\myparagraphnp}[1]{\textbf{\textit{#1}}}
\newcommand{\mysubparagraph}[1]{\textit{#1.}}
\DeclareMathOperator{\diagprod}{\mathbf{diagprod}}
\newcommand{\outprod}[1]{#1^{\vphantom{*}}#1^*}

\definecolor{ltyellow}{rgb}{1, 1, 0.9}
\usepackage{soul}
\sethlcolor{ltyellow}
\newcommand{\warn}[1]{\emph{\hl{#1}}}

\cftsetindents{part}{0em}{2em}
\cftsetindents{section}{1em}{3em}
\makeatletter
\renewcommand{\@pnumwidth}{3em}  
\makeatother

\newcommand{\XTrace}{\textsc{XTrace}\xspace}
\newcommand{\XNysTrace}{\textsc{XNysTrace}\xspace}
\newcommand{\XSymTrace}{\textsc{XSymTrace}\xspace}
\newcommand{\XDiag}{\textsc{XDiag}\xspace}
\newcommand{\XNysDiag}{\textsc{XNysDiag}\xspace}
\newcommand{\XRowNorm}{\textsc{XRowNorm}\xspace}
\newcommand{\XSymRowNorm}{\textsc{XSymRowNorm}\xspace}

\newcommand{\HutchPP}{\textsc{Hutch}\textnormal{\texttt{++}}\xspace}
\newcommand{\DiagPP}{\textsc{Diag}\textnormal{\texttt{++}}\xspace}
\newcommand{\hutchppest}{\hat{\tr}_{\mathrm{H}\textnormal{\texttt{++}}}}
\newcommand{\udiagppest}{\hat{\diag}_{\mathrm{UD}\textnormal{\texttt{++}}}}
\newcommand{\diagppest}{\hat{\diag}_{\mathrm{D}\textnormal{\texttt{++}}}}
\newcommand{\MSE}{\mathrm{MSE}}
\newcommand{\Err}{\mathrm{Err}}
\newcommand{\RowNormName}{Johnson--Lindenstrauss\xspace}
\newcommand{\QiQi}{\mat{Q}_{(i)}^{\vphantom{*}}\mat{Q}_{(i)}^*}
\newcommand{\Bias}{\mathrm{Bias}}
\DeclareMathOperator{\vecop}{\mathbf{vec}}

\newcommand{\RPCholesky}{\textsc{RPCholesky}\index{randomly pivoted Cholesky}\xspace}
\newcommand{\RejectionSampleSubmatrix}{\textsc{RejectionSampleSubmatrix}\index{RejectionSampleSubmatrix@\textsc{RejectionSampleSubmatrix}}\xspace}
\newcommand{\RPCholeskyRejection}{\textsc{RejectionRPCholesky}\index{RejectionRPCholesky@\textsc{RejectionRPCholesky}}\xspace}
\newcommand{\RPQRRejection}{\textsc{RejectionRPQR}\index{RejectionRPQR@\textsc{RejectionRPQR}}\xspace}
\newcommand{\RPQR}{RPQR\index{randomly pivoted QR@randomly pivoted \QR}\xspace}
\newcommand{\CURnoindex}{\textsf{CUR}\xspace}
\newcommand{\CURindex}[1]{\index{CUR approximation@\textsf{CUR} approximation#1}}
\newcommand{\CUR}{\textsf{CUR}\index{CUR approximation@\textsf{CUR} approximation}\xspace}
\newcommand{\RPCURtwo}{\textsc{RPCUR2}\index{RPCUR2 algorithm@\textsc{RPCUR2} algorithm}\xspace}
\newcommand{\RPCURLev}{\textsc{RPCURLev}\index{RPCURLev algorithm@\textsc{RPCURLev} algorithm}\xspace}
\newcommand{\MDCUR}{\textsc{MDCUR}\index{MDCUR algorithm@\textsc{MDCUR} algorithm}\xspace}
\newcommand{\RKHS}{\set{H}}
\newcommand{\Ltwo}{\set{L}^2}

\newcommand{\Hilb}{\set{K}}
\DeclareMathOperator{\cond}{cond}
\newcommand{\TransferenceOfAlgorithms}{\actionbox{\textbf{Gram correspondence: Transference of algorithms.} Every algorithm producing a projection approximation to a general matrix has an analogous algorithm that produces a Nystr\"om approximation to a psd matrix and vice versa.}}
\newcommand{\NP}{\textsf{NP}\xspace}
\newcommand{\kDPP}[1]{\operatorname{DPP}_{#1}}

\newcommand{\simiid}{\stackrel{\text{iid}}{\sim}}
\newcommand{\SkLUPP}{\textsc{SkLUPP}\index{SkLUPP algorithm@\textsc{SkLUPP} algorithm}\xspace}

\DeclareMathOperator{\SD}{SD}
\DeclareMathOperator{\Jack}{Jack}
\DeclareMathOperator*{\minimize}{minimize}

\newif\iffull
\fullfalse

\DeclareMathOperator{\fl}{fl}
\DeclareMathOperator{\err}{err}
\DeclareMathOperator{\BE}{BE}

\makeatletter
\def\@@showbox#1{%
  \color{red}%
  \@tempdima\dp#1\unvbox#1\dp#1\@tempdima}
\let\showbox\@@showbox
\makeatother

\usepackage{xurl}    

\begin{document}

\title{Make the Most of What You Have: \\ Resource-Efficient Randomized Algorithms for \\Matrix Computations}
\author{Ethan N. Epperly}

\degreeaward{Doctor of Philosophy

\vspace{1em}

\textbf{arXiv Version \\ Lightly Edited from Official Version \\ Date: December 17, 2025}}                 
\university{California Institute of Technology}    
\address{Pasadena, California}                     
\unilogo{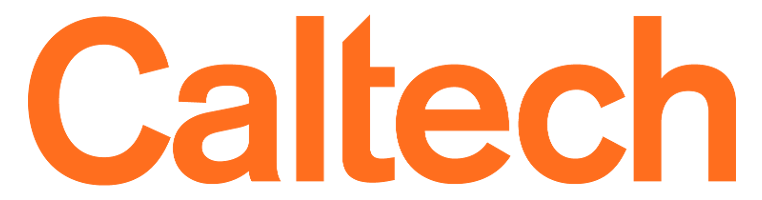}                                 
\copyyear{2025}  
\defenddate{May 1, 2025}          

\orcid{0000-0003-0712-8296}

\rightsstatement{All rights reserved}

\maketitle[logo]

\begin{acknowledgements}
    This PhD thesis is the product of the kindness and mentorship provided to me by many people over many years.
    There are many to thank (and, I'm sure, some who I have regretfully forgotten to acknowledge).

    First, I must thank my advisor Joel A.\ Tropp, from whom I have learned so much---about mathematics, about research, and about writing.
    I thank Joel for always pushing me to dream bigger, setting my sights on more challenging problems.
    Thanks also for your detailed feedback on this thesis.

    My thanks extends to my other mentors throughout my journey into academic research: Franceso Paesani and Jordi Cirera at UCSD; Don Ward, Bernice Mills, Ryan Sills, and Jonathan Hu at Sandia National Laboratories; Shivkumar Chandrasekaran at UCSB; Andrew Barker at Lawrence Livermore National Laboratory; and Lin Lin at Lawrence Berkeley National Laboratory.
    I also give thanks to Alf Morales and Tim Shepodd for giving a 15 year-old kid an early shot at pursuing his dream to be a scientist.

    I thank the remaining members of my committee, Profs.\ Venkat Chandrasekaran, Lin Lin, and Franca Hoffmann, for you time in sitting on my dissertation committee and for all I have gotten to learn from you during my PhD.

    I feel truly blessed to belong to research community filled with such kind, humble, and thoughtful people.
    Getting to collaborate with others is truly my favorite part of research.
    I extend my thanks to the many friends, collaborators, and mentors I have done research with, including Alex Barnett, Chris Cama\~no, Shivkumar Chandrasekaran, Yifan Chen, Lieven De Lathauwer, Mateo D\'iaz, Zhiyan Ding, Gil Goldshlager, Zachary Frangella, Nithin Govindarajan, Anne Greenbaum, Alex Hsu, Daniel Kressner, Lin Lin, Maike Meier, Eliza O'Reilly, Kevin Miller, Raphael Meyer, Elvira Moreno, Chris Musco, Yuji Nakatsukasa, Taejun Park, Akash Rao, Heather Wilber, and Ruizhe Zhang.
    I also thank Deeksha Adil, Haoxuan Chen, Tyler Chen, Alice Cortinovis, Jorge Garza-Vargas, Eitan Levin, and Aleksandros Sobczyk for many insightful conversations.
    The research process is inherently collaborative, and this thesis is suffused with the creativity of these people.

    Special thanks are merited to Rob Webber, in particular, who joined Caltech as a postdoctoral scholar during the second year of my PhD.
    I truly cannot imagine what my PhD would have been without your collaboration, your mentorship, and your friendship.

    I am extremely grateful to Chris Musco for his invitation to visit NYU in the spring of 2024.
    It was a very enjoyable and productive visit, and I learned a lot from collaborating with you.

    My warmest thanks go to administrative staff at Caltech, most especially Jolene Brink and Maria Lopez.
    Thanks to all the help you have provided to me and other students in our department.

    My PhD experience truly would not have been what it's been without the Department of Energy Computational Science Graduate Fellowship and the community of people I have gotten to meet and interact with through it.
    I send my thanks to Lindsey Eilts and all the folks at the Krell Institute who manage the fellowship and provide so much support to us fellows.

    I give my deepest gratitude to my family.
    To my parents Meg and Tom, who have always loved and supported me throughout my life and nurtured my love of learning.
    All of that time playing board games, those long conversations about what I was learning in school, and those \emph{Science on Saturday} talks left a big impression on me.
    To my brother Aidan, it is truly a gift to have a close friend who loves talking math just as much as I do (possibly even more).
    To my grandparents Kathie, Bill, Christine, John, Bob, and Sarah, for all the wisdom and life experience you have shared.

    Doing a PhD is a long journey, and I truly believe I could not have made it without the love and support of my fianc\'ee, Sierra Williams, DVM.
    The level of dedication and intensity you bring to everything you do is an inspiration to me, and you have always been there for me when I needed it.
    This love and thanks extends to our pets: dogs Hulk and Finn and three-legged turtle Shelly.
    Nothing brought me joy throughout my PhD quite like going on a walk with the dogs or sitting down all five of us to watch an episode of Survivor.

    \vfill

    Thanks to the many groups who have funded my PhD research: namely, the U.S. Department of Energy, Office of Science, Office of Advanced Scientific Computing Research, Department of Energy Computational Science Graduate Fellowship under Award Number DE-SC0021110, and under aegis of Joel Tropp, Office of Naval Research BRC awards N00014-18-1-2363 and N00014-24-1-2223, NSF FRG 1952777, and the Caltech Carver Mead New Adventures Fund.
\end{acknowledgements}

\begin{abstract}
   In recent years, randomized algorithms have established themselves as fundamental tools in computational linear algebra, with applications in scientific computing, machine learning, and quantum information science.
   Many randomized matrix algorithms proceed by first collecting information about a matrix and then processing that data to perform some computational task.
   This thesis addresses the following question: How can one design algorithms that use this information as efficiently as possible, reliably achieving the greatest possible speed and accuracy for a limited data budget?
   This question is timely, as randomized algorithms are increasingly being deployed in production software and in applications where accuracy and reliability is critical.

   The first part of this thesis focuses on the problem of low-rank approximation for positive-semidefinite matrices, motivated by applications to accelerating kernel and Gaussian process machine learning methods.
   Here, the goal is to compute an accurate approximation to a matrix after accessing as few entries of the matrix as possible.
   This part of the thesis explores the randomly pivoted Cholesky (\RPCholesky) algorithm for this task, which achieves a level of speed and reliability greater than other methods for the same problem.

   The second part of this thesis considers the task of estimating attributes of an implicit matrix accessible only by matrix--vector products, motivated by applications in quantum physics, network science, and machine learning.
   This thesis describes the \emph{leave-one-out approach} to developing matrix attribute estimation algorithms, and develops optimized trace, diagonal, and row-norm estimation algorithms for this computational model.

   The third part of this thesis considers randomized algorithms for overdetermined linear least squares problems, which arise in statistics and machine learning.
   Randomized algorithms for linear-least squares problems are asymptotically faster than any known deterministic algorithm, but recent work of [Meier et al., SIMAX `24] raised questions about the accuracy of these methods when implemented in floating point arithmetic.
   This thesis shows these issues are resolvable by developing fast randomized least-squares problem achieving backward stability, the gold-standard accuracy and stability guarantee for a numerical algorithm.
\end{abstract}

\begin{publishedcontent}[iknowwhattodo]
\nocite{CETW25w,EM23aw,ETW24aw,ETW24w,ET24w,Epp24aw,EMN24w,DEF+23w}
\end{publishedcontent}

{\sloppy
\tableofcontents
\listoffigures
\listoftables
\listofprograms
}

\mainmatter

\chapter{Introduction}

\epigraph{Our experience suggests that many practitioners of scientific computing view randomized algorithms as a desperate and final resort. Let us address this concern immediately. Classical Monte Carlo methods are highly sensitive to the random number generator and typically produce output with low and uncertain accuracy. In contrast, the algorithms discussed herein are relatively insensitive to the quality of randomness and produce highly accurate results.}{Nathan Halko, Per-Gunnar Martinsson, and Joel Tropp, \emph{Finding structure with randomness: Probabilistic algorithms for constructing matrix decompositions} \cite{HMT11}}

Randomness has played a fundamental role in digital computation throughout its history, dating back to von Neumann and Ulam's groundbreaking work in the 1940s using Monte Carlo to study neutron diffusion problem.
(See \cite{Eck87} for discussion of this history.)
Despite this, randomized methods have also been traditionally viewed with great skepticism by many practitioners of scientific computing.
The belief was that randomized methods could only yield low-accuracy solutions to computational problems.
As such, they were only useful as a ``desperate and final resort'', typically for problems of such high dimensions that traditional methods were doomed to failure.

In modern scientific computing, this conventional wisdom has been completely upended. 
Some of the most exciting progress has been for problems in computational linear algebra, the subject of this thesis.
Randomized algorithms have established them as the most effective methods for linear algebra problems such as matrix low-rank approximation \cite{HMT11,MM15,Nak20,TW23} and highly overdetermined linear least-squares problems \cite{RT08,AMT10}.
Randomized algorithms can produce solutions of high-accuracy (as accurate as the numerical precision, in some case) and with negligible probability of failure.

The modern field of randomized matrix computations has its origins around the turn of the millennium \cite{FKV98,PTRV98}, and its historical roots stretch into the twentieth century \cite{Gir87,KW92}.
As such, this is not a new subject, and many surveys exist to document the field's successes \cite{HMT11,Woo14a,DM16,KV17a,MT20a,TW23,MDM+23a,KT24}.
Despite this multi-decade history, this thesis and its aims remain timely.
We are in the midst of a migration of the randomized matrix algorithm from its natural habitat in the SIAM journal or computer science theory conference to the realm of production software.
Now, more than ever, we need randomized algorithms that are ready for deployment in practice, making maximally efficient use of computational resources and achieving the same level of accuracy and reliability as classical deterministic methods, while maintaining speed.
The goal of this thesis is to describe algorithms developed over the course of my PhD research that meet these criteria.

Many randomized algorithms work by first \emph{collecting information} about a matrix, then using this information to produce an (approximate) solution to a linear algebra problem.
A unifying theme of this thesis is the design of algorithms that \emph{make the most of what you have}---that is, algorithms that attempt to collect as little information as possible about the matrix and use the collected information in a maximally efficient way to achieve the greatest possible accuracy subject to a limited computational budget.
This thesis is divided into three parts, each of which uses the ``make the most of what you have'' principle in a different way.

\section{Part I: Random pivoting}

\vspace{0.5em}

\infoblurb{Entries of a positive-semidefinite matrix $\mat{A}$.}{Produce a \emph{low-rank} approximation $\Ahat\approx \mat{A}$.}

Kernel and Gaussian process methods in machine learning \cite{SS02,RW05} remain some of the most effective tools for scientific machine learning \cite{MCRR20,BDHO24}.\index{kernel method}\index{Gaussian process!type of machine learning method}
However, these methods face a fundamental limitation called the curse of kernelization: \index{curse of kernelization} \index{kernel method!computational difficulties}
\actionbox{\textbf{Curse of kernelization.} Direct implementations of kernel and Gaussian process methods on $n$ data points requires forming and manipulating an $n\times n$ positive-semidefinite kernel matrix $\mat{A}$.\index{kernel matrix}
Doing so requires $\order(n^2)$ storage and $\order(n^3)$ operations.}
The curse of kernelization makes direct implementation of kernel and Gaussian process methods infeasible for large data sets.

Randomization offers a path forward: Using randomized algorithms, we can compute a \emph{low-rank approximation} $\Ahat \approx \mat{A}$ to the kernel matrix.\index{kernel matrix}\index{low-rank approximation!use cases}
By using the low-rank approximation in place of the kernel matrix, the computational cost of kernel and Gaussian process methods can be substantially reduced. 

In order to maintain computational efficiency, kernel matrix low-rank approximation algorithms must produce an approximation to the matrix $\mat{A}$ after only accessing a small number of entries of the matrix; after all, reading the full matrix even a single time incurs an expensive cost of $\order(n^2)$ operations.
Is it possible to produce an accurate approximation to a matrix after looking at a fraction of its entries?
And if so, what is most economical algorithm?

This thesis advocates for the \emph{randomly pivoted Cholesky} algorithm as an answer to these questions. \index{randomly pivoted Cholesky}
The randomly pivoted Cholesky\index{randomly pivoted Cholesky} algorithm produces a \emph{near-optimal} rank-$k$ approximation to any positive-semidefinite matrix using just $(k+1)n$ entry accesses.
Compared to other algorithms for kernel matrix low-rank approximation, randomly pivoted Cholesky\index{randomly pivoted Cholesky} is either faster or more accurate when applied to challenging instances.

\Cref{part:random-pivoting} of this thesis introduces the randomly pivoted Cholesky\index{randomly pivoted Cholesky} algorithm, describes its properties, explains how it can be used to accelerate kernel and Gaussian process machine learning algorithms, and compares it to alternative methods.
After a thorough investigation of randomly pivoted Cholesky\index{randomly pivoted Cholesky} and positive-semidefinite matrix low-rank approximation, this part also discuss extensions of the random pivoting approach to approximating general, rectangular matrices.


\section{Part II: Leave-one-out randomized matrix algorithms}

\vspace{0.5em}

\infoblurb{Matrix--vector products $\mat{B}\vec{\omega}_1,\ldots,\mat{B}\vec{\omega}_s$.}{Estimate attributes associated to the matrix $\mat{B}$ such as its trace, diagonal, or row norms.}

In fields ranging from quantum physics to network science to machine learning, we work with matrices $\mat{B}$ that are accessible only indirectly.
In such settings, we cannot read the entries $b_{ij}$ directly.
Rather, we can access the matrix through matrix--vector products: Given a vector $\vec{\omega}$, we may access $\mat{B}\vec{\omega}$.\index{matvec model}\index{implicit matrix}

Under this implicit access model, linear algebra problems that would ordinarily be trivial become interesting.
For instance, given access to an implicit matrix $\mat{B}$, how can we compute its trace?
Its diagonal?
Its row norms?
For all of these questions, we seek algorithms that are maximally efficient, in the sense of requiring as few matrix--vector products to compute these \emph{matrix attributes} to a desired accuracy level.

\Cref{part:loo} describes the \emph{leave-one-out\index{leave-one-out randomized algorithm} approach} to estimating attributes of an implicit matrix.
We will then use the leave-one-out\index{leave-one-out randomized algorithm} approach to derive optimized estimators for the trace, diagonal, and row norms of both general and positive-semidefinite matrices.
As additional applications of the leave-one-out\index{leave-one-out randomized algorithm} approach, we will develop error estimates for randomized low-rank approximations like the randomized SVD and variance estimation techniques for assessing the quality of more general quantities computed by randomized matrix algorithms.

\section{Part III: Sketching, solvers, and stability}

\vspace{0.5em}

\infoblurb{A ``sketch'' of a tall matrix $\mat{B}\in\real^{m\times n}$ with $m\gg n$, defined to be the product $\mat{S}^*\mat{B}$ of $\mat{B}$ with a wide matrix $\mat{S}^* \in \real^{d\times m}$.}{Solve the overdetermined linear least-squares problem $\vec{x} = \argmin_{\vec{z} \in \real^n} \norm{\vec{c} - \mat{B}\vec{z}}$ to high numerical accuracy.}

One of the major success stories of randomized algorithms in linear algebra is sketch-and-precondition algorithm \cite{RT08,AMT10} for solving over-determined linear least-squares problems \index{least squares} \index{sketch-and-precondition|(}
\begin{equation} \label{eq:ls-intro}
    \vec{x} = \argmin_{\vec{z} \in \real^n} \norm{\vec{c} - \mat{B}\vec{z}} \quad \text{where } \mat{B} \in \real^{m\times n}, \: \vec{c} \in \real^m.
\end{equation}
Here, $\norm{\cdot}$ denotes the vector $\ell_2$ norm, and we consider the \emph{highly overdetermined setting} where $m\gg n$.
The idea of sketch-and-precondition is to preprocess the matrix $\mat{B}$ by applying a randomized dimensionality reduction map or ``sketching matrix'' $\mat{S}^* \in \real^{d\times m}$, where the \emph{embedding dimension} $d$ is a small multiple of the number of columns $n$ in $\mat{B}$, e.g., $d=2n$. \index{sketching matrix}
We then compute a \QR decomposition of the sketched matrix $\mat{S}^*\mat{B} = \mat{Q}\mat{R}$ and solve the least-squares problem \cref{eq:ls-intro} using a Krylov iterative method with $\mat{R}$ as preconditioner.
With an appropriate choice of sketching matrix, this procedure reliably reduces the condition number of $\mat{B}$ below an absolute constant, leading to convergence to high accuracy in a small number of Krylov iterations.
As a result, sketch-and-precondition runs in roughly $\order(mn + n^3)$ operations, a dramatic improvement to the $\order(mn^2)$ cost of classical algorithms for least squares.

The recent paper \cite{MNTW24} cast doubts on the \emph{numerical stability} of this algorithm. \index{sketch-and-precondition!numerical stability}
This paper's experiments show that, when implemented using ordinary double-precision floating point arithmetic, some versions of sketch-and-precondition produce errors that exceed classical methods by orders of magnitude.
Does this mean the sketch-and-precondition algorithm is too unreliable for practical use?

\Cref{part:sketching} critically investigates this question.
We will see that, by using careful initialization and just a single step of \emph{iterative refinement}, sketch-and-precondition-type algorithms can be made just as accurate as classical direct methods, achieving \emph{backward stability}, the gold-standard notion of accuracy in numerical analysis. \index{backward stability} \index{numerical stability!backward stability}

\section{About this thesis}

I have taken non-traditional approach to writing this thesis.
As with many theses, this thesis is primarily based on research papers written during my PhD.
However, rather then collecting these papers and editing them, this thesis contains a new treatment of this material.
As such, this thesis is composed of original writing and constitutes a deeper dive and re-examination of my PhD research.

Rather than trying to include all of my PhD research, I have elected to focus on the three topics described above: random pivoting, leave-one-out randomized matrix algorithms, and numerically stable randomized least-squares algorithms.
My goal was to provide an approachable introduction to these three areas of my research that should be accessible to younger researchers.
As such, I have included a significant amount of discussion of background material and related work, and the end of each part of the thesis describes open problems.
This thesis also contains several new results and algorithms that have not previously been published.

\begin{figure}
    \centering
    \includegraphics[width=0.7\linewidth]{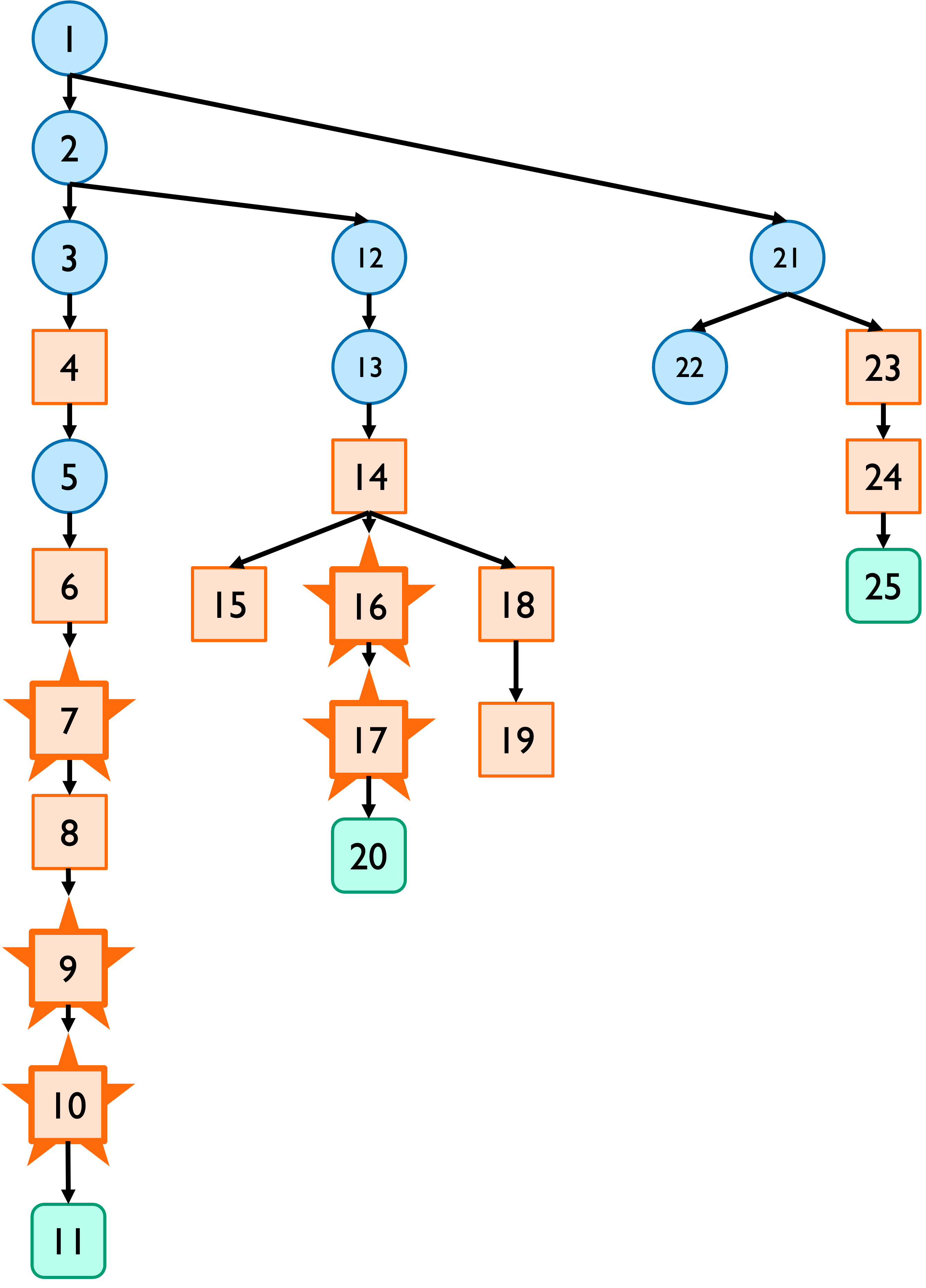}
    \caption[List of chapters in this thesis and dependencies between them]{List of chapters in this thesis and dependencies between them.
    Blue circles indicate chapters whose content is primarily introductory or expository, orange squares indicate sections primarily containing research from my PhD, and green squircles indicate open questions.
    Starred sections contain new research that has not previously been published.}
    \label{fig:roadmap}
\end{figure}

The contents of this thesis are diagrammed in \cref{fig:roadmap}, including dependencies between chapters.
For readers interested primarily in seeing my original research, I have used symbols to distinguish between chapters based on whether they contain primarily exposition, research, or open problems.
I have endeavored to have all of the chapters, including the expository ones, present a unique perspective on the material.
Chapters containing significant research content that is newly presented in this thesis are marked using a star.

Given my focus on a subset of topics, much of my PhD research is not discussed in this thesis, including my work on quantum eigenvalue algorithms \cite{ELN22,DELZ24a}, row action methods for large least-squares problems \cite{EGW25}, superfast non-uniform Fourier transform inversion \cite{WEB25}, tensor network methods \cite{CET25}, and uniqueness of tensor decompositions \cite{GEL22}.
All of these topics are near and dear to my heart, and it is will sadness that I do not include them in the present document.

Taking inspiration from Nick Trefethen's essay \emph{Ten Digit Algorithms} \cite[\S4]{Tre05}, I have elected to present algorithms using code segments rather than pseudocode. 
As a numerical linear algebraist at heart, the programming language will be MATLAB; translation to Python or Julia for those so-inclined will hopefully not be difficult.
These code segments and code files to reproduce most of the numerical experiments in this thesis can be found at
\actionbox{\url{https://github.com/eepperly/Ethan-Epperly-Thesis}}

The version of this thesis that you are reading now was lightly edited from the version that was submitted to the Caltech library, largely to correct typos or other small writing problems.

\section{Notation}

This thesis will work over the field $\field$ of either real or complex numbers, $\field = \real$ or $\field = \complex$.
The imaginary unit is denoted $\iu$, and the complex conjugate is $\overline{x + \iu y} = x - \iu y$.
The nonnegative reals are denoted $\real_+$.

\myparagraph{Matrices and vectors}
The symbol ${}^*$ denotes the adjoint, which reduces to the transpose when $\field = \real$.
We identify the set of vectors $\field^n$ with the set of $n\times 1$ matrices $\field^{n\times 1}$; row vectors $\vec{b}^* \in \field^{1\times n}$ will be denoted as adjoints of column vectors $\vec{b} \in\field^n$.

Matrices will be denoted by bold capital letters (e.g., $\mat{B},\mat{\Omega}$), and vectors will be denoted by bold lowercase letters (e.g., $\vec{b}$, $\vec{\omega}$).
Entries of matrices and vectors are denoted using the corresponding italic letter, e.g., $b_i$ is the $i$th entry of $\vec{b}$ and $\omega_{ij}$ is the $ij$th entry of $\mat{\Omega}$.
Columns of matrices are denoted using the corresponding bold lowercase vector, e.g., $\vec{\omega}_i$ is the $i$th column of $\mat{\Omega}$.
A matrix $\mat{\Omega}$ with its $i$th column deleted is denoted $\mat{\Omega}_{-i}$ (similarly, $\mat{\Omega}_{-ij}$ for two deletions).

The identity matrix is $\Id\in\real^{n\times n}$; its columns are the standard basis vectors $\evec_i \in \real^n$.
The vector of all ones is $\onevec = \sum_{i=1}^n \evec_i$, and the matrix of zeros is $\mat{0} \in \real^{m\times n}$.
We will use subscripts $\Id_n,\onevec_n,\mat{0}_{m\times n}$ when necessary to indicate the dimensions.

For our purposes, a matrix $\mat{A} \in \field^{n\times n}$ is \emph{positive semidefinite} (psd) if it is Hermitian ($\mat{A} = \mat{A}^*$) and satisfies the condition $\vec{x}^*\mat{A}\vec{x} \ge 0$ for all vectors $\vec{x} \in \field^n$.
The concept of a ``real nonsymmetric psd matrix'', valuable in other contexts, will not be considered in this thesis.
The psd order $\succeq$ is defined so that $\mat{A} \succeq \mat{H}$ whenever $\mat{A} - \mat{H}$ is psd.

Throughout this thesis, $\mat{B}\in\field^{m\times n}$ will denote a general matrix, which may be square or rectangular depending on the context.
The matrix $\mat{A}\in\field^{n\times n}$ will always denote a Hermitian matrix, often psd. \index{positive-semidefinite (psd) matrix}

Sets of integer indices, general arbitrary sets, and vector spaces will be denoted using sans serif font (e.g., $\set{S}$, $\set{V}$).
General inner products, when they arise, will be conjugate linear in the \warn{first} coordinate.

We permit matrices $\mat{B} \in \field^{\set{S}\times \set{T}}$ and vectors $\vec{x} \in \field^{\set{S}}$ to be indexed by arbitrary finite sets $\set{S}$ and $\set{T}$.
Given bivariate and univariate functions $\xi : \set{X} \times \set{Y} \to \field$ and $f : \set{X} \to \field$ and subsets $\set{S}\subseteq \set{X}$ and $\set{T} \subseteq \set{Y}$, $\xi(\set{S},\set{T}) \in \field^{\set{S} \times \set{T}}$ and $f(\set{S}) \in \field^{\set{S}}$ denote a matrix and vector with elements $\xi(s,t)$ and $f(s)$ for each $s\in\set{S}$ and $t\in\set{T}$.

We use MATLAB notation to index matrices.
So, e.g.,
\begin{itemize}
    \item $\mat{B}(i,j) = b_{ij}$ is the $ij$th entry of $\mat{B}$,
    \item $\mat{B}(:,i) =\vec{b}_i$ is the $i$th column of $\mat{B}$, and
    \item $\mat{B}(:, \set{S}) = (\vec{b}_i : i \in \set{S})$ is the submatrix of $\mat{B}$ indexed by the set $\set{S}$.
\end{itemize} 
We adopt the convention that submatrices $\mat{B}(\set{S},\set{T}) \in \field^{\set{S}\times \set{T}}$ are indexed by the sets $\set{S}$ and $\set{T}$ rather than by the sets $\{1,\ldots,|\set{S}|\}$ and $\{1,\ldots,|\set{T}|\}$.
So $b_{22}$ is still the $(2,2)$ entry of the matrix $\mat{B}(\{2,3\},\{2,5\})$.

For a square matrix $\mat{B}\in\field^{n\times n}$, the trace is $\tr(\mat{B}) = \sum_{i=1}^n b_{ii} \in \field$ and the diagonal is $\diag(\mat{B}) = (b_{ii} : 1\le i\le n) \in \field^n$.
The diagonal elements of a matrix product $\mat{F}^*\mat{G}$ are denoted  \index{diagprod@$\diagprod$}
\begin{equation*}
    \diagprod(\mat{F},\mat{G}) \coloneqq \diag(\mat{F}^*\mat{G}).
\end{equation*}
Given a vector $\vec{b}\in\field^n$, the diagonal matrix with elements $b_i$ is denoted $\Diag(\vec{b})$.
The squared row norms and squared column norms of a matrix $\mat{B}\in\field^{m\times n}$ are \index{squared row or column norms}
\begin{align*}
    \srn(\mat{B}) &\coloneqq \bigl(\norm{\mat{B}(i,:)}^2 : 1\le i \le m\bigr)& &= \diag(\mat{B}\mat{B}^*) &&\in \field^m, \\
    \scn(\mat{B}) &\coloneqq \bigl(\norm{\mat{B}(:,i)}^2 : 1\le i \le n\bigr)& &= \diag(\mat{B}^*\mat{B}) &&\in \field^n.
\end{align*}

The entrywise product of vectors $\vec{a},\vec{b} \in \field^n$ is $\vec{a} \odot \vec{b} = (a_i \cdot b_i : 1\le i \le n)$.
Nonlinear operations are applied to vectors entrywise.
For instance, $\overline{\vec{a}} = (\overline{a_i} : 1\le i\le n)$ is the entrywise complex conjugate, and $|\vec{a}|^2 = (|a_i|^2 : 1\le i\le n)$ is the entrywise squared modulus.

The $\ell_2$ norm of a vector or spectral norm of a matrix will be denoted $\norm{\cdot}$.
We will also make use of the Frobenius and trace norms, denoted $\norm{\cdot}_{\mathrm{F}}$ and $\norm{\cdot}_*$.
Schatten $p$-norms are denoted $\norm{\cdot}_{\set{S}_p}$.

We will make frequent use of matrix decompositions.
Given a tall matrix $\mat{B} \in \field^{m\times n}$, an economy-size \QR decomposition is a factorization $\mat{B} = \mat{Q}\mat{R}$ into a matrix $\mat{Q} \in \field^{m\times n}$ with orthonormal columns and an upper triangular matrix $\mat{R} \in \field^{n\times n}$. \index{QR decomposition@\QR decomposition!definition}
An economy-size SVD is a factorization $\mat{B} = \mat{U}\mat{\Sigma}\mat{V}^*$ into a matrix $\mat{U} \in \field^{m\times n}$ with orthonormal columns, a nonnegative diagonal matrix $\mat{\Sigma} \in \real_+^{n\times n}$, and a unitary matrix $\mat{V} \in \field^{n\times n}$.
When not otherwise stated, all \QR decompositions and SVDs will be assumed to be economy-size.
Given a positive definite matrix $\mat{A}$, its (full) Cholesky decomposition refers to one of two matrix decompositions, the \emph{lower triangular Cholesky decomposition} $\mat{A} = \mat{L}\mat{L}^*$ or the \emph{upper triangular Cholesky decomposition} $\mat{A} = \mat{R}^*\mat{R}$. \index{Cholesky decomposition!definition}
Both conventions will be convenient for use in different places of this thesis.

We will use $\orth(\mat{B})$ to denote a matrix whose columns form an orthonormal basis for the subspace $\range(\mat{B})$.
When necessary to yield an unambiguous interpretation, $\orth$ of a full column-rank matrix $\mat{B}$ will be given by an economy-size \QR decomposition where the triangular factor has positive diagonal entries.

As usual, $\rank(\mat{B})$ denotes the dimension of the range of $\mat{B}$.
We say a matrix is \emph{rank-$k$} if $\rank(\mat{B}) \le k$.
A rank-$k$ matrix $\mat{B} \in \field^{m\times n}$ is described by its \emph{compact SVD}
\begin{equation*}
    \mat{B} = \mat{U} \mat{\Sigma} \mat{V}^*
\end{equation*}
where $\mat{U} \in \field^{m\times \rank(\mat{B})}$ and $\mat{V} \in \field^{n\times \rank(\mat{B})}$ have orthonormal columns and $\mat{\Sigma} \in \real_+^{\rank(\mat{B})\times \rank(\mat{B})}$ is a diagonal matrix listing the \emph{first $k$} singular values.
Similarly, a rank-$k$ Hermitian matrix $\mat{A} \in \field^{n\times n}$ has a compact eigendecomposition
\begin{equation*}
    \mat{A} = \mat{U}\mat{D}\mat{U}^*
\end{equation*}
where $\mat{U} \in \field^{n\times \rank(\mat{A})}$ has orthonormal columns and $\mat{D} \in \real^{\rank(\mat{A})\times \rank(\mat{A})}$ is a diagonal matrix listing the first $k$ eigenvalues.

\myparagraph{Probability theory}
Probabilities and expectations are denoted $\prob$ and $\expect$.
Nonlinear operations bind before expectations, so that $\expect x^2 \coloneqq \expect(x^2)$, \emph{not} $\expect x^2 = (\expect x)^2$.
The covariance is 
\begin{equation*}
    \Cov(a,b) \coloneqq \expect[(a - \expect a)\overline{(b - \expect b)}],
\end{equation*}
and the variance is $\Var(a) \coloneqq \Var(a,a)$.
Note that, unlike other sesquilinear forms in this thesis, the covariance is conjugate-linear in its \warn{second} argument.

We write $x_1,x_2,\ldots\sim x$ and variations thereof to indicate that $x_1,x_2,\ldots$ are random variables with the same distribution as $x$.
The phrase ``independent and identically distributed'' carries its usual abbreviation \emph{iid}, and we write $x_1,x_2,\ldots\stackrel{\text{iid}}{\sim} x$ when $x_1,x_2,\ldots$ are iid copies of $x$.

Given a nonnegative \emph{weight vector} $\vec{w} \in \real_+^n$, $i \sim \vec{w}$ denotes a random integer $i \in \{1,\ldots,n\}$ selected with probability $\prob \{ i = j \} = w_j / \sum_{k=1}^n w_k$.
We \warn{do not assume the normalization $\sum_{i=1}^n w_i = 1$} when writing $i \sim \vec{w}$.

The uniform distributions on a set $\set{S}$ are denoted $\Unif \set{S}$.
The unit sphere of a vector (sub)space $\set{U}$ is $\sphere(\set{U})$, and $\Unif \sphere(\set{U})$ denotes the (Haar) uniform distribution on the sphere.
A Gaussian distribution over the field $\field$ with mean $\vec{m}$ and covariance matrix\index{covariance matrix} $\mat{\Sigma}$ is denoted $\Normal_\field(\vec{m},\mat{\Sigma})$. \index{Gaussian distribution}
If $\mat{\Sigma}$ is nonsingular, its probability density is
\begin{equation*}
    p(\vec{z}) = \frac{1}{(2\pi)^{\alpha/2} \det(\mat{\Sigma})^{1/2}} \exp \left( - \frac{\vec{z}^*\mat{\Sigma}^{-1}\vec{z}}{2}\right).
\end{equation*}
A random vector $\vec{\omega} \sim \Normal_\field(\vec{m},\mat{\Sigma})$ has jointly Gaussian entries with mean $\expect[\omega_i] = m_i$ and $\Cov(\omega_i,\omega_j) = \sigma_{ij}$.

\myparagraph{Rounding error analysis}
In \cref{part:sketching} of this thesis, we will investigate the \emph{numerical stability} of randomized least-squares solvers, which will require special notation.
This notation will be introduced in \cref{sec:stability-notation}.

\chapter{Low-rank approximation foundations} \label{ch:lra}

\epigraph{Despite decades of research on Lanczos methods, the theory for the randomized algorithms is more complete and provides strong guarantees of excellent accuracy, whether or not there exist any gaps between the singular values.}{Huamin Li, George C.\ Linderman, Arthur Szlam, Kelly P.\ Stanton, Yuval Kluger, and Mark Tygert, \textit{Algorithm 971: An Implementation of a Randomized Algorithm for Principal Component Analysis} \cite{LLS+17}}

This thesis is not \emph{only} about randomized algorithms for low-rank approximation, but low-rank approximation will play an important role throughout our discussion.
This introductory chapter reviews two types of low-rank approximation, projection approximation  and Nystr\"om approximation.

\myparagraph{Sources}
This chapter is introductory and is not based on any particular research article.
The last section is adapted from the blog post \cite{Epp24c}.

\myparagraph{Outline}
\Cref{sec:lra-basics} begins by describing (approximately) low-rank matrices and how they can be represented.
We then go on to discuss two types of low-rank approximations: projection approximations (\cref{sec:projection-approximation}) and Nystr\"om approximations (\cref{sec:nystrom}).
Examples of the former type of approximation is given by the randomized SVD (\cref{sec:rsvd}), which can be improved using subspace iteration (\cref{sec:rsi}).
\Cref{sec:gram-correspondence} concludes by describing a connection between projection approximations and Nystr\"om approximations I have termed the \emph{Gram correspondence}.
The Gram correspondence will be an important tool throughout the first part of this thesis.

\section{Low-rank approximation basics} \label{sec:lra-basics}

The term ``low-rank matrix'' is an informal one.
Colloquially, we say a matrix $\mat{B} \in \field^{m\times n}$ is low-rank if $(\rank \mat{B})$ is much smaller than the dimensions $m$ and $n$.
Low-rank matrices are convenient to work with computationally because a rank-$k$ may be decomposed as an outer product
\begin{equation} \label{eq:general-low-rank-factored}
    \mat{B} = \mat{F}\mat{G}^*
\end{equation}
of thin matrices $\mat{F} \in \field^{m\times k}$ and $\mat{G} \in \field^{n\times k}$.
Consequently, a low-rank matrix $\mat{B}$ can be represented by storing the factor matrices $\mat{F}$ and $\mat{G}$, resulting in a significant reduction in storage cost $(m+n)k \ll mn$.

The starting point of this part of the thesis is a surprising and useful observation: Among the numerous matrices appearing in applications, many of them possess the property of being \emph{well-approximated by a low-rank matrix}.\index{low-rank approximation!representation as a matrix factorization}
For such a matrix $\mat{B} \in \field^{m\times n}$, it can be efficiently represented and processed by storing a low-rank approximation $\Bhat\approx \mat{B}$ in factored form $\Bhat = \mat{F}\mat{G}^*$.
Some explanations for the ubiquity of approximately low-rank matrices in applications are provided in \cite{BT19,UT19}.

Given a factored low-rank approximation $\Bhat = \mat{F}\mat{G}^*$ to a matrix $\mat{B}$, many computational tasks become easy:\index{low-rank approximation!reasons to compute}
\begin{itemize}
    \item \textbf{\textit{Matrix--vector products:}} Matrix--vector products $\vec{z} \mapsto \Bhat\vec{z}$ can be computed in $\order((m+n)k)$ operations by evaluating the product as $\Bhat\vec{z} = \mat{F}(\mat{G}^*\vec{z})$.
    \item \textbf{\textit{Singular values and vectors:}} The low-rank factorization $\Bhat = \mat{F}\mat{G}^*$ can be upgraded to a compact SVD $\Bhat = \mat{U}\mat{\Sigma}\mat{V}^*$ in $\order((m+n)k^2)$ operations.
    From this decomposition, one can obtain the singular values and singular vectors of $\Bhat$, and one can easily compute other quantities such as unitarily invariant norms and projectors onto singular subspaces.
    \item \textbf{\textit{Linear systems and least-squares:}}
    Shifted linear systems of the form
    \begin{equation*}
        (\Bhat + \mu \Id) \vec{x} = \vec{c} \quad \text{for } \mu \in \complex \setminus \{0\}
    \end{equation*}
    and regularized least-squares problems 
    \begin{equation*}
        \vec{x} = \argmin_{\vec{z} \in \field^n} \norm{\Bhat \vec{z} - \vec{c}}^2 + \mu \, \norm{\vec{z}}^2 \quad \text{for } \mu > 0
    \end{equation*}
    can both be solved in $\order((m+n)k^2)$ operations (by the Sherman--Morrison--Woodbury identity in the former case and a compact SVD in the latter case).
\end{itemize}
Thus, given both the great prevalence of approximately low-rank matrices in applications and the great computational benefits of exploiting low-rank structure, the matrix low-rank approximation problem is of fundamental interest in computational mathematics.

The central theoretical result in low-rank approximation is the Schmidt--Mirsky--Eckart--Young theorem, which characterizes the optimal low-rank approximation \cite[Thm.~IV.4.18]{SS90}.\index{low-rank approximation!optimality results}\index{Eckart--Young theorem}

\begin{fact}[Optimal low-rank approximations] \label{fact:eckart-young}
    Let $\mat{B} = \mat{U}\mat{\Sigma}\mat{V}^* \in \field^{m\times n}$ be a matrix and its SVD, and let $0\le k \le \min(m,n)$ be an integer.
    Measured with respect to any unitarily invariant norm $\uinorm{\cdot}$, an optimal low-rank approximation
    \begin{equation*}
        \uinorm{\mat{B} - \Bhat} = \min_{\rank (\mat{C}) \le k} \uinorm{\mat{B} - \mat{C}}
    \end{equation*}
    is given by truncating the SVD to level $k$:
    \begin{equation*}
        \lowrank{\mat{B}}_k \coloneqq \mat{U}(:,1:k) \mat{\Sigma}(1:k,1:k)\mat{V}(:,1:k)^*.
    \end{equation*}
    In particular, the matrix $\lowrank{\mat{B}}_k$ is an optimal low-rank approximation with respect to the trace, Frobenius, and spectral norms.
    With respect to any of these norms, the best approximation is unique if and only if $\sigma_k(\mat{B}) > \sigma_{k+1}(\mat{B})$.
\end{fact}

Throughout this thesis, $\lowrank{\mat{B}}_k$ will denote \warn{any} optimal low-rank approximation in the sense of this theorem.
As this result highlights, the best approximation may not be unique.

As 
\cref{fact:eckart-young} demonstrates, an optimal rank-$k$ approximation to matrix $\mat{B}$ can be computed by forming the SVD and truncating to level $k$.
However, this approach is computationally expensive, since computing even an economy-size SVD of a matrix $\mat{B}$ requires $\order(mn \min\{m,n\})$ operations.
Therefore, it is natural to seek faster methods for obtaining a \emph{near-optimal} low-rank approximation.

Before moving on, we should say a few words about low-rank approximation of a psd matrix $\mat{A}$. 
When approximating a psd matrix, it is natural to use a low-rank approximation $\Ahat \approx \mat{A}$ that is also psd.
The factors $\mat{F}$ and $\mat{G}$ for a psd low-rank matrix $\Ahat$ can always be taken to be the same, yielding the symmetric decomposition
\begin{equation} \label{eq:psd-low-rank-factored}
    \Ahat = \mat{F}\mat{F}^*.
\end{equation}
The eigendecomposition of a psd matrix induces an SVD, so an optimal rank-$k$ approximation to a psd matrix can be obtained by truncating its eigendecomposition.\index{positive-semidefinite low-rank approximation}

\section{Projection approximation} \label{sec:projection-approximation}\index{projection approximation!definition|(}

How might we construct high-quality low-rank approximations without relying on the SVD?
A natural strategy emerges from re-examining the optimal rank-$k$ approximation $\lowrank{\mat{B}}_k$.
In \cref{fact:eckart-young}, we computed the optimal rank-$k$ approximation by truncating the SVD $\mat{B} = \mat{U}\mat{\Sigma}\mat{V}^*$.
Alternately, we can compute the optimal rank-$k$ approximation using the \emph{projection formula}:
\begin{equation} \label{eq:projection-formula}
    \lowrank{\mat{B}}_k = \mat{Q}\mat{Q}^* \mat{B} \quad \text{for } \mat{Q} = \mat{U}(:,1:k)
\end{equation}
Here, we have projected $\mat{B}$ onto the span of its $k$ dominant left singular vectors $\mat{U}(:,1:k)$.
We remind the reader that $\mat{Q}\mat{Q}^*$ acts as an orthoprojector onto the range of a matrix $\mat{Q}$ with orthonormal columns.

The projection formula \cref{eq:projection-formula} motivates the use of projections as general strategy for constructing low-rank approximations.

\begin{definition}[Projection approximation] \label{def:projection-approximation}\index{projection approximation!definition}
    Let $\mat{\Omega} \in \field^{n\times k}$ be a test matrix, and probe the matrix $\mat{B}$ by computing the product $\mat{B}\mat{\Omega}$.
    The \emph{projection approximation} to $\mat{B}$ with \emph{test matrix} $\mat{\Omega}$ is
    \begin{equation*}
        \Bhat \coloneqq \mat{\Pi}_{\mat{B}\mat{\Omega}} \mat{B} = \mat{Q}(\mat{Q}^*\mat{B}) \quad \text{with } \mat{Q} = \orth(\mat{B}\mat{\Omega}).
    \end{equation*}
\end{definition}

As we will see, the class of projection approximations contains many good rank-$k$ approximations to a matrix $\mat{B}$.
Indeed, it is simple to verify that the optimal approximation is, itself, a projection approximation with test matrix $\mat{\Omega} = \mat{V}(:,1:k)$.

The motivation behind projection approximation is that the range of $\mat{B}\mat{\Omega}$ serves as a good, but cheap to compute, proxy for the span of the dominant left singular vectors of $\mat{B}$.
To see this, first expand $\mat{B}$ via its SVD $$\mat{B} = \sum_{i=1}^{\min(m,n)} \sigma_i^{\vphantom{*}}(\mat{B}) \vec{u}_i^{\vphantom{*}} \vec{v}_i^*$$
and consider the product
\begin{equation} \label{eq:B-Omega-projection}
    \mat{B}\mat{\Omega} = \sum_{i=1}^{\min(m,n)}\sigma_i^{\vphantom{*}}(\mat{B}) \vec{u}_i^{\vphantom{*}} (\vec{v}_i^*\mat{\Omega}).
\end{equation}
The influence of each left singular vector $\vec{u}_i$ is controlled by the size of the singular value $\sigma_i(\mat{B})$ and the vector--matrix product $\vec{v}_i^*\mat{\Omega}$.
In particular, \warn{provided that $\vec{v}_i^*\mat{\Omega}$ is not small for each $1\le i\le k$}, the product $\mat{B}\mat{\Omega}$ will have large components in the directions of all of the dominant left singular vectors $\{ \vec{u}_i : 1\le i \le k \}$.
Conversely, the subdominant left singular vectors $\{ \vec{u}_i : i > k \}$ are scaled by smaller singular values $\sigma_i(\mat{B}) \le \sigma_k(\mat{B})$, so $\mat{B}\mat{\Omega}$ will have smaller components in these directions.
Thus, projecting onto the range of $\mat{B}\mat{\Omega}$ provides a computationally cheap alternative to projecting onto the span of the dominant singular vectors.
\index{projection approximation!definition|)}

To assess the cost of projection approximation, we can count the number of \emph{matrix--vector products} (\emph{matvecs}) needed to compute one.
Forming the product $\mat{B}\mat{\Omega}$ requires $k$ matvecs with $\mat{B}$, one with each column of $\mat{\Omega}$:
\begin{equation*}
    \mat{B}\mat{\Omega} = \begin{bmatrix} \mat{B}\vec{\omega}_1 & \cdots & \mat{B}\vec{\omega}_k \end{bmatrix}.
\end{equation*}
Building the second product $\mat{Q}^*\mat{B}$ expends $k$ matvecs \warn{with $\mat{B}^*$}:
\begin{equation*}
    \mat{Q}^*\mat{B} = (\mat{B}^*\mat{Q})^* = \begin{bmatrix} \mat{B}^*\vec{q}_1 & \cdots & \mat{B}^*\vec{q}_k \end{bmatrix}^*.
\end{equation*}
Therefore, computing a projection approximation consists of $k$ matvecs with $\mat{B}$, $k$ matvecs with $\mat{B}^*$, and $\order(mk^2)$ additional arithmetic operations to evaluate $\mat{Q} = \orth(\mat{B}\mat{\Omega})$ (via economy-size \QR decomposition).

Projection approximations satisfy a number of enjoyable properties: \index{projection approximation!properties|(}

\begin{proposition}[Properties of projection approximations] \label{prop:projection-properties}
    Let $\mat{B} \in \field^{m\times n}$ and $\mat{\Omega} \in \field^{n\times k}$ be matrices, and consider the projection approximation $\mat{\Pi}_{\mat{B}\mat{\Omega}} \mat{B}$.
    Then
    \begin{enumerate}[label=(\alph*)]
        \item \label{item:projection-orthogonality} \textbf{\textit{Columnwise orthogonality.}}
        The approximation $\mat{\Pi}_{\mat{B}\mat{\Omega}}\mat{B}$ and its residual $\mat{B} - \mat{\Pi}_{\mat{B}\mat{\Omega}}\mat{B}$ are columnwise orthogonal,
        \begin{equation*}
            (\mat{\Pi}_{\mat{B}\mat{\Omega}}\mat{B})^*(\mat{B} - \mat{\Pi}_{\mat{B}\mat{\Omega}}\mat{B}) = (\mat{B} - \mat{\Pi}_{\mat{B}\mat{\Omega}}\mat{B})^* (\mat{\Pi}_{\mat{B}\mat{\Omega}}\mat{B}) = \mat{0}.
        \end{equation*}
        Consequently,
        \begin{align*}
            (\mat{B} - \mat{\Pi}_{\mat{B}\mat{\Omega}}\mat{B})^*(\mat{B} - \mat{\Pi}_{\mat{B}\mat{\Omega}}\mat{B}) &= \mat{B}^*\mat{B} - (\mat{\Pi}_{\mat{B}\mat{\Omega}}\mat{B})^*(\mat{\Pi}_{\mat{B}\mat{\Omega}}\mat{B}) \\ 
            &= \mat{B}^*(\Id - \mat{\Pi}_{\mat{B}\mat{\Omega}}\mat{B})\mat{B}.
        \end{align*}
        \item \textbf{\textit{Invariance.}} \label{item:projection-invariance} The projection approximation is invariant under right multiplication of the test matrix $\mat{\Omega}$ by a nonsingular matrix $\mat{T}$, $\mat{\Pi}_{\mat{B}\mat{\Omega}\mat{T}}\mat{B} = \mat{\Pi}_{\mat{B}\mat{\Omega}} \mat{B}$.
        In particular, $\mat{\Pi}_{\mat{B}\mat{\Omega}} \mat{B}$ is invariant under reordering of the columns of $\mat{\Omega}$.
        \item \textbf{\textit{Monotonicity.}} \label{item:projection-monotonicity}
        The singular values of the approximation $\mat{\Pi}_{\mat{B}\mat{\Omega}}\mat{B}$ are monotone increasing under extension of the test matrix:
        \begin{equation*}
            \vec{\sigma}(\mat{\Pi}_{\mat{B}\mat{\Omega}'}\mat{B}) \ge \vec{\sigma}(\mat{\Pi}_{\mat{B}\mat{\Omega}}\mat{B}) \quad \text{for } \mat{\Omega}' = \flatonebytwo{\mat{\Omega}}{\mat{\Gamma}}.
        \end{equation*}
        Similarly, the singular values of the residual $\mat{B} - \mat{\Pi}_{\mat{B}\mat{\Omega}}\mat{B}$ are monotone decreasing:
        \begin{equation*}
            \vec{\sigma}(\mat{B} - \mat{\Pi}_{\mat{B}\mat{\Omega}'}\mat{B}) \le \vec{\sigma}(\mat{B} - \mat{\Pi}_{\mat{B}\mat{\Omega}}\mat{B}) \quad \text{for } \mat{\Omega}' = \flatonebytwo{\mat{\Omega}}{\mat{\Gamma}}.
        \end{equation*}
        Consequently, $\uinorm{\mat{B} - \mat{\Pi}_{\mat{B}\mat{\Omega}'}\mat{B}} \le \uinorm{\mat{B} - \mat{\Pi}_{\mat{B}\mat{\Omega}}\mat{B}}$ for any unitarily invariant matrix norm $\uinorm{\cdot}$, such as the trace, Frobenius, and spectral norms.
        \item \label{item:projection-optimality} \textbf{\textit{Optimality.}}
        The projection approximation achieves the lowest Frobenius norm error among all approximations $\mat{C}\approx \mat{B}$ satisfying $\range(\mat{C}) \subseteq \range(\mat{B}\mat{\Omega})$:
        \begin{equation*}
            \norm{\mat{B} - \mat{\Pi}_{\mat{B}\mat{\Omega}} \mat{B}}_{\mathrm{F}} = \min_{\range(\mat{C}) \subseteq \range(\mat{B}\mat{\Omega})} \norm{\mat{B} - \mat{C}}_{\mathrm{F}}.
        \end{equation*}
    \end{enumerate}
\end{proposition}

These properties are all more-or-less standard.
We omit the proof.
\index{projection approximation!properties|)}

\index{projection approximation!right vs.\ left|(}
\begin{remark}[Left versus right] \label{rem:left-right-projection} 
    Projection approximations approximate a matrix $\mat{B}$ by multiplying by an orthoprojector $\mat{\Pi}_{\mat{B}\mat{\Omega}}$ on the left. 
    One can also define \emph{right projection approximations} which apply a projector on the right
    \begin{equation*}
        \mat{B} \mat{\Pi}_{\mat{B}^*\mat{\Psi}} \quad \text{for } \mat{\Psi} \in \field^{m\times k}.
    \end{equation*}
    Left and right projection approximations are formally equivalent, as a right projection approximation of $\mat{B}$ is the adjoint of a left projection approximation to $\mat{B}^*$.
    We will use right projection approximations in \cref{ch:row-norm}, and we will discuss two-sided projection approximations in \cref{ch:cur}.
\end{remark}
\index{projection approximation!right vs.\ left|)}

\section{The randomized SVD} \label{sec:rsvd}

The randomized SVD is a simple and popular method for low-rank approximation, and its output is a projection approximation.\index{randomized SVD|(}
At its simplest level, the randomized SVD approximation can be computed in four steps:
\begin{enumerate}
    \item Generate a random matrix $\mat{\Omega} \in \field^{n\times k}$, constructed without looking at the matrix $\mat{B}$.
    (E.g., $\mat{\Omega}$ could be a standard Gaussian matrix.)
    \item Compute the matrix product $\mat{Y} = \mat{B}\mat{\Omega}$.
    \item Form an orthogonal basis $\mat{Q} = \orth(\mat{Y})$ for the column span of $\mat{Y}$, say by \QR factorization $\mat{Y} = \mat{Q}\mat{R}$.
    \item Evaluate the projection $\mat{C} = \mat{B}^* \mat{Q}$, defining the low-rank approximation $\Bhat = \mat{Q}\mat{Q}^* \mat{B} = \mat{Q}\mat{C}^*$ represented by factors $\mat{Q}$ and $\mat{C}$.
\end{enumerate}
Written as so, the name ``randomized SVD'' is a misnomer for this algorithm, as it does not output a low-rank approximation in SVD form.
If desired, one can ``upgrade'' the approximation $\Bhat = \mat{Q}\mat{C}^*$ to compact SVD form:
\begin{enumerate}\setcounter{enumi}{4}
    \item Compute an economy-size SVD $\mat{C}^* = \mat{W}\widehat{\mat{\Sigma}}\hatbold{V}^*$ and set $\hatbold{U} \coloneqq \mat{Q}\mat{W}$.
    The low-rank approximation $\Bhat = \hatbold{U}\widehat{\mat{\Sigma}}\hatbold{V}^*$ is now expressed in compact SVD form, described by factors $\hatbold{U}$, $\widehat{\mat{\Sigma}}$, and $\hatbold{V}$.
\end{enumerate}
We recognize the output $\Bhat = \mat{Q}\mat{C}^* = \hatbold{U}\widehat{\mat{\Sigma}}\smash{\hatbold{V}}^*$ of the randomized SVD as the projection approximation of $\mat{B}$ generated by test matrix $\mat{\Omega}$.
The cost of steps 1--5 is $k$ matvecs with $\mat{B}$ and $k$ matvecs with $\mat{B}^*$, plus an additional $\order(k^2(m+n))$ operations to compute a \QR decomposition of $\mat{Y}$ and compute an SVD of $\mat{C}^*$.
Code for the randomized SVD is provided in \cref{prog:rsvd}.
In this thesis, we will use the name ``randomized SVD'' to refer to the low-rank approximation $\Bhat$, regardless of whether it is represented as $\Bhat = \mat{Q}\mat{C}^*$ or $\Bhat = \hatbold{U}\widehat{\mat{\Sigma}}\hatbold{V}^*$.

\myprogram{Implementation of the randomized SVD with a (standard) Gaussian test matrix.}{}{rsvd}

The randomized SVD algorithm in its modern form was introduced in the famous paper of Halko, Martinsson, and Tropp \cite{HMT11}.
See \cite[\S2]{HMT11} and \cite[\S3]{TW23} for a discussion of the history of this algorithm, including earlier references featuring algorithms similar to the modern randomized SVD.

The randomized SVD produces an approximation comparable with the best rank-$r$ approximation, where $r$ is smaller than the parameter $k$ used in the randomized SVD algorithm.
Here is an example result \cite[Thm.~8.7]{TW23} (see also \cite[\S10]{HMT11}): 

\index{randomized SVD!a priori error bounds@\emph{a priori} error bounds|(}
\begin{fact}[Randomized SVD error analysis] \label{fact:rsvd-error}
    Let $\mat{\Omega} \in \real^{n\times k}$ be populated with independent \warn{real} standard Gaussian entries, and let $\Bhat$ denote the output of the randomized SVD algorithm. 
    Then 
    \begin{align}
        \expect \norm{\mat{B} - \Bhat}_{\mathrm{F}}^2 &\le \min_{r\le k-2} \frac{k-1}{k-r-1}\norm{\mat{B} - \lowrank{\mat{B}}_r}_{\mathrm{F}}^2, \label{eq:rsvd-fro}\\ 
        \expect \norm{\mat{B} - \Bhat}^2 &\le \min_{r\le k-2} \frac{k+r-1}{k-r-1}\left(\norm{\mat{B} - \lowrank{\mat{B}}_r}^2 + \frac{\e^2}{k-r} \norm{\mat{B} - \lowrank{\mat{B}}_r}_{\mathrm{F}}^2\right). \label{eq:rsvd-spec}
    \end{align}
\end{fact}

The first bound \cref{eq:rsvd-fro} shows that the expected Frobenius error of the randomized SVD is no worse with that of the best rank-$r$ approximation for every $r\le k-2$, up to a prefactor $f(k,r)$ depending on $k$ and $r$.
The second bound \cref{eq:rsvd-spec} demonstrates that the \warn{spectral-norm} error of the randomized SVD still depends on the \warn{Frobenius norm} of the best rank-$r$ approximation
\begin{equation} \label{eq:frob-tail}
    \norm{\mat{B} - \lowrank{\mat{B}}_r}_{\mathrm{F}}^2 = \sum_{i=r+1}^{\min(m,n)} \sigma_i^2(\mat{B}).
\end{equation}
Thus, we see that spectral-norm accuracy of the randomized SVD depends on the entire tail of singular values.
Compare with the optimal approximation, whose error is just the $(r+1)$st singular value: $$\norm{\mat{B} - \lowrank{\mat{B}}_r} = \sigma_{r+1}(\mat{B}).$$
The spectral-norm bound \cref{eq:rsvd-spec} demonstrate that the randomized SVD produces a fairly coarse approximation to a matrix $\mat{B}$ when its singular values decay at a slow rate.
This coarse approximation is nonetheless useful for many purposes.\index{randomized SVD|)}
\index{randomized SVD!a priori error bounds@\emph{a priori} error bounds|)}

\section{Randomized subspace iteration} \label{sec:rsi}
\index{randomized subspace iteration|(}
The randomized SVD can be improved by using powering to build a better test matrix $\mat{\Omega}$.
Fix a number $q\ge 0$ of powering steps, and assume at first that $q$ is even.
To build $\mat{\Omega}$, generate a random matrix $\mat{\Gamma} \in \field^{n\times k}$ and apply powering 
\begin{equation*}
    \mat{\Omega} = (\mat{B}^*\mat{B})^{q/2} \mat{\Gamma}.
\end{equation*}
Now, form the projection approximation $\Bhat \coloneqq \mat{\Pi}_{\mat{B}\mat{\Omega}}\mat{B} = \mat{\Pi}_{\mat{B}(\mat{B}^*\mat{B})^{q/2}\mat{\Gamma}} \mat{B}$, represented as either $\Bhat = \mat{Q}\mat{C}^*$ or $\Bhat = \hatbold{U}\widehat{\mat{\Sigma}}\smash{\hatbold{V}}^*$.
We emphasize that the matrix $\mat{\Omega}$ should be formed by successive matrix multiplications 
\begin{equation} \label{eq:powered-Om}
    \mat{\Omega} = \mat{B}^*(\mat{B}(\mat{B}^*(\mat{B} \cdots (\mat{B} \mat{\Gamma}) \cdots))),
\end{equation}
not by explicitly forming and powering the matrix $\mat{B}^*\mat{B}$.
When $q$ is odd, the test matrix is instead defined as
\begin{equation*}
    \mat{\Omega} = (\mat{B}^*\mat{B})^{(q-1)/2} \mat{B}^*\mat{\Gamma} \quad \text{for random } \mat{\Gamma} \in \field^{m\times k}.
\end{equation*}

This algorithm for computing low-rank approximations by powered test matrices is called \emph{randomized subspace iteration} or the \emph{randomized SVD with subspace iteration} \cite{RST10,HMT11,Gu15a,TW23}.
Code is given in \cref{prog:rsi}.
The cost of randomized subspace iteration is $(q+2)k$ matvecs, split roughly evenly between matvecs with $\mat{B}$ and $\mat{B}^*$, plus $\order(k^2(m+n))$ additional operations for {\QR}s and SVDs.

\myprogram{Implementation of the randomized SVD with subspace iteration.}{Warning: This code can be numerically unstable when the matrix $\mat{B}$ has rapidly decaying singular values of the number of subspace iteration steps $q$ is large.}{rsi}

\begin{remark}[Counting subspace iteration steps]
    We have chosen to count the number of steps of subspace iteration $q$ by the total numbers of matrix products with $\mat{B}$ or $\mat{B}^*$ required to form $\mat{\Omega}$.
    Be warned!
    Some authors use a different convention, counting number of multiplications with $\mat{B}^*\mat{B}$.
\end{remark}

To gain intuition for why subspace iteration helps, consider an expansion $\mat{B}\mat{\Omega}$ analogous to \cref{eq:B-Omega-projection}:
\begin{equation*}
    \mat{B}{\mat{\Omega}} = \sum_{i=1}^{\min(m,n)} \sigma_i(\mat{B})^{q+1} \vec{u}_i (\vec{v}_i^*\mat{\Gamma}) \quad \text{for $q$ even}.
\end{equation*}
We see that subspace iteration has the effect of powering the singular values, boosting the gap between the ``large'' singular values $\sigma_1(\mat{B}),\ldots,\sigma_k(\mat{B})$ and the small singular values $\sigma_i(\mat{B})$ for $i > k$.
The name ``subspace iteration'' is derived from the \emph{power iteration}\index{power iteration} method, which computes the dominant eigenvector or singular vector of a matrix by repeatedly multiplying by a matrix $\mat{B}$ (and possibly its adjoint $\mat{B}^*$).
The process \cref{eq:powered-Om} is referred to as \emph{subspace iteration} because powering is performed on a matrix rather than a single vector.
In subspace iteration, the object of interest is not really the powered matrix $\mat{\Omega}$ itself but the subspace $\range(\mat{\Omega})$.

Subspace iteration, even with a random initialization, is a classical approach in matrix computations \cite[Ch.~14]{Par98}.
The modern randomized algorithms literature has sharpened our understanding of subspace iteration by providing sharp probabilistic analysis and emphasizing the computational benefits of using a large block size $k$ with a small number $q$ of subspace iteration steps.

\index{randomized subspace iteration!stable implementation by reorthogonalization|(}
The basic implementation of subspace iteration we've described can be numerically unstable, as the powered matrix $\mat{\Omega}$ given by \cref{eq:powered-Om} can become rank-deficient up to numerical precision. 
This issue can be addressed by using intermediate orthogonalization during the powering process.
That is, instead of using the plain iteration
\begin{equation*}
    \mat{\Omega} \gets \mat{B}^*(\mat{B}\mat{\Omega}) \quad \text{repeated $q/2$ times},
\end{equation*}
as in \cref{prog:rsi}, orthonormalize after each step:
\begin{equation*}
    \mat{\Omega} \gets \orth(\mat{B}^*(\mat{B}\mat{\Omega})) \quad \text{repeated $q/2$ times}.
\end{equation*}
(Being even more aggressive, one could even use the update rule $\mat{\Omega} \gets \orth(\mat{B}^*\orth(\mat{B}\mat{\Omega}))$.)
\warn{In exact arithmetic}, subspace iteration with and without intermediate orthogonalization produce the same projection approximation as output, in view of the invariance property \cref{prop:projection-properties}\ref{item:projection-invariance}.
In floating-point arithmetic, intermediate orthogonalization can significantly improve the quality of the computed projection approximation.
\index{randomized subspace iteration!stable implementation by reorthogonalization|)}

Error bounds for randomized subspace iteration, analogous to \cref{fact:rsvd-error}, are well-established.
See \cite{HMT11,Gu15a,TW23}.
\index{randomized subspace iteration|)}

\index{randomized SVD!with block Krylov iteration|(}
\index{randomized block Krylov iteration|(}
\begin{remark}[Block Krylov iteration] \label{rem:rbki}
    An even more powerful type of projection approximation is given by block Krylov iteration, where one uses the entire family of powered approximations to define the test matrix $\mat{\Omega}$, i.e.,
    \begin{equation*}
        \mat{\Omega} = \begin{bmatrix} \mat{\Gamma} & (\mat{B}^*\mat{B}) \mat{\Gamma} & \cdots (\mat{B}^*\mat{B})^{q/2} \mat{\Gamma}\end{bmatrix} \quad \text{for $q$ even}.
    \end{equation*}
    References on randomized block Krylov iteration include \cite{RST10,MM15,TW23}.
    The weaker approximations produced by the randomized SVD, possibly with a few steps of subspace iteration, will suffice for the purpose of this thesis.
\end{remark}
\index{randomized block Krylov iteration|)}
\index{randomized SVD!with block Krylov iteration|)}

\section{Nystr\"om approximation} \label{sec:nystrom}

\index{Nystr\"om approximation|(}
To approximate a psd matrix $\mat{A} \in \field^{n\times n}$, we have access to a more efficient class of low-rank approximations known as Nystr\"om approximations.

\begin{definition}[Nystr\"om approximation] \label{def:nystrom-approx}
    Let $\mat{A} \in \field^{n\times n}$ be a psd matrix, and let $\mat{\Omega} \in \field^{n\times k}$ be a test matrix.
    The Nystr\"om approximation to $\mat{A}$ is
    \begin{equation} \label{eq:nys}
        \Ahat = \mat{A}\langle \mat{\Omega}\rangle = \mat{A} \mat{\Omega} (\mat{\Omega}^*\mat{A}\mat{\Omega})^\dagger (\mat{A}\mat{\Omega})^*.
    \end{equation}
    Here, ${}^\dagger$ is the Moore--Penrose pseudoinverse, which agrees with the ordinary inverse for nonsingular matrices.
\end{definition}
Observe that the only access to the matrix $\mat{A}$ needed to compute the Nystr\"om approximation $\mat{A}\langle \mat{\Omega}\rangle$ is the ability to form the single matrix product
\begin{equation} \label{eq:nys-Y}
    \mat{Y} \coloneqq \mat{A}\mat{\Omega}.
\end{equation}
From $\mat{Y}$, the approximation $\Ahat$ may be assembled using the formula
\begin{equation*}
    \Ahat = \mat{Y}(\mat{\Omega}^*\mat{Y})^\dagger \mat{Y}^*.
\end{equation*}
This ``single-pass'' property of the Nystr\"om approximation is an advantage over projection approximations, which require two passes over the matrix $\mat{B}$ (one to compute $\mat{B}\mat{\Omega}$, a second to compute $\mat{Q}^*\mat{B}$).

To motivate the form of the Nystr\"om approximation \cref{eq:nys}, observe that any approximation $\Ahat$ satisfying $\range(\Ahat) = \range(\mat{A}\mat{\Omega})$ must take the form
\begin{equation*}
    \Ahat = (\mat{A}\mat{\Omega}) \mat{M} (\mat{A}\mat{\Omega})^*.
\end{equation*}
The choice $\mat{M} = (\mat{\Omega}^*\mat{A}\mat{\Omega})^\dagger$ enforces that the matrix $\mat{A}$ and the approximation $\Ahat$ agree when multiplied by $\mat{\Omega}$, $\Ahat \mat{\Omega} = \mat{A}\mat{\Omega}$.
Indeed, provided $\mat{\Omega}^*\mat{A}\mat{\Omega}$ is nonsingular, $\mat{M} = (\mat{\Omega}^*\mat{A}\mat{\Omega})^\dagger = (\mat{\Omega}^*\mat{A}\mat{\Omega})^{-1}$ is the \emph{unique} choice of $\mat{M}$ satisfying this condition.\index{interpolation!property of Nystr\"om approximation}

\index{Nystr\"om approximation!properties|(}
The randomized Nystr\"om approximation enjoys many nice properties: 

\begin{proposition}[Properties of the Nystr\"om approximtion] \label{prop:nystrom-properties}
    Let $\mat{A} \in \field^{n\times n}$ be psd, let $\mat{\Omega} \in \field^{n\times k}$ be a matrix, and consider the projection approximation $\mat{A}\langle \mat{\Omega}\rangle$.
    Then
    \begin{enumerate}[label=(\alph*)]
        \item \label{item:nystrom-psd} \textbf{\textit{Psd.}} The Nystr\"om approximation $\mat{A}\langle\mat{\Omega}\rangle$ and its residual $\mat{A} - \mat{A}\langle\mat{\Omega}\rangle$ are psd.
        \item \label{item:nystrom-invariance} \textbf{\textit{Invariance.}} The Nystr\"om approximation is invariant under right multiplication of the test matrix $\mat{\Omega}$ by a nonsingular matrix $\mat{T}$, $\mat{A}\langle \mat{\Omega}\mat{T}\rangle = \mat{A}\langle \mat{\Omega}\rangle$.
        In particular, $\mat{A}\langle \mat{\Omega}\rangle$ is invariant to reordering of the columns of $\mat{\Omega}$.
        \item \label{item:nystrom-monotonicity} \textbf{\textit{Monotonicity.}}
        The Nystr\"om approximation is monotone increasing with respect to the psd order under enlargement of the matrix $\mat{\Omega}$
        \begin{equation*}
            \mat{A} \left\langle [\mat{\Omega} \:\:\: \mat{\Gamma}] \right\rangle \succeq \mat{A}\langle \mat{\Omega} \rangle.
        \end{equation*}
        Recall $\succeq$ denotes the psd order on Hermitian matrices.
        Consequently,
        \begin{equation*}
            \uinorm{\mat{A} - \mat{A} \left\langle[\mat{\Omega} \:\:\: \mat{\Gamma}]\right\rangle} \le \uinorm{\mat{A} - \mat{A} \langle \mat{\Omega}\rangle}
        \end{equation*}
        for any unitarily invariant norm $\uinorm{\cdot}$, such as the trace, Frobenius, or spectral norms. 
        \item \label{item:nystrom-interpolatory} \textbf{\textit{Interpolatory.}} The matrix $\mat{A}$ and its Nystr\"om approximation $\mat{A}\langle \mat{\Omega}\rangle$ have the same action on $\mat{\Omega}$. That is, $\mat{A}\cdot \mat{\Omega} = \mat{A}\langle \mat{\Omega}\rangle \cdot \mat{\Omega}$.\index{interpolation!property of Nystr\"om approximation}
        \item \label{item:nystrom-optimality} \textbf{\textit{Optimality.}}
        Among all Hermitian approximations $\mat{M}$ satisfying $\range(\mat{M}) \subseteq \range(\mat{A}\mat{\Omega})$ with a psd residual $\mat{A} - \mat{M}$, the residual $\mat{A} - \mat{A}\langle \mat{\Omega}\rangle$ is the smallest in the psd order:
        \begin{equation*}
            \mat{A} - \Ahat \preceq \mat{A} - \mat{M}.
        \end{equation*}
        Consequently, $\uinorm{\mat{A} - \mat{A}\langle \mat{\Omega}\rangle} \le \uinorm{\mat{A} - \mat{M}}$ for all such $\mat{M}$ and for any unitarily invariant norm $\uinorm{\cdot}$.
    \end{enumerate}
\end{proposition}

These properties have all been used, implicitly or explicitly, throughout the literature.
The first four properties all have relatively straightforward proofs, and the last property relies on variational properties of Schur complements \cite[Thm.~5.3]{And05}; see \cite{Epp22a} for an exposition.
\index{Nystr\"om approximation!properties|)}

\index{Nystr\"om approximation!of Hermitian indefinite matrices|(}
\begin{remark}[Hermitian indefinite matrices] \label{rem:indefinite}
    While it is most naturally justified for psd matrices, the Nystr\"om approximation can also be used to approximate any Hermitian matrix $\mat{A}$.
    A potential issue is that, for an indefinite matrix $\mat{A}$ (i.e., a matrix for which neither $\mat{A}$ or $-\mat{A}$ is psd), the core matrix $\mat{\Omega}^*\mat{A}\mat{\Omega}$ can be singular (or nearly singular) even if the input matrix $\mat{A}$ and test matrix $\mat{\Omega}$ are both full-rank.
    This near-singularity issue can cause degradations in the accuracy of the approximation $\mat{A}\langle \mat{\Omega}\rangle \approx \mat{A}$.
    This issue is addressed by Nakatsukasa and Park \cite{NP23}, who suggest \emph{oversampling} by using a larger test matrix of size $\mat{\Omega} \in \field^{n\times ck}$, then regularizing the Nystr\"om approximation by truncating the small eigenvalues of the core matrix
    \begin{equation*}
        \Ahat \coloneqq \mat{A}\mat{\Omega} \lowrank{\mat{\Omega}^*\mat{A}\mat{\Omega}}_{k}^\dagger (\mat{A}\mat{\Omega}^*).
    \end{equation*}
    They suggest possible values $c = 1.5$ and $c = 2$.
    This strategy appears to resolve the numerical issues of applying the ordinary, un-regularized Nystr\"om approximation to Hermitian indefinite matrices, though only partial theoretical explanation for the success of this strategy is available \cite[Thm.~3.1]{NP23}.
    A disadvantage of this strategy it increases the cost of computing a rank-$k$ Nystr\"om approximation to $ck$ matvecs, comparable to the cost of a randomized SVD ($2k$ matvecs) or generalized Nystr\"om approximation ($\approx 2.5k$ matvecs); both of these approximations $\Bhat \approx \mat{A}$ are not Hermitian by default, but can be made Hermitian by symmetrizing $\Ahat \coloneqq (\Bhat + \Bhat^*)/2$.
    I believe there is more left to be understood about what the ``right'' algorithm is for approximation of Hermitian indefinite matrices.
\end{remark}
\index{Nystr\"om approximation!of Hermitian indefinite matrices|)}

\index{Nystr\"om approximation!with subspace iteration|(}
\index{Nystr\"om approximation!single-pass|(}
\myparagraph{Choice of test matrix}
Analogous to the randomized SVD, a simple way of invoking the Nystr\"om approximation is to choose $\mat{\Omega}$ to be a random matrix independent from $\mat{A}$ such as a standard Gaussian matrix.
We will call this version of Nystr\"om approximation a \emph{single-pass Nystr\"om approximation}, since it requires only one pass over the matrix to compute.
Alternatively, one can use subspace iteration \cite{RST10,HMT11,Gu15a}
\begin{equation*}
    \mat{\Omega} = \mat{A}^q \mat{\Gamma}
\end{equation*}
or block Krylov iteration \cite{RST10,MM15,TW23}
\begin{equation*}
    \mat{\Omega} = \begin{bmatrix}
        \mat{\Gamma} & \mat{A}\mat{\Gamma} & \cdots & \mat{A}^q \mat{\Gamma}
    \end{bmatrix},
\end{equation*}
both of which require additional passes over the matrix.
All of the warnings about numerical stability and reorthogonalization from \cref{sec:rsi} remain in force when randomized subspace iteration (or block Krylov iteration) is combined with Nystr\"om approximation.
\Cref{part:random-pivoting} of this thesis will explore a special class of Nystr\"om approximations where the matrix $\mat{\Omega}$ is a column submatrix of the identity matrix.
\index{Nystr\"om approximation!with subspace iteration|)}
\index{Nystr\"om approximation!single-pass|)}

\index{Nystr\"om approximation!stable implementation by shifting|(}
\myparagraph{Stable implementation}
Computing the Nystr\"om approximation must be done with care to ensure accurate results in floating-point arithmetic. 
Here, we present a variant of the stable Nystr\"om implementation developed in \cite{TYUC17b} (based on ideas introduced in \cite{LLS+17}).
The idea is to compute a Nystr\"om approximation
\begin{equation*}
    \Ahat_\mu \coloneqq (\mat{A} + \mu \Id) \langle \mat{\Omega}\rangle
\end{equation*}
of a shifted matrix $\mat{A} + \mu \Id$, where $\mu$ is a small shift parameter.
Begin by computing the matrix $\mat{Y}$ in \cref{eq:nys-Y}, and define the shift
\begin{subequations} \label{eq:stable-nystrom}
\begin{equation} \label{eq:nys-mu}
    \mu \coloneqq n^{-1/2}\norm{\mat{Y}}_{\mathrm{F}}u.
\end{equation}
Here, $u$ denotes the unit roundoff, of size $u \approx 10^{-16}$ in double-precision arithmetic.
The shift \cref{eq:nys-mu} differs from the one proposed in \cite{TYUC17b} and was introduced in \cite{ETW24} to obtain a stable shift of the minimum possible size while avoiding computing the spectral norm of $\mat{Y}$.
Now, apply the shift to $\mat{Y}$, obtaining
\begin{equation}
    \mat{Y}_\mu \coloneqq \mat{Y} + \nu \mat{\Omega},
\end{equation}
and form the matrix
\begin{equation}
    \mat{H} \coloneqq \mat{\Omega}^*\mat{Y}_\mu.
\end{equation}
Next, compute a Cholesky decomposition
\begin{equation} \label{eq:nys-cholesky}
    \mat{H} = \mat{R}^*\mat{R}
\end{equation}
and use triangular substitution to form
\begin{equation}
    \mat{F} \coloneqq \mat{Y} \mat{R}^{-1}.
\end{equation}
The factor matrix $\mat{F}$ gives rise to the shifted Nystr\"om approximation $\Ahat_\mu = \mat{F}\mat{F}^*$.
Code is provided in \cref{prog:nystrom}.
In order to reuse this code later in the thesis, we have written it to use a matrix $\mat{\Omega} \sim \Unif \{\pm 1\}^{n\times k}$ of random $\pm 1$ values and to return the shift $\mu$, test matrix $\mat{\Omega}$, and Cholesky factor $\mat{R}$ as optional outputs.
\end{subequations}

\myprogram{Stable implementation of the single-pass Nystr\"om approximation.
Low-rank approximation is outputted in the form $\mat{F}\mat{F}^*$ and computed using shifting.}{The \texttt{random\_signs} subroutine is defined in \cref{prog:random_signs}.}{nystrom}

For many purposes, the shifted Nystr\"om approximation $\Ahat_\mu = \mat{F}\mat{F}^*$ is a perfectly good substitute for the unshifted approximation $\Ahat$, as the shift parameter $\mu$ is tiny---on the order of the unit roundoff.
To achieve the most accurate results, however, we can attempt to correct for the shift.
For some problems, we will have means to correct the shift ``exactly''; see \cref{sec:xnystrace-implementation} for an example of such a scenario.
Alternately, we can upgrade the outer product representation $\Ahat_\mu = \mat{F}\mat{F}^*$ to an eigendecomposition representation $\Ahat_\mu = \mat{U}\mat{D}_\mu\mat{U}^*$ via economy-size SVD $\mat{F} = \mat{U}\mat{\Sigma}\mat{V}^*$ with $\mat{D}_\mu \coloneqq \mat{\Sigma}^2$.
Then, define the shift-corrected Nystr\"om approximation
\begin{equation*}
    \Ahat_{\mathrm{SC}} = \mat{U} \mat{D} \mat{U}^* \quad \text{for } \mat{D} = \max \{ \mat{D}_\mu - \mu \Id, \mat{0} \}.
\end{equation*}
The maximum is taken entrywise.
The shift-corrected Nystr\"om approximation $\Ahat_{\mathrm{SC}}$ is not exactly equal to the original Nystr\"om approximation $\Ahat$, but it is typically more accurate then the uncorrected approximation $\Ahat_\mu$.
Code for the shift-corrected Nystr\"om approximation (with a standard Gaussian test matrix $\mat{\Omega}$) is provided in \cref{prog:nystrom_shiftcor}.

\myprogram{Stable implementation of the single-pass Nystr\"om approximation $\mat{U}\mat{D}\mat{U}^*$ using shift correction.}{The \texttt{nystrom} subroutine is defined in \cref{prog:nystrom}.}{nystrom_shiftcor}
\index{Nystr\"om approximation!stable implementation by shifting|)}

\index{randomized SVD!comparison with Nystr\"om approximation|(}
\index{Nystr\"om approximation!comparison with projection approximation|(}
\index{projection approximation!comparison with Nystr\"om approximation|(}
\myparagraph{Nystr\"om versus projection approximation}
Measured using the Frobenius norm, the projection approximation $\mat{\Pi}_{\mat{A}\mat{\Omega}} \mat{A}$ is more accurate than the Nystr\"om approximation $\mat{A}\langle \mat{\Omega}\rangle$, 
\begin{equation}
    \norm{\mat{A} - \mat{\Pi}_{\mat{A}\mat{\Omega}} \mat{A}}_{\mathrm{F}} \le \norm{\mat{A} - \mat{A}\langle \mat{\Omega}\rangle}_{\mathrm{F}}.
\end{equation}
This conclusion follows from the optimality property \cref{prop:projection-properties}\ref{item:projection-optimality} of the projection approximations.
However, this comparison is not really a fair one since the Nystr\"om approximation $\mat{A}\langle \mat{\Omega}\rangle$ can be computed using $k$ matvecs and a single pass over the matrix, whereas the projection approximation $\mat{\Pi}_{\mat{A}\mat{\Omega}} \mat{A}$ requires $2k$ matvecs and two passes.
Therefore, the fair comparison is between the projection approximation $\mat{\Pi}_{\mat{A}\mat{\Omega}} \mat{A}$ and Nystr\"om approximation computed with one step of subspace iteration $\mat{A}\langle \mat{A}\mat{\Omega}\rangle$.
In this case, a result of Tropp and Webber \cite[Lem.~5.2]{TW23} ensures the Nystr\"om approximation is always more accurate:
\begin{equation}
    \uinorm{\mat{A} - \mat{A}\langle \mat{A}\mat{\Omega}\rangle} \le \uinorm{\mat{A} - \mat{\Pi}_{\mat{A}\mat{\Omega}} \mat{A}}
\end{equation}
for any unitarily invariant norm.

\index{Nystr\"om approximation!a priori error analysis@\emph{a priori} error analysis|(}
\myparagraph{Error analysis}
Analogous to \cref{fact:rsvd-error}, we have error bounds for the single-pass Nystr\"om approximation \cite[Cor.~8.8]{TW23}:

\begin{fact}[Single-pass Nystr\"om error analysis] \label{fact:nystrom-error}
    Let $\mat{\Omega} \in \real^{n\times k}$ be populated with independent (real) standard Gaussian entries, and let $\hatbold{A}$ denote the single-pass Nystr\"om approximation with test matrix $\mat{\Omega}$.
    Then 
    \begin{align}
        \expect \norm{\mat{A} - \Ahat}_* &\le \min_{r\le k-2} \frac{k-1}{k-r-1}\norm{\mat{A} - \lowrank{\mat{A}}_r}_*, \label{eq:nys-trace-norm} \\
        \left( \expect \norm{\mat{A} - \Ahat}^2_{\mathrm{F}} \right)^{1/2} &\le \min_{r\le k-4} \frac{k-2}{k-r-1}\left( \norm{\mat{A} - \lowrank{\mat{A}}_r}_{\mathrm{F}} + \frac{\norm{\mat{A} - \lowrank{\mat{A}}_r}_*}{\sqrt{k-r}} \right), \label{eq:nys-fro} \\
        \left(\expect \norm{\mat{A} - \Ahat}^2\right)^{1/2} &\le \min_{r\le k-4} \frac{k+r-1}{k-r-3} \left( \norm{\mat{A} - \lowrank{\mat{A}}_r} + \frac{\sqrt{3} \e^2}{k-r}\norm{\mat{A} - \lowrank{\mat{A}}_r}_* \right).  \label{eq:nys-spec}
    \end{align}
\end{fact}

With these error bounds, we can again consider the question: Which is more accurate, the randomized SVD (a projection approximation) or the single-pass randomized Nystr\"om approximation?
Recall that the error of the randomized SVD, measured in both the Frobenius and spectral norms, depends on the \emph{Frobenius norm} of the best rank-$r$ approximation \cref{eq:frob-tail}.
By contrast, the error of the single-pass Nystr\"om approximation---measured in either the trace, Frobenius, or spectral norm---depends on the \emph{trace norm} of the best rank-$r$ approximation
\begin{equation*}
    \norm{\mat{A} - \lowrank{\mat{A}}_r}_* = \sum_{i=r+1}^{n} \lambda_i(\mat{A}).
\end{equation*}
The trace norm is always larger than the Frobenius norm.
Therefore, when single-pass Nystr\"om approximation and the randomized SVD are both used to produce a rank-$k$ approximation to a psd matrix $\mat{A}$, the randomized SVD is usually more accurate.
However, measured in matvecs, randomized Nystr\"om approximation ($k$ matvecs) is cheaper than the randomized SVD ($2k$ matvecs).
Thus, for a fixed budget of $s$ matvecs, one can either compute a rank-$s$ randomized Nystr\"om approximation or a rank-$(s/2)$ randomized SVD; because of the higher approximation rank, the former is often preferable to the latter.
The consequences of this comparison will be explored in \cref{part:loo} of this thesis.
\index{Nystr\"om approximation!a priori error analysis@\emph{a priori} error analysis|)}
\index{randomized SVD!comparison with Nystr\"om approximation|)}

\index{Gram correspondence|(}
\section{The Gram correspondence} \label{sec:gram-correspondence}

The \emph{Gram correspondence} is a powerful fact that links projection approximations and Nystr\"om approximations.

\begin{theorem}[Gram correspondence] \label{thm:gram-correspondence}
    Let $\mat{A} \in \field^{n\times n}$ be a psd matrix, and let $\mat{B} \in \field^{m\times n}$ be \emph{any} matrix for which $\mat{A} = \mat{B}^*\mat{B}$.
    For a test matrix $\mat{\Omega} \in \field^{n\times k}$, instantiate the projection approximation $\Bhat \coloneqq \mat{\Pi}_{\mat{B}\mat{\Omega}} \mat{B}$ and Nystr\"om approximation $\Ahat \coloneqq \mat{A}\langle \mat{\Omega}\rangle$.
    These approximations satisfy the following relation: $\Ahat = \Bhat^*\Bhat$.
\end{theorem}

\begin{proof}
    Since $\mat{A} = \mat{B}^*\mat{B}$, the Nystr\"om approximation is 
    \begin{equation*}
        \mat{A}\mat{\Omega} = \mat{B}^*[(\mat{B}\mat{\Omega}) ((\mat{B}\mat{\Omega})^*(\mat{B}\mat{\Omega}))^\dagger (\mat{B}\mat{\Omega})^*] \mat{B}
    \end{equation*}
    Observe that the bracketed matrix is a formula for the projector $\mat{\Pi}_{\mat{B}\mat{\Omega}}$, which equals its square.
    Therefore,
    \begin{equation*}
        \mat{A}\mat{\Omega} = \mat{B}^*\mat{\Pi}_{\mat{B}\mat{\Omega}}^2 \mat{B} = (\mat{\Pi}_{\mat{B}\mat{\Omega}}\mat{B})^*(\mat{\Pi}_{\mat{B}\mat{\Omega}}\mat{B}) = \Bhat^*\Bhat.
    \end{equation*}
    We have obtained the advertised conclusion.
\end{proof}

The Gram correspondence has an interesting history, which we will describe later in this section.
Early references in the randomized matrix computations literature include \cite{BW09a,Git11}.

The Gram correspondence can be stated more concisely by using the following definitions.

\begin{definition}[Gram matrix and Gram square root] \label{def:gram}
    Let $\mat{A} \in \field^{n\times n}$ be a psd matrix.
    \warn{Any} matrix $\mat{B} \in \field^{m\times n}$ for which $\mat{B}^*\mat{B} = \mat{A}$ is called a \textit{Gram square root} of $\mat{A}$. \index{Gram square root}
    Similarly, the matrix $\mat{A} = \mat{B}^*\mat{B}$ is called the \emph{Gram matrix} of $\mat{B}$. \index{Gram matrix}
\end{definition}

The Gram matrix, $\mat{A} = \mat{B}^*\mat{B} = (\vec{b}_i^*\vec{b}_j^{\vphantom{*}})_{1\le i,j\le n}$, named after J{\o}rgen Pedersen Gram of Gram--Schmidt fame, tabulates the pairwise inner products of the columns of a matrix $\mat{B}$.
Every psd matrix $\mat{A}$ has many Gram square roots, including the output $\mat{R}$ of a (pivoted) Cholesky factorization and the positive-semidefinite matrix square root $\mat{A}^{1/2}$.

Using \cref{def:gram}, the Gram correspondence may be rewritten:

\actionbox{\textbf{Gram correspondence (rephrased).} If $\mat{B}$ is a Gram square root of $\mat{A}$, the projection approximation $\Bhat \coloneqq \mat{\Pi}_{\mat{B}\mat{\Omega}} \mat{B}$ is \warn{a} Gram square root of the Nystr\"om approximation $\Ahat \coloneqq \mat{A}\langle \mat{\Omega}\rangle$.
(Just as well, the Nystr\"om approximation $\Ahat$ is the Gram matrix of the projection approximation $\Bhat$.)}

Several important observations are special cases of the Gram correspondence, including the equivalence of single-pass Nystr\"om approximation and the randomized SVD and the equivalence of column-pivoted \QR and Cholesky decompositions.
The latter equivalence will play an important role in \cref{part:random-pivoting} of this thesis; see \cref{ch:low-rank-general}.

\index{Gram correspondence!transference of algorithms|(}
The Gram correspondence has important implications for algorithm design.
At a high-level, the principle is as follows:
\TransferenceOfAlgorithms
An elementary example of a pair of algorithms are the randomized SVD for approximating a general matrix and the single-pass Nystr\"om approximation for approximating a psd matrix.
More sophisticated examples of this principle will be explored in \cref{part:random-pivoting} of thesis, most particularly in \cref{ch:low-rank-general}.
\index{Gram correspondence!transference of algorithms|)}

\index{Gram correspondence!transference of error analysis|(}
The Gram correspondence also has a consequence for error analysis of algorithms:

\begin{corollary}[Gram correspondence: Transference of error bounds] \label{cor:gram-transference}
    Let $\mat{B} \in \field^{m\times n}$ be a Gram square root of a psd matrix $\mat{A} \in \field^{n\times n}$, and fix a test matrix $\mat{\Omega} \in \field^{n\times k}$.
    Then the errors for the Nystr\"om approximation $\Ahat = \mat{A}\langle \mat{\Omega}\rangle$ and projection approximation $\Bhat = \mat{\Pi}_{\mat{B}\mat{\Omega}}\mat{B}$ are related: For any $p\ge 1$, we have
    \begin{equation} \label{eq:transference}
        \norm{\mat{A} - \Ahat}_{\set{S}_p} = \norm{\mat{B} - \Bhat}_{\set{S}_{2p}}^2.
    \end{equation}
    Here, $\norm{\cdot}_{\set{S}_p}$ denotes the Schatten $p$-norm, of which the trace, Frobenius, and spectral norms are the special cases $p=1$, $p=2$, and $p=\infty$.
    The special cases $p=1$ and $p = \infty$ in \cref{eq:transference} yield the useful relations
    \begin{align*}
    \begin{aligned}
        \norm{\mat{A} - \Ahat}_* &= \tr(\mat{A} - \Ahat) 
        &= \norm{\mat{B} - \Bhat}_{\mathrm{F}}^2, \\
        \norm{\mat{A} - \Ahat}_{\hphantom{*}} &&=\norm{\mat{B} - \Bhat}^2.
    \end{aligned}
\end{align*}
\end{corollary}

\begin{proof}
    By the columnwise orthogonality property \cref{prop:projection-properties}\ref{item:projection-orthogonality}, 
    \begin{equation*}
        (\mat{B} - \Bhat)^*(\mat{B} - \Bhat) = \mat{B}^*\mat{B} - \Bhat^*\Bhat = \mat{A} - \Ahat.
    \end{equation*}
    Take Schatten $p$-norms of both sides, and invoke the identity $\norm{\mat{C}^*\mat{C}}_{\set{S}_p} = \norm{\mat{C}}_{\set{S}_{2p}}^2$ to obtain \cref{eq:transference}.
\end{proof}

Thus, for any given (random) test matrix $\mat{\Omega}$, one only has to analyze the error $\norm{\mat{B} - \mat{\Pi}_{\mat{B}\mat{\Omega}}}_{\set{S}_{2p}}$ or $\norm{\mat{A} - \mat{A}\langle \mat{\Omega}\rangle}_{\set{S}_p}$ once, with a bound on the other quantity coming for free.
We will use fact several times in this thesis.
As an example we have already seen, the pair of bounds \cref{eq:rsvd-fro,eq:nys-trace-norm} can be derived from each other in this way.
\index{Gram correspondence!transference of error analysis|)}

\begin{remark}[History]
    The Gram correspondence is implicit in much of the literature on low-rank approximation and pivoted matrix decompositions \cite{Hig90a,BW09a,Git11,GS12,GM13,MW17a,TW23,PBK24}.
    The equivalence between Cholesky and \QR decompositions, a consequence of the Gram correspondence, is classical \cite{Hig90a}.
    Transferring algorithms and analysis between a matrix $\mat{B}$ and its Gram matrix $\mat{A} = \mat{B}^*\mat{B}$ has been a standard technique for column-based matrix decompositions and approximations over many years \cite{Hig90a,BW09a,GS12,CK24}; see \cref{ch:low-rank-general} for examples and discussion.
    A version of the Gram correspondence for the particular Gram square root $\mat{B} \coloneqq \mat{A}^{1/2}$ appears in works of Gittens \cite{Git11,GM13} starting in 2011, and this analytical approach is used to obtain sharp bounds for single-pass Nystr\"om approximations in \cite{TYUC17b}.
    Musco and Woodruff \cite{MW17a} provide a clear statement of the transference of algorithms principle.
    Some aspects of the Gram correspondence are highlighted by Derezi\'nski, Khanna, and Mahoney in \cite[Rem.~2]{DKM20}, who draw attention to algorithmic implications in the work of Belabbas and Wolfe \cite{BW09a}.
    The connection between the randomized SVD and randomized Nystr\"om approximation has been used in a very explicit and direct way in recent papers \cite{TW23,PBK24}.
    In an effort to make \cref{thm:gram-correspondence} and its consequences known beyond the community of experts familiar with it, I described the principle in a general way and suggested the name \emph{Gram correspondence} in the blog post \cite{Epp24c}.
\end{remark}
\index{Nystr\"om approximation|)}
\index{Nystr\"om approximation!comparison with projection approximation|)}
\index{projection approximation!comparison with Nystr\"om approximation|)}
\index{Gram correspondence|)}

\partepigraph{Dedicated to my fianc\'ee Sierra, our dogs Hulk and Finn, and our turtle Shelly.}
\part{Random pivoting} \label{part:random-pivoting}

\chapter{Low-rank approximation of psd matrices}

\epigraph{Symmetric positive definiteness is one of the highest accolades to which a matrix can aspire.}{Nicholas J.\ Higham, \emph{Accuracy and Stability of Numerical Algorithms} \cite[\S10.1]{Hig02}}

The first part of this thesis will largely be concerned with the low-rank approximation of \warn{psd} matrices, though we will return to general matrices in \cref{ch:low-rank-general,ch:cur}.
We will focus on a very limited model of computation where we only have access to a small number of \emph{entries} of the input matrix $\mat{A} \in \field^{n\times n}$.
The main algorithm of this part of the thesis will be \emph{randomly pivoted Cholesky}, which produces near-optimal rank-$k$ approximations to a psd matrix after reading only $(k+1)n$ entries.
As an application of these psd low-rank approximation techniques, we can accelerate computations involving kernel matrices and covariance matrices of Gaussian processes.
These matrices are the core objects in a wide class of machine learning algorithms; see \cref{ch:kernels-gaussian} for an introduction to kernel and Gaussian methods in machine learning.

\myparagraph{Sources}
This chapter largely serves to introduce the psd low-rank approximation problem and summarize the existing literature.
It is a significantly extended version of the literature survey from the following paper:

\fullcite{CETW25}.

\myparagraph{Outline}
\Cref{sec:entry-access,sec:psd-low-rank-approx} discuss the entry access model for matrix computations and the positive-semidefinite low-rank approximation problem.
\Cref{sec:cholesky} introduces pivoted partial Cholesky decompositions and shows how they can be used to compute low-rank approximations to psd matrices; the outputs of pivoted partial Cholesky decompositions are called \emph{column Nystr\"om approximations}, which are discussed in \cref{sec:column-nystrom}.
\Cref{sec:subset-selection} describes subset selection problems in machine learning and computational mathematics and relates them to the psd low-rank approximation task.
\Cref{sec:gaussian-cholesky} describes a connection between Cholesky decomposition, Nystr\"om approximation, and Gaussian random variables.
\Cref{sec:sampling-psd-approximation} concludes with a discussion of sampling methods for psd low-rank approximation, which are some of the main alternatives to the randomly pivoted Cholesky\index{randomly pivoted Cholesky} method.

\index{entry access model|(}
\section{The entry access model} \label{sec:entry-access} 

In computational linear algebra, we usually work with matrices stored directly in memory, with all of the matrix entries immediately available to us to perform whatever operations we so choose.
There are also computational settings where we have much more limited access to the matrix, introducing constraints on algorithm design.
The first part of this thesis will work in one such limited framework, the entry access model.

\actionbox{\textbf{Entry access model.} We are given a matrix $\mat{B} \in \field^{m\times n}$ that may be accessed by requesting individual entries $b_{ij}$.
The (dominant) cost of an algorithm is the total number of entries accessed.}

\Cref{part:loo} of this thesis works in a different computational model, the matvec model, which is a natural model for other linear algebraic computations.

The following definition provides a natural example of a matrix for which the entry access model is appropriate:

\begin{definition}[Function matrix] \label{def:function_matrix}
    Let $\set{D}\subseteq \set{X}$ and $\set{E}\subseteq \set{Y}$ be finite subsets of sets $\set{X}$ and $\set{Y}$, and let $\xi : \set{X} \times \set{Y} \to \field$ be a bivariate function.
    The \emph{function matrix}\index{function matrix} associated with these subsets is the matrix $\mat{B} = \xi(\set{D},\set{E}) \coloneqq (\xi(x,y) : x \in \set{D}, y\in \set{E}) \in \field^{\set{D}\times \set{E}}$.
\end{definition}

A function matrix $\mat{B}$ is described implicitly by the subsets $\set{D}, \set{E}$ and the function $\xi$.
Accessing each entry of $\mat{B}$ requires computing the function $\xi(x,y)$.
Consequently, the entire matrix $\mat{B}$ requires $|\set{D}| \cdot |\set{E}|$ function evaluations, much greater than the number $|\set{D}| + |\set{E}|$ of input elements $\set{D} \cup \set{E}$.
Function matrices, and slight variations thereof, occur in discretizations of integral operators\index{integral operator!discretizations of} in computational physics \cite{SS11a} and in the design of fast algorithms for classical structured matrices\index{rank-structured matrix} \cite{CGS+08,Wil21}.

\index{kernel function|(}\index{kernel matrix|(}
The main motivating example for this part of the thesis will be kernel matrices, a subclass of function matrices.

\begin{definition}[Kernel function and matrix] \label{def:kernel_function}
    Let $\set{X}$ be a set.
    A function $\kappa : \set{X} \times \set{X} \to \field$ is said to be a \emph{(positive definite) kernel function} if the \emph{kernel matrix} $\mat{A} \coloneqq \kappa(\set{D},\set{D})$ is psd for every finite subset $\set{D}\subseteq \set{X}$.
\end{definition}

Kernel functions and kernel matrices are flexible tools that can be used to design algorithms for learning from data, and they are central to the theory of Gaussian processes.
See \cref{ch:kernels-gaussian} for examples of kernel functions and an introduction to kernel and Gaussian process methods in machine learning.

Like other function matrices, an $n\times n$ kernel matrix $\mat{A}$ is defined by a dataset $\set{D}$, say, of size $|\set{D}| = n$, but generating all of the entries of $\mat{A}$ requires $n\times n = n^2$ function evaluations.
Particularly when the base space $\set{X} = \field^d$ is a Euclidean space of large dimension $d\gg 1$, evaluating each entry of $\mat{A}$ is expensive, motivating the search for algorithms for psd low-rank approximation that require a small number of entry evaluations.
\index{kernel function|)}\index{kernel matrix|)}

For function matrices and kernel matrices, the entry access model might not be the most appropriate abstraction for algorithm design, as it ignores the fact that modern computer processers are better at generating these matrices in block.
See \cref{sec:why-blocking} for a discussion of a \emph{submatrix access model} that better captures this phenomenon.

\index{low-rank approximation!impossibility results|(}
\index{entry access model!impossibility results for general low-rank approximation|(}
\subsection{Low-rank approximation in the entry access model}

On its face, it may seem impossible to accurately approximate a matrix from a limited number of entry accesses, at least without prior information.
In general, this intuition is correct.

\begin{proposition}[Impossibility of general matrix approximation from entry accesses] \label{prop:general-lra-impossibility}
    Let $\uinorm{\cdot}$ denote either the spectral, Frobenius or trace norms (or, indeed, any unitarily invariant matrix norm).
    Consider an algorithm that queries an input matrix $\mat{B}$ in $t$ positions and outputs an approximation $\Bhat$ to $\mat{B}$.
    Any such algorithm applied to a matrix $\mat{B}$ with a single nonzero entry in a random position must produce an approximation of high relative error
    \begin{equation} \label{eq:bad-approximation}
        \frac{\uinorm{\mat{B} - \Bhat}}{\uinorm{\mat{B}}} \ge \frac{1}{2}
    \end{equation}
    with probability at least $1 - t/mn - 1/(mn-t)$.
    In particular, even querying half of the matrix entries ($t = mn/2$) still produces a poor approximation \cref{eq:bad-approximation} with probability at least $1/2 - 2/mn$.
\end{proposition}

This result establishes that even an algorithm that reads a large fraction of a general matrix $\mat{B}$'s entries (say, half), still is prone to producing an approximation of high error for some inputs.
The failure mode is intuitive; if a single large entry is placed in $\mat{B}$ at a random position, no algorithm can be guaranteed to find it without exhuming a large number of entries.
This observation dates back to the earliest days of randomized matrix approximation \cite{FKV98}.

\begin{proof}[Proof of \cref{prop:general-lra-impossibility}]
    Let $\mat{B} \in \{0,1\}^{m\times n}$ be a random matrix constructed by placing a single nonzero entry in a uniformly random position,
    \begin{equation*}
        \mat{B} = \evec_{i_\star}^{\vphantom{*}}\evec_{j_\star}^* \quad \text{for } i_\star\sim \Unif \{1,\ldots,m\}, \: j_\star \sim \Unif \{1,\ldots,n\}.
    \end{equation*}
    Without loss of generality, we take this nonzero entry to have value $1$.
    Consider any deterministic algorithm which queries this matrix $\mat{B}$ at $t$ adaptively chosen positions $(i_1,j_1),\ldots,(i_t,j_t)$.
    Except with probability $t/mn$, the algorithm only queries entries with value $0$.
    Conditional on the locations, $\mat{B}$ is equally likely to be any matrix in the collection $\set{C} \coloneqq \{ \evec_i^{\vphantom{*}} \evec_j^* : (i,j) \notin \set{P} \}$.
    For the output $\Bhat$ of the algorithm, at most one point in $\set{C}$ is within distance $(1/2)\uinorm{\mat{B}}$ of $\mat{B}$.
    Therefore, except with probability $1/(mn-t)$, $\uinorm{\mat{B} - \Bhat}/\uinorm{\mat{B}}\ge 1/2$.
\end{proof}
\index{entry access model!impossibility results for general low-rank approximation|)}
\index{low-rank approximation!impossibility results|)}

Remarkably, the psd low-rank approximation problem is much better behaved in the entry access model than the general low-rank approximation problem, and we will see several examples in this thesis of algorithms that produce near-optimal low-rank approximations to a psd matrix while reading a fraction its entries.

\index{positive-semidefinite low-rank approximation!reasons for tractability in the entry access model|(}
\index{off-diagonal inequality|(}
The main structural property that makes psd low-rank approximation feasible in the entry access model is the off-diagonal inequality \cite[7.1.P1]{HJ12}:

\begin{fact}[Off-diagonal inequality]
    Let $\mat{A} \in \field^{n\times n}$ be a psd matrix, and let $1\le i,j\le n$ be indices.
    Then the magnitude $|a_{ij}|$ of the off-diagonal entry is bounded by the geometric mean of the diagonal entries:
    \begin{equation*}
        |a_{ij}| \le a_{ii}^{1/2} \cdot a_{jj}^{1/2} \le \max(a_{ii},a_{jj}).
    \end{equation*}
    In particular, the largest entry of a psd matrix must occur on its diagonal.
\end{fact}

\begin{proof}
    The proof is standard and beautiful.
    Since $\mat{A}$ is psd, the principal submatrix
    \begin{equation*}
        \mat{A}(\{i,j\},\{i,j\}) = \twobytwo{a_{ii}}{a_{ij}}{\overline{a_{ij}}}{a_{jj}} \quad \text{is psd as well}.
    \end{equation*}
    Therefore, the determinant $\det \mat{A}(\{i,j\},\{i,j\}) = a_{ii} a_{jj} - |a_{ij}|^2 \ge 0$ is nonnegative.
    Rearrange to obtain the stated conclusion.
\end{proof}

The off-diagonal inequality shows that large entries cannot ``hide'' in a psd matrix.
By generating the $n$ entries of the diagonal, one obtains a ``heat map'' of all possible places a large entry of $\mat{A}$ can lie.
In particular, large entries of $\mat{A}$ can only exist in columns of $\vec{a}_i$ containing a large diagonal entry, $a_{ii} \gg 0$.
This observation suggests a strategy for psd low-rank approximation: Extract columns of $\mat{A}$ with large diagonal entries.
This strategy forms the basis for the most effective methods for psd low-rank approximation in the entry access model.
\index{off-diagonal inequality|)}
\index{positive-semidefinite low-rank approximation!reasons for tractability in the entry access model|)}

\index{entry access model!impossibility results for general low-rank approximation|(}
\begin{remark}[Algorithms for general low-rank approximation from entry accesses] \label{rem:sublinear}
    The impossibility result \cref{prop:general-lra-impossibility} has not stopped research into algorithms for low-rank approximation of non-psd matrices in the entry access model.
    In order to approximate a general matrix from a small number of entry accesses, one needs either \emph{additional information} (such as the location of large entries, the norms of columns, etc.) or \emph{additional assumptions} (such as \emph{incoherence}, the property that the information in the matrix is ``evenly spread across the rows/columns'').
    Discussion of such methods is beyond the scope of this thesis; see \cite{CD13,CY25} for more information.
\end{remark}
\index{entry access model|)}
\index{entry access model!impossibility results for general low-rank approximation|)}

\index{positive-semidefinite low-rank approximation!statement of problem|(}
\section{The psd low-rank approximation problem} \label{sec:psd-low-rank-approx}

This thesis will consider the following version of the psd low-rank approximation:

\actionbox{\textbf{Psd low-rank approximation problem (entry access model):} Given a psd matrix $\mat{A} \in \field^{n\times n}$ and a target rank $k\ge 1$, compute \warn{the description} of a \emph{nearly optimal} rank-$k$ psd approximation $\Ahat$ to $\mat{A}$.}

This problem statement contains embedded in it two phrases that require elaboration: ``the description'' and ``nearly optimal''.
Let us begin with the former. 
For most of this thesis, ``the description'' of a psd rank-$k$ approximation $\Ahat = \mat{F}\mat{F}^*$ will be provided by a factor matrix $\mat{F} \in \field^{n\times k}$.\index{low-rank approximation!representation as a matrix factorization}
Our algorithms will be randomized, so the approximation $\Ahat$ will be a random matrix.

\index{r,epsilon,p-approximation@$(r,\varepsilon,p)$-approximation!definition|(}
To substantiate the phrase ``nearly optimal'', we employ the following definition:

\begin{definition}[$(r,\varepsilon,p)$-approximation] \label{def:r-eps-p}
    Let $\mat{B} \in \field^{m\times n}$ be a matrix, let $1\le r\le \min(m,n)$ be an integer, and $p \ge 1$ and $\varepsilon \ge 0$ be real numbers.
    An \emph{$(r,\varepsilon,p)$-approximation} is a random matrix $\Bhat$ for which
    \begin{equation} \label{eq:r-eps-appox}
        \expect \norm{\mat{B} - \Bhat}_{\set{S}_p} \le (1+\varepsilon) \norm{\mat{B} - \lowrank{\mat{B}}_r}_{\set{S}_p}.
    \end{equation}
    Here, $\norm{\mat{C}}_{\set{S}_p} \coloneqq (\sum_i \sigma_i(\mat{C})^p)^{1/p}$ denotes the Schatten $p$-norm.
    When $p=1$, we suppress the last parameter and call $\Bhat$ an \emph{$(r,\varepsilon)$-approximation}.
\end{definition}

This definition motivates the question:
\actionbox{For a particular algorithm, what rank $k$ is needed to guarantee the algorithm returns an $(r,\varepsilon)$-approximation?}
As we will see, the algorithms we consider in this thesis will have a rank $k$ that depends in a nearly optimal way on $r$ and $\varepsilon$ (in a precise sense that will be outlined below in \cref{fact:optimal-nys-approx}).
This will be the sense in which the algorithms considered by this thesis will be ``near-optimal''.

\index{r,epsilon,p-approximation@$(r,\varepsilon,p)$-approximation!reasons for focus on $p=1$|(}
Let us now speak to the choice of $p=1$ in the definition of an $(r,\varepsilon)$-approximation.
The problem of computing an $(r,\varepsilon,p)$ becomes more difficult as $p$ becomes larger.
As a simple example, consider an $n\times n$ matrix with eigenvalues $10$ and $1$, the latter with multiplicity $n-1$.
For this matrix, even the zero matrix is a near-optimal rank-1 approximation when $p=1$.
Indeed, the zero matrix is a $(1,\order(1/n),1)$-approximation, and its accuracy parameter $\varepsilon = \order(1/n)$ vanishes in the limit $n\to\infty$.
However, when $p=\infty$, the zero matrix is merely a $(1,9,\infty)$ approximation, and we must work harder to obtain a small accuracy parameter $\varepsilon$.

The example above demonstrates that the psd low-rank approximation problem gets harder as $p$ gets larger, but why $p=1$ specifically?
Why not $p=2$ or $p=4$?
These questions will be explored later in the thesis, most particularly in \cref{sec:frobenius-psd}.
For now, let us just mention one reason why working with $p=1$ is convenient.
Recall that many of the most effective algorithms for psd low-rank approximation output a Nystr\"om low-rank approximation $\Ahat$; see \cref{def:nystrom-approx}.
For such an approximation, \cref{prop:nystrom-properties}\ref{item:nystrom-psd} ensures that $\mat{0} \preceq \Ahat \preceq \mat{A}$, so the Schatten $1$-norm and trace of the residual matrix $\mat{A} - \Ahat$ coincide:
\begin{equation*}
    \norm{\mat{A} - \Ahat}_{\set{S}_1} = \norm{\mat{A} - \Ahat}_* = \tr(\mat{A} - \Ahat).
\end{equation*}
Thus, \warn{for a Nystr\"om approximation $\Ahat$}, the $(r,\varepsilon)$-approximation condition \cref{eq:r-eps-appox} can be written as 
\begin{equation*}
    \expect \tr(\mat{A} - \Ahat) \le (1+\varepsilon) \tr(\mat{A} - \lowrank{\mat{A}}_r).
\end{equation*}
The trace is a linear functional, which will be very useful during algorithm analysis.
\index{r,epsilon,p-approximation@$(r,\varepsilon,p)$-approximation!reasons for focus on $p=1$|)}

\index{r,epsilon,p-approximation@$(r,\varepsilon,p)$-approximation!high-probability bound|(}
\index{high-probability bounds|(}
\begin{remark}[High probability bounds] \label{rem:r-eps-whp}
    Another standard way of analyzing psd low-rank approximations is to define an \emph{$(r,\varepsilon,p)$ approximation with failure probability $\delta$} as a random matrix $\Bhat \approx \mat{B}$ for which
    \begin{equation} \label{eq:r-eps-whp}
        \norm{\mat{B} - \Bhat}_{\set{S}_p} \le (1+\varepsilon) \norm{\mat{B} - \lowrank{\mat{B}}_r}_{\set{S}_p} \quad \text{with probability at least } 1-\delta.
    \end{equation}
    The high probability guarantee \cref{eq:r-eps-whp} is often considered to be stronger or more desirable than the expectation guarantee \cref{eq:r-eps-appox}, though the guarantees are formally incomparable.
    For the purposes of this thesis, we will consider an algorithm to be theoretically supported if it admits either type of guarantee.
\end{remark}
\index{r,epsilon,p-approximation@$(r,\varepsilon,p)$-approximation!high-probability bound|)}
\index{positive-semidefinite low-rank approximation!statement of problem|)}
\index{r,epsilon,p-approximation@$(r,\varepsilon,p)$-approximation!definition|)}
\index{high-probability bounds|)}

\section{Pivoted partial Cholesky decompositions} \label{sec:cholesky}

A classical approach in computational linear algebra for computing low-rank approximations of a matrix is \emph{partial matrix decomposition}.\index{partial matrix decomposition}
A partial matrix decomposition refers to any standard matrix decomposition (\QR, Cholesky, SVD) where rows or columns of the factor matrices have been deleted or simply not computed.
The most famous example of low-rank approximation by partial matrix decomposition is furnished by the Schmidt--Mirsky--Eckart--Young theorem (\cref{fact:eckart-young}), which states that partial singular value decompositions yield optimal low-rank approximations.\index{Eckart--Young theorem}

\index{Cholesky decomposition!definition|(}
For efficient psd low-rank matrix approximation, we use a partial version of a different matrix decomposition, the \emph{Cholesky decomposition}.
The Cholesky decomposition of a psd matrix $\mat{A}$ represents the matrix as a product $\mat{A} = \mat{L} \mat{L}^*$ of a lower triangular matrix $\mat{L}$ and its adjoint.
The Cholesky decomposition is sometimes written $\mat{A} = \mat{R}^*\mat{R}$ using an \warn{upper triangular} matrix $\mat{R}$.
The two conventions are easily inter-converted, $\mat{R} = \mat{L}^*$.
By truncating the Cholesky decomposition to just the first $k$ columns, we obtain a rank-$k$ approximation $\Ahat^{(k)} \coloneqq \mat{L}(:,1:k)\mat{L}(:,1:k)^*\approx \mat{A}$.
The factored approximation $\Ahat^{(k)} = \mat{L}(:,1:k)\mat{L}(:,1:k)^*$ is known as a \emph{partial Cholesky decomposition of $\mat{A}$}.\index{Cholesky decomposition!partial}
\index{Cholesky decomposition!definition|)}

\index{Gaussian elimination|(}
\index{Cholesky decomposition!computation of|(}
We compute the Cholesky decomposition by Gaussian elimination.
We use $\Ahat^{(i)}$ to denote the approximation produced at step $i$, and we let $\mat{A}^{(i)} \coloneqq \mat{A} - \Ahat^{(i)}$ denote the residual.
Under this notation, the entries of the residual matrix $\mat{A}^{(i)}$ are denoted $\mat{A}^{(i)}(j,k) = a^{(i)}_{jk}$, and $\mat{A}^{(0)}= \mat{A}$ denotes the initial matrix.
With this notation established, the Cholesky decomposition may be described thusly.
For $i=1,\ldots,n$, do the following
\begin{enumerate}
    \item \textbf{Rescale.} Extract the $i$th column $\vec{a}_i^{(i-1)}$ of $\mat{A}^{(i-1)}$ and rescale 
    \begin{equation*}
        \vec{\ell}_i \coloneqq \vec{a}_i^{(i-1)} / \sqrt{a_{ii}^{(i-1)}}.
    \end{equation*}
    The vectors $\vec{\ell}_i$ comprise the columns of the lower triangular matrix $\mat{L}$.
    \item \textbf{Eliminate.} Update the residual $\mat{A}^{(i)} \coloneqq \mat{A}^{(i-1)} - \vec{\ell}_i^{\vphantom{*}}\vec{\ell}_i^*$.
    This step has the effect of zeroing out the matrix in the $i$th row and column.
\end{enumerate}
The procedure succeeds provided all of the diagonal entries $a_{ii}^{(i-1)}$ encountered during factorization are nonzero.
Since the procedure zeros out the $i$th row and column of the matrix at every iteration, and the update rule does not introduce nonzero entries into previously nonzero rows.
The procedure terminates with a decomposition $\mat{A} = \mat{L}\mat{L}^*$.
\index{Gaussian elimination|)}
\index{Cholesky decomposition!computation of|)}




\index{pivoted partial Cholesky decomposition|(}
\subsection{Pivoting, in general}

\index{Gaussian elimination|(}
As a method of low-rank approximation, the standard partial Cholesky decomposition can be ineffective, as it always forms an approximation based on the first $k$ columns which may not span a good low-rank approximation.
The procedure can be modified to eliminate the columns in a general order.
The resulting procedure is called a \emph{pivoted Cholesky decomposition}, and the positions $s_i \in \{1,\ldots,n\}$ that are eliminated at each step $i$ are called \emph{pivot indices}.\index{pivot index}
Concretely, beginning from $\mat{A}^{(0)} = \mat{0}$, do the following for $i=1,\ldots,n$:
\begin{enumerate}
    \item \textbf{Select a pivot.} Choose a \emph{pivot index} $s_i \in \{1,\ldots,n\}$ associated with a nonzero \emph{pivot entry} $a^{(i-1)}_{s_is_i} \ne 0$.
    \item \textbf{Rescale.} Extract and rescale the pivot column
    \begin{equation*}
        \vec{\ell}_i \coloneqq \vec{a}_{s_i}^{(i-1)}/\sqrt{a_{s_is_i}^{(i-1)}}.
    \end{equation*}
    \item \textbf{Eliminate.} Update the residual $\mat{A}^{(i)} \coloneqq \mat{A}^{(i-1)} - \vec{\ell}_i^{\vphantom{*}}\vec{\ell}_i^*$.
    This step has the effect of zeroing out the matrix in the $s_i$th row and column.
\end{enumerate}

The ordered list of pivots $\set{S} = \{s_1,\ldots,s_n\}$ gives rise to a reordering of the rows and columns of the matrix $\mat{A}$.
The matrix $\mat{L}$ produced by this procedure is typically not triangular, but it becomes lower triangular after rearranging its rows $\mat{L}(\set{S},:)$ according to $\set{S}$.
The reordered factor matrix $\mat{L}(\set{S},:)$ is the traditional Cholesky factor of the \warn{symmetrically} reordered psd matrix $\mat{A}(\set{S},\set{S})$.
As with the standard Cholesky decomposition, pivoted partial Cholesky decompositions yield low-rank approximations of the matrix, i.e., $\Ahat = \mat{F}\mat{F}^*$ for $\mat{F} = \mat{L}(:,1:k)$.
\index{Gaussian elimination|)}

\index{pivoted partial Cholesky decomposition!efficient implementation|(}
To make pivoted partial Cholesky decomposition an effective method for low-rank approximation in the entry access model,\index{entry access model} we make a final optimization.
The standard Cholesky procedure overwrites the entire residual matrix at every step, at a cost of $\order(n^2)$ operations.
But this is wasteful, as we only ever need to evaluate the residual in the selected pivot columns.
As a more efficient procedure, we avoid updating the residual explicitly, instead generating columns of $\mat{A}^{(i)}$ as-needed using the formula $\mat{A}^{(i)} = \mat{A} - \Ahat^{(i)}$, where $\Ahat^{(i)}$ denotes the low-rank approximation produced at step $i$.

\myprogram{Pivoted partial Cholesky for computing a low-rank approximation to a psd matrix.}{}{pivpartchol}

Code for this optimized version of the pivoted partial Cholesky decomposition appears as \cref{prog:pivpartchol}.
As with other programs that will be introduced in this part of the thesis, this code interacts with the matriz $\mat{A}$ through a function \texttt{Acol}, defined so that \texttt{Acol(i)} outputs the $i$th column $\vec{a}_i$.
\index{pivoted partial Cholesky decomposition|)}
\index{pivoted partial Cholesky decomposition!efficient implementation|)}

\index{greedy pivoted Cholesky|(}
\subsection{Greedy pivoting} \label{sec:greedy}

The pivoted partial Cholesky decomposition gives us a general procedure for solving psd low-rank approximation problems in the entry access model.
But how should we pick the pivots?
The main algorithm of this part of the thesis, randomly pivoted Cholesky\index{randomly pivoted Cholesky}, uses a \emph{randomized rule} for pivot selection.
Before getting to this method in \cref{ch:rpcholesky}, we review a more classical approach \cite{Hig90a,FS01}: \emph{greedy pivoting} (also known as \emph{diagonal} or \emph{complete pivoting}).

The idea of greedy pivoting is based on the principles we discussed in \cref{sec:entry-access}.
Large entries of the matrix can only lie in columns of the matrix with large diagonal entries.
Therefore, it is sensible to employ a greedy approach: Always choose a maximal diagonal entry of the residual matrix $\mat{A}^{(i)}$ as pivot:
\begin{equation*}
    s_{i+1} \in \argmin_{j} a^{(i)}_{jj}.
\end{equation*}
%
We emphasize that the greedy method always picks the largest diagonal entry of the \warn{current} residual matrix $\mat{A}^{(i)}$, which evolves as the iteration counter $i$ increases.
The residual matrix is zero in the columns of previously selected pivot indices, so greedy selection ensures the same pivot is never selected twice.

\myprogram{Pivoted partial Cholesky decomposition with greedy pivoting to compute a low-rank approximation to a psd matrix.}{}{greedy_chol}

An implementation of greedy pivoted (partial) Cholesky is given in \cref{prog:greedy_chol}.
The main difference with the generic pivoted partial Cholesky decomposition implementation in \cref{prog:pivpartchol} is that we track the diagonal $\diag(\mat{A}^{(i)})$ of the residual matrix.
The diagonal of $\mat{A}$ is provided to the program as an input \texttt{d}, and it is updated every iteration using the identity
\begin{equation*}
    \diag \big(\mat{A}^{(i)}\big) = \diag\left( \mat{A}^{(i-1)} - \vec{f}_i^{\vphantom{*}}\vec{f}_i^*\right) = \diag \big(\mat{A}^{(i-1)}\big) - |\vec{f}_i^{\vphantom{*}}|^2.
\end{equation*}
(Recall that $|\cdot|^2$ denotes the entrywise squared modulus of a vector.)
The greedy pivoted Cholesky algorithm reads $(k+1)n$ entries of the matrix and expends $\order(k^2n)$ operations.
The procedure outputs the factor $\mat{F}$ and the set of pivots $\set{S} = \{s_1,\ldots,s_k\}$.

\index{greedy pivoted Cholesky!history|(}
\begin{remark}[History]
    Low-rank approximation by pivoted partial Cholesky decompositions is classical.
    The use of greedy pivoting for Cholesky decomposition can be traced at least as far back as the work of Lawson and Hanson in 1974 \cite{LH74}.
    The greedy pivoting strategy is also classical, and it goes under the traditional names of \emph{diagonal pivoting} and \emph{complete pivoting}.
    The procedure was incorporated into both \textsf{LINPACK} software package in 1979 \cite{DMBS79} and its replacement \LAPACK \cite{ABB+99}.
    Stability analysis and analysis of the approximation quality was done by Higham \cite{Hig90a}.
    The greedy method received new attention in the kernel machine learning community following the work of Fine and Scheinberg \cite{FS01}.
\end{remark}
\index{greedy pivoted Cholesky!history|)}

\index{greedy pivoted Cholesky!failure modes|(}
\subsection{Failure of greedy pivoting}

The greedy method seems very natural, but it can have significant deficiencies on some examples. Consider, for instance, the matrix
\begin{equation*}
    \mat{A} = \twobytwo{(1+\varepsilon)\Id_{n/10}}{\mat{0}}{\mat{0}}{\onevec_{9n/10}^{\vphantom{*}}\onevec_{9n/10}^*} \quad \text{for } \varepsilon > 0 \text{ small}.
\end{equation*}
Selecting any pivot entry from the (2,2) block produces the rank-one approximation
\begin{equation*}
    \Ahat_{\mathrm{good}} \coloneqq \twobytwo{\mat{0}}{\mat{0}}{\mat{0}}{\onevec_{9n/10}^{\vphantom{*}}\onevec_{9n/10}^*}
\end{equation*}
which achieves a relative trace error of roughly 10\%:
\begin{equation*}
    \frac{\tr(\mat{A} - \Ahat_{\mathrm{good}})}{\tr(\mat{A})} = (1+\order(\varepsilon)) \cdot \frac{1}{10}.
\end{equation*}
However, when the greedy pivoted Cholesky algorithm is run on $\mat{A}$, the diagonal entries in the (1,1) block are slightly larger than the entries in the (2,2) block, so it proceeds by eliminating entries in the (1,1) block one at a time.
In particular, the relative trace error of the output $\Ahat_{\mathrm{greedy}}$ of the greedy method remains above roughly 90\% until over one tenth of the matrix entries have been read:
\begin{equation*}
    \frac{\tr(\mat{A} - \Ahat_{\mathrm{greedy}})}{\tr(\mat{A})} \ge (1-\order(\varepsilon)) \cdot \frac{9}{10} \quad \text{ as long as } k\le \frac{n}{10}.
\end{equation*}
This is a dismal performance for an algorithm; a good choice of pivot will approximate the matrix to $\approx10\%$ relative error in one step, but the greedy pivoted Cholesky method fails to obtain error better than $\approx90\%$ even after reading a tenth of the matrix!
The greedy method can fail in ways that are perhaps even more striking.
The examples are somewhat sophisticated and rely on variations of Kahan's famous matrix;\index{Kahan matrix} we refer the interested reader to \cite{Hig90a} for details.

The greedy pivoting strategy is natural and often works well, but it has a fatal flaw.
The greedy approach focuses entirely on \emph{exploiting} large diagonal entries, but fails to \emph{explore} potentially valuable pivot choices outside of the numerically largest diagonal entry.
This issue is rectified by the randomly pivoted Cholesky\index{randomly pivoted Cholesky} algorithm, which uses randomization to preferentially select large diagonal entries while investigating a broader range of pivot choices.
\index{greedy pivoted Cholesky|)}
\index{greedy pivoted Cholesky!failure modes|)}

\index{column Nystr\"om approximation|(}
\section{Column Nystr\"om approximation} \label{sec:column-nystrom}

We now have a procedure, the pivoted partial Cholesky decomposition, for computing a low-rank approximation to a psd matrix.
But what \emph{is} the output of this procedure?
Is there a formula for it?
What are its properties?

\index{column Nystr\"om approximation!definition|(}
To answer these questions, we begin by asking a more abstract question: \emph{How should we approximate a matrix $\mat{A}$ from a subset of columns $\vec{a}_{s_1},\ldots,\vec{a}_{s_k}$?}
For notational convenience, these columns can be packaged into a submatrix $\mat{A}(:,\set{S})$, indexed by the set $\set{S} = \{ s_1,\ldots,s_k\}$.
Once we know the columns of a psd matrix, we also know its rows $\mat{A}(\set{S},:) = \mat{A}(:,\set{S})^*$.
We can build an approximation $\Ahat \approx \mat{A}$ by interpolating\index{interpolation!property of Nystr\"om approximation} the known rows and columns, resulting in an approximation of the form 
\begin{equation*}
    \Ahat = \mat{A}(:,\set{S}) \mat{W} \mat{A}(\set{S},:) \quad\text{for some $\mat{W}\in\field^{k\times k}$}.
\end{equation*}
It is natural to expect that $\Ahat$ agrees with $\mat{A}$ in the selected columns, 
\begin{equation} \label{eq:nys-interp-condition}
    \Ahat(:,\set{S}) = \mat{A}(:,\set{S}).
\end{equation}
The condition \cref{eq:nys-interp-condition} may be ensured by setting $\mat{W} \coloneqq \mat{A}(\set{S},\set{S})^\dagger$.
If $\mat{A}(\set{S},\set{S})$ is invertible, $\mat{W}$ is the \emph{unique} matrix producing an approximation $\Ahat$ satisfying \cref{eq:nys-interp-condition}.
This reasoning motivates the following definition:

\begin{definition}[Column Nystr\"om approximation]
    Let $\mat{A}\in\field^{n\times n}$ be a psd matrix and let $\set{S} \subseteq \{1,\ldots,n\}$ be a set of indices.
    The \emph{column Nystr\"om approximation} with \emph{pivot set}\index{pivot index!set of} $\set{S}$ is:
    \begin{equation} \label{eq:column-nystrom}
        \mat{A}\langle \set{S}\rangle \coloneqq \mat{A}(:,\set{S})\mat{A}(\set{S},\set{S})^\dagger \mat{A}(\set{S},:).
    \end{equation}
\end{definition}

As the name suggests, a column Nystr\"om approximation $\Ahat = \mat{A}\langle \set{S}\rangle$ is a Nystr\"om approximation $\mat{A}\langle \mat{\Omega}\rangle$ in the sense of \cref{def:nystrom-approx}.
The associated test matrix is $\mat{\Omega} = \Id(:,\set{S})$.
Consequently, column Nystr\"om approximations enjoy all of the properties of Nystr\"om approximations presented in \cref{prop:nystrom-properties}.
\index{column Nystr\"om approximation!definition|)}

\index{column Nystr\"om approximation!computation by pivoted partial Cholesky decomposition|(}
As one might hope, the output of a pivoted partial Cholesky decomposition is a column Nystr\"om approximation.

\begin{fact}[Nystr\"om and Cholesky]
    Let $\Ahat$ be the approximation to $\mat{A}$ produced by the pivoted partial Cholesky algorithm (\cref{prog:pivpartchol}) with pivots $\set{S} = \{s_1,\ldots,s_k\}$.
    Then $\Ahat = \mat{A}\langle \set{S}\rangle$ is the column Nystr\"om approximation with pivot set $\set{S}$.
\end{fact}
\index{column Nystr\"om approximation!computation by pivoted partial Cholesky decomposition|)}

\index{column Nystr\"om approximation!optimal approximation results|(}
\index{low-rank approximation!optimality results|(}
\index{r,epsilon,p-approximation@$(r,\varepsilon,p)$-approximation!optimality results|(}
Having identified the class of column Nystr\"om approximations, it is natural to ask: How accurate can these approximations be?
This question admits a precise answer using the concept of an $(r,\varepsilon)$-approximation.

\begin{fact}[Column Nystr\"om approximation: Approximation quality] \label{fact:optimal-nys-approx}
    Fix parameters $r\ge 1$ and $0 < \varepsilon\le r$.
    For any psd matrix $\mat{A} \in \field^{n\times n}$, there exists a set $\set{S}\subseteq \{1,\ldots,n\}$ of size 
    \begin{equation} \label{eq:best-column-nystrom}
        k = \min \left\{ \left\lceil \frac{r}{\varepsilon} + r - 1 \right\rceil, n\right\}
    \end{equation}
    for which $\mat{A}\langle \set{S}\rangle$ is an $(r,\varepsilon)$-approximation to $\mat{A}$.
    Conversely, there exists a psd matrix $\mat{A}$ for which $k\ge r/\varepsilon$ columns are \emph{necessary} to produce a column Nystr\"om approximation that is an $(r,\varepsilon)$-approximation.
\end{fact}

A version of this result for projection approximation of general matrices was proven by Guruswami and Sinop \cite{GS12}.
Their result extends to $(r,\varepsilon)$-approximation approximations by the Gram correspondence; see \cite{CETW25} for a self-contained proof in the psd setting.
The upper bound \cref{eq:best-column-nystrom} is proven by using DPP sampling;\index{determinantal point process sampling} see \cref{fact:k-dpp-nystrom} below.

\Cref{fact:optimal-nys-approx} shows that $k\approx r/\varepsilon$ columns are necessary for a column Nystr\"om approximation to be an $(r,\varepsilon)$-approximation, at least for a worst-case matrix $\mat{A}$.
Therefore, for most of this thesis, a ``near-optimal'' algorithm for the psd low-rank approximation problem will be one which produces a rank-$k$ column Nystr\"om approximation satisfying the $(r,\varepsilon)$-approximation guarantee where $k$ is nearly equal to $r/\varepsilon$.

Incidentally, we note that there is an even sharper approximation guarantee for the case $r = k$:

\begin{theorem}[Column Nystr\"om approximation: Approximation quality, $r=k$] \label{thm:optimal-nys-approx-2}
    For any psd matrix $\mat{A} \in \field^{n\times n}$, there exists a column subset $\set{S}$ of size $k$ such that $\mat{A}\langle\set{S}\rangle$ is a $(k,k)$-approximation, i.e.,
    \begin{equation*}
        \tr(\mat{A} - \mat{A}\langle \set{S}\rangle) \le (k+1) \tr(\mat{A} - \lowrank{\mat{A}}_k).
    \end{equation*}
    Moreover, for every $\gamma > 0$, there exists a psd matrix $\mat{A} \in \real^{(k+1)\times (k+1)}$ such that 
    \begin{equation*}
        \tr(\mat{A} - \mat{A}\langle \set{S}\rangle) \ge (k+1-\gamma) \tr(\mat{A} - \lowrank{\mat{A}}_k) \quad \text{for every $k$-element subset $\set{S}$}.
    \end{equation*}
\end{theorem}

Observe that the existence result is the $\varepsilon = r$ case of \cref{fact:optimal-nys-approx}.
The lower bound requires a separate argument, which follows from \cite[Prop.~3.3]{DRVW06} and the Gram correspondence.
I provide a proof in \cref{sec:optimal-nys-approx-2}.
\index{column Nystr\"om approximation|)}
\index{column Nystr\"om approximation!optimal approximation results|)}
\index{low-rank approximation!optimality results|)}

\index{subset selection problem|(}
\section{Subset selection problems} \label{sec:subset-selection}

Before moving on, let us draw a connection between the low-rank approximation techniques we have been studying and a different type of computational problem: subset selection.
Subset selection is more of a qualitative problem than a quantitative one: Given a (multi)set $\set{X}$ of $n$ items, we wish to identify a subset $\set{S}\subseteq \set{X}$ of $k\ll n$ \emph{representative items}.
Typically, we want a set this set of representatives to be \emph{diverse}; if one element $x$ is repeated many times in $\set{X}$, only one copy of $x$ is needed in the subset $\set{S}$ as a representative.

There are several applications for subset selection:
\begin{enumerate}
    \item \textbf{\textit{Information retrieval.}} For designing user interfaces to large databases, it can be important to surface a small number of ``recommended'' items.\index{information retrieval}
    Examples include product recommendation \cite{WMG19} and document retrieval \cite{CK06}.
    \item \textbf{\textit{Optimal sensor placement, optimal experimental design, and active learning.}} Given multiple information sources---possible locations to place sensors, possible scientific experiments to run, or possible unlabeled data points to collect labels for---which small subset of sources should I consult to learn the maximum possible amount of information?
    These scenarios are the subjects of the closely related problems of optimal sensor placement \cite{RCV14},\index{sensor placement} optimal experimental design \cite{Puk06},\index{experimental design} and active learning \cite{AKG+14},\index{active learning}
    and each of these problems is an example of a subset selection problem.
    \item \index{genetics|(}\textbf{\textit{Genetics.}} A basic question in biology is to determine a small set of genetic markers that predict an observed trait or that characterize variation in a population.
    This application has served as a main motivation for the development of randomized algorithms for subset selection \cite{PZB+07,MD09a}.\index{genetics|)}
    \item \textbf{\textit{Computational mathematics.}}
    There are many instances of subset selection problems in computational mathematics itself.
    Often, these applications require the selection of a small number of columns from a matrix.
    Examples include rank-structured matrix computations \cite{Mar11,Wil21}\index{rank-structured matrix} tensor network algorithms \cite{OT10,TSL24b},\index{tensor network} and recovery of rational functions from measurements \cite{WDT22}.
\end{enumerate}

\index{subset selection problem!formulation as low-rank approximation problem|(}\index{positive-semidefinite low-rank approximation!connection to subset selection|(}
In order to design and analyze algorithms for the subset selection problem, several different ways of mathematizing the problem have been proposed.
One approach is based on linear algebra. 
We represent the items in $\set{X}$ by columns of a matrix $\mat{B}$ and seek a subset $\set{S}$ of columns that span a good low-rank approximation to $\mat{B}$.
We call this problem the \emph{column subset selection} problem.

There are variants of the column subset selection problem both for general, rectangular matrices and for psd matrices. 
We focus on the latter problem; the former is discussed in \cref{ch:low-rank-general}.
The \emph{psd column subset selection problem} is as follows:

\actionbox{\textbf{Psd column subset selection problem.} Given a psd matrix $\mat{A}$ and a subset size $k$, find a subset $\set{S}$ of $k$ pivots such that the trace-error of the Nystr\"om approximation $\tr(\mat{A} - \mat{A}\langle \set{S}\rangle)$ is as small as possible.}

Phrased in this way, the column subset selection problem seems like a reformulation of the psd low-rank approximation problem.
However, there are reasons to consider these problems as distinct.
For psd low-rank approximation, the output of interest is the low-rank approximation $\Ahat$, possibly generated as a column Nystr\"om approximation $\Ahat = \mat{A}\langle \set{S}\rangle$.
For psd column subset selection, the relevant output is the subset $\set{S}$ itself; the trace-error $\tr(\mat{A} - \mat{A}\langle \set{S}\rangle)$ is useful only instrumentally as a way of measuring subset quality.
Another distinction is that, in many subset selection applications like genetics or product recommendation, one wants \emph{every} element of the subset $\set{S}$ to be ``good''.
(It would be considered a large failure to incorrectly suggest a genetic marker is linked to cancer, for instance.)
By contrast, bad pivots in  low-rank approximation are more of a missed opportunity than a negative:
Bad pivots are not helpful in improving the approximation quality, but they do not hurt it either.
Throughout this part of the thesis, we will use pivoted partial Cholesky decompositions both for low-rank approximation and subset selection.
\index{subset selection problem|)}
\index{subset selection problem!formulation as low-rank approximation problem|)}
\index{positive-semidefinite low-rank approximation!connection to subset selection|)}

\index{Gaussian distribution!connection to Cholesky decomposition and Nystr\"om approximation|(}
\index{pivoted partial Cholesky decomposition!connection to Gaussian conditioning|(}
\index{column Nystr\"om approximation!connection to Gaussian conditioning|(}
\section{Column Nystr\"om approximation and Gaussian random variables} \label{sec:gaussian-cholesky}

\index{Gaussian distribution!conditional|(}
Partial Cholesky decomposition and column Nystr\"om approximation are closely related to conditional distributions of (jointly) Gaussian random variables.
Consider a vector $\vec{z} \sim \Normal_\field(\vec{0},\mat{A})$, and let $\set{S} \subseteq \{1,\ldots,n\}$ be a subset of indices.
What is the distribution of $\vec{z}$ \emph{conditional} on observing the coordinates $\vec{z}(\set{S})$?
This question is answered by the following classical result:

\begin{theorem}[Conditional expectations of Gaussian random vectors] \label{thm:Gaussian-expectations}
    Let $\vec{z} \sim \Normal_\field(\vec{0},\mat{A})$ be a Gaussian random vector with psd covariance matrix $\mat{A}$, and let $\set{S} \subseteq \{1,\ldots,n\}$ be a subset of indices.
    Assume $\mat{A}(\set{S},\set{S})$ is nonsingular.
    Then 
    \begin{align*}
        \expect[\vec{z} \mid \vec{z}(\set{S})] &= \mat{A}(:,\set{S})\mat{A}(\set{S},\set{S})^{-1}\vec{z}(\set{S}) \sim \Normal_\field(\vec{0}, \mat{A}\langle \set{S}\rangle), \\
        \vec{z} \mid \vec{z}(\set{S}) &\sim \Normal_\field(\mat{A}(:,\set{S})\mat{A}(\set{S},\set{S})^{-1}\vec{z}(\set{S}), \mat{A} - \mat{A}\langle \set{S}\rangle).
    \end{align*}
\end{theorem}

This result is carefully developed for real Gaussians in \cite[Ch.~21]{Tro23}, and the extension to the complex case is straightforward.
We recognize the Nystr\"om approximation $\mat{A}\langle \set{S}\rangle$ as the covariance matrix\index{covariance matrix} of the conditional expectation $\expect[\vec{z} \mid \vec{z}(\set{S})]$, and its residual $\mat{A} - \mat{A}\langle \set{S}\rangle$ as the covariance of $\vec{z}$ \emph{conditional} on $\vec{z}(\set{S})$.
\index{Gaussian distribution!conditional|)}

\index{experimental design|(}
\subsection{A simple model of experimental design}

To put this connection between Cholesky factorization, Nystr\"om approximation, and Gaussian random variables into action, consider the following very basic experimental design problem.
A scientist has $n$ experiments she could run, but only the budget to run $k$ of them.
Which experiments should she run to maximize the knowledge she learns?
Equivalently, which experiments should she run to minimize her uncertainty about the outcome of the remaining experiments she did not run?

Suppose that we can model the outcomes of the experiments as Gaussian random variables $\vec{z} \sim \Normal(\vec{m},\mat{A})$ with \warn{known} mean $\vec{m}$ and covariance matrix\index{covariance matrix} $\mat{A}$.
Choosing the optimal set of $k$ experiments amounts to choosing $k$ entries of $\vec{z}$ to observe with the goal of minimizing the sum of the variances of the remaining experimental outcomes, conditional on these measurements:
%
\begin{equation} \label{eq:simplified-experimental-design}
    \minimize_{\set{S} \subseteq \{1,\ldots,n\}} \sum_{j=1}^n \Var(z_j \mid \vec{z}(\set{S})) \quad \text{such that } |\set{S}| = k.
\end{equation}

By \cref{thm:Gaussian-expectations}, the sum of posterior variances is precisely the trace error of the Nystr\"om approximation
\begin{equation*}
    \sum_{j=1}^n \Var(z_j \mid \vec{z}(\set{S})) = \tr(\mat{A} - \mat{A}\langle \set{S}\rangle).
\end{equation*}
Therefore, this experimental design problem is fully equivalent to selecting a column subset generating a good Nystr\"om approximation, measured using the trace error.

\subsection{Connection to Cholesky decomposition}

Let us now apply a pivoted partial Cholesky decomposition to solve the experimental design problem \cref{eq:simplified-experimental-design}.
At each step, we have a subset $\set{S}_i \subseteq \{1,\ldots,n\}$ of experiments we have already decided to run, and we must choose the next experiment.
The diagonal entries of the residual $\mat{A}^{(i)} = \mat{A} - \mat{A}\langle \set{S}_i\rangle$ store the conditional variances
\begin{equation*}
    a^{(i)}_{jj} = \Var(z_j \mid \set{S}_i).
\end{equation*}
The greedy method, introduced in \cref{sec:greedy}, chooses the largest diagonal entry as pivot at each step.
Equivalently, it chooses to run the experiment with the \emph{highest variance}, conditional on the already-run experiments.

This strategy---always run the experiment over which there is the most uncertainty---is very natural.
However, it has a flaw in that it doesn't take into account the \emph{correlations} between experiments.
Consider the following matrix
\begin{equation} \label{eq:greedy-bad-experimental-design}
    \mat{A} = \begin{bmatrix}
        1+2\varepsilon & 0 & 0 & 0 \\
        0 & 1+\varepsilon & 1 & 1 \\
        0 & 1 & 1+\varepsilon & 1 \\
        0 & 1 & 1 & 1+\varepsilon 
    \end{bmatrix} \quad \text{for } \varepsilon > 0 \text{ small}.
\end{equation}
The first experiment has slightly higher variance than the other experiments, so the greedy method will choose to run experiment $1$.
However, the outcomes of experiments $2$, $3$, and $4$ are highly correlated; running any one of these experiments will leave the scientist with tiny uncertainty about the outcome of the other experiments.
This example provides another demonstration of why greed isn't always good for column subset selection, and demonstrates how injecting randomness can help improve column subset selection algorithms.
For this example, just picking an experiment to run at random would give better results than the greedy method 75\% of the time.
As we'll see in the next section and in \cref{ch:rpcholesky}, there are much better algorithms for column subset selection than uniform random selection.
\index{experimental design|)}

\myparagraph{A peek forward}
The connection between Cholesky factorization, Nystr\"om approximation, and Gaussian random variables is a powerful tool.
It forms the basis of \cref{ch:kernels-gaussian,ch:rpcholesky-kernel-gp}, where we will use Nystr\"om approximation to accelerate algorithms for learning from data based on \emph{Gaussian processes}.
\index{Gaussian distribution!connection to Cholesky decomposition and Nystr\"om approximation|)}
\index{pivoted partial Cholesky decomposition!connection to Gaussian conditioning|)}
\index{column Nystr\"om approximation!connection to Gaussian conditioning|)}

\index{positive-semidefinite low-rank approximation!by sampling|(}\index{column Nystr\"om approximation!by sampling|(}
\section{Non-adaptive random sampling methods} \label{sec:sampling-psd-approximation}

Interest in psd low-rank approximation and column subset selection was renewed in the early twenty-first century, driven by efforts to accelerate kernel methods in machine learning \cite{WS00,FS01,DM05}.
In addition to continued focus on deterministic methods like greedy selection (and variations thereof), this wave of interest also spurred the development of randomized methods.

This section will summarize non-adaptive random sampling methods for psd low-rank approximation methods that were developed prior to our work on randomly pivoted Cholesky \cite{CETW25}.
As we will detail, randomly pivoted Cholesky is related to but distinct from the methods presented in this section.

\subsection{Uniform sampling and diagonal-power sampling} \label{sec:uniform-diagonal-power}

\index{uniform sampling for column Nystr\"om approximation|(}
The most basic randomized method for computing a column Nystr\"om approximation is to select the pivot set $\set{S}$ uniformly at random, without replacement.
(Sampling with replacement is fine as well, as duplicated columns have no effect on a column Nystr\"om approximation.)
For many practical problems, uniform sampling produces approximations of high-enough quality, though it can product significantly worse approximations on other problems.

One class of problems for which uniform sampling is poorly suited are problems with diagonal entries that span a wide range of magnitudes.\index{diagonal-power sampling for column Nystr\"om approximation|(}
To obtain better results in this setting, we should adapt the sampling distribution to the size of the diagonal entries.
Choosing a power $p > 0$, we may draw pivots $s_1,\ldots,s_k$ iid from the \emph{diagonal-power sampling} distribution
\begin{equation} \label{eq:diagonal-power}
    s_1, \ldots, s_k \stackrel{\text{iid}}{\sim} \diag(\mat{A})^p.
\end{equation}
The diagonal $\diag(\mat{A})\ge 0$ is entrywise nonnegative since $\mat{A}$ is psd, and the power $p$ is applied to the vector $\diag(\mat{A})$ entrywise.
Recall that we write $s \sim \vec{w}$ to denote a sample $\prob \{s = j\} = w_j / \sum_{k=1}^n w_k$ from any \warn{unnormalized} weight vector $\vec{w} \in \real_+^n$.
The power $p$ in diagonal-power sampling controls the amount of ``greediness'' of the procedure; large values of $p$ lead to sampling the large elements with high probability, and smaller values of $p$ lead to a more uniform distribution.

In 2005, Drineas and Mahoney \cite{DM05} proposed the diagonal-power sampling distribution with $p=2$.
The choice $p=2$ was motivated by a line of work initiated by Frieze, Kannan, and Vempala \cite{FKV98}, who proved results for approximating a general matrix $\mat{B} \in \field^{m\times n}$ by projecting onto selected columns $\mat{B}(:,s_i)$ sampled iid from the squared column-norm distribution\index{squared column-norm sampling for column projection approximation}
\begin{equation*}
    \vec{s}_1,\ldots,\vec{s}_k \sim \scn(\mat{B}).
\end{equation*}
Recall that $\scn(\mat{B}) \in \real_+^n$ denotes the squared column norms of $\mat{B}$.\index{squared row or column norms}
Computing the full column norms in the entry access model is expensive, so Drineas and Mahoney suggested sampling from the squared diagonal entries as an alternative.

\index{Gram correspondence!transference of algorithms|(}
The Gram correspondence (\cref{sec:gram-correspondence}) suggests a different value for $p$.
If we treat $\mat{A} = \mat{B}^*\mat{B}$ as the Gram matrix for a general matrix $\mat{B}$, computing a low-rank approximation by projecting onto a subset of $\mat{B}$'s columns is equivalent to computing a Nystr\"om approximation of $\mat{A}$, and sampling the squared column norms of $\mat{B}$ is equivalent to the diagonal-power sampling rule \cref{eq:diagonal-power} with $p=1$, in view of the identity $\diag(\mat{A}) = \scn(\mat{B})$.
As such, the power $p=1$ could be regarded as the more ``natural'' power for psd low-rank approximation (insofar as squared column norm sampling is the ``natural'' approach for general matrix approximation).
Code for $p=1$ diagonal-power sampling is provided in \cref{prog:diag_sample_nys}.\index{Gram correspondence!transference of algorithms|)}

\myprogram{Diagonal-power sampling with power $p=1$ for computing a Nystr\"om approximation to a psd matrix.}{}{diag_sample_nys}

The choice of power $p$ is often moot for the basic diagonal-power sampling scheme because kernel matrices and covariance matrices in machine learning are often normalized to have a constant diagonal.
However, the choice of $p$ will have more influence when we consider adaptive procedures like \RPCholesky.

Here is an error analysis for diagonal-power analysis when $p=1$.

\index{uniform sampling for column Nystr\"om approximation!theoretical analysis|(}\index{diagonal-power sampling for column Nystr\"om approximation!theoretical analysis|(}
\begin{fact}[Diagonal-power sampling, $p=1$]
    Let $\mat{A} \in \field^{n\times n}$ be a psd matrix, fix $r\ge 1,\varepsilon > 0$, and introduce the relative error of the best rank-$r$ approximation:
    \begin{equation*}
        \eta \coloneqq \frac{\tr(\mat{A} - \lowrank{\mat{A}}_r)}{\tr(\mat{A})}.
    \end{equation*}
    Diagonal-power sampling Nystr\"om approximation with $p=1$ produces a $(r,\varepsilon)$-approximation provided that the pivot set has size
    \begin{equation*}
        k\ge \frac{r-1}{\varepsilon \eta} + \frac{1}{\varepsilon}.
    \end{equation*}
    Moreover, the $k = \order(r/\varepsilon \eta)$ complexity is necessary for a worst-case $\mat{A}$ matrix.
\end{fact}
\index{uniform sampling for column Nystr\"om approximation!theoretical analysis|)}
\index{diagonal-power sampling for column Nystr\"om approximation!theoretical analysis|)}

This result is \cite[Thm.~C.3]{CETW25}, which adapts \cite[eq.~(4)]{FKV98} using the Gram correspondence.
The statement in \cite{CETW25} focuses on the real case, but the proof transfers without issue to complex numbers.

\index{uniform sampling for column Nystr\"om approximation!failure modes|(}
\index{diagonal-power sampling for column Nystr\"om approximation!failure modes|(}
This result demonstrates fundamental limitations of the uniform and diagonal-power sampling approaches, which require $\order(1/\eta)$ pivots to produce a column Nystr\"om approximations of relative error $\eta$.
We will substantially improve on this result with randomly pivoted Cholesky\index{randomly pivoted Cholesky}, which has a much lower cost of $\order(\log(1/\eta))$.
\Cref{tab:rpc-comparison} below presents experiments demonstrating this failure mode for the uniform and diagonal-power sampling methods.\index{uniform sampling for column Nystr\"om approximation!limitations|)}\index{diagonal-power sampling for column Nystr\"om approximation!limitations|)}\index{diagonal-power sampling for column Nystr\"om approximation|)}\index{uniform sampling for column Nystr\"om approximation|)}

\index{ridge leverage scores|(}
\subsection{Ridge leverage score sampling}

\index{ridge regularization!for linear least squares|(}
To motivate the ridge leverage score sampling approach, let us take a brief digression to the subject of \emph{ridge-regularized linear regression}.
Consider the task of fitting a \warn{conjugate linear} mapping $\vec{\chi} \mapsto \vec{\chi}^*\vec{\beta}$ from $\field^n$ to $\field$ from input--output pairs $(\vec{\chi}^{(1)},y_1),\ldots,(\vec{\chi}^{(m)},y_m)\in\field^n\times \field$.
Assemble the inputs $\vec{\chi}^{(i)}$ as \warn{rows} of a matrix $\mat{X}$, defined as
\begin{equation*}
    \mat{X}(i,:) = (\vec{\chi}^{(i)})^* \quad \text{for } i=1,\ldots,m,
\end{equation*}
and collect the outputs into a vector $\vec{y} \in \field^m$.
One natural approach to fitting a linear model is \emph{ridge-regularized linear regression}, which chooses the coefficients $\vec{\beta} \in \field^n$ as the solution to an optimization problem 
\begin{equation} \label{eq:opt-ridge-reg-coeffs}
    \vec{\beta} = \argmin_{\vec{\beta} \in \field^n}\,  \norm{\mat{X}\vec{\beta} - \vec{y}}^2 + \lambda \, \norm{\vec{\beta}}^2.
\end{equation}
The \emph{ridge parameter} $\lambda \ge 0$ sets the amount of regularization.
For $\lambda = 0$, the coefficients $\vec{\beta}$ are taken to be the limiting value of \cref{eq:opt-ridge-reg-coeffs} as $\lambda \downarrow 0$.
The solution $\vec{\beta}$ to \cref{eq:opt-ridge-reg-coeffs} is given by the formulas
\begin{align}
    \vec{\beta} &= (\mat{X}^*\mat{X} + \lambda \Id)^\dagger\mat{X}^*\vec{y} \label{eq:ridge-normal} \\
    &= \mat{X}^*(\mat{X}\mat{X}^* + \lambda \Id)^\dagger \vec{y}.  \label{eq:ridge-adj-normal} 
\end{align}
The solution formulas \cref{eq:ridge-normal,eq:ridge-adj-normal} are known as the \emph{normal equations} and \emph{adjoint normal equations} for the ridge regression problem \cref{eq:opt-ridge-reg-coeffs}.\index{normal equations!for ridge-regularized least squares}
\index{ridge regularization!for linear least squares|)}

The coefficients $\vec{\beta}$ give rise to the predicted values
\begin{equation} \label{eq:predicted-vals-ridge}
    \vec{\hat{y}} = \mat{X}\vec{\beta}.
\end{equation}
for the input data elements $\vec{\chi}^{(1)},\ldots,\vec{\chi}^{(m)}$.
The ridge leverage scores measure the sensitivity of the predictions $\vec{\hat{y}}$ to the data $\vec{y}$:

\begin{definition}[(Ridge) leverage scores of a general matrix] \label{def:ridge-general}
    Let $\lambda \ge 0$.
    The \emph{$\lambda$-ridge leverage scores}\index{ridge leverage scores!definition for general matrices} $\vec{\ell}^\lambda$ of a matrix $\mat{X} \in \field^{m\times n}$ are
    \begin{equation} \label{eq:ridge-lev-def}
        \vec{\ell}^\lambda \coloneqq \left( \frac{\partial \hat{y}_i}{\partial y_i} : 1\le i \le m \right) = \diag(\mat{X}\mat{X}^*(\mat{X}\mat{X}^* + \lambda \Id)^\dagger).
    \end{equation}
    The \emph{leverage scores} $\vec{\ell}$ are defined as the $0$-ridge leverage scores. \index{leverage scores!definition}
\end{definition}

The characterization of the ridge leverage scores as a matrix diagonal follows from the definition \cref{eq:predicted-vals-ridge} of the predicted values and the adjoint normal equations \cref{eq:ridge-adj-normal}.

\index{leverage scores|(}
While they are not our main focus for now, the leverage scores are an important object in randomized matrix computations.
The leverage scores can be computed as the squared row norms of an orthonormal basis matrix for $\mat{X}$.
That is,
\begin{equation*}
\vec{\ell} = \srn(\Orth(\mat{X})).
\end{equation*}
The $i$th leverage score is a measure of how ``important'' row $i$ is to the matrix $\mat{X}$.
These scores have a decades-long history in statistics, where they are used to quantify sensitivity of a regression model to changes in output values $\vec{y}$ \cite[\S3.3.3]{JWHT21}.\index{leverage scores|)}

The equation \cref{eq:ridge-lev-def} shows that the ridge leverage scores of a matrix $\mat{X}$ depend only on the matrix $\mat{A}\coloneqq \mat{X}\mat{X}^*$, which is the Gram matrix \warn{of $\mat{X}^*$}.\index{Gram matrix}
This motivates the following definition.

\begin{definition}[Ridge leverage scores of a psd matrix]  \label{def:ridge-psd}
    Let $\lambda \ge 0$ be a number and $\mat{A} \in \field^{n\times n}$ be a \warn{psd} matrix.
    The \emph{$\lambda$-ridge leverage scores}\index{ridge leverage scores!definition for psd matrices} of $\mat{A}$ are $\vec{\ell}^\lambda \coloneqq \diag(\mat{A}(\mat{A}+\lambda \Id)^\dagger)$, and the \emph{$\lambda$-effective dimension}\index{effective dimension!definition} of $\mat{A}$ is $\mathrm{d}_{\mathrm{eff}}(\lambda) \coloneqq \sum_{i=1}^n \ell^\lambda_i$.
\end{definition}

In principle, this definition could be ambiguous since a psd matrix $\mat{A}$ has two sets of ridge leverage scores: its leverage scores as a general matrix under \cref{def:ridge-general} and its leverage scores as a psd matrix using \cref{def:ridge-psd}.
For our purposes, however, the intended meaning should always be clear, and the ridge leverage scores of a matrix that is stated to be psd will always be given by \cref{def:ridge-psd}.
The ridge leverage scores were originally proposed by Alaoui and Mahoney \cite{AM15}.

The effective dimension $d_{\mathrm{eff}}(\lambda)$ a continuous proxy for the rank of a psd matrix $\mat{A}$, where eigenvalues of $\mat{A}$ that are much smaller than level $\lambda$ are treated as negligible.
The $0$-effective dimension is the algebraic rank, $d_{\mathrm{eff}}(0) = \rank \mat{A}$, and the effective dimension decreases to zero as $\lambda \uparrow +\infty$.

\index{ridge leverage scores!for computing a column Nystr\"om approximation|(}
\index{ridge leverage scores!theoretical results|(}
The ridge leverage scores give natural sampling probabilities for selecting columns for Nystr\"om approximation.
We have the following result, slightly simplified from \cite[Thm.~3]{MM17a}:

\begin{fact}[Ridge leverage score sampling: Spectral norm] \label{fact:rls-spectral}
    Let $\mat{A}\in\field^{n\times n}$ be a psd matrix, $\lambda > 0$ be a ridge parameter, $\vec{\ell}^\lambda$ be the ridge leverage scores, and $\delta \in (0,1)$ be a specified failure probability.
    Let $\smash{\widehat{\vec{\ell}}}^\lambda$ be over-approximations to the ridge leverage scores
    \begin{equation} \label{eq:approx-rls}
        \vec{\ell}^\lambda \le \smash{\widehat{\vec{\ell}}}^\lambda\le c\vec{\ell}^\lambda \quad \text{with } c\ge 1,
    \end{equation}
    and define sampling probabilities
    \begin{equation} \label{eq:rls-probs}
        p_i \coloneqq \max \left\{1, 16 \hat{\ell}_i^\lambda \log\left(\delta^{-1}\sum_{j=1}^n \hat{\ell}_j^\lambda\right) \right\}.
    \end{equation}
    Define pivots $\set{S}$ by including each $1\le i \le n$ in $\set{S}$ independently with probability $p_i$.
    With probability at least $1-\delta$, the pivot set is not too large $|\set{S}| \le 32c \, \mathrm{d}_{\mathrm{eff}}(\lambda) \log(c\, \mathrm{d}_{\mathrm{eff}}(\lambda)/\delta)$\index{effective dimension!in the analysis of ridge leverage score sampling} and
    \begin{equation} \label{eq:rls-spectral}
        \mat{A}\langle \set{S}\rangle \preceq \mat{A} \preceq \mat{A}\langle \set{S}\rangle + \lambda \Id. 
    \end{equation}
\end{fact}

Observe that the lower bound $\mat{A}\langle \set{S}\rangle \preceq \mat{A}$ in \cref{eq:rls-spectral} is true for any Nystr\"om approximation, in view of \cref{prop:nystrom-properties}\ref{item:nystrom-psd}.
This result shows that, if we sample $\order(\mathrm{d}_{\mathrm{eff}}(\lambda) \log \mathrm{d}_{\mathrm{eff}}(\lambda))$ pivots $\set{S}$ using the ridge leverage score (RLS) distribution, then we get a matching upper bound up to additive error $\lambda \Id$.\index{effective dimension!in the analysis of ridge leverage score sampling}
The bound \cref{eq:rls-spectral} is very strong; in particular, it applies the \warn{spectral norm} error bound $\norm{\mat{A} - \mat{A}\langle \set{S}\rangle} \le \lambda$.\index{spectral-norm error bounds for column Nystr\"om approximation}

Musco and Musco also showed bounds for the trace norm.

\begin{fact}[Ridge leverage score sampling: Trace norm]
    Let $\mat{A} \in \field^{n\times n}$ be a psd matrix, $\varepsilon > 0$ be a real number, and $1\le r \le n$ be an integer.
    Set $\lambda \coloneqq (2\varepsilon/r) \tr(\mat{A} - \lowrank{\mat{A}}_r)$, and suppose we have approximate ridge leverage scores $\smash{\widehat{\vec{\ell}}}^\lambda$ satisfying \cref{eq:approx-rls}  and defining sampling probabilities $\vec{p}$ by \cref{eq:rls-probs}.
    Define pivots $\set{S}$ by including each $1\le i \le n$ in $\set{S}$ independently with probability $p_i$.
    With probability at least $1-\delta$, the pivot set is not too large $|\set{S}| = \order(\tfrac{r}{\varepsilon} \log(\tfrac{r}{\varepsilon\delta}))$ and $\mat{A}\langle \set{S}\rangle$ and $\mat{A}\langle \set{S}\rangle$ is an $(r,\varepsilon)$-approximation with failure probability $\delta$ (as in \cref{rem:r-eps-whp}).
\end{fact}
\index{ridge leverage scores!theoretical results|)}

\index{ridge leverage scores!sampling algorithms|(}
Musco and Musco proposed the recursive RLS (RRLS) algorithm for performing approximate ridge leverage score sampling \cite{MM17a}.
They also develop versions of their algorithm that produce a set of pivots $\set{S}$ of a prescribed size $k$ and provide an automatic mechanism for selecting the hyperparameter $\lambda$.
The cost of the algorithm is $\order(nk)$ entry evaluations and $\order(nk^2)$ additional arithmetic operations.
MATLAB and Python implementations of RRLS are available \cite{Van19}.
Alternative algorithms for approximate RLS sampling are SQUEAK \cite{CLV17} and BLESS \cite{RCCR19}.
\index{ridge leverage scores!sampling algorithms|)}

\index{ridge leverage scores!limitations|(}
RLS sampling is a mathematically elegant strategy for psd column subset selection, and it is amazing that is even \emph{possible} to perform approximate RLS sampling using a small number of entry evaluations.
Still, there are reasons to continue searching for a more performant algorithm.
First, in our empirical testing, the available RLS sampling implementations require roughly $2kn$ to $3kn$ entry evaluations to produce a rank-$k$ approximation \cite[\S2.4]{CETW25}; pivoted Cholesky-based approaches typically take just $(k+1)n$ evaluations.
Second, empirical and theoretical analysis suggests the constants in the $\order$-notation for RLS sampling are moderately large, even when RLS sampling is performed exactly; for some examples, RLS sampling can be many orders of magnitude less accurate than alternative approaches for producing approximations of a given rank $k$; see \cref{tab:rpc-comparison}.
Finally, RLS sampling can require $\order(r\log r)$ columns to produce a low-rank approximation comparable with the best rank-$r$ approximation.
As the following example shows, this logarithmic overhead is a real property of the algorithm, not an artifact of the analysis.

\index{coupon collector problem|(}
\begin{example}[Collecting coupons]
    Fix parameter $r$ and consider the psd matrix
    \begin{equation} \label{eq:coupon-collector-psd}
        \mat{A} = \begin{bmatrix} \outprod{\onevec_{n/r}} & \mat{0} & \mat{0} \cdots & \mat{0} \\
        \mat{0} & \outprod{\onevec_{n/r}} & \mat{0} & \cdots & \mat{0} \\
        \mat{0} & \mat{0} & \outprod{\onevec_{n/r}} & \cdots & \mat{0} \\
        \vdots & \vdots & \vdots & \ddots & \vdots  \\
        \mat{0} & \mat{0} & \mat{0} & \cdots & \outprod{\onevec_{n/r}}
        \end{bmatrix} \in \field^{n\times n}.
    \end{equation}
    This matrix is block diagonal, with $r$ equally sized diagonal blocks of all ones.
    This matrix has rank $r$, and achieving a column Nystr\"om approximation comparable to the best rank-$r$ approximation requires selecting a pivot in each block.
    Regardless of the ridge parameter $\lambda \ge 0$, the ridge leverage scores are constant $\vec{\ell}^\lambda = c(\lambda) \onevec$, so ridge leverage score sampling coincides with iid uniform sampling.
    This is an example of the well-known coupon collector problem, so it takes $\Theta(r\log r)$ pivots to attain a pivot in each block with high probability; see \cite[\S3.6]{MR95} for an introduction to the coupon collector problem.
\end{example}
\index{ridge leverage scores!limitations|)}
\index{coupon collector problem|)}
\index{ridge leverage scores|)}\index{ridge leverage scores!for computing a column Nystr\"om approximation|)}

\index{determinantal point process sampling|(}
\index{determinantal point process sampling!for computing a column Nystr\"om approximation|(}
\subsection{Determinantal point process sampling}

The example \cref{eq:coupon-collector-psd} presents a challenging example for any column subset selection method based on drawing pivots iid from \emph{any} distribution.
For this reason, we generically expect iid sampling methods to require $k = \order(r\log r)$ columns to produce an approximation comparable with the best rank-$r$ approximation.
To obtain better results, we can move beyond iid sampling.
Pivoted partial Choleksy methods\index{pivoted partial Cholesky decomposition} are one type of non-iid selection strategy, as these methods select columns one of a time in an \emph{adaptive} way.
Determinantal point process (DPP) sampling is an alternate approach that draws a random sample from a \emph{joint distribution} on all possible subsets of $k$ pivots.

To motivate DPP sampling, let us first motivate why the determinant is a useful metric in measuring the quality of a set $\set{S}$ of pivots.
The following result is informative:

\index{pivot entries!relation to determinants|(}
\begin{proposition}[Determinants and pivot entries] \label{prop:det-pivot}
    Let $\mat{A} \in\field^{n\times n}$ be a psd matrix and let $\set{S} = \{s_1,\ldots,s_k\}$ be a set of pivots.
    Consider the partial Cholesky decomposition with this pivot set, and introduce the residuals $\mat{A}^{(i)} \coloneqq \mat{A} - \mat{A}\langle \{s_1,\ldots,s_i\}\rangle$.
    We have the identity
    \begin{equation*}
        \det[\mat{A}(\set{S},\set{S})] = a^{(0)}_{s_1s_1} a^{(1)}_{s_2s_2}  a^{(2)}_{s_3s_3} \cdots a^{(k-1)}_{s_ks_k}.
    \end{equation*}
\end{proposition}

\begin{proof}
    The matrix $\mat{A}(\set{S},\set{S})$ has Cholesky decomposition
    \begin{equation*}
        \mat{A}(\set{S},\set{S}) = \mat{R}^*\mat{R},
    \end{equation*}
    where the diagonal entries of $\mat{R}$ are $|r_{ii}|^2 = a^{(i-1)}_{s_is_i}$.
    Ergo,
    \begin{equation*}
        \det[\mat{A}(\set{S},\set{S})] = |\det(\mat{R})|^2 = \prod_{i=1}^k |r_{ii}|^2  = a^{(0)}_{s_1s_1} a^{(1)}_{s_2s_2}  a^{(2)}_{s_3s_3} \cdots a^{(k-1)}_{s_ks_k}.
    \end{equation*}
    The desired result is proven.
\end{proof}

We recognize the determinant of the submatrix $\det[\mat{A}(\set{S},\set{S})]$ as the product of all the pivot diagonal entries of the matrix $\mat{A}$ during pivoted Cholesky decomposition.
As such, choosing a pivot set $\set{S}$ yielding a submatrix $\mat{A}(\set{S},\set{S})$ of large determinant corresponds to selecting an order for pivoted partial Cholesky decomposition in which all of the pivot diagonal entries are \emph{simultaneously} large.
\index{pivot entries!relation to determinants|)}

This result suggests a computational strategy of selecting the pivot set $\set{S}$ by maximizing the determinant $\det[\mat{A}(\set{S},\set{S})]$.
Unfortunately, the problem of finding the largest-determinant submatrix of a matrix is \NP-hard \cite{CM09}.
To circumvent this intractability result, we can instead use the determinants to define a \emph{sampling distribution}, motivating the following definition.

\index{determinantal point process sampling!definition|(}
\begin{definition}[Fixed-size DPP] \label{def:dpp}
    Let $\mat{A} \in \field^{n\times n}$ be a psd matrix, and fix a number $1\le k \le n$.
    A \emph{determinantal point process of fixed size $k$} or \emph{$k$-DPP} is a random subset $\set{S} \subseteq \{1,\ldots,n\}$ of size $k$ with distribution
    \begin{equation*}
        \prob\{\set{S} = \set{T}\} = \frac{\det[\mat{A}(\set{T},\set{T})]}{\sum_{|\set{R}| = k} \det[\mat{A}(\set{R},\set{R})]} \quad \text{for each subset } \set{T}\subseteq \{1,\ldots,n\} \text{ of $k$ elements}.
    \end{equation*}
    For a $k$-DPP $\set{S}$, we write $\set{S} \sim \kDPP{k}(\mat{A})$.
\end{definition}
\index{determinantal point process sampling!definition|)}

\index{determinantal point process sampling!history|(}
Determinantal point processes have a storied history.
DPPs were originally introduced in 1975 by Macchi \cite{Mac75} in the context of quantum statistical physics, and they have been studied by mathematicians for decades.
There are many subvarietals of DPPs.
In particular, the fixed-size DPPs in \cref{def:dpp} were introduced in 2011 by Kulesza and Taskar as a more useful construction for use in machine learning \cite{KT11}.
Further, Kulesza and Taskar's definition is equivalent to the \emph{volume sampling distribution} introduced five years earlier by Deshpande, Rademacher, Vempala, and Wang \cite{DRVW06} applied to a Gram square root of $\mat{A}$.
Motivated by subset selection problems in machine learning, interest in DPPs has exploded over the past two decades \cite{KT12}.
More recently, DPPs have been deployed as a tool in matrix computations \cite{DM21a}.
See the surveys \cite{KT12,DM21a} for more on DPPs.
\index{determinantal point process sampling!history|)}

\index{determinantal point process sampling!for column Nystr\"om approximation|(}\index{determinantal point process sampling!theoretical results|(}
Pivot sets selected from a $k$-DPP produce excellent Nystr\"om approximations.
We have the following result:

\begin{fact}[$k$-DPPs for Nystr\"om approximation] \label{fact:k-dpp-nystrom}
    Let $\mat{A} \in \field^{n\times n}$ be a psd matrix, and let $1\le k \le n$.
    Draw a pivot set $\set{S} \sim \kDPP{k}(\mat{A})$.
    The trace error of the Nystr\"om approximation admits the exact formula
    \begin{equation*}
        \expect \tr(\mat{A} - \mat{A}\langle \set{S}\rangle) = (k+1) \frac{\e_{k+1}(\lambda_1(\mat{A}),\ldots,\lambda_n(\mat{A}))}{\e_k(\lambda_1(\mat{A}),\ldots,\lambda_n(\mat{A}))}.
    \end{equation*}
    Here, $\e_j$ is the $j$th elementary symmetric polynomial.
    For any $0 \le r\le k$, this error may be bounded as
    \begin{equation*}
        \expect \tr(\mat{A} - \mat{A}\langle \set{S}\rangle) \le \frac{k+1}{k-r+1} \tr(\mat{A} - \lowrank{\mat{A}}_r).
    \end{equation*}
    Consequently, $k$-DPP sampling produces an $(r,\varepsilon)$-approximation for any $k$ satisfying $k\ge r/\varepsilon+r-1$ for any $\varepsilon \in (0,r]$.
\end{fact}

An analog of this general bound was first established by Guruswami and Sinop \cite{GS12} for column projection approximations computed by volume sampling.
Using the Gram correspondence, it was transplanted to Nystr\"om approximation by $k$-DPP sampling in \cite{DM21a}, yielding \cref{fact:k-dpp-nystrom}.
The case $r=k$ was proven first by \cite{DRVW06} (for column projection approximations and volume sampling) and \cite{BW09} (for Nystr\"om approximation by $k$-DPPs).

\index{determinantal point process sampling!algorithms|(}\index{determinantal point process sampling!limitations|(}
\Cref{fact:k-dpp-nystrom} achieves the strongest theoretical bounds (and indeed, the strongest \emph{existence results}) for any column Nystr\"om method.\index{determinantal point process sampling!theoretical results|)}
But turning DPP sampling into a computational strategy is a challenging enterprise.
Standard algorithms for sampling a general $k$-DPP expend $\order(n^3)$ operations and require a full decomposition of the matrix (see, e.g., \cite[\S3]{KT11}).
Responding to this computational bottleneck, researchers investigated Markov chain Monte Carlo algorithms for \emph{approximate} $k$-DPP sampling \cite{AOR16,RO19}; even with efficient implementations, these algorithms are significantly more expensive than pivoted Cholesky methods.
(The most recent algorithms \cite{ALV22} do achieve the same \warn{asymptotic} runtime as pivoted Cholesky algorithms, but they are complicated, and I am not aware of any evaluations of this approach for use in practice.)
An alternate line of work has investigated exact $k$-DPP sampling algorithms that do not require reading the entire input matrix \cite{Der19,DCV19,CDV20}.
These algorithms are impressive, but they are more expensive in both runtime and number of entry evaluations than pivoted Cholesky methods.
In my experience with the $k$-DPP software package \cite{GPBV19}, $k$-DPP sampling has been significantly slower and more resource intensive than other methods for Nystr\"om approximation.
In addition, I have found this software tends to fail by throwing exceptions on challenging problems.
Therefore, while $k$-DPPs are mathematically beautiful and set the mathematical standard to which other methods are compared, I have not found $k$-DPP sampling to be a competitive approach to large-scale psd matrix approximation with the available algorithms and software implementations.
\index{determinantal point process sampling!algorithms|)}\index{determinantal point process sampling!limitations|)}
\index{positive-semidefinite low-rank approximation!by sampling|)}\index{column Nystr\"om approximation!by sampling|)}
\index{determinantal point process sampling|)}
\index{determinantal point process sampling!for computing a column Nystr\"om approximation|)}


\chapter{Randomly pivoted Cholesky} \label{ch:rpcholesky}

\epigraph{So first she tasted the porridge of the Great, Huge Bear, and that was too hot for her. And then she tasted the porridge of the Middle Bear, and that was too cold for her. And then she went to the porridge of the Little, Small, Wee Bear, and tasted that; and that was neither too 
nor too cold, but just right, and she liked it so well that she ate it all up.}{\emph{The Story of the Three Bears} (1905)}

\index{randomly pivoted Cholesky|(}
\index{uniform sampling for column Nystr\"om approximation|(}\index{greedy pivoted Cholesky|(}
Last chapter, we introduced the pivoted partial Cholesky decomposition as a way of computing a low-rank approximation to a psd matrix, and we saw two extreme strategies, uniform sampling and greedy selection, for selecting the pivots.
The greedy strategy always selects the largest diagonal entry of the residual as pivot, and the uniform strategy selects pivots at random without using any information about the size of the diagonal entries.
This chapter will strike a balance between these approaches, selecting a random pivot at each iteration using sampling probabilities weighted by the diagonal entries.
The resulting algorithm is called \emph{randomly pivoted Cholesky} (\RPCholesky).
\index{uniform sampling for column Nystr\"om approximation|)}\index{greedy pivoted Cholesky|)}

\myparagraph{Sources}
This paper is based on the randomly pivoted Cholesky\index{randomly pivoted Cholesky} paper:

\fullcite{CETW25}.

\myparagraph{Outline}
\Cref{sec:alg-impl} introduces the randomly pivoted Cholesky\index{randomly pivoted Cholesky} algorithm and discusses its implementation.
\Cref{sec:rpc-experiments} provides numerical experiments, and \cref{sec:rpc-analysis} discusses analysis.
We conclude by discussing an extension to Gibbs \RPCholesky (\cref{sec:gibbs}) and connections between \RPCholesky and DPPs (\cref{sec:dpp-connections}).

\index{randomly pivoted Cholesky!description of algorithm|(}
\section{Algorithm and implementation} \label{sec:alg-impl}

Randomly pivoted Cholesky (\RPCholesky) is an algorithm for low-rank approximation of psd matrices. 
It employs a pivoted partial Cholesky decomposition, but is distinguished from other pivoted Cholesky methods by using the diagonal of the residual matrix at each iteration as a \emph{sampling distribution} to select the next pivot.
Conceptually, \RPCholesky executes the following iteration. With initial residual $\mat{A}^{(0)} \coloneqq \mat{A}$, do the following for $i=1,2,\ldots,k$:
\begin{enumerate}
    \item \textbf{Draw a random pivot.} Draw a \emph{random} pivot index
    \begin{equation*}
        s_i \sim \diag(\mat{A}^{(i-1)}).
    \end{equation*}
    This random pivoting strategy preferentially selects larger diagonal entries as pivots, but has a nonzero probability of selecting any nonzero diagonal entry as pivot.
    (Recall that we have defined $s\sim\vec{w}$ as a random index sampled from the \warn{unnormalized} weight vector $\vec{w}$.
    That is, $\prob \{s  = j \} = w_j / \sum_{i=1}^n w_i$.)
    \item \textbf{Rescale.} Extract and rescale the pivot column
    \begin{equation*}
        \vec{f}_i \coloneqq \vec{a}_{s_i}^{(i-1)}/\sqrt{a_{s_is_i}^{(i-1)}}.
    \end{equation*}
    \item \textbf{Eliminate.} Update the residual $\mat{A}^{(i)} \coloneqq \mat{A}^{(i-1)} - \vec{f}_i^{\vphantom{*}}\vec{f}_i^*$, zeroing out the matrix in the $s_i$th row and column.
\end{enumerate}
The outputs of \RPCholesky are a factor matrix $\mat{F} \in \field^{n\times k}$ defining a low-rank approximation $\mat{F}\mat{F}^*\approx \mat{A}$ and a set $\set{S} = \{s_1,\ldots,s_k\}$ of pivot entries.

In practice, we do not update the entire residual matrix at every iteration, instead tracking the diagonal $\diag(\mat{A}^{(i)})$ and generating entries from $\mat{A}^{(i)}$ as needed using the formula $\mat{A}^{(i)} = \mat{A} - \mat{F}(:,1:i)\mat{F}(:,1:i)^*$, as in \cref{prog:greedy_chol}.
With these optimizations, \RPCholesky reads $(k+1)n$ entries of the matrix and expends $\order(k^2n)$ arithmetic operations.\index{randomly pivoted Cholesky!computational cost}
See \cref{prog:rpcholesky} for \RPCholesky code.
More efficient block versions of \RPCholesky will be developed in \cref{ch:blocked}.
\index{randomly pivoted Cholesky!description of algorithm|)}

\myprogram{Randomly pivoted Cholesky for psd low-rank approximation and column subset selection.}{}{rpcholesky}

\index{greedy pivoted Cholesky!failure modes|(}
The \RPCholesky algorithm can be seen as a midpoint between greedy selection and uniform sampling.
The greedy method is based on entirely \emph{exploiting} large diagonal entries, without \emph{exploring} smaller diagonal entries as possible pivots.
Conversely, uniform sampling randomly explores the set of all pivots but does not exploit information about the size of the diagonal entries.
\RPCholesky combines both exploration and exploitation, making it more robust than either strategy individually.

\index{randomly pivoted Cholesky!connection to Gaussian distribution|(}\index{Gaussian distribution!connection to Cholesky decomposition and Nystr\"om approximation|(}
Another way of interpreting the \RPCholesky method comes from the Gaussian framing developed in \cref{sec:gaussian-cholesky}.
Introduce a Gaussian vector $\vec{z} \sim \Normal_\field(\mat{0},\mat{A})$ with covariance matrix $\mat{A}$. 
As in \cref{sec:gaussian-cholesky}, we adopt an \emph{experimental design} perspective where we wish to identify a set of coordinates $\set{S}$ that minimizes the sum of conditional variances
\begin{equation*}
    \tr(\mat{A} - \mat{A}\langle \set{S}\rangle) = \sum_{j=1}^n \Var(z_j \mid \vec{z}(\set{S}))
\end{equation*}
for the unseen coordinates.
The greedy method builds the subset $\set{S}$ one pivot at a time, always choosing the \emph{maximum-variance} coordinate as pivot
\begin{equation*}
    s_{i+1} \in \argmax_{1\le j \le n} \Var(z_j \mid \vec{z}(\{s_1,\ldots,s_i\}) ).
\end{equation*}
The \RPCholesky method instead uses the variances as a \emph{sampling distribution}:
\begin{equation*}
    \prob \{s_{i+1} = j\} = \frac{\Var(z_j \mid \vec{z}(\{s_1,\ldots,s_i\}) )}{\sum_{\ell=1}^n \Var(z_\ell \mid \vec{z}(\set{S}))}.
\end{equation*}
Randomizing the procedure in this way balances the need to explore and exploit, and makes the method robust to bad instances like \cref{eq:greedy-bad-experimental-design} where the greedy method can be fooled by tiny differences in variances between coordinates.
\index{randomly pivoted Cholesky!connection to Gaussian distribution|)}\index{greedy pivoted Cholesky!failure modes|)}\index{Gaussian distribution!connection to Cholesky decomposition and Nystr\"om approximation|)}

\myparagraph{Software}
A high-performance implementation of a version of \RPCholesky is under development in the \RandLAPACK software project \cite{MDM+23a}.
At present it may be found at 
\actionbox{\url{https://github.com/BallisticLA/RandLAPACK/blob/main/RandLAPACK/comps/rl_rpchol.hh}}

\index{randomly pivoted Cholesky!numerical results|(}\index{greedy pivoted Cholesky!numerical results|(}\index{uniformly sampling for column Nystr\"om approximation!numerical results|(}\index{ridge leverage scores!numerical results|(}
\section{Experiments} \label{sec:rpc-experiments}

We will see several comparisons of \RPCholesky with alternative methods throughout this part of the thesis.
Here, we provide one initial set of data on \RPCholesky's performance, reproduced from \cite[Tab.~1]{CETW25}.
Here, we evaluate the relative trace error of \RPCholesky against other column selection methods for approximation of psd kernel matrices associated with 20 datasets.\index{kernel matrix}
(See \cref{ch:kernels-gaussian} for an introduction to kernel matrices.)
In addition to \RPCholesky, we test greedy selection and uniform and ridge leverage score (RLS) sampling.
We omit DPP sampling because of its high computational cost.
The Schmidt--Mirsky--Eckart--Young-optimal approximation error is shown for reference.
As the optimal approximation is typically not a column Nystr\"om approximation, it provides a lower limit on the best-possible approximation error for a column Nystr\"om approximation that \warn{typically cannot be attained}.

\begin{table}[t]
    \centering

    \caption[Comparison of uniform sampling, ridge leverage score sampling, greedy selection, and \RPCholesky for Nystr\"om approximation of kernel matrices]{Relative trace error of rank-$1000$ Nystr\"om approximation of psd kernel matrices by four alternative column selection methods: uniform and ridge leverage score (RLS) sampling, greedy selection, and \RPCholesky.
    Each trace error is computed as a median of ten trials.
    Error for the optimal rank-$1000$ approximation is shown for reference.
    Data is taken from \cite[Tab.~1]{CETW25}.}
    \label{tab:rpc-comparison}

    \begin{tabular}{rcccccc}
    \toprule
                                    & Uniform & RLS     & Greedy  & \textsc{RPChol}              & Optimal \\ \midrule
    \texttt{sensit\_vehicle}        & 1.57e-1 & 1.40e-1 & 2.07e-1 & \textbf{1.37e-1}             & 8.77e-2 \\
    \texttt{yolanda}                & 1.46e-1 & 1.39e-1 & 2.08e-1 & \textbf{1.36e-1}             & 8.41e-2 \\
    \texttt{YearPredictionMSD}      & 1.30e-1 & 1.20e-1 & 1.73e-1 & \textbf{1.16e-1}             & 6.79e-2 \\
    \texttt{w8a}                    & 1.42e-1 & 1.17e-1 & 1.91e-1 & \textbf{1.05e-1}             & 6.09e-2 \\
    \texttt{MNIST}                  & 1.21e-1 & 1.10e-1 & 1.67e-1 & \textbf{1.06e-1}             & 5.83e-2 \\
    \texttt{jannis}                 & 1.11e-1 & 1.10e-1 & 1.28e-1 & \textbf{1.09e-1}             & 5.45e-2 \\
    \texttt{HIGGS}                  & 6.73e-2 & 6.31e-2 & 8.51e-2 & \textbf{6.07e-2}             & 2.96e-2 \\
    \texttt{connect\_4}             & 5.81e-2 & 5.07e-2 & 6.48e-2 & \textbf{4.81e-2}             & 2.25e-2 \\
    \texttt{volkert}                & 5.35e-2 & 4.42e-2 & 5.92e-2 & \textbf{4.17e-2}             & 1.99e-2 \\
    \texttt{creditcard}             & 4.83e-2 & 3.71e-2 & 5.77e-2 & \textbf{3.01e-2}             & 1.31e-2 \\
    \texttt{Medical\_Appointment}   & 1.74e-2 & 1.43e-2 & 1.92e-2 & \textbf{1.29e-2}             & 4.59e-3 \\
    \texttt{sensorless}             & 1.20e-2 & 7.78e-3 & 8.70e-3 & \textbf{5.80e-3}             & 2.11e-3 \\
    \texttt{ACSIncome}              & 9.93e-3 & 5.55e-3 & 8.35e-3 & \textbf{4.02e-3}             & 1.27e-3 \\
    \texttt{Airlines\_DepDelay\_1M} & 4.19e-3 & 2.37e-3 & 2.64e-3 & \textbf{1.78e-3}             & 5.08e-4 \\
    \texttt{covtype\_binary}        & 9.10e-3 & 2.12e-3 & 1.41e-3 & \textbf{1.04e-3}             & 2.97e-4 \\
    \texttt{diamonds}               & 1.31e-3 & 2.40e-4 & 1.12e-4 & \textbf{5.85e-5}             & 1.30e-5 \\
    \texttt{hls4ml\_lhc\_jets\_hlf} & 3.78e-4 & 7.30e-5 & 6.68e-5 & \textbf{4.38e-5}             & 1.04e-5 \\
    \texttt{ijcnn1}                 & 3.91e-5 & 3.09e-5 & 2.86e-5 & \textbf{2.13e-5}             & 4.67e-6 \\
    \texttt{cod\_rna}               & 6.00e-4 & 1.36e-5 & 9.19e-6 & \textbf{5.05e-6}             & 9.86e-7 \\
    \texttt{COMET\_MC\_SAMPLE}      & 3.65e-3 & 2.44e-7 & 1.2e-10 & \textbf{4.3e-11}             & 3.5e-12 \\ \bottomrule
    \end{tabular}
\end{table}

Results are shown in \cref{tab:rpc-comparison}.
The performance of \RPCholesky is uniformly good, achieving the lowest trace error on all twenty examples.
We also see that \RPCholesky consistently achieves error close to the Schmidt--Mirsky--Eckart--Young optimal low-rank approximation.
These experiments confirm that \RPCholesky is among the best available approaches for constructing a column Nystr\"om approximation to a psd matrix.
We will see further examples of \RPCholesky's success throughout this part of the thesis.
\index{randomly pivoted Cholesky!numerical results|)}\index{greedy pivoted Cholesky!numerical results|)}\index{uniformly sampling for column Nystr\"om approximation!numerical results|)}\index{ridge leverage scores!numerical results|)}

\index{randomly pivoted Cholesky!theoretical results|(}
\section{Error analysis} \label{sec:rpc-analysis}

As we have seen and will continue to see, the \RPCholesky method consistently produces near-optimal low-rank approximations, often more accurate than competing methods by orders of magnitude. 
This excellent performance can be supported by theoretical analysis.
This section will present analysis of the approximation of \RPCholesky in the trace norm and a weak (but still useful!) result for the spectral norm.
It will then provide proofs of these results.

\index{randomly pivoted Cholesky!trace-norm bounds|(}
\subsection{Trace-norm bounds}

Our first bounds characterize the number of iterations needed to produce a good low-rank approximation to a psd matrix when measured using the trace norm.

\begin{theorem}[Randomly pivoted Cholesky: Trace norm] \label{thm:rpcholesky-trace}
    Let $r\ge 1$, $\varepsilon > 0$, and $\mat{A}\in\field^{n\times n}$ be a psd matrix.
    Introduce the relative error of the best rank-$r$ approximation
    \begin{equation*}
        \eta \coloneqq \tr(\mat{A} - \lowrank{\mat{A}}_r) / \tr(\mat{A}).
    \end{equation*}
    Randomly pivoted Cholesky\index{randomly pivoted Cholesky} produces an $(r,\varepsilon)$-approximation (\cref{def:r-eps-p}) to $\mat{A}$ as long as the number of steps satisfies
    \begin{equation} \label{eq:rpcholesky-trace-k}
        k \ge \frac{r}{\varepsilon} + r \log \left( \frac{1}{\varepsilon \eta} \right).
    \end{equation}
\end{theorem}

We observe that \RPCholesky achieves theoretical guarantees that are nearly optimal within the class of Nystr\"om approximations. (Recall from \cref{fact:optimal-nys-approx} that at least $k\ge r/\varepsilon$ columns are needed to produce an $(r,\varepsilon)$-approximation to a worst-case psd matrix.)

The dependence of the number of columns $k$ on the relative error $\eta$ is significantly improved for \RPCholesky over diagonal sampling (\cref{sec:uniform-diagonal-power}), which requires $\Theta(r/\eta)$ column accesses to produce an approximation comparable with the best rank-$r$ approximation for some input matrices.
\RPCholesky yields an \emph{exponential improvement} in the dependence on $1/\eta$.

The logarithmic factor $\log(1/\eta)$ is, for the purposes of practical computation, a modest constant.
Due to numerical errors, the relative error $\eta$ is effectively bounded from below by the \emph{unit roundoff}\index{unit roundoff} $u$, which captures the size of rounding errors ($u\approx 10^{-16}$ in double precision).
Thus, in double precision, $\log(1/\eta)\lessapprox 37$.

\begin{remark}[What if the relative error is small?] \label{rem:relative-error}
    One somewhat unappealing feature of this result is that when $\mat{A}$ is rank-$r$, the relative error is $\eta = 0$ and the right-hand side of \cref{eq:rpcholesky-trace-k} becomes infinite.
    Here, the bound \cref{eq:rpcholesky-trace-k} badly mischaracterizes the actual behavior of the \RPCholesky algorithm, which recovers a rank-$r$ matrix $\mat{A}$ \warn{with zero error} after precisely $k=\rank \mat{A}$ steps.
    Developing improved bounds for small $\eta$ is an open problem; see \cref{sec:better-rpcholesky-bounds} for discussion.
\end{remark}
\index{randomly pivoted Cholesky!trace-norm bounds|)}

\index{randomly pivoted Cholesky!spectral-norm bounds|(}
\index{spectral-norm error bounds for column Nystr\"om approximation|(}
\subsection{Weak spectral norm-type bounds}

It is natural to desire error bounds for \RPCholesky that hold in the spectral norm.
Perhaps, similar to \cref{fact:rls-spectral}, we could show that \RPCholesky achieves spectral-norm error $\lambda$ in roughly $d_{\mathrm{eff}}(\lambda)$ steps, where $d_{\mathrm{eff}}(\lambda)$ denotes the $\lambda$-effective dimension (\cref{def:ridge-psd}).
While such bounds are not yet known (see \cref{sec:better-rpcholesky-bounds}), we can establish bounds on the \warn{spectral norm of the expected error}:

\begin{theorem}[Randomly pivoted Cholesky: Spectral norm]\label{thm:rpcholesky-spectral}
    Fix parameters $b > 0$ and $\varepsilon > 0$.
    For any psd matrix $\mat{A}$, the $k$-step residual $\mat{A}^{(k)}$ of \RPCholesky satisfies 
    \begin{equation*}
        \norm{\, \expect\mat{A}^{(k)}} \le b + \varepsilon\tr(\mat{A} - \lowrank{\mat{A}}_r)
    \end{equation*}
    provided the number of steps $k$ satisfies
    \begin{equation*}
        k \ge \frac{1}{\varepsilon} + r \log \left( \frac{\norm{\mat{A}}}{b} \right).
    \end{equation*}
\end{theorem}

This result has a significant limitation in that the expectation occurs \warn{inside} the norm.
Pulling the expectation inside the norm can be done at a great cost.
Indeed, for any random psd matrix $\mat{X}$, we have the bound
\begin{equation} \label{eq:expected-norm-norm-expected}
    \norm{\expect \mat{X}} \le \expect \norm{\mat{X}} \le \expect \tr(\mat{X}) = \tr(\expect \mat{X}) \le n\cdot\norm{\expect \mat{X}}.
\end{equation}
The first inequality is Jensen's.
The factor $n$ in the upper bound $\expect\norm{\mat{X}} \le n \cdot \norm{\expect \mat{X}}$ is sharp, as demonstrated by the matrix $\mat{X} = \outprod{\evec_i}$ for $i\sim\Unif\{1,\ldots,n\}$.

The main use of \cref{thm:rpcholesky-spectral} is to bound quadratic forms $\vec{x}^*\mat{A}^{(k)}\vec{x}$ in the residual matrix:
\begin{equation*}
    \expect[\vec{x}^*\mat{A}^{(k)}\vec{x}] = \vec{x}^*[\expect \mat{A}^{(k)}]\vec{x} \le \norm{\vec{x}}^2 \cdot \big\|\expect \mat{A}^{(k)}\big\|.
\end{equation*}
We will use this type of bound to analyze quadrature methods in \cref{sec:quadrature}.

Ignoring the (significant) norm of expectation vs.\ expectation of norm issue, how good is \cref{thm:rpcholesky-spectral}?
To get some insight into this question, consider a matrix with rapidly polynomially decaying eigenvalues $\lambda_j(\mat{A}) = j^{-q}$ for a \warn{fixed} parameter $q > 1$.
The sum of tail eigenvalues is
\begin{equation*}
    \sum_{j=r+1}^n \lambda_j(\mat{A}) \le \sum_{j=r+1}^\infty j^{-q} = \order(r^{1-q}).
\end{equation*}
Setting $\varepsilon = \order(1/r)$ and $b = (r+1)^{-q}/2$, we obtain
\begin{equation*}
    \norm{\expect \mat{A}^{(k)}} \le \lambda_{r+1}(\mat{A}) = (r+1)^{-q} \quad \text{after } k = \order(r \log r) \text{ steps}.
\end{equation*}
\RPCholesky requires at most $\order(r\log r)$ steps for the spectral norm of the expected residual to drop below the spectral norm error of the best rank-$r$ approximation---not bad!

\index{randomly pivoted Cholesky!trace-norm bounds|(}\index{expected residual function $\mat{\Phi}(\mat{A}) = \mat{A} - \mat{A}^2/\tr(\mat{A})$|(}
\subsection{Proofs}

Let us begin with the proof of \cref{thm:rpcholesky-trace,thm:rpcholesky-spectral}.
As a first step, we introduce the \emph{expected residual function}
\begin{equation*}
    \mat{\Phi}(\mat{A}) \coloneqq \expect [\mat{A}^{(1)}]
\end{equation*}
which measures the expected value of the residual matrix $\mat{A}^{(1)}$ after one step of \RPCholesky applied to the input matrix $\mat{A}$.
By direct computation, we observe
\begin{equation*}
    \mat{\Phi}(\mat{A}) = \sum_{j=1}^n \left( \mat{A} - \frac{\outprod{\vec{a}_j}}{a_{jj}} \right) \cdot \prob \{ s_1 = j\} = \mat{A} - \sum_{j=1}^n \frac{\outprod{\vec{a}_j}}{a_{jj}} \frac{a_{jj}}{\tr(\mat{A})} = \mat{A} - \frac{\mat{A}^2}{\tr(\mat{A})}.
\end{equation*}
The map $\mat{\Phi}$ enjoys a number of properties.

\begin{proposition}[Expected residual function] \label{prop:expected-residual-function}
    The expected residual function $\mat{\Phi}$ satisfies the following properties:
    \begin{enumerate}[label=(\alph*)]
        \item \textbf{\textit{Unitarily covariant:}} For a unitary matrix $\mat{U} \in \field^{n\times n}$ and a psd matrix $\mat{A}$, $\mat{\Phi}(\mat{U}\mat{A}\mat{U}^*) = \mat{U}\mat{\Phi}(\mat{A})\mat{U}^*$. \label{item:Phi-unitarily-covariant}
        \item \textbf{\textit{Positive:}} For any psd matrix $\mat{A}$, the expected residual $\mat{\Phi}(\mat{A})$ is psd. \label{item:Phi-positive}
        \item \textbf{\textit{Concave:}} The map $\mat{\Phi}$ is concave with respect to the psd order. \label{item:Phi-concave}
        That is, for psd matrices $\mat{A},\mat{H}$, 
        \begin{equation}
            \mat{\Phi}(\theta \mat{A} + (1-\theta) \mat{H})\succeq  \theta \mat{\Phi}(\mat{A}) + (1-\theta) \mat{\Phi}(\mat{H}) \quad \text{for each }\theta \in[0,1]. \label{eq:residual-concave}
        \end{equation}
        \item \textbf{\textit{Monotone:}} The map $\mat{\Phi}$ is monotone with respect to the psd order. \label{item:Phi-monotone}
        That is, for psd matrices $\mat{A},\mat{H}$, 
        \begin{equation} \label{eq:residual-monotone}
            \mat{\Phi}(\mat{A}+\mat{H}) \succeq \mat{\Phi}(\mat{A}).
        \end{equation}
    \end{enumerate}
\end{proposition}

\begin{proof}
    Throughout this proof, let $\mat{A}$ and $\mat{H}$ denote psd matrices.

    The unitary covariance property \ref{item:Phi-unitarily-covariant} is immediate from the formula $\mat{\Phi}(\mat{A}) = \mat{A} - \mat{A}^2/\tr(\mat{A})$.

    The positivity property \ref{item:Phi-positive} follows by observing $\mat{\Phi}(\mat{A})$ has nonnegative eigenvalues 
    \begin{equation*}
        \lambda_i(\mat{\Phi}(\mat{A})) = \lambda_i(\mat{A}) - \frac{\lambda_i(\mat{A})^2}{\tr(\mat{A})} = \lambda_i(\mat{A}) \left( 1 - \frac{\lambda_i(\mat{A})}{\sum_{j=1}^n \lambda_j(\mat{A})} \right) \ge 0
    \end{equation*}
    for each $i=1,2,\ldots,n$.
    
    The concavity claim \ref{item:Phi-concave} follows from direct computation: For $\theta \in [0,1]$ and $\overline{\theta} = 1-\theta$, we have
    \begin{equation*}
        \mat{\Phi}(\theta \mat{A} + \overline{\theta} \mat{H}) - \theta \mat{\Phi}(\mat{A}) - \overline{\theta} \mat{\Phi}(\mat{H}) = \frac{\theta\overline{\theta}}{\theta \tr(\mat{A}) + \overline{\theta}\tr(\mat{H})} \left( \sqrt{\frac{\tr(\mat{H})}{\tr(\mat{A})}} \mat{A} - \sqrt{\frac{\tr(\mat{A})}{\tr(\mat{H})}}\mat{H}\right)^2,
    \end{equation*}
    which is manifestly psd.
    
    To show the monotonicity property \ref{item:Phi-monotone}, first observe that $\mat{\Phi}$ is positive homogeneous $\mat{\Phi}(\alpha \cdot \mat{A}) = \alpha \cdot \mat{\Phi}(\mat{A})$ for $\alpha \ge 0$.
    Consequently, by \cref{eq:residual-concave},
    \begin{equation*}
         \mat{\Phi}(\mat{A} + \mat{H}) = 2\mat{\Phi}(0.5\mat{A} + 0.5\mat{H}) \succeq \mat{\Phi}(\mat{A}) + \mat{\Phi}(\mat{H}) \succeq \mat{\Phi}(\mat{A}).
    \end{equation*}
    The last inequality holds since $\mat{\Phi}(\mat{H})$ is psd.
    We have established \cref{eq:residual-monotone}.
\end{proof}

Because of the monotonicity and concavity of the expected residual function, we can bound the $i$-step residual $\mat{A}^{(i)}$ by iterating the expected one-step residual function:

\begin{proposition}[Multistep residual] \label{prop:multistep-residual}
    Let $\mat{A}^{(i)}$ denote the residual of the matrix $\mat{A}$ after running $i$ steps of \RPCholesky.
    Then
    \begin{equation*}
        \expect[\mat{A}^{(i)}] \preceq \mat{\Phi}^i(\mat{A}).
    \end{equation*}
    Here, $\mat{\Phi}^i$ denotes the $i$-fold composition of $\mat{\Phi}$.
\end{proposition}

\begin{proof}
    First, apply an appropriate matrix version of Jensen's inequality \cite[Thm.~4.16]{Car10} conditionally on the $(i-1)$-step residual $\mat{A}^{(i-1)}$ to obtain
    \begin{equation*}
        \expect[\mat{A}^{(i)}] = \expect[\expect[\mat{A}^{(i)} \mid \mat{A}^{(i-1)}]] = \expect[\mat{\Phi}(\mat{A}^{(i-1)})] \preceq \mat{\Phi}(\expect [\mat{A}^{(i-1)}]).
    \end{equation*}
    Here, we used concavity \cref{eq:residual-concave} of $\mat{\Phi}$.
    Next, we iterate this inequality using the monotonicity property \cref{eq:residual-monotone}:
    \begin{equation*}
        \expect[\mat{A}^{(i)}] \preceq \mat{\Phi}(\expect [\mat{A}^{(i-1)}]) \preceq \mat{\Phi}^2(\expect [\mat{A}^{(i-2)}]) \preceq \cdots  \preceq \mat{\Phi}^{i-1}(\expect [\mat{A}^{(1)}]) \preceq \mat{\Phi}^i(\mat{A}).
    \end{equation*}
    The desired claim is established.
\end{proof}

\Cref{prop:multistep-residual} shows that the eigenvalues of $\expect[\mat{A}^{(k)}]$ are bounded above by the eigenvalues of the $k$-fold composition of the $\mat{\Phi}$ transformation to the matrix $\mat{A}$.
\index{expected residual function $\mat{\Phi}(\mat{A}) = \mat{A} - \mat{A}^2/\tr(\mat{A})$|)}
Introduce the \warn{vector} function $\vec{\phi} : \real^n_+ \to \real_n^+$
\begin{equation*}
    \vec{\phi}(\vec{\alpha}) = \vec{\alpha} - \frac{\vec{\alpha}^2}{\sum_{j=1}^n \alpha_j}.
\end{equation*}
Here, $\vec{\alpha}^2$ denotes the entrywise square.
The function $\vec{\phi}$ tracks the eigenvalues after applying the one-step residual function $\mat{\Phi}$:
\begin{equation*}
    \vec{\lambda}(\mat{\Phi}(\mat{A})) = \vec{\phi}(\vec{\lambda}(\mat{A})) \quad \text{for a psd matrix} \mat{A}.
\end{equation*}
Consequently,
\begin{equation} \label{eq:eigenvalue-recurrence-vec}
    \vec{\lambda}(\expect[\mat{A}^{(i)}]) \le \vec{\lambda}(\mat{\Phi}^i(\mat{A})) = \vec{\phi}^i(\vec{\lambda}(\mat{A})).
\end{equation}
In this display, $\vec{\phi}^k$ is the $i$-fold composition of $\vec{\phi}$, and the inequality holds entrywise.
The inequality in \cref{eq:eigenvalue-recurrence-vec} follows by \cref{prop:multistep-residual} and Weyl monotononicity principle \cite[Cor.~9.9]{Tro22}:
\begin{equation*}
    \mat{A}\preceq \mat{H} \implies \vec{\lambda}(\mat{A})\le\vec{\lambda}(\mat{H}).
\end{equation*}
The equality in \cref{eq:eigenvalue-recurrence-vec} is unitary covariance (\cref{prop:expected-residual-function}\ref{item:Phi-unitarily-covariant}).
We now provide a proof for \cref{thm:rpcholesky-trace,thm:rpcholesky-spectral}.

\begin{proof}[Proof of \cref{thm:rpcholesky-trace,thm:rpcholesky-spectral}]
    Introduce $\vec{\lambda}^{(i)} \coloneqq \vec{\phi}^i(\vec{\lambda}(\mat{A}))$ to track the right-hand side of \cref{eq:eigenvalue-recurrence-vec}.
    Since $\vec{\lambda}^{(i)}$ is defined as the iteration of the $\vec{\phi}$ map, it obeys the recurrence
    \begin{equation} \label{eq:eigenvalue-recurrence-vec-2}
        \vec{\lambda}^{(i+1)} = \vec{\phi}(\vec{\lambda}^{(i)}) = \vec{\lambda}^{(i)} - \frac{\big[\vec{\lambda}^{(i)}\big]^2}{\sum_{j=1}^n \lambda^{(i)}_j}. 
    \end{equation}
    We make two observations.
    First, we the vectors $\vec{\lambda}^{(i)}$ are entrywise decreasing: For each $i$, it holds that $\vec{\lambda}^{(i+1)}\le \vec{\lambda}^{(i)}$. 
    This observation is immediate from the recurrence \cref{eq:eigenvalue-recurrence-vec-2}.
    Second, the entries of each vector $\vec{\lambda}^{(i)}$ are sorted in nonincreasing order for each $i$.
    
    Let us verify this second observation by induction.
    For the base case $i=0$, note that $\vec{\lambda}^{(0)} = \vec{\lambda}(\mat{A})$ is sorted in nonincreasing order by definition of the $\vec{\lambda}(\cdot)$ function.
    Next, inductively suppose that $\vec{\lambda}^{(i-1)}$ is sorted in nonincreasing order.
    Then for $j\ge 1$,
    \begin{align*}
        \lambda^{(i)}_{j+1} - \lambda^{(i)}_j &= \left[\lambda^{(i-1)}_{j+1} - \lambda^{(i-1)}_j\right] - \frac{\left[\big(\lambda^{(i-1)}_{j+1}\big)^2 - \big(\lambda^{(i-1)}_j\big)^2\right]}{\sum_{\ell=1}^n \lambda^{(i-1)}_\ell} \\
        &=\left[\lambda^{(i-1)}_{j+1} - \lambda^{(i-1)}_j\right] \left(1 - \frac{\lambda^{(i-1)}_{j+1} + \lambda^{(i-1)}_j}{\sum_{\ell=1}^n \lambda^{(i-1)}_\ell} \right) \ge 0.
    \end{align*}
    This completes the inductive argument, showing that the entries of $\vec{\lambda}^{(i)}$ are sorted in nonincreasing order for every $i$.

    Having established some basic properties of the recurrence \cref{eq:eigenvalue-recurrence-vec-2}, we now shall reason about the average
    \begin{equation*}
        a^{(i)} \coloneqq \frac{1}{\ell} \sum_{j=1}^\ell \lambda^{(i)}_j
    \end{equation*}
    of the top $\ell$ eigenvalues of $\vec{\lambda}^{(i)}$.
    We take the parameter $\ell$ to be between $1$ and $r$, and we will primarily be interested in the edge cases $\ell = 1$ (top eigenvalue) and $\ell = r$ (top $r$ eigenvalues).
    By averaging the first $\ell$ entries of the recurrence \cref{eq:eigenvalue-recurrence-vec-2}, we obtain a recurrence for $a^{(i)}$:
    \begin{equation} \label{eq:a-recurrence}
        a^{(i+1)} = a^{(i)} - \frac{\ell^{-1}\sum_{j=1}^\ell \left[\lambda^{(i)}_j\right]^2}{ \sum_{j=1}^n \lambda^{(i)}_j}.
    \end{equation}
    We may bound the dominator of the second term as
    \begin{equation*}
        \sum_{j=1}^n \lambda^{(i)}_j = \sum_{j=1}^r \lambda_j^{(i)} + \sum_{j=r+1}^n \lambda_j^{(i)} \le r a^{(i)} + \sum_{j=r+1}^n \lambda_j(\mat{A}).
    \end{equation*}
    Here, we use the fact that the numbers $\lambda^{(i)}_j$ are nonincreasing in both $i$ and $j$ together with the initial condition $\lambda^{(0)}_j = \lambda_j(\mat{A})$.
    To bound the numerator of the second term of \cref{eq:a-recurrence}, we use Jensen's inequality:
    \begin{equation*}
        \frac{1}{\ell} \sum_{j=1}^\ell \left[\lambda^{(i)}_j\right]^2 \ge \left[\frac{1}{\ell} \sum_{j=1}^\ell \lambda^{(i)}_j\right]^2 = \left[a^{(i)}\right]^2.
    \end{equation*}
    Substituting the two previous displays into \cref{eq:a-recurrence} yields
    \begin{equation} \label{eq:a-recurrence-2}
        a^{(i+1)} \le a^{(i)} - \frac{ \left[a^{(i)}\right]^2}{ra^{(i)} + \sum_{j=r+1}^n \lambda_j}.
    \end{equation}

    To bound the solution to the recurrence \cref{eq:a-recurrence-2}, we compare to an ODE model
    \begin{equation} \label{eq:rpc-ode}
        \frac{\d}{\d t} x(t) = - \frac{x(t)^2}{rx(t) + \sum_{j=r+1}^n \lambda_j} \quad \text{with } x(0) = \frac{1}{\ell} \sum_{j=1}^\ell \lambda_j.
    \end{equation}
    For each $i$, we have $a^{(i)} \le x(i)$ because $x \mapsto -x^2/(rx + \sum_{j=r+1}^n \lambda_j)$ is decreasing on $\real_+$.
    Fixing a level $\gamma > 0$, we may use separation of variables to solve \cref{eq:rpc-ode} and obtain the time $t_\star$ at which $x(t_\star) = \gamma$:
    \begin{equation} \label{eq:rpc-ode-solve} 
        \begin{split}
        t_\star &= \int_\gamma^{\ell^{-1} \sum_{j=1}^\ell \lambda_j} \frac{rx + \sum_{j=r+1}^n \lambda_j}{x^2} \, \d x \\
        &= \frac{\sum_{j=r+1}^n \lambda_j}{\gamma} - \frac{\ell \sum_{j=r+1}^n\lambda_j}{\sum_{j=1}^\ell \lambda_j} + r\log \left(  \frac{\sum_{j=1}^\ell \lambda_j}{\ell \gamma} \right) \\
        &\le \frac{\tr(\mat{A} - \lowrank{\mat{A}}_r)}{\gamma} + r\log \left(  \frac{\sum_{j=1}^\ell \lambda_j}{\ell \gamma} \right).
        \end{split}
    \end{equation}
    In the last inequality, we recall that $\tr(\mat{A} - \lowrank{\mat{A}}_r) = \sum_{j=r+1}^n \lambda_j$.

    We can use \cref{eq:rpc-ode-solve} in two ways.
    First consider the case when $\ell = r$ and bound
    \begin{align*}
        \expect \tr(\mat{A}^{(k)}) &= \sum_{j=1}^r \lambda_j(\expect[\mat{A}^{(k)}]) + \sum_{j=r+1}^n \lambda_j(\expect[\mat{A}^{(k)}]) \\
        &\le \sum_{j=1}^r\lambda^{(k)}_j + \sum_{j=r+1}^n \lambda^{(k)}_j \le r a^{(i)} + \tr(\mat{A} - \lowrank{\mat{A}}_r).
    \end{align*}
    The first inequality uses \cref{eq:eigenvalue-recurrence-vec} and instates the definition of $\lambda_j^{(i)}$, and the second inequality uses the definition of $a^{(i)}$, the fact that $\lambda^{(i)}_j$ is nonincreasing in $i$, and the boundary condition $\lambda^{(0)}_j = \lambda_j(\mat{A})$.
    Apply \cref{eq:rpc-ode-solve} with $\ell = r$ and $\gamma \coloneqq \varepsilon/r \cdot \tr(\mat{A} - \lowrank{\mat{A}}_r)$ to conclude that $\expect \tr(\mat{A}^{(k)}) \le (1+\varepsilon)\tr(\mat{A} - \lowrank{\mat{A}}_r)$ when
    \begin{equation*}
        k \ge \frac{r}{\varepsilon} + r \log \left( \frac{\tr(\lowrank{\mat{A}}_r)}{\varepsilon \tr(\mat{A} - \lowrank{\mat{A}}_r)} \right).
    \end{equation*}
    This result is stronger than as stated in \cref{thm:rpcholesky-trace}.
    Second, consider the case $\ell = 1$ and set $\gamma \coloneqq b + \varepsilon \tr(\mat{A} - \lowrank{\mat{A}}_r)$.
    Then \cref{eq:eigenvalue-recurrence-vec,eq:rpc-ode-solve} imply that
    \begin{equation*}
        \norm{\expect(\mat{A}^{(k)})} = \lambda_1(\expect[\mat{A}^{(k)}]) \le \lambda_1^{(k)} \le b + \varepsilon \tr(\mat{A} - \lowrank{\mat{A}}_r)
    \end{equation*}
    when
    \begin{equation*}
        k \ge \frac{\tr(\mat{A} - \lowrank{\mat{A}}_r)}{b + \varepsilon \tr(\mat{A} - \lowrank{\mat{A}}_r)} + r \log \left( \frac{\norm{\mat{A}}}{b+\tr(\mat{A} - \lowrank{\mat{A}}_r)} \right).
    \end{equation*}
    This is stronger than the result of \cref{thm:rpcholesky-spectral}.
\end{proof}

\index{randomly pivoted Cholesky!theoretical results|)}
\index{randomly pivoted Cholesky!spectral-norm bounds|)}
\index{randomly pivoted Cholesky!trace-norm bounds|)}
\index{spectral-norm error bounds for column Nystr\"om approximation|)}

\index{randomly pivoted Cholesky!Gibbs variant|(}\index{Gibbs randomly pivoted Cholesky|(}
\section{Extension: Gibbs \RPCholesky} \label{sec:gibbs}

The uniform, \RPCholesky, and greedy pivoting strategies can be unified into a common framework.
The Gibbs \RPCholesky method \cite[\S2.3.4]{CETW25} selects a random pivot at each iteration according to the rule
\begin{equation*}
    s_{i+1} \sim \diag(\mat{A}^{(i)})^p.
\end{equation*}
We remind the reader that $s \sim \vec{w}$ denotes a random sample from the (unnormalized) weight vector $\vec{w}$.
The power $p \in (0,\infty)$ controls the level of ``greediness''.
The extreme case $p\downarrow 0$ corresponds to uniform sampling\index{uniform sampling for column Nystr\"om approximation} from the support of $\diag(\mat{A}^{(i)})$; the other extreme $p\uparrow + \infty$ selects the largest diagonal entries, with exact ties broken uniformly at random.
We recognize the cases $p=0$ and $p=\infty$ are variants of uniform and greedy selection.\index{greedy pivoted Cholesky}
\RPCholesky sits at the intermediate value $p = 1$.

As the name suggests, the Gibbs \RPCholesky method samples a pivot index from a Gibbs distribution $\prob \{ s_{i+1} = j \} = \exp(-\beta v_j)$ with ``energies'' $v_j = -\log a^{(i)}_{jj}$ and ``inverse-temperature'' $\beta = p$.
For this reason, the symbol $\beta$ has been used to refer to the power $p$ in previous literature \cite{Ste24,CETW25}.
Under the interpretation of $\beta = p$ as inverse-temperature, the uniform, greedy, and \RPCholesky strategies can be analogized to the three bowls of porridge in the story of the three bears---too hot, too cold, and just right.
The recent paper \cite{DPPL24} by Dong, Pan, Phan, and Lei has explored an alternate definition of the energies $v_j = -a^{(i)}_{jj}$, leading to sampling probabilities $s_{i+1}\sim \exp(\beta \diag(\mat{A}^{(i)}))$.

Our original paper \cite{CETW25} introduced the Gibbs \RPCholesky algorithm but did not provide any numerical experiments.
Stefan Steinerberger \cite{Ste24} took up the task of empirically evaluating the Gibbs \RPCholesky method.
Steineberger's experiments consider several examples, some of which benefit from higher $p$ and others which favor lower $p$.
However, except on specially constructed examples, the differences between different methods are not dramatic (particularly away from the extreme values $p \in \{0,\infty\}$).
Steinerberger's paper also includes theoretical results, including a theorem showing that a single step of the $p=2$ Gibbs \RPCholesky reduces the squared Frobenius norm of a psd matrix by at least a factor $1-1/n$.

My preliminary conclusion from Steinerberger and Dong et al.'s investigations is that, for non-pathological matrices, any sensible random pivoting selection \warn{that incorporates information from the diagonal of the residual matrix} should yield decent performance for low-rank approximation.
To achieve the best possible performance, one can treat the power $p$ (or inverse-temperature $\beta$ for the Dong et al.\ scheme) as a hyperparameter and test multiple $p$ or $\beta$ values.

\index{Gibbs randomly pivoted Cholesky!failure mode for $p>1$|(}
Even though Gibbs \RPCholesky seems to work reliably on most matrices with $p > 1$, significant failures can occur on synthetic, worst-case examples.
Fix $p > 1$, and consider the matrix
\begin{equation*}
    \mat{A} = \twobytwo{n^{1/p}\Id_{100}}{\mat{0}}{\mat{0}}{\outprod{\onevec_{n-100}}}.
\end{equation*}
The optimal rank-one approximation to this matrix is
\begin{equation*}
    \lowrank{\mat{A}}_1 = \twobytwo{\mat{0}}{\mat{0}}{\mat{0}}{\outprod{\onevec_{n-100}}},
\end{equation*}
and this approximation is computed by a partial pivoted Cholesky decomposition in one step if any pivot $s \in \{101,102,\ldots,n\}$ is selected.
For large $n$, this optimal rank-one approximation to this matrix has vanishingly small relative error
\begin{equation*}
    \frac{\tr(\mat{A} - \lowrank{\mat{A}}_1)}{\tr(\mat{A})} = \frac{100n^{1/p}}{n + 100(n^{1/p}-1)\rceil} = (100 + o(1)) n^{1/p-1}.
\end{equation*}
and \RPCholesky produces this optimal approximation in one step with near-certainty:
\begin{equation*}
    \prob\{ s \ge 101 \} = \frac{n-100}{n+100(n^{1/p}-1)} = 1 - o(1) \quad \text{when } s \sim \diag(\mat{A}).
\end{equation*}
With with Gibbs \RPCholesky with power $p$, a bad pivot is selected with roughly 99\% probability:
\begin{equation*}
    \prob\{ s \le 100 \} = \frac{100n}{(n-100)+100n} = \frac{100}{101} + o(1).
\end{equation*}
For large $n$, it takes Gibbs \RPCholesky with power $p$ roughly 50 steps to compute a good low-rank approximation to this matrix $\mat{A}$, where \RPCholesky ($p=1$) produces an excellent approximation with near certainty in a single step.
Given bad examples like these, we believe the power $p=1$ is a sensible default for general-purpose use, particular if only a single $p$ value is to be used. 
\index{randomly pivoted Cholesky!Gibbs variant|)}\index{Gibbs randomly pivoted Cholesky|)}\index{Gibbs randomly pivoted Cholesky!failure mode for $p>1$|)}

\index{randomly pivoted Cholesky!connection to determinantal point process sampling|(}\index{determinantal point process sampling!connection to randomly pivoted Cholesky|(}
\section{Connection to determinantal point processes} \label{sec:dpp-connections}

The \RPCholesky algorithm has a number of connections to (fixed-size) determinantal point processes (DPPs, \cref{def:dpp}).
This section reviews these connections.

\myparagraph{\RPCholesky as iterative 1-DPP sampling}
The first interpretation is the most trivial, but still yields some insight.
A single step of \RPCholesky performs one step of diagonal sampling on the current residual matrix $\mat{A}^{(i)}$, which coincides with the \emph{1-DPP distribution} (since the determinant of a $1\times 1$ matrix---i.e., a number---is just itself).
As a consequence of this interpretation, we can derive a single-step \RPCholesky error bound by invoking \cref{fact:k-dpp-nystrom}:
\begin{equation*}
    \expect \tr(\mat{A}^{(1)}) \le 2\tr(\mat{A} - \lowrank{\mat{A}}_1).
\end{equation*}

\myparagraph{\RPCholesky as iterative conditional DPP sampling}
Another, much more powerful, connection between \RPCholesky and DPPs interprets the \RPCholesky procedure as iteratively performing \emph{conditional} DPP sampling.

\begin{proposition}[\RPCholesky as conditional DPP sampling]
    Suppose we have run the \RPCholesky algorithm for $i$ steps on psd matrix $\mat{A}$, selecting pivots $s_1,\ldots,s_i$.
    The $(i+1)$st pivot $s_{i+1} \sim \diag(\mat{A}^{(i)})$ satisfies
    \begin{equation*}
        \prob \{ s_{i+1} = j \mid s_1,\ldots, s_i\} = \prob \big\{ \set{T} = \{s_1,\ldots,s_{i+1}\} \, \big| \, \set{T} \supseteq \{s_1,\ldots,s_i\} \big\},
    \end{equation*}
    where $\set{T} \sim \kDPP{i+1}(\mat{A})$.
\end{proposition}

We see that the $(i+1)$st step of \RPCholesky can be seen as sampling from an $(i+1)$-DPP, \emph{conditional} on the already-selected pivots belonging to that DPP.
This result follows directly from the definition of fixed-size DPPs and \cref{prop:det-pivot}.

\index{determinantal point process sampling!algorithms|(}
This connection between DPPs and \RPCholesky suggests an \emph{algorithm} for approximate $k$-DPP sampling.
First, run \RPCholesky for $k$ steps, producing pivots $s_1,\ldots,s_k$.
Then, for steps $t=0,1,\ldots,T$, select a random pivot index $j \sim \Unif\{1,\ldots,k\}$, evict pivot $s_j$ from the pivot set, and resample $s_j$ using a single step of the \RPCholesky procedure
\begin{equation*}
    s_j \sim \diag(\mat{A} - \mat{A}\langle s_1,\ldots,s_{j-1},s_{j+1},\ldots,s_k\rangle).
\end{equation*}
The reason behind this procedure's success is intuitive.
Initially, we generate the pivots $s_1,\ldots,s_k$ sequentially, and they do not follow the $k$-DPP distribution.
At each step, we freeze all but one pivot $s_j$ and resample $s_j$ from the $k$-DPP distribution conditional on the other pivots $\{s_i : i\ne j\}$.
By repeating this procedure enough time, the distribution of the pivots converges to the $k$-DPP distribution.

A sampling algorithm for a multivariate distribution of this type (freeze all but one coordinate and sample the conditional distribution) is known as a \emph{Gibbs Markov chain Monte Carlo (MCMC) sampler}.
This Gibbs DPP sampling procedure was proposed by Rezaei and Oveis Gharan \cite{RO19} to sample from a generalization of $k$-DPPs that can be defined on general, possibly \emph{continuous} state spaces.
Fortunately, ordinary $k$-DPPs on the finite set $\{1,\ldots,n\}$ are contained in this general setting.
Rezaei and Oveis Gharan's theoretical bounds show that the \RPCholesky-based Gibbs MCMC sampler  converges in at most $T = \order(k^5\log k)$ steps, and their numerical results suggest a more modest $T = \order(k^2)$ steps suffice.
To implement their algorithm efficiently, one should use Cholesky downdating techniques to handle the pivot evictions; we omit the details.

To raise a potentially provocative point, even Rezaei and Oveis Gharan's empirical results suggest that sampling a $k$-DPP using this algorithm requires \warn{quadratically} more work than executing \RPCholesky.
It is natural to ask: Is that work worth it? 
Is the output of \RPCholesky substantially worse than a sample from a $k$-DPP for practical purposes?
We leave these questions to the reader to ponder.
\index{determinantal point process sampling!algorithms|)}

\index{determinantal point process sampling!projection DPPs|(}
\myparagraph{Projection DPPs}
There is an important class of matrices for which \RPCholesky produces \warn{exact} samples from a $k$-DPP.
We make the following definition.

\begin{definition}[Projection DPP] \label{def:projection-dpp}
    Let $\mat{A}$ be an orthorpojector of \warn{rank exactly $k$}.
    Then the $k$-DPP $\set{S} \sim \kDPP{k}(\mat{A})$ is referred to as a \emph{projection DPP}.
\end{definition}

We emphasize that, in a projection DPP, the rank of the orthoprojector $\mat{A}$ must be equal to the size $k$ of the DPP $\set{S}$.
Projection DPPs are fundamental to both the theory and applications of DPPs.
Theoretically, every ($k$-)DPP can be realized as a mixture of projection DPPs \cite[p.~213]{HKPV06}, which gives rise to procedures for sampling from $k$-DPPs.\index{determinantal point process sampling!algorithms}

Remarkably, the \RPCholesky algorith generates exact samples from a projection DPP when applied to an orthoprojector.

\begin{proposition}[\RPCholesky samples projection DPPs] \label{prop:rpcholesky-projection-dpp}
    Let $\mat{A}$ be an orthoprojector of rank exactly $k$.
    The pivot set $\set{S}$ produced by $k$ steps of \RPCholesky applied to $\mat{A}$ is a sample from the projection DPP $\set{S} \sim \kDPP{k}(\mat{A})$.
\end{proposition}

This proposition appears in its essence in Gillenwater's thesis \cite[sec.~2.2.4]{Gil14a} and a very clear version of the \RPCholesky pseudocode (more-or-less the same as our \cref{prog:rpcholesky}) appears in the work of Poulson \cite{Pou20}.
We emphasize that Poulson only used the \RPCholesky procedure to sample from projection DPPs: In his work, the input matrix $\mat{A}$ is always an orthoprojector, he runs for exactly $\rank \mat{A}$ steps, and the output is the pivot set $\set{S}$.
The factor matrix $\mat{F}$, critical to use cases of \RPCholesky for low-rank approximation, is discarded in Poulson's work.
\index{randomly pivoted Cholesky|)}\index{randomly pivoted Cholesky!connection to determinantal point process sampling|)}\index{determinantal point process sampling!connection to randomly pivoted Cholesky|)}
\index{determinantal point process sampling!projection DPPs|)}

\chapter{Kernels and Gaussian processes} \label{ch:kernels-gaussian}

\epigraph{Although William of Occam first wielded his famous razor against the superfluous elaborations of his Scholastic predecessors, his principle of parsimony has since been incorporated into the methodology of experimental science in the following form: given two explanations of the data, all other things being equal, the simpler explanation is preferable.
This principle is very much alive today in the emerging science of machine learning, whose expressed goal is often to discover the simplest hypothesis that is consistent with the sample data.}{Anselm Blumer, Andrzej Ehrenfeucht, David Haussler, and Manfred K.\ Warmuth, \emph{Occam's razor} \cite{BEHW87}}

The twin theories of reproducing kernel Hilbert spaces and Gaussian processes are among the most beautiful topics in applicable mathematics, and they form the basis for a class of simple and user-friendly machine learning methods.
In this section, we review these theories and how they are used in machine learning.
In the following chapter, we explain how \RPCholesky can be used to accelerate these learning algorithms and provide computational experiments.

\myparagraph{Sources}
This section serves to provide an introduction to kernel methods and Gaussian processes.
It reflects my own approach to this material and is not explicitly based on any of particular published work, of mine or others.
The content of this section is influenced by Houman Owhadi's excellent class ``Stochastic processes and regression'', which I took in Winter of 2022.
Useful references on this material include the books \cite{SS02,RW05} and the survey \cite{KHSS18}.

\myparagraph{Outline}
\Cref{sec:rkhs} introduces the theory of reproducing kernel Hilbert spaces, and \cref{sec:kernel-interpolation} uses this formalism to derive the \emph{kernel interpolation} method for fitting data.
\Cref{sec:gp} describes the parallel theory of Gaussian processes and reformulates kernel interpolation as Gaussian process interpolation.
\Cref{sec:krr_gpr} concludes by discussing regularized kernel and Gaussian process methods.

\index{reproducing kernel Hilbert space|(}
\section{Reproducing kernel Hilbert spaces} \label{sec:rkhs}

\index{L2, function space@$\Ltwo$, function space|(}
One of the great annoyances in mathematical analysis is the fact that functions $u\in\Ltwo(\real^d)$ do not have definite values at any given point. 
Indeed, the space $\Ltwo(\real^d)$ is formally defined to consist of all square-integrable functions on $\real^d$, \warn{modulo almost everywhere equality} (with respect to the Lebesgue measure). 
So for any given $u\in\Ltwo(\real^d)$, there are elements of the equivalence class of $u$ taking every possible value at any given point $x$, since the singleton set $\{x\}$ has measure zero.

If we restrict ourselves to \emph{continuous} functions $u\in\set{C}(\real^d)\cap \Ltwo(\real^d)$, things are a bit better, since $u(\vec{x})$ has a defined value at each point $\vec{x} \in \real^d$.
However, even for continuous $\Ltwo$ functions, function values are not \emph{stable}: For any $\varepsilon>0$, any point $\vec{x}\in\real^d$, and any value $\alpha\in\field$, there exists a continuous perturbation $e$ of norm $\norm{e}_{\Ltwo(\real)^d}\le\varepsilon$ such that $(u+e)(\vec{x})=\alpha$.
The conclusion is that the value of any
$\Ltwo$ function at any point can be changed to any value by an arbitrarily small perturbation.

\index{reproducing kernel Hilbert space!definition|(}
These limitations of the $\Ltwo$ space motivates the concept of a reproducing kernel Hilbert space (RKHS), which may be informally defined as follows.\index{L2, function space@$\Ltwo$, function space|)}

\actionbox{An RKHS is a Hilbert space of functions that have definite values at every point, and these values are stable under small perturbations.}

Formally, we have the following definition:

\begin{definition}[Reproducing kernel Hilbert space]
    Let $\set{X}$ be a set.
    A \emph{reproducing kernel Hilbert space} (RKHS) on $\set{X}$ is a Hilbert space of functions $f:\set{X}\to\field$ over which the evaluation map $x\mapsto f(x)$ is a bounded linear functional for each $x$.
    That is, for each $x\in\set{X}$, there exists a prefactor $c(x)\in\real_+$ such that
    \begin{equation*}
        |f(x)|\le c(x) \norm{f}_\RKHS \quad \text{for each } f \in \RKHS.
    \end{equation*}
\end{definition}

The definition of the RKHS encompasses (and is equivalent to) the notion of stability of function values under small perturbations. 
Indeed, for $e \in \RKHS$ in an RKHS $\RKHS$,
\begin{equation*}
    |(f+e)(x) - f(x)| = |e(x)| \le c(x) \norm{e}_\RKHS.    
\end{equation*}
As we will see, the seemingly innocuous property of point evaluation being a bounded linear functional has surprising and far-reaching consequences.
\index{reproducing kernel Hilbert space!definition|)}

\index{reproducing kernel Hilbert space!examples|(}
Reproducing kernel Hilbert spaces (RKHSs) may initially appear exotic, but they encompass many of the fundamental function spaces in mathematical analysis.

\index{Sobolev space $\set{H}^1([0,1])$!is a reproducing kernel Hilbert space|(}
\begin{example}[A Sobolev space] \label{ex:sobolev}
    Consider, for instance, the \emph{Sobolev space} $\set{H}^1([0,1])$, which consists of all square-integrable functions with a square-integrable (weak) derivative.
    Its norm is
    \begin{equation*}
        \norm{f}_{\set{H}^1([0,1])}^2 = \int_0^1 \left( |f(x)|^2 + |f'(x)|^2
        \right) \, \d x.
    \end{equation*}
    This space is an RKHS.
    To verify this, first consider a \warn{continuous} function $f \in \set{H}^1([0,1]) \cap \set{C}([0,1])$.
    By the intermediate value theorem, there exists a value $x_\star$ at which $f$ achieves its mean value
    \begin{equation*}
        f(x_\star) = \int_0^1 f(x) \, \d x.
    \end{equation*}
    By the Cauchy--Schwarz inequality,
    \begin{equation*}
        |f(x_\star)|^2 = \left|\int_0^1 f(x) \, \d x\right|^2 \le \int_0^1 |f(x)|^2 \, \d x \cdot \int_0^1 1 \, \d x = \int_0^1 |f(x)|^2 \, \d x.
    \end{equation*}
    For any two values $x\le y$ in $[0,1]$,
    \begin{equation*}
        |f(y) - f(x)|^2 = \left| \int_x^y f'(a) \, \d a\right|^2 \le \int_x^y |f'(a)|^2 \, \d a \cdot  \int_x^y 1 \, \d a \le \int_0^1 |f'(a)|^2 \, \d a.
    \end{equation*}
    The same conclusion holds if $y\le x$.
    Employing the two previous displays, we conclude that for any $x \in [0,1]$, we have
    \begin{multline*}
        |f(x)|^2 \le 2|f(x_\star)|^2 + 2|f(x) - f(x_\star)|^2 \\ \le 2\int_0^1 \left(|f(x)|^2 + |f'(x)|^2\right) \, \d x = 2\norm{f}_{\set{H}^1([0,1])}^2.
    \end{multline*}
    We have established the RKHS property with $c(x) \le \sqrt{2}$ for all $x \in [0,1]$.
    This conclusion extends to all $f \in \set{H}^1([0,1])$ by density of $\set{C}([0,1])$ in $\set{H}^1([0,1])$.
    We conclude that $\set{H}^1([0,1])$ is an RKHS.
\end{example}
\index{reproducing kernel Hilbert space!examples|)}\index{Sobolev space $\set{H}^1([0,1])$!is a reproducing kernel Hilbert space|)}

This example is illustrative.
It suggests that functions in an RKHS should be smoother (in the sense of possessing well behaved derivatives) than typical functions in the space $\Ltwo(\mu)$.
This intuition will prove valuable.

\index{reproducing kernel|(}
\subsection{The reproducing kernel}

Whenever one has a linear functional over a Hilbert space, is is generally worth invoking the Riesz representation theorem to see if there are any interesting consequences.
In the cases of RKHSs, this is most certainly the case.

\index{reproducing kernel!construction|(}
For each point $x \in \set{X}$, the Riesz representation theorem furnishes a function $\kappa(\cdot,x) \in \RKHS$ for which
\begin{equation} \label{eq:reproducing-property}
    f(x) = \langle \kappa(\cdot,x), f \rangle_\RKHS.
\end{equation}
We have denoted this function by $\kappa(\cdot,x)$ to indicate that applying this construction at each $x \in \set{X}$ generates a parametric class of univariate functions $(\kappa(\cdot,x) : x \in \set{X})$.
Just as well, this class of univariate functions comprise a bivariate function $\kappa : \set{X} \times \set{X} \to \field$.
The property \cref{eq:reproducing-property} is referred to as the \emph{reproducing property} of the function $\kappa$.
Let me emphasize to the reader that, under the standing conventions, inner products are conjugate linear in their \warn{first} coordinate.
\index{reproducing kernel!construction|)}

\index{reproducing kernel!properties|(}
An especially interesting thing happens if we invoke the reproducing property \cref{eq:reproducing-property} for the special choice $f = \kappa(\cdot,x')$, which gives
\begin{equation} \label{eq:kernel-to-ip}
    \kappa(x,x') = \langle \kappa(\cdot,x), \kappa(\cdot,x') \rangle_\RKHS.
\end{equation}
We see that the bivariate function $\kappa$ tabulates the pairwise inner products of its univariate restrictions $\kappa(\cdot,x)$.
By Hermiticity of the $\RKHS$-inner product, we conclude that
\begin{equation*}
    \kappa(x,x') = \langle \kappa(\cdot,x), \kappa(\cdot,x') \rangle_\RKHS = \overline{\langle \kappa(\cdot,x'), \kappa(\cdot,x) \rangle_\RKHS} = \overline{\kappa(x',x)}.
\end{equation*}
The function $\kappa$ is Hermitian, i.e., conjugate symmetric.
Additionally, given any (finite) set of points $\set{D} \subseteq \set{X}$, the function matrix  $\kappa(\set{D},\set{D})$ is psd (see \cref{def:function_matrix}).
To see thus, observe that for any $\vec{c} \in \field^{\set{D}}$, 
\begin{multline*}
    \vec{c}^*\kappa(\set{D},\set{D})\vec{c} = \sum_{x,x' \in \set{D}} \overline{c_x} c_{x'} \, \kappa(x,x') = \sum_{x,x' \in \set{D}} \overline{c_x} c_{x'} \, \langle \kappa(\cdot,x), \kappa(\cdot,x') \rangle_\RKHS \\
    = \left\langle \sum_{x \in \set{D}} c_{x} \kappa(\cdot,x), \sum_{x' \in \set{D}} c_{x'} \kappa(\cdot,x') \right\rangle \ge 0.
\end{multline*}
The inequality is positive definiteness of the $\RKHS$-inner product.
We conclude that $\kappa(\set{D},\set{D})$ is psd.
Consequently, the function $\kappa$ is a positive-definite kernel function (\cref{def:kernel_function}).\index{kernel function}
Thus, a more appropriate name for the function $\kappa$ is the \emph{reproducing kernel} of $\RKHS$.
\index{reproducing kernel!properties|)}

\begin{figure}[t]
    \centering
    \includegraphics[width=0.9\linewidth]{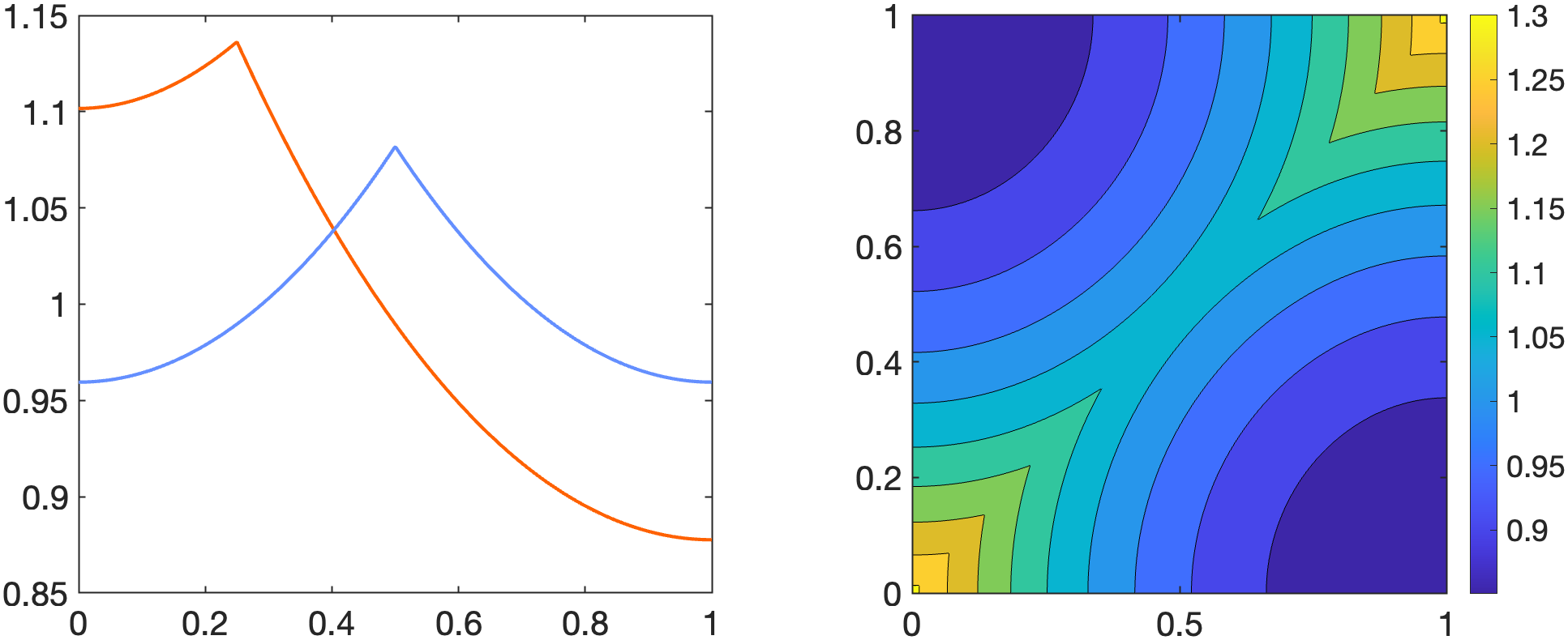}
    \caption[Visualizations of Sobolev kernel]{\emph{Left:} Univariate restrictions $\kappa(\cdot,a)$ of Sobolev kernel \cref{eq:sobolev-kernel} for $a=0.25$ and $a=0.5$.
    \emph{Right:} Contour plot of Sobolev kernel \cref{eq:sobolev-kernel}.}
    \label{fig:sobolev-kernel}
\end{figure}

\index{Sobolev space $\set{H}^1([0,1])$!reproducing kernel|(}
\begin{example}[Reproducing kernel for $\set{H}^1$] \label{ex:sobolev-rkhs-2}
    We return to the case of $\set{H}^1([0,1])$, introduced in \cref{ex:sobolev}.
    Since $\set{H}^1([0,1])$ is an RKHS, it has a reproducing kernel $\kappa : [0,1]^2 \to \field$ satisfying the property that
    \begin{equation} \label{eq:sobolev-reproducing}
        f(x) = \int_0^1 \left( \overline{\kappa(\cdot,a)} f(a) + \overline{\partial_a \kappa(\cdot,a)}f'(a)\right) \, \d a
    \end{equation}
    for any $f \in \set{H}^1([0,1])$.
    We can recover the value of $f$ at any point just by integrating against it and its derivative. 
    Neat!

    A formula for $\kappa$ can be found by integrating the second term of \cref{eq:sobolev-reproducing} by parts and solving an ODE boundary value problem \cite[\S2.11]{Sch07}.
    The resulting formula is
    \begin{equation} \label{eq:sobolev-kernel}
        \kappa(x,x') = \frac{\cosh(\min(x,x')) \cosh(1-\max(x,x'))}{\sinh(1)}.
    \end{equation}
    The property \cref{eq:sobolev-reproducing} can be verified analytically or checked for example functions $f$ and values $x$ using a symbolic computing environment like Mathematica.
    An illustration of this kernel is provided in \cref{fig:sobolev-kernel}.
\end{example}

\index{Sobolev space $\set{H}^1([0,1])$!reproducing kernel|)}

\index{reproducing kernel|)}

\index{reproducing kernel Hilbert space!existence|(}
\subsection{Spaces from kernels}

So far, we have defined the notion of an RKHS and showed that every RKHS $\RKHS$ has a reproducing kernel $\kappa$ (that is, a positive-definite kernel function satisfying the reproducing property \cref{eq:reproducing-property}).
Often in practice, we would like to in the reverse direction: Begin with a set $\set{X}$, endow it with a positive-definite kernel function $\kappa$, and obtain an RKHS $\RKHS$ for which $\kappa$ is the kernel.\index{Moore--Aronszajn theorem|(}
This construction is possible in view of the Moore--Aronszajn theorem \cite[p.~3.44]{Aro50}:

\begin{fact}[RKHSs from kernels] \label{fact:moore-aronszajn}
    Every positive-definite kernel function $\kappa$ on $\set{X}$ gives rise to a unique $\RKHS$ over which $\kappa$ is the reproducing kernel.
\end{fact}

\begin{proof}[Proof sketch]
    Notice that the span of univariate restrictions of the kernel $\kappa(\cdot,x)$ form an inner product space, where the inner product is defined by the relation \cref{eq:kernel-to-ip} and extended to linear combinations by sesquilinearity.
    This space may then be upgraded to a Hilbert space by taking the closure.
\end{proof}

\index{kernel function!examples|(}
This result speaks to how kernels are \emph{used} in practice to learn from data.
One begins with a dataset $\set{D}$ belonging to some ambient space $\set{X}$, typically $\set{X} = \real^d$.
Then one selects a kernel for this data, typically from some standard functional form, such as the square-exponential kernel
\begin{equation} \label{eq:square-exponential-kernel}
    \kappa_{\mathrm{se}}(\vec{x},\vec{x}' \mid \sigma) = \exp \left( - \frac{\uinorm{\vec{x} - \vec{x}'}^2}{2\sigma^2} \right)
\end{equation}
or the Laplace kernel
\begin{equation} \label{eq:laplace-kernel}
    \kappa_{\mathrm{Lap}}(\vec{x},\vec{x}' \mid \sigma) = \exp \left( - \frac{\uinorm{\vec{x} - \vec{x}'}}{\sigma} \right).
\end{equation}
Here, $\uinorm{\cdot}$ is a norm, typically the $\ell_1$ or $\ell_2$ norm.
Hyperparameters for these kernel families, such as the bandwidth $\sigma$ in \cref{eq:square-exponential-kernel,eq:laplace-kernel}, are picked either by ad-hoc techniques like the median heuristic \cite{GJK18} or by systematic procedures like cross validation \cite{LLJ+20}.
One then uses this selected kernel $\kappa$ to perform some data analysis task like interpolation (\cref{sec:kernel-interpolation}) or regression (\cref{sec:krr_gpr}).
\Cref{fact:moore-aronszajn} guarantees the existence of an RKHS $\RKHS$ for which $\kappa$ is the kernel, and the fine details of the Hilbert space $\RKHS$ are typically not needed.
\index{reproducing kernel Hilbert space!existence|)}\index{kernel function!examples|)}\index{Moore--Aronszajn theorem|)}

\index{reproducing kernel!interpretation as nonlinear inner product|(}
\subsection{The kernel as a nonlinear inner product}

There is another interpretation of kernel functions and RKHSs that can be useful.
The kernel $\kappa$ can be seen as defining a sort of nonlinear inner product on a general space $\set{X}$.
More precisely, there exists a \emph{feature map}\index{feature map} $\phi : \set{X} \to \Hilb$ mapping $\set{X}$ into a Hilbert space $\set{K}$  over which $\kappa$ coincides with the inner product
\begin{equation*}
    \kappa(x,x') = \langle \phi(x),\phi(x') \rangle_\Hilb.
\end{equation*}
The space $\Hilb$ is called the \emph{feature space}.\index{feature space}
The feature map\index{feature map} and the feature space\index{feature space} are not unique, but a natural choice of feature space\index{feature space} is given by the RKHS $\RKHS$ itself.
With this choice, the mapping $\phi : x \mapsto \kappa(\cdot,x)$ constitutes a feature map,\index{feature map} in view of \cref{eq:kernel-to-ip}.

\index{reproducing kernel Hilbert space|)}\index{reproducing kernel!interpretation as nonlinear inner product|)}

\index{kernel interpolation|(}
\section{Kernel interpolation} \label{sec:kernel-interpolation}

Having established the RKHS formalism, let us see how it can be used to learn from data.
Consider the task of learning a functional relation $g : \set{X} \to \field$ from input--output pairs $(x_1,y_{x_1}),\ldots,(x_n,y_{x_n}) \in \set{X} \times \field$.
For convenience, we may package the inputs into a multiset $\set{D} = \{ x_1,\ldots,x_n \}$ and the outputs into a vector $\vec{y} = (y_x : x \in \set{D}) \in \field^{\set{D}}$.
(In the case where $\set{D}$ has repeated elements $x_i = x_j$ for $i \ne j$, we abuse notation and permit $y_{x_i}$ and $y_{x_j}$ to be different values.)

If we assume that the output values $y_x$ are provided to us without noise, it makes sense to seek a model $g : \set{X} \to \field$ that \emph{interpolates} the data at the provided points:
\begin{equation*}
    g(x) = y_x \quad \text{for all } x \in \set{D}.
\end{equation*}
The interpolation condition may be written more concisely as $g(\set{D}) = \vec{y}$.

\index{minimum-norm solution to a linear system of equations!in infinite dimensions|(}\index{kernel interpolation!as minimum-norm interpolation|(}
There are infinitely many functions $g \in \RKHS$ interpolating the data.
Of these, it is natural to select the interpolating function $g$ of \emph{minimum norm}.
Indeed, if we think of the RKHS norm as a measure of the smoothness (as in \cref{ex:sobolev}) or complexity of a function, seeking the minimum-norm interpolant can be thought of as finding the interpolating function of minimum complexity (cf.\ the Occam's razor quote at the beginning of the chapter).

Finding the interpolating function of minimum norm is an optimization problem over the typically infinite dimensional space $\RKHS$:
\begin{equation} \label{eq:minimum-norm-interpolant}
    g = \argmin_{\substack{g \in \RKHS \\ g(\set{D}) = \vec{y}}} \, \norm{g}_\RKHS.
\end{equation}
Remarkably, this optimization problem has a closed form solution.

\begin{theorem}[Kernel interpolation: Solution formula] \label{thm:kernel-interpolation}
    Let $\set{D}$ be a finite set of points, and let $\vec{y} \in \field^{\set{D}}$ be output values.
    Assume $\kappa(\set{D},\set{D})$ is nonsingular.
    Then the optimization problem \cref{eq:minimum-norm-interpolant} has the following unique solution
    \begin{equation*}
        g = \sum_{x\in\set{D}} \kappa(\cdot,x) \beta_{x} \quad \text{where } \vec{\beta} = \kappa(\set{D},\set{D})^{-1}\vec{y}.
    \end{equation*}
\end{theorem}
\index{minimum-norm solution to a linear system of equations!in infinite dimensions|)}\index{kernel interpolation!as minimum-norm interpolation|)}

Theorems of this type, which show that an infinite-dimensional optimization problem over an RKHS has a finitely parametrized solution consisting of a linear combination of kernel functions, are known as \emph{representer theorems} \cite[\S4.2]{SS02}. 
For completeness, and because the proof is beautiful and revealing, we shall provide a proof of this representer theorem in the rest of this section.
Code for kernel interpolation is provided in \cref{prog:kernel_interp}.

\myprogram{Code to compute the kernel interpolant through data $\vec{y}$.}{}{kernel_interp}

\index{kernel method!implementation|(}
\begin{remark}[Interface for kernel methods]
    For the programs in this thesis, we use a common interface to implement kernel methods.
    Sets of $n$ data points $\set{D} \subseteq \field^d$ are collected as rows of an $n\times d$ matrix \texttt{D}.
    A univariate functions \texttt{g} can evaluated on a set of inputs $g(\set{D})$ as \texttt{g(D)}.
    Similarly, the evaluation $\kappa(\set{D},\set{E})$ of a bivariate function $\kappa$ can be evaluated as \texttt{kappa(D,E)}.
\end{remark}
\index{kernel method!implementation|)}

\index{minimum-norm solution to a linear system of equations!in finite dimensions|(}
\subsection{Underdetermined systems of linear equations}

To motivate the proof of \cref{thm:kernel-interpolation}, we shall begin with a review of the theory of underdetermined systems of linear equations in finite dimensions.

The linear least-squares problem 
\begin{equation*}
    \vec{x} = \argmin_{\vec{x} \in \field^n}\, \norm{\mat{B} \vec{x} - \vec{c}} \quad \text{for } \mat{B} \in\field^{m\times n}, \: \vec{c} \in \field^m
\end{equation*}
is well-known, as are the normal equations\index{normal equations!for linear least squares} characterizing its solution:
\begin{equation*}
    \mat{B}^*\mat{B}\vec{x} = \mat{B}^*\vec{c}.
\end{equation*}
There is a parallel theory, usually not covered in introductory linear algebra classes, for the minimum norm solution to an underdetermined system of equations
\begin{equation} \label{eq:underdetermined-linsys}
    \vec{x}_\star = \argmin_{\substack{\vec{x} \in \field^m \\ \mat{B}^*\vec{x} = \vec{c}}} \,\norm{\vec{x}} \quad \text{for } \mat{B} \in \field^{m\times n}, \vec{c} \in \field^n. 
\end{equation}
Here is the result.

\begin{theorem}[Minimum-norm solution to an underdetermined system]  \label{thm:underdetermined-linsys}
    Assume $\mat{B}$ has full column rank.
    Then $\vec{x} = \mat{B}(\mat{B}^*\mat{B})^{-1}\vec{c} = (\mat{B}^\dagger)^*\vec{c}$ is the unique solution to \cref{eq:underdetermined-linsys}.\index{normal equations!for underdetermined linear systems}
\end{theorem}

\begin{proof}
    Decompose any solution $\vec{x}$ to $\mat{B}^*\vec{x} = \vec{c}$ as an orthogonal sum $\vec{x} = \vec{x}_\star + \vec{x}_\perp$ of a component $\vec{x}_\star \in \range(\mat{B})$ and a component $\vec{x}_\perp \in \range(\mat{B})^\perp$.
    The orthogonal complement of $\range(\mat{B})$ is the nullspace of $\mat{B}^*$, so 
    \begin{equation*}
        \vec{c} = \mat{B}^*\vec{x} = \mat{B}^*(\vec{x}_\star + \vec{x}_\perp) = \mat{B}^*\vec{x}_\star + \vec{0} = \mat{B}^*\vec{x}_\star.
    \end{equation*}
    Consequently, we see that $\vec{x}_\star$ is also a solution of $\mat{B}^*\vec{x} = \vec{c}$.
    Moreover, the norm of $\vec{x}_\star$ is no larger than the norm of $\vec{x}$:
    \begin{equation} \label{eq:min-norm-derivation}
        \norm{\vec{x}}^2 = \norm{\vec{x}_\star}^2 + \norm{\vec{x}_\perp}^2 \ge \norm{\vec{x}_\star}^2.
    \end{equation}
    Equality holds in \cref{eq:min-norm-derivation} if and only if $\vec{x}_\perp = \vec{0}$.
    Since $\vec{x}_\star \in \range(\mat{B})$, we may write $\vec{x}_\star = \mat{B}\vec{y}_\star$, so that
    \begin{equation*}
        \vec{c} = \mat{B}^*\vec{x}_\star = (\mat{B}^*\mat{B})\vec{y}_\star.
    \end{equation*}
    Since $\mat{B}$ is full-rank, the matrix $\mat{B}^*\mat{B}$ is invertible, so $\vec{y}_\star = (\mat{B}^*\mat{B})^{-1} \vec{c}$.
    We conclude that $\vec{x}_\star = \mat{B}(\mat{B}^*\mat{B})^{-1} \vec{c}$ is the unique minimal-norm solution to $\mat{B}^*\vec{x} = \vec{c}$.
\end{proof}
\index{minimum-norm solution to a linear system of equations!in finite dimensions|)}

\index{minimum-norm solution to a linear system of equations!in infinite dimensions|(}
\subsection{Underdetermined systems of linear equations in a Hilbert space}

With the appropriate apparatus, the derivation for the minimum norm solution of a finite-dimensional system of linear equations translates effortlessly to infinite dimensions.
We make the following definition \cite{Ste98a,TT15}:

\begin{definition}[Quasimatrix]
    Let $\Hilb$ be a Hilbert space.
    An $\Hilb\times n$ \emph{quasimatrix}\index{quasimatrix} $F$ is a collection $F = (f_i : 1\le i \le n)$ of elements $f_i \in \Hilb$.
    For a matrix $\mat{H} \in \field^{n\times p}$, the product $F \mat{H}$ is
    \begin{equation*}
        F\mat{H} = \left( \sum_{i=1}^n f_i h_{ij} : 1\le j \le p \right).
    \end{equation*}
    The matrix--vector product $F\vec{h}$ is defined analogously.
    For an $\Hilb \times m$ quasimatrix\index{quasimatrix} $G$, the product $G^*F$ is 
    \begin{equation*}
        G^*F = \left(\langle g_i,f_j \rangle_\Hilb : 1\le i \le m, 1\le j \le n\right) \in \field^{m \times n}.
    \end{equation*}
    Similarly, $G^*f = (\langle g_i,f \rangle : 1\le i \le m)$ for $f \in \Hilb$.
\end{definition}

With this definition, \cref{thm:underdetermined-linsys} holds for underdetermined systems of equations over a Hilbert space with the same proof.

\begin{theorem}[Minimum-norm solution to an underdetermined system over a Hilbert space] \label{thm:underdetermined-linsys-hilb}
    Let $B$ be an $\Hilb \times n$ quasimatrix\index{quasimatrix} and let $\vec{c} \in \field^n$.
    Assume that $B$ has linearly dependent ``columns'' $(b_i : 1\le i \le n)$ or, equivalently, that $B^*B$ has full rank.
    Then the system of linear systems $B^*x = \vec{c}$ has a unique solution $x \in \Hilb$ of minimum norm.
    This solution is $x = B(B^*B)^{-1}\vec{c}$.\index{normal equations!for underdetermined linear systems}
\end{theorem}

\index{kernel interpolation!as minimum-norm interpolation|(}
This result immediately implies the representer theorem for kernel interpolation (\cref{thm:kernel-interpolation}).
To see this, introduce the quasimatrix\index{quasimatrix}
\begin{equation*}
    \kappa(\cdot,\set{D}) \coloneqq (\kappa(\cdot,x) : x\in\set{D}) \in \field^{\RKHS \times \set{D}}
\end{equation*}
and its adjoint
\begin{equation} \label{eq:kernel-fun-rows}
    \kappa(\set{D},\cdot) \coloneqq \kappa(\cdot,\set{D})^*.
\end{equation}

\begin{proof}[Proof of \cref{thm:kernel-interpolation}]
    By the reproducing property and \cref{eq:kernel-fun-rows},
    \begin{equation*}
        g(\set{D}) = (g(x) : x\in \set{D}) = (\langle \kappa(\cdot,x),g\rangle_\RKHS:x\in\set{D}) = \kappa(\set{D},\cdot)g.
    \end{equation*}
    Consequently, the minimum-norm interpolation problem \cref{eq:minimum-norm-interpolant} is equivalent to finding the minimum norm solution of the linear system $\kappa(\set{D},\cdot) g = \vec{y}$.
    By \cref{thm:underdetermined-linsys-hilb}, the solution is
    \begin{equation*}
        g = \kappa(\cdot,\set{D})[\kappa(\set{D},\cdot)\kappa(\cdot,\set{D})]^{-1}\vec{y}.
    \end{equation*}
    To complete the proof, observe that
    \begin{equation} \label{eq:kDkD}
        \begin{split}
        \kappa(\set{D},\cdot)\kappa(\cdot,\set{D}) &= \left( \langle \kappa(x,\cdot), \kappa(x',\cdot) \rangle_\RKHS : x,x'\in\set{D} \right) \\
        &= \left( \kappa(x,x') : x,x'\in\set{D} \right) = \kappa(\set{D},\set{D}).
        \end{split}
    \end{equation}
    We conclude that $g = \kappa(\cdot,\set{D})\kappa(\set{D},\set{D})^{-1}\vec{y}$, as desired.
\end{proof}
\index{kernel interpolation!as minimum-norm interpolation|)}

\index{kernel interpolation!when kernel matrix is rank-deficient|(}
\begin{remark}[Rank-deficient kernel matrix]
    If the kernel matrix $\kappa$ is rank-deficient, then we have the following alternative version of \cref{thm:kernel-interpolation}:
    \begin{theorem}[Kernel interpolation: General case for solution formula] \label{thm:kernel-interpolation-general}
        Let $\set{D}$ be a finite set of points and let $\vec{y} \in \field^{\set{D}}$ be output values.
        Then there is a unique minimum norm solution $g \in \RKHS$ to the least-squares problem 
        \begin{equation*}
            \min_{g \in \RKHS} \, \norm{\vec{y} - g(\set{D})},
        \end{equation*}
        and it is given by 
        \begin{equation} \label{eq:kernel-interp-equiv}
            g = \kappa(\cdot,\set{D}) \kappa(\set{D},\set{D})^\dagger \vec{y}.
        \end{equation}
        
    \end{theorem}
    Going forward, we shall refer to the function $g$ furnished by this theorem as the kernel interpolant, even when $\kappa(\set{D},\set{D})$ is rank-deficient.
    In the rank-deficient case, the function $g$ will not interpolate the data unless $\vec{y} \in \range(\kappa(\set{D},\set{D}))$.
\end{remark}
\index{minimum-norm solution to a linear system of equations!in infinite dimensions|)}\index{kernel interpolation!when kernel matrix is rank-deficient|)}

\index{kernel interpolation!connection to Nystr\"om approximation|(}
\subsection{Kernel interpolation, Nystr\"om approximation, and error bounds}

Kernel interpolation may be applied to any output data $\vec{y}\in \field^{\set{D}}$, but something special happens when the data $\vec{y} = f(\set{D})$ are outputs of a function $f \in \RKHS$.
In this case, the kernel interpolant, now denoted $\hat{f}_{\set{D}}$, takes the form
\begin{equation} \label{eq:kernel-interp-nystrom}
    \hat{f}_{\set{D}} = \kappa(\cdot,\set{D}) \kappa(\set{D},\set{D})^\dagger f(\set{D}) = \kappa(\cdot,\set{D}) \kappa(\set{D},\set{D})^\dagger \kappa(\set{D},\cdot) f.
\end{equation}
We have used the version \cref{eq:kernel-interp-equiv} of the conclusion of \cref{thm:kernel-interpolation}.
Observe the formal similarity of the right-hand side of \cref{eq:kernel-interp-nystrom} to the column Nystr\"om approximation \cref{eq:column-nystrom} of a psd matrix.
Drawing out this formal connection between kernel interpolation and Nystr\"om approximation will give us a powerful theoretical framework. 
We will develop and use this theory here to prove error bounds, and this theory will play a major role when we extend \RPCholesky to infinite dimensions in \cref{ch:infinite}.\index{infinite-dimensional randomly pivoted Cholesky}

\index{Nystr\"om approximation!of a kernel function|(}\index{kernel function!Nystr\"om approximation|(}
\begin{definition}[Nystr\"om approximation of a kernel function] \label{def:nystrom-kernel}
    Let $\kappa$ be a positive-definite kernel function on $\set{X}$ and let $\set{D}\subseteq \set{X}$ be a finite subset.
    The \emph{Nystr\"om approximation} to $\kappa$ induced by $\set{D}$ is
    \begin{equation*}
        \hat{\kappa}_{\set{D}}(x,x') \coloneqq \kappa(x,\set{D})\kappa(\set{D},\set{D})^\dagger \kappa(\set{D},x').
    \end{equation*}
    The \emph{residual kernel} is $\kappa_{\set{D}} \coloneqq \kappa - \hat{\kappa}_{\set{D}}$.
\end{definition}

Just as a column Nystr\"om approximation to a psd matrix and its residual are psd matrices, the Nystr\"om approximation and its residual are both positive-definite kernels.
The reader is reminded that a ``positive-definite'' kernel function is more analogous to a positive \warn{semi}definite matrix.

\begin{proposition}[Nystr\"om approximation of kernels] \label{prop:nystrom-kernel}
    Let $\kappa$ be a positive-definite kernel function on $\set{X}$, and let $\set{D} \subseteq \set{X}$ be subset.
    The Nystr\"om approximate kernel $\hat{\kappa}_{\set{D}}$ and its residual $\kappa_{\set{D}} = \kappa - \hat{\kappa}_{\set{D}}$ are both positive-definite kernels on $\set{X}$.
\end{proposition}

\begin{proof}
    Let $\set{E} \subseteq \set{X}$ and form the matrix $\mat{A} \coloneqq \kappa(\set{D}\cup\set{E},\set{D}\cup\set{E}) \in \field^{(\set{D}\cup\set{E})\times (\set{D}\cup\set{E})}$.
    Introduce the matrix Nystr\"om approximation $\Ahat \coloneqq \mat{A}\langle \set{D}\rangle$.
    One readily sees that
    \begin{equation*}
        \hat{\kappa}_{\set{D}}(\set{E},\set{E}) = \Ahat(\set{E},\set{E}) \quad \text{and} \quad \kappa_{\set{D}}(\set{E},\set{E}) = (\mat{A}-\Ahat)(\set{E},\set{E}).
    \end{equation*}
    By \cref{prop:nystrom-properties}, these matrices are both psd.
    We conclude that $\hat{\kappa}_{\set{D}}$ and $\kappa_{\set{D}}$ are positive-definite kernels.
\end{proof}

The following result gives a natural representation for the kernel interpolant.

\begin{lemma}[Reproducing representation for kernel interpolant] \label{lem:reproducing-kernel-interpolation}
    Let $\RKHS$ be an RKHS on $\set{X}$, $\set{D} \subseteq \set{X}$ be a finite set of points, and $f \in \RKHS$.
    The kernel interpolant $\hat{f}_{\set{D}}$ through $f$ at $\set{D}$ admits the reproducing representation
    \begin{equation*}
        \hat{f}_{\set{D}}(x) = \langle \hat{\kappa}_{\set{D}}(\cdot,x), f\rangle_\RKHS \quad \text{for each $x \in \set{X}$.}
    \end{equation*}
\end{lemma}

\begin{proof}
    Using quasimatrices\index{quasimatrix}, this result is easy.
    We compute
    \begin{equation*}
        \langle \hat{\kappa}_{\set{D}}(\cdot,x), f\rangle = (\kappa(\cdot,\set{D})\kappa(\set{D},\set{D})^\dagger \kappa(\set{D},x))^* f = \kappa(x,\set{D}) \kappa(\set{D},\set{D})^\dagger \kappa(\set{D},\cdot) f.
    \end{equation*}
    By \cref{eq:kernel-interp-nystrom}, this expression equals $\hat{f}_{\set{D}}(x)$.    
\end{proof}

We are tantalizingly close to proving error bounds for kernel interpolation.
We need one final result.

\begin{lemma}[Self inner product of residual and Nystr\"om kernels] \label{lem:residual-kernel-self}
    With the setting of \cref{lem:reproducing-kernel-interpolation}, we have
    \begin{equation*}
        \langle \hat{\kappa}_{\set{D}}(\cdot,x),\hat{\kappa}_{\set{D}}(\cdot,x)\rangle_\RKHS = \hat{\kappa}_{\set{D}}(x,x) \quad \text{and} \quad \langle \kappa_{\set{D}}(\cdot,x),\kappa_{\set{D}}(\cdot,x)\rangle_\RKHS = \kappa_{\set{D}}(x,x).
    \end{equation*}
\end{lemma} 

\begin{proof}
    To prove the first identity, we employ the quasimatrix\index{quasimatrix} formalism:
    \begin{align*}
        \langle \hat{\kappa}_{\set{D}}(\cdot,x),\hat{\kappa}_{\set{D}}(\cdot,x)\rangle_\RKHS &= (\kappa(\cdot,\set{D})\kappa(\set{D},\set{D})^\dagger \kappa(\set{D},x))^* (\kappa(\cdot,\set{D})\kappa(\set{D},\set{D})^\dagger \kappa(\set{D},x)) \\
        &= \kappa(x,\set{D})\kappa(\set{D},\set{D})^\dagger\kappa(\set{D},\cdot)\kappa(\cdot,\set{D})\kappa(\set{D},\set{D})^\dagger \kappa(\set{D},x) \\
        &= \kappa(x,\set{D})\kappa(\set{D},\set{D})^\dagger\kappa(\set{D},\set{D})\kappa(\set{D},\set{D})^\dagger \kappa(\set{D},x) \\
        &= \kappa(x,\set{D})\kappa(\set{D},\set{D})^\dagger \kappa(\set{D},x) = \hat{\kappa}_{\set{D}}(x,x).
    \end{align*}
    The third line is \cref{eq:kDkD}, and the fourth line is the identity $\mat{M}^\dagger\mat{M}\mat{M}^\dagger = \mat{M}^\dagger$ for the pseudoinverse.
    The second identity follows from the first:
    \begin{align*}
        \langle \kappa_{\set{D}}(\cdot,x),\kappa_{\set{D}}(\cdot,x)\rangle_\RKHS &= \langle \kappa(\cdot,x) - \hat{\kappa}_{\set{D}}(\cdot,x),\kappa(\cdot,x) - \hat{\kappa}_{\set{D}}(\cdot,x)\rangle_\RKHS \\
        &= \kappa(x,x) - 2\hat{\kappa}_{\set{D}}(x,x) + \langle \hat{\kappa}_{\set{D}}(\cdot,x),\hat{\kappa}_{\set{D}}(\cdot,x)\rangle_\RKHS = \kappa_{\set{D}}(x,x).
    \end{align*}
    The final equality is $\langle \hat{\kappa}_{\set{D}}(\cdot,x),\hat{\kappa}_{\set{D}}(\cdot,x)\rangle_\RKHS = \hat{\kappa}_{\set{D}}(x,x)$ and the definition of $\kappa_{\set{D}}$.
\end{proof}

\index{kernel interpolation!pointwise error bounds|(}
With \cref{lem:reproducing-kernel-interpolation,lem:residual-kernel-self} in hand, error bounds for kernel interpolation follow effortlessly.

\begin{theorem}[Kernel interpolation: Pointwise error] \label{thm:kernel-interpolation-error}
    Instate the setting and assumptions of \cref{lem:reproducing-kernel-interpolation}.
    Then the kernel interpolant satisfies the error bound
    \begin{equation} \label{eq:kernel-interpolation-pointwise-bound}
        |f(x) - \hat{f}_{\set{D}}(x)|^2 \le \kappa_{\set{D}}(x,x) \cdot \norm{f}_\RKHS^2.
    \end{equation}
    That is, the squared error at $x$ is bounded in terms of the diagonal entry of the residual kernel $\kappa_{\set{D}}$ at $x$.
    At each $x$, the bound is attained.
\end{theorem}

\begin{proof}
    By \cref{lem:reproducing-kernel-interpolation}, the error at $x$ is
    \begin{equation*}
        f(x) - \hat{f}_{\set{D}}(x) = \langle f, \kappa_{\set{D}}(\cdot,x)\rangle_\RKHS.
    \end{equation*}
    Take absolute values, and bound via Cauchy--Schwarz:
    \begin{equation*}
        |f(x) - \hat{f}_{\set{D}}(x)|^2 \le \langle \kappa_{\set{D}}(\cdot,x), \kappa_{\set{D}}(\cdot,x) \rangle_\RKHS \cdot \norm{f}_\RKHS^2 = \kappa_{\set{D}}(x,x) \cdot \norm{f}_\RKHS^2.
    \end{equation*}
    The final equality is \cref{lem:residual-kernel-self}.
    The bound is attained by setting $f = \kappa_{\set{D}}(\cdot,x)$.
\end{proof}

\index{Nystr\"om approximation!of a kernel function|)}\index{kernel function!Nystr\"om approximation|)}\index{kernel interpolation|)}\index{kernel interpolation!connection to Nystr\"om approximation|)}\index{kernel interpolation!pointwise error bounds|)}

\section{Gaussian processes} \label{sec:gp}

We now discuss the theory of Gaussian processes and their applications to learning from data, which beautifully parallels the theory of RKHSs.

\index{Gaussian process!definition|(}
\begin{definition}[Gaussian process]
    Let $\set{X}$ be a set.
    A \emph{Gaussian process} $g$ on $\set{X}$ with \emph{mean function} $m : \set{X} \to \field$ and \emph{covariance function}\index{covariance function} $\kappa : \set{X}\times \set{X} \to \field$, written $g \sim \GP(m,\kappa)$, is a random function $g$ such that, for every finite subset $\set{D} \subseteq \set{X}$, the function values $g(\set{D})$ obey a Gaussian distribution
    \begin{equation*}
        g(\set{D}) \sim \Normal_\field(m(\set{D}),\kappa(\set{D},\set{D})).
    \end{equation*}
    If $m$ is identically zero, we say that $g$ is \emph{centered} and write $g \sim \GP(\kappa)$.
\end{definition}

Naturally, the covariance function $\kappa$ is required to have the property that $\kappa(\set{D},\set{D})$ is psd for every $\set{D}$; that is, $\kappa$ must be a positive-definite kernel.
The term Gaussian process is frequently abbreviated ``GP''.
\index{Gaussian process!definition|)}

\index{Gaussian process!existence|(}
\begin{remark}[Existence of GPs]
The existence of a GP with any specified mean function and any positive-definite kernel as covariance function is ensured by the Kolmogorov extension theorem \cite[Thm.~14.36]{Kle13}.
Unfortunately, this basic existence result is not sufficient to ensure that statements like ``$g \sim \GP(m,\kappa)$ is continuous'' can be assigned definite probabilities (i.e., the event ``$g$ is continuous'' may not be measurable).
Fortunately, in most cases of practical significance, there are results that are powerful enough both to ensure existence of a GP and that global properties like continuity and differentiability have definite probabilities.
We will not discuss these more nuanced issues of GP theory in this thesis.
\end{remark}
\index{Gaussian process!existence|)}

\index{Gaussian process!type of machine learning method|(}\index{Gaussian process interpolation|(}
We can also use Gaussian processes to design algorithms for learning from data.
The simplest method is \emph{Gaussian process interpolation}.
Suppose we are given a (finite) data set $\set{D} \subseteq \set{X}$ and corresponding labels $\vec{y} \in \field^\set{D}$, and we want to learn a functional relationship $g : \set{D} \to \field$.
We adopt a Bayesian perspective. 
Begin with a \emph{prior}\index{prior distribution} that the ``true'' functional relationship $p : \set{D} \to \field$ is a draw from a Gaussian process, $p \sim \GP(\kappa)$; we use a centered prior for simplicity.
We then model the data $\vec{y}$ as measurements of this function $\vec{y} = p(\set{D})$.
We assume, for now, that the measurements are obtained without noise.
Under this model, the \emph{conditional distribution}\index{Gaussian process!conditioning|(} of $p$ is given by the following result.

\index{reproducing kernel Hilbert space!relation to Gaussian process formalism|(}\index{Gaussian process!relation to reproducing kernel Hilbert space formalism|(}
\begin{theorem}[Conditioning a Gaussian process] \label{thm:gp-conditioning}
    Let $p \sim \GP(\kappa)$ be a GP on a base space $\set{X}$, let $\set{D}\subseteq \set{X}$ be a finite subset, and assume $\kappa(\set{D},\set{D})$ is nonsingular.
    Then
    \begin{equation} \label{eq:gp-posterior}
        p \mid \{p(\set{D}) = \vec{y}\} \sim \GP(g, \kappa_{\set{D}}) \quad \text{for } g = \kappa(\cdot,\set{D})\kappa(\set{D},\set{D})^{-1}\vec{y}. 
    \end{equation}
    In particular, the mean function is the kernel interpolant \cref{eq:kernel-interp-equiv} and the covariance function is the residual kernel (\cref{def:nystrom-kernel}).
\end{theorem}
\index{reproducing kernel Hilbert space!relation to Gaussian process formalism|)}\index{Gaussian process!relation to reproducing kernel Hilbert space formalism|)}

Equation \cref{eq:gp-posterior} characterizes the \emph{posterior} distribution\index{posterior distribution} of the Gaussian process conditional on observing the data, $p(\set{D}) = \vec{y}$.
The mean $g : \set{X} \to \field$ of the posterior provides a model of the functional relationship $\set{X} \to \field$ that interpolates the data $(x,y_x)$, and the covariance function\index{covariance function} $\kappa_{\set{D}}$ captures the remaining uncertainty.
In particular, the posterior variance is
\begin{equation} \label{eq:posterior-variance}
    \Var(p(x) \mid p(\set{D}) = \vec{y}) = \kappa_{\set{D}}(x,x).
\end{equation}
Observe that this posterior variance agrees with the pointwise error bound \cref{eq:kernel-interpolation-pointwise-bound} when $\norm{f}_\RKHS = 1$.\index{kernel interpolation!pointwise error bounds}
We see that the diagonal entries of the residual kernel capture the uncertainty in both the kernel interpolation and Gaussian process interpolation settings.
A notable feature of GPs is that the posterior variance \cref{eq:kernel-interpolation-pointwise-bound} depends only on the location of the input data $\set{D}$, not the values of the output data $\vec{y}$.

The proof of \cref{thm:gp-conditioning} is standard, so we omit it; see \cite[Ch.~21]{Tro23} for an introduction to conditioning results for Gaussian random variables.\index{Gaussian process!conditioning|)}

\begin{figure}
    \centering
    \includegraphics[width=0.48\linewidth]{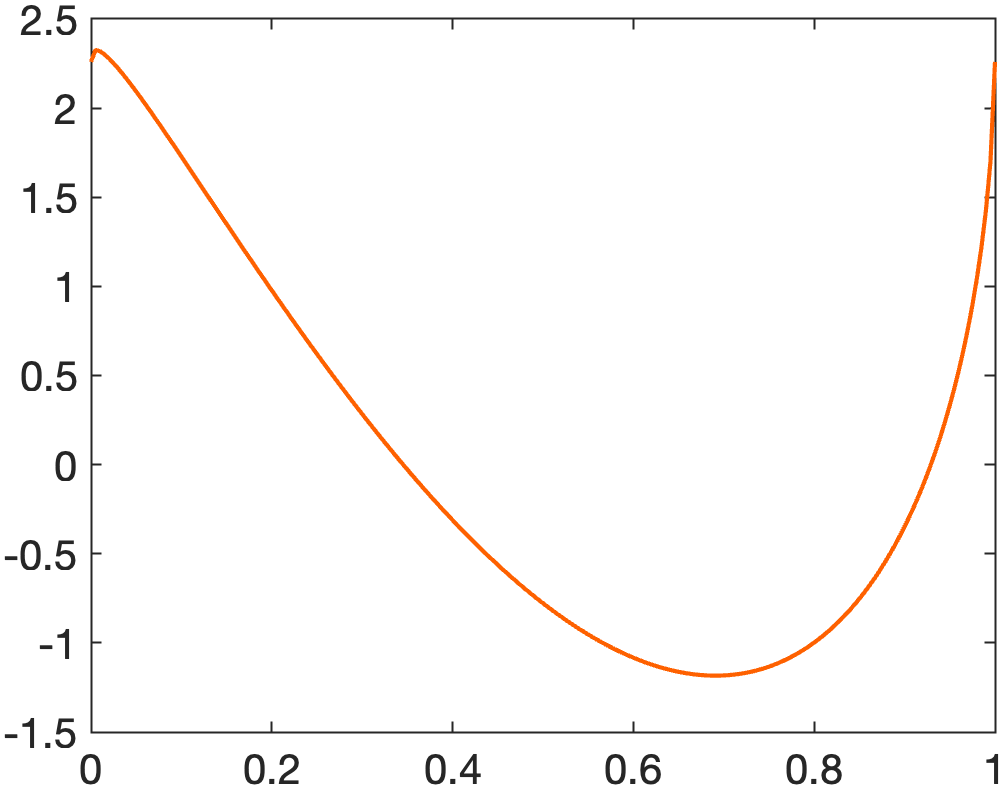}
    \includegraphics[width=0.48\linewidth]{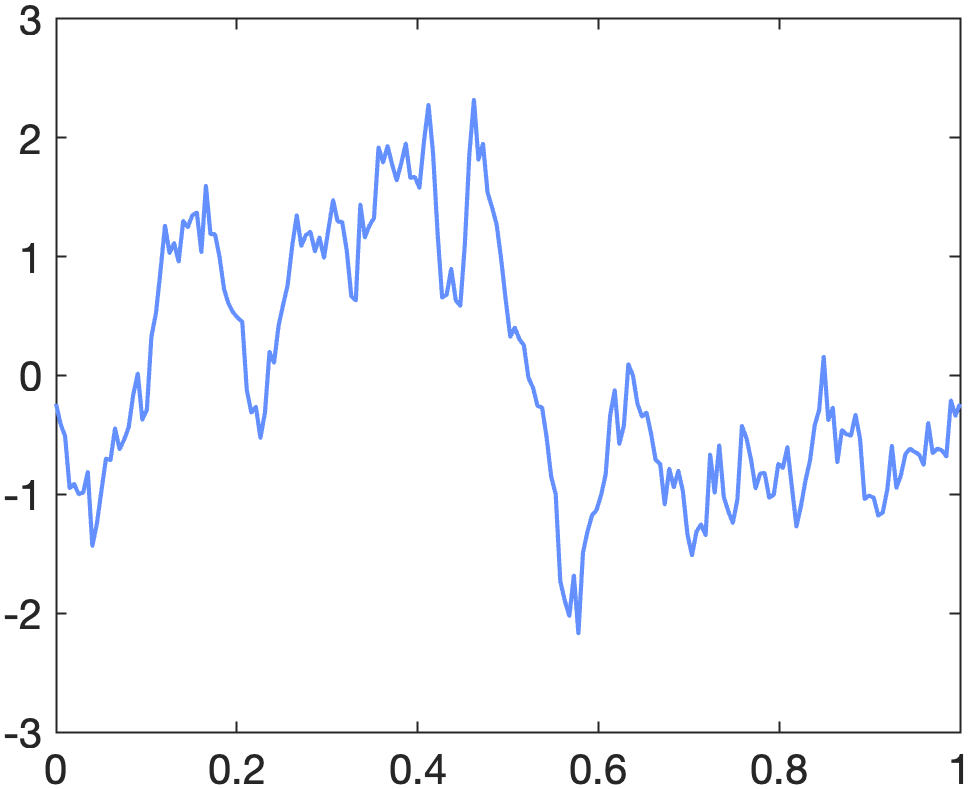}
    \caption[Comparison of RKHS function and draw from a Gaussian process with the same positive-definite kernel]{RKHS function (\emph{left}) and single draw of a Gaussian process (\emph{right}) for the same positive-definite kernel $\kappa$.
    The draw from the Gaussian process is observed to be much ``rougher'' than the RKHS function.}
    \label{fig:gp-vs-rkhs}
\end{figure}

\index{reproducing kernel Hilbert space!relation to Gaussian process formalism|(}\index{Gaussian process!relation to reproducing kernel Hilbert space formalism|(}
\begin{remark}[Gaussian processes don't lie in the RKHS]
    Under fairly general conditions, draws $g \sim \GP(\kappa)$ \warn{do not} belong to the associated RKHS $\RKHS$ for which $\kappa$ is the kernel (with 100\% probability)!
    A visual illustration is provided in \cref{fig:gp-vs-rkhs}, which shows a function $f$ from the \warn{periodic} Sobolev space $\set{H}^1_{\mathrm{per}}([0,1])$ and a draw $g \sim \GP(\kappa)$ from a GP whose covariance function $\kappa$ is the kernel for $\set{H}^1_{\mathrm{per}}([0,1])$.
    The GP realization $g \notin \set{H}^1_{\mathrm{per}}([0,1])$ is observed to be much rougher than the RKHS function $f \in \set{H}^1_{\mathrm{per}}([0,1])$.
\end{remark}\index{Gaussian process interpolation|)}\index{reproducing kernel Hilbert space!relation to Gaussian process formalism|)}\index{Gaussian process!relation to reproducing kernel Hilbert space formalism|)}

\section{Kernel ridge regression and Gaussian process regression} \label{sec:krr_gpr}

Gaussian process interpolation and kernel interpolation are two names and two interpretations for the same methodology.
As the names suggest, this method \warn{interpolates} the data (i.e., $g(x) = y_x$ for every $x \in \set{D}$).  
Conventional wisdom from early machine learning practice and statistical learning theory suggest interpolation can lead to \emph{overfitting},\index{overfitting} yielding a model that fits the data but fails to generalize \cite[\S\S1.1, 3.2, \& 5.5]{Bis06}.
Recent machine learning practice has challenged this conventional wisdom, suggesting that overfitting may not be an issue for many settings in modern machine learning and that one can interpolate the training data without fear.
Many theoretical explanations for this phenomenon of \emph{benign overfitting} have been proposed (e.g., \cite{BLLT20,LR20,Bel21a,CCBG22}).\index{benign overfitting}

\index{Gaussian process regression!reasons for regularization|(}\index{Gaussian process regression|(}\index{kernel ridge regression!reasons for regularization|(}\index{kernel ridge regression|(}
Still, we may still have reasons to want to \emph{regularize} kernel or GP fitting methods.
First, for problems in lower dimensions, overfitting with kernel and Gaussian process methods may be a serious issue.
Second, our data may be provided to us with noise, which regularization may help to mitigate.
Third, regularizing the problem makes it better conditioned and easier to solve with preconditioned iterative methods; see \cref{sec:full-data-preconditioning}.
Finally, even if one wants to interpolate the data, the kernel matrix $\kappa(\set{D},\set{D})$ can sometimes be rank-deficient, at least up to the resolution of floating-point errors. 
For such a rank-deficient problems, adding even a small amount of regularization (say, at the level of the machine precision) may be necessary to obtain meaningful results.\index{Gaussian process regression!reasons for regularization|)}\index{kernel ridge regression!reasons for regularization|)}

Just as before, one can develop regularized data fitting methods in either the RKHS or GP formalisms.
The resulting method will be called \emph{kernel ridge regression} (KRR) or \emph{Gaussian process regression} (GPR).
As a change of pace, we will begin with the GP approach and then develop the corresponding RKHS perspective.\index{kernel ridge regression|(}

\subsection{Gaussian process regression}

To incorporate regularization into a GP data-fitting pipeline, we assume the following model.
We begin from the same assumption the true underlying relationship is drawn from a \emph{prior distribution}\index{prior distribution}
\begin{subequations} \label{eq:gpr-data-model}
\begin{equation}
    p \sim \GP(\kappa)
\end{equation}
We now assume that the data $\vec{y} \in \field^{\set{D}}$ is provided to us corrupted by noise:
\begin{equation} \label{eq:noise-model}
    \vec{y} = p(\set{D}) + \vec{\varepsilon} \quad \text{where } \vec{\varepsilon} \sim \Normal_\field(0,\lambda \, \Id).
\end{equation}
\end{subequations}
The parameter $\lambda \ge 0$ sets the variance of the noise $\vec{\varepsilon} \in \field^{\set{D}}$, which is assumed to have iid entries $(\varepsilon_x : x \in \set{D})$, independent of the GP $g$.
Observe that $\lambda = 0$ recovers the data model for GP interpolation.
The posterior\index{posterior distribution} of $p$ given the data is as follows:
\begin{theorem}[Gaussian process regression] \label{thm:gpr}
    Let $p \sim \GP(\kappa)$ be a GP on a base space $\set{X}$, let $\set{D}\subseteq \set{X}$ be a finite subset, and assume $\lambda > 0$.
    Under the noise model \cref{eq:noise-model}, the posterior is
    \begin{equation*}
        p \mid \vec{y} \sim \GP(g, \kappa^\lambda_{\set{D}}), 
    \end{equation*}
    where
    \begin{equation*}
        g = \kappa(\cdot,\set{D})[\kappa(\set{D},\set{D}) + \lambda \Id]^{-1}\vec{y}
    \end{equation*}
    and
    \begin{equation*}
        \kappa^\lambda_{\set{D}}(x,x') = \kappa(x,x') - \kappa(x,\set{D})[\kappa(\set{D},\set{D}) + \lambda \Id]^{-1} \kappa(\set{D},x') \quad \text{for } x,x' \in \set{X}.
    \end{equation*}
\end{theorem}

\begin{proof}
    We can employ a formal device to reduce the analysis of GPR to the result \cref{thm:gp-conditioning} we already know for Gaussian process interpolation.\index{Gaussian process interpolation}
    Begin by introducing a copy $\set{D}'$ of the set $\set{D}$ and defining $\set{X}' \coloneqq \set{X} \cup \set{D}'$.
    Now, package the data $(p,\vec{y})$ into a Gaussian process $p'$ on $\set{X}'$ by setting $p' = p$ on $\set{X}$ and $p' = \vec{y}$ on $\set{D}'$.
    This GP is centered and its covariance function $\kappa'$ is characterized by the formula 
    \begin{equation} \label{eq:covariance-function-extended}
        \kappa'(\set{D}'\cup \set{E},\set{D}'\cup \set{E}) = \twobytwo{\kappa(\set{D},\set{D}) + \lambda \Id}{\kappa(\set{D},\set{E})}{\kappa(\set{E},\set{D})}{\kappa(\set{E},\set{E})} \quad \text{for } \set{E} \subseteq \set{X} \text{ finite}.
    \end{equation}
    Let us justify this statement block-by-block.
    The $(1,1)$-block is the covariance of the output data $\vec{y}$, and it equals $\kappa(\set{D},\set{D}) + \lambda \Id$ by \cref{eq:noise-model}.
    The $(1,2)$- and $(2,1)$-blocks are the cross-covariances between $\vec{y}$ and $p(\set{E})$, which equal the cross-covariances between $p(\set{D})$ and $p(\set{E})$ since the noise $\vec{\varepsilon}$ is independent of $p$. 
    The $(2,2)$-block is the covariance of $p(\set{E})$, which is $\kappa(\set{E},\set{E})$ by definition of the covariance function.

    The Gaussian process regression posterior conditions $p$ on the output values $\vec{y}$.
    Since $p$ and $p'$ agree on $\set{X}$, conditioning $p$ on $\vec{y}$ is equivalent to conditioning $p'$ on the event $p'(\set{D}') = \vec{y}$, i.e., interpolating the values $\vec{y}$ on the copied space $\set{D}'$.
    In this sense, the Gaussian process \emph{regression} prior on $\set{D}$ is equivalent to the Gaussian process \emph{interpolation} prior on $\set{D}'$.
    Therefore, by \cref{thm:gp-conditioning}, the posterior for Gaussian process regression is 
    \begin{equation*}
        p(\set{E}) \mid \vec{y} \sim \GP(\kappa'(\set{E},\set{D}')\kappa'(\set{D}',\set{D}')^{-1} \vec{y}, \kappa'_{\set{D}'}(\set{E},\set{E})) \quad \text{for } \set{E} \subseteq \set{X} \text{ finite}.
    \end{equation*}
    Invoking the formula \cref{eq:covariance-function-extended} for $\kappa'$ completes the proof.
\end{proof}

\begin{figure}
    \centering
    \includegraphics[height=1.4in]{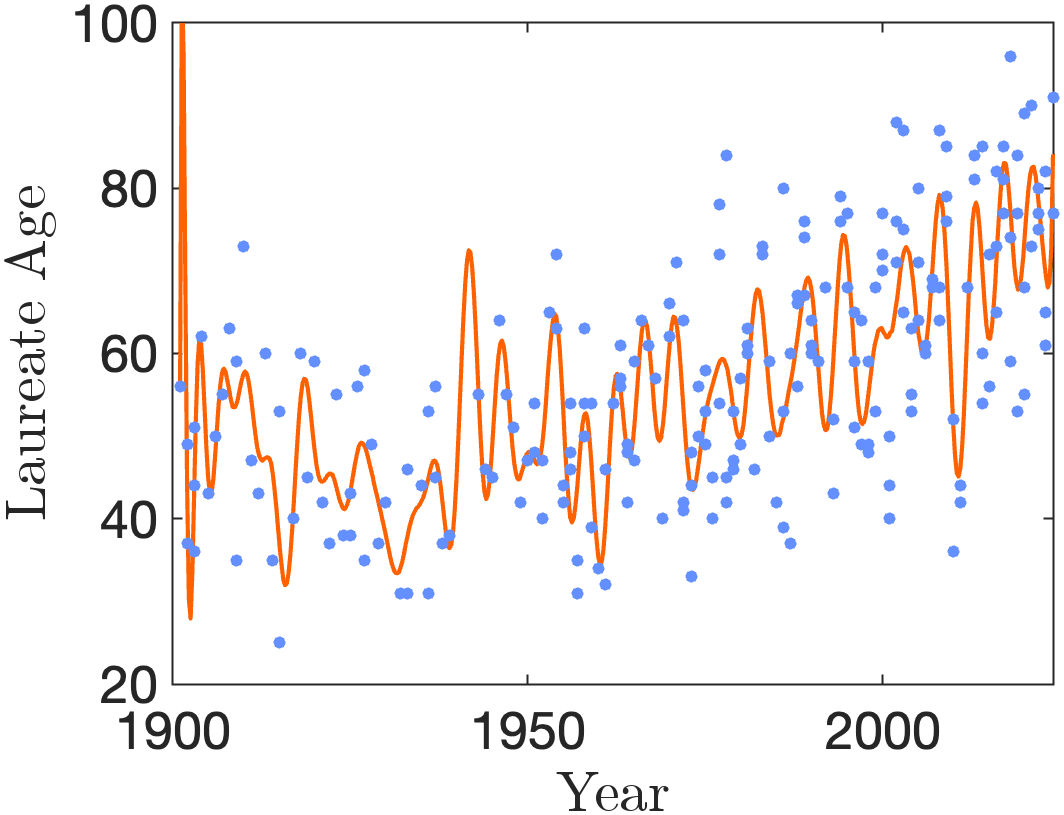}
    \includegraphics[height=1.4in]{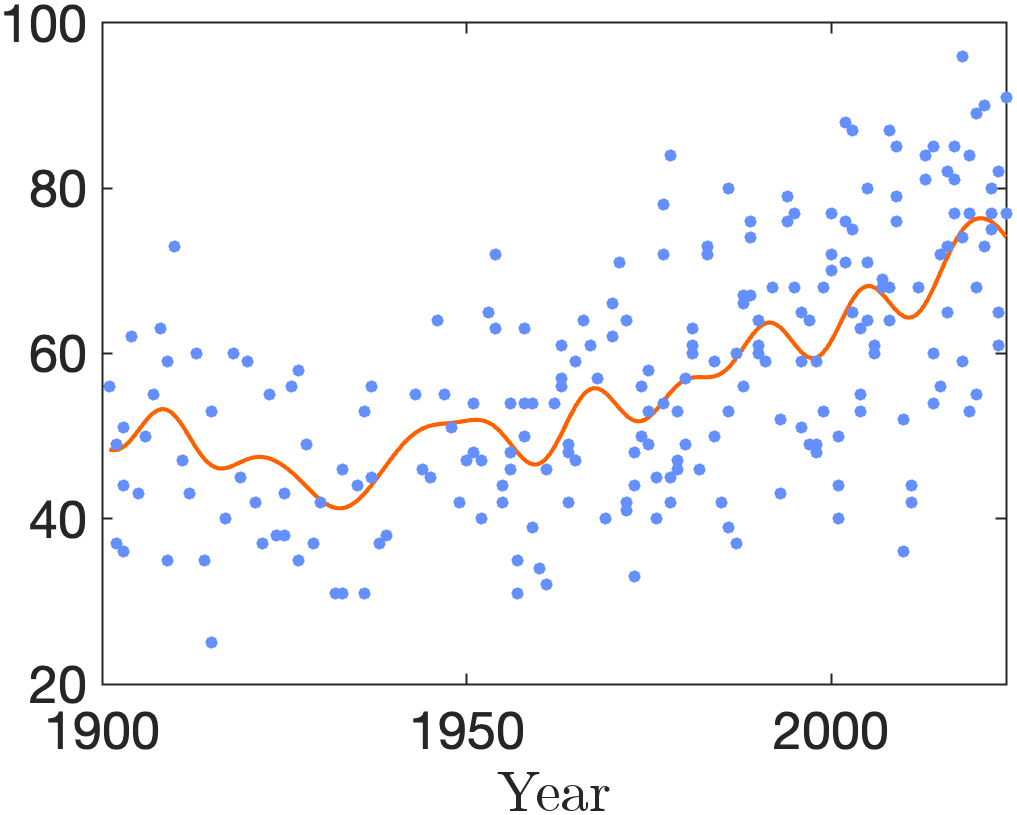}
    \includegraphics[height=1.4in]{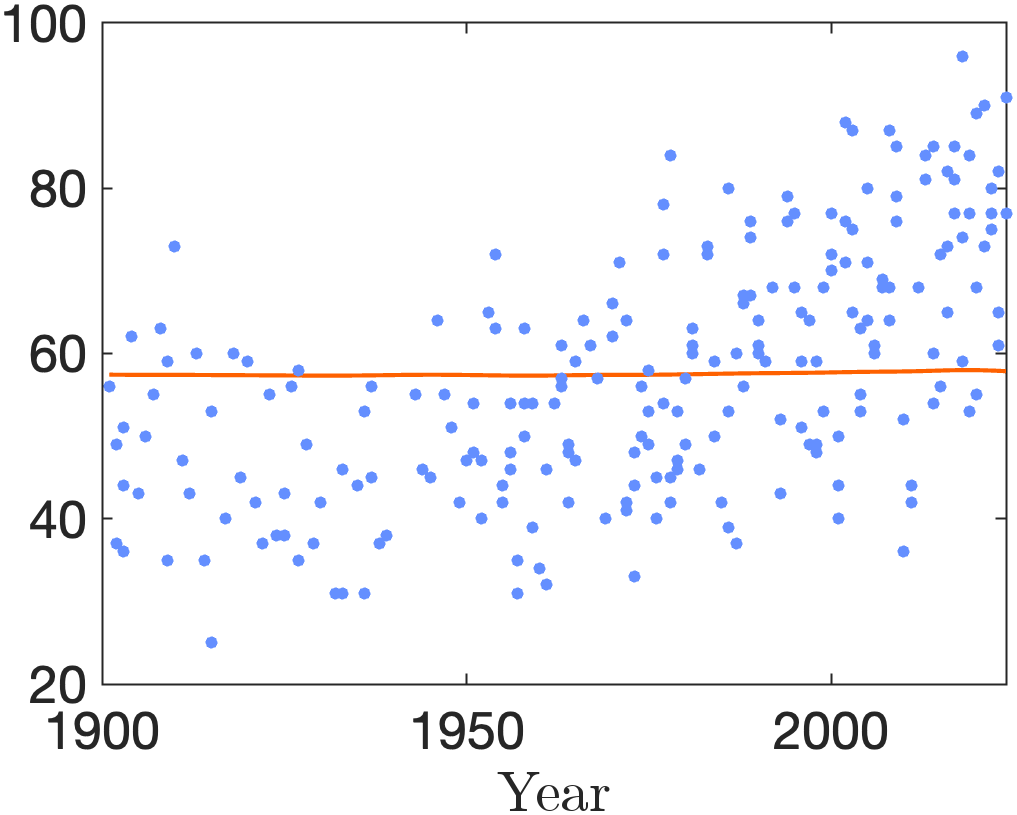}
    
    \makebox[0.32\textwidth]{\centering \qquad$\lambda = 10^{-13}$}%
    \makebox[0.32\textwidth]{\centering \qquad$\lambda = 1$}%
    \makebox[0.32\textwidth]{\centering \qquad$\lambda = 100$}%
    
    \caption[Fit of Nobel laureate ages by prize year using Gaussian process regression with three levels of regularization]{Fitting of Nobel laureate ages by prize year using GPR for three values of the regularization $\lambda = 10^{-13}$ (\emph{left}), $\lambda = 1$ (\emph{middle}), and $\lambda = 100$ (\emph{right}); further details are in the text.}
    \label{fig:krr-reg}
\end{figure}

\index{Gaussian process regression!parameter selection|(}
For low-dimensional fitting tasks with noisy data, choosing the correct regularization $\lambda$ is critical to the success of GPR.
An illustration is provided in \cref{fig:krr-reg}, which applies GPR to fitting the ages of Nobel laureates in Physics versus the year of their prize.
To apply GPR, we first center the ages in the dataset to have mean zero and employ a square-exponential kernel \cref{eq:square-exponential-kernel} with bandwidth $\sigma = 5$.
The left panel uses almost no regularization, and the resulting model is wildly oscillatory.
The right panel uses too much regularization, and the resulting model is nearly constant and thus uninformative.
The middle panel sits between these extremes and produces an interesting model that captures trends in the data.
In particular, this model demonstrates that the average age of Nobel laureates was relatively stable in the twentieth century and has since been rising.
Even this goldilocks models is perhaps more ``wiggly'' then one would like.

\begin{remark}[Parameter selection] \label{rem:mle}
    The Gaussian process framework also provides a natural prescription for picking hyperparameters\index{hyperparameter} such as the regularization parameter $\lambda$ and the bandwidth $\sigma$ for a covariance function such as \cref{eq:square-exponential-kernel} or \cref{eq:laplace-kernel}.
    In general, assume the kernel $\kappa = \kappa_{\vec{\theta}}$ is parametrized by hyperparameters $\vec{\theta}$.
    \warn{Under the data generation model} \cref{eq:gpr-data-model}, the likelihood\index{likelihood} (probability density) of observing the data $\vec{y}$ is
    \begin{equation*}
        L(\vec{y};\vec{\theta},\lambda) = \frac{1}{(2\pi)^{d/2}\det(\kappa_{\vec{\theta}}(\set{D},\set{D}) + \lambda \Id)^{1/2}} \exp \left( - \frac{\vec{y}^*[\kappa_{\vec{\theta}}(\set{D},\set{D}) + \lambda \Id]^{-1}\vec{y}}{2} \right).
    \end{equation*}
    Up to an affine transformation, the log-likelihood\index{log-likelihood} is 
    \begin{equation*}
        \ell(\vec{y};\vec{\theta},\lambda) = -\vec{y}^*[\kappa_{\vec{\theta}}(\set{D},\set{D}) + \lambda \Id]^{-1}\vec{y}-\log\det[\kappa_{\vec{\theta}}(\set{D},\set{D}) + \lambda \Id].
    \end{equation*}
    One can then pick the hyperparameters $\vec{\theta}$ and $\lambda$ to maximize the log-likelihood by an iterative procedure like gradient descent.
    The computation of the log-determinant\index{log-determinant} and its gradients can be a nontrivial problem, for which algorithms based on randomized trace estimation\index{trace estimation} can be employed; see \cref{part:loo}.
\end{remark}\index{Gaussian process!type of machine learning method|)}\index{Gaussian process regression|)}\index{Gaussian process regression!parameter selection|)}

\index{kernel ridge regression|(}
\subsection{Kernel ridge regression}

The GPR method can also be derived within the RKHS formalism.
The kernel interpretation of GPR is referred to as \emph{kernel ridge regression} (KRR).

To derive the KRR problem, we begin by formulating a least-squares regression problem over the RKHS $\RKHS$.
Given data $\vec{y} \in \field^{\set{D}}$, the least-squares fit
\begin{equation*}
    \min_{g \in \RKHS}\, \norm{\vec{y} - g(\set{D})}^2.
\end{equation*}
has infinitely many solutions, all of which interpolate the data (at least as well as possible; see \cref{thm:kernel-interpolation-general}).
To make the procedure more resilient to noise, we can reformulate this problem by adding a ridge penalty term:
\begin{equation} \label{eq:krr-def}
    \hat{g} = \argmin_{g \in \RKHS}\, \norm{\vec{y} - g(\set{D})}^2 + \lambda \, \norm{g}_\RKHS^2 \quad \text{for } \lambda > 0.
\end{equation}
The ridge penalty $\lambda \, \norm{g}_\RKHS^2$ serves to penalize ``roughness'' or ``complexity'', as measured by the RKHS norm.
The (unique) solution $\hat{g}$ is the KRR regression model, which provides a model of the input output relation $\set{X} \to \field$.

\index{kernel ridge regression!solution formula|(}
\begin{theorem}[Kernel ridge regression: Solution formula] \label{thm:krr-solution}
    Let $\RKHS$ be an RKHS and let $\lambda > 0$.
    The KRR problem \cref{eq:krr-def} has a unique solution, which satisfies 
    \begin{equation*}
        g = \kappa(\cdot,\set{D})[\kappa(\set{D},\set{D}) + \lambda \Id]^{-1}\vec{y}.
    \end{equation*}
    Equivalently,
    \begin{equation*} \label{eq:krr-opt-interp}
        g = \sum_{x \in \set{D}} \kappa(\cdot,x) \beta_x \quad \text{where } \vec{\beta} = \argmin_{\vec{\beta} \in \field^{\set{D}}} \, \norm{\kappa(\set{D},\set{D}) \vec{\beta} - \vec{y}}^2 + \lambda \, \vec{\beta}^*\kappa(\set{D},\set{D})\vec{\beta}.
    \end{equation*}
\end{theorem}

\begin{proof}[Proof sketch]
    First, argue that the solution $g$ lies in the column span of the quasimatrix\index{quasimatrix} $\kappa(\cdot,\set{D})$, arguing similarly to \cref{thm:underdetermined-linsys}.
    Then, use the ansatz $g = \kappa(\cdot,\set{D})\vec{\beta}$ and solve for $\vec{\beta}$.
\end{proof}\index{kernel ridge regression!solution formula|)}

Code for kernel ridge regression appears in \cref{prog:krr}.

\myprogram{Kernel ridge regression for data fitting.}{}{krr}

\chapter{Accelerating kernel and Gaussian process methods by subset selection and column Nystr\"om approximation} \label{ch:rpcholesky-kernel-gp}

\index{randomly pivoted Cholesky!for accelerating kernel methods|(}

\epigraph{Fortunately, there is a path forward. To implement kernel methods, we simply need to approximate the kernel matrix\ldots  Even a poor approximation of the kernel can suffice to achieve near-optimal performance, both in theory and in practice.}{Per-Gunnar Martinsson and Joel A.\ Tropp, \textit{Randomized numerical linear algebra: Foundations and algorithms} \cite{MT20a}}

In last chapter, we saw how to use the theories of reproducing kernel Hilbert spaces (RKHSs) and Gaussian processes (GPs) to develop methods for learning from data.
Throughout this discussion, we had little to say about the \emph{computational cost} of these methods.
When run on a dataset $\set{D}$ of size $|\set{D}| = n$, direct implementation of all the methods in the previous chapter require forming, storing, and factorizing the $n\times n$ kernel matrix $\kappa(\set{D},\set{D})$.
This incurs a heavy computational cost, sometimes called the \emph{curse of kernelization} \cite{WCV12}:\index{curse of kernelization}\index{kernel method!computational difficulties}

\actionbox{\textbf{Curse of kernelization.}\index{curse of kernelization}\index{kernel method!computational difficulties} The cost of implementing kernel interpolation, kernel ridge regression, and many other kernel methods on $n$ data points \warn{using standard direct linear algebra methods} requires $\order(n^2)$ storage and $\order(n^3)$ operations.}

In this chapter, we will see how \RPCholesky\index{randomly pivoted Cholesky} and other psd column subset selection methods can be used to accelerate kernel methods, resulting in faster algorithms.

We say up front that there is no free lunch here.
\RPCholesky\index{randomly pivoted Cholesky} and other psd low-rank approximation are effective for kernel problems when the kernel matrix $\kappa(\set{D},\set{D})$ is well-approximated by a low-rank matrix.
Fortunately, many kernel matrices possess this property, so the approaches described in this chapter have wide---though not universal---applicability.

\myparagraph{Sources}
\Cref{sec:full-data-preconditioning} is adapted from the paper

\fullcite{DEF+23}

\iffull\else\Cref{sec:restricted-krr} is adapted from the original \RPCholesky paper\fi 

\fullcite{CETW25}.

The material in \cref{sec:active-krr} is new.

\myparagraph{Outline}
\Cref{sec:full-data-preconditioning,sec:restricted-krr,sec:active-krr} present three approaches to accelerating KRR (equivalently, GPR) using \RPCholesky or other psd column subset selection algorithms; these sections are ordered from most expensive and most accurate to least expensive and most approximate.
\Cref{sec:full-data-preconditioning} discusses using \RPCholesky to precondition the KRR linear systems; this approach leads to full accuracy, but requires $\order(n^2)$ operations even under favorable conditions.
\Cref{sec:restricted-krr} discusses the restricted KRR problem, a cheaper approximate version of KRR with a reduced cost of $\order(k^2n)$, where $k$ is a user-tunable subset size.
\Cref{sec:active-krr} discusses using \RPCholesky for KRR in the setting of \emph{active learning problems}.

\index{column Nystr\"om preconditioning|(}
\section{Column Nystr\"om preconditioning}\label{sec:full-data-preconditioning}

Our first way of using psd low-rank approximation methods to accelerate kernel methods is via \emph{preconditioning}.
Let us focus on KRR.
Per \cref{thm:krr-solution}, the optimal coefficients $\vec{\beta} \in \field^n$ for KRR are the solution to a linear system\index{kernel ridge regression!solution formula}
\begin{equation} \label{eq:krr-lin-sys}
    (\mat{A} + \lambda \Id) \vec{\beta} = \vec{y} \quad \text{for } \mat{A} \coloneqq \kappa(\set{D},\set{D}).
\end{equation}
When $n$ is large (say, so large that one cannot fit the entire matrix $\mat{A}$ in memory at once), it is natural to use an iterative method like conjugate gradient\index{conjugate gradient} \cite[\S6.7]{Saa03} to solve \cref{eq:krr-lin-sys}.
However, if the matrix $\mat{A} + \lambda \Id$ is ill-conditioned, the convergence of iterative methods will be slow.
To improve convergence, we can form a column Nystr\"om approximation $\Ahat \approx \mat{A}$ and use this approximation to \emph{precondition} the linear system, resulting in faster convergence.\index{preconditioning}


\index{column Nystr\"om preconditioning!implementation|(}
Let $\Ahat = \mat{F}\mat{F}^*$ be a column Nystr\"om approximation computed by any method; we recommend \RPCholesky\index{randomly pivoted Cholesky} (or, more precisely, a fast implementation of \RPCholesky,\index{randomly pivoted Cholesky!blocked versions} see \cref{ch:blocked}) for use in practice.
Define the Nystr\"om preconditioner $\mat{P} \coloneqq \Ahat + \lambda \Id$.
In practice, it is grossly inefficient to store $\mat{P}$ directly.
Instead, we compute an (economy-size) SVD $\mat{F} = \mat{U}\mat{\Sigma} \mat{V}^*$, from which we can apply the action of $\mat{P}$ and its inverse via the formulas
\begin{align}
    \mat{P}\vec{z} &= \mat{U}(\mat{\Sigma}^2(\mat{U}^*\vec{z})) + \lambda \vec{z}, \nonumber \\
    \mat{P}^{-1}\vec{z} &= \mat{U}((\mat{\Sigma}^2 + \lambda \Id)^{-1}(\mat{U}^*\vec{z})) + \lambda^{-1}(\vec{z} - \mat{U}(\mat{U}^*\vec{z})). \label{eq:Pinv-formula}
\end{align}
To solve \cref{eq:krr-lin-sys}, we run preconditioned conjugate gradient\index{conjugate gradient} (PCG) with preconditioner $\mat{P}$.
The matrix $\mat{A}+\lambda \Id$ is applied via the formula $(\mat{A} + \lambda \Id)\vec{z} = \mat{A}\vec{z} + \lambda \vec{z}$, and the inverse-preconditioner is applied via the formula \cref{eq:Pinv-formula}.
See \cref{prog:rpcholesky_precon} for an implementation.
(Note that we use the faster \emph{accelerated} version of \RPCholesky algorithm\index{accelerated randomly pivoted Cholesky} for our implementation; see \cref{sec:acc-rpcholesky}.)

\myprogram{\RPCholesky-preconditioned conjugate gradient for solving KRR problems.}{Subroutine \texttt{mypcg} is provided in \cref{prog:mypcg}.}{rpcholesky_precon}

\index{column Nystr\"om preconditioning!computational cost|(}
\myparagraph{Computational cost}
The computational cost of KRR with column Nystr\"om preconditioning consists of generating the entries of the kernel matrix $\mat{A} = \kappa(\set{D},\set{D})$, forming the preconditioner $\mat{P}$, and performing iterations with preconditioned conjugate gradient.\index{conjugate gradient}
We analyze the cost of two variants, a high-memory version where the kernel matrix $\mat{A}$ is formed once and stored, and a low-memory version where the kernel matrix is regenerated each PCG iteration.
We denote by ``$\mathrm{niter}$'' the number of PCG iterations and $t_\kappa\ge 1$ the number of operations required for a single kernel function evaluation.

First, suppose that we generate and store the whole kernel matrix once, requiring $\order(n^2)$ memory and $\order(t_\kappa n^2)$ time.
Afterwards, computing and factorizing $\mat{F} = \mat{U}\mat{\Sigma}\mat{V}^*$ using most column Nystr\"om methods (\RPCholesky, RLS sampling via the RRLS algorithm, uniform sampling, greedy selection) requires $\order(k^2n)$ operations.
Each PCG iteration consists of one matvec with $\mat{A} + \lambda \Id$, costing $\order(n^2)$ operations, and one invocation of the primitive \cref{eq:Pinv-formula}, costing $\order(kn)$ operations.
The total cost is thus $\order((t_\kappa + \mathrm{niter})n^2 + k^2n)$ operations.

\index{column Nystr\"om preconditioning!low-memory version|(}
Second, suppose we regenerate entries of kernel matrix on an as-needed basis.
The storage costs are now dominated by storing the factor $\mat{U}$ needed for the inverse-preconditioner operation \cref{eq:Pinv-formula}, requring $\order(kn)$ memory.
Forming the preconditioner requires $\order(t_\kappa kn + k^2n)$ operations, and each PCG iteration requires $\order(t_\kappa n^2 + kn)$ operations, since we must regenerate the matrix at each iteration.
The total cost is thus $\order(\mathrm{niter}\cdot t_\kappa n^2 + k^2n)$ operations.

\begin{table}[t]
    \centering
    \begin{tabular}{ccc} \toprule
       Implementation  & High-Memory    & Low-Memory \\\midrule
       Storage         & $\order(n^2)$  & $\order(kn)$ \\
       Runtime         & $\order((t_\kappa + \mathrm{niter})n^2 + k^2n)$ & $\order(\mathrm{niter}\cdot t_\kappa n^2 + k^2n)$ \\\bottomrule
    \end{tabular}
    \caption[Storage and runtime costs for low-memory and high-memory column Nystr\"om preconditioned kernel ridge regression]{Storage and runtime costs for low-memory and high-memory column Nystr\"om preconditioned kernel ridge regression. Here, $\mathrm{niter}$ is the number of PCG iterations and $t_\kappa\ge 1$ the number of operations required for a single kernel function evaluation.}
    \label{tab:nystrom-pcg}
\end{table}

\Cref{tab:nystrom-pcg} compares both implementations.
These implementations represent a classic time--space tradeoff. 
Regenerating the kernel matrix makes each PCG iteration significantly more expensive but also cuts the storage costs dramatically.
\index{column Nystr\"om preconditioning!computational cost|)}\index{column Nystr\"om preconditioning!implementation|)}\index{column Nystr\"om preconditioning!low-memory version|)}

\index{column Nystr\"om preconditioning!theoretical results|(}
\myparagraph{Analysis}
How many PCG iterations are required for \RPCholesky-preconditioned KRR to converge?
The following result provides a partial answer:

\index{tail rank|(}
\begin{theorem}[\RPCholesky preconditioning] \label{thm:rpcholesky-preconditioning}
    Fix $\lambda > 0$ and psd matrix $\mat{A} \in \field^{n\times n}$.
    Introduce the \emph{tail rank}
    \begin{equation*}
        \mathrm{d}_{\mathrm{tail}}(\lambda) \coloneqq \min \left\{ r \ge 0 : \sum_{i=r+1}^n \lambda_i(\mat{A}) \le \lambda \right\},
    \end{equation*}
    execute \RPCholesky for $k\ge \mathrm{d}_{\mathrm{tail}}(\lambda)(1 + \log(\tr(\mat{A})/\lambda))$ steps to produce low-rank approximation $\Ahat$, and instantiate the \RPCholesky preconditioner $\mat{P} \coloneqq \Ahat + \lambda \Id$.
    With 90\% probability, the preconditioned condition number is controlled as
    \begin{equation*}
        \cond(\mat{P}^{-1/2}(\mat{A} + \lambda \Id)\mat{P}^{-1/2}) \le 30.
    \end{equation*}
    Consequently, PCG produces a solution $\vec{\beta}^{(\mathrm{niter})}$ satisfying the guarantee 
    \begin{equation*}
        \norm{\smash{\vec{\beta}^{(\mathrm{niter})} - \vec{\beta}}}_{\mat{A} + \lambda \Id} \le \varepsilon \cdot \norm{\smash{\vec{\beta}}}_{\mat{A} + \lambda \Id}
    \end{equation*}
    after at most $\mathrm{niter} \le \lceil 6 \log(2/\varepsilon)\rceil$ steps.
    Here, $\norm{\vec{z}}_{\mat{M}} \coloneqq (\vec{z}^*\mat{M}\vec{z})^{1/2}$ denotes the norm associated to a positive definite matrix $\mat{M}$.
\end{theorem}

A slight strengthening of this result and its proof appear in \cite[Thm.~2.2]{DEF+23}.

\index{effective dimension!vs.\ tail rank|(}
We conjecture that \cref{thm:rpcholesky-preconditioning} holds with the tail dimension $\mathrm{d}_{\mathrm{tail}}(\lambda)$ replaced by the effective dimension\index{effective dimension!conjectured results for randomly pivoted Cholesky} $\mathrm{d}_{\mathrm{eff}}(\lambda)$, defined in \cref{def:ridge-psd}.
To get a sense of the difference between $\mathrm{d}_{\mathrm{tail}}(\lambda)$ and $\mathrm{d}_{\mathrm{eff}}(\lambda)$, suppose the dimension $n$ is very large, and consider a matrix with polynomially decaying spectrum $\lambda_i(\mat{A}) = i^{-p}$, where $p > 1$ is \warn{fixed}.
The tail rank and effective dimensions are $\mathrm{d}_{\mathrm{tail}}(\lambda) = \Theta(\lambda^{-1/(p-1)})$ and $\mathrm{d}_{\mathrm{eff}}(\lambda) = \Theta(\lambda^{-1/p})$.
For a small power, say $p = 2$, these parameters are $\mathrm{d}_{\mathrm{tail}}(\lambda) = \Theta(\lambda^{-1})$ and $\mathrm{d}_{\mathrm{tail}}(\lambda) = \Theta(\lambda^{-1/2})$, leading to a dramatic difference when $\lambda \ll 1$.
For large $p$ (or high $\lambda$), the distinction between tail rank and effective dimension is less significant.
Ultimately, the ``missing link'' that would allow us to replace $\mathrm{d}_{\mathrm{tail}}(\lambda)$ in \cref{thm:rpcholesky-preconditioning} by $\mathrm{d}_{\mathrm{eff}}(\lambda)$ is better bounds for \RPCholesky in the spectral norm; see \cref{sec:better-rpcholesky-bounds} for discussion on the types of bounds I conjecture.
\index{column Nystr\"om preconditioning!theoretical results|)}\index{tail rank|)}\index{effective dimension!vs.\ tail rank|)}

\index{column Nystr\"om preconditioning!parameter selection|(}
\myparagraph{Parameter selection} There are many approaches to choosing parameters for \RPCholesky preconditioning.
One approach, useful in memory-constrained settings, is to just choose the parameter $k$ as large as memory will allow.
A second rule of thumb, guided more by runtime considerations, is to set $k \sim n^{1/2}$, so that the $\order(k^2n)$ cost of forming the preconditioner is comparable to the $\order(n^2)$ cost of a single conjugate gradient step.
Since many CG steps are typically required, I would recommend values of $k \approx 3\sqrt{n}$ or $k\approx 10\sqrt{n}$.
Last, one can use the residual trace to set the approximation rank, running \RPCholesky until the trace of the residual matrix falls below a tolerance.
The parameter $\lambda$ provides a natural guide for the scale of the tolerance.
\index{column Nystr\"om preconditioning!parameter selection|)}

\index{column Nystr\"om preconditioning!numerical results|(}\index{greedy pivoted Cholesky!numerical results|(}\index{uniform sampling for column Nystr\"om approximation!numerical results|(}\index{randomly pivoted Cholesky!numerical results|(}
\myparagraph{Experiment} 
Here, we report an experiment performed by myself and coauthors in \cite[Fig.~6]{DEF+23}.
Here, we apply several column Nystr\"om preconditioners to predict the highest occupied molecule orbital energy for $n=5\times 10^4$ data points from the \texttt{QM9} dataset \cite{RvBR12,RDRv14}.\index{chemistry}
We set the approximation rank to be $k\coloneqq 10^3$.
For each choice of KRR coefficients $\vec{\beta}$, we report the relative residual
\begin{equation*}
    \norm{(\kappa(\set{D},\set{D}) + \lambda \Id)\vec{\beta} - \vec{y}}/\norm{\vec{y}}
\end{equation*}
and the symmetric mean average percentage error\index{symmetric mean average percentage error}
\begin{equation*}
    \operatorname{SMAPE}(\vec{y}^{(\mathrm{test})},\vec{\hat{y}}^{(\mathrm{test})}) = \frac{1}{m} \sum_{i=1}^m \frac{|y_i^{(\mathrm{test})} - \hat{y}_i^{(\mathrm{test})}|}{(|y_i^{(\mathrm{test})}|+|\hat{y}_i^{(\mathrm{test})}|)/2},
\end{equation*}
computed between the test data $\vec{y}^{(\mathrm{test})}$ and the predicted values $\vec{\hat{y}}^{(\mathrm{test})} = \kappa(\set{D}^{(\mathrm{test})},\set{D})\vec{\beta}$ for the test data.

\begin{figure}
    \centering
    \includegraphics[width=0.49\linewidth]{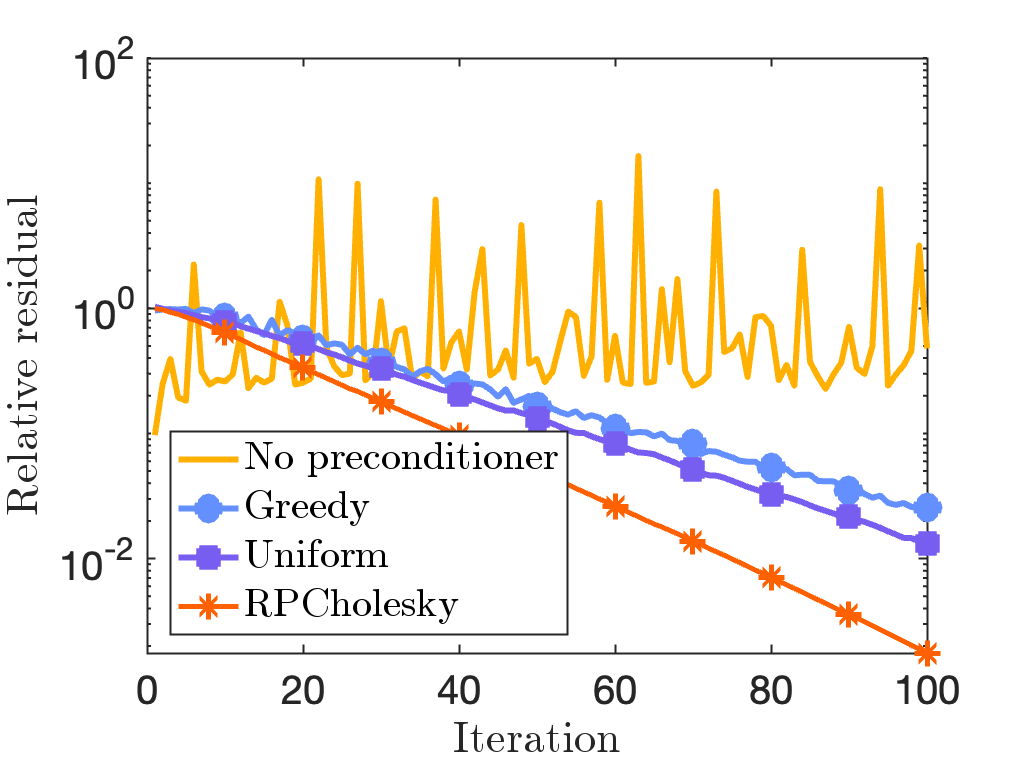}
    \includegraphics[width=0.49\linewidth]{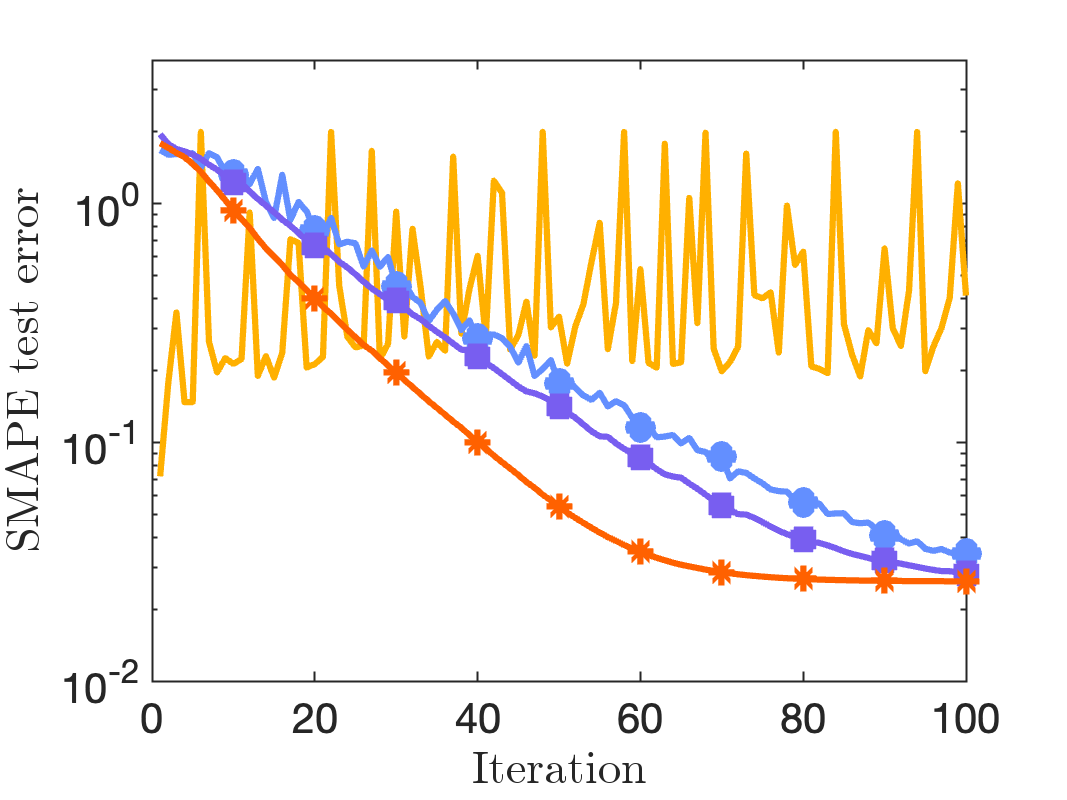}
    
    \caption[Comparison of relative residual and test error for several column Nystr\"om preconditioners]{Relative residual (\emph{left}) and SMAPE test error (\emph{right}) for column Nystr\"om-preconditioned conjugate gradient with greedy (blue circles), uniform (purple squares), \RPCholesky (orange asterisks), and no preconditioning (yellow).}
    \label{fig:rpcholesky-preconditioning}
\end{figure}

Results are shown in \cref{fig:rpcholesky-preconditioning}.
We see that, on this example, \RPCholesky outperforms the other two column Nystr\"om preconditioning strategies, achieving the lowest possible test error in about 60 iterations versus 100 iterations for the other two methods.
More than anything, this example illustrates that \emph{some} preconditioner is absolutely necessary to solve this problem using an iterative method, as the method fails to converge at all without preconditioning.

\Cref{fig:rpcholesky-preconditioning} shows just a single example of the success of Nystr\"om preconditioning. 
See \cite[\S2.1]{DEF+23} for many more experiments, including performance plots which demonstrate that \RPCholesky preconditioning is generally the most effective column Nystr\"om preconditioning method among available strategies on a testbed of examples.\index{column Nystr\"om preconditioning!numerical results|)}\index{greedy pivoted Cholesky!numerical results|)}\index{uniform sampling for column Nystr\"om approximation!numerical results|)}\index{randomly pivoted Cholesky!numerical results|)}\index{column Nystr\"om preconditioning|)}

\index{restricted kernel ridge regression|(}
\section{Restricted kernel ridge regression} \label{sec:restricted-krr}

Even in the most optimistic setting, column Nystr\"om preconditioning for KRR still requires work \emph{quadratic} in the data size $n$.
In this section, we will develop methods with a runtime that is, in principle, \emph{linear} in size of the data.
More precisely, we will describe \emph{restricted kernel ridge regression}, an approximate form of KRR that requires at most $\order(kn)$ storage and $\order(k^2n)$ time, where $1\le k\le n$ is a tunable parameter.
Larger values of $k$ typically leads to more accurate results, at the cost of being more expensive.

\index{restricted kernel ridge regression!derivation and interpretation|(}
\subsection{Description of restricted KRR method}

We have largely described kernel and GP fitting algorithms in abstract terms, as minimum norm interpolants in a Hilbert space, solutions to infinite-dimensional regularized least-squares problems, or as conditional expectations of Gaussian processes.
More prosaically, kernel interpolation is just interpolation of scattered data $(x,y_x)$ by a linear combination of functions $\kappa(\cdot,x)$ for $x \in \set{D}$:
\begin{subequations} \label{eq:krr-full}
\begin{equation} \label{eq:krr-full-fun}
    g = \sum_{x \in \set{D}} \kappa(\cdot,x) \beta_x.
\end{equation}
For KRR, we fit the data by a function of this form that minimizes a regularized least-squares loss.
More precisely, following \cref{thm:krr-solution}, $\vec{\beta}$ is chosen to satisfy
\begin{equation} \label{eq:krr-full-coeffs}
    \vec{\beta} = \argmin_{\vec{\beta} \in \field^{\set{D}}}\, \norm{\kappa(\set{D},\set{D})\vec{\beta}  - \vec{y}}^2 + \lambda \, \vec{\beta}^*\kappa(\set{D},\set{D})\vec{\beta}.
\end{equation}
\end{subequations}

To reduce computational cost and obtain a model $\hat{g} : \set{D} \to \field$ with fewer parameters, it is natural to consider a restricted version of the full optimization problem \cref{eq:krr-full} where $\vec{\beta} \in \field^{\set{D}}$ is only permitted to be nonzero in $k$ selected positions $\set{S} \subseteq \set{D}$.
Denoting $\smash{\widehat{\vec{\beta}}} = \vec{\beta}(\set{S})$, we have the following \emph{restricted} version of the KRR problem:
\begin{subequations} \label{eq:krr-restricted}
\begin{equation} \label{eq:krr-restricted-fun}
    \hat{g} = \sum_{x \in \set{S}} \kappa(\cdot,x) \hat{\beta}_x
\end{equation}
with
\begin{equation} \label{eq:krr-restricted-coeffs}
    \smash{\widehat{\vec{\beta}}} = \argmin_{\vec{\beta} \in \field^{\set{S}}} \, \norm{\kappa(\set{D},\set{S})\smash{\widehat{\vec{\beta}}}  - \vec{y}}^2 + \lambda \, \smash{\smash{\widehat{\vec{\beta}}}}^*\kappa(\set{S},\set{S})\smash{\widehat{\vec{\beta}}}.
\end{equation}
\end{subequations}
We call \cref{eq:krr-restricted} the \emph{restricted KRR method}.
Its solution satisfies the \emph{normal equations}\index{normal equations!for restricted kernel ridge regression}
\begin{equation} \label{eq:krr-restricted-normal-eq}
    [\kappa(\set{S},\set{D})\kappa(\set{D},\set{S}) + \lambda \kappa(\set{S},\set{S})] \smash{\widehat{\vec{\beta}}} = \kappa(\set{S},\set{D})\vec{y}.
\end{equation}
The landmark set $\set{S}$ is fixed during restricted KRR, though it can (and should!) be adaptively selected based on the data $\set{D}$ using a procedure such as \RPCholesky.
The landmarks $\set{S}$ can be chosen by a psd column selection algorithm like \RPCholesky.\index{randomly pivoted Cholesky}
To facilitate comparison, we call \cref{eq:krr-full} the \emph{full-data KRR method}.\index{full-data kernel ridge regression}

Observe that the restricted KRR problem \cref{eq:krr-restricted-coeffs} is a regularized least-squares problem of dimension $n\times k$, where $k = |\set{S}|$.
Its solution defines a regression function $\hat{g}$ that can evaluated at a point using only $k$ kernel function evaluations; compare to the $n$ kernel function evaluations required to evaluate the full KRR model \cref{eq:krr-full-fun}.

\begin{remark}[History and terminology]
    Methods for simplifying a regression problem by using only a subset of basis functions date back at least to work on scattered data interpolation in the 1980s \cite[Ch.~7]{Wah90}.
    The modern literature on restricted KRR began in the GP community with the work of Smola and Bartlett \cite{SB00}, who coined the name \emph{sparse Gaussian process regression} (SGPR).\index{sparse Gaussian process regression}
    The method reemerged in the kernel literature in the work of Rudi, Camoriano, and Rosasco \cite{RCR15} as the \emph{Nystr\"om method}.
    
    I find the terms ``sparse Gaussian process'' and ``Nystr\"om method'' both to be potentially misleading.
    In sparse Gaussian process regression, it is the \emph{coefficient vector} $\vec{\beta}$ that is sparse, not the kernel matrix or the target solution $\hat{g}$.
    Confusingly, methodologies where one designs a kernel with compact support \cite[Ch.~9]{Wen04} or zero out small entries in the kernel matrix, resulting in a sparse kernel matrix, are also common \cite{FGN06}.
    As we have seen already and will continue to see in this thesis, there are several non-equivalent ways of using Nystr\"om approximation to accelerate the solution of KRR problems, so the term ``Nystr\"om method'' can also be ambiguous.
    For these reasons. I use the term ``restricted KRR'' for \cref{eq:krr-restricted}, which my co-authors and I introduced in \cite{DEF+23}.
\end{remark}

\subsection{Restricted KRR as Nystr\"om approximation of the kernel}

So far, we have presented restricted KRR as an ad hoc solution to reduce computational costs.
We can also give interpretations of this method in both the RKHS and GP formalisms.
We look at the RKHS formalism first.

Recall that a subset $\set{S} \subseteq \set{D}$ induces a Nystr\"om approximation 
\begin{equation*}
    \kappa(\set{D},\set{S}) \kappa(\set{S},\set{S})^\dagger\kappa(\set{S},\set{D}) 
\end{equation*}
to the kernel matrix $\kappa(\set{D},\set{D})$ \emph{and} a Nystr\"om approximation of the entire kernel function $\kappa$:\index{Nystr\"om approximation!of a kernel function}
\begin{equation*}
    \hat{\kappa}_{\set{S}}(x,x') = \kappa(x,\set{S}) \kappa(\set{S},\set{S})^\dagger \kappa(\set{S},x') \quad \text{for } x,x' \in \set{X}.
\end{equation*}
This observation leads to an interpretation of restricted KRR, which was suggested to me by Yifan Chen:

\begin{theorem}[Restricted KRR as Nystr\"om approximation of the kernel] \label{thm:yifan}
    Let $\set{S} \subseteq \set{D} \subseteq \set{X}$ be nested finite subsets of a base space $\set{X}$, and let $\vec{y} \in \field^{\set{D}}$ be data.
    Assume the regularization parameter $\lambda >0$ is \warn{positive} and the kernel matrix $\kappa(\set{S},\set{S})$ is \warn{invertible}.
    The restricted KRR function $\hat{g}$ given by \cref{eq:krr-restricted-fun} is the output of the full-data KRR method applied to the data $\vec{y} \in \field^{\set{D}}$ with the Nystr\"om approximate kernel:
    \begin{equation} \label{eq:yifan}
        \hat{g} = \hat{\kappa}_{\set{S}}(\cdot,\set{D}) [\hat{\kappa}_{\set{S}}(\set{D},\set{D}) + \lambda \Id]^{-1} \vec{y}.
    \end{equation}
\end{theorem}

\index{ridge regression identity|(}\index{nonsymmetric ridge regression identity|(}
To prove this result, we will employ the following lemma:

\begin{lemma}[Nonsymmetric ridge regression identity] \label{lem:nonsymmetric-ridge}
    Let $\mat{B},\mat{C} \in \field^{n\times k}$ be matrices, and choose $\lambda \in \field$ to ensure that the matrix inverses below are well-defined.
    Then
    \begin{equation*}
        (\mat{B}^*\mat{C} + \lambda \Id)^{-1} \mat{B}^* = \mat{B}^* (\mat{C}\mat{B}^* + \lambda \Id)^{-1}.
    \end{equation*}
\end{lemma}

The case $\mat{B} = \mat{C}$, yields the identity $(\mat{B}^*\mat{B} + \lambda \Id)^{-1}\mat{B}^* = \mat{B}^*(\mat{B}\mat{B}^* + \lambda \Id)$, useful in the analysis of ridge-regularized linear regression,\index{ridge regularization!for linear least squares} which may be called the (symmetric) ridge regression identity.
The formula is easily checked by multiplying to clear the inverses and checking that the left- and right-hand sides agree.
It may also be derived using the Sherman--Morrison--Woodbury formula.

     
\index{ridge regression identity|)}\index{nonsymmetric ridge regression identity|)}

\begin{proof}[Proof of \cref{thm:yifan}]
    Let $h$ denote the right-hand side of \cref{eq:yifan}.
    Our goal is to show $\hat{g} = h$.
    We prove the case $\lambda > 0$.
    
    Using the definition of the Nystr\"om-approximate kernel, we may write
    \begin{equation*}
        h = \kappa(\cdot,\set{S})\kappa(\set{S},\set{S})^{-1} \kappa(\set{S},\set{D})[\kappa(\set{D},\set{S})\kappa(\set{S},\set{S})^{-1}\kappa(\set{S},\set{D}) + \lambda \Id]^\dagger \vec{y}.
    \end{equation*}
    Introduce $\mat{B} = \kappa(\set{D},\set{S})\kappa(\set{S},\set{S})^{-1}$ and $\mat{C} = \kappa(\set{D},\set{S})$, and invoke the nonsymmetric ridge regression identity\index{ridge regression identity}\index{nonsymmetric ridge regression identity} to obtain 
    \begin{align*}
        h &= \kappa(\cdot,\set{S}) [\kappa(\set{S},\set{S})^{-1}\kappa(\set{S},\set{D})\kappa(\set{D},\set{S}) + \lambda \Id]^{-1}\kappa(\set{S},\set{S})^{-1}\kappa(\set{D},\set{S}) \vec{y} \\
        &= \kappa(\cdot,\set{S}) [\kappa(\set{S},\set{D})\kappa(\set{D},\set{S}) + \lambda \kappa(\set{S},\set{S})]^{-1}\kappa(\set{D},\set{S}) \vec{y}.
    \end{align*}
    Utilize the normal equations\index{normal equations!for kernel ridge regression} \cref{eq:krr-restricted-normal-eq} to complete the proof.
\end{proof}

\index{sparse Gaussian process regression|(}
\subsection{Sparse Gaussian process regression}

The GP framework provides alternative perspectives on the restricted KRR method.
In their original work, Smola and Bartlett \cite{SB00} introduce SGPR as a way of restricting the GPR problem to a subset of basis functions to increase computational efficiency, similar to our original derivation of restricted KRR above.

One may also develop probabilistic interpretations of restricted KRR.
The easiest interpretation is that restricted KRR is equivalent to ordinary Gaussian process regression on the \emph{restricted prior}
\begin{equation*}
    \hat{p} \coloneqq \expect[p \mid p(\set{S})].
\end{equation*}
The restricted prior $\hat{p}$ averages out all the randomness in the initial prior $p$ except at the landmark points $\set{S}$.
The equivalence of this approach with restricted KRR follows from \cref{thm:yifan,thm:Gaussian-expectations}.
There is also a more sophisticated probabilistic interpretation due to Titsias \cite{Tit09}.

\index{restricted kernel ridge regression!derivation and interpretation|)}\index{sparse Gaussian process regression|)}

\index{restricted kernel ridge regression!with randomly pivoted Cholesky|(}\index{restricted kernel ridge regression!numerical results|(}\index{uniform sampling for column Nystr\"om approximation!numerical results|(}\index{greedy pivoted Cholesky!numerical results|(}\index{randomly pivoted Cholesky!numerical results|(}
\subsection{Experiment}

Experiments for restricted KRR with \RPCholesky and other column Nystr\"om methods on a scientific dataset are provided in \cite[\S4.2]{CETW25}.
Here, we provide a more conceptual experiment to illustrate the benefits of \RPCholesky as a method for selecting the subset $\set{S}$ for use with restricted KRR.

\begin{figure}
    \centering
    \includegraphics[width=0.41\linewidth]{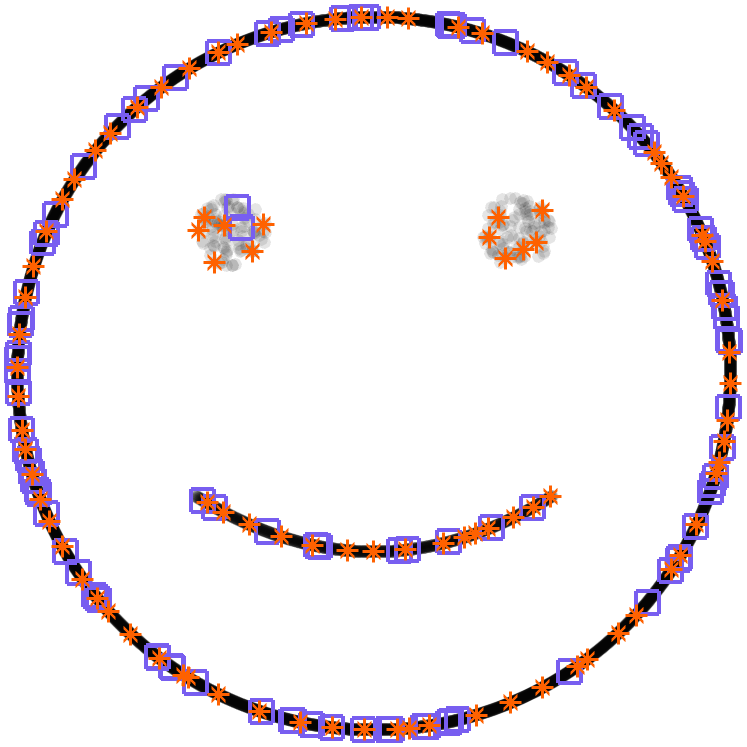}
    \includegraphics[width=0.58\linewidth]{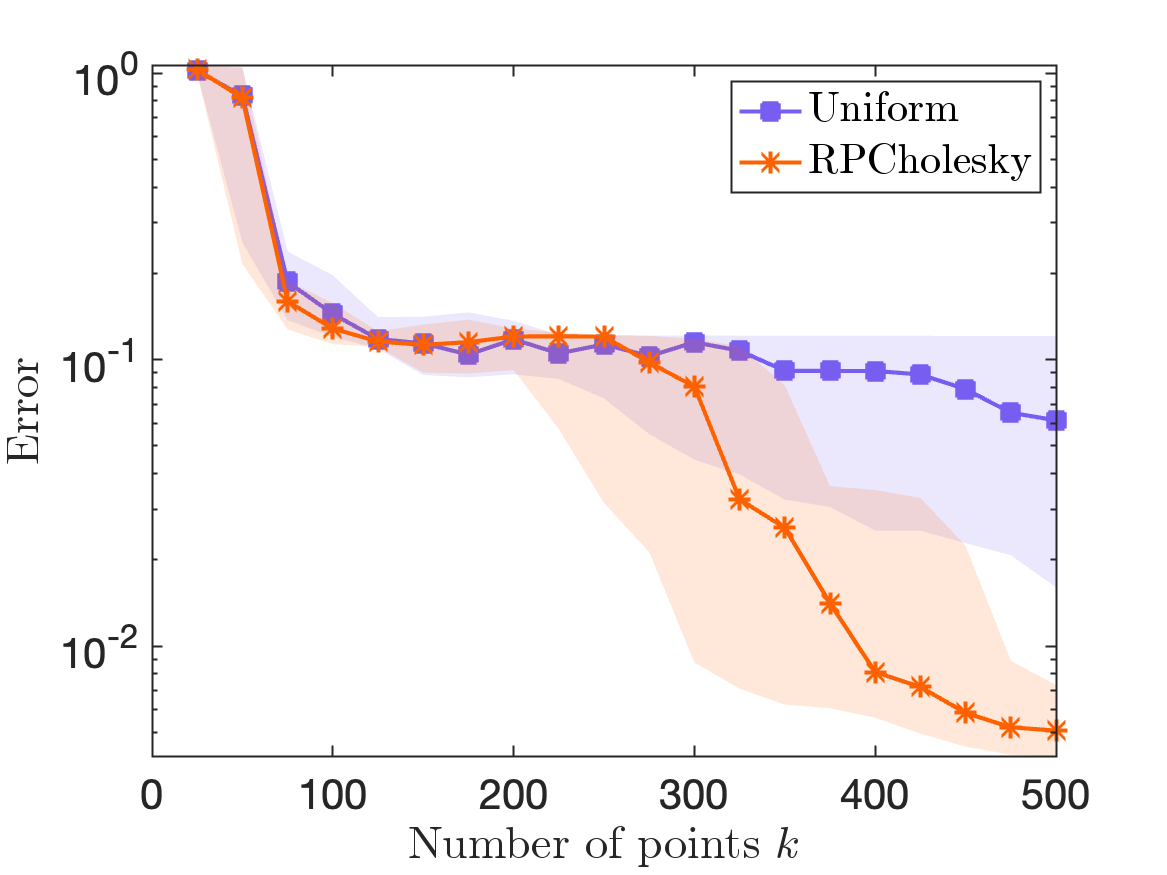}
    \caption[Test errors for restricted KRR with landmarks selected by \RPCholesky and uniformly selected points]{\emph{Left}: Error \cref{eq:error-restricted-krr} for restricted KRR with \RPCholesky-selected points (orange asterisks) and uniformly selected points (purple squares).
    Lines show the median of 100 trials, and shaded regions show the 10\% and 90\% quantiles.
    \emph{Right:} Sets of $k=100$ points selected either by \RPCholesky or uniformly at random, overlayed on the data $\set{D}$ (black translucent circles).
    \RPCholesky is seen to produce an even coverage of the data, whereas uniform sampling misses the right eye.}
    \label{fig:restricted-krr}
\end{figure}

We consider the task of fitting a function $f: \real^2\to\real$ from data $\set{D} \subseteq \real^2$ and $\vec{y} \coloneqq f(\set{D}) \in \real^{\set{D}}$.
For concreteness, we choose the function to be $f(\vec{x}) = \sin((x_1+x_2)/10)$ and $\set{D}$ to consist of $n=10^4$ points in $\real^2$ forming a smiley face with radius $10$; see the right panel of \cref{fig:restricted-krr} for illustration.
We use the Mat\'ern-5/2 kernel with bandwidth $\sigma=40$, and we set the regularization to $\lambda \coloneqq 0$.

Results for restricted KRR with both \RPCholesky-selected points and uniformly selected data points are shown in the right panel of \cref{fig:restricted-krr}.
We see that the error, measured as the largest difference between the true function $f$ and the model $g$ on the dataset $\set{D}$,
\begin{equation} \label{eq:error-restricted-krr}
    \mathrm{error} \coloneqq \max_{\vec{x}\in\set{D}} | f(\vec{x}) - g(\vec{x})|,
\end{equation}
is substantially lower for \RPCholesky restricted KRR than uniform restricted KRR.
An explanation for why can be seen in the left panel of \cref{fig:restricted-krr}, which shows $k=100$ points selected either by \RPCholesky or uniformly at random.
Each eye comprises just $n^{1/2} \ll n$ data points, which are easy for uniform sampling to miss.
\RPCholesky, by contrast, selects points that are well-distributed across the whole data set, including both eyes.

\index{restricted kernel ridge regression|)}\index{restricted kernel ridge regression!with randomly pivoted Cholesky|)}\index{restricted kernel ridge regression!numerical results|)}\index{uniform sampling for column Nystr\"om approximation!numerical results|)}\index{greedy pivoted Cholesky!numerical results|)}\index{randomly pivoted Cholesky!numerical results|)}

\index{active learning|(}\index{active learning!for kernel regression|(}\index{kernel ridge regression!active learning|(}

\section{Active learning for kernel interpolation and ridge regression} \label{sec:active-krr}

So far, we have seen two approaches to that use \RPCholesky and other subset selection algorithms to improve the efficiency of KRR.
First, in \cref{sec:full-data-preconditioning}, we employed column Nystr\"om approximations to precondition the full-data KRR problem.
Under ideal conditions, the cost of this approach was as low as $\order(n^2)$ operations.
Next, in \cref{sec:restricted-krr}, we used subset selection algorithms to \emph{restrict} the KRR problem to a set of $k$ selected centers, at a cost of $\order(k^2n)$ operations.
In both of these sections, we assumed access to a dataset of $n$ data points $\set{D} \subseteq \set{X}$ \emph{with labels} $\vec{y} \in \field^n$.
In this section, we will consider an \emph{active learning setting} where we are given just the unlabeled training data $\set{D} \subseteq \set{X}$ and must choose what data to label.

For now, we suppose the set of data points $\set{D}$ to potentially label is \warn{finite}.
Unlike our other settings, the assumption that $\set{D}$ is finite is non-vacuous; one can imagine applications where there is an infinite set of points that we could label. 
We will return to this infinite setting in \cref{ch:infinite}.
Given the set $\set{D}$, we can select $k$ points to label by running \RPCholesky (or another subset selection algorithm).
After we have a subset $\set{S} \subseteq \set{D}$, we are free to fit whatever model we want to the data $\set{S}$ and labels $\vec{y} \in \field^{\set{S}}$.
In this framework, it is most natural to use a kernel or GP method, like kernel interpolation or KRR, to fit the data.
One could also combine kernel-based subset selection with any type of machine learning method such as an artificial neural network.
Code for active KRR with \RPCholesky is provided in \cref{prog:rpcholesky_active_krr}, and theoretical guarantees are provided later in \cref{cor:kernel-interpolation-rpcholesky}.

\myprogram{Active kernel ridge regression with data points selected by \RPCholesky.}{Subroutines \texttt{krr} and \texttt{acc\_rpcholesky} are provided in \cref{prog:acc_rpcholesky,prog:krr}.}{rpcholesky_active_krr}

\index{active learning!definition of|(}
\begin{remark}[Is this active learning?]
    The term ``active learning'' is contested in the literature.
    Miller and Calder suggest drawing a distinction between \emph{coreset methods}, which ``leverage the geometry of the underlying dataset\ldots but not the set of labels observed at labeled points during the active learning process'', and \emph{active learning methods}, which use both the geometry of the data set $\set{D}$ \emph{and} the observed labels $y_x$ at queried points $x \in \set{D}$ \cite[\S3.1]{MC23}.
    Other researchers do not make this distinction and refer to any machine learning that chooses points to label as an active learning method \cite{CP19,MMWY22}.
    We shall adopt the second definition in this work, though our methods would be described as a ``coreset method'' in Miller and Calder's taxonomy.
\end{remark}
\index{active learning!definition of|)}

\index{randomly pivoted Cholesky!numerical results|(}\index{uniform sampling for column Nystr\"om approximation!numerical results|(}
To test \RPCholesky as a method for active learning, we apply active kernel interpolation to $n=10^5$ randomly selected points from the \texttt{SUSY} dataset \cite{BSW14}; we hold out an additional $n$ randomly selected points as a test set.
Data is standardized.
We compare active kernel interpolation with subsets selected uniformly at random to those selected by \RPCholesky.
As the kernel, we use a Laplace kernel with $\ell_2$ distances and bandwidth $\sigma = 4$.
As a baseline, we also report the errors with restricted KRR (with regularization parameter $\lambda = 0$) using the same subsets $\set{S}$.

\begin{figure}
    \centering
    \includegraphics[width=0.7\linewidth]{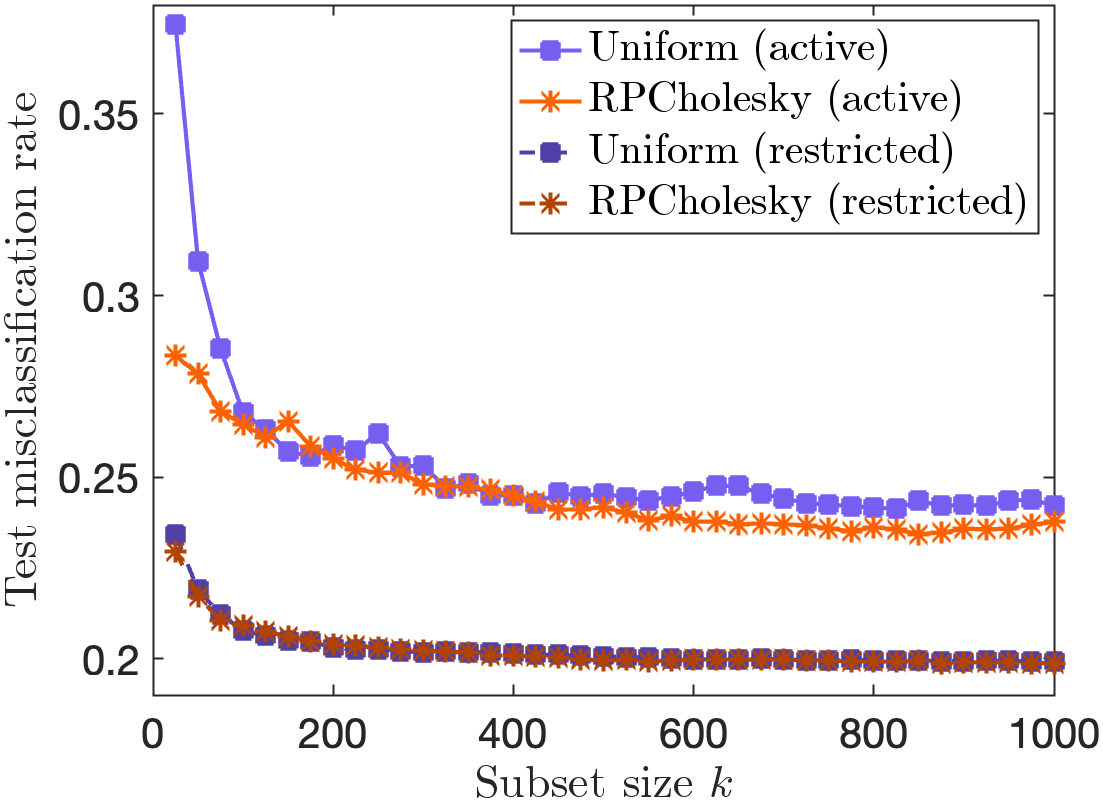}
    \caption[Comparison of active and restricted KRR with landmarks selected by \RPCholesky and uniformly at random]{Comparison of active (line solid lines) and restricted (dark dashed lines) kernel interpolation with subset $\set{S}$ selected by \RPCholesky (orange asterisks) or uniformly at random (purple squares).}
    \label{fig:active-krr}
\end{figure}

Results are shown in \cref{fig:active-krr}.
We see that \RPCholesky active kernel interpolation outperforms active kernel interpolation with uniformly selected points for most values of $k$; the benefits are pronounced when the number of labeled points is small, say, $k < 100$.
The baselines results for restricted KRR show much smaller errors than the active learning methods, demonstrating that it is valuable to use label information $y_x$ at each point $x \in \set{D}$ if it is available.

\index{kernel ridge regression|)}\index{active learning|(}\index{active learning!for kernel regression|)}\index{kernel ridge regression!active learning|)}\index{randomly pivoted Cholesky!numerical results|)}\index{uniform sampling for column Nystr\"om approximation!numerical results|)}

\index{randomly pivoted Cholesky!for accelerating kernel methods|)}

\chapter{To infinite dimensions} \label{ch:infinite}

\epigraph{One should never try to prove anything that is not almost obvious.}{A sentiment attributed to Alexandre Grothendieck by Allyn Jackson in \emph{As if summoned from the void, the life of Alexandre Grothendieck} \cite{Jac04}}

We have seen the \RPCholesky algorithm is effective at producing low-rank approximations to a psd matrix.
In the previous chapter, we applied this algorithm to the kernel matrix associated with a finite set of data points, and we used the resulting low-rank approximations and subsets of landmark points to accelerate kernel machine learning algorithms.

In \cref{ch:kernels-gaussian}, we worked with kernels $\kappa : \set{X} \times \set{X}\to \field$ defined over an infinite base space $\set{X}$, but we only applied \RPCholesky to kernel matrices $\kappa(\set{D},\set{D})$ associated with \warn{finite} subsets $\set{D} \subseteq \set{X}$.
It is natural to ask: \emph{Is there an infinite-dimensional version of \RPCholesky that we can use to construct low-rank approximations to the infinite kernel function $\kappa$ and to identify landmark points $\set{S}\subseteq\set{X}$ for the infinite set $\set{X}$?}
This chapter answers this question in the affirmative.

\myparagraph{Sources}
This chapter is based on the following paper:

\fullcite{EM23a}.

The material has been significantly expanded, including a much lengthier introduction to the functional analysis setting, the new \cref{thm:cont-rpcholesky-trace} (generalizing \cref{thm:rpcholesky-trace}), and a new application to active kernel interpolation in \cref{sec:infinite-active-kernel-interp}.

\myparagraph{Outline}
\Cref{sec:infinite-setting} sets the stage by introducing appropriate operators and function spaces for infinite-dimensional low-rank approximation.
\Cref{sec:infinite-rpcholesky} presents the infinite-dimensional \RPCholesky algorithm, and \cref{sec:rpc-rejection} discusses implementations of this procedure using rejection sampling.
\Cref{sec:infinite-active-kernel-interp,sec:quadrature} provide applications to active learning and quadrature.

\section{Mathematical setting} \label{sec:infinite-setting}

In this section, we develop a functional analysis setting for low-rank approximation of kernels and subset selections on general, possibly infinite sets, following the treatments given in \cite{HB04,Bac17,BBC19,EM23a}.
Inspired by the Grothendieck quote at the start of this chapter, our goal will be to develop a function space setting where the infinite dimensional extension of \RPCholesky and its analysis is a direct translation of the ordinary matrix setting.

\index{reproducing kernel Hilbert space!relationship to $\Ltwo$ space|(}\index{L2, function space@$\Ltwo$, function space|(}
\myparagraph{Function spaces}
Our setting will be a topological space $\set{X}$ endowed with a Borel measure $\mu$, upon which we define an RKHS $\RKHS$ of functions $f : \set{X} \to \field$.
We let $\kappa : \set{X}\times \set{X} \to \field$ denote the reproducing kernel of $\set{H}$, which we assume is continuous.

This set $\set{X}$ supports two spaces of functions, the RKHS $\RKHS$ and the space of square integrable functions $\Ltwo(\mu)$.
The relation between these spaces will be crucial to our development.
To ensure these two spaces place nicely with each other, we make the assumption that the kernel $\kappa$ is integrable along the diagonal:
\begin{equation*}
    \int_{\set{X}} \kappa(x,x) \, \d\mu(x) < +\infty
\end{equation*}
and that $\set{H}$ is dense in $\Ltwo(\mu)$.
These assumptions imply that $\set{H}$ is compactly embedded in $\Ltwo(\mu)$ \cite[Prop.~2]{HB04}.

\index{integral operator!for a reproducing kernel Hilbert space|(}\index{reproducing kernel Hilbert space!integral operator|(}
\myparagraph{Discovering the integral operator}
Since $\RKHS$ is compactly embedded in $\Ltwo(\mu)$, the identity mapping $\iota f \coloneqq f$ defines a compact linear map $\iota : \RKHS \to \Ltwo(\mu)$.
The maps $\iota$ ``forgets'' that a function $f$ belongs to $\RKHS$ and treats it as a function in $\Ltwo(\mu)$.
Introduce the symbol $A \coloneqq \iota^*$ for its adjoint.
The adjoint $A : \Ltwo(\mu) \to \RKHS$ maps a function $u\in \Ltwo(\mu)$ (which, in general, can be quite rough) to a ``smooth'' function $A u \in \set{H}$.
What could this mysterious operator $A$ be?

We shall discover the correct answer by a formal calculation.
Let $u \in \Ltwo(\mu)$ and $f \in \RKHS$ be arbitrary, and compute the inner product
\begin{equation*}
    \langle \iota f, u \rangle_{\Ltwo(\mu)} = \int_{\set{X}} \overline{f(x)} u(x) \, \d \mu(x).
\end{equation*}
Recall that the adjoint is defined via the relation $\langle \iota f,u \rangle_{\Ltwo(\mu)} = \langle f,\iota^* u\rangle_\RKHS$, so we must find a way of introducing the $\RKHS$ inner product.
To do so, employ the reproducing property $\langle f, \kappa(\cdot,x) \rangle_\RKHS = \overline{f(x)}$ to obtain
\begin{equation*}
    \langle \iota f, u \rangle_{\Ltwo(\mu)} = \int_{\set{X}} \langle f, \kappa(\cdot,x) \rangle_\RKHS u(x)  \, \d \mu(x).
\end{equation*}
Pull the integral over $x$ into the inner product, yielding
\begin{equation*}
    \langle \iota f, u \rangle_{\Ltwo(\mu)} = \left\langle f, \int_{\set{X}} \kappa(\cdot,x )u(x) \, \d \mu(x) \right\rangle_\RKHS \eqqcolon \langle f,\iota^* u\rangle_\RKHS.
\end{equation*}
We conclude that the adjoint $A = \iota^*$ is the integral operator
\begin{equation} \label{eq:integral-operator}
    Au = \int_{\set{X}} \kappa(\cdot,x )u(x) \, \d \mu(x).
\end{equation}
A rigorous version of this argument is provided in \cite[Prop.~2]{HB04}.

\myparagraph{Properties of the integral operator}
The integral operator $A$, defined by \cref{eq:integral-operator}, can also be seen as an operator on $\Ltwo(\mu)$, which we denote $A_{\Ltwo(\mu) \to \Ltwo(\mu)}$.
Using the ``forgetting'' map $\iota : \RKHS \to \Ltwo(\mu)$, this operator can be written
\begin{equation*}
    A_{\Ltwo(\mu) \to \Ltwo(\mu)} = A_{\Ltwo(\mu) \to \RKHS} \cdot \iota = \iota\iota^*.
\end{equation*}
Thus, since the inclusion operator $\iota$ is compact (since $\RKHS$ is compactly embedded in $\Ltwo(\mu)$), we conclude that $A_{\Ltwo(\mu) \to \Ltwo(\mu)}$ is compact and psd.

To this point, we have been very careful about domains and codomains of linear mappings; we will now permit ourselves to be more lax and will use $A$ to refer to the transformation \cref{eq:integral-operator}, whatever its domain and codomain are.
Since $A$ is a compact psd operator on $\Ltwo(\mu)$, it admits a spectral decomposition
\begin{equation} \label{eq:operator-eval-decomp}
    A = \sum_{i=1}^\infty \lambda_i^{\vphantom{*}} e_i^{\vphantom{*}} e_i^*,
\end{equation}
where $\lambda_1\ge \lambda_2 \ge \cdots \ge 0$ are the eigenvalues of $A$ and $\{e_i\}$ form an orthonormal basis of $\Ltwo(\mu)$.
(We will tacitly assume that $\Ltwo(\mu)$ is infinite-dimensional for this chapter, though the formalism all carries through for finite-dimensional spaces as well.)
Here, we let 
\begin{equation*}
    u^* \coloneqq \langle u,\cdot \rangle_{\Ltwo(\mu)} : \Ltwo(\mu) \to \field
\end{equation*}
denote the adjoint of a function $u \in \Ltwo(\mu)$, identified with the linear transformation $\alpha \mapsto \alpha u$ from $\field \to \Ltwo(\mu)$.

In addition to being compact, the operator $A$ on $\Ltwo(\mu)$ is trace-class.
Its trace is
\begin{equation} \label{eq:trace-equals-integral}
    \tr(A) = \int_{\set{X}} \kappa(x,x) \, \d\mu.
\end{equation}
This conclusion may be derived formally as follows.
Observe that an arbitrary function in $\Ltwo(\mu)$ can be represented a limit $u = \lim_{n\to\infty} (\sum_{i=1}^n e_i^{\vphantom{*}}e_i^*)u$.
Therefore, we compute
\begin{align*}
    \int_{\set{X}} \kappa(x,x) \, \d\mu 
    &= \int_{\set{X}} \sum_{i=1}^\infty e_i^{\vphantom{*}}(x) (e_i^*\kappa(\cdot,x)) \, \d\mu(x) \\
    &= \int_{\set{X}} \sum_{i=1}^\infty \left( \int_{\set{X}} \kappa(x',x) \overline{e_i(x')} \, \d\mu(x')\right) e_i(x) \, \d\mu(x) \\
    &= \sum_{i=1}^\infty \int_{\set{X}} \left( \int_{\set{X}} \kappa(x',x) e_i(x) \, \d\mu(x)\right) \overline{e_i(x')} \, \d\mu(x') \\
    &= \sum_{i=1}^\infty \langle Ae_i, e_i\rangle_{\Ltwo(\mu)} = \sum_{i=1}^\infty \lambda_i = \tr(A).
\end{align*}
A rigorous proof of \cref{eq:trace-equals-integral} is provided in \cite[Thm.~3.10]{SS12}.

From the eigendecomposition \cref{eq:operator-eval-decomp}, we can define the operator square root
\begin{equation*}
    A^{1/2} \coloneqq \sum_{i=1}^\infty \lambda_i^{1/2} e_i^{\vphantom{*}} e_i^* : \Ltwo(\mu) \to \Ltwo(\mu).
\end{equation*}
In fact, the range of the linear transformation $A^{1/2}$ is $\RKHS$ and the map $A^{1/2} : \Ltwo(\mu) \to \RKHS$ is an isometric embedding.
To see this, first take $u \in \RKHS$ and compute
\begin{equation*}
    \norm{u}_{\Ltwo(\mu)}^2 = \langle u,u\rangle_{\Ltwo(\mu)} = \langle Au,u\rangle_\RKHS = \langle A^{1/2}u,A^{1/2}u\rangle_\RKHS = \norm{\smash{A^{1/2}u}}_\RKHS^2.
\end{equation*}
The second equality is an invocation of the adjoint relation $\langle u, f\rangle_{\Ltwo(\mu)} = \langle Au, f\rangle_\RKHS$.
Since $\RKHS$ is dense in $\Ltwo(\mu)$, we conclude that $A^{1/2} : \Ltwo(\mu) \to \RKHS$ is an isometric embedding.

Since $A^{1/2} : \Ltwo(\mu) \to \RKHS$ is an isometric embedding, it follows that $\{A^{1/2}e_i\}_{i=1}^\infty = \{ \lambda_i^{1/2} e_i \}_{i=1}^\infty$ is an orthonormal system in $\RKHS$.
In fact, it is an orthonormal basis, and $A^{1/2}$ is an isometry \cite[Thm.~3.1]{SS12}.

The expression \cref{eq:operator-eval-decomp} constitutes an eigenvalue decomposition for the integral operator $A$ associated with the kernel $\kappa$, but what about $\kappa$ itself?
Does it have such an eigendecomposition?
Proceeding formally again, we compute
\begin{multline*}
    \int_{\set{X}} \kappa(x, x') u(x') \, \d\mu(x') = Au(x) = \sum_{i=1}^\infty \lambda_i^{\vphantom{*}} e_i(x) \langle e_i,u\rangle_{\Ltwo(\mu)} \\
    = \sum_{i=1}^\infty \int_{\set{X}} \lambda_i e_i(x) \overline{e_i(x')} u(x') \, \d\mu(x')
    = \int_{\set{X}} \left( \sum_{i=1}^\infty \lambda_i e_i(x) \overline{e_i(x')}\right) u(x') \, \d\mu(x').
\end{multline*}
Since this equation holds for all $x$ and all $u$, it is natural to conjecture that
\begin{equation*}
    \kappa(x,x') = \sum_{i=1}^\infty \lambda_i e_i(x) \overline{e_i(x')}.
\end{equation*}
Indeed, this \emph{Mercer decomposition}\index{Mercer decomposition} holds true, and the convergence holds pointwise \cite[Thm.~3.1]{SS12}.\index{reproducing kernel Hilbert space!relationship to $\Ltwo$ space|)}\index{L2, function space@$\Ltwo$, function space|)}\index{integral operator!for a reproducing kernel Hilbert space|)}\index{reproducing kernel Hilbert space!integral operator|)}

\index{positive-semidefinite low-rank approximation!in infinite dimensions|(}\index{subset selection problem!in infinite dimensions|(}
\section{Infinite-dimensional psd low-rank approximation and \RPCholesky} \label{sec:infinite-rpcholesky}

Having established the proper infinite-dimensional setting, we can now formulate the appropriate low-rank approximation problem and devise an infinite-dimensional version of randomly pivoted Cholesky\index{randomly pivoted Cholesky} to solve it.

\index{Nystr\"om approximation!of a kernel function|(}
As we saw in \cref{sec:kernel-interpolation}, a subset $\set{S} \subseteq \set{X}$ induces a rank-$k$ approximation
\begin{equation*}
    \hat{\kappa}_{\set{S}}(x,x') \coloneqq \kappa(x,\set{S})\kappa(\set{S},\set{S})^\dagger \kappa(\set{S},x')
\end{equation*}
to the kernel function $\kappa$.
It also defines the residual
\begin{equation*}
    \kappa_{\set{S}} = \kappa - \hat{\kappa}_{\set{S}}.
\end{equation*}
These objects suggest an infinite-dimensional version of the psd low-rank approximation and column subset selection problems:

\actionbox{\textbf{Infinite-dimensional column subset selection problem.} Find a subset $\set{S} \subseteq \set{X}$ of size $k$ that approximately minimizes the trace error
\begin{equation*}
    \int_{\set{X}} \kappa_{\set{S}}(x,x) \, \d\mu(x) = \int_{\set{X}} \left[ \kappa(x,x) - \kappa(x,\set{S})\kappa(\set{S},\set{S})^\dagger \kappa(\set{S},x) \right] \, \d\mu(x).
\end{equation*}
}
\index{Nystr\"om approximation!of a kernel function|)}

\index{Nystr\"om approximation!of an integral operator|(}\index{integral operator!Nystr\"om approximation of|(}
Using the formalism developed in the previous section, we can recast this problem in terms of integral operators.
Indeed, introduce the Nystr\"om-approximate operator
\begin{equation} \label{eq:nystrom-integral-operator}
    \hat{A}_{\set{S}}u \coloneqq \int_{\set{X}} \hat{\kappa}_{\set{S}}(\cdot,x) u(x) \, \d \mu(x).
\end{equation}
Since the residual kernel $\kappa_{\set{S}}(x,x) = \kappa(x,x) - \hat{\kappa}_{\set{S}}(x,x)$ is psd and integrable along the diagonal (\cref{prop:nystrom-kernel}), the residual operator $A - \hat{A}_{\set{S}}$ is psd, and the trace error is
\begin{equation} \label{eq:trace-error-residual-kernel}
    \tr(A - \hat{A}_{\set{S}}) = \int_{\set{X}} \kappa_{\set{S}}(x,x) \, \d\mu(x).
\end{equation}
Thus, choosing a subset $\set{S}$ to control the integrated diagonal of the residual $\kappa_{\set{S}}$ is equivalent to construction a Nystr\"om approximation of the integral operator $A$ with a small trace error $\tr(A - \hat{A}_{\set{S}})$.\index{Nystr\"om approximation!of an integral operator|)}\index{integral operator!Nystr\"om approximation of|)}\index{positive-semidefinite low-rank approximation!in infinite dimensions|)}\index{subset selection problem!in infinite dimensions|)}

\index{infinite-dimensional randomly pivoted Cholesky|(}
Given this setup, the infinite-dimensional version of \RPCholesky is natural.

\actionbox{\textbf{\RPCholesky on a general space.} Begin from the empty subset $\set{S}_0 \coloneqq \emptyset$, iterate for $i=0,1,2\ldots,k-1$:
\begin{enumerate}
    \item \textit{\textbf{Sample}} a random pivot
    \begin{equation} \label{eq:infinite-rpc-sample}
        s_{i+1} \sim \kappa_{\set{S}_i}(x,x) \, \d \mu(x)
    \end{equation}
    from the diagonal of the residual kernel $\kappa_{\set{S}_i}(x,x)$.
    Here, $s \sim \nu$ to denotes a random element sampled from an \warn{unnormalized} finite measure $\nu$. 
    That is, $\prob \{ s \in \set{B} \} = \nu(\set{B}) / \nu(\set{X})$ for every measurable set $\set{B}$.
    \item \textit{\textbf{Induct}} the sampled pivot into the landmark set, $\set{S}_{i+1} \coloneqq \set{S}_i \cup \{s_{i+1}\}$.
\end{enumerate}}

This procedure is easy to state, but it is unclear how to implement the sampling step \cref{eq:infinite-rpc-sample}.
We will turn to the issue of implementation in \cref{sec:rpc-rejection}.

\begin{figure}[t]
    \centering
    \includegraphics[width=0.48\linewidth]{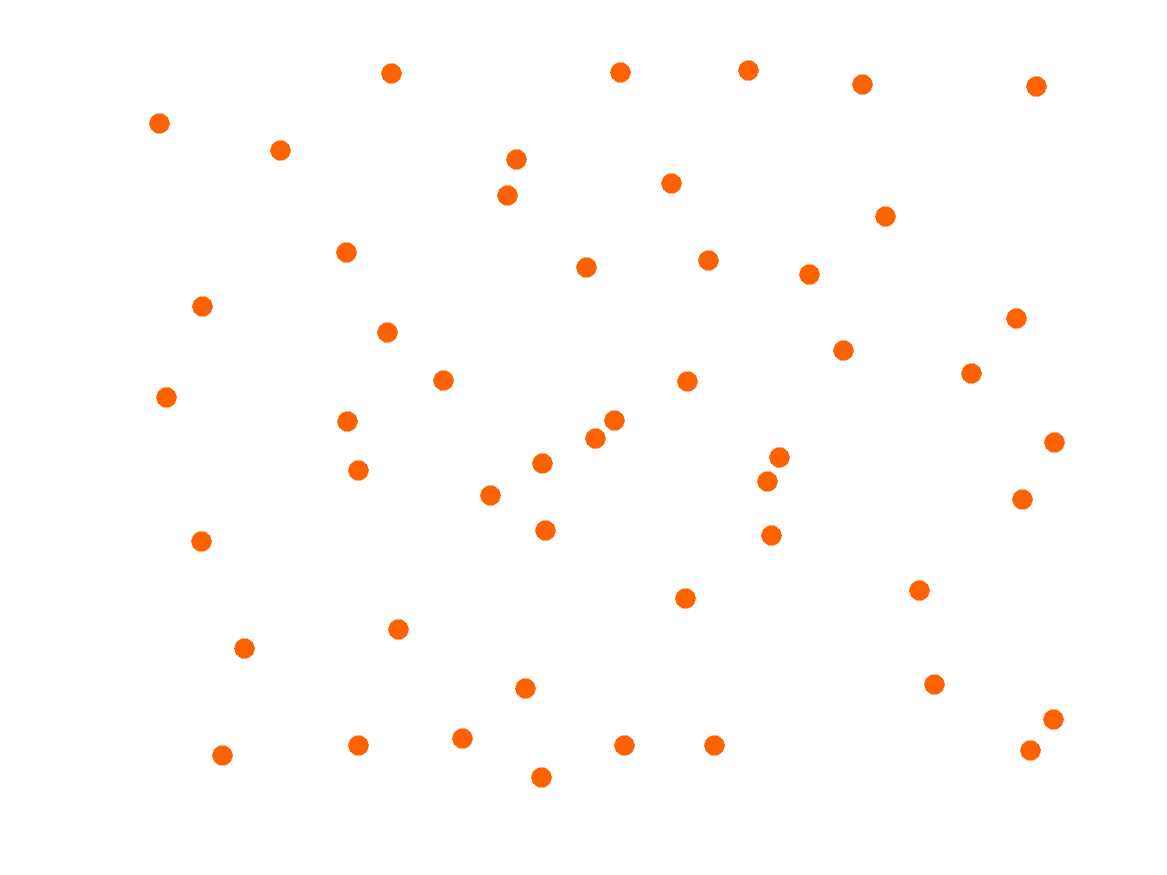}
    \includegraphics[width=0.48\linewidth]{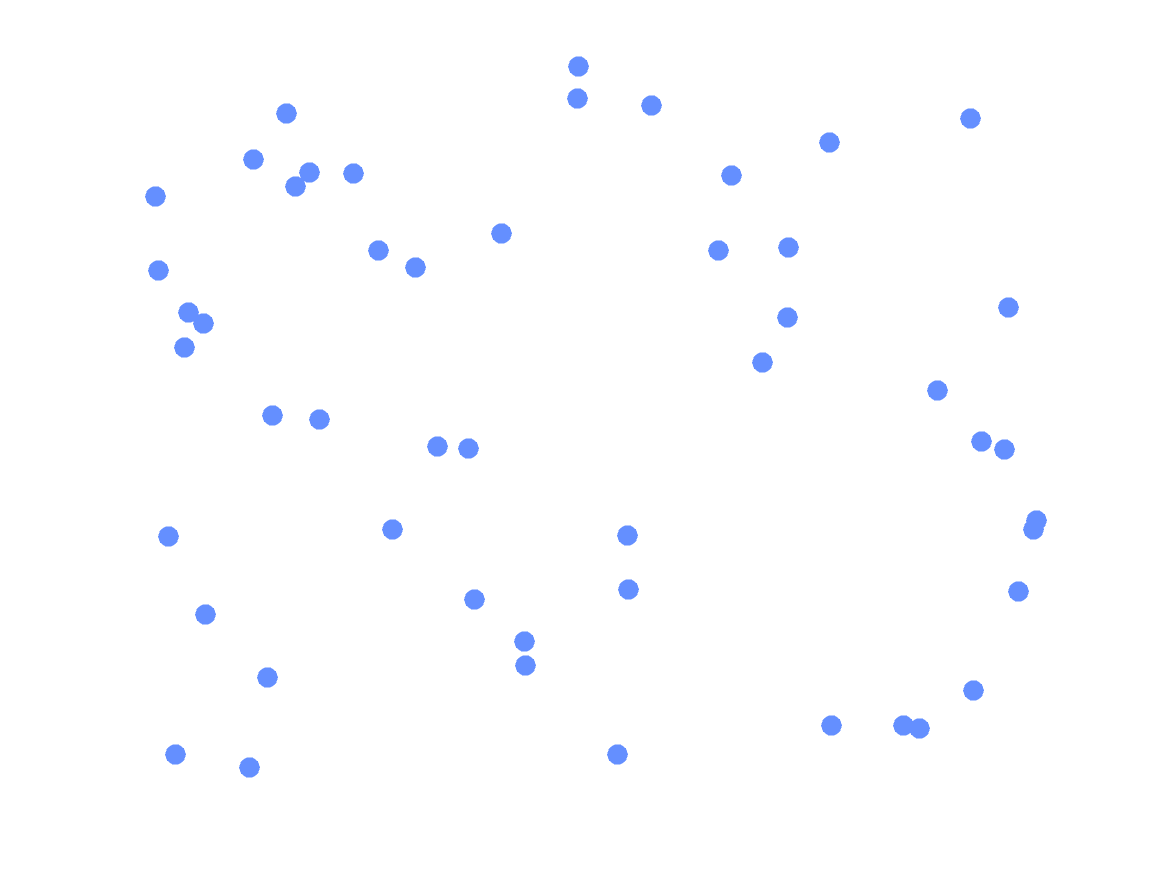}
    \caption[Illustration of random points in the unit square selected by continuous \RPCholesky and uniformly at random]{50 random points in the unit square generated by \RPCholesky (\emph{left}) and uniform sampling (\emph{right}). Points sampled using \RPCholesky are more spread out, demonstrating the repulsive nature of the \RPCholesky sampling process.
    \RPCholesky is implemented using the uniform measure $\mu = \Unif([0,1]^2)$ and the Laplace kernel with bandwidth $10$.}
    \label{fig:continuous-rpcholesky}
\end{figure}

An illustration of the landmarks selected by \RPCholesky on a continuous space is given in \cref{fig:continuous-rpcholesky}.
We see that the points produced by \RPCholesky are more well-distributed across the region, whereas points generated by uniform sampling are more clustered.
This example demonstrates that the \RPCholesky sampling process is repulsive; once a landmark is selected, it is less likely that nearby points will be selected as landmarks.

\index{infinite-dimensional randomly pivoted Cholesky!theoretical analysis|(}
The error bounds from \cref{sec:rpc-analysis} (\cref{thm:rpcholesky-trace,thm:rpcholesky-spectral}) generalize to infinite-dimensional \RPCholesky in the natural way.
The proofs are the same, modulo appropriate details to handle the infinite-dimensional setting; we omit the details.

\begin{theorem}[Continuous \RPCholesky: Trace error] \label{thm:cont-rpcholesky-trace}
    Instate the prevailing notation and assumptions, let $r\ge 1$ and $\varepsilon > 0$, and introduce the relative error of the best rank-$r$ approximation
    \begin{equation*}
        \eta \coloneqq \frac{\sum_{i=r+1}^\infty \lambda_i}{\sum_{i=1}^\infty \lambda_i}.
    \end{equation*}
    The continuous \RPCholesky method is guaranteed to produce a pivot set $\set{S}$ achieving the guarantee
    \begin{equation*}
        \expect \tr(A - \hat{A}_{\set{S}}) \le (1+\varepsilon) \sum_{i=r+1}^\infty \lambda_i
    \end{equation*}
    after running for any number of steps $k$ satisfying
    \begin{equation*}
        k \ge \frac{r}{\varepsilon} + r \log \left( \frac{1}{\varepsilon \eta} \right).
    \end{equation*}
\end{theorem}

\begin{theorem}[Continuous \RPCholesky: Spectral norm of expected error] \label{thm:cont-rpcholesky-spectral}
    Instate the prevailing notation and assumptions, let $r\ge 1$ and $\varepsilon > 0$, and introduce the relative error of the best rank-$r$ approximation
    \begin{equation*}
        \eta \coloneqq \frac{\sum_{i=r+1}^\infty \lambda_i}{\sum_{i=1}^\infty \lambda_i}.
    \end{equation*}
    The continuous \RPCholesky method is guaranteed to produce a pivot set $\set{S}$ achieving the guarantee
    \begin{equation*}
        \lambda_{\mathrm{max}}(\expect(A - \hat{A}_{\set{S}})) \le b + \varepsilon \sum_{i=r+1}^\infty \lambda_i
    \end{equation*}
    after running for any number of steps $k$ satisfying
    \begin{equation*}
        k \ge \frac{1}{\varepsilon} + r \log \left( \frac{\lambda_1}{b} \right).
    \end{equation*}
    Here, $\lambda_{\mathrm{max}}(\cdot)$ denotes the largest eigenvalue of a Hermitian operator on $\Ltwo(\mu)$.
    In particular, for such $k$ and any function $g \in \Ltwo(\mu)$, 
    \begin{equation*}
        \expect \langle g, (A - \hat{A}_{\set{S}})g \rangle_{\Ltwo(\mu)} \le \left( b + \varepsilon \sum_{i=r+1}^\infty \lambda_i \right) \norm{g}_{\Ltwo(\mu)}^2.
    \end{equation*}
\end{theorem}

\Cref{thm:cont-rpcholesky-spectral} originally appeared as \cite[Thm.~7]{EM23a}.
The infinite-dimensional trace error bound \cref{thm:cont-rpcholesky-trace} has not appeared in print before; it is an immediate transplant of the finite-dimensional result \cref{thm:rpcholesky-trace}, originally proven in \cite[Thm.~1.1]{CETW25}, to the infinite setting.\index{infinite-dimensional randomly pivoted Cholesky!theoretical analysis|)}

\section{Implementing \RPCholesky using rejection sampling} \label{sec:rpc-rejection}

\index{infinite-dimensional randomly pivoted Cholesky!implementation by rejection sampling|(}\index{randomly pivoted Cholesky!implementation by rejection sampling|(}
How could we implement the \RPCholesky method in infinite dimensions? 
Nominally, storing even a single ``column'' $\kappa(\cdot,s)$ of the ``infinite kernel matrix'' $\kappa$ requires infinite memory, which makes a direct implementation of infinite \RPCholesky similar to \cref{prog:rpcholesky} impossible in general.
\index{rejection sampling!for general distributions|(}
We circumvent this limitation by implementing \RPCholesky in infinite dimensions using an alternative approach based on \emph{rejection sampling}.\index{infinite-dimensional randomly pivoted Cholesky!implementation by rejection sampling|)}\index{randomly pivoted Cholesky!implementation by rejection sampling|)}

\subsection{Rejection sampling in general}

\myprogram{General implementation of rejection sampling from a target distribution $\tau(x) \, \d\mu(x)$ using a proposal distribution $\pi(x) \, \d\mu(x)$.}{}{rejection_sample}

Rejection sampling is a standard technique in probabilistic computing, originally due to von Neumann \cite{von51}.
See \cite[\S2.2]{Liu04} for a contemporary introduction to rejection sampling.
Suppose that we are interested in sampling $s \sim \tau(x) \, \d\mu(x)$ from an (unnormalized) target distribution $\tau(x) \, \d\mu(x)$.
Assume the target distribution is difficult to sample from directly, but that we have easy access to samples from an (unnormalized) \emph{proposal distribution}\index{proposal distribution} $p \sim \pi(x) \, \d\mu(x)$.
Suppose additionally that we know that the proposal distribution bounds the target distribution pointwise, $0\le \tau \le \pi$.\index{target distribution}
With these preliminaries, the rejection sampling procedure is as follows.
\actionbox{\textbf{Rejection sampling.} Repeat until termination:
\begin{enumerate}
    \item \textbf{\textit{Generate}} a sample $p \sim \pi(x) \, \d\mu(x)$ from the proposal distribution.\index{proposal distribution}
    \item \textbf{\textit{Accept/reject?}} With probability $\tau(p)/\pi(p)$, \emph{accept} by setting $s \coloneqq p$ and exiting this loop.
    Otherwise, \emph{reject} the sample and repeat step 1.
\end{enumerate}}
A program implementing rejection sampling is given as \cref{prog:rejection_sample}.

The correctness of the rejection sampling procedure is established by the following standard result, in essence due to von Neumann \cite{von51}.

\begin{fact}[Rejection sampling]
    The rejection sampling procedure works:
    Provided $\tau > 0$ on a set of positive measure, the procedure terminates with probability $1$, and it outputs a sample $s \sim \tau(x) \, \d\mu(x)$.
\end{fact}

The proof is short and simple; see \cite[\S2.2]{Liu04} for a proof when $\mu$ is the Lebesgue measure on $\real^n$. (The proof is no different for general probability spaces.)\index{rejection sampling!for general distributions|)}

\index{infinite-dimensional randomly pivoted Cholesky!implementation by rejection sampling|(}\index{randomly pivoted Cholesky!implementation by rejection sampling|(}\index{rejection sampling!for implementing randomly pivoted Cholesky|(}
\subsection{Rejection sampling for \RPCholesky}

We can use rejection sampling to perform the diagonal sampling step for \RPCholesky.
We assume access to two primitive operations:
\begin{enumerate}
    \item \textbf{\textit{Kernel evaluations.}} We can evaluate the $\kappa(x,x')$ at any pair of points $x,x'\in\set{X}$.
    \item \textbf{\textit{Sampling from the diagonal.}} We have access to an efficient procedure for sampling $p \sim \kappa(x,x) \, \d\mu(x)$.
\end{enumerate}
In typical machine learning applications, these primitives are easily accessible, as (1) the kernel function $\kappa$ is often given by an explicit functional form and (2) the kernel function often satisfies $\kappa(x,x) = 1$ for all $x$, so sampling from the diagonal amounts to uniform sampling from $\mu$.
As we will see, access to these two weak primitives will suffice to implement \RPCholesky in infinite dimensions.

Consider step $i$ of \RPCholesky.
We seek a sample $s_{i+1} \sim \kappa_{\set{S}_i}(x,x)\,\d\mu(x)$ from the target distribution\index{target distribution} $\kappa_{\set{S}_i}(x,x)\,\d\mu(x)$.
To do so, we will employ rejection sampling with $\kappa(x,x) \, \d\mu(x)$ as the proposal distribution.\index{proposal distribution}
These distributions form a valid pair for rejection sampling in view of the inequality
\begin{equation*}
    0 \le \kappa_{\set{S}_i}(x,x) \le \kappa(x,x) \quad \text{for every } x \in \set{X},
\end{equation*}
which holds due to \cref{prop:nystrom-kernel}.

\myprogram{Rejection sampling-based implementation of \RPCholesky on general spaces.}{}{rejection_rpcholesky}

Since we have efficient access to proposals $p \sim \kappa(x,x) \, \d\mu(x)$, implementing this procedure just requires the ability to evaluate the ratio
\begin{equation*}
    \mathrm{ratio}_i(x) \coloneqq \frac{\kappa_{\set{S}_i}(x,x)}{\kappa(x,x)} = \frac{\kappa(x,x) - \hat{\kappa}_{\set{S}_i}(x,x)}{\kappa(x,x)}.
\end{equation*}
To do so, we maintain a Cholesky decomposition $\kappa(\set{S}_i,\set{S}_i) = \mat{R}_i^*\mat{R}_i^{\vphantom{*}}$ of the currently selected kernel submatrix and evaluate $\mathrm{ratio}_i(x)$ using the formula
\begin{equation} \label{eq:rpcholesky-rejection-ratio}
    \mathrm{ratio}_i(x) = \frac{\kappa(x,x) - \norm{\vec{c}_i(x)}^2}{\kappa(x,x)} \quad \text{where } \vec{c}_i(x) \coloneqq \mat{R}_i^{-*} \kappa(\set{S}_i,x).
\end{equation}
Once a new pivot $s_{i+1}$ is accepted, the Cholesky factor can be updated using the formula
\begin{equation*}
    \mat{R}_{i+1} \coloneqq \twobytwo{\mat{R}_i}{\vec{c}_i(s_{i+1})}{\mat{0}}{(\kappa(s_{i+1},s_{i+1}) - \norm{\vec{c}_i(s_{i+1})}^2)^{1/2}} .
\end{equation*}
This rejection sampling-based implementation of \RPCholesky is called \RPCholeskyRejection. 
An implementation is given in \cref{prog:rejection_rpcholesky}.

\index{infinite-dimensional randomly pivoted Cholesky!computational cost|(}
\subsection{Computational cost}

The computational cost of \RPCholeskyRejection can be characterized in multiple ways.
Each execution of the rejection sampling loop requires at most $k$ kernel function evaluations and $\order(k^2)$ operations to perform the triangular solve in \cref{eq:rpcholesky-rejection-ratio}.
Thus, $b$ proposals of \RPCholeskyRejection expends at most $bk$ kernel function evaluations and $\order(bk^2)$ additional arithmetic operations, where $k \le b$ denotes the number of accepted proposals.
An end-to-end cost estimate is provided by the following result:

\begin{proposition}[\RPCholeskyRejection] \label{prop:rpcholesky-rejection}
    \Cref{prog:rejection_rpcholesky} is correct: It produces exact samples from the continuous \RPCholesky algorithm.
    To measure the computational cost, introduce the relative trace error of the best rank-$(i-1)$ approximation:
    \begin{equation*}
        \eta_i \coloneqq \frac{\sum_{j = i}^\infty \lambda_i}{\sum_{j = 1}^\infty \lambda_i}.
    \end{equation*}
    The expected computational cost of \cref{prog:rejection_rpcholesky} is at most
    \begin{align*}
        &\order \left( \sum_{i=1}^k \frac{1}{\eta_i} \right) \le \order\left( \frac{k}{\eta_k} \right)&&\text{rejection sampling steps}, \\
        &\order \left( \sum_{i=1}^k \frac{i}{\eta_i} \right) \le \order\left( \frac{k^2}{\eta_k} \right)&&\text{kernel function evaluations}, \\
        &\order \left( \sum_{i=1}^k \frac{i^2}{\eta_i} \right) \le \order\left( \frac{k^3}{\eta_k} \right)&&\text{additional arithmetic operations}.
    \end{align*}
\end{proposition}

This result is \cite[Thm.~4]{EM23a}.
The main idea is that, to sample $s_{i+1}$, the probability of each proposal being accepted is $\tr(A - \hat{A}_{\{s_1,\ldots,s_i\}})/\tr(\mat{A}) \ge \eta_i$.
Ergo, the expected number of proposals is a geometric random variable with expectation $1/\eta_i$, and the result follows.

\index{curse of smoothness|(}\index{spectral decay|(}
This result demonstrates that rejection sampling-based \RPCholesky suffers from a \emph{curse of smoothness}: the runtime of the algorithm is higher when the eigenvalues $\lambda_i$ decrease at a rapid rate.
The reason for this curse is simple.
When the eigenvalues decay rapidly, the Nystr\"om approximations $\hat{\kappa}_{\set{S}_i}$ become increasingly close to the kernel $\kappa$.
As such, the acceptance probability \cref{eq:rpcholesky-rejection-ratio} goes down, requiring more proposals to achieve each acceptance.

Elvira Moreno and I proposed a potentially costly fix for the curse of smoothness in \cite[Alg.~4]{EM23a}.
Here is the idea:
When the algorithm begins rejecting proposals at a high rate (say, 10 rejections in a row), solve a global minimization problem
\begin{equation*}
    \alpha \gets \max_{x \in \set{X}} \mathrm{ratio}_i(x),
\end{equation*}
where $\mathrm{ratio}_i(x)$ was defined in \cref{eq:rpcholesky-rejection-ratio}.
Then, run rejection sampling with the rescaled proposal distribution $\alpha^{-1} \kappa(x,x) \, \d\mu(x)$.
This rescaling does not effect the distribution of proposals $p \sim \alpha^{-1} \kappa(x,x) \, \d\mu(x)$, but the acceptance probabilities are scaled by a factor $\alpha^{-1}$.
Since $\alpha \le 1$, this has the effect of boosting the acceptance probabilities.
This procedure can be very effective; our paper \cite{EM23a} reports speedups of $39\times$ over the non-accelerated version.
However, the accelerated \RPCholeskyRejection procedure requires solving a global nonconvex optimization problem over the domain $\set{X}$, which can be computationally costly.
The development of tractable rejection sampling schemes for \RPCholesky that avoid nonconvex optimization is an interesting direction for future research.\index{infinite-dimensional randomly pivoted Cholesky!implementation by rejection sampling|)}\index{randomly pivoted Cholesky!implementation by rejection sampling|)}\index{rejection sampling!for implementing randomly pivoted Cholesky|)}\index{infinite-dimensional randomly pivoted Cholesky!computational cost|)}\index{curse of smoothness|)}\index{spectral decay|)}


\index{infinite-dimensional randomly pivoted Cholesky!applied to active kernel interpolation|(}\index{active learning!for kernel regression|(}\index{kernel ridge regression!active learning|(}\index{kernel interpolation!active learning|(}
\section{Application: Active kernel interpolation} \label{sec:infinite-active-kernel-interp}

Consider the active regression task of approximating a function $f \in \RKHS$ whose value we may be query at $k$ locations $\set{S} = \{s_1,\ldots,s_k\} \subseteq \set{X}$, which we are free to choose.
A natural approach to this problem is to select the locations $\set{S}$ using \RPCholesky, then interpolate the function values at those points using kernel interpolation.
Throughout this section, we let $\hat{f}_{\set{S}}$ denote the kernel interpolant through the points $\set{S}$.

\index{kernel interpolation!mean-squared error bounds|(}
The success of this active kernel interpolation strategy can be analyzed by combining the following known result for kernel interpolation with the trace error bound \cref{thm:cont-rpcholesky-trace} for \RPCholesky:

\begin{proposition}[Kernel interpolation: mean-squared error] \label{fact:kernel-interpolation-error}
    Let $f \in \RKHS$ be a function and let $\hat{f}_{\set{S}}$ be its kernel interpolant through the points $\set{S}$.
    We have the error bound
    \begin{equation*}
        \norm{\smash{f - \hat{f}_{\set{S}}}}_{\Ltwo(\mu)}^2 \le \norm{f}_\RKHS^2 \cdot \tr(A - \hat{A}_{\set{S}}).
    \end{equation*}
\end{proposition}

\begin{proof}
    The bound follows by integrating the pointwise bound \cref{thm:kernel-interpolation-error} and using the identity \cref{eq:trace-error-residual-kernel} for $\tr(A - \hat{A}_\set{S})$.
\end{proof}

Combining with \cref{thm:cont-rpcholesky-trace} yields the following immediate corollary:

\begin{corollary}[Kernel interpolation with \RPCholesky] \label{cor:kernel-interpolation-rpcholesky}
    With the setting and value of $k$ from \cref{thm:cont-rpcholesky-trace}, a set $\set{S}$ picked by $k$ steps of \RPCholesky satisfies the bound
    \begin{equation*}
        \expect \norm{\smash{f - \hat{f}_{\set{S}}}}_{\Ltwo(\mu)}^2 \le \norm{f}_\RKHS^2 \cdot (1+\varepsilon) \sum_{i=r+1}^\infty \lambda_i.
    \end{equation*}
    Moreover, we have an estimate that is uniform over all functions in the RKHS:
    \begin{equation*}
        \expect \sup \left\{ \frac{\norm{\smash{f - \hat{f}_{\set{S}}}}_{\Ltwo(\mu)}^2}{\norm{f}_\RKHS^2} : f \in \RKHS \setminus \{0\} \right\} \le (1+\varepsilon) \sum_{i=r+1}^\infty \lambda_i.
    \end{equation*}
\end{corollary}
\index{kernel interpolation!mean-squared error bounds|)}

\index{randomly pivoted Cholesky!numerical results|(}
\begin{figure}
    \centering
    \includegraphics[height=2in]{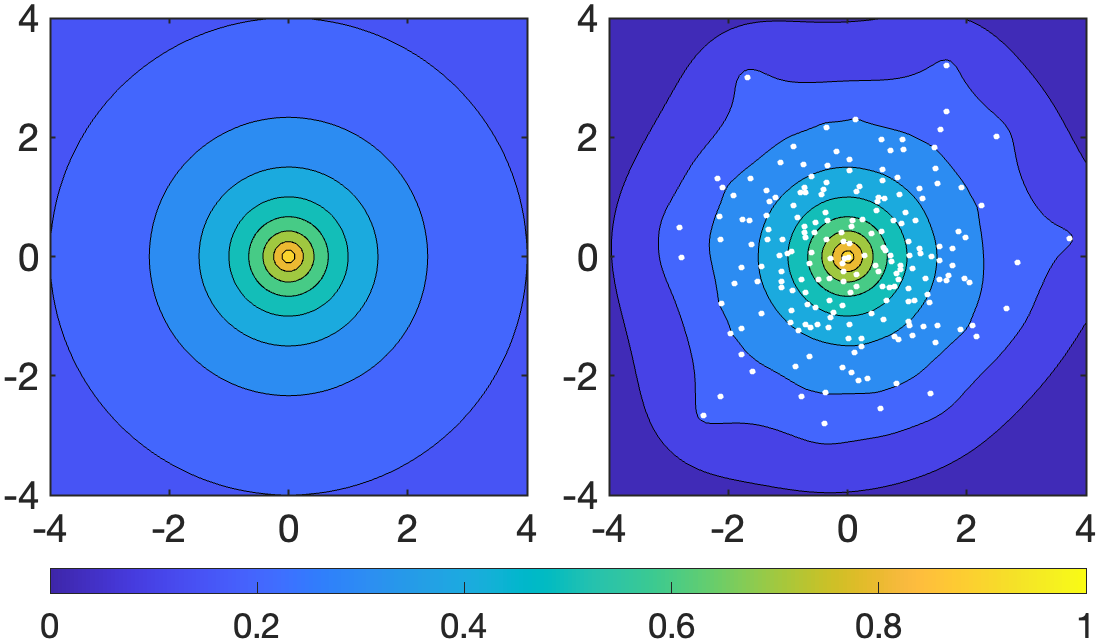}
    \includegraphics[height=2in]{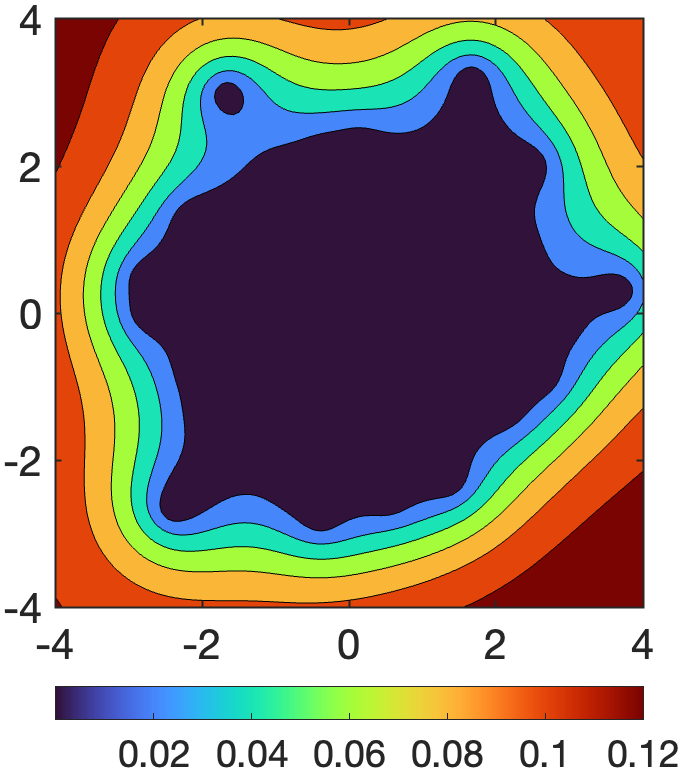}
    
    \includegraphics[width=0.75\linewidth]{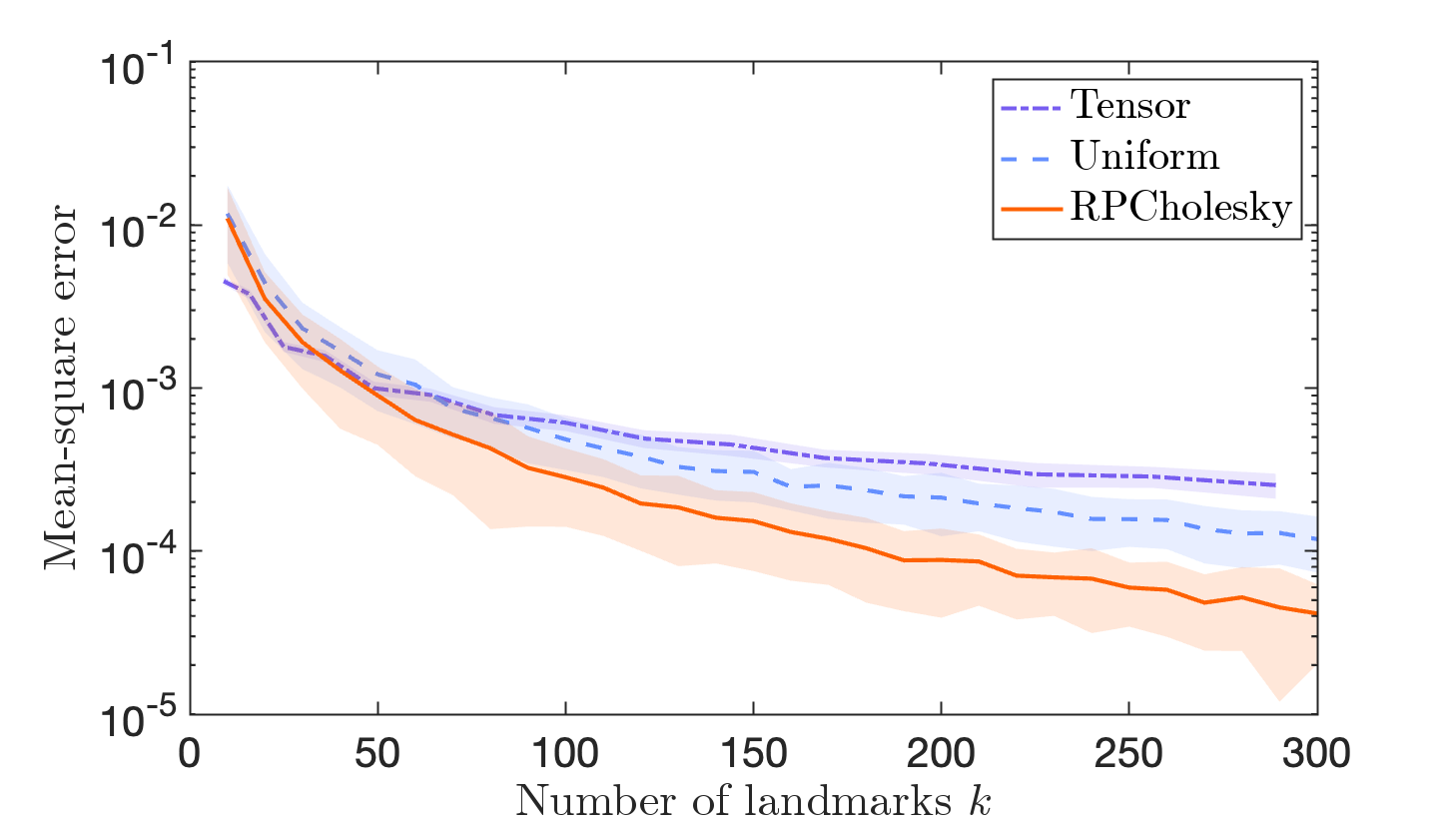}
    
    \caption[Active kernel interpolation using continuous \RPCholesky, iid random selection, and a tensor-product grid]{\emph{Top:} Function $f$ (\emph{left}), kernel interpolant $\hat{f}_{\set{S}}$ with \RPCholesky (\emph{middle}), and error $|f - \hat{f}_{\set{S}}|$ (\emph{right}).
    Landmarks selected by \RPCholesky are overlayed as white dots in the middle plot.
    Function $f$, kernel $\kappa$, and measure $\mu$ are described in the text.
    \emph{Bottom:} mean-squared error of kernel interpolation through three different sets of nodes.
    Lines show a mean of 100 trials, and shaded regions show one standard deviation.}
    \label{fig:active-kernel-interpolation}
\end{figure}

An example of the active interpolation methdology is shown in \cref{fig:active-kernel-interpolation}.
Here, we apply kernel interpolation to approximate the function
\begin{equation*}
    f(\vec{x}) = (1+\norm{\vec{x}})^{-1}
\end{equation*}
on the domain $\set{X} \coloneqq \real^2$.
We use the standard Gaussian measure $\mu = \Normal_\real(\vec{0},\Id_2)$ and set $\kappa$ to be the Laplace kernel with bandwidth $1$.

The top panels in \cref{fig:active-kernel-interpolation} illustrate the approximation produced by active kernel interpolation with \RPCholesky-selected landmarks.
The left panel plots the target function $f$ over to the finite window $[-4,4]^2$.
The middle panel shows the kernel interpolant $\hat{f}_{\set{S}}$ computing by interpolating through $200$ points selected by \RPCholesky, shown as white dots.
The right panel shows the absolute error $|f - \hat{f}_{\set{S}}|$.
\RPCholesky kernel interpolation produces a small error near the origin, with growing errors farther away.
This behavior is expected, as the measure $\mu$ places most of its mass near the center, leading most of the landmarks to be placed there.
Fortunately, higher error in regions of small measure is usually acceptable in applications.

The bottom panel of \cref{fig:active-kernel-interpolation} compares the mean-squared error of three different sets of kernel interpolation through three sets of landmarks: uniformly random points $s_1,\ldots,s_k \stackrel{\text{iid}}{\sim} \mu$, points selected by \RPCholesky, and points from a tensor product grid $\set{S} = \set{E} \times \set{E}$ where $\set{E} = \{ \Phi^{-1}((2i-1)/2\sqrt{k}) : i = 1,\ldots,k\}$.
Here, $\Phi$ denotes the cdf of the standard Gaussian distribution $\Normal_\real(0,1)$.
We estimate the mean-squared error $\norm{f - \smash{\hat{f}_{\set{S}}}}_{\Ltwo(\mu)}^2$ using the Monte Carlo estimator $s^{-1} \sum_{i=1}^s (f(x_i) - \hat{f}_{\set{S}}(x_i))^2$ where $x_1,\ldots,x_s \stackrel{\text{iid}}{\sim} \mu$; we use $s=10^3$ for these experiments.
For each value of $k$, we display the mean of $100$ trials; error bars show one standard deviation.
The tensor product nodes are deterministic, so the error bars capture only the Monte Carlo fluctuations in the estimated mean-squared error, which are observed to be small.
The \RPCholesky method achieves the lowest error of these three point sets.
\index{randomly pivoted Cholesky!numerical results|)}

\begin{remark}[Beyond kernel interpolation]
    The continuous \RPCholesky method provides a general strategy for selecting points to label during active learning.
    While we have used kernel interpolation as the learning method in this section, one can also combine \RPCholesky with any other functional regression technique.
    Kernel ridge regression is one natural choice, but one can also use nonlinear approximation methods such as artificial neural networks.
\end{remark}

\index{active learning!for kernel regression|)}\index{kernel ridge regression!active learning|)}\index{infinite-dimensional randomly pivoted Cholesky!applied to active kernel interpolation|)}\index{kernel interpolation!active learning|)}

\index{kernel quadrature|(}
\section{Application: Kernel quadrature} \label{sec:quadrature}

Another task we can apply \RPCholesky to is numerical quadrature.
Given a weight function $u \in \Ltwo(\mu)$, the task is to devise weights $\vec{w} \in \field^k$ and nodes $s_1,\ldots,s_k \in \set{X}$ so that
\begin{equation*}
    \int_{\set{X}} \overline{u(x)}f(x) \, \d\mu(x) \approx \sum_{i=1}^k \overline{w_i} f(s_i) \quad \text{for all } f \in \RKHS.
\end{equation*}
To simplify notation in later parts of this section, I have made the unorthodox choice to complex conjugate both the weight function $u$ and the quadrature weights $w_i$ in this expression.
The quality of a quadrature scheme (i.e., both weights and nodes) for a given function $f$ is quantified by the error:
\begin{equation} \label{eq:quadrature-error}
    \Err(\set{S},\vec{w},f) \coloneqq \left| \int_{\set{X}} \overline{u(x)}f(x) \, \d\mu(x) - \sum_{i=1}^k \overline{w_i} f(s_i) \right|.
\end{equation}
We will use \RPCholesky to pick the nodes $\set{S} = \{s_1,\ldots,s_k\}$, and there is a standard procedure for picking the weights $\vec{w}$.
We begin by reviewing these ideal weights and then prove error bounds for optimally weighted \RPCholesky quadrature.

\index{kernel quadrature!characterizations of ideal weights|(}
\myparagraph{The ideal weights}
There are five equivalent ways of characterizing the ideal weights, which I summarized in \cite{Epp23b}.
We quickly derive the ideal weights and review these interpretations, as each sheds different light on the kernel quadrature problem.

\mysubparagraph{Characterization \#1: Minimizing the worst-case error}
Under this interpretation, we choose $\vec{w}$ to minimize the worst-case value of the error \cref{eq:quadrature-error} over all $f \in \RKHS$ with $\norm{f}_\RKHS \le 1$.
To derive the worst-case error, use the integral operator $A$ and the reproducing property to rewrite the error \cref{eq:quadrature-error} as
\begin{equation} \label{eq:quadrature-error-computation}
    \Err(\set{S},\vec{w},f) = \left| \langle u,f \rangle_{\Ltwo(\mu)} - \sum_{i=1}^k \overline{w_i} \langle \kappa(\cdot,s_i), f \rangle_\RKHS \right| = \left| \langle Au - \kappa(\cdot,\set{S})\vec{w},f\rangle_\RKHS \right|.
\end{equation}
By the Cauchy--Schwarz inequality, the maximum value of this quantity over all $f$ with $\norm{f}_\RKHS \le 1$ is
\begin{equation} \label{eq:worst-case-error-quadrature}
    \Err(\set{S},\vec{w}) \coloneqq \max_{\norm{f}_\RKHS \le 1} \Err(\set{S},\vec{w},f) = \norm{Au - \kappa(\cdot,\set{S})\vec{w}}_\RKHS.
\end{equation}
Minimizing this quantity over $\vec{w}$ is a linear least-squares problem.
The normal equations characterizing the minimizer $\vec{w}_\star$ are
\begin{equation} \label{eq:ideal-weights-system}
    \kappa(\set{S},\set{S}) \vec{w}_\star = Au(\set{S}) = \left( \int_{\set{X}} \kappa(s_i,x) u(x) \, \d\mu(x) : 1\le i\le k \right).
\end{equation}
We call the solution to this system the \emph{ideal weights} $\vec{w}_\star$.

\mysubparagraph{Characterization \#2: Exactness}
A classical way of deriving a quadrature rule on $k$ nodes is to require that it integrate $k$ chosen functions exactly.
Here, the natural choice for $k$ functions are the kernel functions $f_i \coloneqq \kappa(\cdot,s_i)$.
Enforcing exactness $\Err(\set{S},\vec{w},f_i) = 0$ on these $f_i$ yields the same system \cref{eq:ideal-weights-system} characterizing the ideal weights.

\mysubparagraph{Characterization \#3: Interpolate and integrate}
Another classical approach to deriving quadrature rules is to interpolate the function at the nodes, then integrate the interpolant.
Using kernel interpolation as our method of interpolation, this approach yields ideally weighted kernel quadrature:
\begin{equation*}
    \sum_{i=1}^k \overline{\vec{w}_\star(i)} f(s_i) = \int_{\set{X}} \overline{u(x)}  \hat{f}_{\set{S}}(x) \, \d\mu(x).
\end{equation*}

\mysubparagraph{Characterization \#4: Conditional expectation}
Our last two interpretations will drawn on a Gaussian process frame.
Imagine the function $f \sim \GP(\kappa)$ is drawn from a Gaussian process with covariance function $\kappa$.
Under this model, the integral $\int_{\set{X}} f(x) \overline{u(x)} \, \d\mu(x)$ is a random quantity, and we can obtain a quadrature rule by evaluating the conditional expectation of this random quantity conditional on the value $f(s_i)$  of the function at queried positions $s_i$.
Evaluating this conditional expectation also yields kernel quadrature with the ideal weights:
\begin{equation*}
    \sum_{i=1}^k \overline{\vec{w}_\star(i)} f(s_i) = \expect \left[ \int_{\set{X}} \overline{u(x)} f(x) \, \d\mu(x) \, \middle| \, f(\set{S}) \right].
\end{equation*}
Based on this interpretation, ideally weighted kernel quadrature is sometimes called \emph{Bayesian quadrature}.

\mysubparagraph{Characterization \#5: Minimizing the mean-squared error}
As a last notion of how the ``ideal weights'' could be designed, we could also choose the weights to minimize the \emph{mean-squared error} of the integral for a function $f \sim \GP(\kappa)$ drawn from the Gaussian process $\GP(\kappa)$:
\begin{equation*}
    \MSE(\set{S},\vec{w}) \coloneqq \expect \left|\int_{\set{X}} \overline{u(x)} f(x) \, \d\mu(x) - \sum_{i=1}^k \overline{w_i} f(s_i) \right|^2.
\end{equation*}
%
We compute this mean-squared error formally.
Begin by re-expressing the mean-squared error using the $\RKHS$-inner product, following \cref{eq:quadrature-error-computation}:
\begin{equation*}
    \MSE(\set{S},\vec{w}) = \left| \langle Au - \kappa(\cdot,\set{S})\vec{w},f\rangle_\RKHS \right|^2
\end{equation*}
Now, employ a Karhunen-Lo\`eve expansion $f = \sum_{i=1}^\infty z_i (\lambda_i^{1/2} e_i)$.
The functions $\lambda_i^{1/2} e_i$ (defined above in \cref{eq:operator-eval-decomp}) form an orthonormal basis for $\RKHS$, and the random variables $z_1,z_2,\ldots \stackrel{\text{iid}}{\sim} \Normal_\field(0,1)$ are iid and Gaussian.
Similarly, we may decompose $Au - \kappa(\cdot,\set{S})\vec{w} = \sum_{i=1}^\infty b_i (\lambda_i^{1/2}e_i)$ in this orthonormal basis is well.
The sequence $(b_i)_{i=1}^\infty$ is square summable, and it admits the Pythagorean identity $\sum_{i=1}^\infty |b_i|^2 = \norm{Au - \kappa(\cdot,\set{S})\vec{w}}_{\RKHS}^2$.
By orthonormality of $\{\lambda_i^{1/2} e_i\}_{i=1}^\infty$, we compute
\begin{align*}
    \MSE(\set{S},\vec{w}) &= \expect \left| \langle Au - \kappa(\cdot,\set{S})\vec{w},f\rangle_\RKHS \right|^2 = \expect \left| \sum_{i=1}^\infty z_i b_i \right|^2\\ &= \sum_{i=1}^\infty |b_i|^2 =  \norm{Au - \kappa(\cdot,\set{S})\vec{w}}_{\RKHS}^2.
\end{align*}
The mean-squared error $\MSE(\set{S},\vec{w})$ for a function $f \sim \GP(\kappa)$ is equal to square the worst-case error $\Err(\set{S},\vec{w})$ over all $f \in \RKHS$ with $\norm{f}_\RKHS \le 1$, in view of \cref{eq:worst-case-error-quadrature}.
Thus, the same set of ideal weights minimizes both the mean-squared error and worst-case error.
\index{kernel quadrature!characterizations of ideal weights|)}

\index{kernel quadrature!error bounds|(}
\myparagraph{Error bounds for ideally weighted kernel quadrature}
The fact that five different definitions of the ideal weights all coincide is a powerful demonstration of the robustness of the RKHS and Gaussian process formalisms.
Having thoroughly convinced ourselves that the ideal weights given by \cref{eq:ideal-weights-system} are, in five senses, the natural way of weighting quadrature rules on RKHSs (or with respect to Gaussian processes), let us now analysis the error for ideally weighted kernel quadrature.

Introduce the error 
\begin{equation*}
    \Err(\set{S}) \coloneqq\Err(\set{S},\vec{w}_\star)
\end{equation*}
for ideally weighted kernel quadrature.
This quantity expresses both the worst-case quadrature error over a function in the unit ball of $\RKHS$ and the root-mean-squared error of a function drawn from the corresponding Gaussian process (see interpretations \#1 and \#5 above).
Combing formulas \cref{eq:worst-case-error-quadrature,eq:ideal-weights-system}, this quantity may be written
\begin{equation} \label{eq:ideally-weighted-worst-case}
    \Err(\set{S})^2 = \norm{Au - \kappa(\cdot,\set{S})\kappa(\set{S},\set{S})^{-1}Au(\set{S})}_\RKHS^2.
\end{equation}
(We assume here and throughout this section that $\kappa(\set{S},\set{S})$ is nonsingular.)
Introducing the definition of $A$, we obtain
\begin{align*}
    \Err(\set{S})^2 &= \norm{\int_{\set{X}} \kappa(\cdot,x) u(x) \, \d\mu(x)  - \int_{\set{X}} \kappa(\cdot,\set{S})\kappa(\set{S},\set{S})^{-1}\kappa(\set{S},x) u(x) \, \d\mu(x)}_\RKHS^2 \\
    &= \norm{(A - \hat{A}_\set{S}) u }_\RKHS^2.
\end{align*}
In the second line, we recognize the Nystr\"om approximate integral operator $\hat{A}_\set{S}$, defined above in \cref{eq:nystrom-integral-operator}.
In \cref{sec:infinite-setting}, we saw that $A^{1/2}$ is an isometry and a bijection between $\Ltwo(\mu)$ and $\RKHS$.
As such, it has an inverse, which we will denote $A^{-1/2}$, that is an isometry and bijection from $\RKHS$ to $\Ltwo(\mu)$.
This transformations acts according to the rule
\begin{equation*}
    A^{-1/2}f = \sum_{i=1}^\infty \lambda_i^{-1/2} e_i^{\vphantom{*}} (e_i^*f) \quad \text{for } f \in \RKHS.
\end{equation*}
Additionally, for $f \in \RKHS$ for which the sequence $\{\lambda_i^{-1} (e_i^*f)\}$ is square-summable we may also define a map $A^{-1}$ by the formula
\begin{equation*}
    A^{-1}f = \sum_{i=1}^\infty \lambda_i^{-1} e_i^{\vphantom{*}} (e_i^*f).
\end{equation*}
Using these linear maps, we may bound
\begin{align*}
    \Err(\set{S})^2 &= \norm{A^{-1/2}(A - \hat{A}_\set{S}) u }_{\Ltwo(\mu)}^2 \\
    &= \langle (A - \hat{A}_\set{S})A^{-1}(A - \hat{A}_\set{S})u, u\rangle_{\Ltwo(\mu)}  \le \langle (A - \hat{A}_\set{S})u, u\rangle_{\Ltwo(\mu)}.
\end{align*}
The last inequality is valid because $A - \hat{A}_{\set{S}}$ and $\hat{A}_{\set{S}}$ are both psd operators.
We have established the following property:

\begin{proposition}[Kernel quadrature error]
    The worst-case quadrature error $\Err(\set{S})$ defined in \cref{eq:ideally-weighted-worst-case} admits the bound
    \begin{equation*}
        \Err(\set{S})^2 \le \langle (A - \hat{A}_\set{S})u, u\rangle_{\Ltwo(\mu)}.
    \end{equation*}
\end{proposition}

\index{kernel quadrature!with randomly pivoted Cholesky|(}
Combining with \cref{thm:cont-rpcholesky-spectral} immediately yields a result for \RPCholesky kernel quadrature with ideal weights.

\begin{corollary}[Kernel quadrature with \RPCholesky]
    Instate the prevailing notation and assumptions, and let $r \ge 0$ and $\varepsilon > 0$.
    Then
    \begin{equation*}
        \expect \Err(\set{S})^2 \le \varepsilon \cdot \norm{u}_{\Ltwo(\mu)}^2 \cdot \sum_{r+1}^\infty \lambda_r
    \end{equation*}
    where $\set{S}$ is selected by running $k$ steps of \RPCholesky where
    \begin{equation*}
        k \ge r \log \left( \frac{2\lambda_1}{\varepsilon \sum_{i=r+1}^\infty \lambda_i} \right) + \frac{2}{\varepsilon}.
    \end{equation*}
\end{corollary}\index{kernel quadrature!error bounds|)}

\myprogram{Program to compute the ideal kernel quadrature weights for computing integrals $\int_\set{X} f(x) \overline{u(x)} \, \d\mu(x)$ of functions $f$ drawn from an RKHS.}{}{kernel_quad_wts}
\index{kernel quadrature!with randomly pivoted Cholesky|)}

\index{kernel quadrature!computational considerations|(}
\myparagraph{Computational considerations}
Implementing ideally weighted kernel quadrature requires evaluating the integrals
\begin{equation} \label{eq:kq-integrals}
    \int_{\set{X}} \kappa(s_i,x) u(x) \, \d\mu(x) \quad \text{for } i=1,2,\ldots,k
\end{equation}
that compose the right-hand side of the system \cref{eq:ideal-weights-system}.
Many schemes for quadrature on RKHSs have been proposed that avoid the computation of these exact integrals \cite{DM22,SDM22,HOL22}, but the rate of convergence is reduced.
There is no obvious way to obtain \emph{spectrally accurate} kernel quadrature schemes without high-accuracy evaluation of the integrals \cref{eq:kq-integrals}.\index{kernel quadrature!computational considerations|)}

\index{randomly pivoted Cholesky!numerical results|(}\index{kernel quadrature!with randomly pivoted Cholesky|(}
\myparagraph{Example}
As an illustration, we consider a simple problem in two dimensions.
The measure is $\mu = \Normal(\vec{0},\Id_2)$, the kernel is square-exponential with bandwidth $1$, the weight function is $u\equiv 1$, and the function is $f(\vec{x}) = \cos(\norm{\vec{x}})$.
The integral is known exactly: 
\begin{equation*}
    \int_{\set{X}} f(\vec{x}) \, \d\mu(\vec{x}) = 1 - \sqrt{2} \mathrm{F}(1/\sqrt{2}) = 0.27522154099292\ldots
\end{equation*}
Here, $\mathrm{F}$ is the Dawnson integral  (\texttt{DawsonF} in Mathematica).
The values of $Au$ are also known exactly:
\begin{equation*}
    Au(\vec{x}) = \frac{1}{2} \exp \left( - \frac{\norm{\vec{x}}^2}{4} \right).
\end{equation*}

\begin{figure}
    \centering
    \includegraphics[width=0.48\linewidth]{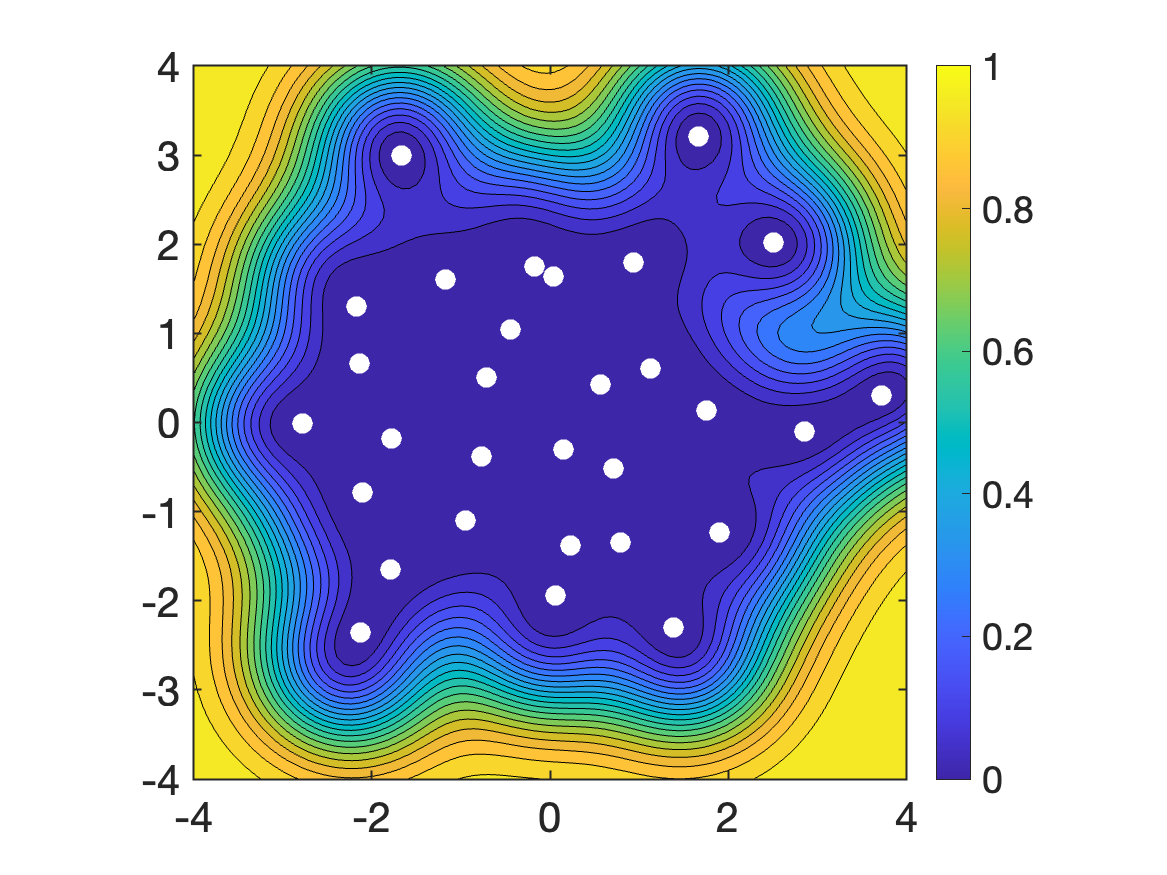}
    \includegraphics[width=0.48\linewidth]{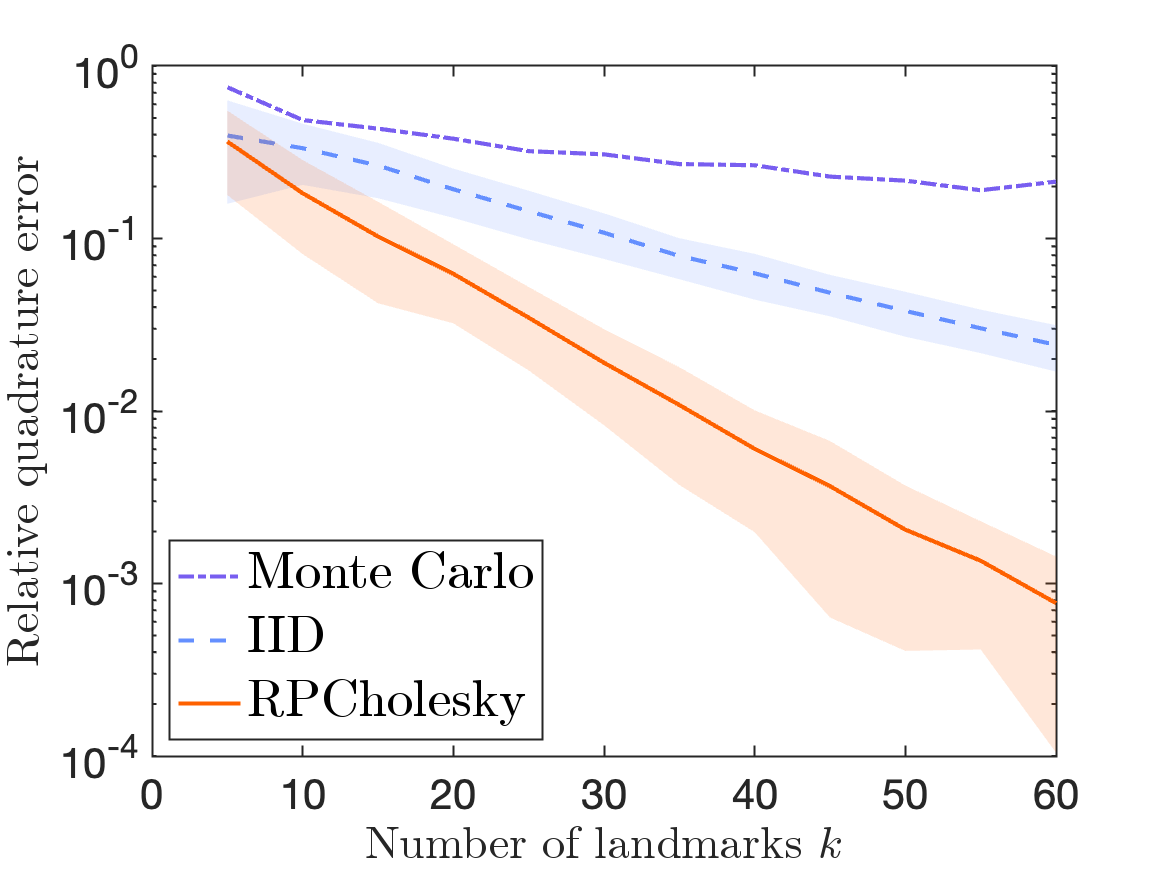}
    \caption[Kernel quadrature using nodes selected by continuous \RPCholesky and iid random selection]{\emph{Left:} Quadrature nodes selected by \RPCholesky ($k=30$, white). 
    Shading shows the diagonal elements of the residual kernel $\kappa_{\set{S}}(\vec{x},\vec{x})$.
    \emph{Right:} Error in evaluating integral $\int_{\set{X}} f(\vec{x}) \, \d\mu(\vec{x})$ using kernel quadrature for different numbers of nodes $k$, selected by \RPCholesky or by iid sampling.
    We also show the error of the simple Monte Carlo estimate $\int_{\set{X}} f(\vec{x}) \, \d\mu(\vec{x}) \approx k^{-1} \sum_{i=1}^k f(\vec{x}_i)$ for $\vec{x}_1,\ldots,\vec{x}_k \stackrel{\text{iid}}{\sim} \mu$.
    The function $f$, kernel $\kappa$, and measure $\mu$ are described in the text.
    Lines show a mean of 100 trials, and shaded regions show one standard deviation.}
    \label{fig:kernel-quadrature}
\end{figure}

\Cref{fig:kernel-quadrature} shows the results.
The left panel shows $k=30$ nodes selected by \RPCholesky, with the shading indicating the diagonal values $\kappa_{\set{S}}(\vec{x},\vec{x})$ of the residual kernel.
The right panel shows the relative error for the outputs of ideally weighted kernel quadrature with $5\le k \le 60$ nodes selected with \RPCholesky and by iid sampling.
We also compare the simple Monte Carlo integral estimate $\int_{\set{X}} f(\vec{x}) \, \d\mu(\vec{x}) \approx k^{-1} \sum_{i=1}^k f(\vec{x}_i)$ for $\vec{x}_1,\ldots,\vec{x}_k \stackrel{\text{iid}}{\sim} \mu$.
The rate of convergence for \RPCholesky kernel quadrature is faster than both other methods.

See \cite{EM23a} for more experiments with \RPCholesky quadrature, accuracy and timing comparisons of \RPCholesky quadrature to more methods for kernel quadrature \cite{BBC20,DM22,SDM22,HOL22}, and an application to chemistry.\index{chemistry}

\index{kernel quadrature|)}\index{kernel quadrature!with randomly pivoted Cholesky|)}\index{randomly pivoted Cholesky!numerical results|)}\index{infinite-dimensional randomly pivoted Cholesky|)}

\chapter{Blocked algorithms} \label{ch:blocked}

\epigraph{Since all machines have stores of finite size often divided up into high speed and auxiliary sections, storage considerations often have a vitally important part to play.}{Jim Wilkinson, \emph{The use of iterative methods for finding the latent roots and vectors of matrices} \cite{Wil55}}

\index{randomly pivoted Cholesky!blocked versions|(}\index{blocked algorithm!reasons for blocking|(}

When designing algorithms, we often operate under the convenient fiction that all memory elements can be accessed equally quickly and that the runtime of an algorithm is roughly proportional to the total number of arithmetic operations.
Unfortunately, the performance of algorithms in practice is more complicated than this simple model suggests.
In particular, matrix algorithms tend to be more performant when data is processed in \emph{blocks} rather than column-by-column or entry-by-entry.
This chapter responds to this reality by discussing improved variants of the \RPCholesky algorithm that are more hardware-efficient and utilize block matrix computations.\index{randomly pivoted Cholesky!blocked versions|)}

\myparagraph{Sources}
This chapter is based on the paper

\fullcite{ETW24a},

though some of the content appears in the original \RPCholesky paper

\fullcite{CETW25}.

It expands on these references by providing a more detailed discussion of the benefits of block matrix computations (\cref{sec:why-blocking}) and a more thorough discussion of the RBRP Cholesky method of \cite{DCMP23} (\cref{sec:rbrp-rpcholesky}).

\myparagraph{Notation}
Throughout this chapter, we will be interested in computing a pivoted partial Cholesky decomposition of a psd matrix $\mat{A} \in \field^{n\times n}$ of size $n$.
The pivots will be denoted $s_1,\ldots,s_k$, and $\mat{A}^{(i)}$ denotes the residual of the decomposition after the first $i$ elimination steps:
\begin{equation} \label{eq:block-chapter-residual}
    \mat{A}^{(i)} = \mat{A} - \mat{A}\langle \{s_1,\ldots,s_i\} \rangle.
\end{equation}
The matrix $\mat{F} \in \field^{n\times k}$ denotes the factor matrix computed by the pivoted partial Cholesky algorithm (\cref{prog:pivpartchol}), which satisfies $\mat{A}\langle \{s_1,\ldots,s_i\}\rangle = \mat{F}(:,1:i)\mat{F}(:,1:i)^*$ for every $i$.

\myparagraph{Outline}
\Cref{sec:why-blocking} discusses the benefits of blocking for matrix computations with a focus on kernel computations in particular.
The next three sections present three blocked versions of the \RPCholesky algorithm, referred to as \emph{block \RPCholesky} (\cref{sec:block-rpcholesky}), \emph{RBRP Cholesky} (\cref{sec:rbrp-rpcholesky}), and \emph{accelerated \RPCholesky} (\cref{sec:acc-rpcholesky}).
A comparison of these methods appears in \cref{sec:block-rpcholesky-compare}.

\index{memory transfer|(}
\section{Why blocking?} \label{sec:why-blocking}

Why are blocked algorithms faster than unblocked algorithms? 
The answer boils down to \emph{memory transfers}.
Applied mathematicians and computer scientists tend to analyze numerical algorithms by counting the number of arithmetic operations.
However, on modern computing architectures, movement of data between random access memory and the cache can become the rate-limiting factor in the speed of algorithms, exceeding the cost of floating point operations \cite[pp.~13, 581--582]{DFF+03}.

\subsection{A simplified model for kernel matrix computations}

To get a sense of why memory transfer costs can make blocked algorithms faster than unblocked methods, consider the following simplified scenario.
Suppose we wish to generate columns of a kernel matrix\index{kernel matrix} $\mat{A} = \kappa(\set{D},\set{D})$, where $\set{D}\subseteq \set{X}$ is a set of $n$ data points taking values in some set $\set{X}$.
We consider an extremely simplified computing machine with a memory hierarchy consisting of two layers; the data elements $x \in \set{D}$ are stored in \emph{memory}.
In order to perform computations on data elements, they must be moved from memory into the \emph{cache}.\index{cache}
However, the cache has a limit size and is only able to store $1\ll t\ll n$ elements of $\set{X}$ at once.

First, suppose we generate columns of $\mat{A}$ one at a time, as is done in the standard \RPCholesky algorithm.
To generate each column $\kappa(\set{D},x)$ requires bringing each entry of $x'\in\set{D}$ into memory one at a time and evaluating $\kappa(x',x)$. 
Since the cache has $t \ll n$ entries, most entries of $\set{D}$ must be brought back from memory to the cache\index{cache} for each column generation.
Thus, generating $k$ columns requires $\order(kn)$ memory transfers.

Now contrast this with a strategy where we generate columns in \emph{blocks} $\set{B}$ of size $|\set{B}| = b\le t-1$.
To evaluate a block of columns $\kappa(\set{D},\set{B})$, we can now generate the matrix row-wise $\kappa(x',\set{B})$, moving each $x'\in\set{D}$ from memory to the cache\index{cache} in sequence.
The consequence is we can generate \emph{all $b$} columns at once using just $n$ memory transfers.
Generating $k\ge b$ columns in blocks of size $b$ requires $\order(nk/b)$ memory transfers, a factor $b$ smaller than the one-column-at-a-time approach.
Thus, for an algorithm whose runtime is dominated by memory movement, the blocked algorithm can be significantly faster.

\begin{remark}[Cache and memory: Models and reality]
    The model we have used in this section is a simplified version of the \emph{ideal cache model}\index{ideal cache model} \cite{FLPR99}.
    The memory architecture for a modern computer is significantly more complicated than the basic description given here.
    In particular, memory hierarchies are divided into many more than two levels and data is moved between levels in chunks of fixed size called \emph{cache lines}.\index{cache}
    Given this complexity, we advise against using the model described above to make quantitative predictions about the runtime of algorithms (i.e., it is not typically true that generating columns in blocks of size $b$ is precisely $b$ times faster than generating one-at-a-time).
    Rather, our point should be taken as qualitative and conceptual.
\end{remark}

\subsection{The submatrix access model}

In view of the discussion in the previous subsection, we see the entry access model described in \cref{sec:entry-access} may be too simplistic an abstraction to describe the behavior of kernel matrix algorithms.
A better abstraction, introduced in \cite{ETW24a}, is the \emph{submatrix access model}:

\actionbox{\textbf{Submatrix access model.} We are given a matrix $\mat{B} \in \field^{m\times n}$ that may be efficiently accessed by submatrices $\mat{B}(\set{S},\set{T})$ of small size $|\set{S}| \cdot |\set{T}| \ll mn$.
The access is most efficient when both $\set{S}$ and $\set{T}$ have size $\gg 1$.}
The goal of this section will be to develop versions of \RPCholesky that are efficient in the submatrix access model.

\index{matrix--matrix operations|(}
\subsection{Matrix--matrix arithmetic}

Memory transfers also help explain why block matrix operations like matrix--matrix multiplication are fast on modern computers.
The core idea is that, with careful implementation to manage data movement, computing a matrix--matrix product $\mat{B}\mat{C}$ all at once requires fewer memory transfers than computing the product one column $\mat{B}\vec{c}_i$ at a time.


Achieving high performance for block matrix arithmetic on modern computer architectures demands significant programming effort.
Fortunately, much of this work has already been done, with hardware vendors providing optimized implementations of matrix multiplication as part of the Level-3 Basic Linear Algebra Subprograms (\BLASthree).
Modern linear algebra libraries such as \LAPACK are built to exploit these routines, computing matrix factorizations (\LU, Cholesky, \QR, etc.) through sequences of matrix multiplications implemented via \BLASthree.
For the algorithm designer or numerical programmer, the takeaway is clear: To maximize performance, it is advantageous to reorganize computations to rely on matrix--matrix operations rather than on matrix--vector or vector--vector operations.

\index{blocked algorithm!reasons for blocking|)}\index{memory transfer|)}\index{matrix--matrix operations|)}\index{randomly pivoted Cholesky!blocked versions|(}\index{blocked algorithm!blocked versions of randomly pivoted Cholesky|(}\index{block randomly pivoted Cholesky|(}
\section{Algorithm 1: Block \RPCholesky} \label{sec:block-rpcholesky}

Given currently selected pivots $s_1,\ldots,s_i$, the standard \RPCholesky algorithm draws a single next pivot at each iteration using the sampling rule
\begin{equation*}
    s_{i+1} \sim \diag(\mat{A}^{(i)}),
\end{equation*}
and it computes the next column $\vec{f}_{i+1}$ of the factor matrix $\mat{F}$ by the formula
\begin{equation} \label{eq:rpcholesky-new-col}
    \vec{f}_{i+1} \coloneqq \frac{\vec{a}_{s_{i+1}} - \mat{F}(:,1:i)\mat{F}(s_{i+1},1:i)^*}{(a_{s_{i+1}s_{i+1}} - \mat{F}(s_{i+1},1:i)\mat{F}(s_{i+1},1:i)^*)^{1/2}}.
\end{equation}
(Recall the definition of $\mat{A}^{(i)}$ in \cref{eq:block-chapter-residual}.)
Each iteration exposes a single column $\vec{a}_{s_i}$ of the matrix $\mat{A}$ and evaluates the matrix--vector product $\mat{F}(:,1:i)\mat{F}(s_{i+1},1:i)^*$.

To create a block variant of this algorithm, we can draw multiple random pivots at once
\begin{subequations} \label{eq:block-rpcholesky}
\begin{equation} \label{eq:block-rpcholesky-pivots}
    s_{i+1},\ldots,s_{i+b} \simiid \diag(\mat{A}^{(i)}),
\end{equation}
where $b \ge 1$ is a user-specified block size. 
Setting $\set{S}' \coloneqq \{s_{i+1},\ldots,s_{i+b}\}$, a block of $b$ new columns of the factor matrix can be generating using the formula
\begin{equation} \label{eq:block-rpcholesky-new-cols}
    \mat{F}(:,i+1:i+b) \coloneqq (\mat{A}(:,\set{S}') - \mat{F}(:,1:i)\mat{F}(\set{S}',1:i)^*)\mat{R}^{-1},
\end{equation}
where
\begin{equation} \label{eq:block-rpcholesky-cholesky}
    \mat{R}^*\mat{R} = \mat{A}(\set{S}',\set{S}') - \mat{F}(\set{S}',1:i)\mat{F}(\set{S}',1:i)^* = \mat{A}^{(i)}(\set{S}',\set{S}').
\end{equation}
\end{subequations}
is an upper-triangular Cholesky factorization.
Equation \cref{eq:block-rpcholesky-new-cols} is the block Cholesky analog of \cref{eq:rpcholesky-new-col}.
The update steps \cref{eq:block-rpcholesky} generate $b$ columns of the kernel matrix, and they employ block matrix--matrix operations to update the factor matrix $\mat{F}$.
We call the procedure \cref{eq:block-rpcholesky} the \emph{block \RPCholesky algorithm}.
Code is provided in \cref{prog:block_rpcholesky}.

\myprogram{Block \RPCholesky algorithm for psd low-rank approximation and column subset selection.}{Subroutine \texttt{sqrownorms} is defined in \cref{prog:sqrownorms}.}{block_rpcholesky}

\index{redundant pivot problem|(}\index{pivots!redundant|(}
\myparagraph{The redundant pivot problem}
The block \RPCholesky algorithm is simple and performant.
However, it can suffer from the \emph{redundant pivot problem}: The selected pivots $\set{S}' = \{s_{i+1},\ldots,s_{i+b}\}$ can contain similar information to each other, resulting in a submatrix $\mat{A}^{(i)}(\set{S}',\set{S}')$ that is ill-conditioned.
Redundant pivots can hurt the performance of the algorithm in two ways: 
\begin{enumerate}
    \item \textbf{Cholesky decomposition failure.} On some examples, the matrix $\mat{A}^{(i)}(\set{S}',\set{S}')$ can be (numerically) rank-deficient, causing the Cholesky decomposition \cref{eq:block-rpcholesky-cholesky} to fail.\index{Cholesky decomposition!failure when matrix is rank-deficient}
    \item \textbf{Lower quality approximation/pivot set.} For a fixed rank $k$, a low-rank approximation with redundant pivots is less accurate than an approximation where each pivot is distinct. 
    Additionally, for subset selection problems, the computational task is to select a \emph{diverse} set of representatives for a data set; selecting redundant pivots is contrary to this goal.
\end{enumerate}

The first issue can be addressed by adding a shift on the order of the machine precision to the matrix $\mat{A}^{(i)}(\set{S}',\set{S}')$ in \cref{eq:block-rpcholesky-cholesky} or applying \RPCholesky to a shifted version $\mat{A} + \mu \Id$ of the matrix.
My coauthors and I have experimented with different values for the shift parameter across different papers and different versions of our software, with the goal of finding the minimum possible shift that guarantees convergence of the algorithm.
Our current recommendation is to use a shift of $4 \max\{a_{ii}\} u$, where $u$ is the unit roundoff\index{unit roundoff} ($u\approx 10^{-16}$ in double precision).

The second issue may or may not be serious depending on the application.
We emphasize that redundant pivots do not \emph{hurt} the quality of the approximation (at least in exact arithmetic); by \cref{prop:nystrom-properties}\ref{item:nystrom-monotonicity}, the quality of a low-rank approximation $\mat{A}\langle \set{S} \cup \set{S}_{\mathrm{redundant}} \rangle$ with redundant pivots is never worse than the low-rank approximation $\mat{A}\langle \set{S} \rangle$ without redundant pivots.
Rather, each redundant pivot is a wasted opportunity to find a better pivot that better represents the data set.
The next two sections outline improved versions of block \RPCholesky that \emph{subsample} the set of iid pivots \cref{eq:block-rpcholesky-pivots} to remove redundant pivots.\index{redundant pivot problem|)}\index{pivots!redundant|)}

\index{block randomly pivoted Cholesky!theoretical analysis|(}
\myparagraph{Analysis}
The block \RPCholesky bound satisfies the following guarantee, which we take from \cite[Thm.~4.1]{ETW24a}.

\begin{theorem}[Block randomly pivoted Cholesky] \label{thm:block-rpcholesky}
    Let $r\ge 1$ be an integer, $\varepsilon > 0$ be a real number, and $\mat{A}\in\field^{n\times n}$ be a psd matrix.
    Introduce the relative error of the best rank-$r$ approximation:
    \begin{equation*}
        \eta \coloneqq \tr(\mat{A} - \lowrank{\mat{A}}_r) / \tr(\mat{A}).
    \end{equation*}
    Fix a block size $b$ that divides the rank $k$.
    Block randomly pivoted Cholesky produces an $(r,\varepsilon)$-approximation (\cref{eq:r-eps-appox}) provided
    \begin{equation} \label{eq:block-rpcholesky-trace-k}
        k \ge \frac{r}{\varepsilon} + (r+b) \log \left( \frac{1}{\varepsilon \eta} \right).
    \end{equation}
\end{theorem}

The guarantee \cref{eq:block-rpcholesky-trace-k} for block \RPCholesky is very similar to the guarantee \cref{eq:rpcholesky-trace-k} for standard \RPCholesky; the only difference between the bounds is that the $(r+b)$ factor in the block \RPCholesky bound is replaced by $r$ in the standard \RPCholesky bound.
The implication is that the \warn{worst-case guarantees} for standard and block \RPCholesky are similar for any block size $1 \le b \le \order(r)$.
However, while the theoretical bounds for standard and block \RPCholesky are different, the performance of these methods differs greatly in practice.
On some instances, standard \RPCholesky produce far smaller errors than block \RPCholesky.
(On the other hand, block \RPCholesky can be much faster than standard \RPCholesky, so it can accommodate a larger value of $k$ given a fixed runtime budget.)

We refer the interested reader to \cite[Thm.~4.1]{ETW24a} for a proof of this result, which proceeds by a comparison with the accelerated \RPCholesky algorithm (\cref{sec:acc-rpcholesky}).
See also the discussion around \cite[Thm.~4.3]{ETW24a} for a discussion of the error bounds for block \RPCholesky that can be inferred from the earlier works of Deshpande, Rademacher, Vempala, and Wang \cite{DRVW06}.\index{block randomly pivoted Cholesky!theoretical analysis|)}\index{block randomly pivoted Cholesky|)}

\index{robust blockwise random pivoting!for positive-semidefinite low-rank approximation|(}
\section{Algorithm 2: RBRP Cholesky} \label{sec:rbrp-rpcholesky}

The first approach to improving on block \RPCholesky is \emph{robust blockwise random pivoting} (RBRP), which was proposed by Dong, Chen, Martinsson, and Pearce \cite{DCMP23}.
This approach was developed concurrently with and released prior to the accelerated \RPCholesky algorithm, which we describe in the next section.
The RBRP strategy was originally developed for use with the randomly pivoted \QR\index{randomly pivoted QR@randomly pivoted \QR} algorithm (\cref{ch:low-rank-general}), but the idea works equally well when combined with \RPCholesky.
We describe the \RPCholesky variant here.

The RBRP algorithm is delightfully simple.
Just like ordinary block \RPCholesky, we begin by drawing a set of iid random pivots
\begin{equation*}
    \set{S}' = \{s_{i+1}',\ldots,s_{i+b}'\} \quad \text{with}\quad s_1',\ldots,s_b' \simiid \diag(\mat{A}^{(i)}).
\end{equation*}
Then, form the submatrix 
\begin{equation*}
    \mat{H} \coloneqq \mat{A}^{(i)}(\set{S}',\set{S}') = \mat{A}(\set{S}',\set{S}') - \mat{F}(\set{S}',1:i)\mat{F}(\set{S}',1:i)^*.
\end{equation*}
To filter any redundant pivots, we apply a step Dong et al.\ call \emph{robust blockwise filtering},\index{robust blockwise filtering} which runs Cholesky with greedy pivoting\index{greedy pivoted Cholesky} (\cref{prog:greedy_chol}) on $\mat{H}$ until the trace has been reduced by a factor $\tau \in (0,1)$.
(Dong et al.\ suggest $\tau = 1/b$.)
Letting $\set{T} \subseteq \set{S}'$ denote the set of pivots selected by the greedy method, we induct the subselected pivots as new pivots:
\begin{equation*}
    s_{i+1} = s_{i_1}', s_{i+2} = s_{i_2}',\ldots,s_{i+\ell} \coloneqq s_{i_{\ell}'} \quad \text{where }\set{T} = \{s_{i_1}',\ldots,s_{i_{\ell}}'\}.
\end{equation*}
Observe that this procedure inducts a random number of pivots $1 \le \ell \coloneqq |\set{T}| \le b$ at each step.
Code for the robust blockwise filtering and RBRP Cholesky procedures appear in \cref{prog:robust_block_filter,prog:rbrp_chol}.
At present, theoretical analysis for RBRP Cholesky is unavailable.

\myprogram{Implementation of robust blockwise filtering, a subroutine for the RBRP Cholesky algorithm (\cref{prog:rbrp_chol}).}{}{robust_block_filter}

\myprogram{RBRP Cholesky algorithm for psd low-rank approximation and column subset selection.}{Subroutines \texttt{robust\_block\_filter} and \texttt{sqrownorms} ared defined in \cref{prog:robust_block_filter,prog:sqrownorms}.}{rbrp_chol}\index{robust blockwise random pivoting!for positive-semidefinite low-rank approximation|)}

\index{accelerated randomly pivoted Cholesky|(}\index{rejection sampling!for implementing randomly pivoted Cholesky|(}
\section{Algorithm 3: Accelerated \RPCholesky} \label{sec:acc-rpcholesky}

\emph{Accelerated \RPCholesky} is another blocked variant of \RPCholesky that solves the redundant pivot problem\index{redundant pivot problem} \cite{ETW24a}.
It is conceptually similar to the RBRP Cholesky algorithm but uses rejection sampling in place of robust blockwise filtering.\index{robust blockwise filtering}
Compared to block \RPCholesky and RBRP Cholesky, accelerated \RPCholesky has the virtue that it produces \warn{the same random output as standard \textsc{RPCholesky}} (\cref{prog:rpcholesky}).
Thus, it inherits \RPCholesky's theoretical guarantees (\cref{sec:rpc-analysis}), and it can be used for applications where the precise distribution of the \RPCholesky pivots is important, such as fixed-size DPP sampling (\cref{sec:dpp-connections}).\index{determinantal point process sampling!connection to randomly pivoted Cholesky}\index{determinantal point process sampling!algorithms}
We will discuss the relative strengths of the three algorithms in \cref{sec:block-rpcholesky-compare,sec:block-rpcholesky-compare-discussion}.

Recall from \cref{sec:rpc-rejection} that rejection sampling allows one to sample from a (potentially complicated) \emph{target distribution}\index{target distribution} using samples from a simpler \emph{proposal distribution}.\index{proposal distribution}
In that section, we used the diagonal of a kernel function $\kappa$ as a proposal distribution, with the \RPCholesky sampling distribution as the target.
Here, we will use the same principles in a somewhat different way.

Suppose we have sampled an initial batch of pivots $s_1,\ldots,s_i$, giving rise to the residual matrix $\mat{A}^{(i)}$.
We are interested in sampling a new set of pivots $s_{i+1},\ldots,s_{i+\ell}$, each of which has distribution
\begin{equation*}
    s_{i+j+1} \sim \diag(\mat{A}^{(i+j)}) \quad \text{for } j = 0,\ldots,\ell-1.
\end{equation*}
We shall sample these distributions using rejection sampling, with $\diag(\mat{A}^{(i)})$ serving as the sampling distribution.

\myprogram{Accelerated \RPCholesky method for psd low-rank approximation and column subset selection.}{Subroutine \texttt{rejection\_sample\_submatrix} is defined in \cref{prog:rejection_sample_submatrix}.}{acc_rpcholesky}

More precisely, we do the following.
As in the other algorithms, we draw a block of proposal pivots $\set{S}' = \{s_1',\ldots,s_b'\}$ iid from the diagonal of the residual matrix:
\begin{equation*}
    s_1',\ldots,s_b'\simiid \diag(\mat{A}^{(i)}).
\end{equation*}
We emphasize that the sampling must be performed \warn{with replacement}.
(The other algorithms can be implemented with or without replacement.)
We now use this block of pivots to perform $b$ proposals of rejection sampling.
Beginning with $\ell = 0$, we perform the following step for each $j=1,\ldots,b$:
\begin{equation*}
    \text{With probability } \frac{\mat{A}^{(i+\ell)}(s_j',s_j')}{\mat{A}^{(i)}(s_j',s_j')}, \text{ \emph{accept} by setting } \ell\gets \ell+1 \text{ and }s_{i+\ell}\gets s_j'.
\end{equation*}
The output of this rejection sampling loop is a collection of new pivots $s_{i+1},\ldots,s_{i+\ell}$ of random size $1\le \ell\le b$.
(We have $\ell\ge 1$ since the first pivot is always accepted.)
Having selected new pivots, we then generate columns $\mat{F}(:,i+1:i+\ell)$ using the block update formula \cref{eq:block-rpcholesky-new-cols}.
In the next round, we use the diagonal of $\mat{A}^{(i+\ell)}$ as the new proposal distribution for selecting new pivots.
The accelerated \RPCholesky method is shown in \cref{prog:acc_rpcholesky}.

The accelerated \RPCholesky method uses rejection sampling to simulate the performance of the original \RPCholesky method while taking advantage of blockwise matrix operations.
We emphasize that the output of accelerated \RPCholesky has the \warn{same random distribution} as the output of the original \RPCholesky algorithm.

\myparagraph{\RejectionSampleSubmatrix subroutine}
For computational efficiency and extensibility, the rejection sampling steps in the accelerated \RPCholesky program are encapsulated in a subroutine called \RejectionSampleSubmatrix, shown in \cref{prog:rejection_sample_submatrix}.
This subroutine takes as input a submatrix $\mat{H} = \mat{A}^{(i)}(\set{S}',\set{S}') \in \field^{\set{S}'\times \set{S}'}$ and proposal distribution subvector $\vec{u} \ge \diag(\mat{H}) \in \real_+^{\set{S}'}$.
The vector $\vec{u}$ stores the entries of the proposal distribution for the proposed pivots $\set{S}'$.
In accelerated \RPCholesky, the proposal weights are $\diag(\mat{A}^{(i)}(\set{S}',\set{S}'))$.
One can also \RejectionSampleSubmatrix to develop a block version of the infinite-dimensional \RPCholesky sampler (\cref{prog:rejection_rpcholesky}), in which case the proposal distribution subvector is $\vec{u} = \diag(\kappa(\set{S}',\set{S}'))$.

This subroutine works by performing a Cholesky decomposition of $\mat{H}$ in place, with a probabilistic decision made at each step whether to eliminate entry $j$ or ignore it entirely.
For $j=1,\ldots,b$, do the following:
\begin{enumerate}
    \item \textbf{Accept/reject.} With probability $\mat{H}(s_j',s_j')/\vec{u}(s_j')$ accept and induct $s_j'$ as a new pivot $\set{T} \gets \set{T} \cup \{s_j'\}$ and go to step 2.
    Otherwise, reject and skip step 2.
    \item \textbf{Update and eliminate.} Perform a step of Cholesky decomposition to eliminate $s_j'$.
    Specifically, set
    \begin{equation*}
        \mat{H}(\{s_j',\ldots,s_b'\},s_j') \gets \mat{H}(\{s_j',\ldots,s_b'\},s_j') / \mat{H}(s_j',s_j')^{1/2},
    \end{equation*}
    then update 
    \begin{equation*}
        \mat{H}(\set{S}'_{j+1:b},\set{S}'_{j+1:b}) \gets \mat{H}(\set{S}'_{j+1:b},\set{S}'_{j+1:b}) - \mat{H}(\set{S}'_{j+1:b},s_j)\mat{H}(\set{S}'_{j+1:b},s_j)^*
    \end{equation*}
    with $\set{S}'_{j+1:b} \coloneqq \{s_{j+1}',\ldots,s_b'\}$.
\end{enumerate}
This procedure terminates with $\set{T}$ storing the set of accepted pivots and the lower triangular portion $\mat{L}$ of $\mat{H}(\set{T},\set{T})$ containing a Cholesky factorization of the $\set{T}$-submatrix of the (unmodified) input matrix $\mat{H}(\set{T},\set{T})$.
Code is provided in \cref{prog:rejection_sample_submatrix}.

\myprogram{Sample a set of \RPCholesky pivots using rejection sampling.}{\iffull\else This subroutine is used in \cref{prog:acc_rpcholesky}.\fi}{rejection_sample_submatrix}

\index{accelerated randomly pivoted Cholesky!comparison to \RPCholeskyRejection|(}\index{RejectionRPCholesky@\RPCholesky!comparision to \RPCholeskyRejection|(}
\myparagraph{Comparison of \RPCholeskyRejection and accelerated \RPCholesky}
We now have two different strategies for implementing \RPCholesky using rejection sampling, \RPCholeskyRejection and accelerated \RPCholesky.
A comparison is warranted.

We originally introduced \RPCholeskyRejection for sampling from an infinite space $\set{X}$ endowed with a positive-definite kernel function $\kappa$, but the algorithm also works for implementing \RPCholesky on a finite psd matrix $\mat{A}$.
\RPCholeskyRejection uses the same proposal distribution $\diag(\mat{A})$ throughout the algorithm, and it forms only the pivot set $\set{S}$, not the factor matrix $\mat{F}$.
Assuming $\diag(\mat{A}) = \onevec$ is the vector of ones, the \RejectionSampleSubmatrix runs in roughly
\begin{equation*}
    \order(k^3/\eta_k) \text{ operations},
\end{equation*}
where $\eta_k = \tr(\mat{A} - \lowrank{\mat{A}}_{k-1})/\tr(\mat{A})$ is the relative error of the best rank-$(k-1)$ approximation.
(Here, we assume that a random integer between $1\le i \le n$ can be generated in $\order(1)$ operations, which is true in the WordRAM model of computation.)
The runtime of the algorithm is \emph{independent} of the dimension $n$, but it suffers from the \emph{curse of smoothness},\index{curse of smoothness} with the algorithm slowing when the eigenvalues of $\mat{A}$ decay rapidly.
The reason for the curse of smoothness\index{curse of smoothness} is that the proposal distribution\index{proposal distribution} $\diag(\mat{A})$ remains static throughout the algorithm, and the probability of each proposal being accepted becomes smaller and smaller as the diagonal $\diag(\mat{A}^{(i)})$ of the residual shrinks.
An additional weakness of \RejectionSampleSubmatrix is that it does not generate the factor matrix $\mat{F}$, although it can be assembled in $\order(k^2n)$ operations if desired.

The accelerated \RPCholesky defeats the curse of smoothness\index{curse of smoothness} by updating the proposal distribution\index{proposal distribution} after every $b$ proposals. 
The cost of the algorithm is at most
\begin{equation*}
    \order(k^2n + kb^2) \text{ operations and } (k+1)n + kb^2 \text{ entry accesses},
\end{equation*}
only slightly more than standard \RPCholesky provided the block size $b$ is not large.
Moreover, accelerated \RPCholesky is designed so that most of the operations are performed in a block-wise fashion, making it 5$\times$ to 40$\times$ faster than standard \RPCholesky in practice.

Ultimately, \RejectionSampleSubmatrix (or a block version of it) is a useful algorithm for certain scenarios, namely for applying \RPCholesky to matrices with slow spectral decay or for extremely large or infinite-dimensional problems.
Accelerated \RPCholesky, by contrast, is excellent for \emph{general-purpose use}, and I recommend the algorithm for deployment in practice.\index{accelerated randomly pivoted Cholesky!comparison to \RPCholeskyRejection|)}\index{RejectionRPCholesky@\RPCholesky!comparision to \RPCholeskyRejection|)}

\index{accelerated randomly pivoted Cholesky!choice of block size|(}
\myparagraph{Picking the block size}
An advantage of the accelerated \RPCholesky algorithm over the block \RPCholesky algorithm is that it can accommodate a very large block size while maintaining approximation quality, as the rejection sampling step helps filter any redundant pivots.
When the number of pivots $k$ is set by the user in advance, I recommend $b=k/10$ or even $b=k/3$ as a good default value for the block size.

One can also set the block size automatically to avoid additional algorithmic parameters.
Here is one procedure I have used in my code.
Initialize the block size $b$ at some default value.
Each round of the algorithm, time both the rejection step (i.e., sampling proposals $\set{S}'$, forming the submatrix $\mat{A}^{(i)}(\set{S}',\set{S}')$, and running \RejectionSampleSubmatrix) and the processing step (reading columns $\mat{A}(:,\set{T})$ and updating $\mat{F}$).
Let $t_{\mathrm{r}}$ and $t_{\mathrm{p}}$ be the respective runtimes.
We may assume for simplicity that the runtimes of these steps satisfy the proportionality relations
\begin{equation} \label{eq:time-proportionality}
    t_{\mathrm{r}}\approx c_{\mathrm{r}}b^2 \quad \text{and} \quad t_{\mathrm{p}} \approx c_{\mathrm{p}}b.
\end{equation}
for appropriate constants $c_{\mathrm{r}}$ and $c_{\mathrm{p}}$.
It is reasonable want the rejection step to comprise a modest fraction of the total runtime, say, 20\%.
Setting $t_{\mathrm{r}} = 0.25t_{\mathrm{p}}$ and solving \cref{eq:time-proportionality} suggests choosing the next block size as $b = t_{\mathrm{p}}/(4t_{\mathrm{r}})$.
To make this procedure robust, we make sure that the block size does not change too much from iteration to iteration, and weset a maximum block size $b_{\mathrm{max}}$ (e.g., $b_{\mathrm{max}} = k/3$ if $k$ is known in advance).
The resulting update rule is
\begin{equation*}
    b \gets \max\left\{ \min\left\{\frac{t_{\mathrm{p}}}{4t_{\mathrm{r}}} \cdot b,\left\lceil 1.5b\right\rceil,b_{\mathrm{max}}\right\}, \left\lceil\frac{b}{3}\right\rceil \right\}.
\end{equation*}\index{accelerated randomly pivoted Cholesky!choice of block size|)}

\index{accelerated randomly pivoted Cholesky!theoretical analysis|(}
\myparagraph{Analysis}
One can show that the accelerated \RPCholesky algorithm achieves the $(r,\varepsilon)$-guarantee using a limited number of \warn{proposals}.
Specifically, we have the following result \cite[Thm.~4.1]{ETW24a}:

\begin{theorem}[Accelerated randomly pivoted Cholesky] \label{thm:acc-rpcholesky}
    Let $r\ge 1$ be an integer, $\varepsilon > 0$ be a real number, and $\mat{A} \in \field^{n\times n}$ be a psd matrix.
    Introduce the relative error of the best rank-$r$ approximation:
    \begin{equation*}
        \eta \coloneqq \tr(\mat{A} - \lowrank{\mat{A}}_r) / \tr(\mat{A}).
    \end{equation*}
    Fix a number of rounds $t$ and a number of proposals per round $b$.
    Accelerated randomly pivoted Cholesky produces an $(r,\varepsilon)$-approximation to $\mat{A}$ provided that
    \begin{equation} \label{eq:acc-rpcholesky-trace-k}
        bt \ge \frac{r}{\varepsilon} + (r+b) \log \left( \frac{1}{\varepsilon \eta} \right).
    \end{equation}
\end{theorem}

We omit the proof, which is a more involved version of the standard \RPCholesky analysis (\cref{thm:rpcholesky-trace}).
The similarity between the block \RPCholesky guarantee \cref{eq:block-rpcholesky-trace-k} and the accelerated \RPCholesky guarantee \cref{eq:acc-rpcholesky-trace-k} is no accident. 
\Cref{thm:block-rpcholesky} follows as a corollary of \cref{thm:acc-rpcholesky} under the principle that ``Accepting more pivots can only help the approximation quality.''\index{accelerated randomly pivoted Cholesky|)}\index{rejection sampling!for implementing randomly pivoted Cholesky|)}\index{accelerated randomly pivoted Cholesky!theoretical analysis|)}

\index{accelerated randomly pivoted Cholesky!numerical results|(}\index{randomly pivoted Cholesky!numerical results|(}\index{block randomly pivoted Cholesky!numerical results|(}
\section{Experiments} \label{sec:block-rpcholesky-compare}

\begin{figure}
    \centering
    \includegraphics[width=0.49\linewidth]{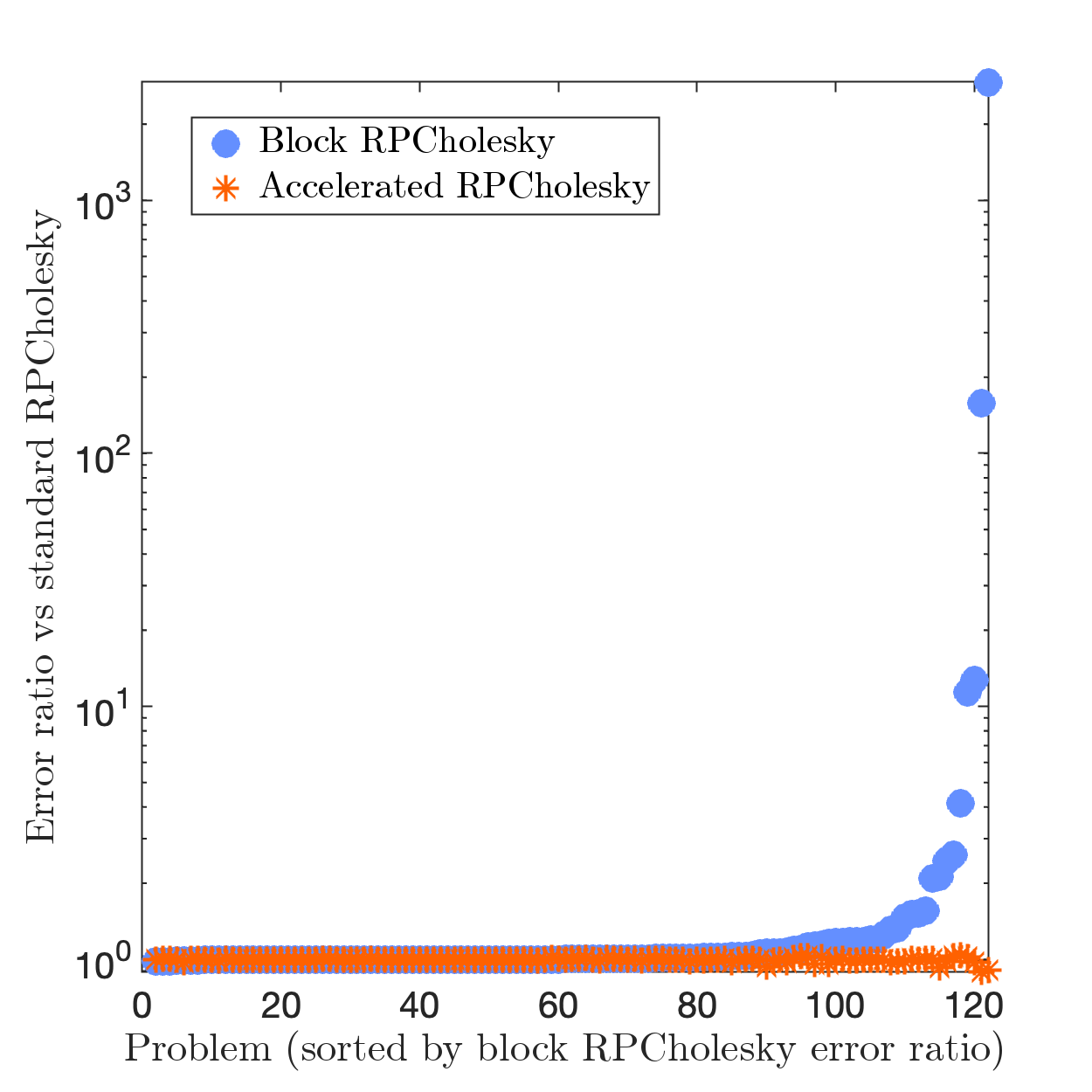}
    \includegraphics[width=0.49\linewidth]{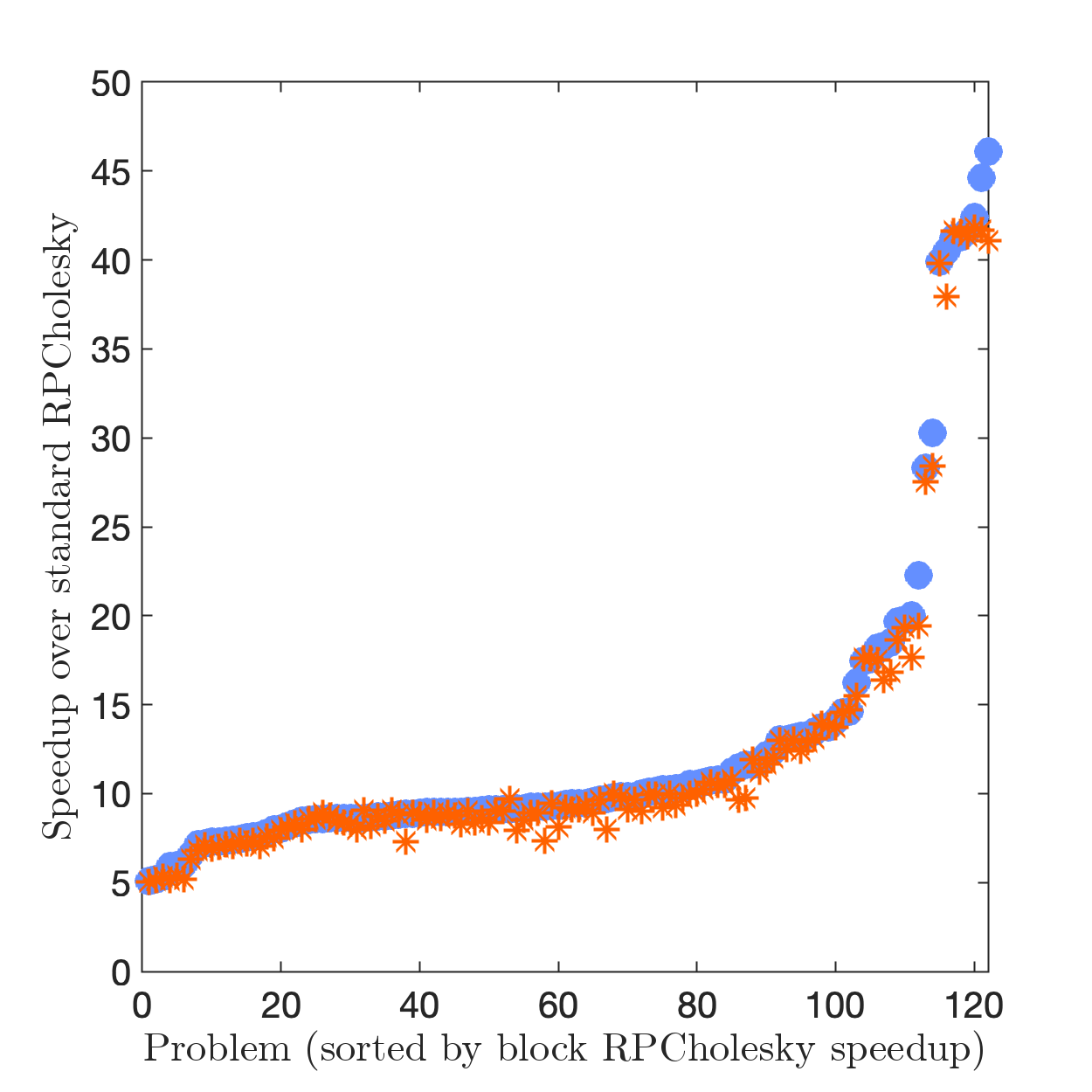}
    \caption[Comparison of speed and accuracy of blocked \RPCholesky algorithms]{Error ratio (\cref{eq:error-ratio-block}, \emph{left}) and speedup over standard \RPCholesky (\emph{right}) for block \RPCholesky and accelerated \RPCholesky on a testbed of 125 psd kernel matrices.
    The data is taken from \cite[Fig.~3]    {ETW24a}.}
    \label{fig:block-vs-acc-rpc}
\end{figure}

In \cite[\S3.1]{ETW24a}, my coauthors and I compared the error and runtime of block \RPCholesky and accelerated \RPCholesky for rank-$1000$ approximation of 125 test matrices with sizes between $4\times 10^4$ and $10^5$, including both synthetic and real data.
I reproduce the data from those experiments here in \cref{fig:block-vs-acc-rpc}.
Standard \RPCholesky provides a baseline for both speed and accuracy, and \cref{fig:block-vs-acc-rpc} plots the speedup and error ratio 
\begin{equation} \label{eq:error-ratio-block}
    \text{error ratio} \coloneqq \frac{\text{algorithm error}}{\text{standard \RPCholesky error}}
\end{equation}
for each matrix in the test suite.
The runtimes in this data were computed using the Python code employed in \cite{CETW25,ETW24a}, not the MATLAB code provided in this thesis.
The block size was set to $b= 150$.

We see that accelerated and block \RPCholesky achieve comparable speed to each other, with both methods achieving speedups of $5\times$ to $40\times$ over standard \RPCholesky.
These experiments demonstrate the large benefits of blocking for the speed of \RPCholesky-type low-rank approximation algorithms.

These results show that the accuracy of block \RPCholesky and accelerated \RPCholesky differs significantly.
On all test matrices, the accelerated \RPCholesky method achieves the same error as standard \RPCholesky, which makes sense because both methods produce the same random distribution of outputs.
On 100 out of 125 examples, the error of block \RPCholesky is similar to the other methods, with an error ratio between 1.00 and 1.16.
Because of the redundant pivot problem, block \RPCholesky suffers significantly higher errors (up to $3000\times$) than standard \RPCholesky.

\begin{table}[t]
    \centering
    \begin{tabular}{ccc}
        \toprule
         & Relative trace error & Runtime (sec) \\
        \midrule
        Block \RPCholesky & 3.89e-04 $\pm$ 1.40e-04 & 7.00 $\pm$ 0.63  \\
        Accelerated \RPCholesky	& 4.85e-07 $\pm$ 3.60e-08	& 7.43 $\pm$ 1.06 \\
        RBRP Cholesky &	6.62e-07 $\pm$ 1.35e-07 &	7.53 $\pm$ 0.80 \\
        \bottomrule
    \end{tabular}
    \caption[Runtime and relative error of three blocked \RPCholesky algorithms on a difficult synthetic test matrix]{Runtime and relative error for rank-$1000$ approximation produced by three blocked \RPCholesky algorithms for the synthetic test matrix of dimension $n=10^5$ from \cite[Fig.~1]{ETW24a}. 
    (All methods use block size $b=120$.)
    The table shows the mean and standard deviation computed over $100$ trials, and data is taken from \cite[Tab.~3]{ETW24a}.}
    \label{tab:rbrp}
\end{table}

\index{robust blockwise randomly pivoting!for positive-semidefinite low-rank approximation|(}
Our paper \cite{ETW24a} does not provide detailed experiments comparing RBRP Cholesky to the other blocked \RPCholesky algorithms, but it does compare the algorithms on a single challenging synthetic  matrix.
We reproduce this table as \cref{tab:rbrp}.
We see that the performance of RBRP Cholesky and accelerated \RPCholesky are similar, with accelerated \RPCholesky being slightly more accurate and faster.
(Notably, the output of accelerated \RPCholesky is also \emph{less variable} than RBRP Cholesky, which is curious because the accelerated \RPCholesky algorithm uses ``more randomness'' than the RBRP Cholesky does.)
The small advantage of accelerated \RPCholesky over RBRP Cholesky on this example pales in comparison to the large improvement in accuracy both methods have over block \RPCholesky.\index{accelerated randomly pivoted Cholesky!numerical results|)}\index{randomly pivoted Cholesky!numerical results|)}\index{block randomly pivoted Cholesky!numerical results|)}

\index{accelerated randomly pivoted Cholesky|(}\index{block randomly pivoted Cholesky|(}
\section{Comparison of three algorithms} \label{sec:block-rpcholesky-compare-discussion}

All three blocked \RPCholesky algorithms have merits for some use cases.

The block \RPCholesky is simple and fast, and it performs operations on chunks of data of \emph{fixed} size $n\times b$.
As such, block \RPCholesky algorithm is ideal for computing hardware like GPUs that are optimized for repetitive operations on data buffers of fixed size.
Additionally, the block \RPCholesky algorithm may be appropriate for applications where its subpar accuracy can be compensated for by taking a larger approximation rank $k$.

Both accelerated \RPCholesky and RBRP Cholesky algorithms both yield much better low-rank approximations on tough examples than block \RPCholesky.
At least for the CPU experiments reported in \cref{sec:block-rpcholesky-compare}, the cost of using these more-accurate methods is minimal, as these methods achieve nearly the same speedup as block \RPCholesky does.
For this reason, I would recommend accelerated either \RPCholesky or RBRP Cholesky as a natural choice for most applications of \RPCholesky.
The merits of these methods over block \RPCholesky are particularly pronounced when used for subset selection\index{subset selection} problems, where it is important not to select redundant pivots.\index{redundant pivot problem}\index{pivots!redundant}

The differences between RBRP Cholesky and accelerated \RPCholesky are fairly minor, but there are some reasons to prefer accelerated \RPCholesky:
\begin{enumerate}
    \item Accelerated \RPCholesky produces the same random pivot distribution as standard \RPCholesky.
    As such, it inherits \RPCholesky's theoretical guarantees (\cref{sec:rpc-analysis}), and it can be used in applications where the precise distribution of pivots is important (e.g., sampling from a projection DPP, \cref{sec:dpp-connections}).
    \item In my computational experience, accelerated \RPCholesky sometimes produces modestly better low-rank approximations on some examples (like the one in \cref{tab:rbrp}).
\end{enumerate}
Notwithstanding these differences, accelerated \RPCholesky and RBRP Cholesky tend to perform similarly in practice, and both are natural choices for deployment in software.

\index{randomly pivoted Cholesky!blocked versions|)}\index{blocked algorithm!blocked versions of randomly pivoted Cholesky|)}\index{robust blockwise randomly pivoting!for positive-semidefinite low-rank approximation|)}\index{accelerated randomly pivoted Cholesky|)}\index{block randomly pivoted Cholesky|)}

\chapter{\texorpdfstring{Randomly pivoted \QR: Low-rank approximation of general matrices}{Randomly pivoted QR: Low-rank approximation of general matrices}} \label{ch:low-rank-general}

\epigraph{But the idea of using actual columns and rows of the matrix $\mat{A}$ [for low-rank approximation] can be highly attractive.
Those vectors have meaning.
They are often sparse and/or nonnegative and they
reflect useful properties that we wish to preserve in approximating $\mat{A}$.}{Gilbert Strang and Cleve Moler, \textit{\textsf{LU} and \textsf{CR} Elimination} \cite[\S6]{SM22}}

So far, this thesis has focused on low-rank approximation of psd matrices.
Algorithms for this task can achieve remarkable results, producing a near-optimal low-rank approximation to \emph{any} psd matrix after reading a fraction of its entries.
The algorithms we have seen accomplish this goal by forming a column Nystr\"om approximation to a matrix $\mat{A}$ using a judiciously chosen subset of columns.
Among the available algorithms for psd low-rank approximation from limited entry evaluations, \RPCholesky and its variants are among the fastest and most reliable.

The impossibility result \cref{prop:general-lra-impossibility} shows that producing accurate approximations to a general matrix is impossible from a small budget of entry evaluations without additional information.\index{low-rank approximation!impossibility results}
Therefore, to produce a low-rank approximation to a general matrix, we are generally interested in algorithms that read the full input matrix, opening up a much larger design space for algorithms.
In particular, natural algorithms include the randomized SVD\index{randomized SVD} and its variants (\cref{sec:rsvd,sec:rsi}), which multiply the input matrix $\mat{B}$ with a random test matrix $\mat{\Omega}$.

However, in some applications, we require low-rank approximations to $\mat{B}$ that are spanned by a subset of its columns.
This chapter will review \emph{randomly pivoted \QR}\index{randomly pivoted QR@randomly pivoted \QR} algorithms for constructing such low-rank approximations and describe their connection to \RPCholesky and other methods for psd low-rank approximation.

\myparagraph{Sources}
The algorithm we call randomly pivoted \QR algorithm was originally proposed as adaptive sampling by Deshpande, Rademacher, Vempala, and Wang \cite{DRVW06,DV06}.
Our treatment follows the papers

\fullcite{CETW25}

and

\fullcite{ETW24a}.

These papers revisited the randomly pivoted \QR\index{randomly pivoted QR@randomly pivoted \QR} algorithm, though their focus was more on \RPCholesky.
This chapter's treatment of randomly pivoted \QR\index{randomly pivoted QR@randomly pivoted \QR} is substantially expanded from those papers.
It contains new research including discussion of numerical stability issues.

\myparagraph{Outline}
\Cref{sec:column-projection} introduces column projection approximations for approximating a matrix using a subset of their columns and discusses pivoted partial \QR\index{randomly pivoted QR@randomly pivoted \QR} algorithms for computing them.
\Cref{sec:rpqr} presents randomly pivoted \QR\index{randomly pivoted QR@randomly pivoted \QR}, and \cref{sec:acc-rpqr} discuss fast implementations using rejection sampling.
Related work on sketchy pivoting is discussed in \cref{sec:sketchy-pivoting}, and numerical experiments are provided in \cref{sec:rpqr-experiments}.
Finally, \cref{sec:rpcholesky-rpqr-history} concludes by discussing the history and connections between randomly pivoted \QR\index{randomly pivoted QR@randomly pivoted \QR} and randomly pivoted Cholesky\index{randomly pivoted Cholesky}.

\section{Low-rank approximation via column selection} \label{sec:column-projection}

Recall from \cref{ch:lra} that many of the commonly used low-rank approximations for general matrices are \emph{projection approximations}.
Given a matrix $\mat{B} \in \field^{m\times n}$ and a test matrix $\mat{\Omega}$, the projection approximation to $\mat{B}$ is 
\begin{equation*}
    \Bhat \coloneqq \mat{\Pi}_{\mat{B}\mat{\Omega}} \mat{B} = \mat{Q}\mat{Q}^*\mat{B} \quad \text{where } \mat{Q} = \orth(\mat{B}\mat{\Omega}).
\end{equation*}
Throughout this chapter, $\mat{\Pi}_{\mat{F}}$ denotes the orthoprojector onto $\range(\mat{F})$.

\index{column projection approximation!definition|(}
Just as we obtained column \emph{Nystr\"om} approximations from choosing the test matrix $\mat{\Omega} = \Id(:,\set{S})$, we may obtain column \emph{projection} approximations in the same way:

\begin{definition}[Column projection approximation]
    Let $\mat{B} \in \field^{m\times n}$ be a matrix and let $\set{S} \subseteq \{1,\ldots,n\}$ be a set of \emph{pivot indices}.
    The \emph{column projection approximation} to $\mat{B}$ induced by $\set{S}$ is
    \begin{equation*}
        \Bhat \coloneqq \mat{\Pi}_{\mat{B}(:,\set{S})} \mat{B}.
    \end{equation*}
\end{definition}

Naturally, the column projection approximation $\mat{\Pi}_{\mat{B}(:,\set{S})} \mat{B}$ is the projection approximation (\cref{def:projection-approximation}) to $\mat{B}$ with test matrix $\mat{\Omega} = \Id(:,\set{S})$.
As column projection approximations are projection approximations, they enjoy the properties listed in \cref{prop:projection-properties}.\index{column projection approximation!definition|)}

\index{column projection approximation!computation by pivoted partial \QR decomposition|(}\index{pivoted partial .QR decomposition@pivoted partial \QR decomposition|(}
\subsection{Pivoted partial \QR decomposition}

Just as column Nystr\"om approximations can be computed by pivoted Cholesky decompositions,\index{pivoted partial Cholesky decomposition} column projection approximations can be computed using pivoted \QR decompositions.
We represent the low-rank approximation by a factorization $\Bhat = \mat{Q}\mat{F}^*$, where $\mat{Q}\in\field^{m\times k}$ has orthonormal columns and $\mat{F} \in \field^{n\times k}$ is a general matrix. 
Beginning from initial residual $\mat{B}^{(0)} \coloneqq \mat{B}$ , do the following for $i=0,1,2,\ldots,k-1$:
\begin{itemize}
    \item \textbf{Select a pivot.} Choose a \emph{pivot index} $s_i \in \{1,\ldots,n\}$ associated with a nonzero pivot column $\vec{b}^{(i)}_{s_i}\ne \vec{0}$.
    \item \textbf{Update approximation.} Set $\vec{q}_i \coloneqq \vec{b}_{s_i}^{(i)} / \norm{\smash{\vec{b}_{s_i}^{(i)}}}$ and $\vec{f}_i \coloneqq \mat{B}^*\vec{q}_i$.
    Observe that $\vec{q}_i^{\vphantom{*}}\vec{f}_i^* = \outprod{\vec{q}_i}\mat{B}^{(i)}$ is the orthogonal projection of the residual matrix $\mat{B}^{(i)}$ onto its $s_i$th column.
    \item \textbf{Update residual.} Define $\mat{B}^{(i+1)} \coloneqq \mat{B}^{(i)} - \vec{q}_i^{\vphantom{*}}\vec{f}_i^* = (\Id - \outprod{\vec{q}_i})\mat{B}^{(i)}$.
\end{itemize}
As written, this procedure is effectively a Gram--Schmidt orthogonalization:\index{Gram--Schmidt orthogonalization} We successively select a pivot column and orthonormalize the remaining columns against it.
\warn{Neglecting rounding errors},\index{Gram--Schmidt orthogonalization!numerical stability issues} the matrix $\mat{Q}$ produced by this procedure has orthonormal columns.
See \cref{prog:pivpartqr} for an implementation.

\myprogram{Pivoted partial \QR decomposition based on modified Gram--Schmidt orthogonalization.}{}{pivpartqr}

\index{Gram--Schmidt orthogonalization!numerical stability issues|(}\index{Gram--Schmidt orthogonalization!classical and modified Gram--Schmidt|(}
\begin{remark}[Numerical stability of Gram--Schmidt] \label{rem:gram-schmidt}
    The implementation of the pivoted partial \QR algorithm described above and presented in \cref{prog:pivpartqr} is an instance of the \warn{modified} Gram--Schmidt algorithm.
    The classical and modified Gram--Schmidt algorithms differ in the sequence in which they perform orthogonalization: Modified Gram--Schmidt orthogonalizes the entire $\mat{B}$ matrix against the selected pivot column, whereas classical Gram--Schmidt orthogonalizes only on an as-needed basis, forming columns of $\mat{B}^{(i)}$ on an as-needed basis using the formula $\vec{b}_s^{(i)} = \vec{b}_s^{\vphantom{*}} - \sum_{j=1}^i \outprod{\vec{q}_j}\vec{b}_s^{\vphantom{*}}$.
    The modified Gram--Schmidt method has significant numerical benefits over classical Gram--Schmidt; see \cite[\S19.8]{Hig02} for discussion and analysis.
    Even using this modified Gram--Schmidt procedure, the columns of the matrix $\mat{Q}$ produced by this procedure can lose orthogonality in finite-precision arithmetic.
    To ensure one obtains a $\mat{Q}$ matrix with numerically orthonormal columns, one can perform pivoted partial \QR using Householder reflectors\index{Householder QR factorization@Householder \QR factorization!for pivoted partial \QR decomposition}; see \cite[\S5.1]{GV13}.
\end{remark}\index{Gram--Schmidt orthogonalization!numerical stability issues|)}\index{Gram--Schmidt orthogonalization!classical and modified Gram--Schmidt|)}

The pivoted partial \QR algorithm computes a column projection approximation:

\begin{proposition}[Pivoted partial \QR computes a column projection approximation]
    Pivoted partial \QR decomposition (\cref{prog:pivpartqr}) with input matrix $\mat{B}$ and pivot set $\set{S}$ computes a column projection approximation $\mat{\Pi}_{\mat{B}(:,\set{S})} \mat{B} = \mat{Q}\mat{F}^*$.
\end{proposition}
\index{column projection approximation!computation by pivoted partial \QR decomposition|)}\index{pivoted partial .QR decomposition@pivoted partial \QR decomposition|)}

\index{column projection approximation!representation as a matrix factorization|(}
\subsection{Two factorizations of a column projection approximation}

In applications, there are two common ways of representing a column projection approximation in factored form.
The first way is a pivoted partial \QR decomposition:
\begin{equation*}
    \mat{\Pi}_{\mat{B}(:,\set{S})} \mat{B} = \mat{Q}\mat{F}^* \quad \text{for } \mat{Q} = \orth(\mat{B}(:,\set{S}))\text{ and } \mat{F}^* = \mat{Q}^*\mat{B}.
\end{equation*}
The pivoted partial \QR algorithm (\cref{prog:pivpartqr}) returns this factorization.

\index{interpolative decomposition|(}
Another convenient way of representing the projection approximation is the \emph{interpolative decomposition} (ID\index{interpolative decomposition}, \cite[\S13]{MT20a}), which takes the form
\begin{equation} \label{eq:exact-id}
    \mat{\Pi}_{\mat{B}(:,\set{S})} \mat{B} = \mat{B}(:,\set{S})\mat{W}^* \quad \text{where } \mat{W}^* = \mat{B}(:,\set{S})^\dagger \mat{B}.
\end{equation}
The \emph{interpolation matrix}\index{interpolative matrix} $\mat{W}$ approximates the columns of $\mat{B}$ as linear combinations of the pivot columns $\mat{B}(:,\set{S})$.
The pivot columns are represented perfectly by the interpolative decomposition\index{interpolative decomposition}, and, provided $\mat{B}(:,\set{S})$ is full rank, the pivot columns of $\mat{W}^*$ are $\mat{W}^*(:,\set{S}) = \Id$.
The interpolative decomposition\index{interpolative decomposition} is a key ingredient for rank-structured matrix computations\index{rank-structured matrix} \cite{Mar11,Wil21} and tensor network algorithms\index{tensor network} \cite{OT10,TSL24b}.

The pivoted partial \QR decomposition, as produced by many of the methods in this chapter, can easily be converted to an interpolative decomposition\index{interpolative decomposition}.
Indeed, $\mat{F}^*(:,\set{S}) = \mat{F}(\set{S},:)^*$ is an upper-triangular matrix, so
\begin{equation*}
    \mat{\Pi}_{\mat{B}(:,\set{S})} \mat{B} = \mat{Q}\mat{F}^* = \mat{Q} \mat{F}(\set{S},:)^* \mat{F}(\set{S},:)^{-*} \mat{F}^* = \mat{B}(:,\set{S}) [\mat{F}(\set{S},:)^{-*} \mat{F}^*].
\end{equation*}
Ergo, the interpolation matrix is $\mat{W} = \mat{F}\mat{F}(\set{S},:)^{-1}$, which can be computed in $\order(nk^2)$ operations using triangular solves.

To facilitate discussion of other algorithms later in this chapter, let us also introduce the notion of an \emph{inexact} interpolative decomposition\index{interpolative decomposition}.

\index{interpolative decomposition!exact vs.\ inexact|(}
\begin{definition}[Exact and inexact interpolative decomposition] \label{def:inexact-id}
    An \emph{exact} interpolative decomposition\index{interpolative decomposition} is a decomposition of the form \cref{eq:exact-id}.
    An \emph{inexact interpolative decomposition\index{interpolative decomposition}} is any low-rank approximation of the form
    \begin{equation*}
        \Bhat = \mat{B}(:,\set{S}) \smash{\hatbold{W}}^*.
    \end{equation*}
    We say a randomized algorithm for computing an inexact ID\index{interpolative decomposition} is \emph{near-optimal} if it is guaranteed to produce an output satisfying
    \begin{equation*}
        \norm{\mat{B} - \mat{B}(:,\set{S}) \smash{\hatbold{W}}^*}_{\mathrm{F}}^2 \le \mathrm{const} \cdot \norm{\mat{B} - \mat{\Pi}_{\mat{B}(:,\set{S})} \mat{B}}_{\mathrm{F}}^2
    \end{equation*}
    with 90\% probability for \warn{any} matrix $\mat{B}$ and pivot set $\set{S}$.
\end{definition}\index{column projection approximation!representation as a matrix factorization|)}\index{interpolative decomposition|)}\index{interpolative decomposition!exact vs.\ inexact|)}

\index{Gram correspondence!equivalence between Cholesky and \QR|(}\index{QR decomposition@\QR decomposition!equivalence with Cholesky decomposition|(}\index{Cholesky decomposition!equivalence with QR decomposition@equivalence with \QR decomposition|(}\index{column Nystr\"om approximation!equivalence with column projection approximation|(}\index{column projection approximation@equivalence with column Nystr\"om approximation|(}
\subsection{The Gram correspondence: Equivalence between Cholesky and \QR}

The Gram correspondence\index{Gram correspondence} (\cref{thm:gram-correspondence}) establishes a link between column projection approximation and column Nystr\"om approximation and, consequently, a link between the pivoted partial \QR and Cholesky algorithms.
The general Gram correspondence theorem (\cref{thm:gram-correspondence}) yields the following corollary:

\begin{corollary}[\QR and Cholesky] \label{cor:qrcholesky}
    Let $\mat{A}$ be a psd matrix and let $\mat{B}$ be \warn{any} Gram square root\index{Gram square root} of $\mat{A}$ (that is, $\mat{A} = \mat{B}^*\mat{B}$).
    Let $\mat{Q}$ and $\mat{F}$ be factors produced by pivoted partial \QR applied to $\mat{B}$ with pivot set $\set{S}$.
    The following conclusions hold.
    \begin{enumerate}[label=(\alph*)]
    \item \textbf{Same factor.} The pivoted partial Cholesky on $\mat{A}$ with pivot set $\set{S}$ produces the factor matrix $\mat{F}$, up to a possible scaling of the columns by elements $\omega \in \field$ of unit modulus.
    \item \textbf{Connection between low-rank approximations.} 
    The low-rank approximation $\Bhat = \mat{Q}\mat{F}^*$ produced by pivoted partial \QR is a Gram square root\index{Gram square root} of the approximation $\Ahat = \mat{F}\mat{F}^*$ produced by pivoted partial Cholesky,
    \begin{equation*}
        \Ahat = \Bhat^*\Bhat,
    \end{equation*}
    \item \textbf{Same errors.}
    For any $p\ge 1$, the approximation errors are related
    \begin{equation*}
        \norm{\mat{A} - \Ahat}_{\set{S}_p} = \norm{\mat{B} - \Bhat}_{\set{S}_{2p}}^2.
    \end{equation*}
    Here, $\norm{\cdot}_{\set{S}_p}$ denotes the Schatten $p$-norm.
    \item \textbf{Connection between residuals.}
    At every step $i$, the $i$-step \QR residual $\mat{B}^{(i)} \coloneqq \mat{B} - \mat{Q}(:,1:i)\mat{F}(:,1:i)^*$ is a Gram square root\index{Gram square root} of the $i$-step Cholesky residual $\mat{A}^{(i)} \coloneqq \mat{A} - \mat{F}(:,1:i)\mat{F}(:,1:i)^*$.
    \item \textbf{Diagonals and squared column norms.}\index{squared row or column norms}
    At each step $i$, the squared column norms of $\mat{B}^{(i)}$ are the diagonal entries of $\mat{A}^{(i)}$. \label{item:qrcholesky-diag-scn}
    \end{enumerate}
\end{corollary}

As we shall see, this result allows us to to seamlessly convert Cholesky-based algorithms for psd low-rank approximation to \QR-based algorithms for computing column projection approximations to a general matrix.
We restate the transference of algorithms principle below:\index{Gram correspondence!transference of algorithms}
\TransferenceOfAlgorithms

\index{Gram correspondence!equivalence between Cholesky and \QR|)}\index{QR decomposition@\QR decomposition!equivalence with Cholesky decomposition|)}\index{Cholesky decomposition!equivalence with QR decomposition@equivalence with \QR decomposition|)}

\index{randomly pivoted QR@randomly pivoted \QR|(}
\section{Randomly pivoted \QR} \label{sec:rpqr}

Under the transference of algorithms,\index{Gram correspondence!transference of algorithms} the \RPCholesky procedure has an analog for computing projection approximations to a general matrix using \QR decomposition.
We will call this algorithm \emph{randomly pivoted \QR} (\RPQR) \cite{DRVW06,DV06,CETW25}; its history will be discussed later in \cref{sec:rpcholesky-rpqr-history}.

The randomly pivoted \QR algorithm is straighforward: Execute a pivoted partial \QR decomposition, drawing a random pivot column at each iteration
\begin{equation*}
    s_{i+1} \sim \scn(\mat{B}^{(i)}).
\end{equation*}
sampled according to the \warn{squared} column norms\index{squared row or column norms} of the current residual matrix $\mat{B}^{(i)}$.
The squared column norm distribution is the natural analog of the diagonal sampling $s_{i+1} \sim \diag(\mat{A}^{(i)})$ used in \RPCholesky in view of \cref{cor:qrcholesky}\ref{item:qrcholesky-diag-scn}.
See \cref{prog:rpqr} for an implementation.

\myprogram{A modified Gram--Schmidt-based implementation of the randomly pivoted \QR algorithm for computing a column projection approximation.}{Subroutine \texttt{sqcolnorms} is defined in \cref{prog:sqcolnorms}.}{rpqr}

\index{randomly pivoted QR@randomly pivoted \QR!theoretical results|(}
\subsection{Theoretical results}

Under the Gram correspondence, theoretical results for \RPCholesky immediately lead to results for \RPQR.
In particular, we have the following:

\begin{corollary}[Randomly pivoted \QR] \label{cor:rpqr-error}
    Let $\mat{B} \in \field^{m\times n}$ be a matrix, and fix $r\ge 1$ and $\varepsilon \ge 0$.
    Introduce the squared relative error of the best rank-$r$ approximation:
    \begin{equation} \label{eq:rpqr-eta}
        \eta \coloneqq \frac{\norm{\mat{B} - \lowrank{\mat{B}}_r}_{\mathrm{F}}^2}{\norm{\mat{B}}_{\mathrm{F}}^2}.
    \end{equation}
    Randomly pivoted \QR produces an approximation $\Bhat$ satisfying
    \begin{equation*}
        \expect \norm{\mat{B} - \Bhat}_{\mathrm{F}}^2 \le (1+\varepsilon) \norm{\mat{B} - \lowrank{\mat{B}_r}}_{\mathrm{F}}^2.
    \end{equation*}
    provided the number of steps satisfies
    \begin{equation} \label{eq:rpqr-frob-k}
        k \ge \frac{r}{\varepsilon} + r \log \left( \frac{1}{\varepsilon \eta} \right).
    \end{equation}
\end{corollary}

More concisely, \RPQR produces an $(r,\varepsilon/2,2)$-approximation after \cref{eq:rpqr-frob-k} steps.
(Recall \cref{def:r-eps-p} and the Lyapunov inequality $\expect \norm{\mat{B} - \Bhat}_{\mathrm{F}} \le \smash{(\expect \norm{\mat{B} - \Bhat}_{\mathrm{F}}^2)^{1/2}}$.)\index{randomly pivoted QR@randomly pivoted \QR!theoretical results|)}

\index{randomly pivoted QR@randomly pivoted \QR!implementation|(}
\subsection{Implementation}

We will not spend too much time on the implementation of the standard \RPQR algorithm, as we mainly advocate the \emph{accelerated} \RPQR algorithm\index{accelerated randomly pivoted QR@accelerated randomly pivoted \QR} for practical computations, which will be introduced in the next section.
Therefore, we describe only a few implementation details for \RPQR.

The implementation of \RPQR in \cref{prog:rpqr} is based on a modified Gram--Schmidt procedure;\index{Gram--Schmidt orthogonalization!classical and modified Gram--Schmidt} see \cref{rem:gram-schmidt} for discussion.
This implementation is numerically stable\index{Gram--Schmidt orthogonalization!numerical stability issues} enough for most use cases, but it may be worth using an implementation based on Householder reflectors\index{Householder QR factorization@Householder \QR factorization!for pivoted partial \QR decomposition} if obtaining a $\mat{Q}$ matrix that has orthonormal columns up to machine accuracy is necessary.

\index{blocked algorithm!blocked versions of randomly pivoted QR@blocked versions of randomly pivoted \QR|(}
A second potential implementation issue is hardware efficiency.
As we described in the previous chapter, the fastest matrix algorithms are based on block matrix computations; the implementation of \RPQR provided in \cref{prog:rpqr} is inherently sequential, selecting a column and orthogonalizing the entire matrix against it at every iteration.
This deficit will be addressed using rejection sampling with the accelerated \RPQR algorithm,\index{accelerated randomly pivoted QR@accelerated randomly pivoted \QR} but there are other fixes that have been proposed for pivoted \QR decomposition with deterministic greedy pivoting that are worth mentioning.
For fully deterministic implementations, the state of the art is provided by the \LAPACK routine \texttt{xGEQP3}.
Quoting from the \LAPACK manual, this routine ``only updates one column and one row of the rest of the matrix (information necessary for the next pivoting phase) and delays the update of the rest of the matrix until a block of columns has been processed'' \cite[\S2.4.2.3]{ABB+99}.
This modification improves the efficiency of pivoted \QR decompositions significantly and reorganizes computations so a significant fraction of them utilize block matrix operations.\index{blocked algorithm!blocked versions of randomly pivoted QR@blocked versions of randomly pivoted \QR|)}
A second more recent idea is to use randomized dimensionality reduction, applying the slow sequential pivoted \QR algorithm to the matrix $\mat{B}$ after it has been compressed using randomized dimensionality reduction; see \cite{MBM+23} for details.\index{randomized CholeskyQR}
\index{randomly pivoted QR@randomly pivoted \QR!implementation|)}

\index{block randomly pivoted QR@block randomly pivoted \QR|(}\index{randomly pivoted QR@randomly pivoted \QR!blocked versions|(}
\subsection{Block randomly pivoted \QR}

For historical reasons, we also mention that there is a straightforward block implementation of \RPQR, analogous to block \RPCholesky (\cref{sec:block-rpcholesky}).
For most use cases, robust block random pivoting\index{robust blockwise random pivoting!for general low-rank approximation} \cite{DCMP23} or accelerated \RPQR (introduced next section) provide a better alternative.\index{block randomly pivoted QR@block randomly pivoted \QR|)}

\index{accelerated randomly pivoted QR@accelerated randomly pivoted \QR|(}\index{rejection sampling!for implementing randomly pivoted QR@for implementing randomly pivoted \QR|(}
\section{Accelerated randomly pivoted \QR} \label{sec:acc-rpqr}

The accelerated \RPCholesky method (introduced in \cref{sec:acc-rpcholesky}) uses rejection sampling to simulate the performance of the original \RPCholesky algorithm.
It is much faster than the standard \RPCholesky algorithm in practice due to runtime efficiencies from block-wise computations, and it produces the \warn{same random output} as ordinary \RPCholesky.
Under the Gram correspondence,\index{Gram correspondence!transference of algorithms} the accelerated \RPCholesky algorithm has an analog, accelerated \RPQR, that simulates the performance of the \RPQR algorithm while running much raster due to block-matrix arithmetic.

\myparagraph{Description of algorithm}
One ``round'' of the accelerated \RPQR algorithm may be described as follows.
Suppose we have already generated pivots $s_1,\ldots,s_i$ sampled from the \RPQR distribution, and, for each $j$, denote the residual as
\begin{equation*}
    \mat{B}^{(j)} \coloneqq \mat{B} - \mat{Q}^{(j)}\big(\mat{Q}^{(j)}\big)^*\mat{B} \quad \text{for }\mat{Q}^{(j)} = \orth(\mat{B}(:,\{s_1,\ldots,s_j\})).
\end{equation*}
We wish to generate new pivots $s_{i+1},\ldots,s_{i+\ell}$ generated from the \RPQR distribution
\begin{equation*}
    s_{i+j} \sim \scn(\mat{B}^{(i+j)}).
\end{equation*}
To do so,we fix a block size $b\ge 1$ and draw a collection $\set{S}'=\{s_1',\ldots,s_b'\}$ iid \warn{with replacement} from the squared column norm distribution of $\mat{B}^{(i)}$
\begin{equation*}
    s_1',\ldots,s_b'\simiid \scn(\mat{B}^{(i)}).
\end{equation*}
Then, we filter the proposals $\set{T}\subseteq\set{S}'$ using rejection sampling, which can be accomplished using the \RejectionSampleSubmatrix subroutine (\cref{prog:rejection_sample_submatrix}) on $\mat{H} = [\mat{B}^{(i)}(:,\set{S}')]^*\mat{B}^{(i)}(:,\set{S}')$ with $\vec{u} = \diag(\mat{H})$; see \cref{sec:acc-rpcholesky} for details on \RejectionSampleSubmatrix.
Once the new pivots $s_{i+1},\ldots,s_{i+\ell}$ are inducted, we update $\mat{Q}^{(i+\ell)}$, columns of the factor matrix $\mat{F}$, and residual $\mat{B}^{(i+\ell)}$ as needed.

Below, we describe two implementations of accelerated \RPQR using different types of orthogonalization.
The first implementation (\cref{prog:acc_rpqr_bgs}) is based on block Gram--Schmidt orthogonalization,\index{Gram--Schmidt orthogonalization!in blocks} and the second (\cref{prog:acc_rpqr}) employs Householder reflectors.\index{Householder QR factorization@Householder \QR factorization!for pivoted partial \QR decomposition}

\index{Gram--Schmidt orthogonalization!in blocks|(}
\myprogram{Block Gram--Schmidt-based implementation of accelerated \RPQR algorithm for computing a column projection approximation to a general matrix.}{Subroutines \texttt{rejection\_sample\_submatrix} and \texttt{sqcolnorms} are provided in \cref{prog:rejection_sample_submatrix,prog:sqcolnorms}.}{acc_rpqr_bgs}

\myparagraph{Implementation \#1: Block Gram--Schmidt}
Our first implementation of accelerated \RPQR is shown in \cref{prog:acc_rpqr_bgs} and uses block Gram--Schmidt\index{Gram--Schmidt orthogonalization!in blocks} to perform orthogonalization.
It is a fairly direct extension of the modified Gram--Schmidt implementation of standard \RPQR in \cref{prog:rpqr}.
Throughout the algorithm, \texttt{B} stores the current residual, the columns of the left orthonormal factor $\mat{Q}$ are obtained by computing a \QR decomposition of selected columns of \texttt{B}, and the factor matrix is computed using the handy formula $\mat{F}(:,i+1:i+\ell) = (\mat{B}^{(i)})^*\mat{Q}(:,i+1:i+\ell)$.
Following best practices \cite{GLRE05}, we orthonormalize the residual matrix $\mat{B}$ against the newly orthonormalized columns $\mat{Q}(:,i+1:i+\ell)$ \warn{twice} to improve numerical stability.
This implementation has worse stability properties than the Householder reflector-based implementation in \cref{prog:acc_rpqr}.\index{Gram--Schmidt orthogonalization!numerical stability issues}\index{Gram--Schmidt orthogonalization!in blocks|)}

\index{Householder QR factorization@Householder \QR factorization!for pivoted partial \QR decomposition|(}
\myparagraph{Implementation \#2: Householder \QR}
A more numerically robust implementation of accelerated \RPQR can be obtained using Householder reflectors.
Rather than storing the orthonormal factor $\mat{Q}$ directly, we represent $\mat{Q}$ implicitly as the first $k$ columns of a product $\mat{H}_1\cdots \mat{H}_k$ of \emph{Householder reflector matrices} $\mat{H}_i$, each of which takes the form $\mat{H}_i = \Id - 2\outprod{\vec{u}_i}$ for a unit vector $\vec{u}_i$.
This representation remains storage-efficient, as we only need to store the unit vectors $\vec{u}_i$.
See \cite[\S\S5.1 \& 5.2]{GV13} for more on Householder reflectors and \QR factorization. 

The advantage of a Householder reflector-based implementation over a Gram--Schmidt-based implementation is that Householder \QR methods are guaranteed to produce a $\mat{Q}$ matrix whose columns are orthonormal up to machine precision \cite[\S19.3]{Hig02}.\index{Gram--Schmidt orthogonalization!numerical stability issues}
Gram--Schmidt-type methods do not possess such a guarantee \cite[\S19.8]{Hig02}.

\myprogram{Householder reflector-based implementation of accelerated \RPQR algorithm for computing a column projection approximation to a general matrix.}{Subroutines \texttt{rejection\_sample\_submatrix}, \texttt{hhqr}, \texttt{apply\_Qt}, \texttt{get\_Q}, and \texttt{sqcolnorms} are provided in \cref{prog:rejection_sample_submatrix,prog:hhqr,prog:apply_Qt,prog:get_Q,prog:sqcolnorms}.}{acc_rpqr}

\Cref{prog:acc_rpqr} provides a Householder-reflector based implementation of accelerated \RPQR.
This implementation uses compact representations of the Householder \QR decomposition and carefully updates a \QR decomposition of the matrix $\mat{B}(:,\set{S})$ throughout the algorithm, even as the pivot set $\set{S}$ increases in size.
See \cref{app:incremental-qr} for discussion.\index{incremental QR decomposition!incremental \QR decomposition}\index{QR decomposition@\QR decomposition!incremental}\index{randomly pivoted QR@randomly pivoted \QR!blocked versions|)}\index{accelerated randomly pivoted QR@accelerated randomly pivoted \QR|)}\index{rejection sampling!for implementing randomly pivoted QR@for implementing randomly pivoted \QR|)}\index{Householder QR factorization@Householder \QR factorization!for pivoted partial \QR decomposition|)}

\index{randomly pivoted QR@randomly pivoted \QR|)}

\index{sketchy pivoting|(}
\section{Related work: Sketchy pivoting} \label{sec:sketchy-pivoting}

Per-Gunnar Martinsson and collaborators have pioneered a different family of methods for column subset selection called \emph{sketchy pivoting methods} \cite{VM17,DM23,DCMP23}.
For completeness, we provide a brief introduction to the idea.

\subsection{Fast random embeddings}

The starting point for these methods is a type of matrix that we will call a \emph{fast random embedding} or \emph{fast sketching matrix}.\index{sketching matrix}\index{subspace embedding!fast constructions}
Informally, a fast random embedding $\mat{S} \in \field^{m\times d}$ is a random matrix satisfying the following two properties:
\begin{itemize}
    \item \warn{Fix} an arbitrary $k$-dimensional subspace $\set{U}\subseteq\field^m$.
    With high probability, $\mat{S}^*$ preserves the lengths of \warn{all} vectors in $\set{U}$ up to a constant factor
    \begin{equation} \label{eq}
        \mathrm{c}_1 \norm{\vec{x}} \le \norm{\mat{S}^*\vec{x}} \le \mathrm{c}_2 \norm{\vec{x}} \quad \text{for all } \vec{x} \in \set{U}.
    \end{equation}
    Here, the $\mathrm{c}_1,\mathrm{c}_2 \approx 1$ are absolute constants.
    The embedding dimension $d$ should be nearly proportional to the subspace dimension, $d = \order(k)$ or $d = \order(k \log k)$.
    \item For any matrix $\mat{B} \in \field^{m\times n}$, the product $\mat{S}^*\mat{B}$ can be computed in at most $\order(mn \log m)$ operations.
\end{itemize} 
Sketching will play a major role in \cref{part:sketching} of this thesis; see that section for formal definitions, constructions of fast random embeddings, and guidance on which type of embedding to use.

\subsection{Sketching pivoting}

Suppose we are interested in selecting $k$ columns of a matrix $\mat{B}$, where the parameter $k$ is \warn{known to the algorithm in advance}.
Using the sketchy pivoting approach, we first draw a fast random embedding of size, say $d = 2k$, and compute $\mat{S}^*\mat{B}$.
Then, to select a subset of columns for $\mat{B}$, we apply some traditional column subset selection algorithm, such as a pivoted matrix decomposition, to the \emph{sketched matrix} $\mat{S}^*\mat{B}$.
The key feature of sketchy pivoting is that the traditional matrix decompositions needs only be applied to the sketched matrix $\mat{S}^*\mat{B}$, which is much smaller than the original matrix $\mat{B}$.

One major strength of sketchy pivoting is that the randomized embedding $\mat{S}^*$ has the effect of \emph{regularizing} the matrix $\mat{B}$.
Empirically, one finds that one obtains a nearly optimal set of columns even when one uses a poor column subset selection algorithm on $\mat{S}^*\mat{B}$ such as \LU with partial pivoting (applied to the adjoint $(\mat{S}^*\mat{B})^*$) \cite{DM23a}.
We will abbreviate the sketchy pivoting with partial-pivoted \LU scheme as \SkLUPP.
With an appropriate type of embedding (namely, a sparse sign embedding\index{sparse sign embedding} with constant sparsity parameter $\zeta$; see \cref{sec:sse}), \SkLUPP requires just $\order(mn + k^2n)$ operations \warn{to select the column set $\set{S}$}, faster than \RPQR and its variants.

\subsection{Quickly computing an inexact interpolative decomposition}

\index{interpolative decomposition!exact vs.\ inexact|(}\index{oversampled sketchy interpolative decomposition (OSID)|(}\index{interpolative decomposition!oversampled sketchy|(}
To upgrade a skeleton set $\set{S}$ computed by a sketchy pivoting method to an \emph{exact} interpolative decomposition\index{interpolative decomposition} (in the sense of \cref{def:inexact-id}) requires $\order(kmn)$ operations, the same asymptotic cost as \RPQR and its variants.
However, a \emph{nearly optimal inexact} ID\index{interpolative decomposition} (also as in \cref{def:inexact-id}) can be obtained rapidly by applying another level of sketching.
The idea is due to \cite{DCMP23} and is called \emph{oversampled sketchy interpolative decomposition}.
The idea is simple and clever: Generate another random embedding $\mat{\Phi} \in \field^{m\times d}$ and solve the sketched least-squares problem
\begin{equation*}
    \hatbold{W} = \argmin_{\hatbold{W}} \norm{\mat{\Phi}^*\mat{B} - (\mat{\Phi}^*\mat{B}(:,\set{S})) \smash{\hatbold{W}}^*}_{\mathrm{F}} = [(\mat{\Phi}^*\mat{B}(:,\set{S}))^\dagger (\mat{\Phi}^*\mat{B})]^*.
\end{equation*}
OSID is essentially an instantiation of the \emph{sketch-and-solve}\index{sketch-and-solve} method for computing an approximate least-squares solution, which will be discussed in \cref{sec:sketch-and-solve}.
Theoretical guarantees for sketch-and-solve are discussed in \cref{sec:sketch-and-solve,app:sketch-and-solve-analysis}, which confirm that---with appropriate parameter choices---OSID compiutes a nearly optimal inexact ID\index{interpolative decomposition}.
The cost of sketchy pivoting with OSID may be as low as $\order(mn+k^2n)$ operations.
\index{interpolative decomposition!exact vs.\ inexact|)}\index{oversampled sketchy interpolative decomposition (OSID)|)}\index{interpolative decomposition!oversampled sketchy|)}

\subsection{\RPQR vs.\ sketchy pivoting}

There are three main advantages of \RPQR-based approaches over sketchy pivoting.
First, and least importantly in my opinion, \RPQR-based approaches come with theoretical guarantees of performance.
At present, sketchy pivoting lacks such theoretical guarantees, thought it has been put through extensive numerical testing.
Second, \RPQR produces an exact ID\index{interpolative decomposition}, whereas the fast $\order(mn + k^2n)$ version of \SkLUPP produces an inexact ID\index{interpolative decomposition}.
One can use the column set from \SkLUPP to compute an exact ID\index{interpolative decomposition}, but then the total asymptotic cost of $\order(kmn)$ is the same as accelerated \RPQR.
Finally, to be most efficient, sketchy pivoting methods require (an upper bound on) the column subset size $k$ to be known in advance.
To determine the number of columns $k$ at runtime may require periodically regenerating the embedding matrix $\mat{S}^*$ and recomputing $\mat{S}^*\mat{B}$ with a larger value of $d$.
By contrast, \RPQR methods can be run for a general number of steps $k$ and stopped at whatever iteration one pleases.

Notwithstanding these limitations, sketchy pivoting methods are asymptotically faster for subset selection and inexact ID\index{interpolative decomposition} than \RPQR methods.
\iffull
\else
I regard both techniques (and adaptive randomized pivoting methods \cite{CK24}) as effective, general methods for column subset selection and ID\index{interpolative decomposition} computation.
\fi

\index{sketchy pivoting|)}

\index{randomly pivoted QR@randomly pivoted \QR!numericalresults|(}
\section{Experiments} \label{sec:rpqr-experiments}

Throughout this section, we use the Householder \QR-based implementation of accelerated \RPQR (\cref{prog:acc_rpqr}) in all experiments.
We also compare pivoted partial \QR with greedy pivoting (often simply called \emph{column-pivoted \QR}, CPQR in the numerical literature)\index{greedy pivoted QR@greedy pivoted \QR} and \SkLUPP.
(Similarly to greedy pivoted Cholesky, CPQR selects the largest-norm column as pivot at every step of \QR decomposition.)
For \SkLUPP, we compute the interpolative decomposition\index{interpolative decomposition} exactly (i.e., not with OSID\index{oversampled sketchy interpolative decomposition (OSID)}\index{interpolative decomposition!oversampled sketchy}) and choose the sketching matrix to be a sparse sign embedding\index{sparse sign embedding} with embedding dimension $d=2k$; see \cref{sec:sse}.
We use the test matrix $\mat{B} \in \field^{n\times n}$ of dimension $n=2500$ with entries
\begin{equation} \label{eq:rpqr-test-matrix}
    b_{ij} = \frac{1}{\norm{\vec{x}_i - \vec{y}_j}} \quad \text{for every } i,j = 1,\ldots,n.
\end{equation}
Here, points $\{\vec{x}_i\}$ and $\{\vec{y}_j\}$ are equispaced Cartesian grids on $[0,1]\times [0,1]$ and $[1,2]\times [0,1]$, respectively.
Matrices similar to this appear as discretizations of integral operators,\index{integral operator!discretizations of} and significant work has gone into developing ways of approximating them (e.g., by proxy point methods\index{proxy point method} \cite{YXY20}).

\begin{figure}
    \centering
    \includegraphics[width=0.7\linewidth]{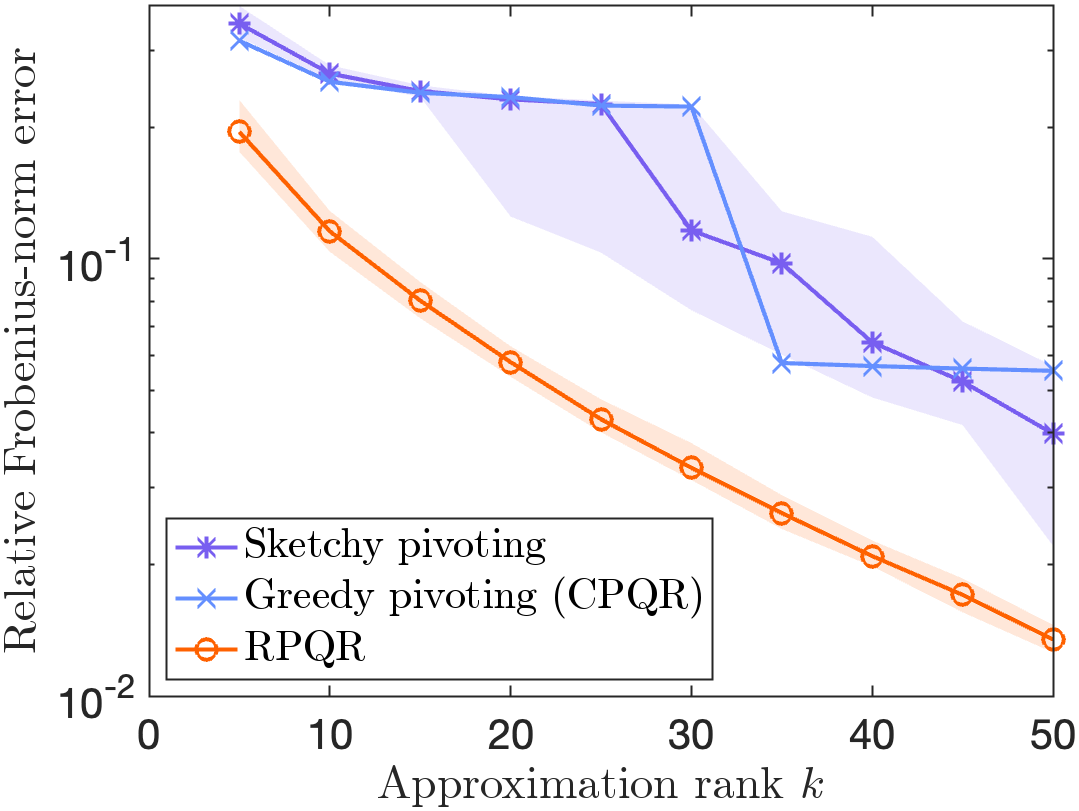}
    \caption[Comparison of randomly pivoted \QR, greedy pivoted \QR, and sketchy pivoting]{Accuracy of \RPQR (orange circles), greedy pivoted \QR (purple asterisks), and sketchy pivoting method (\SkLUPP, blue crosses) for computing a low-rank approximation to the matrix \cref{eq:rpqr-test-matrix}. 
    Lines shown median of 100 trials, and shaded regions show 10\% and 90\% quantiles.}
    \label{fig:rpqr-accuracy}
\end{figure}

\myparagraph{Accuracy}
\Cref{fig:rpqr-accuracy} presents an accuracy comparison of \RPQR, greedy pivoted \QR,\index{greedy pivoted QR@greedy pivoted \QR} and \SkLUPP on this matrix.
We see that \RPQR consistently achieves smaller errors than the other methods on this matrix, demonstrating the virtues of incorporating randomness into pivot selection for computing a pivoted partial \QR decomposition.
At rank $k=50$ on this example, \RPQR is $3.0\times$ more accurate than \SkLUPP and $4.2\times$ more accurate than greedy pivoted \QR.

\begin{figure}
    \centering
    \includegraphics[width=0.7\linewidth]{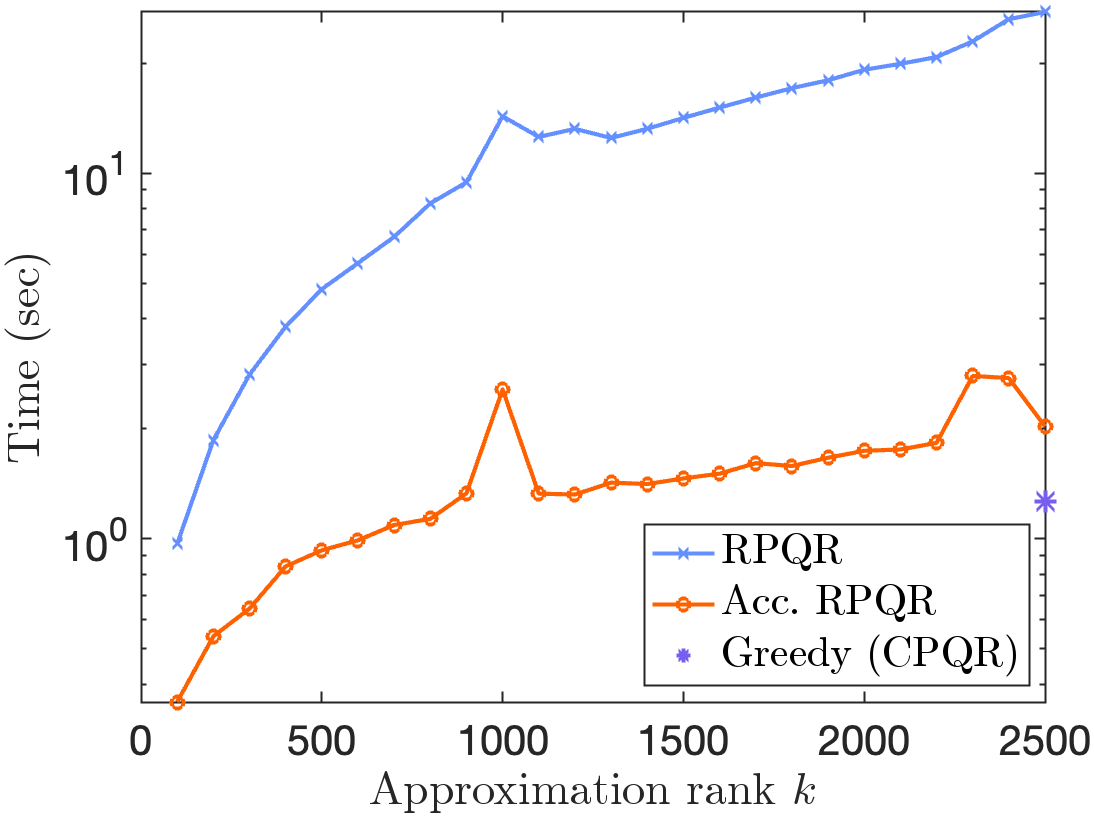}
    \caption[Runtime for standard randomly pivoted \QR, accelerated randomly pivoted \QR, and MATLAB's greedy pivoted \QR for full pivoted \QR decomposition]{Runtime for standard \RPQR (blue crosses) and accelerated \RPQR (orange circles) for different approximation ranks $k$ (single execution).
    The runtime for MATLAB's (full) greedy pivoted \QR (i.e., \texttt{[Q,R,P] = qr(B)}) is shown as a purple asterisk.
    }
    \label{fig:rpqr-time}
\end{figure}

\myparagraph{Speed}
\Cref{fig:rpqr-time} compares the runtime of the ordinary and accelerated \RPQR\index{accelerated randomly pivoted QR@accelerated randomly pivoted \QR} algorithms (\cref{prog:rpqr,prog:acc_rpqr}).
For accelerated \RPQR,\index{accelerated randomly pivoted QR@accelerated randomly pivoted \QR} we set the block size to $b\coloneqq \max \{ k/2, 200\}$.
We see that accelerated \RPQR\index{accelerated randomly pivoted QR@accelerated randomly pivoted \QR} is $13.5\times$ faster than standard \RPQR when we reach a rank of $k=n$, at which point we have computed a full pivoted \QR decomposition of the entire matrix.
At this point, we can compare the performance of accelerated \RPQR\index{accelerated randomly pivoted QR@accelerated randomly pivoted \QR} to MATLAB's greedy column-pivoted \QR decomposition,\index{greedy pivoted QR@greedy pivoted \QR} executed as \texttt{[Q,R,P] = qr(B)}.
We omit comparisons with sketchy pivoting as it is not appropriate for computing a full pivoted \QR decomposition.
We see that MATLAB's built-in column-pivoted \QR decomposition subroutine, calling appropriate \LAPACK routines written in Fortran, is only $1.6\times$ faster than accelerated \RPQR,\index{accelerated randomly pivoted QR@accelerated randomly pivoted \QR} written in MATLAB.
I consider these timing results to be promising, and it suggests that an implementation of accelerated \RPQR in a low-level programming language might be faster than existing \LAPACK routines for column-pivoted \QR.\index{randomly pivoted QR@randomly pivoted \QR!numericalresults|)}

\myparagraph{More experiments}
Many further demonstration of the virtues of random pivoting for approximation of general matrices appear in the paper \cite{DCMP23}, including comparisons of the blocked and unblocked \RPQR algorithm and the RBRP \QR algorithm.\index{robust blockwise random pivoting!for general low-rank approximation}

\index{randomly pivoted Cholesky!history|(}\index{randomly pivoted QR@randomly pivoted \QR!history|(}
\section{\RPCholesky and \RPQR: History} \label{sec:rpcholesky-rpqr-history}

In our presentation, we introduced \RPCholesky first and used the transference of algorithms principle to derive the \RPQR algorithm.
Historically, it went the other way around.
Here is the story, as best I can tell it.

\index{squared column-norm sampling for column projection approximation|(}
\myparagraph{Squared column norm sampling}
In the late 1990s, theoretical computer scientists and researchers in adjacent fields became interested in using randomization to accelerate matrix computations, motivated by applications in data analysis \cite{PTRV98,FKV98}.
As a fast way of constructing a low-rank approximation to a large matrix, Frieze, Kannan, and Vempala suggested approximating a matrix by low-rank approximations spanned by columns $s_1,\ldots,s_k \stackrel{\text{iid}}{\sim} \scn(\mat{B})$ drawn \warn{iid} from the squared column norm distribution \cite{FKV98}.\index{squared row or column norms}
Their main result \cite[Thm.~2]{FKV98} establishes the guarantee
\begin{equation*}
    \expect \norm{\mat{B} - \mat{\Pi}_{\mat{B}(:,\{s_1,\ldots,s_k\})} \mat{B}}_{\mathrm{F}}^2 \le \norm{\mat{B} - \lowrank{\mat{B}}_r}_{\mathrm{F}}^2 + \frac{r}{k} \norm{\mat{B}}_{\mathrm{F}}^2 \quad \text{for } s_1,\ldots,s_k \stackrel{\text{iid}}{\sim} \scn(\mat{B}).
\end{equation*}
Here, $r\ge 1$ is any fixed target rank.
Introducing the squared relative-error of the best rank-$r$ approximation $\eta$ from \cref{eq:rpqr-eta}, we see that $k = \order(r/\eta)$ columns are sufficient to obtain a low-rank approximation comparable to the best rank-$r$ approximation to $\mat{B}$.
This analysis is tight:
\warn{For a worst-case input matrix}, $k=\Omega(r/\eta)$ columns from the squared column norm distribution are necessary to obtain this guarantee; see \cite[Thm.~C.3(b)]{CETW25}.\index{squared column-norm sampling for column projection approximation|)}

\myparagraph{Randomly pivoted \QR introduced as ``adaptive sampling''}
The Frieze--Kannan--Vempala algorithm suggests a natural improvement by iteration.
Draw a single column (or a small batch) from the squared column norm distribution, subtract off the projection of the matrix onto these column(s), and repeat. 
This refinement was proposed in 2006 by Deshpande, Rademacher, Vempala, and Wang \cite{DRVW06} and expanded on later that year by Deshpande and Vempala \cite{DV06}; these authors called their method \emph{adaptive sampling}.
The algorithms we have called \RPQR and block \RPQR\index{block randomly pivoted QR@block randomly pivoted \QR} are essentially the same as Deshpande et al.'s adaptive sampling algorithm, up to implementation details.
The papers \cite{DRVW06,DV06} also proposed \emph{volume sampling},\index{volume sampling} which draws a subset $\set{S}$ of $k$ columns of $\mat{B}$ with probability proportional to the squared volume $\det[\mat{B}(:,\set{S})^*\mat{B}(:,\set{S})]$; in modern language, the $k$-volume sampling distribution is the $k$-DPP distribution\index{determinantal point process sampling!connection to volume sampling} on the Gram matrix\index{Gram matrix} $\mat{B}^*\mat{B}$.

Desphande and coauthors combine their three primitives \RPQR, block \RPQR,\index{block randomly pivoted QR@block randomly pivoted \QR} and volume sampling\index{volume sampling} in creative ways to prove existence results and algorithms for computing column projection approximations.
Underlying these combinations are two main theoretical results, which we now summarize.
The first result \cite[Thm.~1.2]{DRVW06} shows that, with a sufficiently large block size, block \RPQR\index{block randomly pivoted QR@block randomly pivoted \QR} produces accurate low-rank approximations.

\begin{fact}[Block \RPQR: Large block size] \label{fact:block-rpqr}
    Let $\mat{B} \in \field^{m\times n}$ be a matrix, and fix target rank $r\ge 1$ and accuracy parameter $\varepsilon \in (0,1)$.
    Set the block size of block \RPQR to be $b\ge r/\varepsilon$, and execute the procedure for $t$ (block) steps.
    Then block \RPQR outputs a low-rank approximation $\Bhat$ satisfying
    \begin{equation*}
        \expect \norm{\mat{B} - \Bhat}_{\mathrm{F}}^2 \le (1-\varepsilon)^{-1} \norm{\smash{\mat{B} - \lowrank{\mat{B}}_r}}_{\mathrm{F}}^2 + \varepsilon^t \norm{\mat{B}}_{\mathrm{F}}^2.
    \end{equation*}
\end{fact}

This result establishes the qualitative conclusion that block \RPQR, \warn{with a sufficiently large block size}, produces an approximation comparable to the best rank-$r$ approximation after drawing $k = \order(r \log(1/\eta))$ columns, where $\eta$ is defined in \cref{eq:rpqr-eta}.
This result demonstrates an \emph{exponential separation} between iid squared column norm sampling\index{squared column norm sampling for column projection approximation} and block \RPQR\index{block randomly pivoted QR@block randomly pivoted \QR} in the parameter $\eta$.

The second result provides a sharp bound on column projection approximations built from the volume sampling\index{volume sampling} distribution and a very crude bound for the plain \RPQR algorithm:

\begin{fact}[\RPQR and volume sampling]
    Let $\mat{B} \in \field^{m\times n}$ and construct $k$-column projection approximations $\Bhat_{\mathrm{vol}}$ and $\Bhat_{\RPQR}$ using volume sampling and \RPQR (with block size $1$), respectively.
    Then
    \begin{align}
        \expect \norm{\mat{B} - \Bhat_{\mathrm{vol}}}_{\mathrm{F}}^2 &\le (k+1) \norm{\mat{B} - \lowrank{\mat{B}}_k}_{\mathrm{F}}^2, \label{eq:dv-vs} \\
        \expect \norm{\mat{B} - \Bhat_{\RPQR}}_{\mathrm{F}}^2 &\le (k+1)! \norm{\mat{B} - \lowrank{\mat{B}}_k}_{\mathrm{F}}^2. \label{eq:dv-rpqr}
    \end{align}
\end{fact}

The volume sampling\index{volume sampling} bound is \cite[Thm~1.3]{DRVW06}, and the \RPQR bound is \cite[Prop.~2]{DV06}.
The volume sampling bound\index{volume sampling} \cref{eq:dv-vs} is optimal in the sense that no method can achieve a prefactor smaller than $k+1$ (which follows by \cref{thm:optimal-nys-approx-2} and the Gram correspondence\index{Gram correspondence}).
By contrast, the bound \cref{eq:dv-rpqr} is quite weak, and seems to suggest that the unblocked versions of \RPQR produce approximations of extremely low quality.
Using the weak bound \cref{eq:dv-rpqr}, Deshpande and Vempala produce algorithms with accuracy guarantees by combining $r$ steps of \RPQR with multiple rounds of block \RPQR\index{block randomly pivoted QR@block randomly pivoted \QR} with a large block size $b = \Omega(r)$.
The very weak bound \cref{eq:dv-vs} is the only analysis of \RPQR with block size $b=1$ established by Deshpande et al.

\myparagraph{Empirical evaluation of adaptive sampling/\RPQR}
In the years following the publication of Deshpande and coauthors work, there were various efforts to analyze the empirical performance of adaptive sampling/\RPQR-type algorithms and compare them to alternatives.
There was particular interest in using variants of the adaptive sampling idea for Nystr\"om approximation of psd matrices.
A comparison of column Nystr\"om methods by Kumar, Mohri, and Talwalkar concluded \cite[p.~989]{KMT12} that adaptive sampling ``requires a full pass through [the kernel matrix] $\mat{K}$ at each iteration and is thus inefficient for large $\mat{K}$''.
As solution, the papers \cite{KMT12,WZ13} proposed cheaper, approximate versions of adaptive sampling that were more tractable.
Despite its appealing properties and attempts to address the method's weaknesses, the adaptive sampling algorithm has not seen wide use for applied computation in the 2010s and early 2020s; indeed, the algorithm is not mentioned in a recent comparison of popular methods for column selection \cite{DM23a}.

\myparagraph{Randomly pivoted Cholesky}
In 2017, Musco and Woodruff considered the problem of computing a low-rank approximation to a psd matrix from a small number of entry accesses \cite{MW17a}.
Their goal was to produce a low-rank approximation $\Ahat$ to $\mat{A}$ that is competitive with the best rank-$r$ approximation when measured in the \warn{Frobenius norm}:\index{Frobenius-norm positive-semidefinite low-rank approximation}
\begin{equation} \label{eq:frob-psd-lra}
    \norm{\mat{A} - \Ahat}_{\mathrm{F}}^2 \le (1+\varepsilon) \norm{\mat{A} - \lowrank{\mat{A}}_r}_{\mathrm{F}}^2 \quad \text{with } \rank(\Ahat) \le \poly(r,1/\varepsilon).
\end{equation}
Producing a low-rank approximation of this quality from a small number of entry accesses is a stringent requirement; see \cref{sec:frobenius-psd} for more discussion.\index{Frobenius-norm positive-semidefinite low-rank approximation}

To motivate why psd low-rank approximation is even \emph{possible} without reading the whole input matrix, Musco and Woodruff use what we have called the transference of algorithms principle\index{Gram correspondence!transference of algorithms} to extend adaptive sampling to an algorithm for psd low-rank approximation.
They write:
\begin{quote}
    Since $\mat{A}^{1/2}\mat{A}^{1/2} = \mat{A}$, the entry $a_{ij}$ is just the dot product between the $i$th and $j$th columns of $\mat{A}^{1/2}$. So with $\mat{A}$ in hand, the dot products have been `precomputed' and [adaptive sampling] yields a low-rank approximation algorithm for $\mat{A}^{1/2}$ running in just $n\cdot \poly(k/\varepsilon)$ time. Note that, aligning with our initial intuition that reading the diagonal entries of $\mat{A}$ is necessary to avoid the $\operatorname{nnz}(\mat{A})$ time lower bound for general matrices [i.e., \cref{prop:general-lra-impossibility}],\index{low-rank approximation!impossibility results} the diagonal entries of $\mat{A}$ are the column norms of $\mat{A}^{1/2}$, and hence their values are critical to computing the adaptive sampling probabilities.
\end{quote}
The algorithm Musco and Woodruff describe is, in its essence, \RPCholesky.

The discussion of \RPCholesky in Musco and Woodruff's work is brief, serving to motivate the development of more sophisticated algorithms with stronger theoretical guarantees.\index{Frobenius-norm positive-semidefinite low-rank approximation}
The original paper of Musco and Woodruff does not provide pseudocode for \RPCholesky, document an implementation of the algorithm, or report any numerical experiments.
(I have since learned from Cameron and Christopher Musco that numerical experiments were done and not published.)
Also, given the brief treatment, some the subtleties of the adaptive sampling in Desphande and coauthor's work go unmentioned, such as the specific mixtures of unblocked and blocked adaptive sampling that Deshpande and Vempala use to obtain algorithms with relative error guarantees \cite[\S3.2]{DV06}.
To the best of my knowledge, Musco and Woodruff's work is the sole mention of the \RPCholesky algorithm for psd low-rank approximation in the literature prior to our work.
(As mentioned in \cref{sec:dpp-connections}, Poulson \cite{Pou20} also uses the same computational steps as \RPCholesky for projection DPP sampling.)

Our paper \cite{CETW25}, originally released as a preprint in 2022, revisited the \RPCholesky algorithm, reinterpreted the procedure as a partial Cholesky decomposition, and suggested the name \emph{randomly pivoted Cholesky}.
To simply the algorithm and remove the large block size limitation in \cref{fact:block-rpqr}, we established new theoretical results ensuring that \RPCholesky produces near-optimal low-rank approximations with block size $b=1$.
(Our subsequent work \cite{ETW24a} extends this to general block sizes.)
This paper provided numerical experiments comparing \RPCholesky to other algorithms for column Nystr\"om and using \RPCholesky to accelerate kernel computations in scientific machine learning.

\myparagraph{From adaptive sampling to randomly pivoted \QR}
After we released our paper on \RPCholesky, we were surprised by the amount of interest our colleagues had in the \QR version of the procedure, even though it requires significantly more computational work than \RPCholesky ($\order(kn^2)$ operations versus $\order(k^2n)$ operations).
This feedback, as well as input from anonymous referees during peer review, led us to devote more attention to the column projection approximations for general matrices in revisions of our \RPCholesky manuscript \cite{CETW25} and in its followup \cite{ETW24a}; I have tried to add further details in this thesis.

Let me summarize our contributions and proposals for \RPQR.
First, we reinterpreted the adaptive sampling procedure of Deshpande and coauthors as a pivoted partial \QR decomposition with a random pivoting rule, and suggested the name \emph{randomly pivoted \QR}.
We hope the new name will make this algorithm more transparent to researchers in numerical analysis and scientific computing and help to differentiate between the \QR and Cholesky versions of the adaptive sampling idea.
Second, we proved new theoretical results \cite[Cor.~5.2]{CETW25} and \cite[Cor.~5.2]{ETW24a}.
These results demonstrate that \RPCholesky and \RPQR produce approximations comparable to the best rank-$r$ approximation when the block size $b$ is set to any value between $1$ and $\order(r)$.
In particular, these results justify the use of the algorithm without the multistep procedures and large block sizes used in Deshpande et al.'s existing work.
Third, we developed accelerated block implementations of \RPQR (\cref{sec:acc-rpqr}).
Finally, in this thesis, I discussed stable numerical implementations using modified Gram--Schmidt\index{Gram--Schmidt orthogonalization!classical and modified Gram--Schmidt} and Householder reflectors;\index{Householder QR factorization@Householder \QR factorization!for pivoted partial \QR decomposition} related ideas appear in \cite{CK24}.
\index{randomly pivoted Cholesky!history|)}\index{randomly pivoted QR@randomly pivoted \QR!history|)}

\chapter{\texorpdfstring{\CURnoindex decompositions}{CUR decompositions}} \label{ch:cur}

\epigraph{Thus, $\mat{C}$ and/or $\mat{R}$ can be used in place of the eigencolumns and eigenrows, but since they consist of actual data elements they will be interpretable in terms of the field from which the data are drawn (to the extent that the original data points and/or features are interpretable).}{Michael W.\ Mahoney and Petros Drineas, \emph{\CUR matrix decompositions for improved data analysis} \cite{MD09a}}

A \CUR decomposition or \CUR approximation refers to a low-rank approximation of a matrix $\mat{B}$ of the form
\begin{equation*}
    \mat{B} \approx \mat{C} \mat{U} \mat{R} \quad \text{where } \mat{C} = \mat{B}(:,\set{S}) \text{ and } \mat{R} = \mat{B}(\set{T},:).
\end{equation*}
The matrix $\mat{B}$ is approximated by a low-rank approximation spanned by a subset of both its rows and columns.\CURindex{!reasons for computing|(}
The \CUR decomposition is valued for its increased \emph{interpretability} over other types of low-rank approximations like the SVD \cite{MD09a}.
The \CUR decomposition is also favored in certain contexts because the $\mat{C}$ and $\mat{R}$ matrices inherit structural properties like sparsity from the matrix $\mat{B}$ or because, in some cases, \CUR decompositions can be constructed without looking at the entire input matrix (see \cref{rem:sublinear}).\CURindex{!reasons for computing|)}
This section reviews the theory of \CUR decompositions and discussed algorithms for computing them based on random pivoting.

\myparagraph{Sources}
This chapter presents original research that has not previously been published.
The references \cite{MD09a,PN24} have been helpful in shaping this chapter.

\myparagraph{Outline}
The \CUR decomposition does not fit into our existing taxonomy of randomized low-rank approximation algorithms from \cref{ch:lra}, and it is not evident how we should pick the middle factor $\mat{U}$.
To remedy this deficit, we shall begin by briefly reviewing two additional types of low-rank approximation: two-sided projection approximations (\cref{sec:two-sided}) and generalized Nystr\"om approximations (\cref{sec:gen-nys}).
These two classes yield two \warn{inequivalent} types of \CUR approximations, which we will refer to as \CUR projection approximations and \CUR cross approximations.
\Cref{sec:cur-numerics} discusses numerically stable representation of \CUR approximations, and \cref{sec:cur-algs} discusses random pivoting algorithms for computing \CUR decompositions.
Finally, \cref{sec:md-cur} reviews the Mahoney--Drineas algorithm for computing \CUR decompositions and \cref{sec:cur-experiments} presents a numerical comparison.

\index{two-sided projection approximation|(}
\section{Two-sided and \CUR projection approximation} \label{sec:two-sided}

Given a test matrix $\mat{\Omega} \in \field^{n\times k}$, the projection approximation is the Frobenius-norm optimal approximation to a matrix $\mat{B} \in \field^{m\times n}$ in the column span of $\mat{B}\mat{\Omega}$.
Given a second test matrix $\mat{\Psi} \in \field^{m\times \ell}$, it is natural to define a \emph{two-sided} projection approximation as the Frobenius-norm optimal approximation to $\mat{B}$ that both lies in the column span of $\mat{B}\mat{\Omega}$ and lies in the row span of $\mat{\Psi}^*\mat{B}$.\index{two-sided projection approximation!optimality|(}
To this end, we have the following result \cite[\S4]{Ste99}.

\begin{proposition}[Optimal two-sided approximation] \label{prop:optimal-two-sided}
    Fix matrices $\mat{B} \in \field^{m\times n}$, $\mat{Y} \in \field^{m\times k}$, and $\mat{Z} \in \field^{n\times \ell}$.
    Then the approximation
    \begin{equation*}
        \Bhat \coloneqq \mat{\Pi}_{\mat{Y}} \mat{B}\mat{\Pi}_{\mat{Z}} = \mat{Y}\mat{U}\mat{Z}^* \quad \text{for } \mat{U} \coloneqq \mat{Y}^\dagger \mat{B}\mat{Z}^{\dagger *}
    \end{equation*}
    is the unique Frobenius-norm optimal approximation to $\mat{B}$ whose columns are spanned by the columns of $\mat{Y}$ and whose rows are spanned by the rows of $\mat{Z}^*$.
\end{proposition}

\begin{proof}
    Consider the space of all approximations
    \begin{equation*}
        \set{S} \coloneqq \{ \mat{Y} \mat{T} \mat{Z}^* : \mat{T} \in \field^{k\times \ell} \}.
    \end{equation*}
    Geometrically, we are interested in characterizing the unique point $\mat{B}_\star \in \set{S}$ that is closest to $\mat{B}$.
    The vectorization operation $\vecop : (\field^{m\times n},\norm{\cdot}_{\mathrm{F}}) \to (\field^{mn},\norm{\cdot})$ is a linear isometry, so we can equivalently find the point $\vecop(\mat{B}_\star)$ closest to $\vecop(\mat{B})$ in
    \begin{equation*}
        \vecop(\set{S}) = \{ \vecop(\mat{Y} \mat{T} \mat{Z}^*) : \mat{T}\in\field^{k\times \ell} \} = \{ (\mat{\overline{Z}} \otimes \mat{Y}) \vec{t} : \vec{t} = \vecop(\mat{T}) \in \field^{k\ell} \}.
    \end{equation*}
    The second equality is the vec--Kronecker product identity.
    Thus, we have
    \begin{align*}
        \vecop(\mat{B}_\star) &= \mat{\Pi}_{\mat{\overline{Z}} \otimes \mat{Y}}(\vecop(\mat{B})) = (\mat{\overline{Z}} \otimes \mat{Y})(\mat{\overline{Z}} \otimes \mat{Y})^\dagger \vecop(\mat{B}) \\
        &= (\mat{\overline{Z}} \otimes \mat{Y})(\mat{\overline{Z}}^\dagger \otimes \mat{Y}^\dagger) \vecop(\mat{B}) = \vecop(\mat{Y}\mat{Y}^\dagger\mat{B}\mat{Z}^{\dagger *}\mat{Z}^*).
    \end{align*}
    In the penultimate identity, we used the fact that the pseudoinverse of a Kronecker product is the Kroncker product of the pseudoinverses, which follows by the SVD and the mixed product property for the Kronecker product.
    In the last equality, we applied the vec--Kronecker product identity twice.
    Ergo, $\mat{B}_\star = \mat{Y}\mat{Y}^\dagger\mat{B}\mat{Z}^{\dagger *}\mat{Z}^*$, as promised.
\end{proof}\index{two-sided projection approximation!optimality|)}

This result motivates the following definition:

\begin{definition}[Two-sided projection approximation and \CUR projection approximation]
    Given test matrices $\mat{\Omega} \in \field^{n\times k}$ and $\mat{\Psi} \in \field^{m\times \ell}$, the \emph{two-sided projection approximation} to $\mat{B} \in \field^{m\times n}$ is
    \begin{equation*}
        \Bhat \coloneqq \mat{\Pi}_{\mat{B}\mat{\Omega}}\mat{B}\mat{\Pi}_{\mat{B}^*\mat{\Psi}}.
    \end{equation*}
    If $\mat{\Omega} = \Id(:,\set{S})$ and $\mat{\Psi} = \Id(:,\set{T})$ are column submatrices of the identity matrix, the resulting then the two-sided projection approximation is called a \emph{CUR projection approximation}.\CURindex{!CUR projection approximation@\textsf{CUR} projection approximation}
    It takes the form
    \begin{equation*}
        \Bhat = \mat{C}\mat{U}\mat{R} \quad \text{where $\mat{C} = \mat{B}(:,\set{S})$, $\mat{R} \coloneqq \mat{B}(\set{T},:)$, and $\mat{U} = \mat{C}^\dagger \mat{B} \mat{R}^\dagger$}.
    \end{equation*}
\end{definition}

This definition provides the first natural choice $\mat{U} = \mat{C}^\dagger \mat{B}\mat{R}^\dagger$ for the core matrix in a \CUR decomposition.
As a consequence of \cref{prop:optimal-two-sided}, this type of approximation yields the smallest Frobenius norm error for any middle factor in a \CUR approximation.
In view of this property, Park and Nakatsukasa \cite{PN24} call this type of \CUR approximation a \emph{\CUR best approximation}.\CURindex{!CUR projection approximation@\textsf{CUR} projection approximation}
For this thesis, we will use the term \emph{\CUR projection approximation} to emphasize the connection with (two-sided) projection approximations.
It is important to emphasize that a \CUR projection approximation achieves the lowest possible Frobenius norm error \warn{for a given choice of $\mat{C}$ and $\mat{R}$}; with a bad choice of columns and rows, the accuracy of this ``best approximation'' can be quite poor.

\begin{remark}[Decomposition or approximation]
    In the literature, the term ``decomposition'' is often used to describe the factored approximation $\Bhat = \mat{C}\mat{U}\mat{R} \approx \mat{B}$.
    I prefer the term \CUR approximation, as one typically has $\Bhat\ne \mat{B}$ in applications.
    Conditions for exactness of the \CUR approximation are provided in \cite{HH20}.
\end{remark}

Before going forward, we catalog the following useful and standard (for example, see \cite{MD09a,SE16,DM23a,CK24}) result:

\begin{proposition}[One-sided to two-sided] \label{prop:one-to-two}
    Fix matrices $\mat{B} \in \field^{m\times n}$, $\mat{Y} \in \field^{m\times k}$, and $\mat{Z} \in \field^{n\times \ell}$.
    In any unitarily invariant norm $\uinorm{\cdot}$, we have
    \begin{equation*}
        \uinorm{\mat{B} - \mat{\Pi}_{\mat{B}\mat{\Omega}}\mat{B}\mat{\Pi}_{\mat{B}^*\mat{\Psi}}} \le \uinorm{\mat{B} - \mat{\Pi}_{\mat{B}\mat{\Omega}}\mat{B}} + \uinorm{\mat{B} - \mat{B}\mat{\Pi}_{\mat{B}^*\mat{\Psi}}}. 
    \end{equation*}
\end{proposition}

\begin{proof}
    Apply the triangle inequality:
    \begin{equation*}
        \uinorm{\mat{B} - \mat{\Pi}_{\mat{B}\mat{\Omega}}\mat{B}\mat{\Pi}_{\mat{B}^*\mat{\Psi}}} \le \uinorm{\mat{B} - \mat{\Pi}_{\mat{B}\mat{\Omega}}\mat{B}} + \uinorm{\mat{\Pi}_{\mat{B}\mat{\Omega}}(\mat{B} - \mat{B}\mat{\Pi}_{\mat{B}^*\mat{\Psi}})}.
    \end{equation*}
    By the operator ideal property $\uinorm{\mat{F}\mat{G}} \le \norm{\mat{F}}\cdot \uinorm{\mat{G}}$, the second term may be bounded as
    \begin{equation*}
        \uinorm{\mat{\Pi}_{\mat{B}\mat{\Omega}}(\mat{B} - \mat{B}\mat{\Pi}_{\mat{B}^*\mat{\Psi}})} \le \norm{\mat{\Pi}_{\mat{B}\mat{\Omega}}} \cdot \uinorm{\mat{B} - \mat{B}\mat{\Pi}_{\mat{B}^*\mat{\Psi}}} = \uinorm{\mat{B} - \mat{B}\mat{\Pi}_{\mat{B}^*\mat{\Psi}}}.
    \end{equation*}
    This completes the proof.
\end{proof}
\index{two-sided projection approximation|)}

\index{generalized Nystr\"om approximation|(}
\section{Generalized Nystr\"om approximation and \CUR cross approximation} \label{sec:gen-nys}

The generalized Nystr\"om approximation is a powerful low-rank approximation format that encompasses both projection approximations and Nystr\"om approximations as special cases.
They also lead to another useful, inequalivant way of constructing the middle factor for \CUR approximations.

\begin{definition}[Generalized Nystr\"om approximation and \CUR cross approximation]
    Given test matrices $\mat{\Omega} \in \field^{n\times k}$ and $\mat{\Psi} \in \field^{m\times \ell}$, the \emph{generalized Nystr\"om approximation} to $\mat{B} \in \field^{m\times n}$ is
    \begin{equation*}
        \mat{B}\langle \mat{\Omega},\mat{\Psi}\rangle \coloneqq (\mat{B}\mat{\Omega})(\mat{\Psi}^*\mat{B}\mat{\Omega})^\dagger (\mat{\Psi}^*\mat{B}).
    \end{equation*}
    If $\mat{\Omega} = \Id(:,\set{S})$ and $\mat{\Psi} = \Id(:,\set{T})$ are column submatrices of the identity matrix, then the generalized Nystr\"om approximation is called a \emph{\CUR cross approximation}.\CURindex{!CUR cross approximation@\textsf{CUR} cross approximation}
    It takes the form
    \begin{equation*}
        \mat{B}\langle \set{S},\set{T} \rangle \coloneqq \mat{C}\mat{U}\mat{R} \quad \text{where $\mat{C} = \mat{B}(:,\set{S})$, $\mat{R} \coloneqq \mat{B}(\set{T},:)$, and $\mat{U} = \mat{B}(\set{T},\set{S})^\dagger$}.
    \end{equation*}
\end{definition}

Generalized Nystr\"om approximation includes projection approximation and classical Nystr\"om approximation as the special cases $\mat{\Psi} = \Id$ and $\mat{\Psi} = \mat{\Omega}$.

\index{generalized Nystr\"om approximation!interpolatory property|(}
The generalized Nystr\"om approximation is \emph{interpolatory}.
Consider the case when $k \le \ell$ and $\mat{\Psi}^*\mat{B}\mat{\Omega}$ is full-rank.
Then
\begin{equation*}
    \mat{B}\langle \mat{\Omega},\mat{\Psi}\rangle \cdot \mat{\Omega} = \mat{B}\mat{\Omega}(\mat{\Psi}^*\mat{B}\mat{\Omega})^\dagger (\mat{\Psi}^*\mat{B}\mat{\Omega}) = \mat{B}\mat{\Omega}.
\end{equation*}
Consequently, if $|\set{S}| \le |\set{T}|$ and $\mat{B}(\set{S},\set{T})$ has full rank, then $\Bhat \coloneqq \mat{B}\langle \set{S},\set{T}\rangle$ agrees with the matrix $\mat{B}$ in the selected columns: $\mat{B}(:,\set{S}) = \Bhat(:,\set{S})$.
Analogous statements hold in the opposite case ($k\ge \ell$).\index{generalized Nystr\"om approximation!interpolatory property|)}

\CUR projection approximation\CURindex{!CUR projection approximation@\textsf{CUR} projection approximation} $\mat{U} = \mat{C}^\dagger \mat{B}\mat{R}^\dagger$ and \CUR cross approximation\CURindex{!CUR cross approximation@\textsf{CUR} cross approximation} $\mat{U} = \mat{B}(\set{S},\set{T})^\dagger$ are both natural.
The former enjoys the optimality result \cref{prop:optimal-two-sided}, and the latter satisfies the interpolatory property just discussed.
The \CUR cross approximation is particularly attractive for approximating very large matrices, owing to the following property: \warn{Once the index sets $\set{S}$ and $\set{T}$ have been identified}, forming the \CUR cross approximation only requires accessing $mk+n\ell$ entries of $\mat{B}$.
Thus, the \CUR cross approximation forms the basis for low-rank approximation algorithms for general matrices that avoid reading the entire matrix (\cref{rem:sublinear}).

\begin{figure}
    \centering
    \includegraphics[width=0.7\linewidth]{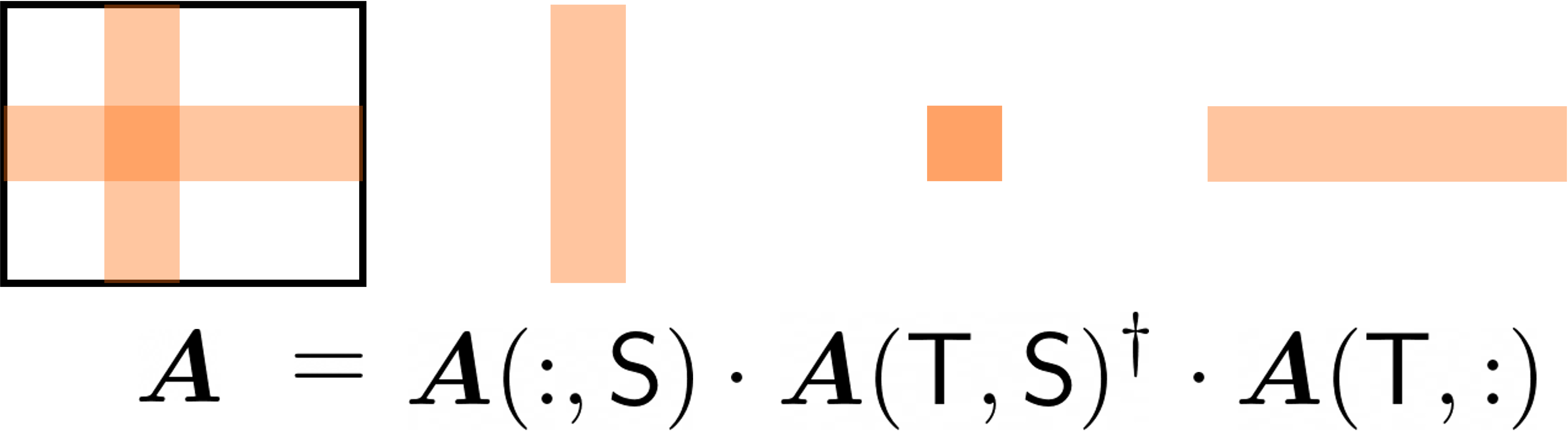}
    \caption{Diagram of cross approximation, showing how the entries selected to form the approximation form a cross.}
    \label{fig:cross}
\end{figure}

\index{generalized Nystr\"om approximation!history|(}
\begin{remark}[History and terminology]
    The generalized Nystr\"om approximation has an interesting history.
    The first references in the modern randomized linear algebra literature are the seemingly independent works of Woolfe, Liberty, Rokhlin, and Tygert \cite{WLRT08} and of Clarkson and Woodruff \cite{CW09}.
    An algebraically equivalent type of approximation appears in \cite{TYUC17}, and they suggest an implementation with improved stability properties.
    The name \emph{generalized Nystr\"om approximation} was suggested by Yuji Nakatsukasa \cite{Nak20}, who investigated different stable implementations that zero out small singular values in the core matrix $\mat{\Psi}^*\mat{B}\mat{\Omega}$.
    Per-Gunnar Martinsson and Alex Townsend trace the origin of this type of approximation all the way back to the work of Wedderburn \cite[p.~69]{Wed34} in 1934, whose writings contain the formula for a rank $k=\ell=1$ generalized Nystr\"om approximation.
    The generalization to $k=\ell > 1$ appears 45 years later in \cite{CF79}; see the paper \cite{CFG95} for more discussion.
    An important difference between the premodern usage of generalized Nystr\"om approximation and its contemporary usage are the settings of the parameters $k$ and $\ell$; the early literature exclusively focused on the case $k = \ell$, whereas the modern literature has exposed benefits of oversampling (e.g., $\ell = \lceil1.5k\rceil$ or $k = \lceil1.5\ell\rceil$).
    
    The name \CUR cross approximation\CURindex{!CUR cross approximation@\textsf{CUR} cross approximation} is taken from \cite{PN24}.
    The name \emph{cross approximation} for the \CUR cross approximation is classical \cite{BR03}, particularly when the row and column subset have the same size $|\set{S}| = |\set{T}|$.
    The name refers to the fact that the decomposition can be computed from only a subset of columns and rows of the matrix.
    If these subsets are contiguous, they form the shape of a cross on the matrix \cite{GT01}; see \cref{fig:cross} for illustration.
    Earlier literature also uses the term \emph{pseudoskeleton approximations} \cite{GTZ97,GZT97,GT01}.
\end{remark}
\index{generalized Nystr\"om approximation!history|)}

\CURindex{!weighting|(}\CURindex{!CUR cross approximation@\textsf{CUR} cross approximation|(}
\myparagraph{Weighted \CUR cross approximations}
One can also consider weighted versions of the \CUR cross approximation.

\begin{definition}[Weighted \CUR cross approximation]
    Let $\set{S} \subseteq \{1,\ldots,n\}$ be a set of $k$ indices and $\set{T} \subseteq \{1,\ldots,m\}$ be a subset of $\ell$ indices, and let $\mat{W}_1 \in \field^{k\times k}$ and $\mat{W}_2 \in \field^{\ell\times \ell}$ be positive \warn{definite}.
    The weighted \CUR cross approximation is 
    \begin{align*}
        \Bhat &\coloneqq \mat{B}(:,\set{S})\mat{W}_1 (\mat{W}_2\mat{B}(\set{T},\set{S})\mat{W}_1)^\dagger \mat{W}_2\mat{B}(\set{T},:)\\
        &=\mat{B}(:,\set{S})\mat{U}\mat{B}(\set{T},:) &&\text{for } \mat{U} = \mat{W}_1 (\mat{W}_2\mat{B}(\set{T},\set{S})\mat{W}_1)^\dagger \mat{W}_2.
    \end{align*}
\end{definition}

If the matrix $\mat{B}(\set{T},\set{S})$ has full rank, then only one of the weight matrices plays a role in the approximation.
Indeed, if $\ell \ge k$ and $\mat{B}(\set{T},\set{S})$ has full column-rank, then 
\begin{equation*}
    (\mat{W}_2\mat{B}(\set{T},\set{S})\mat{W}_1)^\dagger = \mat{W}_1^{-1} \cdot (\mat{W}_2\mat{B}(\set{T},\set{S}))^\dagger,
\end{equation*}
and the core matrix only depends on $\mat{W}_2$:
\begin{equation*}
    \mat{U} = (\mat{W}_2\mat{B}(\set{T},\set{S}))^\dagger \mat{W}_2.
\end{equation*}
A similar result holds if $k\ge \ell$.
In particular, if $k=\ell$ and $\mat{B}(\set{T},\set{S})$ is full-rank, then $\Bhat$ depends on neither weight matrix, and $\Bhat$ is just the ordinary \CUR cross approximation.

The justification for considering weighted \CUR approximations is provided by the following result.

\begin{theorem}[The value of weighting for \CUR approximations] \label{thm:weighted-cur}
    Fix $k\ge 1$ and any $m$ that is divisible by $k+1$.
    There exists a matrix $\mat{B} \in \real^{m\times 2k}$ and a column subset $\set{S} = \{1,\ldots,2k\}$ of $k$ elements such that the squared \warn{unweighted} \CUR cross approximation error is at least
    \begin{equation*}
        \norm{\mat{B} - \mat{B}(:,\set{S})\mat{B}(\set{T},\set{S})^\dagger\mat{B}(\set{T},:)}_{\mathrm{F}}^2 \ge 1.5 \norm{\mat{B} - \mat{B}(:,\set{S})\mat{B}(:,\set{S})^\dagger\mat{B}(\set{T},:)}_{\mathrm{F}}^2
    \end{equation*}
    unless the row subset $\set{T}$ comprises at least
    \begin{equation*}
        \ell \ge (3-2\sqrt{2})m \approx 0.172\cdot m \text{ elements}.
    \end{equation*}
    Conversely, for (diagonally) \warn{weighted} \CUR cross approximations, a row subset of size
    \begin{equation} \label{eq:cur-upper-bound}
        \ell = \order(k \log k + k/\varepsilon)
    \end{equation}
    is sufficient to obtain the guarantee
    \begin{equation*}
        \norm{\mat{B} - \mat{B}(:,\set{S})(\mat{W}_2\mat{B}(\set{T},\set{S}))^\dagger\mat{W}_2\mat{B}(\set{T},:)}_{\mathrm{F}}^2 \le(1+\varepsilon) \norm{\mat{B} - \mat{B}(:,\set{S})\mat{B}(:,\set{S})^\dagger\mat{B}(\set{T},:)}_{\mathrm{F}}^2.
    \end{equation*}
\end{theorem}

The proof is somewhat lengthy, so we defer it to \cref{sec:weighted-cur-proof}.
I believe that the dependence of $\ell$in \cref{eq:cur-upper-bound} can be reduced to $\ell = \order(k/\varepsilon)$ by extending the results of \cite{CP19} to least-squares problems with multiple right-hand sides.

This result establishes a fundamental separation between the approximation capabilities of weighted and unweighted \CUR approximations.
In short, the conclusion of \cref{thm:weighted-cur} may be summarized as follows.
\actionbox{Weighted \CUR cross approximations achieve $(1+\varepsilon)$-type approximation guarantees using a row-subset of size $\ell = \order(k \log k + k/\varepsilon)$, whereas unweighted approximations a constant fraction of the matrix rows $\ell = \Omega(m)$ to achieve such a guarantee.}
In practice, the examples which demonstrate a separation between the weighted and unweighted approximations are somewhat pathological, and unweighted \CUR cross approximations are usually fine for applications.\index{generalized Nystr\"om approximation|)}\CURindex{!weighting|)}\CURindex{!CUR cross approximation@\textsf{CUR} cross approximation|)}

\CURindex{!stable representation|(}\index{interpolatie decomposition!numerical stability|(}
\section{Numerically stable representations} \label{sec:cur-numerics}


Standard implementations of \CUR approximations suffer from issues of numerical stability, as Martinsson and Tropp explain \cite[\S13.1]{MT20a}:
\begin{quote}
    A disadvantage of the [\CUR cross approximation] is that, when the singular values of $\mat{B}$ decay rapidly, the factorization [$\mat{B}\langle \set{S},\set{T} \rangle = \mat{C}\mat{U}\mat{R}$] is typically numerically ill-conditioned.
    The reason is that, whenever the factorization is a good representation of $\mat{B}$, the singular values of $\mat{B}(\set{T},\set{S})$ should approximate the $k$ dominant singular values of $\mat{B}$, so the singular values of $\mat{U}$ end up approximating the inverses of these singular values. This means that $\mat{U}$ will have elements of magnitude $1/\sigma_k$, which is clearly undesirable when $\sigma_k$ is small. In contrast, the ID\index{interpolative decomposition} is numerically benign.
\end{quote}
(We have modified this quote to have consistent notation with the rest of this chapter.)
The issues identified in this quote show that working with the \warn{standard representation} $\Bhat = \mat{C}\mat{U}\mat{R}$ of the \CUR approximation can be numerically problematic.\index{interpolatie decomposition!numerical stability|)}

To fix this instability, the authors of \cite{ADM+15} propose representing the \CUR projection approximation using the factorization
\begin{equation} \label{eq:previous-stable-cur}
    \Bhat = \mat{Q}_1\mat{\Xi}\mat{Q}_2^*,
\end{equation}
where $\mat{Q}_1 \coloneqq \orth(\mat{B}(:,\set{S}))$, $\mat{Q}_2 \coloneqq \orth(\mat{B}(\set{T},:)^*)$, and $\mat{\Xi} \coloneqq \mat{Q}_1^*\mat{B}\mat{Q}_2$.
While this format is stable, it no longer represents the low-rank approximation as a \CUR factorization $\Bhat = \mat{C}\mat{U}\mat{R}$.
In particular, the factors $\mat{Q}_1$ and $\mat{Q}_2$ no longer inherit structural properties of $\mat{B}$ such as sparsity, and $\mat{Q}_1$ and $\mat{Q}_2$ lack the interpretability benefits of approximating $\mat{B}$ using columns and rows.

To remedy the numerical stability issues with \CUR approximation while maintaining a factored representation using selected columns and rows, I propose the following new way to stably represent the \CUR approximation.
The observation is that it can be more numerically stable to work with an \emph{implicit representation} of the core matrix $\mat{U}$.
Specifically, we represent $\mat{U}$ as a product
\begin{equation} \label{eq:cur-U-factored}
    \mat{U} = \mat{P}^{-1} \cdot \mat{G}
\end{equation}
of a \emph{well-conditioned} matrix $\mat{G}$ and a (possibly ill-conditioned) upper-triangular matrix $\mat{P}$.
This representation can be computed in the following ways:
\begin{enumerate}
    \item For \CUR projection approximations, first use a stable algorithm to compute a \warn{row} interpolative decomposition\index{interpolative decomposition}
    \begin{equation*}
        \mat{B} \approx \mat{W}\mat{B}(\set{T},:) \quad \text{for } \mat{W} = \mat{B} \mat{B}(\set{T},:)^\dagger .
    \end{equation*}
    To compute a stable representation \cref{eq:cur-U-factored}, we compute a \QR decomposition of the column submatrix $\mat{B}(:,\set{S})$:
    \begin{equation*}
        \mat{B}(:,\set{S}) = \mat{Q}\mat{P}
    \end{equation*}
    and define $\mat{G} \coloneqq \mat{Q}^*\mat{W}$.
    The core matrix $\mat{U}$ admits the factorization \cref{eq:cur-U-factored}.
    \item For (weighted) \CUR cross approximations with $\ell \ge k$ and weight matrix $\mat{W}_2$, we compute a \QR decomposition
    \begin{equation*}
        \mat{W}_2 \mat{B}(\set{T},\set{S}) = \mat{Q}\mat{P}.
    \end{equation*}
    and define $\mat{G} \coloneqq \mat{Q}^*\mat{W}_2$.
    The core matrix $\mat{U}$ admits the factorization \cref{eq:cur-U-factored}.
\end{enumerate}
As we will see in \cref{sec:cur-experiments}, working with $\mat{G}$ and $\mat{P}$ matrices, rather than the standard core matrix $\mat{U}$, can lead to significant stability improvements to the \CUR approximation while maintaining the core essence of the decomposition $\mat{B} \approx \mat{C}\mat{U}\mat{R}$.
We emphasize that, to reap the numerical stability benefits of this approach, it is critical to store the triangular matrix $\mat{P}$, \warn{not its inverse $\mat{P}^{-1}$}.

\myparagraph{Motivation for improved stability}
To motivate why the representation \cref{eq:cur-U-factored} is more stable, first consider a square matrix $\mat{B} \in \field^{n\times n}$, and work with the extreme case when $\set{S} = \set{T} = \{1,\ldots,n\}$.
The standard \CUR decomposition $\mat{B} = \mat{C}\mat{U}\mat{R}$ has factors $\mat{C} = \mat{R} = \mat{B}$ and $\mat{U} = \mat{B}^{-1}$; in particular, forming $\mat{U}$ as a dense array requires explicit computation of the inverse $\mat{B}^{-1}$.
As can be seen in the following MATLAB code segment, the reconstruction error with this \CUR decomposition can be large:
\begin{lstlisting}
B = rand_with_evals(logspace(0,-15,1000)); % Condition num = 1e15
norm(B - B*inv(B)*B)/norm(B)               % Outputs 9.7e-3\end{lstlisting}
If we instead compute a \QR decomposition $\mat{B} = \mat{Q}\mat{P}$ and represent the inverse as $\mat{B}^{-1} = \mat{P}^{-1}\mat{Q}^*$, the reconstruction error drops to the level of the machine precision:
\begin{lstlisting}
[Q,P] = qr(B);
norm(B - B*(P\(Q'*B))/norm(B)             % Outputs 1.6e-15
\end{lstlisting}
The stability difference between these two approaches is consequence of the fact that, as a way of solving the system $\mat{B}\vec{X} = \vec{Y}$, the procedure \texttt{inv(B)*y} is \warn{forward but not backward stable}\index{numerical stability!backward stability}\index{numerical stability!foward stability} \cite{DT12}, whereas \texttt{P\textbackslash{}(Q'*B)} is backward stable.
Notions of forward and backward stability are discussed a great deal in \cref{part:sketching} of this thesis.

The above example demonstrates the ineffectiveness of explicit inversion in storing and representing the core matrix $\mat{U}$ in a \CUR decomposition in extreme case when $\set{S}= \set{T} = \{1,\ldots,n\}$.
The same principle---that explicit computation of the inverse of a matrix is bad---motivates the use of the implicitly represented core matrix $\mat{U} = \mat{P}^{-1} \cdot \mat{G}$ for storing general \CUR approximations.

\myparagraph{Using the stable representation}
The stable representation \cref{eq:cur-U-factored} of the core matrix allows the \CUR decomposition to be used in a numerically stable way for downstream tasks.
The approximation $\Bhat$ can be stably reconstituted as a dense array by the right--left evaluation $\mat{C}(\mat{P}^{-1}(\mat{G}\mat{R}))$ or the left--right evaluation $((\mat{C}\mat{P}^{-1})\mat{G})\mat{R}$; the mixed evaluation order $\mat{C}((\mat{P}^{-1}\mat{G})\mat{R})$ should be avoided, however.
The left--right and right--left evaluation orders can also be used to efficiently and stably generate submatrices or individual entries of $\Bhat$.
One can also compute matrix products as 
\begin{equation*}
    \Bhat\cdot \mat{X} = \mat{C}(\mat{P}^{-1}(\mat{G}(\mat{R}\mat{X}))) \quad \text{or} \quad \mat{Y}^*\Bhat = (((\mat{Y}^*\mat{C})\mat{P}^{-1})\mat{G})\mat{R}.
\end{equation*}

Storing a \CUR decomposition with stable representation \cref{eq:cur-U-factored} is economical,
as it only requires retaining the two submatrices $\mat{C} = \mat{B}(:,\set{S})$ and $\mat{R} = \mat{B}(\set{T},:)$, plus $\order(k\ell)$ storage for the factored core matrix \cref{eq:cur-U-factored}.
If the matrix is structured, storing $\mat{C}$ and $\mat{R}$ may be significantly more efficient then storing general matrices $\mat{Q}_1$ and $\mat{Q}_2$ in the representation \cref{eq:previous-stable-cur}.
In particular, when working with a function matrix $\mat{B} = \xi(\set{X},\set{Y})$ (\cref{def:function_matrix}), the entries of $\mat{C} = \mat{B}(:,\set{S})$ and $\mat{R} = \mat{B}(\set{T},:)$ need not ever be stored all at once, as they can be generated on-the-fly using the data $\set{X}$ and $\set{Y}$ and the function $\xi$.

\begin{remark}[Comparison with \cite{PN24}]
    Park and Nakatsukasa \cite{PN24} also consider the numerical stability of \CUR approximation, focusing on unweighted \CUR cross approximations.
    These authors are principally concerned with developing stable algorithms for \emph{reconstructing} a \CUR cross approximation as a dense array, which differs from our emphasis on finding a stable way to \emph{store} a \CUR approximation while maintaining the traditional $\mat{C}\mat{U}\mat{R}$ factorization.
    
    To convert the \CUR approximation $\Bhat = \mat{C}\mat{U}\mat{R}$ to a dense matrix, Park and Nakatsukasa compute an (economy-size) SVD $\mat{B}(\set{T},\set{S}) = \mat{W}\mat{\Sigma}\mat{V}^*$ and reconstruction $\Bhat$ using the formula $\Bhat = (\mat{C}\mat{V}\mat{\Sigma}^{-1})(\mat{W}^*\mat{R})$.
    (They also mention the possibility of using a \QR decomposition in place of the SVD, which agrees with our approach.)
    Park and Nakatsukasa provide a formal proof of stability for their approach, and discuss a variant where the pseudoinverse $\mat{B}(\set{T},\set{S})^\dagger$ is regularized by setting the small singular values of $\mat{B}(\set{T},\set{S})$ to zero.

    This work extends the insights of Nakatsuksa and Park by proposing the $\mat{C}(\mat{P}^{-1}\mat{G})\mat{R}$ format as a stable way of representing a \CUR decomposition of \emph{any} type.
    In particular, we show how to use this format to stably represent weighted \CUR cross approximations and \CUR projection approximations, which is beyond the scope of \cite{PN24}.
    My hope is that this new work makes clear that, when appropriate care for numerical stability is taken, the \CUR approximation \emph{storage format} is safe for use in general-purpose computations, resolving the issues identified by Martinsson and Tropp in the beginning of this section.
\end{remark}
\CURindex{!stable representation|)} 

\CURindex{!algorithms|(}

\section{Algorithms} \label{sec:cur-algs}

As \CUR decompositions require identifying subsets $\set{S}$ and $\set{T}$ of the columns and rows of the matrix $\mat{B}$ that approximately span the matrix.
As demonstrated throughout this part of the thesis, random pivoting algorithms are well-suited to exactly this task.
There are a wide array of possible strategies for computing \CUR decompositions using the random pivoting approach.
This thesis discusses two possible approaches, which I will call \RPCURtwo and \RPCURLev.
The former approach separately applies \RPQR to $\mat{B}$ and $\mat{B}^*$, and the latter approach applies \RPQR to obtain a low-rank approximation $\mat{B} \approx \mat{Q}\mat{F}^*$, then applies leverage score sampling\index{leverage scores!as a sampling distribution for computing \textsf{CUR} approximations} to $\mat{Q}$.

\subsection{\RPCURtwo: Running \RPQR twice}
Our first strategy for computing a \CUR decomposition is to simply run \RPQR twice, once on $\mat{B}$ to obtain a column set $\set{S}$ and once on $\mat{B}^*$ to obtain a row set $\set{T}$.
One may be then combine these subsets to create a \CUR decomposition of $\mat{B}$.
While \RPCURtwo can be used to create either type of \CUR decomposition, it is most natural when used for computing a \CUR projection approximation.\CURindex{!CUR projection approximation@\textsf{CUR} projection approximation}
Code is provided in \cref{prog:rpcur2}.

\myprogram{Implementation of \RPCURtwo for computing a \CUR projection approximation.}{Subroutine \texttt{rpqr} is provided in \cref{prog:rpqr}.}{rpcur2}

\index{RPCUR2 algorithm@\textsc{RPCUR2} algorithm!implementation|(}
\myparagraph{Implementation}
\RPCURtwo requires two separate executions of \RPQR, plus additional post-processing to construct the core matrix.
The total runtime is $\order(mn(k+\ell))$ operations.
For runtime speed and numerical stability, I use the Householder-based\index{Householder \QR factorization!for pivoted partial \QR decomposition} accelerated \RPQR\index{accelerated randomly pivoted QR@accelerated randomly pivoted \QR} implementation (\cref{prog:acc_rpqr}) in my code.
(I have found examples where \RPCURtwo fails catastrophically when using the modified Gram--Schmidt implementation in \cref{prog:acc_rpqr_bgs}.)\index{RPCUR2 algorithm@\textsc{RPCUR2} algorithm!implementation|)}

\index{RPCUR2 algorithm@\textsc{RPCUR2} algorithm!theoretical results|(}
\myparagraph{Analysis}
\emph{A priori} error bounds for the \RPCURtwo algorithm can easily be inferred from the results for \RPQR (\cref{cor:rpqr-error}) and \cref{prop:one-to-two}.

\begin{corollary}[\RPCURtwo]
    Let $\mat{B} \in \field^{m\times n}$ be a matrix, fix $r\ge 1$, and introduce the squared relative error of the best rank-$r$ approximation:
    \begin{equation*}
        \eta \coloneqq \frac{\norm{\mat{B} - \lowrank{\mat{B}}_r}_{\mathrm{F}}^2}{\norm{\mat{B}}_{\mathrm{F}}^2}.
    \end{equation*}
    \RPCURtwo (\cref{prog:rpcur2}) produces an approximation $\Bhat$ satisfying
    \begin{equation*}
        \expect \norm{\mat{B} - \Bhat}_{\mathrm{F}}^2 \le 4 \norm{\mat{B} - \lowrank{\mat{B}}_r}_{\mathrm{F}}^2,
    \end{equation*}
    provided the parameters $k$ and $\ell$ satisfy $k,\ell \ge r\log(\e/\eta)$.
\end{corollary}\index{RPCUR2 algorithm@\textsc{RPCUR2} algorithm!theoretical results|)}

\index{RPCUR2 algorithm@\textsc{RPCUR2} algorithm!parameter selection|(}
\myparagraph{Automatic rank determination}
If one wishes to compute a \CUR decomposition with error controlled by a tolerance $\tau$ in some \warn{unitarily invariant} norm, one can achieve this goal by running each \RPQR step with a tolerance of $\tau/2$, in view of \cref{prop:one-to-two}.\index{RPCUR2 algorithm@\textsc{RPCUR2} algorithm!parameter selection|)}

\index{leverage scores!as a sampling distribution for computing \textsf{CUR} approximations|(}\index{RPCURLev algorithm@\textsc{RPCURLev} algorithm|(}
\subsection{\RPCURLev: Combining \RPQR and leverage score sampling}

Our second algorithm for computing a \CUR decomposition combines \RPQR to compute the column set $\set{S}$ with leverage score sampling (\cref{def:ridge-general}) to determine the pivot set $\set{T}$.
It most naturally outputs a (diagonally) \warn{weighted} \CUR cross approximation,\CURindex{!CUR cross approximation@\textsf{CUR} cross approximation}\CURindex{weighting} but it can also be used with to compute a \CUR projection approximation or unweighted \CUR cross approximation if desired.

The \RPCURLev algorithm proceeds as follows.
First, run \RPQR to obtain a column projection approximation $\mat{B} \approx \mat{Q}\mat{F}^*$.
Recall that, for the matrix $\mat{Q}$ with orthonormal columns, its leverage scores $\vec{\lambda} = \srn(\mat{Q})$ are its squared row norms.\index{squared row or column norms}
We then form $\set{T}$ by sampling iid from the leverage score distribution
\begin{equation*}
    \set{T} = \{t_1,\ldots,t_\ell\} \quad \text{where } t_1,\ldots,t_\ell \simiid \vec{\lambda}.
\end{equation*}
We emphasize sampling is to be done \warn{with replacement}.
For reasons that will soon become more clear, it is most natural to output a weighted \CUR decomposition with weight matrix $\mat{W}_2 = \Diag(\vec{\lambda}(\set{T}))^{-1/2}$.

\index{RPCURLev algorithm@\textsc{RPCURLev} algorithm!implementation|(}
\myprogram{Implementation of \RPCURLev for computing a weighted \CUR cross approximation to a matrix.}{Subroutines \texttt{rpqr} and \texttt{sqrownorms} are provided in \cref{prog:rpqr,prog:sqrownorms}.}{rpcur_lev}

Code is given in \cref{prog:rpcur_lev}.
We highlight that this implementation computes the pseudoinverse numerically using a column-pivoted \QR decomposition, discarding entries from the index set $\set{S}$ to ensure that $\mat{B}(\set{T},\set{S})$ is numerically full-rank.
The use of column pivoting makes this approach robust, even for sparse problems which could have many zero entries.

\myparagraph{Runtime}
\RPCURLev requires running $k$ steps of \RPQR, sampling $\ell$ elements from a weighted probability distribution on $m$ items, and doing post-processing on an $\order(k\ell)$ matrix.
The runtime of \RPCURLev is $\order(mnk + k^2\ell)$.
I use accelerated \RPQR\index{accelerated randomly pivoted \QR} (\cref{prog:acc_rpqr}) in my implementation.
As we will see in \cref{fig:rpcur_sparse}, the \RPCURLev algorithm can be meaningfully faster than the \RPCURtwo algorithm.\index{RPCURLev algorithm@\textsc{RPCURLev} algorithm!implementation|)}

\index{leverage scores!theoretical results|(}\index{RPCURLev algorithm@\textsc{RPCURLev} algorithm!theoretical results|(}
\myparagraph{Analysis}
To analyze \RPCURLev, we can make use of existing results from the leverage score literature.
We shall use the following result:

\begin{fact}[Leverage score sampling for least squares] \label{fact:lev-ls}
    Let $\mat{C} \in \field^{m\times k}$ be a full-rank matrix, and let $\vec{F} \in \field^{m\times n}$.
    Let $\vec{\lambda} = \srn(\orth(\mat{C}))$ be the leverage scores of $\mat{C}$ and sample an index set $\set{T}$ of $\ell$ elements \warn{iid with replacement} from $\vec{\lambda}$.
    Introduce the diagonal weight matrix $\mat{W} \coloneqq \Diag(\vec{\lambda}(\set{T}))^{-1/2}$.
    Provided 
    \begin{equation*}
        \ell \ge \Omega \left( k \log k + \frac{k}{\varepsilon} \right),
    \end{equation*}
    it holds with at least 99\% probability that
    \begin{equation*}
        \norm{\mat{F} - \mat{C}[\mat{W}\mat{C}(\set{T},:)]^\dagger [\mat{W}\mat{F}(\set{T},:)]}_{\mathrm{F}} \le  (1+\varepsilon)\norm{\mat{F} - \mat{C}\mat{C}^\dagger \mat{F}}_{\mathrm{F}}.
    \end{equation*}
\end{fact}

This result is a variant of \cite[Thm.~7.9]{CW17}, with the version stated here following straightforwardly from the proof of that result.\index{leverage scores!theoretical results|)}
This result leads immediately to a error bound for \RPCURLev.

\begin{corollary}[\RPCURLev] \label{cor:rpcur_lev}
    Let $\mat{B} \in \field^{m\times n}$ be a matrix, fix $r\ge 1$ and $\varepsilon \in (0,1)$, and introduce the squared relative error of the best rank-$r$ approximation:
    \begin{equation*}
        \eta \coloneqq \frac{\norm{\mat{B} - \lowrank{\mat{B}}_r}_{\mathrm{F}}^2}{\norm{\mat{B}}_{\mathrm{F}}^2}.
    \end{equation*}
    With 99\% probability, the weighted \CUR cross approximation $\Bhat$ produced by \RPCURtwo (\cref{prog:rpcur2}) produces an approximation satisfying
    \begin{equation*}
        \expect_{\set{S}} \norm{\mat{B} - \Bhat}_{\mathrm{F}}^2\le (1+\varepsilon) \norm{\mat{B} - \lowrank{\mat{B}}_r}_{\mathrm{F}}^2,
    \end{equation*}
    provided the parameters $k$ and $\ell$ satisfy
    \begin{equation*}
        k \ge \frac{2r}{\varepsilon} + r \log \left( \frac{2}{\eta\varepsilon} \right), \quad \ell \ge \Omega\left(k\log k + \frac{k}{\varepsilon}\right)
    \end{equation*}
    Here, $\expect_{\set{S}}$ denotes the expectation over the randomness in the first index set $\set{S}$.
\end{corollary}

The proof is immediate from \cref{cor:rpqr-error,fact:lev-ls}.
This result suggests (correctly) that both oversampling and weighting $\ell > k$ are necessary for \RPCURLev to produce high-quality approximations in general.\index{RPCURLev algorithm@\textsc{RPCURLev} algorithm!theoretical results|)}

\index{RPCURLev algorithm@\textsc{RPCURLev} algorithm!parameter choices|(}
The \RPCURLev procedure represents only one way of combining \RPQR with an extra step of row sampling to produce a (weighted) \CUR cross approximation.
In my experience, it reliably produces high-quality approximations \warn{when implemented with moderately high oversampling} $\ell = \order(k\log k)$; approaches using pivoting \cite{SE16,PN24,CK24} are more effective for small $\ell\approx k$.
In place of leverage score sampling, natural alternatives are volume sampling \cite{DRVW06,DW17,DW18}, adaptive randomized pivoting \cite{CK24}, leveraged volume sampling \cite{DWH18}, the Chen--Price method \cite{CP19}, and pivoting strategies such as the Park--Nakatsukasa algorithm \cite[Alg.~3.2]{PN24}.

\myparagraph{Parameter settings}
On worst-case examples, the logarithmic oversampling $\ell \ge \Omega(k\log k)$ suggested by \cref{cor:rpcur_lev} is necessary for \RPCURLev to succeed.
We recommend $\ell = \lceil 1.5 k \log k\rceil$ as a sensible default value.
If one wants to control the error up to a tolerance $\tau$ in the Frobenius norm, we recommend running \RPQR with a lower tolerance ($\tau/3$, say) to allow for the additional error produced by selection of the row subset $\set{T}$.\index{RPCURLev algorithm@\textsc{RPCURLev} algorithm|)}\index{RPCURLev algorithm@\textsc{RPCURLev} algorithm!parameter choices|)}

\index{MDCUR algorithm@\textsc{MDCUR} algorithm|(}
\section{Related work: Mahoney and Drineas' algorithm} \label{sec:md-cur}

In their influential work on \CUR decompositions for data analysis, Mahoney and Drineas \cite{MD09a} proposed a pure leverage score sampling-based approach to select index sets $\set{S}$ and $\set{T}$ for \CUR approximations; see also the follow-up work \cite{DMM08} which obtains better theoretical results.
Up to some small tweaks, their algorithm works as follows:
\begin{enumerate}
    \item \textbf{SVD.} Compute an SVD $\mat{B} = \mat{U}\mat{\Sigma}\mat{V}^*$.
    \item \textbf{Column sampling.} Choose a parameter $k'\le k$, and sample $k$ indices $\set{S} = \{s_1,\ldots,s_k\}$ iid from the leverage score distribution of $\mat{V}(:,1:k')$:
    \begin{equation*}
        s_1,\ldots,s_k \simiid \srn(\mat{V}(:,1:k')).
    \end{equation*}
    
    \item \textbf{Row sampling.} Choose a parameter $\ell'\le \ell$, and sample $\ell$ indices $\set{T} = \{t_1,\ldots,t_\ell\}$ iid from the leverage score distribution of $\mat{U}(:,1:\ell')$:
    \begin{equation*}
        t_1,\ldots,t_k \simiid \srn(\mat{U}(:,1:\ell')).
    \end{equation*}
    \item \textbf{Output.} Return the \CUR projection approximation with index sets $\set{S}$ and $\set{T}$.
\end{enumerate}
To achieve a $(r,1+\varepsilon,2)$-approximation, Mahoney and Drineas suggest choosing $k' = \ell' = r$ and $k,\ell = \order((r\log r)/\varepsilon^2)$.
I believe the parameters $k,\ell = \order(r/\varepsilon + r\log r)$ should suffice with modern proof techniques, but I have not confirmed this.
In their original algorithm, Drineas and Mahoney suggest including or not including each element $1\le s\le n$ or $1\le t \le m$ in the sets $\set{S}$ or $\set{T}$ independently with some probability so that the \warn{expected} number of accepted indices is $\expect[k],\expect[\ell] = \order((r\log r)/\varepsilon^2)$; we have modified the algorithm here to pick a fixed number of indices.
In our experiments, we set $k' = \lceil k/(2\log k)\rceil$ and $\ell' = \lceil\ell/(2\log \ell)\rceil$.
We refer to this algorithm as Mahoney--Drineas \CUR (\MDCUR) and provide an implementation in \cref{prog:md_cur}.

\myprogram{Implementation of \MDCUR for computing a \CUR projection approximation.}{}{md_cur}

\myparagraph{Runtime}
The dominant cost of the \MDCUR algorithm is the SVD computation.
Using a dense SVD, its cost is $\order(mn\min\{m,n\})$ operations.
We can accelerate the algorithm by using a randomized SVD, reducing the cost to $\order(mn(k'+\ell'))$.
The remaining steps of the algorithm incur a post-processing cost of $\order(mk^2 + n\ell^2)$.

\myparagraph{Further alternatives}
Another approach to computing \CUR decompositions is to use sketchy pivoting, discussed in \cref{sec:sketchy-pivoting} \cite{DM23}.
The comparison of randomly pivoted \CUR methods to sketchy pivoting methods for \CUR is similar to the comparison of \RPQR methods to sketchy pivoting methods for interpolative decomposition.
In particular, sketchy pivoting methods can be faster but generally require an upper bound on the approximation rank $k$ to be known in advance, which is a significant limitation for some applications.
See \cref{sec:sketchy-pivoting} and \cite{DCMP23} for further discussion on the relative merits of sketchy pivoting and random pivoting.

There are many other algorithms for computing \CUR decompositions including squared row and column norm sampling \cite{DKM06b}, adaptive cross approximation \cite{BR03}, adaptive randomized pivoting \cite{CK24}, and more sophisticated multi-stage sampling approaches \cite{BW17}.
Comparison of the random pivoting algorithms against these approaches is a natural subject for future work.

\CURindex{!algorithms|)}\index{MDCUR algorithm@\textsc{MDCUR} algorithm|)}\index{leverage scores!as a sampling distribution for computing \textsf{CUR} approximations|)}

\CURindex{!numerical results|(}\index{MDCUR algorithm@\textsc{MDCUR} algorithm!numerical results|(}\index{RPCUR2 algorithm@\textsc{RPCUR2} algorithm!numerical results|(}\index{RPCURLev algorithm@\textsc{RPCURLev} algorithm!numerical results|(}

\section{Experiments} \label{sec:cur-experiments}

To begin exploring the performance of the two randomly pivoted \CUR algorithms, we present some preliminary tests on three different example matrices.

\myparagraph{Experiment \#1: Function matrix\index{function matrix}}
First, as a simple test of approximation quality and numerical stability, we test the two randomly pivoted \CUR algorithms of the previous section on a simple test matrix $\mat{B} \in \complex^{n\times n}$ with entries
\begin{equation*}
    b_{ij} = \frac{1}{z_i - w_j}
\end{equation*}
The points $z_i$ and $w_j$ are chosen to be equispaced on the complex unit circle $\unitcircle(\complex)$:
\begin{equation*}
    z_i \coloneqq \exp \left(2\pi \iu \cdot\frac{i-1}{2n}\right), \: w_i \coloneqq \exp \left(2\pi \iu \cdot\frac{n+i-1}{2n}\right) \quad \text{for } i=1,2,\ldots,n,
\end{equation*}
and we set $n\coloneqq 1000$.
The points $z_i$ trace the portion of the unit circle in the upper half plane, $\unitcircle(\complex) \cap \{\Im(z) \ge 0\}$, and the $w_j$ trace the lower half, $\unitcircle(\complex) \cap \{\Im(z) \le 0\}$.
Matrices similar to $\mat{B}$ occur in the design of algorithms for structured systems of linear equations and least-squares problems \cite{CGS+08,XXG12,XXCB14,WEB25,BKW25}.\index{rank-structured matrix}

\begin{figure}[t]
    \centering
    \includegraphics[height=2.4in]{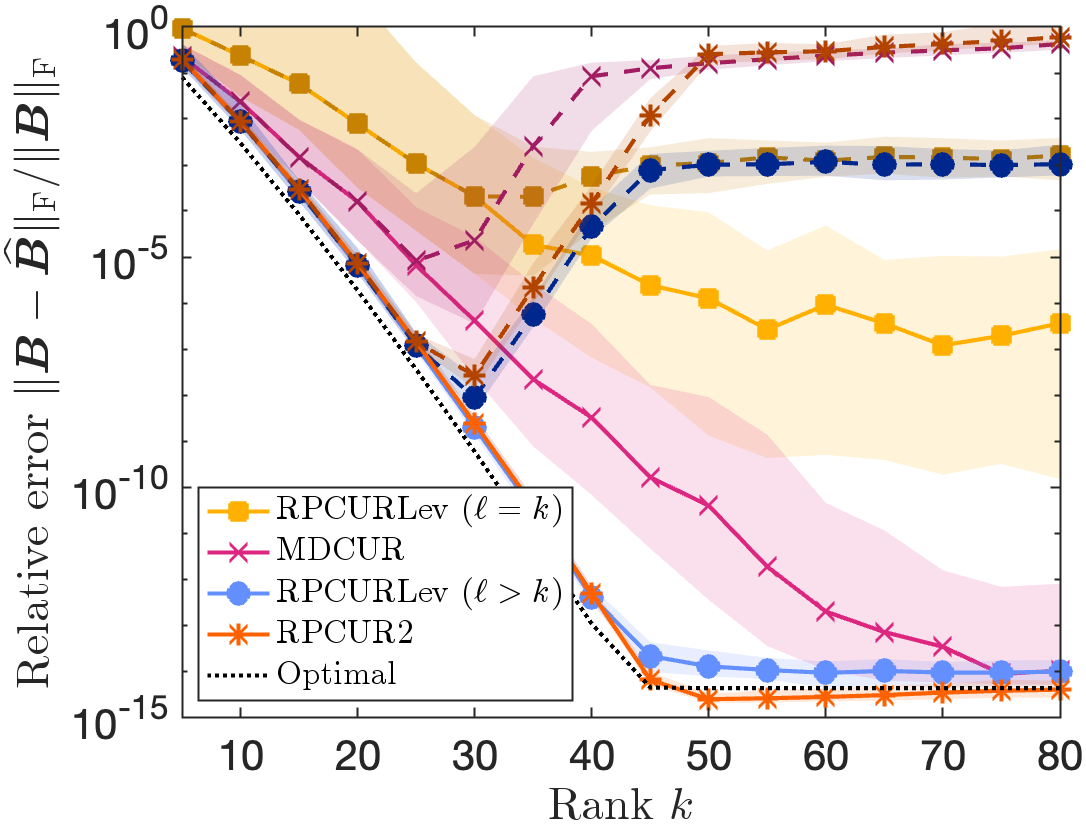}
    \includegraphics[height=2.4in]{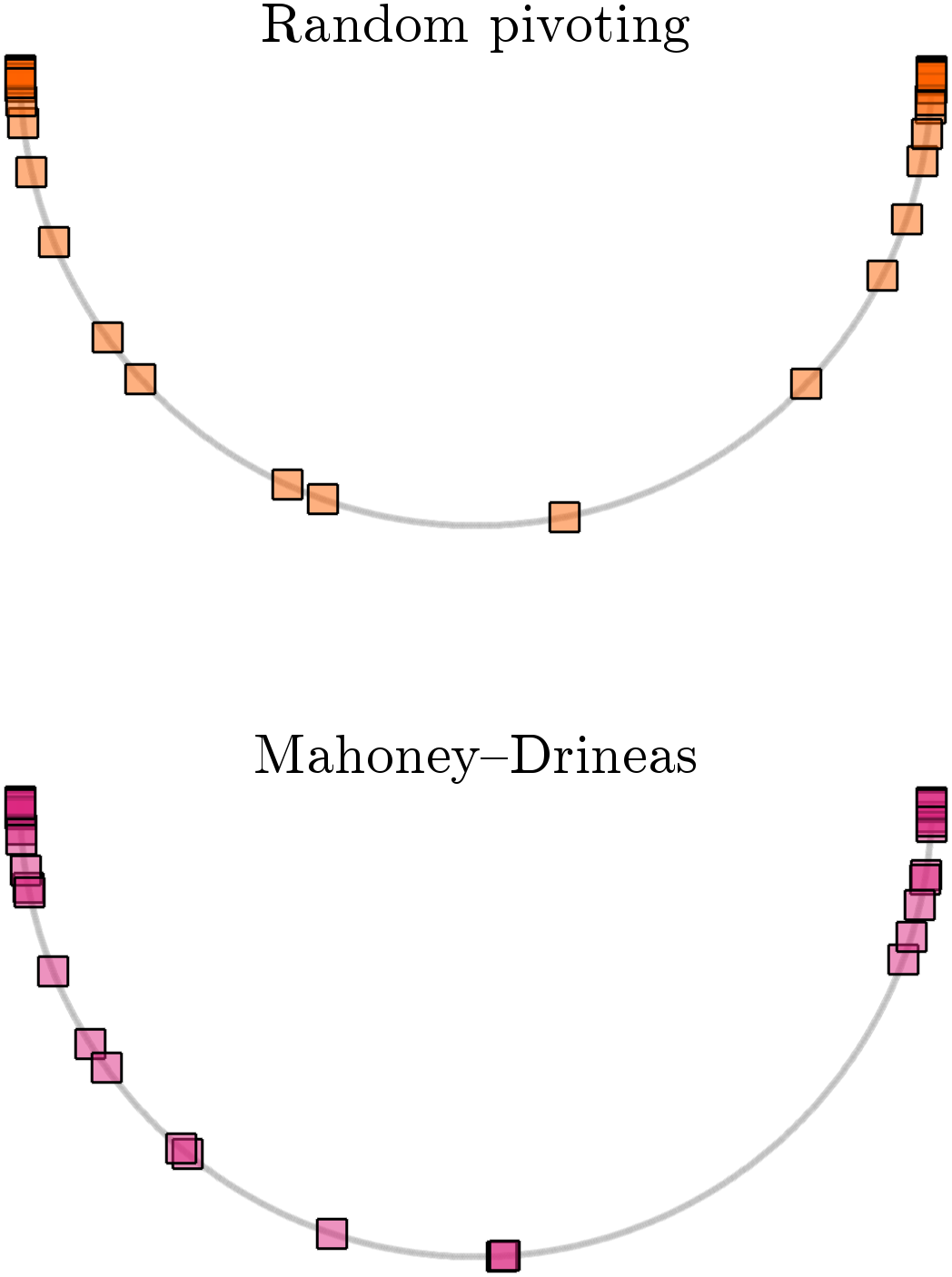}
    \caption[Relative error and chosen landmarks for various \CUR algorithms applied to a function matrix]{\emph{Left:} Relative error $\norm{\mat{B} - \Bhat}_{\mathrm{F}}/\norm{\mat{B}}_{\mathrm{F}}$ computed by the following methods: \RPCURtwo (orange stars), \RPCURLev with (blue circles) and without oversampling (yellow squares), and \MDCUR (pink crosses).
    Solid light lines for \CUR algorithms use the stable representation $\mat{U} = \mat{P}^{-1}\cdot \mat{G}$ from \cref{sec:cur-numerics}, and dashed dark lines store $\mat{U}$ as a dense array.
    The error for the best rank-$k$ approximation $\Bhat = \lowrank{\mat{B}}_k$ is shown as a dotted black line.
    We plot the median of 100 trials, with error bars showing 10\% and 90\% quantiles.
    \emph{Right:} Points $w_j$ in column subset $j \in \set{S}$ selected by random pivoting methods (\emph{top}) and the Mahoney--Drineas algorithm (\emph{bottom}) with $k=30$.}
    \label{fig:rpcur}
\end{figure}

\Cref{fig:rpcur} shows results for \RPCURtwo, \RPCURLev (with oversampling $\ell = \lceil 1.5k\log k \rceil$ and no oversampling $\ell = k$), and \MDCUR.
As points of comparison, I also include the optimal rank-$k$ approximation.
We note the following conclusions:

\begin{itemize}
    \item \textbf{Near-optimality for random pivoting.}
    On this example, \RPCURtwo and \RPCURLev with oversampling both produces approximations of comparable quality to the optimal rank-$k$ approximation.
    \MDCUR, by contrast, converges at a slower rate.
    These results demonstrates the benefits of using random pivoting to select at least one of the row or column subsets.
    \item \textbf{Necessity of oversampling for \RPCURLev.}
    Comparing the accuracy of the \RPCURLev method with and without oversampling, we see that oversampling $\ell > k$ is crucial to maintain accuracy of the procedure.
    \item \textbf{Stability benefits of factored $\mat{U}$ matrix representation \cref{eq:cur-U-factored}.}\CURindex{!stable representation}
    The solid light curves show the results of the \CUR algorithms with the stable representation\CURindex{!stable representation} \cref{eq:cur-U-factored}, which eventually achieve machine accuracy and even slightly beat the optimal rank-$k$ approximation \warn{due to finite-precision effects}. 
    The dark dashed curves show the results of the \CUR algorithms with an explicit $\mat{U}$ matrix.
    When the accuracy is smaller ($\lessapprox u^{1/2}$ for the random pivoting methods), the explicit $\mat{U}$ implementations have similar accuracy to the factored representation.
    However, while the relative error of the stable methods improve past this point, the relative error of the unstable implementations begins to \emph{increase} for large values of $k$.
    These results confirm the significant numerical issues of the standard implementation of the \CUR decomposition for high-accuracy matrix approximation.
    \item \textbf{Visualization of column subsets.}
    Column subsets $\set{S}$ of $k=30$ elements selected by random pivoting methods and the Mahoney--Drineas method are depicted visually in the right panels of \cref{fig:rpcur}.
    Both methods cluster landmarks near the edges of the semicircle, with random pivoting methods producing a more even distribution.
    Because the Mahoney--Drineas method samples points iid, its landmark set contains examples of nearly overlapping landmarks and regions that lack landmarks, which explains the differences in approximation quality with the random pivoting methods.
\end{itemize}

\myparagraph{Example \#2: Sparse matrix}
As a second example, we evaluate the methods on the \texttt{Meszaros/large} matrix from the SuiteSparse collection.
This example was used in the recent paper \cite[\S6.1]{CK24}.
We store this matrix as a dense matrix for the experiments in this section.

\begin{figure}[t]
    \centering
    \includegraphics[width=0.48\textwidth]{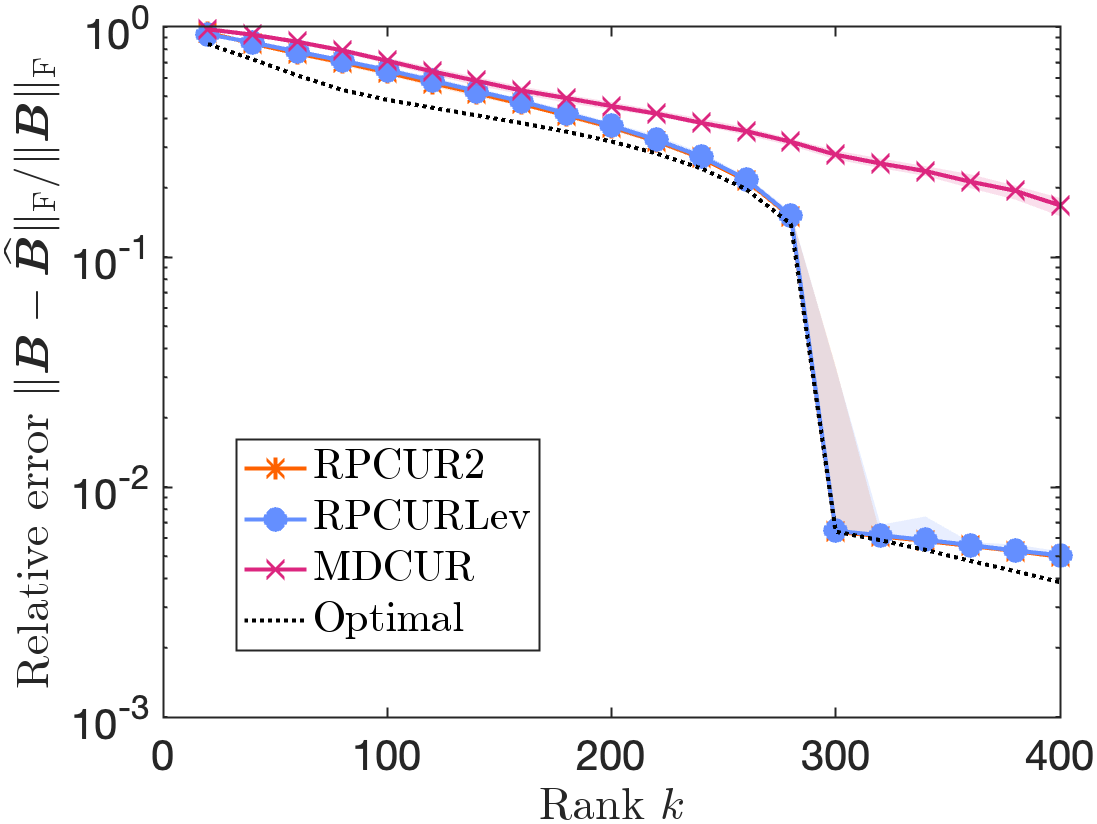}
    \includegraphics[width=0.48\textwidth]{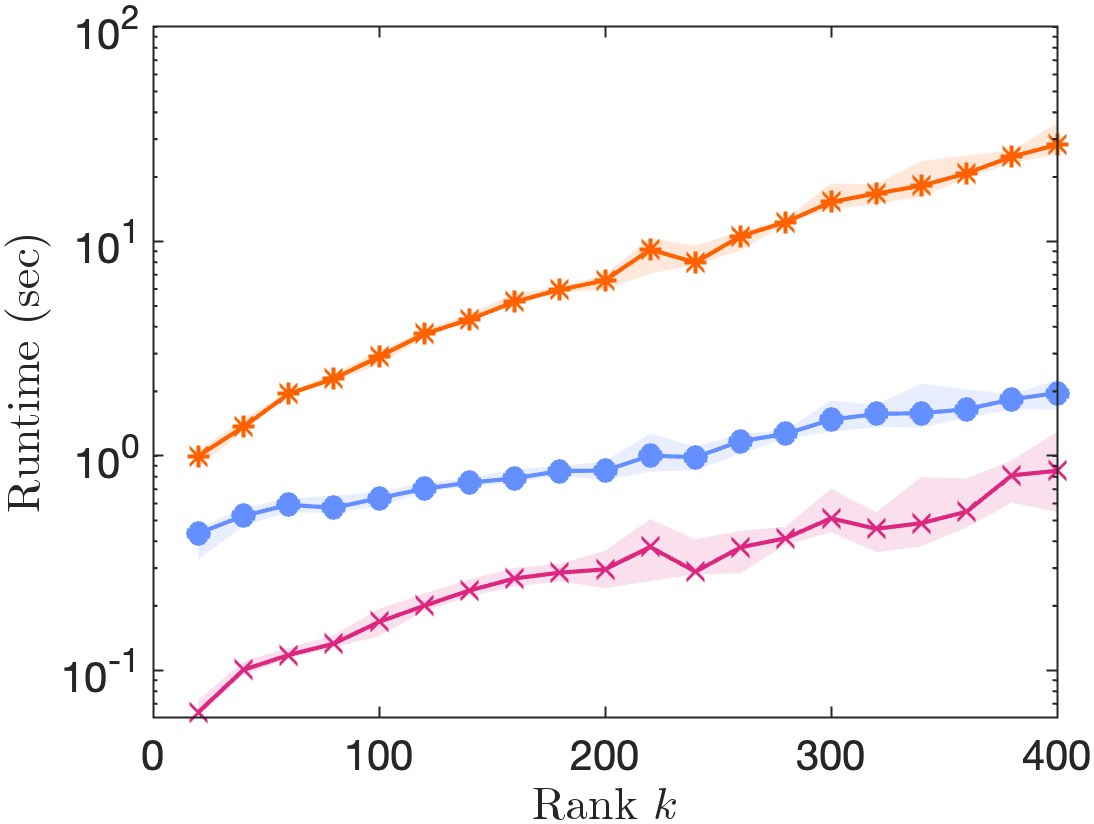}
    \caption[Comparison of speed and accuracy for various \CUR algorithms applied to a sparse matrix]{Relative error (\emph{left}) and runtime (\emph{right}) for \RPCURtwo (orange stars), \RPCURLev (blue circles), and \MDCUR (pink crosses) on \texttt{Meszaros/large} sparse matrix.
    We plot the median of 100 trials, with error bars showing 10\% and 90\% quantiles.}
    \label{fig:rpcur_sparse}
\end{figure}

We evaluate the runtime and accuracy of \RPCURtwo, \RPCURLev (with oversampling $\ell = \lceil 1.5k\log k\rceil$), and \MDCUR on this example.
To make the \MDCUR competitive in terms of runtime, we use the randomized SVD (\cref{prog:rsvd}) with rank $k'$ to approximate the dominant left and right singular vectors of $\mat{B}$.
\RPCURtwo and \MDCUR are implemented without oversampling $k=\ell$.
Due to the sparsity of $\mat{B}$, the column pivoted \QR decomposition in \cref{prog:rpcur_lev} is necessary for the algorithm to succeed.

Results are shown in \cref{fig:rpcur_sparse}.
The \MDCUR algorithm is the fastest, and \RPCURtwo is the slowest.
Both \RPCURtwo and \RPCURLev achieve accuracy comparable to the optimal rank-$k$ approximation, whereas \MDCUR lags behind the other methods, particularly when $k\ge 300$.
These results demonstrate how, on difficult examples, \MDCUR can produce approximations comparable with the best rank-($r<k$) approximation, whereas the random pivoting approximations empirically produce approximations comparable to the best rank-$k$ approximation.
The need for oversampling presents a weakness of pure leverage score sampling based approaches over random pivoting approaches.

\CURindex{!numerical results|)}\index{MDCUR algorithm@\textsc{MDCUR} algorithm!numerical results|)}\index{RPCUR2 algorithm@\textsc{RPCUR2} algorithm!numerical results|)}\index{RPCURLev algorithm@\textsc{RPCURLev} algorithm!numerical results|)}

\chapter{Random pivoting: Open problems}

\index{randomly pivoted Cholesky!open problems|(}
We conclude this part of the thesis with two open problems: improved error bounds for \RPCholesky (\cref{sec:better-rpcholesky-bounds}) and whether there are simple, effective algorithms for approximating a psd matrix from a small number of entry accesses when errors are measured using the Frobenius norm, rather than the trace (\cref{sec:frobenius-psd}).

\index{randomly pivoted Cholesky!theoretical results|(}
\section{Open problem: \RPCholesky error bounds} \label{sec:better-rpcholesky-bounds}

The bound \cref{thm:rpcholesky-trace} already provides a useful description of the performance of \RPCholesky, showing that one achieves an $(r,\varepsilon)$-approximation (\cref{def:r-eps-p}) in 
\begin{equation} \label{eq:rpcholesky-bound-open}
    k \le \frac{r}{\varepsilon} + r \log\left( \frac{1}{\varepsilon\eta} \right) \text{ steps}.
\end{equation}
Here, 
\begin{equation*}
    \eta = \frac{\tr(\mat{A} - \lowrank{\mat{A}}_r)}{\tr(\mat{A})}
\end{equation*}
is the relative error of the best rank-$r$ approximation.
The bound \cref{thm:rpcholesky-spectral} provides weak---but sometimes useful control---related to the spectral norm.
There are a few directions along which one could hope to improve our theoretical understanding of \RPCholesky.

\myparagraph{Bounds in the limit $\eta \downarrow 0$}
As noted in \cref{rem:relative-error}, the bound \cref{eq:rpcholesky-bound-open} degenerates to infinity in the limit when the relative error $\eta \downarrow 0$.
At least in some cases, this bound badly mischaracterizes the performance of the \RPCholesky algorithm, which converges in exactly $r$ steps when applied to a matrix with $\rank \mat{A} = r$.

For cases when $\eta$ is small, \cite[Thm.~5.1]{CETW25} establishes that the $(r,\varepsilon)$-guarantee holds when
\begin{equation*} \label{eq:rpcholesky-quadratic-r}
    k \ge \frac{r}{\varepsilon} + r+ \max \left\{ 0, r^2\log(2) + r\log(1/\varepsilon) \right\}.
\end{equation*}
This bound is nice because it purges the relative error $\eta$ from the bound entirely
On the other hand, $k$ is required to be \warn{quadratic} in $r$ to produce an approximation comparable to the best-$r$ approximation---ouch! 

R.\ J.\ Webber and I conjecture that this bound can be improved as follows.

\begin{conjecture}[Better $\eta$-free \RPCholesky bounds] \label{conj:eta-free-bound}
    \RPCholesky achieves an $(r,\varepsilon)$-approximation to any psd matrix in $r/\varepsilon + \order(r\log(r/\varepsilon))$ operations.
\end{conjecture}

Here is one concrete approach that could resolve this conjecture.
Specifically, we make another conjecture:

\begin{conjecture}[Improved $r$-step bound] \label{conj:r-step-bound}
    For any psd matrix $\mat{A}$, the $r$-step \RPCholesky residual $\mat{A}^{(r)}$ satisfies
    \begin{equation} \label{eq:conj-poly-bound}
        \tr(\mat{A} - \mat{A}^{(r)}) \le C(r) \tr(\mat{A} - \lowrank{\mat{A}}_r) \quad \text{for } C(r) = \poly(r).
    \end{equation}
\end{conjecture}

\Cref{conj:r-step-bound} would imply \cref{conj:eta-free-bound}.
Unfortunately, we are currently exponentially far from proving \cref{conj:r-step-bound}; our best bound \cite[Lem.~5.5]{CETW25} establishes \cref{eq:conj-poly-bound} with $C(r) = 2^r$.

\index{randomly pivoted Cholesky!spectral-norm bounds|(}\index{spectral-norm error bounds for column Nystr\"om approximation|(}
\myparagraph{Spectral-norm bounds}
For reasons that will be described in \cref{sec:frobenius-psd}, it is not possible for \RPCholesky or any psd low-rank approximation algorithm that uses limited entry accesses to produce a relative-error approximation in the spectral norm (that is, an $(r,\varepsilon,\infty)$-approximation) in $\order(n\poly(r/\varepsilon))$ entry access operations.
However, it is reasonable to hope that \RPCholesky produces approximations that are accurate when measured in the spectral norm once the number of steps is comparable to the \emph{effective dimension},\index{effective dimension!conjectured results for randomly pivoted Cholesky} defined in \cref{def:ridge-psd}.
Indeed, we saw earlier in \cref{fact:rls-spectral} that ridge leverage score sampling achieves such a bound.
We conjecture that \RPCholesky does equally well.

\begin{conjecture}[\RPCholesky: Spectral norm bounds]
    Let $\lambda \ge 0$.
    If \RPCholesky is executed for $k = \order(\mathrm{d}_{\mathrm{eff}}(\lambda) \polylog(\mathrm{d}_{\mathrm{eff}}(\lambda),\norm{\mat{A}}/\lambda))$ steps, it produces an index set satisfying\index{effective dimension!conjectured results for randomly pivoted Cholesky}
    \begin{equation*}
        \mat{A}\langle \set{S}\rangle \preceq \mat{A} \preceq \mat{A}\langle \set{S}\rangle + \lambda \Id
    \end{equation*}
    with at least 99\% probability.
    Here, $\polylog$ denotes an unspecified polylogarithmic function.
\end{conjecture}

If such a result were true, it would provide theoretical support for the empirical observation that \RPCholesky performs similarly to or better than RLS sampling for column Nystr\"om preconditioning\index{column Nystr\"om preconditioning} (\cref{sec:full-data-preconditioning}).
To be conservative, I have included polylogarithmic factors in the conjecture, but I believe the conjecture will also be true without them.\index{randomly pivoted Cholesky!spectral-norm bounds|)}\index{spectral-norm error bounds for column Nystr\"om approximation|)}\index{randomly pivoted Cholesky!theoretical results|)}\index{randomly pivoted Cholesky!open problems|)}

\index{Frobenius-norm positive-semidefinite low-rank approximation|(}\index{r,epsilon,p-approximation!reasons for focus on $p=1$|(}
\section{Open problem: Frobenius-norm psd low-rank approximation} \label{sec:frobenius-psd}

The goal of this part of the thesis has been to develop algorithms that evaluate a small number of entries to compute an $(r,\varepsilon,1)$-approximation to a psd matrix.
Recall from \cref{def:r-eps-p} that an $(r,\varepsilon,p)$-approximation is a random matrix $\Ahat$ such that
\begin{equation*}
    \expect \norm{\mat{A} - \Ahat}_{\set{S}_p} \le (1+\varepsilon) \norm{\mat{A} - \lowrank{\mat{A}}_r}_{\set{S}_p}.
\end{equation*}
Recall that $\norm{\cdot}_{\set{S}_p}$ is the Schatten $p$-norm.
It is frame the question: Does the psd low-rank approximation problem admit efficient algorithms when $p > 1$?

Many of the algorithms we have considered so far, including \RPCholesky, uniform column Nystr\"om approximation, greedy pivoted partial Cholesky, and the entire family of Gibbs \RPCholesky methods (\cref{sec:gibbs}), are \emph{diagonal--column access algorithms}.\index{diagonal--column access algorithm}
These methods interact with the matrix $\mat{A}$ only by reading the diagonal and full columns.
As we will see, diagonal--column access algorithms\index{diagonal--column access algorithm} have fundamental limitations for psd low-rank approximation, and a more general class of algorithms are needed order to compute $(r,\varepsilon,p)$-approximations for $p > 1$.

The following result establishes establishes limits on both general and diagonal--column access\index{diagonal--column access algorithm} psd low rank approximation algorithms.

\begin{proposition}[Psd low-rank approximation: Lower bound] \label{prop:psd-low-rank-lower}
    Fix $p \in [1,\infty]$.
    On a worst-case input matrix,
    \begin{itemize}
        \item A diagonal--column access algorithm\index{diagonal--column access algorithm} must read $\Omega(n^{2-1/p})$ entries to guarantee a $(1,\mathrm{c},p)$-approximation.  
        \item Any algorithm must read $\Omega(\max\{n,n^{2-2/p}\})$ entries to guarantee a $(1,\mathrm{c},p)$-approximation.  
    \end{itemize}
    Here, $\mathrm{c} > 0$ is a universal constant.
\end{proposition}

\begin{proof}[Proof sketch]
    Choose $n$ to be large, and construct the matrix 
    \begin{equation*}
        \mat{A} = \mat{P}\twobytwo{\Id_{n-\lceil n^{1/p}\rceil}}{\mat{0}}{\mat{0}}{\outprod{\onevec_{\lceil 2n^{1/p}\rceil}}}\mat{P}^*,
    \end{equation*}
    where $\mat{P}$ is a uniformly random permutation matrix.
    The optimal rank-one approximation to this matrix is 
    \begin{equation*}
        \lowrank{\mat{A}}_1 = \mat{P}\twobytwo{\mat{0}}{\mat{0}}{\mat{0}}{\outprod{\onevec_{\lceil 2n^{1/p}\rceil}}}\mat{P}^*.
    \end{equation*}
    We compute
    \begin{align*}
        \norm{\mat{A}}_{\set{S}_p} &= \left( (n-\lceil 2n^{1/p}\rceil) \cdot 1 + (\lceil 2n^{1/p}\rceil)^p \right)^{1/p} &&\ge (2+o(1))n^{1/p}, \\
        \norm{\mat{A} - \lowrank{\mat{A}}_1}_{\set{S}_p} &= \left( (n-\lceil 2n^{1/p}\rceil) \cdot 1  \right)^{1/p} &&= (1+o(1)) n^{1/p}.
    \end{align*}
    To produce a nontrivial low-rank approximation requires identifying a nonzero off-diagonal entry.
    
    Since each diagonal entry is one, any diagonal--column access algorithm\index{diagonal--column access algorithm} must keep accessing columns until it finds a nonzero off-diagonal entry.
    Each column access has a $\lceil 2n^{1/p} \rceil / n = (2+o(1))n^{1/p-1}$ probability of identifying a nonzero off-diagonal entry, so it takes $\Omega(n^{1-1/p})$ column accesses ($=\Omega(n^{2-1/p})$ entry accesses) to find a nonzero off-diagonal entry with high probability.

    A general psd low-rank approximation algorithm is free to query off-diagonal entries one-by-one.
    Until a nonzero entry is found, each query finds a nonzero off-diagonal entry with probability $\order(n^{2/p-2})$, so it requires $\Omega(n^{2-2/p})$ accesses to find a nonzero off-diagonal entry with high probability.
    The fact that it always requires $\Omega(n)$ accesses to produce a nontrivial low-rank approximation can be established by considering the matrix $\outprod{\evec_i}$ where $i \sim \Unif \{1,\ldots,n\}$.
\end{proof}

\Cref{prop:psd-low-rank-lower} shows that, for \RPCholesky and other diagonal--column access algorithms, approximating matrices in the trace norm is essentially as good as one can expect, at least for an algorithm with an $\order(n)$ runtime.\index{r,epsilon,p-approximation!reasons for focus on $p=1$|)}
However, this result leaves open the possibility of achieving relative-error approximations in the Frobenius norm with $\order(n)$ operations using a more general access pattern.
Indeed, such $\order(n)$ algorithms for producing $(r,\varepsilon,2)$-approximations do exist.
The first such algorithm was discovered by Musco and Woodruff \cite{MW17a}, which produces an $(r,\varepsilon,2)$ in $\tilde{\order}(nr/\varepsilon^{2.5})$ entry accesses; this was improved to $\tilde{\order}(nr/\varepsilon)$ entries in \cite{BCW20}.
Despite these appealingly low entry access counts, the algorithms of \cite{MW17a,BCW20} are complicated and have large prefactor constants.
That leads us to the open question:

\actionbox{Is there a simple, effective, and practically performant algorithm for psd low-rank approximation that achieves relative error guarantees in the Frobenius norm?}

As \cref{prop:psd-low-rank-lower} shows, such an algorithm must use a more sophisticated access model than diagonal and column accesses alone.
Additionally, even basic existence questions are open, as far as I am aware.
For instance, does there exist a set $\set{S}$ of $\poly(r/\varepsilon)$ pivots defining a column Nystr\"om approximation $\mat{A}\langle \set{S}\rangle$ that is an $(r,\varepsilon,2)$ approximation?
Or must we use more sophisticated low-rank approximation formats, such as weighted \CUR cross approximations,\CURindex{!CUR cross approximation@\textsf{CUR} cross approximation} to produce Frobenius-norm relative-error approximations?

\index{Frobenius-norm positive-semidefinite low-rank approximation|)}

\partepigraph{Dedicated to my parents Meg and Tom Epperly.}
\part{Leave-one-out randomized matrix algorithms} \label{part:loo}

\chapter{Matrix attribute estimation problems}

\epigraph{This kind of structure [sparsity] is readily exploited by the iterative methods we shall discuss, for these algorithms use a matrix in the form of a \emph{black box}.
The iterative algorithm requires nothing more than the ability to determine $\mat{A}\vec{x}$ for any $\vec{x}$, which in a computer program will be effected by a procedure whose internal workings need be of no concern to the designer of the iterative algorithm.}{Lloyd N. Trefethen and David Bau, III, \textit{Numerical Linear Algebra} \cite[p.~244]{TB22}}

Frequently in applications, we encounter problems where we wish to learn information about an unknown matrix $\mat{B} \in \complex^{m\times n}$ that can only be accessed by matrix--vector products $\vec{\omega} \mapsto \mat{B}\vec{\omega}$ and, possibly, matrix--vector products with the adjoint $\vec{\omega} \mapsto \mat{B}^*\vec{\omega}$.\index{matrix attribute estimation problem}
As examples, we might be interested in learning the matrix trace $\tr(\mat{B})$, the matrix diagonal $\diag(\mat{B})$, or individual entries $b_{ij}$.

In this part of the thesis, I will present a new approach to designing randomized algorithms for such \emph{matrix attribute estimation} problems called the \emph{leave-one-out\index{leave-one-out randomized algorithm}} approach.
The leave-one-out\index{leave-one-out randomized algorithm} approach leads to fast, resource-efficient algorithms with state-of-the-art accuracy.
After introducing the necessary preliminaries in \cref{ch:fundamental-tools}, we will discuss the leave-one-out\index{leave-one-out randomized algorithm} approach in \cref{ch:loo}.
This initial chapter serves to motivate the problem of estimating attributes of a matrix through matrix--vector products.
We will address the questions ``When are we able to access a matrix $\mat{B}$ only through matrix--vector products?'' and ``Which attributes of $\mat{B}$ do we want to estimate in applications?''

\myparagraph{Sources}
This is an introductory chapter and is not based on any particular research article.
The resources \cite{US17,Pop23} are useful resources for applications of trace and diagonal estimation algorithms, respectively.

\myparagraph{Outline}
\Cref{sec:matvec} discusses matrix attribute estimation problems\index{matrix attribute estimation problem} and the matrix--vector (``matvec'') model for measuring their computational cost.
\Cref{sec:matrix-attribute-examples} discusses examples of matrix attribute estimation problems including trace, diagonal, and row-norm estimation.
It also discusses applications.

\index{matvec model|(}
\section{Matrix attribute estimation and the matvec model} \label{sec:matvec}

Linear algebraic algorithms that access a matrix $\mat{B}$ through matrix--vector products (affectionately abbreviated \emph{matvecs}) are extremely common.
Examples include power iteration\index{power iteration} for computing eigenvalues \cite[\S8.2]{GV13} as well as Krylov iteration\index{Krylov method} for solving linear systems \cite{Saa03}, computing eigenvalues \cite{Saa11}, and applying matrix functions \cite{Che24}. 
For these classical iterative methods, the matrix $\mat{B}$ is often stored explicitly in memory as either a dense or sparse array, and matvec algorithms are used due to their computational efficiency (particularly when $\mat{B}$ is sparse, making $\mat{z}\mapsto \mat{B}\vec{z}$ is cheap to compute).

This thesis focuses on an alternative setting where the matrix $\mat{B}$ is not stored explicitly and $\mat{B}$ can \emph{only} be accessed through matvecs (and possibly matvecs with the adjoint $\mat{B}^*$).
This more restrictive setting is motivated by the following examples:
\begin{itemize}
    \item \textbf{Products of matrices.} The matrix $\mat{B}$ is a product of several other matrices, $\mat{B} = \mat{C}_1 \mat{C}_2 \cdots \mat{C}_\ell$.
    It is expensive to evaluate the matrix $\mat{B}$, as it requires forming the product, but matvecs with $\mat{B}$ can be computed efficiently by multiplying a vector against the $\mat{C}_i$ matrices one at a time: $\mat{B}\vec{\omega} = \mat{C}_1(\mat{C}_2(\cdots (\mat{C}_\ell\vec{\omega}) \cdots )$.
    Important special cases are the Gram matrix $\mat{B} = \mat{C}^*\mat{C}$ and powers $\mat{B} = \mat{C}^k$.
    \item \textbf{Matrix functions.}\index{matrix function} Let $\mat{C}\in\complex^{n\times n}$ be a matrix and $f : \complex \to \complex$ be a function defined on the spectrum of $\mat{C}$.
    Set $\mat{B} \coloneqq f(\mat{C})$ using the standard extension of scalar functions to matrix inputs \cite[\S1.2]{Hig08}.
    Computing $\mat{B} = f(\mat{C})$ explicitly is typically costly, but matvecs $\mat{B}\vec{\omega} = f(\mat{C})\vec{\omega}$ can often be computed efficiently by algorithms such as the Arnoldi method\index{Arnoldi method} or, if $\mat{C}$ is Hermitian, the Lanczos method\index{Lanczos method!for matrix functions} \cite[\S13.2]{Hig08}.
    (In particular, for the special case $f(z) = z^{-1}$, computing products $\mat{C}^{-1}\vec{\omega}$ is equivalent to solving linear systems of the form $\mat{C} \vec{y} = \vec{\omega}$, for which there has been \emph{much} work \cite{Saa03}.)
    \item \textbf{Automatic differentiation.}\index{automatic differentiation} Given a twice differentiable function $f : \real^n \to \real$ and an input $\vec{x} \in\real^n$, automatic differentiation allows one to compute matvecs $D^2 f(\vec{x}) \cdot \vec{\omega}$ with the \emph{Hessian matrix}\index{Hessian matrix} $D^2 f(\vec{x}) \in \real^{n\times n}$ in a small multiple of the runtime required to evaluate $f(\vec{x})$ \cite{BR24}.
    Similarly, given a function $g : \real^n\to \real^m$, matvecs $Dg(\vec{x})\cdot \vec{\omega}$ with the Jacobian matrix\index{Jacobian matrix} $Dg(\vec{x}) \in \real^{m\times n}$ are also efficiently computable by automatic differentiation.
\end{itemize}
In all of these settings, matvecs with the adjoint $\mat{B}^*$ are typically also available.
These examples motivate the \emph{matvec} model:

\actionbox{\textbf{Matvec model.} A matrix $\mat{B} \in \real^{m\times n}$ is provided that can only be accessed through matvecs $\vec{\omega} \mapsto \mat{B}\vec{\omega}$ and adjoint matvecs $\vec{\omega} \mapsto \mat{B}^*\vec{\omega}$.
The cost of an algorithm is measured by the number of matvecs used.}

Under the matvec model, the entire matrix $\mat{B}$ can be recovered one column at a time by computing the matvecs $\mat{B}\evec_1,\ldots,\mat{B}\evec_n$ with each of the standard basis vectors $\evec_i$.
Thus, in the matvec model, any problem can be solved at the ``trivial cost'' of $n$ matvecs.
Thus, in this model, an algorithm is considered ``efficient'' if it beats this trivial cost of $n$ matvecs.
In this part of the thesis, we will seek algorithms that produce approximate solutions but will expend a number of matvecs that is \emph{independent} of the dimensions $m$ and $n$ of the matrix.

Problems that are trivial for a matrix $\mat{B}$ stored in memory become more challenging when considered in the matvec model.
For instance, consider the trace\index{trace estimation}
\begin{equation*}
    \tr(\mat{B}) = \sum_{i=1}^n b_{ii}
\end{equation*}
of a square matrix $\mat{B} \in \real^{n\times n}$.
If $\mat{B}$ is stored in memory, the trace can be computed exactly in $\order(n)$ operations.
However, in the matvec model, computing the trace \warn{exactly}---by either a deterministic or randomized algorithm---incurs the trivial cost of $n$ matvecs.
This observation motivates the study of algorithms for \warn{approximating} the trace in the matvec model.
We will discuss trace estimation in the matvec model throughout this part of the thesis.

Trace estimation is a member of the family of \emph{matrix attribute estimation problems}:\index{matrix attribute estimation problem}
\actionbox{\textbf{Matrix attribute estimation problem.} Let $Q(\mat{B})$ denote some attribute of a matrix such as its trace $\tr(\mat{B})$ or its entry $b_{ij}$ in position $(i,j)$. 
The matrix attribute estimation problem is to compute an estimate $\hat{Q}$ for $Q(\mat{B})$.
The cost an algorithm is measured by the number of matvecs, and the quality of the solution is measured by the error $|\hat{Q} - Q(\mat{B})|$.}
\noindent We have written $Q(\mat{B})$ for a scalar-valued attribute $Q(\mat{B}) \in \field$, but vector-valued or matrix-valued attributes $\vec{q}(\mat{B}) \in \field^d$ or $\mat{Q}(\mat{B}) \in \field^{m'\times n'}$ are fine as well.
In these cases, we will measure the error using an appropriate norm.

In a sense, any linear algebra problem is an example of a matrix attribute estimation problem.\index{matrix attribute estimation problem}
In this thesis, we will use this term more narrowly to describe attributes like the trace that are linear or quadratic functions of the matrix $\mat{B}$ and are computationally challenging only within the matvec model.

\index{matvec model!limitations|(}
As we will see in this part of the thesis, the matvec model provides a convenient framework for designing and analyzing algorithms.
However, like all computational models, the matvec model is idealized.
We note a few weaknesses:
\begin{itemize}
    \item \textbf{Post-processing costs.}
    In the matvec model, we measure the cost of an algorithm only by counting the number of matvecs.
    However, for problems where the matvec operations are relatively cheap, the computational cost can be dominated by post-processing to assemble the collected matvecs into the estimator. 
    Responding to this weakness, this thesis will focus on algorithms with fast post-processing.
    We will even develop improved implementations of one estimator with up to $7\times$ faster post-processing  (\cref{rem:xnystrace-new-impl}).
    \item \textbf{Blocking.} Another deficit of measuring algorithm cost through matvec count is that it ignores the computational benefits of blocking.
    For reasons discussed in \cref{sec:why-blocking}, performing $\ell$ matvecs $\mat{B}\vec{\omega}_1,\ldots,\mat{B}\vec{\omega}_\ell$ individually can be much more expensive than performing these matvecs in a batch,\index{matrix--matrix operations}\index{blocked algorithm!matrix--matrix products} i.e., by computing the matrix--matrix product $\mat{B}[\begin{matrix}
        \vec{\omega}_1 & \cdots & \vec{\omega}_\ell
    \end{matrix}]$.
    The algorithms in this thesis will always compute matvecs in batches, so they remain efficient even accounting for the benefits of blocking.
    
    \item \textbf{Opening the black box.} The matvec model treats the matvec subroutine $\vec{\omega} \mapsto \mat{B}\vec{\omega}$ as a black box.
    But, for matrix functions\index{matrix function} $\mat{B} = f(\mat{C})$, matvecs $f(\mat{C})\vec{\omega}$ are often computed by the Arnoldi\index{Arnoldi method} or Lanczos methods,\index{Lanczos method!for matrix functions} which utilize matvecs with $\mat{C}$.
    In this context, matvecs with $\mat{C}$ are the primitive operation that should be minimized.
    Tyler Chen and collaborators have productively used this observation to develop faster ``Krylov-aware'' matrix attribute estimation algorithms \cite{CH23a,PCM25}.\index{Krylov-aware algorithms}
\end{itemize}\index{matvec model!limitations|)}

Notwithstanding these limitations, we will use the matvec model in this thesis.\index{matvec model|)}

\section{Examples of matrix attribute estimation problems: Trace, diagonal, and row-norm estimation} \label{sec:matrix-attribute-examples}

There are several interesting and practically useful matrix attribute estimation problems,\index{matrix attribute estimation problem} including trace, diagonal, and row-norm estimation.

\index{trace estimation|(}\index{trace estimation!applications of|(}
\subsection{Trace estimation}

The trace estimation problem is to approximate $\tr(\mat{B})$ in the matrix--vector product (matvec) model.\index{matvec model}
Trace estimation has a number of applications, some of which were surveyed in the paper \cite{US17}.
Here is a partial list:

\begin{itemize}
    \item \textbf{Triangle counting.}\index{triangle counting}
    Counting the number of triangles in a graph is a basic problem in network science, and it is necessary to compute the clustering coefficient \cite{AD18}.
    This problem is a trace estimation problem, as the number of triangles in a graph with adjacency matrix $\mat{A}$ is $\tr(\mat{A}^3)/6$.
    This problem is a natural demonstration of the matvec model, as matvecs with $\mat{A}^3$ can be computed by iterated multiplication $\mat{A}(\mat{A}(\mat{A}\vec{\omega}))$.
    \item \textbf{Log-determinant \cite{MMB21,WPH+22}.}
    Estimates of the log-determinant\index{log-determinant} are used in maximum likelihood estimation for Gaussian process\index{Gaussian process!type of machine learning method} methods (\cref{rem:mle}).
    This problem can be viewed as a trace estimation problem, in view of the identity $\log \det(\mat{A}) = \tr \log(\mat{A})$.
    Matvecs with $\log(\mat{A})$ can be computed using the Lanczos method.
    \item \textbf{Continuous normalizing flows \cite{CRBD18,GCB+19,SSK+21,Now22}.}\index{continuous normalizing flow}
    Several machine learning models---such as neural ODEs \cite{CRBD18}, FFJORD \cite{GCB+19}, and diffusion models \cite{SSK+21}---evolve a random initial value $\vec{x}(0)$ under continuous dynamics $\tfrac{\d}{\d t} \vec{x}(t) = f(\vec{x}(t),t)$.
    These algorithms must estimate the instantaneous rate of change of the log-likelihood\index{log-likelihood} of $\vec{x}(t)$, which is the negative trace of the Jacobian\index{Jacobian matrix} $-\tr(\mathrm{D}f)$.
    Matvecs with the Jacobian can be computed using automatic differentiation.
    \item \textbf{Statistical physics.}\index{statistical physics}
    The partition function\index{partition function} of a quantum system with Hamiltonian $\mat{A}$ at inverse-temperature $\beta > 0$ is $z \coloneqq \tr \exp(-\beta \mat{A})$.
    Other thermodynamics quantities can also be expressed as matrix traces; for instance, the average energy is $e \coloneqq z^{-1} \tr(\mat{A} \exp(-\beta \mat{A}))$.
    Trace-exponentials also appear in network science as Estrada index\index{Estrada index} centrality measures \cite{Est22}.
    \item \textbf{Norm estimation.}
    Computing or estimating the norm of a matrix is a ubiquitous task in matrix computations.
    This computation can be seen as a trace estimation problem, as the Schatten $p$-norm\index{Schatten $p$-norm} is $\norm{\mat{B}}_{\set{S}_p} \coloneqq \tr[(\mat{B}^*\mat{B})^{p/2}]$.
    Matvecs with $(\mat{B}^*\mat{B})^{p/2}$ can be computed by iterated multiplication (if $p$ is an even integer) or by the Lanczos method.
\end{itemize}
\index{trace estimation|)}\index{trace estimation!applications of|)}

\index{diagonal estimation|(}\index{diagonal estimation!applications of|(}
\subsection{Diagonal estimation}

Given a square matrix $\mat{B} \in \field^{n\times n}$, the diagonal estimation problem is to estimate $\diag(\mat{B}) \in \field^n$.
There are several applications of diagonal estimation algorithms.

\begin{itemize}
    \item \textbf{Statistics.} The diagonal entries of a covariance matrix\index{covariance matrix} $\mat{A}$ for a random vector $\vec{x}\in\field^n$ are the variances of the individual entries $a_{ii} = \Var(x_i)$.
    Similarly, \warn{for jointly Gaussian random variables}, the \warn{reciprocals} of the diagonal of the inverse-covariance matrix\index{covariance matrix} (precision matrix) $\mat{A}^{-1}$ store the \emph{conditional} variances:
    \begin{equation*}
        \frac{1}{\mat{A}^{-1}(i,i)} = \Var(x_i \mid x_1,\ldots,x_{i-1},x_{i+1},\ldots,x_n).
    \end{equation*}
    \item \textbf{Centrality measures.}\index{centrality measure} Many centrality measures for graphs are defined using the diagonal of functions of the graph adjacency matrix $\mat{A}$ (or related matrices such as the graph Laplacian\index{Laplacian matrix}).
    For instance, the number of triangles incident on vertex $i$ is the $i$th diagonal entry of $\mat{A}^3/2$ \cite{AD18}, and the \emph{subgraph centrality}\index{subgraph centrality} of vertex $i$ is the $i$th diagonal entry of $\exp(\mat{A})$ \cite{Est22}.
    For a survey on centrality measures,\index{centrality measure} see \cite{BJT23}.
    \item \textbf{Optimization.}
    The convergence of gradient descent\index{gradient descent} methods can be slow if the problem is poorly scaled.
    This issue can be remedied by using the diagonal of the Hessian matrix\index{Hessian matrix} to precondition the descent method \cite{YGS+21}.
    Another application of diagonal estimation to semidefinite programming\index{semidefinite programming} appears in \cite{Lin23a}.
\end{itemize}
Other applications of diagonal estimation include electronic structure calculations in materials science \cite{BKS07,LLY+09},\index{materials science} uncertainty quantification for linear regression \cite{Pop23}, low-rank matrix approximation algorithms \cite{FL24}.
\cref{ch:diagonal} discusses diagonal estimation algorithms.

\index{diagonal estimation|)}\index{diagonal estimation!applications of|)}

\index{row-norm estimation|(}\index{row-norm estimation!applications of|(}

\subsection{Row-norm estimation}
The problem of estimating the (squared) row- or column-norms\index{squared row or column norms} of a matrix $\mat{B} \in \field^{m\times n}$ also has many applications.
This problem has received significantly less attention than trace or diagonal estimation.
Assuming matvec access to both $\mat{B}$ and $\mat{B}^*$, the row- and column-norm problems are equivalent, as column-norm estimation on $\mat{B}$ is row-norm estimation on $\mat{B}^*$.
Here are several applications:
\begin{itemize} 
    \item \textbf{Leverage scores.}\index{leverage scores!estimation of}
    As introduced in \cref{def:ridge-general}, the leverage scores of a matrix $\mat{A}$ are the squared row norms of any orthonormal basis matrix $\mat{U}$ for $\range(\mat{A})$.
    Row-norm estimation algorithms play a crucial role in fast algorithms for leverage score estimation \cite{DMMW12}.
    Related ideas are used to estimate effective resistances in networks \cite{SS08}.
    \item \textbf{Matrix computations.} Many standard linear algebra algorithms require computation of row norms, and many of these algorithms are robust to those row norms being computed approximately.
    For such algorithms, row-norm estimates can allow for significant acceleration (see, e.g., \cite{FL24}).
    Many of the codes in this thesis require row or column norms, and can also be accelerated using row-norm estimates (see \cref{rem:faster-loo}).
    \item \textbf{Electronic structure theory.}\index{material science}\index{electronic structure theory}
    In electronic structure theory calculations, the electron density can be obtained from the row norms of an implicit matrix, making it possible to accelerate the computation using row-norm estimation algorithms.
    This application appears in the thesis of Aleksandros Sobczyk \cite{Sob24a}.
\end{itemize}
Row-norm estimation also has close connections to the diagonal estimation problem (\cref{sec:square-root-trick}), and many of the most effective algorithms for diagonal estimation for \warn{psd} matrices proceed through row-norm estimation \cite{MMB21,Lin23a,FL24}.
See \cref{ch:row-norm} for more discussion of and algorithms for row-norm estimation.

\index{row-norm estimation|)}\index{row-norm estimation!applications of|)}

\chapter{Fundamental tools: Low-rank and Monte Carlo approximations} \label{ch:fundamental-tools}

\epigraph{What we are calling a quadratic trace estimator is often called the Hutchinson's trace estimator, especially when $\vec{v}$ is chosen uniformly from the set of vectors with entries $\pm n^{-1/2}$. However, \cite{Hut89} was not the first use of quadratic trace estimators for the task of approximating the trace of an implicit matrix; \cite{Hut89} itself cites \cite{Gir87} which addresses the same task by using samples of $\vec{v}$ drawn uniformly from the unit hypersphere. Algorithms based on the use of random vectors date back at least to the mid 1970s \cite{ABKS75,WW76,WW77,Dd89}.}{Tyler Chen, \emph{Lanczos-based methods for matrix functions} \cite[\S4.1.1]{Che22}}

This chapter introduces two approaches to matrix attribute estimation---Monte Carlo approximation and low-rank approximation---which serve as building blocks for matrix attribute estimation algorithms.
After introducing these techniques individually, we will see that they become more powerful when used in combination, as demonstrated in the famous \HutchPP algorithm.
This ideas will be refined in the next section, where we will use them in combination with the \emph{leave-one-out\index{leave-one-out randomized algorithm} approach} to the design of randomized matrix algorithm.

\myparagraph{Sources}
The main aims of this section are expository, and it is not based on any particular research article.
\Cref{sec:mc-plus-lra} discusses variance reduction technique for matrix attribute estimation and \HutchPP algorithm from the paper \cite{MMMW21}.
The concept of ``resphering'' matrix attribute estimation algorithms, which we discuss throughout this section, was first introduced in \cite{ETW24}.

\myparagraph{Outline}
\Cref{sec:monte-carlo-matrix-attribute} introduces the Monte Carlo approach to matrix attribute estimation, and it shows how this approach can be used to develop Monte Carlo estimators of the trace, diagonal, and row norms.
This section also discusses choice of random distribution for such algorithms and introduces the idea of ``resphering'' a Monte Carlo estimator for a matrix with a known nullspace.
\Cref{sec:lra} discusses the use of low-rank approximation to estimate matrix attributes, and \cref{sec:mc-plus-lra} describes how low-rank approximation can be used as a \emph{variance reduction} technique to improve the accuracy of Monte Carlo estimators.
This variance technique is exemplified by the \HutchPP algorithm, which is discussed in this section.
We conclude \cref{sec:mc-plus-lra} by discussing weaknesses of the \HutchPP algorithm which will be addressed by the leave-one-out\index{leave-one-out randomized algorithm} approach in \cref{ch:loo}.

\section{Monte Carlo approximation} \label{sec:monte-carlo-matrix-attribute}

Monte Carlo approximation is a simple yet powerful method for constructing randomized estimators of a matrix attribute.
In its basic form, the approach is to design an unbiased estimator $\hat{q}$ for a quantity of interest $q$ (so that $\expect[\hat{q}] = q$), and then average several independent copies of $\hat{q}$ to reduce variance \cite{Liu04}.
More general versions of the method also allow for estimators that are not fully independent or that introduce slight bias.

For matrix computations, many Monte Carlo estimators can be built using \emph{isotropic random vectors}.
\begin{definition}[Isotropic random vector]
    A random vector $\vec{\omega} \in \field^n$ is \emph{isotropic} if it satisfies $\expect[\vec{\omega}^{\vphantom{*}}\vec{\omega}^*] = \Id$.
\end{definition}
Examples of isotropic random vectors include standard Gaussian vectors $\vec{\omega} \sim \Normal_\field(\vec{0},\Id)$ or vectors $\vec{\omega} \sim \Unif(\sqrt{n}\,\sphere(\field^n))$ drawn uniformly from the sphere of radius $\sqrt{n}$.
Both of these constructions are defined in either the field of real or complex numbers, $\field \in \{\real,\complex\}$.
Another popular isotropic random distribution is the random sign distribution $\Unif \{ \pm 1\}^n$.

An isotropic random vector $\vec{\omega}$ gives rise to a rank-one matrix $\vec{\omega}^{\vphantom{*}}\vec{\omega}^*$ that serves as an unbiased estimator for to the identity matrix.
This estimator can be improved by forming an average $s^{-1} \sum_{i=1}^s \vec{\omega}_i^{\vphantom{*}} \vec{\omega}_i^*$ of independent copies $\vec{\omega}_1,\ldots,\vec{\omega}_s\simiid \vec{\omega}$.
Using this observation, one can design unbiased estimates for matrix attributes by introducing a copy of the identity matrix and replacing it by a stochastic approximation. 
In the following examples, we develop several classical Monte Carlo estimators in matrix computations using this perspective.

\subsection{Example: Girard--Hutchinson trace estimator}

As a first example of a Monte Carlo algorithm for matrix attribute estimation, we consider the trace estimation problem.
Let $\mat{B} \in \field^{n\times n}$ be a square matrix.
To estimate $\tr(\mat{B})$, introduce a copy of the identity matrix 
\begin{equation*}
    \tr(\mat{B}) = \tr(\mat{B} \cdot \Id)
\end{equation*}
and use the stochastic approximation $s^{-1} \sum_{i=1}^s \vec{\omega}_i^{\vphantom{*}} \vec{\omega}_i^* \approx \Id$, resulting in the estimator 
\begin{equation*}
    \hat{\tr}_{\mathrm{GH}} \coloneqq \tr\left(\mat{B} \cdot \frac{1}{s} \sum_{i=1}^s \vec{\omega}_i^{\vphantom{*}} \vec{\omega}_i^* \right) \approx \tr(\mat{B}\cdot \Id) = \tr(\mat{B}).
\end{equation*}
Invoking the linearity and the cyclic property of the trace, we can  rewrite the estimator $\hat{\tr}_{\mathrm{GH}}$ in the more computationally useful form
\begin{equation} \label{eq:girard-hutchinson}
    \hat{\tr}_{\mathrm{GH}} = \frac{1}{s} \sum_{i=1}^s \vec{\omega}_i^*\left(\mat{B}\vec{\omega}_i^{\vphantom{*}}\right).
\end{equation}
The formula \cref{eq:girard-hutchinson} clearly demonstrates that $\hat{\tr}_{\mathrm{GH}}$ can be computed using $s$ matvecs.
See \cref{prog:girard_hutchinson} for un-optimized MATLAB code.

\myprogram{Unoptimized MATLAB implementation of the Girard--Hutchinson estimator \cref{eq:girard-hutchinson}.}{The \texttt{random\_signs} subroutine is defined in \cref{prog:random_signs}.}{girard_hutchinson}

Estimators of the form \cref{eq:girard-hutchinson} were used by Girard \cite{Gir87,Gir89} and popularized by Hutchinson \cite{Hut89}, so we will call this approximation the \emph{Girard--Hutchinson estimator}.
Tyler Chen has conducted an in-depth study of the history of trace estimation.
His analysis identifies estimators related to \cref{eq:girard-hutchinson} in work from quantum physics dating back as early as the 1920s \cite[\S4.1.1]{Che22}.

\subsection{Example: BKS diagonal estimator} \label{sec:bks}

Another Monte Carlo algorithm for a matrix attribute estimation problem is the BKS diagonal estimator.
Again, we begin with a square matrix $\mat{B}$ and write $\mat{B} = \mat{B} \cdot \Id$. 
Replacing $\Id$ with the stochastic approximation $s^{-1} \sum_{i=1}^s \vec{\omega}_i^{\vphantom{*}} \vec{\omega}_i^* \approx \Id$ yields the diagonal estimator
\begin{equation} \label{eq:bks-formal}
    \hat{\diag}_{\mathrm{BKS}} \coloneqq \diag\left(\mat{B} \cdot \frac{1}{s} \sum_{i=1}^s \vec{\omega}_i^{\vphantom{*}} \vec{\omega}_i^*\right) \approx \diag(\mat{B}).
\end{equation}
This estimator can be simplified by means of the identity
\begin{equation} \label{eq:diagonal-rank-one}
    \diag(\vec{a}\vec{b}^*) = \vec{a} \odot \overline{\vec{b}} \quad \text{for } \vec{a},\vec{b} \in \field^n.
\end{equation}
Here, $\odot$ denotes the entrywise product of vectors, and $\overline{(\cdot)}$ denotes the entrywise complex conjugate.
Applying \cref{eq:diagonal-rank-one} and the linearity of the diagonal map yields a simpler form for the diagonal estimator \cref{eq:bks-formal}:
\begin{equation} \label{eq:bks}
    \hat{\diag}_{\mathrm{BKS}} = \frac{1}{s} \sum_{i=1}^s (\mat{B} \vec{\omega}_i) \odot \overline{\vec{\omega}_i}.
\end{equation}
Estimators similar to \cref{eq:bks} appear to have originally developed by Bekas, Kokiopoulou, and Saad \cite{BKS07}.
In this work, we will call $\hat{\diag}_{\mathrm{BKS}}$ the \emph{BKS diagonal estimator}.
Code for the BKS diagonal estimator appears in \cref{prog:bks}.
We discuss diagonal estimation more thoroughly in \cref{ch:diagonal}.

\myprogram{Bekas--Kokiopoulou--Saad estimator for the diagonal of a matrix.}{Subroutine \texttt{random\_signs} is provided in \cref{prog:random_signs}.}{bks}

\subsection{Example: The Johnson--Lindenstrauss row-norm estimator}

As a final example, we consider the problem of estimating the row norms of a \warn{possibly rectangular} matrix $\mat{B}\in\field^{m\times n}$.
Row-norm estimation is an interesting problem because the row norms are \emph{nonlinear} functions of the matrix.
We shall focus on developing unbiased estimates for the \emph{squared} row norms of $\mat{B}$, which we denote
\begin{equation}
    \srn(\mat{B}) \coloneqq \left( \norm{\mat{B}(i,:)}^2 : 1\le i \le m \right) \in \real^m_+.
\end{equation}

To fashion an estimator for $\srn(\mat{B})$, we first invoke the identity 
\begin{equation*}
    \srn(\mat{B}) = \diag(\mat{B}\mat{B}^*).
\end{equation*}
Now insert a copy of the identity matrix in the middle of the matrix product $\mat{B} \mat{B}^* = \mat{B} \cdot \Id \cdot \mat{B}^*$, and use the stochastic approximation $s^{-1} \sum_{i=1}^s \vec{\omega}_i^{\vphantom{*}} \vec{\omega}_i^* \approx \Id$, resulting in the Monte Carlo estimator
\begin{equation*}
    \hat{\srn}_{\mathrm{JL}} \coloneqq \diag\left( \mat{B} \cdot \frac{1}{s} \sum_{i=1}^s \vec{\omega}_i^{\vphantom{*}} \vec{\omega}_i^* \cdot \mat{B}^*\right) \approx \srn(\mat{B}).
\end{equation*}
This estimator can more conveniently be written as 
\begin{equation} \label{eq:jl-row-norm}
    \hat{\srn}_{\mathrm{JL}} = \frac{1}{s} \sum_{i=1}^s |\mat{B}\vec{\omega}_i|^2,
\end{equation}
where the function $|\cdot|^2$ denotes the squared modulus, evaluated entrywise for a vector input.
We call this estimator the \emph{Johnson--Lindenstrauss row-norm estimator}, as it can be analyzed using the Johnson--Lindenstauss lemma \cite{JL84}.
Code is provided in \cref{prog:jl_rownorm}.
We discuss row-norm estimation more in \cref{ch:row-norm}.

\myprogram{Johnson--Lindenstrauss estimator for the (squared) row norms of a matrix.}{Subroutines \texttt{random\_signs} and \texttt{sqrownorms} re provided in \cref{prog:random_signs,prog:sqrownorms}.}{jl_rownorm}

\subsection{Which isotropic vector to use?} \label{sec:which-isotropic}

Having seen that isotropic random vectors can be used to build Monte Carlo estimators of matrix attributes, we now discuss the choice of which isotropic distribution to use.

Let us first catalog the popular options, many of which have separate definitions over the fields $\field = \real$ and $\field = \complex$ of real and complex numbers:
\begin{itemize}
    \item \textbf{Gaussian.} A standard Gaussian vector $\vec{\omega} \sim \Normal_\field(\vec{0},\Id)$ is isotropic. 
    Its entries are independent and drawn from the real or complex standard Gaussian distribution.
    (Recall that a complex standard Gaussian random variable $g \sim \Normal_\complex(0,1)$ takes the form $g = (g_1 + \iu g_2) / \sqrt{2}$, where $g_1,g_2 \sim \Normal_\real(0,1)$ are independent real standard Gaussians.)
    \item \textbf{Sphere.} A random vector $\vec{\omega} \sim \Unif(\sqrt{n} \sphere(\field^n))$ is isotropic.
    One can be generated by drawing a standard Gaussian vector $\vec{g} \in \Normal_\field(\vec{0},\Id)$ and scaling it to have constant length $\vec{\omega} \coloneqq \sqrt{n} / \norm{\vec{g}} \cdot \vec{g}$.
    \item \textbf{Random signs.} A vector of random signs $\vec{\omega} \sim \Unif \{\pm 1\}^n$ is isotropic.
    This distribution has very low resource requirements, requiring only $n$ independent uniformly random \emph{bits} to generate.
    The random sign distribution is often called the \emph{Rademacher distribution}.
    \item \textbf{Random phases.} An analog of the random sign distribution for the complex field is a vector of random phases $\vec{\omega} \sim \Unif (\unitcircle(\complex)^n)$, defined to be a random vector whose entries are independent and drawn uniformly from the complex unit circle $\unitcircle(\complex) \coloneqq \{ \phi \in \complex : |\phi| = 1 \}$.
    The random phase distribution is also called the \emph{Steinhaus distribution}.
    
    Note that the random signs and random phases distributions can be unified $\vec{\omega} \sim \Unif (\unitcircle(\field)^n)$ by defining the ``unit circle'' of the real numbers to be $\unitcircle(\real) \coloneqq \{ \phi \in \real : |\phi| = 1\} = \{\pm 1\}$.
    \item \textbf{Random coordinate.} Last, one can also generate an isotropic random vector by drawing a random standard basis element $\vec{\omega} \sim \Unif \{\sqrt{n} \evec_i : 1\le i\le n\}$, appropriately rescaled.
\end{itemize}

Which of these options should one use?
To answer this question with some degree of precision, we will focus on which of these random vectors to use for estimating the trace of a real symmetric matrix $\mat{A} \in \real^{n\times n}$ using the Girard--Hutchinson trace estimator, although the same principles apply for other matrix attribute estimation problems and algorithms as well.

Recall that the Girard--Hutchinson estimator for $\mat{A}$ is defined as 
\begin{equation*}
    \hat{\tr}_{\mathrm{GH}} \coloneqq \frac{1}{m} \sum_{i=1}^m \vec{\omega}_i^* \mat{A}\vec{\omega}_i^{\vphantom{*}} \quad \text{for } \vec{\omega}_1,\ldots,\vec{\omega}_m \stackrel{\text{iid}}{\sim} \vec{\omega}.
\end{equation*}
We will compare the accuracy of this estimator with several random isotropic vectors $\vec{\omega}$ by using the mean-squared error $\expect (\hat{\tr}_{\mathrm{GH}} - \tr(\mat{A}))^2$.
Because the terms $\vec{\omega}_i^* \mat{A}\vec{\omega}_i^{\vphantom{*}}$ are iid and unbiased estimators for $\tr(\mat{A})$, the mean-squared error is
\begin{equation} \label{eq:GH-average-to-simple}
    \expect (\hat{\tr}_{\mathrm{GH}} - \tr(\mat{A}))^2 = \frac{1}{m} \Var(\vec{\omega}^*\mat{A}\vec{\omega}).
\end{equation}
In particular, we see that the mean-squared error decays at a $\order(1/m)$ rate regardless of the choice of isotropic test vector $\vec{\omega}$.
The $\order(1/m)$ convergence rate in the mean-squared error is typical of Monte Carlo methods.
The quality of a distribution for trace estimation is quanitified by the prefactor $\Var(\vec{\omega}^*\mat{A}\vec{\omega})$, which depends on the choice of test vector $\vec{\omega}$.

\begin{fact}[Girard--Hutchinson estimator: Variance formulas] \label{fact:gh-variance}
    Let $\mat{A} \in \real^{n\times n}$ be a \warn{real} symmetric matrix with eigenvalues $\lambda_1,\ldots,\lambda_n \in \real$. 
    Denote the mean eigenvalue $\overline{\lambda}\coloneqq n^{-1}\tr(\mat{A}) = n^{-1} \sum_{i=1}^n \lambda_i$ and the mean diagonal element $\overline{a} \coloneqq n^{-1} \sum_{i=1}^n a_{ii}$.
    The following equations give the variance of the basic Girard--Hutchinson estimate $\vec{\omega}^*\mat{A}\vec{\omega}$ for several choices for the isotropic random vector $\vec{\omega} \in \complex^n$:
    \begin{align*}
        \Var(\vec{\omega}^*\mat{A}\vec{\omega}) &= \mathrm{C}_{\field}\norm{\mat{A}}_{\mathrm{F}}^2 &&= \mathrm{C}_{\field}\sum_{i=1}^n \lambda_i^2&& \text{for } \vec{\omega} \sim \Normal_\field(\vec{0},\Id), \\
        \Var(\vec{\omega}^*\mat{A}\vec{\omega}) &= K_{\field}\norm{\mat{A} - \overline{\lambda} \Id}_{\mathrm{F}}^2 &&= K_{\field}\sum_{i=1}^n (\lambda_i - \overline{\lambda})^2&& \text{for } \vec{\omega} \sim \Unif(\sqrt{n}\sphere(\field^n)), \\
        \Var(\vec{\omega}^*\mat{A}\vec{\omega}) &= \mathrm{C}_\field \sum_{i\ne j} a_{ij}^2&& && \text{for } \vec{\omega} \sim \Unif (\unitcircle(\field)^n), \\
        \Var(\vec{\omega}^*\mat{A}\vec{\omega}) &= n^2\sum_{i=1}^n (a_{ii} - \overline{a})^2 && &&\text{for } \vec{\omega} \sim \Unif \{\sqrt{n} \evec_i\}_{1\le i\le n}.
    \end{align*}
    The prefactors are $\mathrm{C}_\real = 2$, $\mathrm{C}_\complex = 1$, $K_\real = 2n/(n+2)$, and $K_\complex = n/(n+1)$.
    These equalities become upper bounds for a real nonsymmetric matrix.
\end{fact}

These variance formulas are standard, and I have collected the simplest proofs of these formulas I know in \cite{Epp23a}.
The random phase distribution $\Unif (\unitcircle(\field)^n)$ and sphere distribution $\Unif(\sqrt{n}\sphere(\field^n))$ are known to be optimal in certain senses, and these optimality results are also discussed in \cite{Epp23a}.
We may draw several conclusions from these formulas:

\myparagraph{The Gaussian distribution is dominated}
The variance for the Gaussian distribution $\vec{\omega} \sim \Normal_\field(\vec{0},\Id)$ is higher than both the sphere and random sign/phase distributions.
Compared to the sphere distribution, the Gaussian distribution has a variance depending on the \textit{aggregate size} of $\mat{A}$'s eigenvalues, whereas the sphere distribution depends only on the \emph{spread} of $\mat{A}$'s eigenvalues. 
This size/spread distinction can make a big difference for a matrix $\mat{A}$ with large eigenvalues that are tightly clustered.
Compared to the random phase distribution, the variance for the Gaussian distribution depends on Frobenius norm $\norm{\mat{A}}_{\mathrm{F}}^2 = \sum_{i,j} a_{ij}^2$, which reflects the magnitude of all of $\mat{A}$'s entries.
By contrast, the variance of the random sign/phase distribution depends on $\sum_{i\ne j} a_{ij}^2$, which omits the contribution from $\mat{A}$'s diagonal entries.
This property makes the random sign/phase distribution especially effective for matrices with a heavy diagonal.
I describe the case against Gaussians for stochastic trace estimation in more detail and provide numerical evidence in \cite{Epp24b}.

\myparagraph{Sphere vs.\ signs/phases}
The sphere and sign/phase distributions both dominate the Gaussian distribution, but which should be preferred?
Ultimately, both choices of distribution are effective, and we will use both in this thesis.
The random sign/phase distribution has the benefit that it ignores the influence of the diagonal of $\mat{A}$, making it effective for matrices with a heavy diagonal.
The random sphere distribution, however, is known to be (minimax) optimal for matrices $\mat{A}$ drawn from unitarily invariant families, such as the class of all symmetric matrices with $\norm{\mat{A}}_{\mathrm{F}} \le 1$.
This optimality result appears to have first been discovered in unpublished work of Richard Kueng \cite[Prob.~1.23]{Tro20a}; see \cite{Epp23a} for a proof.

\myparagraph{Real vs.\ complex}
The variance formulas in \cref{fact:gh-variance} show that using complex-valued test vectors $\vec{\omega}$ results in a variance about half that of their real-valued counterparts.
However, when applied to real matrices $\mat{A}$, this benefit is typically offset by the increased computational cost of using complex arithmetic.
Further, existing codes for evaluating the matvec operation might not be compatible with complex data.
As such, the reduced variance of complex-valued test vectors is usually not worth the additional cost when the matrix $\mat{A}$ is real-valued.
For a complex-valued matrix $\mat{B} \in \complex^{n\times n}$, however, it is generally preferable to use complex-valued test vectors $\vec{\omega}$.

\myparagraph{Coordinate sampling is often dangerous}
The variance for trace estimation with random coordinate vectors is often much higher than alternate approaches.
This can be seen by a simple ``back of the envelope'' computation: Consider a matrix $\mat{A}$ where all the entries are roughly of unit magnitude $|a_{ij}| \approx 1$, with a similar number of positive and negative entries, so that the average of diagonal entries is small: $\overline{a} \approx 0$.
Then the variance of the Girard--Hutchinson estimator with random coordinate sampling has scaling $\sim n^3$, whereas the other isotropic distributions in \cref{fact:gh-variance} have a variance of $\sim n^2$.
Thus, in view of \cref{eq:GH-average-to-simple}, the coordinate sampling distribution may require up to $n$ matvecs to achieve accuracy comparable to what other estimators can achieve with a single matvec.
This is a dismal state of affairs, particular since the trace estimation problem can be solved exactly in $n$ matvecs.
Random coordinate vectors do have their uses, particularly in computational models more restrictive than the matvec model \cite{BKM22a}, but they are best avoided except in special situations.

\subsection{Resphering: Improved Monte Carlo for rank-deficient matrices} \label{sec:resphering}

So far, we have designed Monte Carlo methods for matrix attribute estimation based on introducing the identity matrix $\mat{B} = \mat{B} \cdot \Id$ and replacing $\Id$ with an unbiased stochastic approximation.
As we will see shortly, we will have many occasions to apply Monte Carlo approximations to a rank-deficient matrix $\mat{B}$ with a known (right or left) nullspace.
In such cases, we can use \emph{resphering} to obtain lower variance Monte Carlo estimators.
Here is the main definition:

\begin{definition}[Isotropic random vector on a subspace]
    Let $\set{U} \subseteq \field^n$ be a subspace over the real or complex numbers, and let $\mat{\Pi}$ denote the orthoprojector on $\set{U}$.
    A random vector $\vec{\nu}$ is said to be \emph{isotropic on $\set{U}$} if $\expect[\vec{\nu}\vec{\nu}^*] = \mat{\Pi}$.
\end{definition}

For any subspace $\set{U}$ of dimension $s$, a uniformly random vector $\vec{\nu} \sim \Unif(\sqrt{s}\sphere(\set{U}))$ from the sphere of radius $\sqrt{s}$ in $\set{U}$ is isotropic on $\set{U}$.
This will be our only example of an isotropic vector on a subspace in this thesis.

To see how isotropic random vectors on a subspace can be employed, let $\mat{B} \in \complex^{n\times n}$ be a square matrix and suppose we have access to a matrix $\mat{Q} \in \complex^{n\times k}$ with orthonormal columns for which $\mat{B}\mat{Q} = \mat{0}$.
To design a trace estimator for $\mat{B}$, write 
\begin{equation*}
    \mat{B} = \mat{B} \cdot \mat{\Pi} \quad \text{for } \mat{\Pi} \coloneqq \Id - \mat{Q}\mat{Q}^*.
\end{equation*}
To approximate $\mat{\Pi}$, generate isotropic random vectors $\vec{\nu}_1,\ldots,\vec{\nu}_s$ on the subspace $\range(\mat{Q})^\perp$ and introduce the stochastic approximation $s^{-1} \sum_{i=1}^s \vec{\nu}_i^{\vphantom{*}}\vec{\nu}_i^*$.
This leads to the \emph{resphered Girard--Hutchinson trace estimator}
\begin{equation*}
    \hat{\tr}_{\mathrm{SGH}} \coloneqq \tr \left( \mat{B} \cdot \frac{1}{s} \sum_{i=1}^s \vec{\nu}_i^{\vphantom{*}}\vec{\nu}_i^*\right) = \frac{1}{s} \sum_{i=1}^s \vec{\nu}_i^* \left(\mat{B}\vec{\nu}_i^{\vphantom{*}}\right).
\end{equation*}

To generate isotropic random vectors 
\begin{equation*}
    \vec{\nu} \sim \Unif\bigl(\sqrt{n-k}\cdot \sphere(\range(\mat{Q})^\perp)\bigr)    
\end{equation*}
on the subspace $\range(\mat{Q})^\perp$ proceed as follows: Draw a random vector $\vec{\omega}$ from the normal or sphere distribution, orthogonalize it against $\mat{Q}$, and rescale to the proper norm
\begin{equation} \label{eq:resphered}
    \vec{\nu} \coloneqq \frac{\sqrt{n-k}}{\norm{(\Id - \mat{Q}\mat{Q}^*)\vec{\omega}}} \cdot (\Id - \mat{Q}\mat{Q}^*)\vec{\omega} = \frac{\sqrt{n-k}}{\norm{\vec{\omega}- \mat{Q}(\mat{Q}^*\vec{\omega})}} \cdot (\vec{\omega}- \mat{Q}(\mat{Q}^*\vec{\omega})).
\end{equation}
We call the process of replacing $\vec{\omega}$ by $\vec{\nu}$ in this way \emph{resphering}.


When applied to a rank-deficient matrix $\mat{B}$ whose nonzero eigenvalues are clustered, the resphered Girard--Hutchinson estimator can lead to significantly lower variance than the standard Girard--Hutchinson estimator.
Here is an example result.

\begin{corollary}[Resphered Girard--Hutchinson estimator: Variance] \label{cor:resphere-gh-variance}
    Let $\mat{A} \in \real^{n\times n}$ be a symmetric matrix with at most $r$ nonzero eigenvalues $\lambda_1,\ldots,\lambda_r$, let $\overline{\lambda} \coloneqq r^{-1} \sum_{i=1}^r \lambda_i$ denote their average, and let $\vec{\nu} \sim \Unif(\sqrt{r} \cdot \sphere(\range(\mat{A})))$ be a test vector.
    Then the resphered Girard--Hutchinson estimator $\vec{\nu}^*\mat{A}\vec{\nu}$ has variance
    \begin{equation*}
        \Var(\vec{\nu}^*\mat{A}\vec{\nu}) = \frac{2r}{r+2} \sum_{i=1}^r (\lambda_i - \overline{\lambda})^2.
    \end{equation*}
    Excepting the trivial case when $\lambda_i = 0$ for all $i$, the variance $\Var(\vec{\nu}^*\mat{A}\vec{\nu})$ is always strictly smaller than for the plain Girard--Hutchinson estimator $\Var(\vec{\omega}^*\mat{A}\vec{\omega})$ with $\vec{\omega} \sim \Unif(\sqrt{n} \cdot \sphere(\real^n))$.
\end{corollary}

The technique now called resphering was first introduced in \cite{ETW24}, where it appeared under the name ``normalization.''
In this thesis, I adopt the more descriptive term ``resphering,'' a name suggested to me by Joel Tropp.

\section{Low-rank approximation} \label{sec:lra}

Low-rank approximation provides another paradigm for solving matrix attribute estimation problems.
As we saw in \cref{ch:lra}, we can cheaply compute near-optimal low-rank approximations to a general matrix $\mat{B}$ using the randomized SVD (\cref{sec:rsvd}) or to a psd matrix $\mat{A}$ using Nystr\"om approximation (\cref{sec:nystrom}).

Low-rank approximation gives a natural approach to any matrix attribute estimation problem.
To estimate an attribute $Q(\mat{B})$ of a matrix $\mat{B} \in \complex^{m\times n}$, simply replace $\mat{B}$ by a low-rank approximation $\Bhat$ and use the $Q(\Bhat)$ as an estimator for $Q(\mat{B})$.
This use of low-rank approximation in this way was proposed by Saibaba, Alexanderian, and Ipsen for trace estimation  \cite{SAI17}.

The quality of the approximation $Q(\Bhat)\approx Q(\mat{B})$ is dictated by the quality of the low-rank approximation $\mat{B} \approx \Bhat$ and, consequently, the rate of singular value decay in the matrix $\mat{B}$ (recall \cref{fact:rsvd-error,fact:nystrom-error}).
For matrices with rapid singular value decay, low-rank approximation-based estimators can produce highly accurate results.
On the other hand, for matrices with slow singular value decay, plain low-rank approximation based estimators are often wholly inaccurate.
This inconsistent performance makes pure low-rank approximation based estimators only situationally useful.

\section{Combining Monte Carlo and low-rank approximation} \label{sec:mc-plus-lra}

The Monte Carlo and low-rank approximation approaches can be combined to achieve better and more consistent results than either approach yields by itself.

Let us illustrate by deriving a simplified version of the \HutchPP algorithm for trace estimation \cite{MMMW21a}.
Suppose we wish to estimate the trace of a matrix $\mat{B} \in \complex^{n\times n}$ and are given a fixed budget of $s$ matvecs to accomplish the task.
Assume, for simplicity, that $s$ is divisible by $3$.
Begin by running the randomized SVD (\cref{sec:rsvd}) to compute a matrix $\mat{Q} \in \complex^{n\times (s/3)}$, which defines a rank-$(s/3)$ approximation $\Bhat = \mat{Q}\mat{Q}^*\mat{B}$ to $\mat{B}$.
By linearity, we may decompose the trace of $\mat{B}$:
\begin{equation*}
    \tr(\mat{B}) = \tr(\Bhat) + \tr(\mat{B} - \Bhat) = \tr(\mat{Q}\mat{Q}^*\mat{B}) + \tr((\Id-\mat{Q}\mat{Q}^*)\mat{B}).
\end{equation*}
The first term, $\tr(\Bhat) = \tr(\mat{Q}\mat{Q}^*\mat{B})$, can be computed exactly by forming $\mat{Q}^*\mat{B}$.
To estimate the second term, we employ a Monte Carlo method, specifically the Girard--Hutchinson estimator:
\begin{equation*}
    \tr(\mat{B} - \Bhat) \approx \frac{1}{s/3} \sum_{i=1}^{s/3} \vec{\gamma}_i^* \left( \left(\Id - \mat{Q}\mat{Q}^*\right)(\mat{B}\vec{\gamma}_i^{\vphantom{*}})\right)
\end{equation*}
Here, $\vec{\gamma}_1,\ldots,\vec{\gamma}_{s/3}$ denote freshly generated isotropic random vectors, independent of each other and the random test matrix used to execute the randomized SVD.
Combining the exact computation of the first term and the stochastic approximation of the second yields the (simplified) \HutchPP estimator
\begin{equation} \label{eq:simple-hutchpp}
    \hat{\tr}_{\mathrm{SH}\texttt{++}} \coloneqq \tr(\mat{Q}(\mat{Q}^*\mat{B})) + \frac{1}{s/3} \sum_{i=1}^{s/3} \vec{\gamma}_i^* \left( \left(\Id - \mat{Q}\mat{Q}^*\right)(\mat{B}\vec{\gamma}_i^{\vphantom{*}})\right) \quad \text{for }\mat{Q} = \orth(\mat{B}\mat{\Omega}).
\end{equation}
The computational cost of the simplified \HutchPP estimator is $s$ matvecs ($2s/3$ with $\mat{B}$ and $s/3$ with $\mat{B}^*$) plus $\order(s^2n)$ additional arithmetic operations.

The full \HutchPP estimator contains two optimizations over the simplified version \cref{eq:simple-hutchpp}:

\begin{enumerate}
    \item In some contexts, matvecs with $\mat{B}^*$ are expensive or are entirely unavailable.
    To avoid matvecs with $\mat{B}^*$, we use the cyclic property of the trace to write
    \begin{equation*}
        \tr(\Bhat) = \tr(\mat{Q}(\mat{Q}^*\mat{B})) = \tr(\mat{Q}^*(\mat{B}\mat{Q})).
    \end{equation*}
    The expression $\tr(\mat{Q}^*(\mat{B}\mat{Q}))$ can be evaluated using only matvecs with $\mat{B}$, removing matvecs with $\mat{B}^*$ from the \HutchPP algorithm entirely.
    \item The matrix $\Id - \mat{Q}\mat{Q}^*$ is an orthoprojector and thus satisfies $\Id - \mat{Q}\mat{Q}^* = (\Id - \mat{Q}\mat{Q}^*)^2$.
    Therefore, the residual trace can be written more symmetrically as
    \begin{equation*}
        \tr(\mat{B} - \Bhat) = \tr((\Id - \mat{Q}\mat{Q}^*) \mat{B}) = \tr((\Id - \mat{Q}\mat{Q}^*)^2 \mat{B}) =  \tr((\Id - \mat{Q}\mat{Q}^*) \mat{B}(\Id - \mat{Q}\mat{Q}^*)).
    \end{equation*}
    The symmetrically projected matrix $(\Id - \mat{Q}\mat{Q}^*) \mat{B}(\Id - \mat{Q}\mat{Q}^*)$ always has a smaller Frobenius norm than the one-sided projection $(\Id - \mat{Q}\mat{Q}^*) \mat{B}$, so we expect a smaller error in applying the Girard--Hutchinson estimator to the former matrix rather than the latter (cf.\ \cref{fact:gh-variance}).
\end{enumerate}
Combining these two optimizations yields the standard \HutchPP estimator
\begin{equation} \label{eq:hutchpp}
    \hutchppest  \coloneqq \tr(\mat{Q}^*(\mat{B}\mat{Q})) + \frac{1}{s/3} \sum_{i=1}^{s/3} \vec{\gamma}_i^* \left( \left(\Id - \mat{Q}\mat{Q}^*\right)\left(\mat{B}\left( \left(\Id - \mat{Q}\mat{Q}^*\right)\vec{\gamma}_i^{\vphantom{*}}\right)\right)\right),
\end{equation}
where $\mat{Q} = \orth(\mat{B}\mat{\Omega})$.
Code is provided in \cref{prog:hutchpp}.

\myprogram{\HutchPP algorithm for trace estimation.}{The \texttt{random\_signs} subroutine is defined in \cref{prog:random_signs}.}{hutchpp}

The combination of low-rank approximation with Monte Carlo for trace estimation and related problems was explored prior to the original \HutchPP paper.
Notable examples include the works of Gambhir, Stathopoulos, and Originos \cite{GSO17} and Lin \cite{Lin17}.
Meyer, Musco, Musco, and Woodruff \cite{MMMW21a} crystallized these ideas in the \HutchPP algorithm, provided mathematical analysis, and established lower bounds on the best-possible accuracy for any trace estimation algorithm.

The \HutchPP algorithm illustrates the use of \emph{variance reduction} in Monte Carlo methods \cite[\S2.3]{Liu04}.
Specialized to trace estimation, the key idea of variance reduction is to choose a matrix $\Bhat$, called a \emph{control variate}, that closely approximates $\mat{B}$ and whose trace can be computed exactly.
Instead of applying a Monte Carlo estimator to estimate $\tr(\mat{B})$ directly, we instead estimate the trace of the residual $\mat{B} - \Bhat$ and add the result to $\tr(\Bhat)$.
This strategy typically yields an estimator with lower variance.
In principle, this variance reduction strategy can be applied using any type of approximation $\Bhat\approx\mat{B}$, but low-rank approximations have proven the most effective for trace estimation so far.

\subsection{Theoretical analysis of \HutchPP}

Error bounds for \HutchPP can be derived by combining error bounds for the randomized SVD (\cref{eq:rsvd-fro}) with variance bounds for the Girard--Hutchinson estimator (\cref{fact:gh-variance}).
We will state and prove such a bound here, as it will serve as a useful comparison for the error bounds for the \XTrace and \XNysTrace estimators developed in the next section; see \cref{sec:trace-error-bounds}.

\begin{theorem}[\HutchPP: mean-squared error] \label{thm:hutchpp-mse}
    Let $\mat{B} \in \real^{n\times n}$ be a \warn{real} matrix and let $\hutchppest $ be the \HutchPP estimator \cref{eq:hutchpp-error} with a \warn{real} standard Gaussian test matrix $\mat{\Omega} \in \real^{n\times (s/3)}$ and iid $\vec{\gamma}_1,\ldots,\vec{\gamma}_{s/3}\in\real^n$ drawn from any of the \warn{real} isotropic distributions from \cref{fact:gh-variance}, except the random coordinate distribution.
    Then
    \begin{equation} \label{eq:hutchpp-error}
        \expect \left( \hutchppest - \tr(\mat{B}) \right)^2 \le \frac{6}{s} \cdot \min_{r \le s/3-2} \frac{s-3}{s-3r-3}\cdot \norm{\mat{B} - \lowrank{\mat{B}}_r}_{\mathrm{F}}^2.
    \end{equation}
\end{theorem}

\begin{proof}
    Using the \emph{reduced matrix}
    \begin{equation*}
        \mat{B}_{\mathrm{red}} = (\Id-\mat{Q}\mat{Q}^*)\mat{B}(\Id-\mat{Q}\mat{Q}^*),
    \end{equation*}
    the error of \HutchPP may be expressed as
    \begin{equation} \label{eq:hutchpp-error-formula}
        \hutchppest - \tr(\mat{B}) = \frac{1}{s/3} \sum_{i=1}^{s/3} \vec{\omega}_i^* \mat{B}_{\mathrm{red}}\vec{\omega}_i^{\vphantom{*}} - \tr(\mat{B}_{\mathrm{red}}) = \frac{1}{s/3} \sum_{i=1}^{s/3} (\vec{\omega}_i^* \mat{B}_{\mathrm{red}}\vec{\omega}_i^{\vphantom{*}} - \tr(\mat{B}_{\mathrm{red}})).
    \end{equation}
    The right-hand side of \cref{eq:hutchpp-error-formula} is a sum of mean-zero random variables that are iid conditional on $\mat{B}_{\mathrm{red}}$.
    Therefore,
    \begin{equation*}
        \expect \left[(\hutchppest - \tr(\mat{B}))^2 \, \middle| \, \mat{B}_{\mathrm{red}} \right] = \frac{1}{s/3} \Var\left(\vec{\omega}_1^*\mat{B}_{\mathrm{red}}\vec{\omega}_1^{\vphantom{*}} \, \middle| \, \mat{B}_{\mathrm{red}}\right).
    \end{equation*}
    Applying \cref{fact:gh-variance}, we obtain
    \begin{equation*}
        \expect \left[(\hutchppest - \tr(\mat{B}))^2 \, \middle| \, \mat{B}_{\mathrm{red}} \right] \le \frac{6}{s} \cdot \norm{\mat{B}_{\mathrm{red}}}_{\mathrm{F}}^2.
    \end{equation*}

    Multiplying by an orthoprojector can only reduce the Frobenius norm.
    Therefore,
    \begin{equation*}
        \norm{\mat{B}_{\mathrm{red}}}_{\mathrm{F}} \le \norm{(\Id - \mat{Q}\mat{Q}^*)\mat{B}}_{\mathrm{F}}.
    \end{equation*}
    Combining the two previous displays and taking the expectation, we obtain
    \begin{equation*}
        \expect (\hutchppest - \tr(\mat{B}))^2 \le \frac{6}{s} \cdot \expect \norm{(\Id - \mat{Q}\mat{Q}^*)\mat{B}}_{\mathrm{F}}^2.
    \end{equation*}
    Invoking the randomized SVD error bound \cref{eq:rsvd-fro} establishes the desired result.
\end{proof}

In \cref{thm:hutchpp-mse}, we assumed the matrix $\mat{B}$ was real to use \cref{fact:gh-variance}.
Versions of this result are easy to derive for complex vectors.

The error bound \cref{eq:hutchpp-error} shows that the mean-squared error of \HutchPP is proportional to $1/s$ times the (squared) Frobenius norm error of the best rank-$r$ approximation, where $r\approx s/3$.
A bound similar to \cref{thm:hutchpp-mse} appears in \cite[Thm.~1.1]{ETW24}; see also \cite{Mey21} and \cite[\S5]{MMMW21} in the \warn{arXiv version} of the \HutchPP paper.

The bound \cref{eq:hutchpp-error} describes the mean-squared error of the \HutchPP method in practice, but it can be a bit imposing.
Therefore, it can be informative to derive simplified versions of the bound.
Choosing $r = s/6-1$ in the minimum \cref{eq:hutchpp-error} yields
\begin{equation} \label{eq:hutchpp-error-1}
    \expect \left( \hutchppest  - \tr(\mat{B}) \right)^2\le \frac{12}{s} \cdot \norm{\mat{B} - \lowrank{\mat{B}}_{s/6-1}}_{\mathrm{F}}^2.
\end{equation}
To further simplify, we can employ a crude bound on the best rank-$r$ approximation error (see, e.g., \cite[Lem.~7]{GSTV07}, \cite[Lem.~13]{MMMW21a}, and \cite[Fact~5.5]{ETW24} for versons of this result).

\begin{fact}[Rank-$r$ approximation error] \label{fact:l1l2}
    For any matrix $\mat{B}\in \field^{m\times n}$ and $r\ge 1$, 
    \begin{equation*}
        \norm{\mat{B} - \lowrank{\mat{B}}_r} \le \frac{\norm{\mat{B}}_*}{r+1}\quad \text{and} \quad \norm{\mat{B} - \lowrank{\mat{B}}_r}_{\mathrm{F}} \le \frac{\norm{\mat{B}}_*}{2\sqrt{r}}.   
    \end{equation*}
\end{fact}

Using \cref{fact:l1l2}, we can further bound \cref{eq:hutchpp-error-1} as
\begin{equation} \label{eq:hutchpp-error-2}
    \expect \left( \hutchppest  - \tr(\mat{B}) \right)^2 \le \frac{72}{s^2} \cdot \norm{\mat{B}}_*^2.
\end{equation}
This result shows that the mean-squared error of \HutchPP is dominated by a quantity that decays at a $\order(1/s^2)$ rate, improving on the Monte Carlo $\order(1/s)$ rate \cref{eq:GH-average-to-simple} of the Girard--Hutchinson estimator.
However, the reader should be aware that this ``$\order(1/s^2)$ convergence rate'' for \HutchPP is easy to misinterpret; see \cref{sec:interpreting-hutchpp-bounds}.

As a final simplification, suppose we apply \HutchPP to a psd matrix $\mat{A}$. 
In this case, the trace norm and the trace are the same, and \cref{eq:hutchpp-error-2} implies that
\begin{equation} \label{eq:hutchpp-error-3}
    \left( \expect \left( \hutchppest  - \tr(\mat{A}) \right)^2\right)^{1/2} \le \varepsilon \tr(\mat{A}) \quad \text{when } s \ge \frac{6\sqrt{2}}{\varepsilon}.
\end{equation}
We obtain a root-mean-squared error of $\varepsilon \tr(\mat{A})$ when the number of matvecs $s$ is $s = \order(1/\varepsilon)$.
Meyer et al.\ show that the $s = \order(1/\varepsilon)$ parameter complexity is optimal for trace estimation; in particular, no algorithm achieves relative error $\varepsilon$ using $s = \order(1/\varepsilon^{0.999})$ matvecs for every input matrix $\mat{A}$ \cite[\S4]{MMMW21a} (see also \cite[\S2.3.5]{Mey24}).

\begin{remark}[High probability bounds]
    It is straightforward to prove error bounds for \HutchPP that control the error $|\hutchppest  - \tr(\mat{B})|$ with high probability by combining high-probability error bounds for the randomized SVD with the Hanson--Wright inequality.
    For the former, see relevant results in \cite{HMT11,MM20,TW23}.
    For the latter, see \cite{RV13,Epp22}.
    Using this approach, Meyer et al.\ established the following result \cite[Thm.~1.1]{MMMW21a}:

    \begin{fact}[\HutchPP: High probability error bound] \label{fact:hpp-hp}
        When applied to a \warn{real} psd matrix $\mat{A}$, \HutchPP with random sign vectors achieves the guarantee 
        \begin{equation} \label{eq:hutchpp-whp}
            |\hutchppest - \tr(\mat{A})| \le \varepsilon \tr(\mat{A}) \quad \text{with probability at least } 1-\delta
        \end{equation}
        using $s = \order(\varepsilon^{-1}\cdot\sqrt{\log(1/\delta)} + \log(1/\delta))$ matvecs.
    \end{fact}

    This result establishes that the probability of failing to produce a trace approximation of relative error $\varepsilon$ decreases exponentially in the number of matvecs $s$.
\end{remark}

\subsection{Improving \HutchPP using resphering}

\myprogram{\HutchPP algorithm with resphering for trace estimation.}{The \texttt{sqcolnorms} subroutine is defined in \cref{prog:sqcolnorms}.}{hutchpp_resphere}

In passing, let us observe that we can enhance the \HutchPP algorithm by using resphering, introduced in \cref{sec:resphering}.
Simply use the resphered Girard--Hutchinson estimator rather than the Girard--Hutchinson estimator to estimate $\tr((\Id - \mat{Q}\mat{Q}^*)\mat{B}(\Id - \mat{Q}\mat{Q}^*))$.
Code is provided in \cref{prog:hutchpp_resphere}.

\begin{figure}
    \centering
    \includegraphics[width=0.7\linewidth]{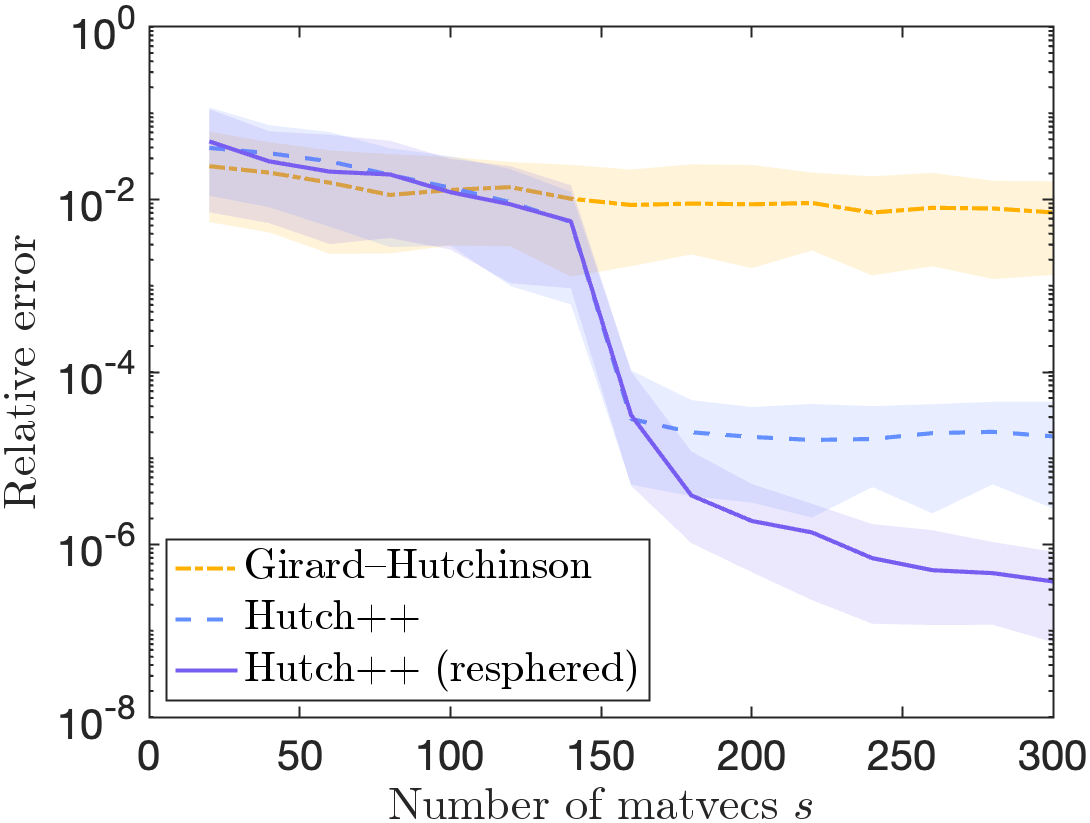}
    \caption[Comparison of Girard--Hutchinson and \HutchPP estimators with and without resphering on a matrix with a ``step'' spectrum]{Comparison of Girard--Hutchinson and \HutchPP estimators with and without resphering on the \texttt{step} matrix \cref{eq:step-hutchpp}.
    Lines show median of 100 trials, and error bars show 10\% and 90\% quantiles.}
    \label{fig:hutchpp-resphere}
\end{figure}

The benefits of resphering are demonstrated in \cref{fig:hutchpp-resphere}, which compares the Girard--Hutchinson estimator to the \HutchPP algorithm with and without resphering.
I use input matrix
\begin{equation} \label{eq:step-hutchpp}
    \texttt{step} \coloneqq \mat{U} \Diag(\underbrace{1,\ldots,1}_{\text{50 times}},\underbrace{10^{-3},\ldots,10^{-3}}_{\text{950 times}}) \mat{U}^* \in \real^{10^3\times 10^3}.
\end{equation}
Here, $\mat{U}$ is a Haar random orthogonal matrix.
The matrix \texttt{step} has $k=50$ large eigenvalues, after which the spectrum is flat.
It takes the \HutchPP methods roughly $s \approx 3k = 150$ matvecs to compute a low-rank approximation capturing these dominant eigenvalues, explaining the substantial drop in the error of the \HutchPP methods at $s=150$ matvecs.
With resphering, the error drops two more orders of magnitude as $s$ is increased beyond $150$, owing to the relatively flat distribution of \warn{nonzero} eigenvalues of the projected matrix $(\Id - \mat{Q}\mat{Q}^*)\texttt{step}(\Id - \mat{Q}\mat{Q})$.
For $s=300$, \HutchPP with resphering is 50$\times$ more accurate than \HutchPP without resphering and 1800$\times$ more accurate than the Girard--Hutchinson estimator.

\subsection{Weaknesses of \HutchPP} \label{sec:hutchpp-weaknesses}

The \HutchPP algorithm is an elegant approach for trace estimation, and it achieves excellent accuracy compared to the Girard--Hutchinson estimator when applied to matrices with rapidly decaying singular values.
Even so, we recognize opportunities for improvement, which will be realized by the \XTrace and \XNysTrace estimators in the next chapter.

\myparagraph{Apportionment}
In its standard form, the \HutchPP algorithm uses a fixed division of $2s/3$ matvecs for low-rank approximation and $s/3$ matvecs for residual trace estimation.
For a given matrix, this apportionment need not be optimal.
For matrices with rapidly decaying singular values, the ideal strategy would be to put \emph{all} matvecs into low-rank approximation; for matrices with slowly decreasing singular values, the opposite strategy is warranted.
The apportionment problem for \HutchPP can be addressed with an adaptive scheme that apportions matvecs between low-rank approximation and residual trace estimation, as is done in \cite{PCK22}.
The \XTrace and \XNysTrace estimators introduced in next section solve the apportionment problem by redesigning the estimator to exploit a leave-one-out\index{leave-one-out randomized algorithm} approach.

\myparagraph{Exchangeability}
The second weakness of \HutchPP is more conceptual.
We have the following property of the optimal (i.e., minimum-variance) trace estimator.

\actionbox{\textbf{Exchangeability principle.} Consider the class of estimators for $\tr(\mat{A})$ computed using matvecs $\mat{A}\vec{\omega}_1,\ldots,\mat{A}\vec{\omega}_k$ with \warn{iid} random vectors $\vec{\omega}_1,\ldots,\vec{\omega}_k$. Among these methods, the minimum-variance unbiased estimator for $\tr(\mat{A})$ is a \emph{permutation-invariant} function of $\vec{\omega}_1,\ldots,\vec{\omega}_k$.}

This principle is due to Halmos \cite{Hal46}, and it is easy to confirm.
Indeed, the variance of any estimator $\hat{\tr}(\vec{\omega}_1,\ldots,\vec{\omega}_k)$ is reduced by averaging over all permutations of the test vectors $\vec{\omega}_1,\ldots,\vec{\omega}_k$:
\begin{align*}
    &\Var\left( \frac{1}{k!} \sum_{\sigma \in \set{S}_k} \hat{\tr}(\vec{\omega}_{\sigma(1)},\ldots,\vec{\omega}_{\sigma(k)}) \right) \\ &\qquad= \frac{1}{(k!)^2} \sum_{\sigma,\sigma' \in \set{S}_k}  \Cov\left( \hat{\tr}(\vec{\omega}_{\sigma(1)},\ldots,\vec{\omega}_{\sigma(k)}),\hat{\tr}(\vec{\omega}_{\sigma'(1)},\ldots,\vec{\omega}_{\sigma'(k)}) \right) \\
    &\qquad\le \frac{1}{(k!)^2} \sum_{\sigma,\sigma' \in \set{S}_k}  \left[\Var\left( \hat{\tr}(\vec{\omega}_{\sigma(1)},\ldots,\vec{\omega}_{\sigma(k)})\right)\Var\left(\hat{\tr}(\vec{\omega}_{\sigma'(1)},\ldots,\vec{\omega}_{\sigma'(k)}) \right)\right]^{1/2} \\
    &\qquad= \Var( \hat{\tr}(\vec{\omega}_{1},\ldots,\vec{\omega}_{k})).
\end{align*}
The first identity is sesquilinearity of the covariance, the inequality is Cauchy--Schwarz for the variance, and the second identity follows from the observation that $\hat{\tr}(\vec{\omega}_{\sigma(1)},\ldots,\vec{\omega}_{\sigma(k)})$ has the same distribution for every permutation $\sigma$.

The \HutchPP estimator is not exchangeable: Half the test vectors are used for low-rank approximation, while the other half are used for residual trace estimation.
As such, there must be an exchangeable estimator with lower variance than \HutchPP.
One way we could fashion such an estimator is by averaging the value of the \HutchPP estimator over all possible splits of test vectors between these two tasks, as in the derivation above.
Unfortunately, computing this symmetrized \HutchPP estimator is computationally infeasible because there are $\binom{2s/3}{s/3} \ge 2^{s/3}$ divisions of $2s/3$ test vectors into two equal groups.
A new, different approach is needed to obtain an exchangeable, variance-reduced trace estimator.
The \XTrace estimator will remedy this shortcoming of \HutchPP by using a leave-one-out\index{leave-one-out randomized algorithm} design, which is exchangeable by construction.

\chapter{The leave-one-out approach and application to trace estimation} \label{ch:loo}

\epigraph{An idea which can be used only once is a trick. If one can use it more than once it becomes a method.}{George P\'olya and Gabor Szeg\H{o}, \textit{Problems and Theorems in Analysis I} \cite[p.~VIII]{PS98}}

In the previous chapter, we developed our basic tools---randomized Monte Carlo and low-rank approximations for matrices---and used them to estimate matrix attributes.
This discussion culminated with the \HutchPP algorithm, which showed how combining Monte Carlo estimation and low-rank approximation can lead to more accurate trace estimates than either approach individually.

In this chapter, we will present the leave-one-out\index{leave-one-out randomized algorithm} approach for matrix attribute estimation, a flexibility methodology that combines Monte Carlo and low-rank approximation in a way that squeezes as much information possible out of every matvec.
When the leave-one-out\index{leave-one-out randomized algorithm} approach is applied to trace estimation, it results in the \XTrace and \XNysTrace algorithms.
Each method is an exchangeable estimator that uses every matvec both for low-rank approximation and for Monte Carlo estimation.
This strategy ameliorates the weaknesses with \HutchPP identified at the end of last chapter.

This chapter presents a tutorial-style introduction to leave-one out randomized matrix algorithm design, focusing on developing trace estimators for general and psd matrices.
The basic idea is simple, but the formulas appearing in the final algorithms can be complicated.
This chapter will attempt to demystify these expressions and to provide a clear example of how to derive a leave-one-out\index{leave-one-out randomized algorithm} randomized matrix algorithm.
Subsequent chapters will present several additional applications of the leave-one-out\index{leave-one-out randomized algorithm} approach to matrix attribute estimation.

\myparagraph{Sources}
Both the leave-one-out\index{leave-one-out randomized algorithm} approach and the \XTrace and \XNysTrace algorithms were developed in the paper:

\fullcite{ETW24}.

This chapter refines the original paper \cite{ETW24} by developing new implementations of the \XNysTrace algorithm with processing costs that are up to $7\times$ faster than the implementations given in the original paper.

\myparagraph{Outline}
\Cref{sec:xtrace-derivation} introduces the leave-one-out\index{leave-one-out randomized algorithm} approach to matrix attribute estimation and uses it to derive the \XTrace algorithm for trace estimation.
To implement this algorithm efficiently, we require a downdating formula for the randomized SVD, which is developed in \cref{sec:rsvd-downdate}.
We devise an efficient \XTrace implementation using this formula in \cref{sec:xtrace-implementation}.
\Cref{sec:xnystrace-derivation} presents \XNysTrace, an improved version of \XTrace for psd matrices.
\Cref{sec:nystrom-downdating,sec:xnystrace-implementation} discusses implementation of \XNysTrace.
\cref{sec:trace-experiments} contains experimental comparison of \XTrace and \XNysTrace to \HutchPP and the Girard--Hutchinson estimators, \cref{sec:xtrace-resphere} introduces resphered versions of \XTrace and \XNysTrace, and \cref{sec:estrada} presents an application of trace estimators to computing the Estrada index of a network.
We conclude in \cref{sec:loo-summary} by summarizing the leave-one-out\index{leave-one-out randomized algorithm} approach to matrix attribute estimation.

\section{\XTrace: The leave-one-out approach} \label{sec:xtrace-derivation}

The leave-one-out\index{leave-one-out randomized algorithm} approach to matrix attribute estimation consists of five steps:
\begin{enumerate}
    \item Compute a low-rank approximation to the input matrix by multiplying it against a collection of test vectors.
    \item Decompose the quantity of interest into known piece depending on the low-rank approximation plus a residual.
    \item Construct a Monte Carlo estimate of the residual using a single random vector.
    \item Downdate the low-rank approximation by recomputing it with a test vector removed, and use the left-out test vector as the random vector in step 3.
    \item Average the estimator from step 4 over all choices of vectors to leave out.
\end{enumerate}
At present, these steps are fairly abstract.
To make this program more complete, we shall derive \XTrace, a leave-one-out\index{leave-one-out randomized algorithm} algorithm for estimating the trace of a general square matrix $\mat{B}$.

\myparagraph{Step 1: Compute a low-rank approximation to the input matrix by multiplying it against a collection of test vectors}
The first step of the leave-one-out\index{leave-one-out randomized algorithm} approach is to construct a low-rank approximation of the matrix $\mat{B}$.
Here, since $\mat{B}$ is a general square matrix, we employ the randomized SVD (\cref{sec:rsvd}).
Let $\mat{\Omega} \in \field^{n\times k}$ be a random matrix with isotropic columns, compute the product $\mat{B}\mat{\Omega}$, and orthogonalize, obtaining $\mat{Q} = \orth(\mat{B}\mat{\Omega})$.
The resulting randomized SVD low-rank approximation is 
\begin{equation*}
    \Bhat \coloneqq \mat{Q}\mat{Q}^*\mat{B}.
\end{equation*}

\myparagraph{Step 2: Decompose the quantity of interest into known piece depending on the low-rank approximation plus a residual}
Next, we decompose the quantity of interest, $\tr(\mat{B})$ in this case, into a known piece depending on the low-rank approximation and an unknown residual term.
The simplest decomposition exploits the linearity of the trace:
\begin{equation*}
    \tr(\mat{B}) = \tr(\Bhat) + \tr(\mat{B} - \Bhat) = \tr(\mat{Q}^*\mat{B}\mat{Q}) + \tr((\Id - \mat{Q}\mat{Q}^*)\mat{B}).
\end{equation*}
We saw a better decomposition from our derivation of the \HutchPP algorithm.
Indeed, using the identity $(\Id - \mat{Q}\mat{Q}^*)^2 = \Id - \mat{Q}\mat{Q}^*$ and the cyclic property of the trace, the residual trace is
\begin{equation*}
    \tr((\Id - \mat{Q}\mat{Q}^*)\mat{B}) = \tr((\Id - \mat{Q}\mat{Q}^*)^2\mat{B}) = \tr((\Id - \mat{Q}\mat{Q}^*)\mat{B}(\Id - \mat{Q}\mat{Q}^*)).
\end{equation*}
Thus, in \XTrace, we shall also use the following decomposition of the trace:
\begin{equation*}
    \tr(\mat{B}) = \tr(\Bhat) + \tr(\mat{B} - \Bhat) = \tr(\mat{Q}^*\mat{B}\mat{Q}) + \tr((\Id - \mat{Q}\mat{Q}^*)\mat{B}(\Id - \mat{Q}\mat{Q}^*)).
\end{equation*}

\myparagraph{Step 3: Construct a Monte Carlo estimate of the residual using a single random vector}
Now, we construct a Monte Carlo approximation to the residual trace $\tr((\Id - \mat{Q}\mat{Q}^*)\mat{B}(\Id - \mat{Q}\mat{Q}^*))$ using a single random vector $\vec{\omega}$.
Assuming $\vec{\omega}$ is isotropic, the natural estimator is the single-vector Girard--Hutchinson estimator
\begin{equation*}
    \vec{\omega}^*(\Id - \mat{Q}\mat{Q}^*)\mat{B}(\Id - \mat{Q}\mat{Q}^*)\vec{\omega} \approx \tr((\Id - \mat{Q}\mat{Q}^*)\mat{B}(\Id - \mat{Q}\mat{Q}^*)).
\end{equation*}
One could also use resphering at this step, but we will postone this refinement to simplify the presentation.
Using the Monte Carlo estimate for the residual, we obtain the following unbiased trace estimate:
\begin{equation} \label{eq:basic-xtrace-estimate}
    \hat{\tr} \coloneqq \tr(\mat{Q}^*\mat{B}\mat{Q}) + \vec{\omega}^*(\Id - \mat{Q}\mat{Q}^*)\mat{B}(\Id - \mat{Q}\mat{Q}^*)\vec{\omega}.
\end{equation}

\myparagraph{Step 4: Downdate the low-rank approximation by recomputing it with a test vector removed, and use the left-out test vector as the random vector in step 3}
So far, the estimate \cref{eq:basic-xtrace-estimate} is nothing special, basically just a lopsided version of \HutchPP where only a single vector $\vec{\omega}$ is used to estimate the residual trace $\tr(\mat{B} - \Bhat)$.
Now, we employ the core device of the leave-one-out approach.
Pick an index $1\le i \le k$, and leave out the $i$th column of the test matrix $\mat{\Omega}$, resulting in a downdated low-rank approximation
\begin{equation*}
    \Bhat_{(i)} \coloneqq \mat{Q}_{(i)}^{\vphantom{*}}\mat{Q}_{(i)}^* \mat{B} \quad \text{where } \mat{Q}_{(i)} \coloneqq \orth(\mat{B}\mat{\Omega}_{-i}).
\end{equation*}
Here and elsewhere, $\mat{\Omega}_{-i}$ denotes $\mat{\Omega}$ without its $i$th column, and $\vec{\omega}_i$ denotes the $i$th column of $\mat{\Omega}$.
The downdated low-rank approximation $\Bhat_{(i)}$ is somewhat less accurate because it has smaller rank $k-1$, but we have freed up an isotropic vector $\vec{\omega}_i$ that is \warn{independent} of the low-rank approximation $\Bhat_{(i)}$ and the orthonormal basis $\mat{Q}_{(i)}$.
Introducing the downdated approximation $\Bhat_{(i)}$ and the left-out vector vector $\vec{\omega}_i$ in the trace estimate \cref{eq:basic-xtrace-estimate} from step 3, we obtain a new basic trace estimate:
\begin{equation} \label{eq:xtrace-i}
    \hat{\tr}_i \coloneqq \tr(\mat{Q}_{(i)}^*\mat{B}\mat{Q}_{(i)}^{\vphantom{*}}) + \vec{\omega}_i^*(\Id - \mat{Q}_{(i)}^{\vphantom{*}}\mat{Q}_{(i)}^*)\mat{B}(\Id - \mat{Q}_{(i)}^{\vphantom{*}}\mat{Q}_{(i)}^*)\vec{\omega}_i^{\vphantom{*}}.
\end{equation}
We call $\hat{\tr}_i$ the \emph{$i$th basic \XTrace estimator}.
It is an unbiased estimate for $\tr(\mat{B})$, owing to the independence of $\vec{\omega}_i$ from $\mat{Q}_{(i)}$.

\myparagraph{Step 5: Average the estimator from step 4 over all choices of vectors to leave out}
Each basic \XTrace estimator $\hat{\tr}_i$ can individually be a poor estimate of the trace, but we can reduce variance by averaging all of them:
\begin{equation} \label{eq:xtrace-definition}
    \hat{\tr}_{\mathrm{X}} \coloneqq \frac{1}{k} \sum_{i=1}^k \hat{\tr}_i.
\end{equation}
We call $\hat{\tr}_{\mathrm{X}}$ the (full) \emph{\XTrace estimator}, short for the e\textbf{X}changeable \textbf{Trace} estimator.

\myparagraph{Discussion}
While it is not obvious yet, the \XTrace estimator \cref{eq:xtrace-definition} can be computed using only $2k$ matvecs ($k$ to compute $\mat{A}\mat{\Omega}$, $k$ to compute $\mat{A}\mat{Q}$).
Equivalently, a fixed budget of $s$ matvecs can accommodate a rank of $k = \lfloor s/2\rfloor$.
(Henceforth, we will assume $s$ is even for simplicity.)
Thus, the \XTrace estimator takes the form
\begin{equation} \label{eq:xtrace}
    \hat{\tr}_{\mathrm{X}} = \frac{1}{s/2} \sum_{i=1}^{s/2} \left[\tr(\mat{Q}_{(i)}^*\mat{B}\mat{Q}_{(i)}^{\vphantom{*}}) + \vec{\omega}_i^*(\Id - \mat{Q}_{(i)}^{\vphantom{*}}\mat{Q}_{(i)}^*)\mat{B}(\Id - \mat{Q}_{(i)}^{\vphantom{*}}\mat{Q}_{(i)}^*)\vec{\omega}_i^{\vphantom{*}}\right].
\end{equation}

Let us compare \XTrace to \HutchPP.
With its budget of $s$ matvecs, \HutchPP dedicates $2s/3$ matvecs to generate a low-rank approximation $\Bhat = \mat{Q}\mat{Q}^*\mat{B}$ and $s/3$ matvecs to forming an estimate of the residual trace $\tr(\mat{B} - \Bhat)$.
\XTrace, in effect, uses the same pool of $s/2$ matvecs both to generate a \emph{family} of low-rank approximations $\Bhat_{(i)}$ and to estimate all of the residual traces $\tr(\mat{B} - \Bhat_{(i)})$.
In this sense, \textit{\XTrace uses all of its multipications with isotropic random test vectors both for low-rank approximation and for Monte Carlo estimation.}
\XTrace also expends an additional $s/2$ matvecs to compute $\mat{B}\mat{Q}$, which is necessary to evaluate the expressions $\tr(\mat{Q}_{(i)}^*\mat{B}\mat{Q}_{(i)}^{\vphantom{*}})$.

\XTrace fixes both weaknesses of \HutchPP that we identified in \cref{sec:hutchpp-weaknesses}.
First, \XTrace is exchangeable, being a symmetric function of all the matvecs $\mat{B}\vec{\omega}_i$ collected in the first phase of the algorithm.
Second, and more importantly, \XTrace is based on a rank-$(s/2)$ approximation to $\mat{B}$, which is more powerful than \HutchPP's rank-$(s/3)$ approximation.
Thus, \XTrace can be significantly more accurate than \HutchPP for matrices $\mat{B}$ with rapid singular value decay, while maintaining accuracy for matrices with slow singular value decay.
Moreover, \XTrace achieves this increased accuracy automatically through an exchangeable, leave-one-out\index{leave-one-out randomized algorithm} design; it does not require manually or adaptively apportioning matvecs between the low-rank approximation and Monte Carlo estimation roles, as in the adaptive \HutchPP algorithm \cite{PCK22}. 
(See \cref{sec:adaptive-hutchpp} for discussion.)

\section{Leave-one-out formula for the randomized SVD} \label{sec:rsvd-downdate}

To implement \XTrace, we need a way of efficiently evaluating the \XTrace estimator \cref{eq:xtrace}.
For the purpose of error estimation (\cref{sec:trace-error-estimation}), we will also need to form each of the individual estimators $\hat{\tr}_i$ defined in \cref{eq:xtrace-i}.
To develop efficient implementations, we can use the following \emph{downdating formula} for the randomized SVD:

\begin{theorem}[Downdating the randomized SVD] \label{thm:rsvd-downdate}
   Let $\mat{B} \in \field^{n\times n}$ and $\mat{\Omega} \in \field^{n\times k}$ be matrices, and assume $\mat{B}\mat{\Omega}$ is full-rank.
   Compute the matrix $\mat{Q}$ defining the randomized SVD approximation $\Bhat = \mat{Q}\mat{Q}^*\mat{B}$ by an economy-size \QR decomposition $\mat{B}\mat{\Omega} = \mat{Q}\mat{R}$.
   The downdated matrices $\mat{Q}_{(i)} = \orth(\mat{B}\mat{\Omega}_{-i})$ admit the representation
   \begin{equation} \label{eq:rsvd-downdate}
       \mat{Q}_{(i)}^{\vphantom{*}}\mat{Q}_{(i)}^* = \mat{Q}(\Id - \vec{s}_i^{\vphantom{*}}\vec{s}_i^*) \mat{Q}^* \quad \text{for } i=1,2,\ldots,k.
   \end{equation}
   The vectors $\vec{s}_i$ are equal to the columns of $\mat{R}^{-*}$, scaled to have unit norm.
\end{theorem}

This result was developed by myself and collaborators in \cite{ET24,ETW24}.

I find this result to be quite surprising. To compute the randomized SVD, we must orthogonalize the matrix $\mat{B}\mat{\Omega}$, which is conventionally accomplished using a \QR decomposition $\mat{B}\mat{\Omega} = \mat{Q}\mat{R}$.
\Cref{thm:rsvd-downdate} demonstrates that the matrix $\mat{R}$, a useless byproduct in most randomized SVD implementations, contains all the information needed to extract all $k$ downdated approximations $\Bhat_{(i)} = \mat{Q}_{(i)}^{\vphantom{*}}\mat{Q}_{(i)}^*\mat{B}$.
These downdated approximations can be represented implicitly using the formula \cref{eq:rsvd-downdate}, and the arithmetic cost is only $\order(k^3)$ operations to compute the $\vec{s}_i$ vectors.
In particular, the cost is independent of both dimensions of the original matrix $\mat{B}$!

\begin{proof}[Proof of \cref{thm:rsvd-downdate}]
    Fix an index $i$.
    To compute $\mat{Q}_{(i)}$, we to form a \QR decomposition of $\mat{B}\mat{\Omega}_{-i}$.
    The decomposition $\mat{B}\mat{\Omega}_{-i} = \mat{Q}\mat{R}_{-i}$ is nearly a \QR decomposition, but the matrix $\mat{R}_{-i}$ is no longer triangular when $i < k$ because its $i$th column has been deleted.
    To restore triangularity, we take a (full) \QR decomposition of $\mat{R}_{-i}$, which we partition as
    \begin{equation*}
        \mat{R}_{-i} = \onebytwo{\tildebold{Q}}{\tildevector{q}}\twobyone{\tildebold{R}}{\mat{0}} \quad \text{for } \tildebold{Q} \in \field^{k\times (k-1)}, \tildebold{R} \in \field^{(k-1)\times (k-1)},\tildevector{q} \in \field^k.
    \end{equation*}
    Using this factorization, we obtain a \QR decomposition of $\mat{B}\mat{\Omega}_{-i}$, namely $$\mat{B}\mat{\Omega}_{-i} = \mat{Q}\mat{R}_{-i} = (\mat{Q}\tildebold{Q}) \tildebold{R}.$$ In particular, $\mat{Q}_{(i)}\coloneqq \mat{Q}\tildebold{Q}$ is an orthonormal basis for the column space of $\mat{B}\mat{\Omega}_{-i}$.
    Therefore, the outer product of $\mat{Q}_{(i)}$ with itself is 
    \begin{equation*}
        \mat{Q}_{(i)}^{\vphantom{*}}\mat{Q}_{(i)}^* = \mat{Q}(\tildebold{Q}\tildebold{Q}^*)\mat{Q}^*.
    \end{equation*}
    Since $\onebytwo{\tildebold{Q}}{\tildevector{q}}$ is unitary, the parenthesized term may be written as 
    \begin{equation*}
        \tildebold{Q}\tildebold{Q}^* = \onebytwo{\tildebold{Q}}{\tildevector{q}}\onebytwo{\tildebold{Q}}{\tildevector{q}}^* - \tildevector{q}\, \tildevector{q}^* = \Id - \tildevector{q}\,\tildevector{q}^*.
    \end{equation*}
    Combining the two previous displays yields
    \begin{equation*}
        \mat{Q}_{(i)}^{\vphantom{*}}\mat{Q}_{(i)}^* = \mat{Q}(\Id - \tildevector{q}\,\tildevector{q}^*)\mat{Q}^*.
    \end{equation*}
    
    To establish the desired result, it remains to show that the unit vector $\tildevector{q}$ is proportional to the $i$th column of $\mat{R}^{-*}$.
    By construction, $\tildevector{q}$ is orthogonal to the column space of $\mat{R}_{-i}$.
    That is, $\mat{R}_{-i}^*\tildevector{q} = \vec{0}$.
    Restoring the deleted $i$th column to $\mat{R}$ and recalling that $\mat{R}$ is nonsingular, we conclude that $\mat{R}^*\tildevector{q} = \alpha \cdot \evec_i$ for some nonzero scalar $\alpha$.
    Therefore, $\tildevector{q} = \alpha \cdot \mat{R}^{-*}\evec_i$ is proportional to the $i$th column of $\mat{R}^{-*}$.
\end{proof}

\section{Implementing \XTrace efficiently} \label{sec:xtrace-implementation}

Having established \cref{thm:rsvd-downdate}, we can use it to calculate the basic \XTrace estimates $\hat{\tr}_i$, defined in \cref{eq:xtrace-i}.
This section begins with a derivation for an efficiently computable formula for the $\hat{\tr}_i$'s, after which implementation is discussed.

\subsection{Formula for the basic \XTrace estimates: Derivation}

The easiest way to derive efficient formulas for a leave-one-out algorithm,\index{leave-one-out randomized algorithm} in my experience, is to proceed methodically and introduce variables to represent intermediate matrices that arise during the derivation.
Let me demonstrate as we derive an efficient \XTrace implementation.

We begin with the randomized SVD step.
First, generate a random matrix $\mat{\Omega} \in \field^{n\times (s/2)}$ with isotropic random columns.
Next, form the product
\begin{equation} \label{eq:xtrace-Y-pre-QR}
    \mat{Y} \coloneqq \mat{B}\mat{\Omega}
\end{equation}
and obtain its (economy-size) \QR decomposition
\begin{equation} \label{eq:xtrace-Y}
    \mat{Y} = \mat{Q}\mat{R}.
\end{equation}
Then, to use the randomized SVD downdating formula \cref{eq:rsvd-downdate}, generate the matrix $\mat{S}$ (with columns $\vec{s}_i$) by building $\mat{R}^{-*}$ and recaling each of its columns to have norm one.
Finally, compute
\begin{equation} \label{eq:xtrace-Z}
    \mat{Z} \coloneqq \mat{B}\mat{Q}.
\end{equation}
The matrices $\mat{Y}$, $\mat{Q}$, $\mat{S}$, and $\mat{Q}$ will be used later in our derivation.

Next, we turn to the basic \XTrace estimators $\hat{\tr}_i$, which were defined in \cref{eq:xtrace-i}.
To use the randomized SVD downdating formula \cref{eq:rsvd-downdate}, invoke the cyclic property of the trace and write:
\begin{equation*}
    \hat{\tr}_i = \tr(\mat{B}\mat{Q}_{(i)}^{\vphantom{*}}\mat{Q}_{(i)}^*) + \vec{\omega}_i^*(\Id - \mat{Q}_{(i)}^{\vphantom{*}}\mat{Q}_{(i)}^*)\mat{B}(\Id - \mat{Q}_{(i)}^{\vphantom{*}}\mat{Q}_{(i)}^*)\vec{\omega}_i^{\vphantom{*}}.
\end{equation*}
Now, invoke the downdating formula \cref{eq:rsvd-downdate}:
\begin{equation} \label{eq:xtrace-downdated-two-terms}
    \hat{\tr}_i = \underbrace{\tr(\mat{B}\mat{Q}^{\vphantom{*}}(\Id - \vec{s}_i^{\vphantom{*}}\vec{s}_i^*)\mat{Q}^*)}_{\textcircled{A}} + \underbrace{\vec{\omega}_i^*(\Id - \mat{Q}(\Id - \vec{s}_i^{\vphantom{*}}\vec{s}_i^*)\mat{Q}^*)\mat{B}(\Id - \mat{Q}^{\vphantom{*}}(\Id - \vec{s}_i^{\vphantom{*}}\vec{s}_i^*)\mat{Q}^*)\vec{\omega}_i^{\vphantom{*}}}_{\textcircled{B}}.
\end{equation}
The result is a sum of two terms, \textcircled{A} and \textcircled{B}, which we will treat separately.

Begin with \textcircled{A}.
Use the cyclic property of the trace and \cref{eq:xtrace-Z} to rewrite
\begin{equation*}
    \textcircled{A} = \tr(\mat{Q}^*\mat{B}\mat{Q} (\Id - \vec{s}_i^{\vphantom{*}}\vec{s}_i^*)) = \tr(\mat{Q}^*\mat{Z} (\Id - \vec{s}_i^{\vphantom{*}}\vec{s}_i^*)).
\end{equation*}
To evaluate this expression, compute
\begin{equation} \label{eq:xtrace-H}
    \mat{H} \coloneqq \mat{Q}^*\mat{Z}.
\end{equation}
The expression \textcircled{A} now simplifies as 
\begin{equation} \label{eq:xtrace-A}
    \textcircled{A} = \tr(\mat{H}) - \vec{s}_i^*\mat{H}\vec{s}_i^{\vphantom{*}}.
\end{equation}
This formula constitutes our final expression for \textcircled{A}.

Now, we treat \textcircled{B}.
Begin by forming the matrix
\begin{equation} \label{eq:xtrace-W}
    \mat{W} \coloneqq \mat{Q}^*\mat{\Omega},
\end{equation}
from which we may define a matrix $\mat{X} \in \field^{(s/2)\times (s/2)}$ with columns
\begin{equation} \label{eq:xtrace-x}
    \vec{x}_i \coloneqq (\Id - \vec{s}_i^{\vphantom{*}}\vec{s}_i^*)\mat{Q}^*\vec{\omega}_i = \vec{w}_i^{\vphantom{*}} - \vec{s}_i^{\vphantom{*}}\cdot\vec{s}_i^*\vec{w}_i^{\vphantom{*}}.
\end{equation}
With this definition, \textcircled{B} can be written as
\begin{align*}
    \textcircled{B} &= \vec{\omega}_i^* (\Id - \mat{Q}\mat{Q}^* + \mat{Q}\vec{s}_i^{\vphantom{*}}\vec{s}_i^*\mat{Q}^*)\mat{B}(\vec{\omega}_i^{\vphantom{*}} - \mat{Q}\vec{x}_i) \\ 
    &= (\vec{\omega}_i^*(\Id - \mat{Q}\mat{Q}^*) + \mat{w}_i^*\vec{s}_i^{\vphantom{*}}\cdot\vec{s}_i^*\mat{Q}^*)(\vec{y}_i^{\vphantom{*}} - \mat{Z}\vec{x}_i).
\end{align*}
Since $\mat{Y} = \mat{B}\mat{\Omega} = \mat{Q}\mat{R}$ is a \QR decomposition,
\begin{equation*}
    (\Id - \mat{Q}\mat{Q}^*)\vec{y}_i = \vec{0} \quad \text{and} \quad \mat{Q}^*\vec{y}_i = \vec{r}_i.
\end{equation*}
Therefore, 
\begin{equation*}
    \textcircled{B} = -\vec{\omega}_i^* (\Id - \mat{Q}\mat{Q}^*)\mat{Z}\vec{x}_i + \mat{w}_i^*\vec{s}_i^{\vphantom{*}}\cdot \vec{s}_i^*(\vec{r}_i^{\vphantom{*}} - \mat{H}\vec{x}_i).
\end{equation*}
To simplify further, introduce and form the matrix 
\begin{equation} \label{eq:xtrace-T}
    \mat{T} \coloneqq \mat{Z}^*\mat{\Omega}.
\end{equation}
Now, we may simplify \textcircled{B} as
\begin{equation*}
    \textcircled{B} = -\vec{t}_i^*\vec{x}_i^{\vphantom{*}} + \vec{w}_i^*\mat{H}\vec{x}_i^{\vphantom{*}} + \mat{w}_i^*\vec{s}_i^{\vphantom{*}}\cdot \vec{s}_i^*(\vec{r}_i^{\vphantom{*}} - \mat{H}\vec{x}_i).
\end{equation*}
Finally, using the definition \cref{eq:xtrace-x} of the vectors $\vec{x}_i$, we simplify
\begin{equation*}
    \textcircled{B} = -\vec{t}_i^*\vec{x}_i^{\vphantom{*}} + \vec{x}_i^*\mat{H}\vec{x}_i^{\vphantom{*}} + \mat{w}_i^*\vec{s}_i^{\vphantom{*}}\cdot \vec{s}_i^*\vec{r}_i^{\vphantom{*}}.
\end{equation*}
Combining this expression for \textcircled{B} with the expression \cref{eq:xtrace-A} for \textcircled{A}, we obtain our final expression for $\hat{\tr}_i$:
\begin{equation} \label{eq:xtrace-i-formula}
    \hat{\tr}_i = \tr(\mat{H}) - \vec{s}_i^*\mat{H}\vec{s}_i^{\vphantom{*}}-\vec{t}_i^*\vec{x}_i^{\vphantom{*}} + \vec{x}_i^*\mat{H}\vec{x}_i^{\vphantom{*}} + \mat{w}_i^*\vec{s}_i^{\vphantom{*}}\cdot \vec{s}_i^*\vec{r}_i^{\vphantom{*}}.
\end{equation}

\subsection{Formula for the basic \XTrace estimates: Implementation}

We now discuss a few implementation details for evaluating the formula \cref{eq:xtrace-i-formula} for the basic trace estimates $\hat{\tr}_i$.

\myprogram{Efficient implementation of the $\diagprod$ operation \cref{eq:diagprod}.}{}{diagprod}
\myparagraph{Diagonal entries of products}
The formula \cref{eq:xtrace-i-formula} contains many expressions such as $\vec{f}_i^*\vec{g}_i^{\vphantom{*}}$, which constitute the diagonal entries of the matrix product $\mat{F}^*\mat{G}$ for $\mat{F},\mat{G} \in \field^{d_1\times d_2}$.
We denote the vector of all $\vec{f}_i^*\vec{g}_i^{\vphantom{*}}$ using the $\diagprod$ operation: 
\begin{equation} \label{eq:diagprod}
    \diagprod(\mat{F},\mat{G}) \coloneqq \diag(\mat{F}^*\mat{G}) = (\vec{f}_i^*\vec{g}_i^{\vphantom{*}} : 1\le i \le d_2).
\end{equation}
One way of evaluating the $\diagprod$ operation would be to form the matrix product $\mat{F}^*\mat{G}$ and extract its diagonal; this approach expends $\order(d_1^2d_2)$ operations.
However, computing these expressions directly is cheaper, requiring just just $\order(d_1d_2)$ operations to evaluate $\vec{f}_i^*\vec{g}_i^{\vphantom{*}}$ for each $i$.
Code is given in \cref{prog:diagprod}

\myparagraph{The vectors $\vec{x}_i$}
The vectors $\vec{x}_i$ defined in \cref{eq:xtrace-x} can be packaged into a matrix $\mat{X}$ with formula
\begin{equation} \label{eq:xtrace-X}
    \mat{X} = \mat{W} - \mat{S} \cdot \Diag(\diagprod(\mat{S},\mat{W})).
\end{equation}
Recall that $\Diag(\vec{a})$ denotes the diagonal matrix with diagonal entries $a_i$.
In MATLAB, expression $\mat{S} \cdot \Diag(\vec{a})$ can be evaluated rapidly as \texttt{S .* a.'}.

\myparagraph{\XTrace implementation}
Using all of the formulas we've developed, the vector $\vec{\hat{\tr}} = (\hat{\tr}_i : 1\le i \le s/2)$ of basic \XTrace estimators \cref{eq:xtrace-i-formula} is
\begin{multline} \label{eq:xtrace-vec}
    \vec{\hat{\tr}} = \tr(\mat{H})\onevec - \diagprod(\mat{S},\mat{H}\mat{S}) - \diagprod(\mat{T},\mat{X})+ \diagprod(\mat{X},\mat{H}\mat{X}) \\ + \diagprod(\mat{W},\mat{S}) \odot \diagprod(\mat{S},\mat{R}).
\end{multline}
The \XTrace estimator \cref{eq:xtrace} is merely the mean of the entries of this vector; the variance of the entries will be used for error estimation in \cref{sec:trace-error-estimation}.
To evaluate $\vec{\hat{\tr}}$, we evaluate the equations \cref{eq:xtrace-Y-pre-QR,eq:xtrace-Y,eq:xtrace-Z,eq:xtrace-H,eq:xtrace-W,eq:xtrace-T,eq:xtrace-X,eq:xtrace-vec}.
Stringing these formulas together produces a daunting and mysterious looking program, but this program is nothing complicated---just matrix algebra.

\myprogram{Efficient implementation of the \XTrace estimator.}{Subroutines \texttt{diagprod}, \texttt{random\_signs}, and \texttt{cnormc} appear in \cref{prog:diagprod,prog:random_signs,prog:cnormc}.}{xtrace}

An implementation of \XTrace is provided in \cref{prog:xtrace}.
This implementation outputs the trace estimator $\hat{\tr}_{\mathrm{X}}$ as first output \texttt{tr}; it also outputs an error estimate \texttt{est}, which will be discussed in \cref{sec:trace-error-estimation}.

\section{\XNysTrace: Trace estimation for psd matrices} \label{sec:xnystrace-derivation}

\XTrace is an effective trace estimator for general square matrices, but it can be made even more efficient for psd matrices.
This section will derive \XNysTrace, an optimized trace estimator for psd matrices using the leave-one-out\index{leave-one-out randomized algorithm} approach.
We uncover the algorithm by following the five-step process for deriving a leave-one-out\index{leave-one-out randomized algorithm} randomized matrix algorithm, similar to our derivation of \XTrace in \cref{sec:xtrace-derivation}.
\Cref{sec:nystrom-downdating,sec:xnystrace-implementation} will discuss efficient implementation.
Throughout this section, and the following sections, $\mat{A} \in \field^{n\times n}$ denotes a psd matrix.

\myparagraph{Step 1: Compute a low-rank approximation to the input matrix by multiplying it against a collection of test vectors}
Since the matrix $\mat{A}$ is psd, we have a more rich set of low-rank approximation algorithms available to us.
Here, we will employ the single pass Nystr\"om approximation 
\begin{equation*}
    \Ahat \coloneqq \mat{A}\langle \mat{\Omega}\rangle
\end{equation*}
associated with a test matrix $\mat{\Omega} \in \field^{n\times s}$ with independent, isotropic random columns (\cref{sec:nystrom}).
An advantage of the single-pass Nystr\"om approximation is that we can form a rank-$s$ approximation using only $s$ matvecs, whereas $s$ matvecs only allow us to obtan a rank-$(s/2)$ approximation with the randomized SVD.

\myparagraph{Step 2: Decompose the quantity of interest into known piece depending on the low-rank approximation plus a residual}
We decompose the trace by employing the simplest possible approach
\begin{equation*}
    \tr(\mat{A}) = \tr(\Ahat) + \tr(\mat{A} - \Ahat).
\end{equation*}

\myparagraph{Step 3: Construct a Monte Carlo estimate of the residual using a single random vector}
To estimate the residual trace $\tr(\mat{A} - \Ahat)$ using a single random vector, we employ the single-vector Girard--Hutchinson estimator
\begin{equation*}
    \vec{\omega}^*(\mat{A}-\Ahat)\vec{\omega}\approx \tr(\mat{A} - \Ahat),
\end{equation*}
which leads to the trace estimate
\begin{equation} \label{eq:xnystrace-before-leave-out}
    \hat{\tr} \coloneqq \tr(\Ahat) + \vec{\omega}^*(\mat{A}-\Ahat)\vec{\omega}. 
\end{equation}

\myparagraph{Step 4: Downdate the low-rank approximation by recomputing it with a test vector removed, and use the left-out test vector as the random vector in step 3}
Now, we invoke the leave-one-out\index{leave-one-out randomized algorithm} trick.
Define the downdated the Nystr\"om approximation by leaving out column $i$ of the test matrix $\mat{\Omega}$,
\begin{equation*}
    \Ahat_{(i)} \coloneqq \mat{A}\langle \mat{\Omega}_{-i} \rangle.
\end{equation*}
Using the left-out vector $\vec{\omega}_i$ as the test vector $\vec{\omega}$ in \cref{eq:xnystrace-before-leave-out} gives the family of basic \XNysTrace estimators
\begin{equation} \label{eq:xnystrace-i}
    \hat{\tr}_i \coloneqq \tr(\Ahat_{(i)}) + \vec{\omega}_i^* (\mat{A} - \Ahat_{(i)})\vec{\omega}_i^{\vphantom{*}} \quad \text{for } i=1,\ldots,s.
\end{equation}

\myparagraph{Step 5: Average the estimator from step 4 over all choices of vectors to leave out}
Averaging the basic \XNysTrace estimators $\hat{\tr}_i$ gives the full \XNysTrace estimator
\begin{equation} \label{eq:xnystrace}
    \hat{\tr}_{\mathrm{XN}} \coloneqq \frac{1}{s} \sum_{i=1}^s \hat{\tr}_i = \frac{1}{s} \sum_{i=1}^s \left[\tr(\Ahat_{(i)}) + \vec{\omega}_i^* (\mat{A} - \Ahat_{(i)})\vec{\omega}_i^{\vphantom{*}}\right].
\end{equation}

\myparagraph{Discussion}
As we will show in \cref{sec:nystrom-downdating,sec:xnystrace-implementation}, the estimator \XNysTrace estimator \cref{eq:xnystrace} can be computed using only the $s$ matvecs composing $\mat{A}\mat{\Omega}$.
\XNysTrace is exchangeable and dedicates all matvecs both to low-rank approximation and to Monte Carlo estimation of the residual trace.

The use of Nystr\"om approximation for variance reduction in trace estimation was originally proposed by Persson, Cortinovis, and Kressner \cite{PCK22}, who developed a \HutchPP-style algorithm using Nys\"trom approximation called \textsc{Nystr\"om}\texttt{++}.
The \XNysTrace algorithm improves on \textsc{Nystr\"om}\texttt{++} by using an exchangeable, variance-reduced design.

The main difference between \HutchPP, \XTrace, and \XNysTrace is the rank of the matrix approximation.
Given a fixed budget of $s$ matvecs, \HutchPP uses a rank-$(s/3)$ approximation, \XTrace uses a rank-$(s/2)$ approximation, and \XNysTrace uses a rank-$s$ approximation.
This disparity results in significant differences among these estimators \warn{when applied to a matrix with rapidly decaying eigenvalues}.
As an example, consider a psd matrix $\mat{A}$ whose eigenvalues decay at an exponential rate $\lambda_i(\mat{A}) \le \alpha^i$ for $\alpha \in (0,1)$.
As we will show in \cref{sec:trace-error-bounds}, \HutchPP, \XTrace, and \XNysTrace satisfy error bounds of the form
\begin{equation} \label{eq:trace-bounds-compare-ch14}
    \begin{split}
	\left(\mathbb{E} \bigl|\hutchppest - \tr (\mat{A}) \bigr|^2\right)^{1 \slash 2} &\le \makebox[\widthof{$\sqrt{s}$}][r]{}\,C_1(\alpha)\, \alpha^{\color{red} s/3}; \\ \left(\mathbb{E} \bigl|\hat{\tr}_{\mathrm{X}} - \tr (\mat{A}) \bigr|^2\right)^{1 \slash 2}&\le \sqrt{s} \,C_2(\alpha)\,\alpha^{\color{red} s/2}; \\
	\left(\mathbb{E} \bigl|\hat{\tr}_{\mathrm{XN}} - \tr (\mat{A}) \bigr|^2\right)^{1 \slash 2}&\le \makebox[\widthof{$\sqrt{s}$}][r]{$s$} \,C_3(\alpha)\, \alpha^{\color{red} s}.
    \end{split}
\end{equation}
Here, $C_i(\alpha)$ denote prefactors depending only on $\alpha$.
For this matrix, the rate of convergence of \XTrace is $3/2\times$ faster than \HutchPP, and the convergence rate of \XNysTrace is $3\times$ faster than \HutchPP.

\section{Leave-one-out formula for randomized Nystr\"om approximation} \label{sec:nystrom-downdating}

Just as the randomized SVD downdating formula (\cref{thm:rsvd-downdate}) was the main ingredient in making an efficient \XTrace implementation, a fast Nystr\"om approximation downdating formula will support a fast \XNysTrace implementation.
This formula will also be useful in deriving other leave-one-out\index{leave-one-out randomized algorithm} algorithms based on randomized Nystr\"om approximation.

\begin{theorem}[Downdating randomized Nystr\"om approximation] \label{thm:nystrom-downdate}
    Let $\mat{A} \in \field^{n\times n}$ be psd and let $\mat{\Omega} \in \field^{n\times k}$ be a matrix, and assume $\mat{A}\mat{\Omega}$ is full-rank.
    Denote $\mat{H} \coloneqq \mat{\Omega}^*(\mat{A}\mat{\Omega})$.
    Then, the downdated Nystr\"om approximation $\Ahat_{(i)} \coloneqq \mat{A}\langle\mat{\Omega}_{-i}\rangle$ has the representation
    \begin{subequations}\label{eq:nystrom-downdate}
    \begin{align} 
        \Ahat_{(i)} &= (\mat{A}\mat{\Omega}) \left( \mat{H}^{-1} - \frac{\mat{H}^{-1}(:,i)\mat{H}^{-1}(i,:)}{\mat{H}^{-1}(i,i)} \right) (\mat{A}\mat{\Omega})^* \label{eq:nystrom-downdate-1}\\
        &= \Ahat - \frac{(\mat{A}\mat{\Omega}\mat{H}^{-1}(:,i))(\mat{A}\mat{\Omega}\mat{H}^{-1}(:,i))^*}{\mat{H}^{-1}(i,i)}.\label{eq:nystrom-downdate-2}
    \end{align}
    \end{subequations}
\end{theorem}

The proof of this result relies on a consequence of the Banachiewicz inversion formula.
Here, we restate the version of this formula given in \cite[eq.~(0.7.2)]{PS05}.
\begin{fact}[Banachiewicz inversion formula] \label{fact:Banachiewicz}
    Let $\mat{P}\in\field^{d_1\times d_1},\mat{Q}\in\field^{d_1\times d_2},\mat{R}\in\field^{d_2\times d_1},\mat{S}\in\field^{d_2\times d_2}$ with $\mat{P}$ nonsingular.
    Then the block matrix 
    \begin{equation*}
        \mat{M} \coloneqq \twobytwo{\mat{P}}{\mat{Q}}{\mat{R}}{\mat{S}}
    \end{equation*}
    is invertible  if and only if the \emph{Schur complement} $\mat{M} / \mat{P} \coloneqq \mat{S} - \mat{R}\mat{P}^{-1}\mat{Q}$ is invertible, in which case 
    \begin{subequations} \label{eq:Banachiewicz}
    \begin{align}
        \mat{M}^{-1} = \twobytwo{\mat{P}}{\mat{Q}}{\mat{R}}{\mat{S}}^{-1} &= \twobytwo{\mat{P}^{-1}+\mat{P}^{-1}\mat{Q}(\mat{M}/\mat{P})^{-1}\mat{R}\mat{P}^{-1}}{-\mat{P}^{-1}\mat{Q}(\mat{M}/\mat{P})^{-1}}{-(\mat{M}/\mat{P})^{-1}\mat{R}\mat{P}^{-1}}{(\mat{M}/\mat{P})^{-1}} \label{eq:Banachiewicz-1} \\
        &= \twobytwo{\mat{P}^{-1}}{\mat{0}}{\mat{0}}{\mat{0}} + \twobyone{-\mat{P}^{-1}\mat{Q}}{\Id} (\mat{M}/\mat{P})^{-1}\onebytwo{-\mat{R}\mat{P}^{-1}}{\Id}. \label{eq:Banachiewicz-2}
    \end{align}
    \end{subequations}
\end{fact}

As a consequence, we obtain a formula relating the inverse of a matrix to the inverse of a submatrix.

\begin{corollary}[Downdating the inverse] \label{cor:inverse-downdating}
    Instate the notation and assumptions of \cref{fact:Banachiewicz}, and denote $\set{E} \coloneqq \{d_1+1,\ldots,d_1+d_2\}$.
    Then
    \begin{equation*}
        \twobytwo{\mat{P}^{-1}}{\mat{0}}{\mat{0}}{\mat{0}} = \mat{M}^{-1} - \mat{M}^{-1}(:,\set{E})[\mat{M}^{-1}(\set{E},\set{E})]^{-1} \mat{M}^{-1}(\set{E},:).
    \end{equation*}
\end{corollary}

\begin{proof}
    Rewrite \cref{eq:Banachiewicz-2} as 
    \begin{equation*}
        \mat{M}^{-1} = \twobytwo{\mat{P}^{-1}}{\mat{0}}{\mat{0}}{\mat{0}} + \twobyone{-\mat{P}^{-1}\mat{Q}(\mat{M}/\mat{P})^{-1}}{(\mat{M}/\mat{P})^{-1}} (\mat{M}/\mat{P}) \onebytwo{-(\mat{M}/\mat{P})^{-1}\mat{R}\mat{P}^{-1}}{(\mat{M}/\mat{P})^{-1}},
    \end{equation*}
    and use \cref{eq:Banachiewicz-1} to recognize the factors of the second term as submatrices of $\mat{M}^{-1}$:
    \begin{equation*}
        \twobyone{-\mat{P}^{-1}\mat{Q}(\mat{M}/\mat{P})^{-1}}{(\mat{M}/\mat{P})^{-1}} = \mat{M}^{-1}(:,\set{E}), \quad \onebytwo{-(\mat{M}/\mat{P})^{-1}\mat{R}\mat{P}^{-1}}{(\mat{M}/\mat{P})^{-1}} = \mat{M}^{-1}(\set{E},:).
    \end{equation*}
    Finally observe that $\mat{M}/\mat{P} = [\mat{M}^{-1}(\set{E},\set{E})]^{-1}$ is the inverse of the block-$(2,2)$ entry of $\mat{M}^{-1}$.
    Combining these observations yields the desired result.
\end{proof}

With this formula and its corollary in hand, the proof of \cref{thm:nystrom-downdate} is immediate.

\begin{proof}[Proof of \cref{thm:nystrom-downdate}]
    The Nystr\"om approximation $\mat{A}\langle \mat{\Omega}\rangle$ is invariant to permutation of the columns of $\mat{\Omega}$ (\cref{prop:nystrom-properties}\ref{item:nystrom-invariance}).
    As such, we can assume without loss of generality that $i = k$ by permuting the $i$th column of $\mat{\Omega}$ to appear last.

    The downdated Nystr\"om approximation $\Ahat_{(k)}$ takes the form
    \begin{equation*}
        \Ahat_{(k)} = \mat{A}\mat{\Omega}_{-k}^{\vphantom{*}} \big(\mat{H}_{(k)}\big)^{-1} \mat{\Omega}_{-k}^*\mat{A} = \mat{A}\mat{\Omega} \twobytwo{\mat{H}_{(k)}^{-1}}{\vec{0}}{\vec{0}^*}{0}\mat{\Omega}^*\mat{A}.
    \end{equation*}
    Invoking \cref{cor:inverse-downdating} and repackaging $\mat{\Omega}^*\mat{A} = (\mat{A}\mat{\Omega})^*$ yields the stated result.
\end{proof}

\subsection{Using the formula \cref{eq:nystrom-downdate}}

The Nystr\"om downdating formula \cref{eq:nystrom-downdate} can be combined with the stable Nystr\"om implementation from \cref{eq:stable-nystrom}.
We treat the shift $\mu$ as zero for the following discussion.
In practice, the shift $\mu$ given by \cref{eq:nys-mu} should be used.

We compute the Nystr\"om approximation in outer product form $\Ahat = \mat{F}\mat{F}^*$ using the stable implementation in \cref{eq:stable-nystrom}.
First, we compute the matrix product $\mat{Y} = \mat{A}\mat{\Omega}$.
Then, we obtain a Cholesky decomposition $\mat{\Omega}^*\mat{A}\mat{\Omega} = \mat{R}^*\mat{R}$.
Finally, we construct the factor matrix $\mat{F} = \mat{Y} \mat{R}^{-1}$.
Using these matrices, the downdating formula \cref{eq:nystrom-downdate} may be written
\begin{subequations} \label{eq:nystrom-downdate-practice}
\begin{equation} \label{eq:nystrom-downdate-practice-1}
    \Ahat_{(i)} = \mat{F}\mat{F}^* - \vec{z}_i^{\vphantom{*}}\vec{z}_i^* \quad \text{for } i =1,2,\ldots,k
\end{equation}
where
\begin{equation} \label{eq:nystrom-downdate-practice-2}
    \mat{Z} \coloneqq \mat{F}\mat{R}^{-*} \cdot \Diag\Bigl(\srn\bigl(\mat{R}^{-1}\bigr)\Bigr)^{-1/2}.
\end{equation}
\end{subequations}
The matrix $\mat{Z}$ contains all of the information needed to do Nystr\"om downdating.

\section{Implementing \XNysTrace efficiently} \label{sec:xnystrace-implementation}

We can use \cref{thm:nystrom-downdate} to form the \XNysTrace estimator rapidly.
This section begins with a derivation of an efficient formula for the basic \XNysTrace estimates \cref{eq:xnystrace-i}.
Then, we will discuss how to implement this formula.

\subsection{Formula for the basic \XNysTrace estimates: Derivation}

Begin by generating a matrix $\mat{\Omega}$ with isotropic columns, compute the Nystr\"om $\Ahat = \mat{F}\mat{F}^*$ via \cref{eq:stable-nystrom}, and form $\mat{Z}$ from \cref{eq:nystrom-downdate-practice-2}.
Then, substitute the Nystr\"om downdating formula \cref{eq:nystrom-downdate-1} in the definition \cref{eq:xnystrace-i} of the basic \XNysTrace estimators:
\begin{equation} \label{eq:nystrom-intermediate}
    \hat{\tr_i} = \tr(\mat{F}\mat{F}^* - \vec{z}_i^{\vphantom{*}}\vec{z}_i^*) + \vec{\omega}_i^*(\mat{A} - \Ahat + \vec{z}_i^{\vphantom{*}}\vec{z}_i^*)\vec{\omega}_i^{\vphantom{*}}.
\end{equation}
By \cref{prop:nystrom-properties}\ref{item:nystrom-interpolatory}, the Nystr\"om approximation $\Ahat$ satisfies the \emph{interpolation condition} $\Ahat \vec{\omega}_i = \mat{A}\vec{\omega}_i$.
Using this interpolation property, the cyclic property of the trace, and the identity $\norm{\mat{F}}_{\mathrm{F}}^2 = \tr(\mat{F}^*\mat{F})$, the formula \cref{eq:nystrom-intermediate} simplifies as
\begin{equation*}
    \hat{\tr}_i = \norm{\mat{F}}_{\mathrm{F}}^2 - \norm{\vec{z}_i}^2 + |\vec{z}_i^*\vec{\omega}_i|^2.
\end{equation*}
More concisely, the vector of trace estimates $\vec{\tr} = (\hat{\tr}_i : 1\le i\le s)$ is 
\begin{equation*}
    \vec{\hat{\tr}} = \norm{\mat{F}}_{\mathrm{F}}^2\cdot\onevec - \srn(\mat{Z}) + |\,{\diagprod(\mat{Z},\mat{\Omega})}\,|^2.
\end{equation*}

\subsection{Formula for the basic \XNysTrace estimates: Implementation}

\myprogram{Efficient and stable implementation of \XNysTrace estimator.}{Subroutines \texttt{nystrom} and \texttt{diagprod} appear in \cref{prog:nystrom,prog:diagprod}.}{xnystrace}

In practice, we use the shift $\mu$ given by formula \cref{eq:nys-mu} to ensure numerical stability and success of the Cholesky decomposition \cref{eq:nys-cholesky}.
As a result, the \XNysTrace estimator produces an unbiased estimate of the \emph{shifted matrix} $\mat{A} + \mu \Id$.
To correct for the shift, we remove the trace of the correction $\tr(\mu \Id) = n\mu$ from each trace estimate, resulting in the alternative formula
\begin{equation*}
    \vec{\hat{\tr}} = \norm{\mat{F}}_{\mathrm{F}}^2\cdot\onevec - \srn(\mat{Z}) + |\,{\diagprod(\mat{Z},\mat{\Omega})}\,|^2 - n\mu \cdot \onevec.
\end{equation*}
We use this formula in our code.
Code for \XNysTrace is provided in \cref{prog:xnystrace}, which outputs the \XNysTrace estimator $\hat{\tr}_{\mathrm{XN}}$ as \texttt{tr} and an error estimate \texttt{est}; see \cref{sec:trace-error-estimation} for discussion.

\begin{remark}[Comparison to \XNysTrace implementation in \cite{ETW24}] \label{rem:xnystrace-new-impl}
    Our original paper on \XNysTrace \cite{ETW24} uses a distnctive implementation based on a \QR decomposition of $\mat{Y}$.
    The implementation in \cref{prog:xnystrace} is significantly faster.
    When applied to a problem of dimension $n = 10^5$ with $s = 10^3$ matvecs, the processing time (that is, the total runtime minus the time required to perform matvecs) for \cref{prog:xnystrace} was 7$\times$ faster than the implementation give in \cite{ETW24}.
\end{remark}

\section{Synthetic Experiments} \label{sec:trace-experiments}

To compare the \XTrace and \XNysTrace algorithms to \HutchPP and the Girard--Hutchinson estimator, we evaluate on four test matrices with different spectra:
\begin{subequations} \label{eq:test-matrices}
\begin{align}
    \texttt{flat} &\coloneqq \mat{U}_1^{\vphantom{*}} \Diag( \texttt{linspace}(1,3,n)) \mat{U}_1^*, \label{eq:flat} \\
    \texttt{poly} &\coloneqq \mat{U}_2^{\vphantom{*}} \Diag(i^{-2} : i=1,\ldots,n)\mat{U}_2^*,\label{eq:poly} \\
    \texttt{exp} &\coloneqq \mat{U}_3^{\vphantom{*}} \Diag( 0.7^i : i=0,\ldots,n-1) \mat{U}_3^*, \label{eq:exp} \\
    \texttt{step} &\coloneqq \mat{U}_4^{\vphantom{*}} \Diag(\underbrace{1,\ldots,1}_{50 \text{ times}},\underbrace{10^{-3},\ldots,10^{-3}}_{n-50 \text{ times}})\mat{U}_4^*. \label{eq:step}
\end{align}
\end{subequations}
Here, $\mat{U}_i \in \real^{n\times n}$ denote Haar-random orthogonal matrices, $\texttt{linspace}(1,3,n)$ denotes equally spaced entries between 1 and 3, and $n\coloneqq 10^3$.
Each matrix is real and positive definite, although the \texttt{exp} matrix is singular up to numerical precision.

\begin{figure}
    \centering
    \includegraphics[width=0.99\linewidth]{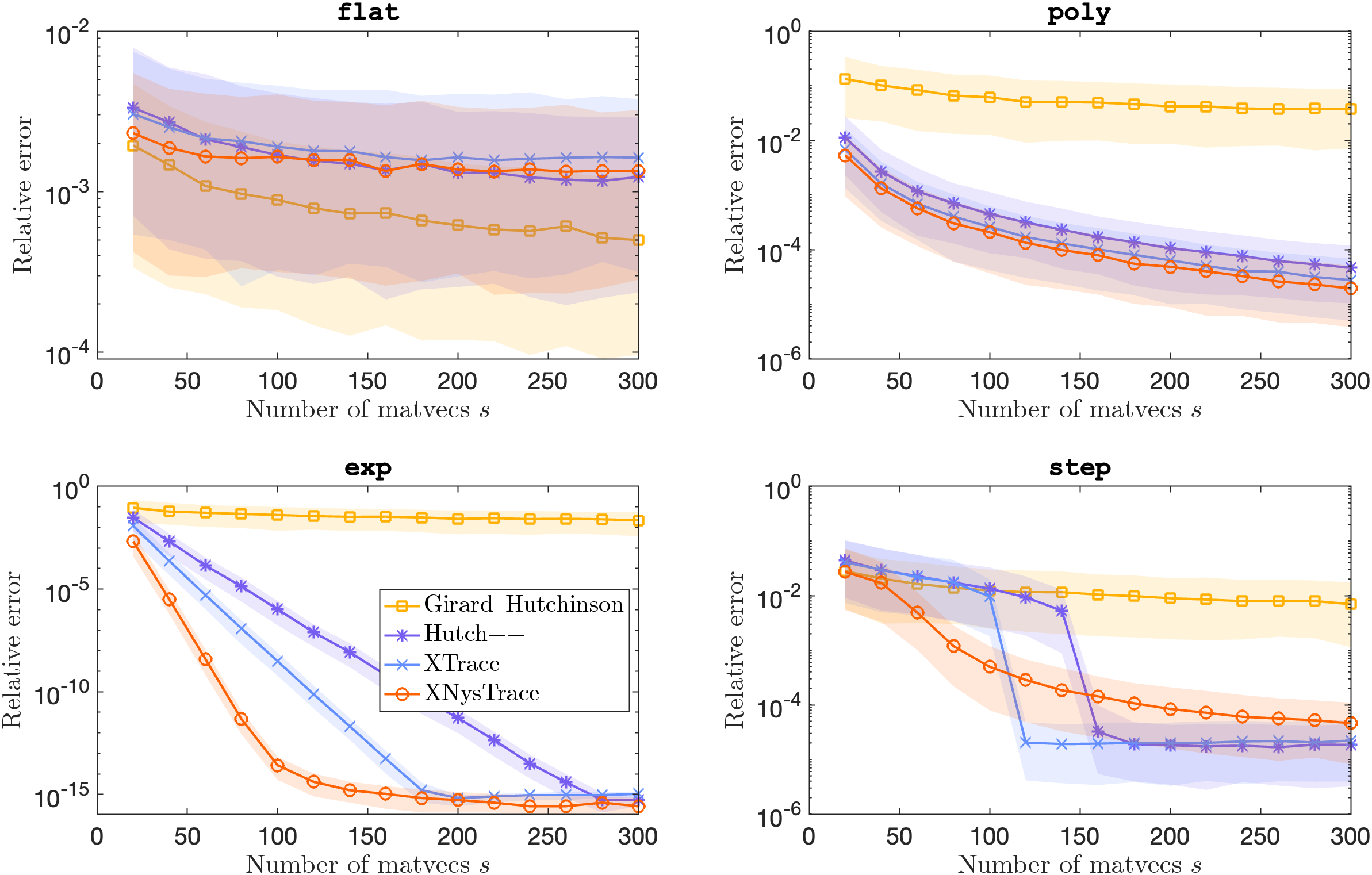}
    \caption[Comparison of Girard--Hutchinson, \HutchPP, \XTrace, and \XNysTrace for estimating the trace of matrices with different spectra]{Relative error of trace estimates by Girard--Hutchinson estimator (yellow squares), \HutchPP (purple asterisks), \XTrace (blue crosses), and \XNysTrace (orange circles) for the four test matrices \cref{eq:test-matrices} for different numbers of matvecs $s$.
    Lines show median of 1000 trials and error bars show 10\% and 90\% quantiles.}
    \label{fig:trace-comparison}
\end{figure}

Results are shown in \cref{fig:trace-comparison}.
Here are my conclusions:

\myparagraph{With spectral decay, \XNysTrace $>$ \XTrace $>$ \HutchPP $\gg$ Girard--Hutchinson}
For matrices with eigenvalues decaying at a steady, and sufficiently rapid, rate (e.g., the \texttt{poly} and \texttt{exp} examples), the ranking of methods is clear: \XNysTrace, \XTrace, \HutchPP, then Girard--Hutchinson.
This disparity is most visible with the \texttt{exp} examples, where \XNysTrace converges at a $3\times$ faster exponential rate than \HutchPP and \XTrace converges $1.5\times$ faster.
The Girard--Hutchinson estimator, by contrast, converges at the much slower Monte Carlo rate $\Theta(s^{-1/2})$.

\myparagraph{Without spectral decay, variance reduction by low-rank approximation doesn't help}
The \HutchPP, \XTrace, and \XNysTrace estimators depend on reducing the variance of the trace estimate by using a low-rank approximation of the matrix as a control variate.
This technique does not provide any benefit, and in fact slightly harms the quality of the estimate, when the matrix has a flat spectrum.
This effect is visible with the \texttt{flat} matrix, where the Girard--Hutchinson estimator beats the \HutchPP, \XTrace, and \XNysTrace estimators by a small multiple.
(Though, resphering can help improve the variance-reduced estimators on problems like this, as in \cref{fig:xtrace-resphere} below.)

\myparagraph{Weakness of \XNysTrace on \texttt{step} matrix}
Perhaps the most interesting example of these four test matrices is the \texttt{step} matrix.
This matrix has $k=50$ large eigenvalues, with the remaining $n-50$ eigenvalues being much smaller.
It takes \HutchPP about $s=3k =150$ matvecs to produce a low-rank approximation capturing these dominant eigenvalues, and it takes \XTrace about $s=2k=100$ matvecs.
After hitting this number of matvecs, the error of \HutchPP and \XTrace drop about three orders of magnitude (and even more with resphering, see \cref{fig:xtrace-resphere} below).
By contrast, \XNysTrace begins reaping the benefits of low-rank approximation at about $s=50$ matvecs.
Yet, the convergence of \XNysTrace for $s>50$ is more gradual than for either \HutchPP and \XTrace.
Consequently, \HutchPP and \XTrace achieve lower error than \XNysTrace for sufficiently large $s$.

An explanation for this behavior is visible in the error bounds for the randomized SVD and the single-pass Nystr\"om approximation provided by \cref{fact:rsvd-error,fact:nystrom-error}.
For a psd matrix, the error of the single-pass Nystr\"om approximation depends on the \emph{sum} of the tail eigenvalues $\sum_{i \gtrapprox k} \lambda_i$, whereas the randomized SVD error depends on the \emph{$\ell_2$ norm} of the tail eigenvalues $(\sum_{i \gtrapprox k} \lambda_i^2)^{1/2}$.
On this example, the sum of tail eigenvalues is larger than the $\ell_2$ norm of the tail eigenvalues by a factor of about $\sqrt{n}$.
Thus, the randomized SVD-based \XTrace and \HutchPP estimators can approximate the \texttt{step} matrix better than the Nystr\"om-based \XNysTrace estimator.

\myparagraphnp{Which method should I use?}
Unfortunately, the pattern of results in \cref{fig:trace-comparison} defies a truly simple conclusion.
Still, I think there are some pretty clear recommendations that can be gleaned from these experiments.

One can imagine two distinct settings for trace estimation.
In the first setting, one is writing general-purpose software, and the trace estimator must be designed to handle arbitrary input matrices.
In the second setting, one is interested in a specific application, and the trace estimator needs only work well for matrices appearing in that application.
For the former setting, I would make the following recommendation:
\actionbox{For general-purpose use, I would recommend \XTrace or, for psd matrices, either \XTrace or \XNysTrace.
Both should be implemented with resphering (\cref{sec:xtrace-resphere}).}
While \XTrace and \XNysTrace are not the best trace estimators for every single problem, the benefits over \HutchPP and the Girard--Hutchinson estimator can be substantial on some problems.
For problems with slow spectral decay (like the \texttt{flat} matrix), the resphering step (\cref{sec:xtrace-resphere}) can substantially improve \XNysTrace.

In a specific application, the choice of trace estimator can be determined by profiling.
As a rule of thumb, \XNysTrace is the best estimator for psd matrices with \emph{consistent} spectral decay, \XTrace is the best estimate for general matrices with at least some singular value decay, and the Girard--Hutchinson estimator is the best estimator (by a small margin) on problems with very little spectral decay.
For additional approaches to trace estimation, see \cref{sec:trace-alternatives}.
The paper \cite{ETW24} for comparisons of \XTrace with the \textsc{Nystr\"om}\texttt{++} and adaptive \HutchPP algorithms of \cite{PCK22}.

\begin{figure}
    \centering
    \includegraphics[width=0.6\linewidth]{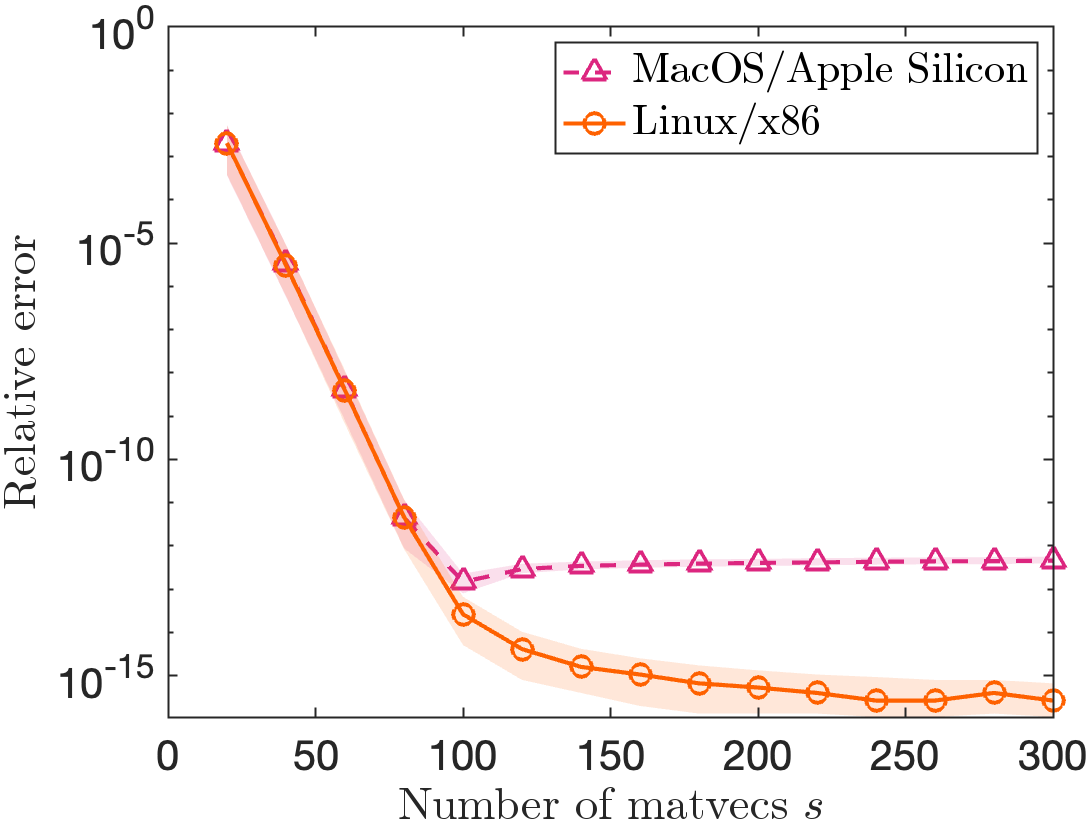}
    \caption[Comparison of \XNysTrace implemented on two different computer architectures. Errors are observed to be significantly higher on an Apple silicon Mac than on a Linux system]{Comparison of the error of \XNysTrace on \texttt{exp} matrix \cref{eq:exp} on Linux x86 machine (MATLAB Online, orange solid circles) and MacOS Apple Silicon (Macbook Pro with M3 ARM chip, pink dashed triangles).
    The error on the Linux system saturates at 2.7e-16 and the Mac system saturates at 4.6e-13, 1700$\times$ higher.
    Lines show median of 1000 trials and error bars show 10\% and 90\% quantiles.}
    \label{fig:platform-dependence}
\end{figure}

\begin{remark}[Platform dependence]
    Curiously, during the numerical experiments for this thesis, I observed that the numerical errors for \XNysTrace were many orders of magnitude higher for matrices than in \cite{ETW24}.
    In \cite{ETW24}, the results show \XNysTrace achieving machine precision, whereas the new results showed the error saturating roughly $10^3\times$ higher.
    Eventually, I isolated the discrepancy to the platform; the original experiments in \cite{ETW24} were performed on an Mac computer with an Intel x86 chip, and the new experiments werre performed on a Mac with an Apple Silicon ARM chip.
    I was able to reproduce the numerical behavior in the original paper \cite{ETW24} on MATLAB Online, which uses x86 Linux machines.
    See \cref{fig:platform-dependence} for a comparison.

    I am not aware for an underlying reason for the platform dependence on the numerical accuracy.
    The distinction between x86 and Apple ARM systems persisted in every numerical experiment I ran across multiple machines, including multiple MATLAB versions up to 2024b.
    To show the best performance for the algorithm, the experiments in \cref{fig:trace-comparison} were performed on MATLAB Online (x86).
\end{remark}

\section{Resphering \XTrace and \XNysTrace} \label{sec:xtrace-resphere}

We can improve the \XTrace and \XNysTrace algorithms by using the resphering technique, discussed at the end of \cref{sec:monte-carlo-matrix-attribute}.

\subsection{Resphering \XTrace}

To resphere \XTrace, we replace $\vec{\omega}_i$ by
\begin{equation*}
    \vec{\nu}_i \coloneqq \sqrt{n-s/2+1}\cdot \frac{(\Id - \outprod{\mat{Q}_{(i)}})\vec{\omega}_i}{\norm{\smash{(\Id - \outprod{\mat{Q}_{(i)}})\vec{\omega}_i}}}.
\end{equation*}
The effect of this substitution is to scale the factor \textcircled{B} in \cref{eq:xtrace-downdated-two-terms} by 
\begin{equation*}
    \alpha_i \coloneqq \frac{n-s/2+1}{\norm{\smash{(\Id - \outprod{\mat{Q}_{(i)}})\vec{\omega}_i}}^2}.
\end{equation*}
Using the $\mat{X}$ matrix defined in \cref{eq:xtrace-X}, we evaluate the denominator to be
\begin{equation*}
    \norm{\smash{(\Id - \outprod{\mat{Q}_{(i)}})\vec{\omega}_i}}^2 = \vec{\omega}_i^*(\Id - \mat{Q}(\Id - \outprod{\vec{s}_i})\mat{Q}^*)\vec{\omega}_i^{\vphantom{*}} = \norm{\vec{\omega}_i}^2 - \norm{\vec{x}_i}^2.
\end{equation*}
Consequently, the vector of scaling factors is
\begin{equation*}
    \vec{\alpha} = \frac{n-s/2+1}{\scn(\mat{\Omega}) - \scn(\mat{X})}.
\end{equation*}
As usual, division is performed elementwise.
With resphering, the vector of \XTrace vectors is given by 
\begin{multline*}
    \vec{\hat{\tr}} = \tr(\mat{H})\onevec - \diagprod(\mat{S},\mat{H}\mat{S}) + \vec{\alpha} \odot [-\diagprod(\mat{T},\mat{X}) \\ + \diagprod(\mat{X},\mat{H}\mat{X}) + \diagprod(\mat{W},\mat{S}) \odot \diagprod(\mat{S},\mat{R})].
\end{multline*}
An implementation of \XTrace with resphering is provided in \cref{prog:xtrace_resphere}.

\myprogram{Efficient implementation of \XTrace algorithm with resphering.}{Subroutines \texttt{diagprod} and \texttt{sqcolnorms} are provided in \cref{prog:diagprod,prog:sqcolnorms}.}{xtrace_resphere}

\subsection{Resphering \XNysTrace}

To resphere \XNysTrace requires a bit more thought. 
Remember that the basic \XNysTrace estimators take the form
\begin{equation*}
    \hat{\tr}_i = \tr(\Ahat_{(i)}) + \vec{\omega}_i^* (\mat{A} - \Ahat^{(i)}) \vec{\omega}_i^{\vphantom{*}} \quad \text{where } \Ahat_{(i)} \coloneqq \mat{A} \langle \mat{\Omega}_{-i} \rangle.
\end{equation*}
To resphere \XNysTrace, we need to identify a matrix $\mat{G}$ such that $(\mat{A} - \Ahat^{(i)})\mat{G} = \mat{0}$.
Such a matrix is furnished the the \emph{interpolatory property} of Nystr\"om approximation (\cref{prop:nystrom-properties}\ref{item:nystrom-interpolatory}), which shows that $\mat{G} = \mat{\Omega}_{-i}$ has this feature:
\begin{equation*}
    (\mat{A} - \Ahat^{(i)})\mat{\Omega}_{-i} = \mat{0}.
\end{equation*}
Thus, we can resphere \XNysTrace by replacing $\vec{\omega}_i$ by 
\begin{equation*}
    \vec{\nu}_i \coloneqq \sqrt{n-s+1}\cdot \frac{(\Id - \outprod{\mat{Q}_{(i)}})\vec{\omega}_i}{\norm{\smash{(\Id - \outprod{\mat{Q}_{(i)}}})\vec{\omega}_i}} \quad \text{for } \mat{Q}_{(i)} \coloneqq \orth(\mat{\Omega}_{-i}).
\end{equation*}

To compute the vectors $\vec{\nu}_i$, we need to orthonormalize the columns of the matrix $\mat{\Omega}$ after every possible column deletion.
Fortunately, the randomized SVD downdating formulas (\cref{thm:rsvd-downdate}) are exactly what we need to perform this computation. 
Begin by computing a \QR decomposition of $\mat{\Omega}$,
\begin{equation} \label{eq:xnystrace-resphere-qr}
    \mat{\Omega} = \mat{Q}\mat{T},
\end{equation}
and form the downdating matrix $\mat{S}$ by normalizing the columns of $\mat{T}^{-*}$.
The downdated orthonormal matrices $\mat{Q}_{(i)}$ then admit the relation
\begin{equation*}
    \outprod{\mat{Q}_{(i)}} = \mat{Q} (\Id - \outprod{\vec{s}_i})\mat{Q}^*.
\end{equation*}

Using this display, a short computation shows the resphered \XNysTrace estimators with shifting are
\begin{subequations} \label{eq:xnystrace-resphere}
\begin{equation}
    \vec{\hat{\tr}} = \norm{\mat{F}}_{\mathrm{F}}^2\cdot\onevec - \srn(\mat{Z}) + \vec{\alpha}\odot|\,{\diagprod(\mat{Z},\mat{\Omega})}\,|^2 - n\mu \cdot \onevec
\end{equation}
where
\begin{equation}
    \vec{\alpha} \coloneqq \frac{n-s+1}{\scn(\mat{\Omega}) - \scn(\mat{X})} \quad \text{with } \mat{X} \coloneqq \mat{T} - \mat{T}^{-*} \cdot \Diag(\scn(\mat{T}^{-*}))^{-1}.
\end{equation}
\end{subequations}

An interesting observation is that the resphered \XNysTrace estimator \cref{eq:xnystrace-resphere} depends only on the triangular factor of the \QR decomposition \cref{eq:xnystrace-resphere-qr}.
As such, we can develop an implementation of the resphered \XNysTrace estimator that avoids \QR decomposition entirely, which is beneficial since \QR decomposition is expensive.
Indeed, we can instead compute a Cholesky decomposition
\begin{equation*}
    \mat{\Omega}^*\mat{\Omega} = \mat{T}^*\mat{T}
\end{equation*}
of the Gram matrix $\mat{\Omega}^*\mat{\Omega}$.
Forming the Gram matrix is highly discouraged as a general practice in matrix computations \cite{Hig22}, but it is benign here because Gaussian random matrices are very well-conditioned (at least if $n\ge s/2$) \cite[\S11.2 and p.~166]{Tro21}.
An implementation of the resphered \XNysTrace estimator is provided in \cref{prog:xnystrace_resphere}.

\myprogram{Efficient implementation of \XNysTrace algorithm with resphering.}{Subroutines \texttt{diagprod}, \texttt{sqcolnorms}, and \texttt{sqrownorms} are provided in \cref{prog:diagprod,prog:sqcolnorms,prog:sqrownorms}.}{xnystrace_resphere}

\begin{remark}[Improved resphered \XNysTrace implementation]
    The publicly available code from \cite{ETW24} uses a \QR-based implementation of the resphered \XNysTrace estimator.
    We have improved it here by introducing the faster Cholesky-based implementation.
\end{remark}

\subsection{Experiments}

\begin{figure}
    \centering
    \includegraphics[width=0.49\linewidth]{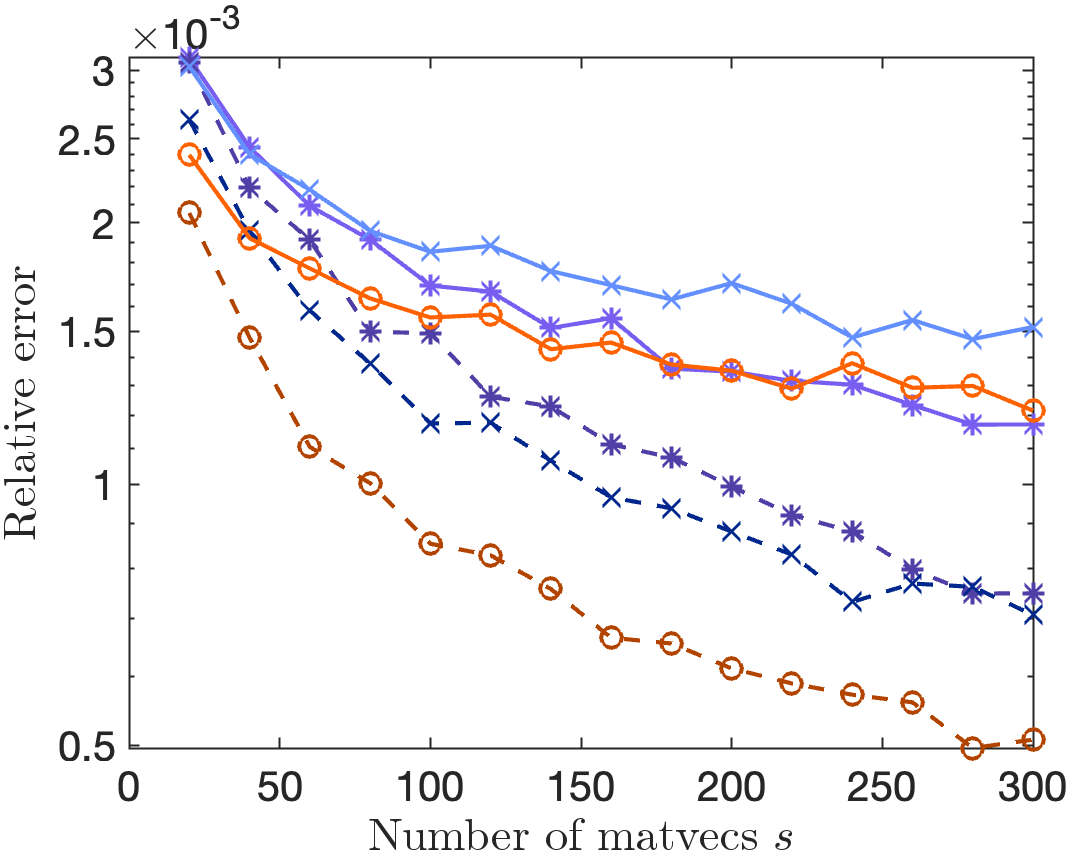}
    \includegraphics[width=0.49\linewidth]{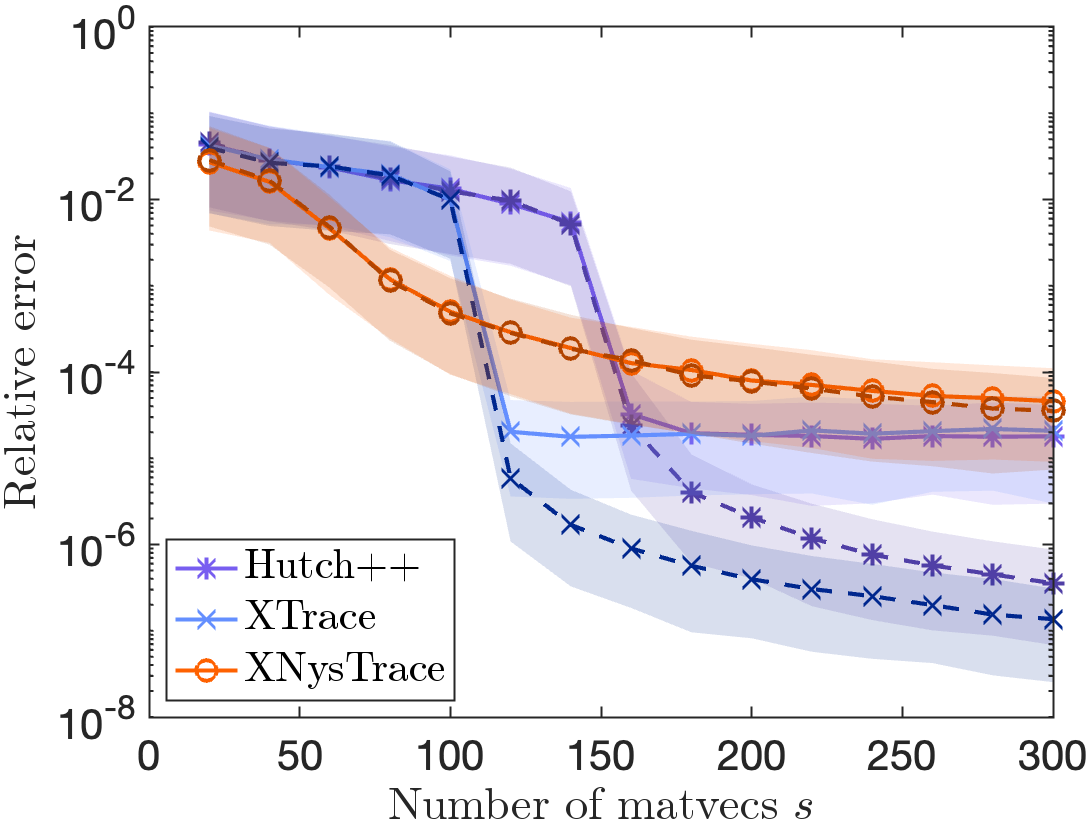}
    \caption[Comparison of \HutchPP, \XTrace, and \XNysTrace trace estimators and their resphered versions on matrices with different spectra]{Comparison of \HutchPP, \XTrace and \XNysTrace algorithms (light solid lines) and their resphered versions (dark dashed lines) on the \texttt{flat} (\cref{eq:flat}, \emph{left}) and \texttt{step} matrices (\cref{eq:step}, \emph{right}).
    Lines show median of 100 trials, and error bars on the right panel show 10\% and 90\% quantiles.
    (Error bars are omitted on the left panel for clarity.)}
    \label{fig:xtrace-resphere}
\end{figure}

\Cref{fig:xtrace-resphere} compares the \HutchPP, \XTrace and \XNysTrace algorithms and their resphered versions on the matrices \texttt{flat} and \texttt{step} defined in \cref{eq:test-matrices}.
We see that resphering significantly improves the performance of all estimators on the \texttt{flat} example and the \XTrace and \HutchPP estimators on the \texttt{step} example.

\section{Application: Estrada index} \label{sec:estrada}


As a running example throughout this part of the thesis, we will apply trace and diagonal estimators to problems in network science.
Given a graph with adjacency matrix $\mat{M}$, the \emph{Estrada index} \cite{Est22} is defined as the trace-exponential of the adjacency matrix:
\begin{equation*}
    \mathrm{estr} \coloneqq \tr(\exp(\mat{M})).
\end{equation*}
The Estrada index is a measure of \emph{centralization} for a graph, that is, how ``clustered'' or ``spread out'' its nodes are.
This quantity is an ideal candidate for variance-reduced trace estimators like \XTrace, as the exponential function promotes spectral decay in the matrix $\mat{A} \coloneqq \exp(\mat{M})$.
The matrix $\mat{A}$ is psd, allowing us to apply estimators designed for psd matrices like \XNysTrace.

\begin{figure}
    \centering
    \includegraphics[width=0.6\linewidth]{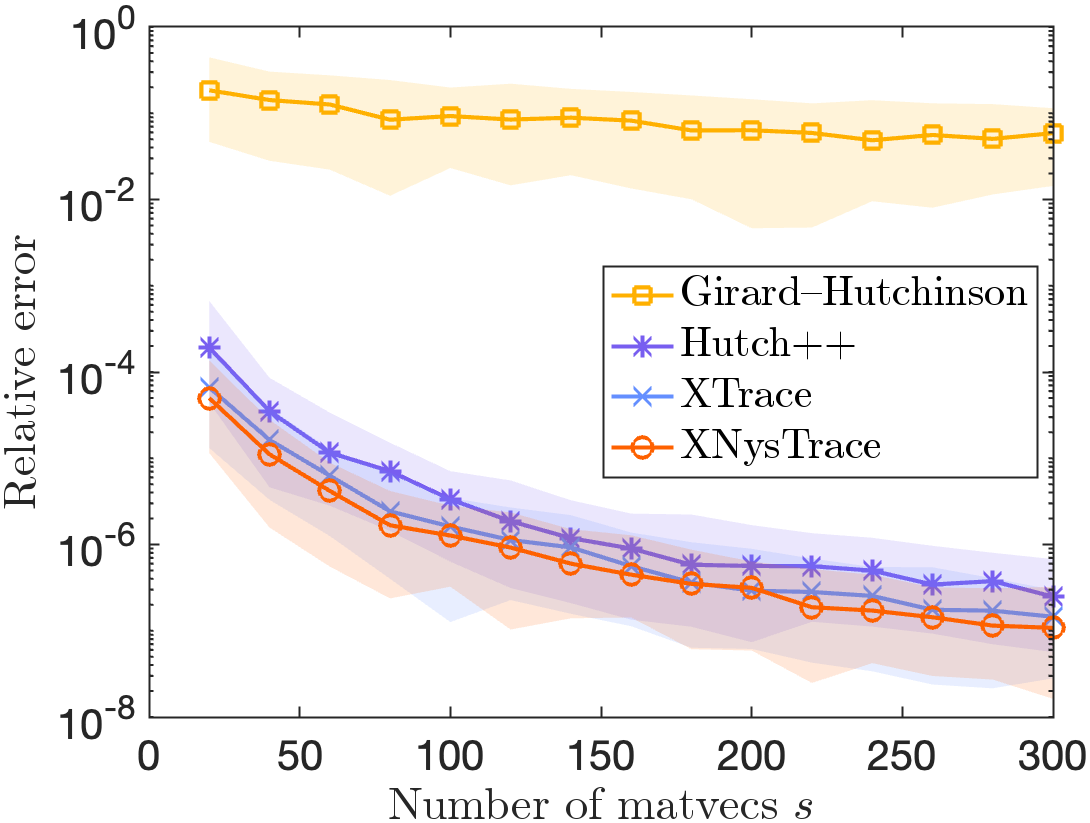}
    \caption[Comparison of Girard--Hutchinson, \HutchPP, \XTrace, and \XNysTrace for estimating the Estrada index]{Relative error of trace estimates by Girard--Hutchinson estimator (yellow squares), \HutchPP (purple asterisks), \XTrace (blue crosses), and \XNysTrace (orange circles) for Estrada index problem as a function of the number of matvecs $s$.
    Lines show median of 100 trials, and error bars show 10\% and 90\% quantiles.}
    \label{fig:estrada}
\end{figure}

\Cref{fig:estrada} provides a demonstration of various trace estimators applied to the Estrada index.
I test on the \texttt{yeast} network $\mat{M} \in \real^{2361\times 2361}$ from \cite{BZC+03}; for this small example, the exact values of the Estrada index can be computed for reference.
To obtain matvecs with $\mat{A} = \exp(\mat{M})$, I employ forty steps of the Lanczos algorithm  \cite[Ch.~6]{Che24}, which is sufficient to compute matvecs to relative error $10^{-11}$.
I implemented all trace estimators with resphering.
The conclusions are much the same as in \cref{fig:trace-comparison,fig:xtrace-resphere}, with ranking from best-to-worst \XNysTrace, \XTrace, \HutchPP, Girard--Hutchinson.
On this example, the difference between \XNysTrace, \XTrace, and \HutchPP is comparatively small, and all of these variance-reduced methods leading to substantial accuracy improvements over the Girard--Hutchinson estimator.

\section{The leave-one-out approach: Summary} \label{sec:loo-summary}

This chapter descrived the leave-one-out\index{leave-one-out randomized algorithm} approach to designing randomized algorithms for matrix attribute estimation.
At a high level, this approach proceeds in five steps:
\begin{enumerate}
    \item Compute a low-rank approximation to the input matrix by multiplying it against a collection of test vectors.
    \item Decompose the quantity of interest into known piece depending on the low-rank approximation plus a residual.
    \item Construct a Monte Carlo estimate of the residual using a single random vector.
    \item Downdate the low-rank approximation by recomputing it with a test vector removed, and use the left-out test vector as the random vector in step 3.
    \item Average the estimator from step 4 over all choices of vectors to leave out.
\end{enumerate}
We employed this approach to derive the \XTrace and \XNysTrace estimators, which improve on \HutchPP by using an exchangeable design.
These algorithms all the matrix--vector products \emph{both} for low-rank approximation \emph{and} residual trace estimation.
We will see more leave-one-out\index{leave-one-out randomized algorithm} randomized matrix algorithms in upcoming chapters.

The leave-one-out\index{leave-one-out randomized algorithm} approach is simple, but deriving efficient implementations using it required some effort.
To derive fast algorithms, our main tools were the downdating formulas for the randomized SVD (\cref{thm:rsvd-downdate}) and randomized Nystr\"om approximation (\cref{thm:nystrom-downdate}).
The resphering technique further improves leave-one-out estimators, at the cost of some additional complexity.

\chapter{More on trace estimation} \label{ch:more-trace}

\epigraph{From 1987 to 2020, an algorithm called Hutchinson’s Estimator was the state-of-the-art for the trace estimation problem, with analysis giving matching upper and lower bounds at $\Theta(1/\varepsilon^2)$ Matrix–Vector products. Notably, this $\Omega(1/\varepsilon^2)$ lower bound was known only for this estimator, and was not known to hold in general \cite{WWZ14} We tried to generalize that lower bound to hold for arbitrary Matrix-Vector algorithms, but that didn't work. While looking into why, we instead found an algorithm that only used $\order(1/\varepsilon)$ Matrix-Vector products, called \HutchPP.}{Raphael A.\ Meyer}

The last chapter introduced the leave-one-out\index{leave-one-out randomized algorithm} approach for designing randomized algorithms for matrix computations.
As a demonstration of that technique, we saw the \XTrace and \XNysTrace algorithms for trace estimation and developed efficient implementations using downdating formulas for the randomized SVD and Nystr\"om approximation.

In this chapter, we will continue our discussion of trace estimation by discussing \emph{a priori} error analysis, \emph{a posteriori} error estimation, adaptive determination of parameters, and alternatives to \XTrace and \XNysTrace.

\myparagraph{Sources}
This chapter is based on the \XTrace paper

\fullcite{ETW24}.

Several aspects of the presentation have been significantly expanded, including the discussion of how to interpret the \XTrace error bounds in \cref{sec:interpreting-hutchpp-bounds}, the discussion of adaptivity in \cref{sec:adaptive-xtrace}, and the discussion of trace estimation alternatives in \cref{sec:trace-alternatives}.

\myparagraph{Outline}
\Cref{sec:trace-error-bounds} presents \emph{a priori} error bounds for \XTrace and \XNysTrace, and \cref{sec:interpreting-hutchpp-bounds} discusses how to interpret them.
\Cref{sec:trace-error-estimation} describes how to estimate the error of \XTrace and \XNysTrace \emph{a posteriori}, and \cref{sec:adaptive-xtrace} explains how these error estimates can be used to implement \XTrace and \XNysTrace adaptively to meet an error tolerance.
Finally, \cref{sec:trace-alternatives} presents alternatives to \XTrace and \XNysTrace for trace estimation.

\section{\emph{A priori} error bounds} \label{sec:trace-error-bounds}

In this section, we introduce \emph{a priori} error bounds for \XTrace and \XNysTrace, analogous to \cref{thm:hutchpp-mse} for \HutchPP.
Versions of these bounds originally appeared as \cite[Thm.~1.1]{ETW24}.
The first result is the following bound for \XTrace:

\begin{theorem}[\XTrace: mean-squared error bound] \label{thm:xtrace-mse}
    Let $\mat{B} \in \real^{n\times n}$ be a \warn{real} square matrix, and consider the \XTrace estimator with test vectors $\vec{\omega}_1,\ldots,\vec{\omega}_{s/2}$ drawn iid from the standard normal $\Normal_\real(\vec{0},\Id_n)$ or sphere $\sqrt{n} \cdot \sphere(\real^n)$.
    This estimator $\hat{\tr}_{\mathrm{X}}$ is an unbiased estimate for the trace $\tr(\mat{B})$ with mean-squared error
    \begin{equation} \label{eq:xtrace-variance}
        \Var(\hat{\tr}_{\mathrm{X}}) \le \min_{r\le s/2-4} \left[\frac{4\e^2 s}{(s/2-r-3)^2} \norm{\mat{B} - \lowrank{\mat{B}}_r}_{\mathrm{F}}^2 + \frac{4s}{s/2-r-3} \norm{\mat{B} - \lowrank{\mat{B}}_r}^2\right].
    \end{equation}
    In particular, the root-mean-squared error decreases at a rate of at least $\order(s^{-1})$:
    \begin{equation} \label{eq:xtrace-1/s}
        [\expect (\hat{\tr}_{\mathrm{X}} - \tr(\mat{B}))^2]^{1/2} \le \frac{\mathrm{const}}{s} \cdot \norm{\mat{B}}_*.
    \end{equation}
\end{theorem}

An analogous result holds for \XNysTrace \cite[Thm.~1.1]{ETW24}:

\begin{theorem}[\XNysTrace: mean-squared error bound] \label{thm:xnystrace-mse}
    Let $\mat{A} \in \real^{n\times n}$ be a \warn{real} psd matrix and consider the \XNysTrace estimator with test vectors $\vec{\omega}_1,\ldots,\vec{\omega}_{s/2}$ drawn iid from the standard normal $\Normal_\real(\vec{0},\Id_n)$ or sphere $\sqrt{n} \cdot \sphere(\real^n)$ distributions.
    Then the \XNysTrace estimator $\hat{\tr}_{\mathrm{XN}}$ provicdes an unbiased estimate for $\tr(\mat{A})$ with root-mean-squared error
    \begin{multline} \label{eq:xnystrace-rmse}
        \Var(\hat{\tr}_{\mathrm{XN}})^{1/2} \le s \min_{r\le s-6} \Bigg[\frac{5\e^2}{(s-r-5)^2} \norm{\mat{A} - \lowrank{\mat{A}}_r}_* \\ + \frac{\sqrt{2}}{(s-r-5)^{3/2}} \norm{\mat{A} - \lowrank{\mat{A}}_r}_{\mathrm{F}} + \frac{\sqrt{8}}{s-r-5} \norm{\mat{A} - \lowrank{\mat{A}}_r}\Bigg].
    \end{multline}
    In particular, the root-mean-squared error decreases at a rate of at least $\order(s^{-1})$:
    \begin{equation} \label{eq:xnystrace-1/s}
        [\expect (\hat{\tr}_{\mathrm{XN}} - \tr(\mat{A})^2]^{1/2} \le \frac{\mathrm{const}}{s} \cdot \tr(\mat{A}).
    \end{equation}
\end{theorem}

These bounds, particularly the \HutchPP-style $\order(1/s)$ ``convergence rates'' of \XTrace \cref{eq:xtrace-1/s} and \XNysTrace \cref{eq:xnystrace-1/s} are easy to misinterpret; we will elaborate more on what these bounds do and do not mean in \cref{sec:interpreting-hutchpp-bounds}.

For now, let us focus on the main error bounds \cref{eq:xtrace-variance,eq:xnystrace-rmse}, which are more informative than the $\order(1/s)$ convergence rate but perhaps harder to interpret.
To help make things simpler, we remind the reader of the comparison \cref{eq:trace-bounds-compare-ch14} above.
Let $\mat{A}$ be a psd matrix with eigenvalues decaying at an exponential rate $\lambda_i(\mat{A}) \le \alpha^i$ for $\alpha \in (0,1)$.
Then the \XTrace (\cref{thm:xtrace-mse}, \XNysTrace (\cref{thm:xnystrace-mse}), and \HutchPP (\cref{thm:hutchpp-mse}) theorems imply that the root-mean-squared errors decay at rates
\begin{equation} \label{eq:trace-exp-convergence}
\begin{split}
	\bigl[\expect \bigl(\hat{\tr}_{\mathrm{H}\texttt{++}} - \tr (\mat{A}) \bigr)^2\bigr]^{1/22} &\le \makebox[\widthof{$s^{1/2}$}][r]{}\,c_1(\alpha)\, \alpha^{\color{red} s/3}; \\ 
    \bigl[\expect \bigl(\hat{\tr}_{\mathrm{X}} - \tr (\mat{A}) \bigr)^2\bigr]^{1/2}&\le s^{1/2} \,c_2(\alpha)\,\alpha^{\color{red} s/2}; \\
	\bigl[\expect \bigl(\hat{\tr}_{\mathrm{XN}} - \tr (\mat{A}) \bigr)^2\bigr]^{1 / 2}&\le \makebox[\widthof{$s^{1/2}$}][r]{$s$} \,c_3(\alpha)\, \alpha^{\color{red} s}.
\end{split}
\end{equation}
Here, the prefactors $c_1(\alpha),c_2(\alpha),c_3(\alpha) > 0$ depend only on $\alpha$.
We see that \XNysTrace and \XTrace converge at $3\times$ and $1.5\times$ the rate of \HutchPP, respectively.
These theoretically predicted convergence rates are mirrored in the \texttt{exp} example in \cref{fig:trace-comparison}.

We will prove the \XTrace bound (\cref{thm:xtrace-mse}) and omit the proof of the \XNysTrace bound (\cref{thm:xnystrace-mse}), which is similar.
The essence of the theorem is contained in the following structural bound.

\begin{lemma}[\XTrace: Structural bound] \label{lem:xtrace-structural}
    Import the setting of \cref{thm:xtrace-mse}.
    The \XTrace estimator is unbiased $\expect[\hat{\tr}_{\mathrm{X}}] = \tr(\mat{B})$ and the variance admits the bound
    \begin{equation*}
        \Var(\hat{\tr}_{\mathrm{X}}) \le \frac{4}{s} \cdot \expect \norm{(\Id - \outprod{\mat{Q}_{(1)}})\mat{B}}_{\mathrm{F}}^2 + 4 \, \expect \norm{(\Id - \outprod{\mat{Q}_{(12)}})\mat{B}}^2.
    \end{equation*}
    Here, $\mat{Q}_{(1)} = \orth(\mat{B}\mat{\Omega}_{-1})$ and $\mat{Q}_{(12)} = \orth(\mat{B} \mat{\Omega}(:,3:s/2))$.
\end{lemma}

Before we prove this lemma, let us use it to establish \cref{thm:xtrace-mse}.

\begin{proof}[Proof of \cref{thm:xtrace-mse}]
    Combining the \XTrace structural bound \cref{lem:xtrace-structural} and the randomized SVD bound (\cref{fact:rsvd-error}) gives
    \begin{multline*}
        \Var(\hat{\tr}_{\mathrm{X}}) \le \frac{4}{s} \min_{r \le s/2-3} \frac{s/2-2}{s/2-r-2} \norm{\mat{B} - \lowrank{\mat{B}}_r}_{\mathrm{F}}^2 \\ + 4 \min_{r\le s/2-4} \frac{s/2+r-3}{s/2-r-3} \left( \norm{\mat{B} - \lowrank{\mat{B}}_r}^2 + \frac{\e^2}{s/2-r-2} \norm{\mat{B} - \lowrank{\mat{B}}_r}_{\mathrm{F}}^2 \right).
    \end{multline*}
    Simplifying, we obtain \cref{eq:xtrace-variance}.
    To prove \cref{eq:xtrace-1/s}, we set $r \coloneqq \lfloor s/4\rfloor - 4$ and invoke \cref{fact:l1l2}.
\end{proof}

We now prove the lemma:

\begin{proof}[Proof of \cref{lem:xtrace-structural}]
    The unbiasedness of the \XTrace estimator follows form its construction; see \cite[p.~18]{ETW24} for a detailed proof.

    By bilinearity of the variance and since every pair $(\hat{\tr}_i,\hat{\tr}_j)$ has the same distribution, we have
    \begin{equation} \label{eq:var-decomp}
        \Var(\hat{\tr}_{\mathrm{X}}) = \frac{2}{s} \cdot \Var(\hat{\tr}_1) + \left(1 - \frac{2}{s}\right)\Cov(\hat{\tr}_1,\hat{\tr}_2).
    \end{equation}
    We will bound the variance $\Var(\hat{\tr}_1)$ and the covariance $\Cov(\hat{\tr}_1,\hat{\tr}_2)$ separately.

    To evaluate $\Var(\hat{\tr}_1)$, use the chain rule for the variance:
    \begin{equation} \label{eq:variance-decomposition-tr1}
        \Var(\hat{\tr}_1) = \expect[\Var(\hat{\tr}_1 \mid \mat{\Omega}_{-1})] + \Var(\expect[\hat{\tr}_1 \mid \mat{\Omega}_{-1}]).
    \end{equation}
    To compute the second term in \cref{eq:variance-decomposition-tr1}, we observe that $\expect[\hat{\tr}_1 \mid \mat{\Omega}_{-1}] = \tr(\mat{B})$.
    That is, even conditional on the other test vectors $\vec{\omega}_2,\ldots,\vec{\omega}_{s/2}$, the first basic \XTrace estimator $\hat{\tr}_1$ is unbiased.
    Thus, the second term in \cref{eq:variance-decomposition-tr1} is zero.
    To evaluate the first term of \cref{eq:variance-decomposition-tr1}, we invoke \cref{fact:gh-variance} to obtain
    \begin{align*}
        \Var(\hat{\tr}_1 \mid \mat{\Omega}_{-1}) &= \Var(\vec{\omega}_i^*(\Id - \outprod{\mat{Q}_{(i)}}) \mat{B}(\Id - \outprod{\mat{Q}_{(i)}})\vec{\omega}_i^{\vphantom{*}} \mid \mat{\Omega}_{-1}) \\ &\le 2 \norm{(\Id - \outprod{\mat{Q}_{(i)}}) \mat{B}(\Id - \outprod{\mat{Q}_{(i)}})}_{\mathrm{F}}^2.
    \end{align*}
    Thus, we have
    \begin{equation} \label{eq:var-bound}
        \Var(\hat{\tr}_1) = \expect[\Var(\hat{\tr}_1 \mid \mat{\Omega}_{-1})] \le 2 \norm{(\Id - \outprod{\mat{Q}_{(i)}}) \mat{B}}_{\mathrm{F}}^2.
    \end{equation}
    We have dropped the second projector $\Id - \outprod{\mat{Q}_{(i)}}$, which can only increase the Frobenius norm.

    Now, we bound $\Cov(\hat{\tr}_1,\hat{\tr}_2)$.
    Introduce
    \begin{equation} \label{eq:xtrace-proof-X-def}
        \mat{X} \coloneqq (\Id - \outprod{\mat{Q}_{(12)}}) \mat{B} (\Id - \outprod{\mat{Q}_{(12)}}).
    \end{equation}
    This matrix is independent of both $\vec{\omega}_1$ and $\vec{\omega}_2$ by design.
    Therefore, we may write
    \begin{align*}
        \Cov(\hat{\tr}_1,\hat{\tr}_2) &= \expect[(\hat{\tr}_1 - \tr(\mat{B}))(\hat{\tr}_2 - \tr(\mat{B}))] \\
        &= \expect[(\hat{\tr}_1 - \tr(\mat{B}) + \tr(\mat{X}) - \vec{\omega}_1^*\mat{X}\vec{\omega}_1^{\vphantom{*}})(\hat{\tr}_2 - \tr(\mat{B}))] \\
        &= \expect[(\hat{\tr}_1 - \tr(\mat{B}) + \tr(\mat{X}) - \vec{\omega}_1^*\mat{X}\vec{\omega}_1^{\vphantom{*}})(\hat{\tr}_2 - \tr(\mat{B}) +  \tr(\mat{X}) - \vec{\omega}_2^*\mat{X}\vec{\omega}_2^{\vphantom{*}})].
    \end{align*}
    The second line follows because $\hat{\tr}_2$ is unbiased approximation to $\mat{B}$, even conditional on $(\vec{\omega}_i : i \ne 2)$.
    The third line follows because the first factor is mean zero, even conditionally on $(\vec{\omega}_i : i \ne 1)$.
    Now, invoke Cauchy--Schwarz and the fact that the two factors are identifically distributed to obtain
    \begin{equation} \label{eq:cov-bound-1}
        \Cov(\hat{\tr}_1,\hat{\tr}_2) \le \expect \left|\hat{\tr}_1 - \tr(\mat{B}) + \tr(\mat{X}) - \vec{\omega}_1^*\mat{X}\vec{\omega}_1^{\vphantom{*}}\right|^2
    \end{equation}
    Conditional on $(\vec{\omega}_i : i \ne 1)$, the terms $\tr(\mat{B})$ and $\tr(\mat{X})$ are constants and the quantity in the absolute value signs is mean-zero.
    Furthermore, the first term of 
    \begin{equation*}
        \hat{\tr}_1 = \tr(\mat{Q}_{(1)}^*\mat{B}\mat{Q}_{(1)}^{\vphantom{*}}) + \vec{\omega}_1^*(\Id - \outprod{\mat{Q}_{(1)}})\mat{B}(\Id - \outprod{\mat{Q}_{(1)}}) \vec{\omega}_1^{\vphantom{*}},
    \end{equation*}
    is also conditionally constant.
    Therefore, we obtain,
    \begin{equation*}
        \Cov(\hat{\tr}_1,\hat{\tr}_2) \le \expect\left[ \Var \left( \vec{\omega}_1^*[(\Id - \outprod{\mat{Q}_{(1)}})\mat{B}(\Id - \outprod{\mat{Q}_{(1)}}) - \mat{X}] \vec{\omega}_1^{\vphantom{*}} \mid \mat{\Omega}_{-1} \right) \right].
    \end{equation*}
    Using \cref{fact:gh-variance}, we bound the right-hand side as
    \begin{equation} \label{eq:cov-bound-2}
        \Cov(\hat{\tr}_1,\hat{\tr}_2) \le 2 \expect \norm{(\Id - \outprod{\mat{Q}_{(1)}})\mat{B}(\Id - \outprod{\mat{Q}_{(1)}}) - \mat{X}}_{\mathrm{F}}^2.
    \end{equation}

    We now do some wrangling to bound the right-hand side of \cref{eq:cov-bound-2}.
    Introduce the projectors $\mat{\Pi}_1 \coloneqq \outprod{\mat{Q}_{(1)}}$ and $\mat{\Pi}_{12} \coloneqq \outprod{\mat{Q}_{(12)}}$, and note the identity
    \begin{equation} \label{eq:projector-product}
        \Id - \mat{\Pi}_1 = (\Id - \Pi_1)(\Id - \Pi_{12}) = (\Id - \Pi_{12})(\Id - \Pi_1).
    \end{equation}
    We compute
    \begin{align*}
        \norm{(\Id - \mat{\Pi}_1)\mat{B}(\Id - \mat{\Pi}_1) - \mat{X}}_{\mathrm{F}}^2
        &= \norm{(\Id - \mat{\Pi}_1)(\Id - \mat{\Pi}_{12})\mat{B}(\Id - \mat{\Pi}_{12})(\Id - \mat{\Pi}_1) - \mat{X}}_{\mathrm{F}}^2\\
        &= \norm{(\Id - \mat{\Pi}_1)\mat{X}(\Id - \mat{\Pi}_1) - \mat{X}}_{\mathrm{F}}^2\\
        &= \norm{(\Id - \mat{\Pi}_1)\mat{X}(\Id - \mat{\Pi}_1) - (\Id - \mat{\Pi}_1)\mat{X}\Id - \mat{\Pi}_1\mat{X}}_{\mathrm{F}}^2 \\
        &= \norm{-(\Id - \mat{\Pi}_1)\mat{X}\mat{\Pi}_1  - \mat{\Pi}_1\mat{X}}_{\mathrm{F}}^2 \\
        &= \norm{(\Id - \mat{\Pi}_1)\mat{X}\mat{\Pi}_1}_{\mathrm{F}}^2 + \norm{\mat{\Pi}_1\mat{X}}_{\mathrm{F}}^2 \\
        &\le \norm{\mat{X}\mat{\Pi}_1}_{\mathrm{F}}^2 + \norm{\mat{\Pi}_1\mat{X}}_{\mathrm{F}}^2 \\
        &= \norm{\mat{X}(\mat{\Pi}_1 - \mat{\Pi}_{12})}_{\mathrm{F}}^2 + \norm{(\mat{\Pi}_1 - \mat{\Pi}_{12})\mat{X}}_{\mathrm{F}}^2 \\
        &= \norm{\mat{X}(\mat{\Pi}_1 - \mat{\Pi}_{12})}^2 + \norm{(\mat{\Pi}_1 - \mat{\Pi}_{12})\mat{X}}^2 \\
        &\le 2 \norm{\mat{X}}^2.
    \end{align*}
    The first line is the identity \cref{eq:projector-product},
    the second line is the definition of $\mat{X}$ \cref{eq:xtrace-proof-X-def},
    the fifth line is the Pythagorean theorem,
    the sixth line is the nonexpansiveness of the orthoprojector $\Id - \mat{\Pi}_1$,
    the seventh line invokes the identity $\mat{X}\mat{\Pi}_{12} = \mat{0}$,
    the eighth line observes that $\mat{X}(\mat{\Pi}_1 - \mat{\Pi}_{12})$ is a rank-one matrix, 
    and the final line holds again because the orthoprojector $\mat{\Pi}_1 - \mat{\Pi}_{12}$ is nonexpansive.
    Combining this result with \cref{eq:cov-bound-2} and using the definition of $\mat{X}$, we obtain
    \begin{equation*}
        \Cov(\hat{\tr}_1,\hat{\tr}_2) \le 4 \expect \norm{(\Id - \outprod{\mat{Q}_{(12)}})\mat{B}(\Id - \outprod{\mat{Q}_{(12)}})}^2.
    \end{equation*}
    Substituting this covariance bound and the variance bound \cref{eq:var-bound} into the variance decomposition \cref{eq:var-decomp} yields the stated result.
\end{proof}

\section{How to interpret the \HutchPP and \XTrace error bounds} \label{sec:interpreting-hutchpp-bounds}

When the \HutchPP paper was released, it generated a large amount of interest, both among specialists in randomized matrix computations and researchers in adjacent areas such as machine learning, quantum information, and scientific computing.
In my view, part of this interest can be attributed to the the \HutchPP paper showing that its method is \emph{optimal}.
Speaking precisely, the paper \cite{MMMW21} shows that \HutchPP is optimal in the following sense:
\actionbox{\textbf{Optimality of \HutchPP.} \HutchPP produces an approximation to the trace of a \warn{psd} matrix $\mat{A}$ with relative error $\varepsilon$ using $s = \order(1/\varepsilon)$ matvecs.
Conversely, there is no algorithm that produces a trace estimate up to relative error $\varepsilon$ \warn{for an arbitrary psd input matrix} in fewer than $s = \order(1/[\varepsilon \log(1/\varepsilon)])$ matvecs.}
This lower bound was later improved to $s = \order(1/\varepsilon)$ \cite{Mey24}, establishing that the $s = \order(1/\varepsilon)$ rate of \HutchPP is optimal \textit{even up to logarithmic factors}.
I consider the matching $\order(1/\varepsilon)$ complexity upper and lower bounds for trace estimation to be one of the most beautiful results in the field of randomized matrix computations.

\begin{figure}[t]
    \centering
    \includegraphics[width=0.95\linewidth]{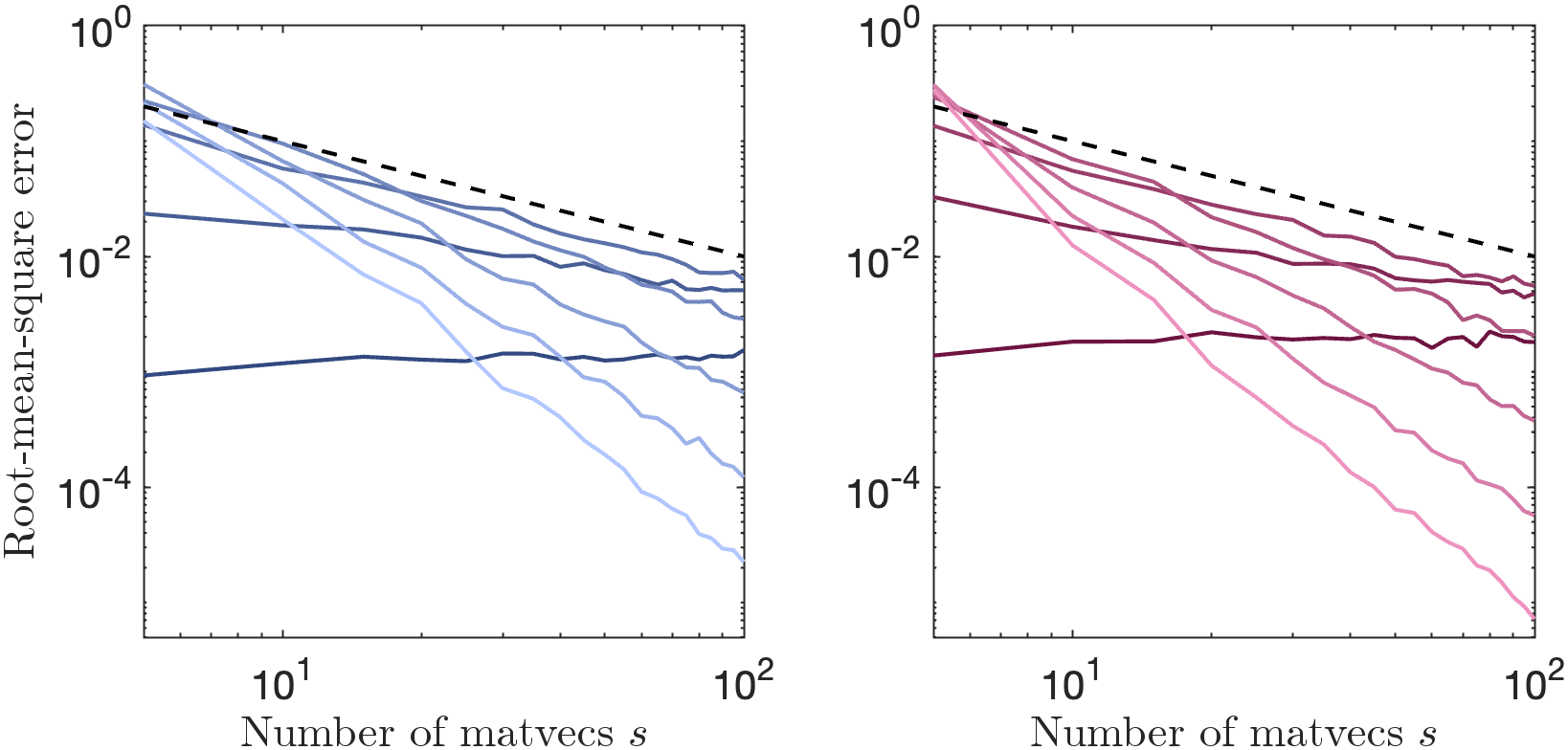}
    \caption[Errors for \HutchPP and \XTrace for matrices with different spectral profiles, compared to the predicted $1/s$ ``convergence rate'']{Root-mean-squared errors for \HutchPP (\emph{left}) and \XTrace (\emph{right}) on the family of psd matrices with eigenvalues $\lambda_i(\mat{A}) = i^{-\alpha}$ for $\alpha=0,0.5,\ldots,3$.
    Higher values of $\alpha$ are indicated using lighter shades, and the the dashed black line shows $\varepsilon = 1/s$.
    Root-mean-squared errors are estimated using 100 trials.}
    \label{fig:hutchpp-xtrace-err-bounds}
\end{figure}

Given these facts, it is very tempting to invert the relation $s \sim 1/\varepsilon$ to claim that the \emph{convergence rate} for \HutchPP is $\varepsilon \sim 1/s$.
Based on this analysis, we would expect that plotting the error versus number of matvecs for \HutchPP on a log--log plot should result in a curve with a slope of $-1$.
\Cref{fig:hutchpp-xtrace-err-bounds} tests this proposition.
Here, we evaluate \HutchPP (\emph{left}) and \XTrace (\emph{right}) on a family of psd matrices of dimension $n = 1000$ with eigenvalues $\lambda_i(\mat{A}) = i^{-\alpha}$ for each fixed $\alpha = 0,0.5,\ldots,3$.
Darker curves show lower values of $\alpha$ and lighter curves show higher values.
We report the average relative error over 100 trials.
For reference, the curve $1/s$ is marked as a black dashed line.
Perhaps surprisingly, the error-versus-matvec curves for both \HutchPP show a number of convergence rates that are both faster \emph{and slower} than the $\mathrm{error} \sim 1/s$ convergence rate we predicted from \HutchPP analysis.
What's going on?

Part of the issue is a cultural difference between how different researchers talk about the complexity or rate of convergence of an algorithm.
\begin{itemize}
    \item Under one view, which might be labeled the \textbf{\textit{computational mathematics view}}, the performance of an algorithm is determined by its \emph{convergence rate}.
    In the context of trace estimate, the convergence rate is the number $\beta$ such that the error $\varepsilon$ scales like $s^{\beta}$ (or better).
    Equivalently, the convergence rate $\beta$ is gives upper bound on the slope of a log--log plot of the error plotted against the number of matvecs.
    Under this interpretation, an algorithm has a convergence rate of $1/s$ only if doubling $s$ halves the error (or better).
    As \cref{fig:hutchpp-xtrace-err-bounds} shows, the convergence rate of \HutchPP and \XTrace can be both slower and faster than $1/s$ for different input matrices.
    \item A second perspective, which we might label as the \textbf{\textit{theoretical computer science view}}, is concerned with the amount of computational resources (in this case, matvecs $s$) required to achieve an of error $\varepsilon$ for a \emph{worst-case input matrix}. 
    Under this viewpoint, the complexity of an algorithm is $s = \order(1/\varepsilon)$ if the error-versus-matvec curve is \emph{bounded above} by a curve of the form ``$\mathrm{const}/s$'' for all $s\le n$.
    It is in this second sense in which \HutchPP is optimal and has a complexity of $s = \order(1/\varepsilon)$.
\end{itemize}
Having understood these two viewpoints, we reconcile the perhaps surprising results in \cref{fig:hutchpp-xtrace-err-bounds} .
The error analysis of \HutchPP and \XTrace does \emph{not} say that their \emph{convergence rate} is $\order(1/s)$, just that the error-versus-matvec curve is bounded above by a curve of the form ``$\mathrm{const}/s$'' for all $s\le n$.
Indeed, this is exactly the case in \cref{fig:hutchpp-xtrace-err-bounds}: The error curves for both \HutchPP and \XTrace are bounded above by $3/s$ for each of the seven test matrices.

This understanding of the \HutchPP and \XTrace error bounds delivers both good news and bad news about how these algorithms perform.

\myparagraph{The good: Spectral convergence}
On the positive side, \HutchPP and \XTrace often converge \emph{much faster} than the $\order(1/s)$ analysis would suggest.
In fact, they converge \emph{spectrally}, with the trace errors decreasing at a rate proportional to the error of the best rank-$\order(s)$ approximation, measured in an appropriate norm.
This phenomenon was demonstrated theoretically in \cref{eq:trace-exp-convergence}, which shows that the variance-reduced trace estimators converge \emph{geometrically} when the eigenvalues decay geometrically.
Spectral convergence is also demonstrated empirically in \cref{fig:hutchpp-xtrace-err-bounds}, where several curves are seen to decrease significantly more rapidly then the dashed $\order(1/s)$ rate.
This spectral convergence means that stochastic trace estimators are very effective at estimating the trace of a matrix with rapidly decaying eigenvalues, such a matrix of the form $\exp(-\beta \mat{A})$ arising from partition function calculations in quantum physics.

\myparagraph{The bad: \HutchPP and \XTrace are not always better than Girard--Hutchinson}
The Girard--Hutchinson \warn{converges at} a rate of $\order(1/s^{1/2})$, in both the ``computational math'' and ``theoretical computer science'' senses.
Given the $s = \order(1/\varepsilon)$ complexity of \HutchPP/\XTrace and the $s = \order(1/\varepsilon^2)$ complex of the Girard--Hutchinson estimator, one might surmise that \HutchPP or \XTrace will always achieve lower error than the Girard--Hutchinson estimator, except possibly for very small values of $s$.
Unfortunately, this is not the case.

\actionbox{For problems with sufficiently slow spectral decay, the Girard--Hutchinson estimator will be (slightly) more accurate than \HutchPP and \XTrace.}

Let us justify this claim by applying both \HutchPP and the Girard--Hutchinson estimator to the identity matrix $\mat{A} = \Id$.
To yield generalizable conclusions, we use \emph{Gaussian} test vectors $\vec{\omega} \sim \Normal_\real(\vec{0},\Id)$---these are a bad choice for computation in practice, but they suffice to illustrate the conceptual point.
By \cref{fact:gh-variance}, the variance of the Girard--Hutchinson estimator for this problem is 
\begin{equation*}
    \Var(\hat{\tr}_{\mathrm{GH}}) = \frac{2}{s} \norm{\Id}_{\mathrm{F}}^2 = \frac{2n}{s}.
\end{equation*}
Consequently, the root-mean-squared \warn{relative} error is 
\begin{equation*}
    \left[\expect \left( \frac{\tr(\Id) - \hat{\tr}_{\mathrm{GH}}}{\tr(\Id)} \right)^2\right]^{1/2} = \frac{\sqrt{2}}{s^{1/2}} \cdot \frac{1}{n^{1/2}}.
\end{equation*}
Now, let us turn to \HutchPP.
For the matrix $\mat{Q}$ produced by the \HutchPP algorithm, the matrix $(\Id-\mat{Q}\mat{Q}^*)\Id(\Id - \mat{Q}\mat{Q}^*)$ is a compression of the identity matrix to a uniformly random subspace of dimension $n - s/3$. 
Thus,
\begin{equation*}
    \Var(\hat{\tr}_{\mathrm{H}\texttt{++}}) = \frac{2}{s/3} \norm{(\Id-\mat{Q}\mat{Q}^*)\Id(\Id - \mat{Q}\mat{Q}^*)}_{\mathrm{F}}^2 = \frac{6(n-s/3)}{s}.
\end{equation*}
Ergo, the root-mean-squared relative error for \HutchPP is
\begin{equation*}
    \left[\expect \left( \frac{\tr(\Id) - \hat{\tr}_{\mathrm{H}\texttt{++}}}{\tr(\Id)} \right)^2\right]^{1/2} = \frac{\sqrt{6}}{s^{1/2}} \cdot \frac{(n-s/3)^{1/2}}{n}.
\end{equation*}
We observe two things: First, the \warn{convergence rate} of \HutchPP on this example is $1/s^{1/2}$, the same as the Girard--Hutchinson estimator.
Second, except for large values $s \approx 3n$, the root-mean-squared error for \HutchPP is larger than the Girard--Hutchinson estimator by about $\sqrt{3}$.
This simple computation demonstrates that, for problems with slow spectral decay, the Girard--Hutchinson estimator can be better than \HutchPP by a constant factor.
The same qualitative conclusions hold for \XTrace and \XNysTrace.

\section{Posterior error estimation} \label{sec:trace-error-estimation}

In the past sections, we studied \emph{a priori} bounds on the error for trace estimation algorithms.
Sharp forms of these error bounds are typically inaccessible in applications, as they depend on the singular values of the matrix $\mat{B}$.
As such, it can be helpful to have \emph{posterior} bounds on the error for trace estimation algorithms.
This section develops posterior error estiamtes for the \XTrace and \XNysTrace algorithms.

Our approach to error estimation is straightforward.
The \XTrace and \XNysTrace algorithms generate a \emph{family} of trace estimates $\hat{\tr}_i$ for $i=1,\ldots,\ell$ with $\ell = s/2$ (\XTrace) or $\ell = s$ (\XNysTrace).
The scaled standard deviation of these estimates serves as a natural posterior estimate for the error $|\hat{\tr}_{\mathrm{X}} - \tr(\mat{B})|$ or $|\hat{\tr}_{\mathrm{XN}} - \tr(\mat{A})|$.
Specifically, we define the error estimate 
\begin{equation} \label{eq:trace-error-estimate}
    \smash{\hat{\Err}}^2 \coloneqq \frac{1}{\ell(\ell-1)} \sum_{i=1}^\ell |\hat{\tr}_i - \hat{\tr}|^2 \quad \text{where } \hat{\tr} = \frac{1}{\ell} \sum_{i=1}^\ell \hat{\tr}_i.
\end{equation}
We have the following result:

\begin{proposition}[\XTrace and \XNysTrace: Posterior error estimate] \label{prop:trace-posterior}
    For both algorithms, the posterior error estimate \cref{eq:trace-error-estimate} satisfies
    \begin{equation*}
        \expect\big[\smash{\hat{\Err}}^2\big] = \frac{1 - \Cor(\hat{\tr}_1,\hat{\tr}_2)}{1+(\ell-1)\Cor(\hat{\tr}_1,\hat{\tr}_2)} \cdot \expect \big( {\tr(\mat{B}) - \hat{\tr}} \big)^2,
    \end{equation*}
    where $\Cor$ is the correlation.
    Moreover, if $\mat{B}$ is a \warn{real} symmetric matrix, the columns $\vec{\omega}_1,\ldots,\vec{\omega}_\ell$ are iid standard Gaussian vectors, and the algorithm is \XTrace, then
    \begin{equation} \label{eq:cor-bound-trace}
        \Cor(\hat{\tr}_1,\hat{\tr}_2) \le 2\frac{\expect \norm{\smash{(\Id - \outprod{\mat{Q}_{(12)}})\mat{B}(\Id - \outprod{\mat{Q}_{(12)}})}}^2\vphantom{\mat{Q}_{(12)}}}{\expect \norm{\smash{(\Id - \outprod{\mat{Q}_{(1)}})\mat{B}(\Id - \outprod{\mat{Q}_{(1)}})}}_{\mathrm{F}}^2}.
    \end{equation}
    Here, $\mat{Q}_{(1)}$ and $\mat{Q}_{(12)}$ were defined in \cref{lem:xtrace-structural}.
\end{proposition}

See \cite[\S5.6]{ETW24} for a proof.

In practice, we find that the \XTrace and \XNysTrace estimates have a small positive correlation.
For example, \cite[\S5.6]{ETW24} reports empirical correlations of at most $0.06$ on the \texttt{exp} example \cref{eq:exp}.
In this typical setting, \cref{prop:trace-posterior} shows that $\hat{\Err}$ provides a slight \emph{underestimate} of the error, on average.

The second result of \cref{prop:trace-posterior} provides some control on the size of these correlations in the case where $\mat{B}$ is real and symmetric and the test vectors are Gaussian.
In particular, if the matrix has a \emph{slow} rate of spectral decay, then the Frobenius norm in the denominator of \cref{eq:cor-bound-trace} denominates the spectral norm in the numerator, so the correlations are small.

\begin{figure}
    \centering
    \includegraphics[width=0.99\linewidth]{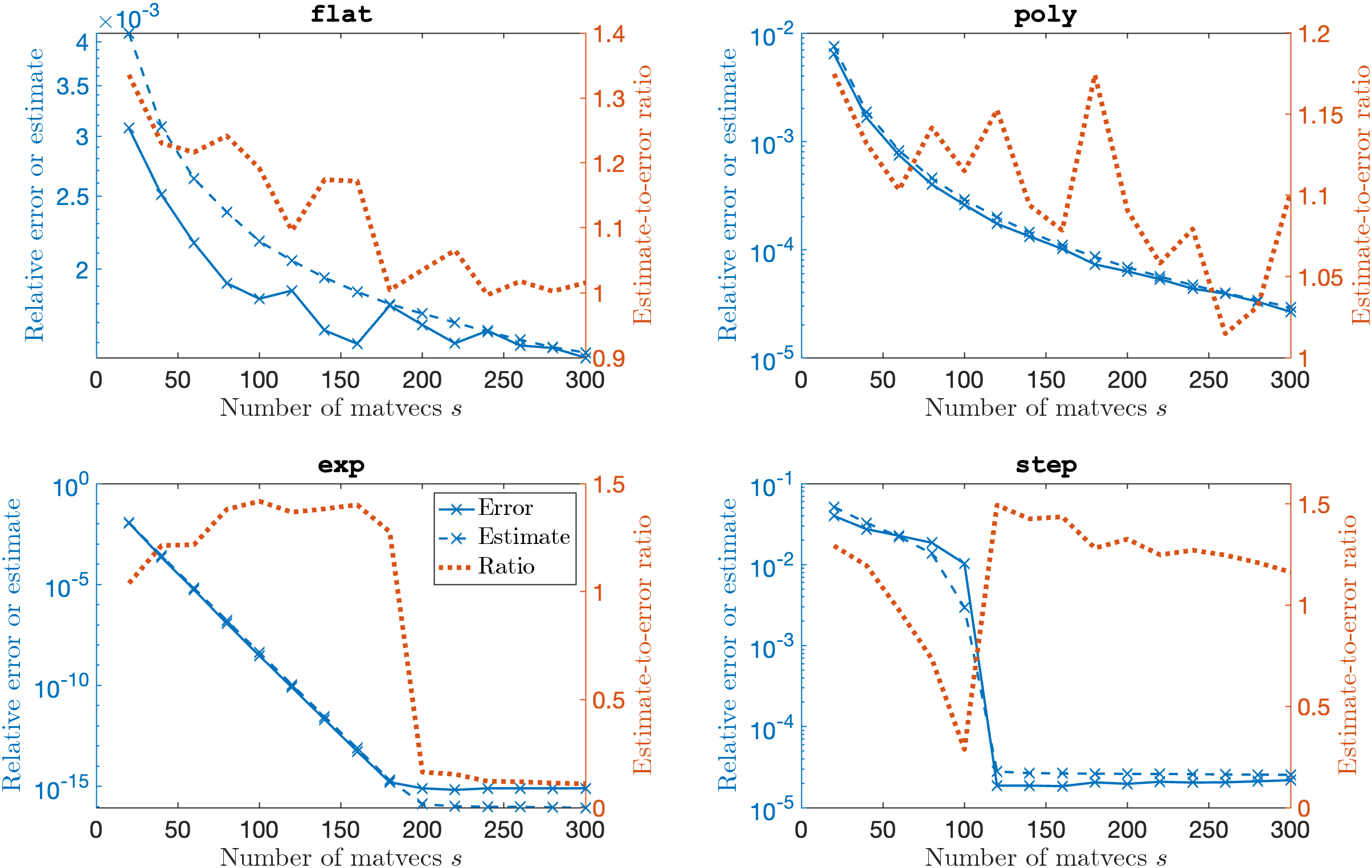}
    \caption[Errors and error estimates for \XTrace applied to matrices with different spectra]{Left-hand axes on each subplot show the error (blue solid crosses) and the median error estimate (blue dashed crosses) for \XTrace on the four test matrices \cref{eq:test-matrices}. 
    Right-hand axes show the ratio of the median error estimate to the median error (red dotted).
    Medians are computed using 1000 trials.}
    \label{fig:trace-error-estimate}
\end{figure}

\Cref{fig:trace-error-estimate} compares the actual error to the error estimate for the \XTrace algorithm applied to the four test matrices \cref{eq:test-matrices}.
Across each test matrix for most numbers of matvecs $s$, the error estimate closely tracks the true error, agreeing up to a factor of $1.5\times$.
There are two situations where the error estimate performs somewhat less effectively:
The first is at the value of $s=100$ matvecs for the step example.
This example confounds the error estimator because there is a large gap between the eigenvalues $\lambda_{s/2}(\mat{A})$ and $\lambda_{s/2+1}(\mat{A})$, causing the basic \XTrace estimators to become more correlated, underestimating the error by a factor of about $0.3$.
The second example occurs on the \texttt{exp} example for large values of $s$.
Here, the trace estimates are so accurate that the size of rounding errors comes into play, and the error estimate is less quantitatively sharp.
This is a non-issue in practice (in double precision computation, at least), because the error estimate correctly diagnoses that the trace has been computed up to ``machine accuracy''.

\section{Adaptivity} \label{sec:adaptive-xtrace}

The \XTrace or \XNysTrace algorithms can be implemented to adaptively determine the number of matvecs required to achieve an error tolerance.
The basic scheme is simple.
Begin by running the \XTrace or \XNysTrace algorithm with some initial number of matvecs $s$ and compute the error estimate \cref{eq:trace-error-estimate}.
If the error estimate is below the tolerance, output the estimated trace.
Otherwise, collect an additional $\Delta s$ matvecs, set $s \gets s + \Delta s$, and update estimates with the newly collected vectors, repeating as many times as necessary to achieve the error tolerance.

A basic version of this scheme with $\Delta s \coloneqq s$ was proposed in \cite[\S3.2]{ETW24}.
With this approach, the number of matvecs is doubled every time the error tolerance is met.
This approach can be wasteful, as it may collect roughly twice as many matvecs as needed.
On the other hand, it maintains the postprocessing cost of $\order(s^2n)$ of the original algorithm, since
\begin{equation*}
    s^2n + (s/2)^2n + (s/4)^2n + \cdots = \order(s^2n).
\end{equation*}

A more economical approach is to use a constant increment $\Delta s$, say, $\Delta s = 10$.
Implemented directly, this approach requires a post-processing cost of $\order(s^3n)$ operations, since
\begin{equation*}
    s^2n + (s-10)^2n + (s-20)^2n + \cdots = \order(s^3n).
\end{equation*}
However, the cost of the algorithm can be reduced to $\order(s^2n)$ by carefully updating the matrix factorizations as new matvecs are collected.
In particular, \XTrace requires an incremental version of the \QR decomposition, discussed in \cref{app:incremental-qr}.
An efficient adaptive version of \XTrace (using the \texttt{scipy} routine \texttt{qr\_update} for the incremental \QR functionality) was developed by Matt Piekenbrock and is implemented in the \texttt{primate} software \cite{Pie24}.

\section{Alternatives to \XTrace and \XNysTrace} \label{sec:trace-alternatives}

\XTrace and \XNysTrace are very effective trace estimators for \emph{truly black-box trace estimation} of matrices enjoying \emph{rapidly decaying singular values}.
However, in many applications, we have additional structure or lack decaying singular values.
In such cases, there may be more effective trace estimators, which we briefly review in this section.

\subsection{Adaptive \HutchPP} \label{sec:adaptive-hutchpp}

The adaptive \HutchPP algorithm was introduced by Persson, Cortinovis, and Kressner \cite{PCK22} as a way of both solving the \emph{apportionment} problem of \HutchPP (how should matvecs be allocated between low-rank approximation and residual trace estimation?) and determining when to stop the algorithm.
The core of their idea is that the number of matvecs needed to perform trace estimation on the residual matrix $(\Id - \mat{Q}\mat{Q}^*)\mat{A}(\Id - \mat{Q}\mat{Q}^*)$ is determined by the Frobenius norm of that matrix (in view of \cref{fact:gh-variance} and related results).
Thus, the number of additional residual trace estimation matvecs can be estimated at runtime using a randomized norm estimator.

The adaptive \HutchPP algorithm gives a natural strategy for improving the \HutchPP algorithm, and it is more straightforward to implement than adaptive versions of \XTrace.
However, in experiments reported in \cite{ETW24}, we find that adaptive \HutchPP produces less accurate results than \XTrace does for a fixed budget of matvecs, with a few exceptions.

\subsection{Stochastic Lanczos quadrature}

Consider the task of computing $\tr f(\mat{A})$ for a Hermitian matrix $\mat{A}$.
The starting point for stochastic Lanczos quadrature is the observation that the trace is an integral
\begin{equation*}
    \tr(f(\mat{A})) = \int_\real f(x) \, \d \mu_{\mat{A}}(x),
\end{equation*}
where $\mu_{\mat{A}} =\sum_{i=1}^n \delta_{\lambda_i(\mat{A})}$ denotes the (unnormalized) \emph{spectral measure} of $\mat{A}$.
Thus, computing trace-functions of $\mat{A}$ can be reformulated as a quadrature problem to approximate integral with respect to the measure $\mu_{\mat{A}}$.

Our first move to solve this problem will be to introduce a stochastic approximation for $\mu_{\mat{A}}$.
For a vector $\vec{\omega} \in \field^n$, we define the \emph{eigenvector spectral measure} of $\mat{A}$ and $\vec{\omega}$ to be
\begin{equation*}
    \mu_{\mat{A},\vec{\omega}} = \sum_{i=1}^n |\vec{u}_i^*\vec{\omega}|^2 \, \delta_{\lambda_i(\mat{A})} \quad \text{where } \mat{A} = \sum_{i=1}^n \lambda_i(\mat{A}) \outprod{\vec{u}_i} \text{ is a spectral decomposition}.
\end{equation*}
The eigenvector spectral measure of $\mat{A}$ and $\vec{\omega}$ has the property that
\begin{equation*}
    \vec{\omega}^*f(\mat{A})\vec{\omega} = \int_\real f(x) \, \d\mu_{\mat{A},\vec{\omega}}(x).
\end{equation*}

Stochastic Lanczos quadrature is based on the following observation: If $\vec{\omega}$ is an \emph{isotropic} random vector, then the \warn{random measure} $\mu_{\mat{A},\vec{\omega}}$ is an unbiased approximation to the spectral measure $\mu_{\mat{A}}$, $\expect[\mu_{\mat{A},\vec{\omega}}] = \mu_{\mat{A}}$.
(In general, defining the expectation of a random measure requires care; no issues arise here because the measures $\mu_{\mat{A},\vec{\omega}}$ and $\mu_{\mat{A}}$ are both supported on the spectrum of $\mat{A}$, a finite set.)
Thus, given approximate eigenvector spectral measures $\hat{\mu}_{\mat{A},\vec{\omega}_i} \approx \mu_{\mat{A},\vec{\omega}_i}$ for iid isotropic vectors $\vec{\omega}_1,\ldots,\vec{\omega}_k$, we can fashion a Monte Carlo approximation to the spectral measure $\mu_{\mat{A}}$ by averaging, $\mu_{\mat{A}} \approx k^{-1} \sum_{i=1}^k \mu_{\mat{A},\vec{\omega}_i}$.

To approximate the eigenvector spectral measure $\mu_{\mat{A},\vec{\omega}}$ of the pair $\mat{A}$ and $\vec{\omega}$, we can deploy the Lanczos algorithm \cite{Lan50,Che24}.
Given an initial vector $\vec{\omega}$, the Lanczos algorithm computes an orthonormal basis $\mat{Q} \in \field^{n\times t}$ for the \emph{Krylov subspace} 
\begin{equation*}
    \set{K}_t(\mat{A},\vec{\omega}) \coloneqq \operatorname{span} \{ \vec{\omega}, \mat{A}\vec{\omega},\mat{A}^2\vec{\omega},\ldots,\mat{A}^{t-1}\vec{\omega} \}
\end{equation*}
and a Hermitian tridiagonal matrix $\mat{T}$ such that
\begin{equation*}
    \mat{Q}^*\mat{A}\mat{Q} = \mat{T},
\end{equation*}
where $\mat{Q}(:,1) = \vec{\omega}/\norm{\vec{\omega}}$. 
We have the following amazing fact:

\begin{fact}[Lanczos and quadrature]
    Introduce $\hat{\mu}_{\mat{A},\vec{\omega}} \coloneqq \norm{\vec{\omega}}^2 \cdot \mu_{\mat{T},\evec_1}$, where $\mat{T}$ is the tridiagonal matrix produced by $t$ steps of the Lanczos procedure.
    Then $\hat{\mu}_{\mat{A},\vec{\omega}}$ is the unique $t$-point Gaussian quadrature measure for the eigenvector spectral measure $\mu_{\mat{A},\vec{\omega}}$.
    That is, 
    \begin{equation*}
        \vec{\omega}^*p(\mat{A})\vec{\omega} = \int_\real p(x) \, \d\mu_{\mat{A},\vec{\omega}}(x) = \int_\real p(x) \, \d\hat{\mu}_{\mat{A},\vec{\omega}}(x)
    \end{equation*}
    for all polynomial $p$ of degree $<2t$.
\end{fact}

See \cite[Ch.~2--3]{Che24}, \cite[\S3]{Tro20a}, and \cite{GM10} for more details on the Lanczos method and the proof of this result.

Stochastic Lanczos quadrature consists of two levels of approximation
\begin{equation*}
    \mu_{\mat{A}} \approx \frac{1}{k} \sum_{i=1}^k \mu_{\mat{A},\vec{\omega}_i}  \approx \hat{\mu}_{\mat{A}} \coloneqq \frac{1}{k} \sum_{i=1}^k \mu_{\mat{T}_i,\evec_1}.
\end{equation*}
First, we approximate the spectral measure $\mu_{\mat{A}}$ by an average of eigenvector spectral measures $\mu_{\mat{A},\vec{\omega}_i}$ generated by iid isotropic vectors $\vec{\omega}_i$.
Second, we approximate each $\mu_{\mat{A},\vec{\omega}_i}$ by a quadrature approximation $\mu_{\mat{T}_i,\evec_1}$ generated from the tridiagonal matrix $\mat{T}_i$ obtained from applying the Lanczos method to each $\vec{\omega}_i$.
Once formed, the \emph{stochastic Lanczos quadrature} approximation
\begin{equation*}
    \hat{\mu}_{\mat{A}} \coloneqq \frac{1}{k} \sum_{i=1}^k \mu_{\mat{T}_i,\evec_1} = \sum_{i=1}^k \sum_{j=1}^t \frac{|\vec{u}_{i,j}^*\vec{\omega}_i|^2}{k} \delta_{\lambda_j(\mat{T}_i)} \approx \mu_{\mat{A}}
\end{equation*}
can be use to compute the trace of \emph{any} trace-function of $\mat{A}$ at a minimal additional cost of $\order(k\cdot t)$ operations:
\begin{equation*}
    \tr(f(\mat{A})) \approx \int_\real f(x) \, \d\hat{\mu}_{\mat{A}}(x) = \sum_{i=1}^k \sum_{j=1}^t \frac{|\vec{u}_{i,j}^*\vec{\omega}_i|^2}{k} f(\lambda_j(\mat{T}_i)).
\end{equation*}
Here, $\mat{T}_i = \sum_{j=1}^t \lambda_j(\mat{T}_i) \outprod{\vec{u}_{i,j}}$ denotes a spectral decomposition of each $\mat{T}_i$.

To refine the stochastic Lanczos quadrature approximation, one must increase both the number of isotropic vectors $k$ and the number of Lanczos steps $t$.
The rate of convergence is Monte Carlo in $k$, and the rate of convergence in $t$ depends on the regularity of the function $f$.
For analytic functions, the convergence is geometric \cite[\S3.2]{CTU25}.

Stochastic Lanczos quadrature has advantages and disadvantages for computing $\tr(f(\mat{A}))$ over variance-reduced trace estimators like \HutchPP and \XTrace.
As a disadvantage, the (root-mean-square) error of stochastic Lanczos quadrature converges at the un-accelerated Monte Carlo rate $\order(k^{-1/2})$, the same as the ordinary non-variance-reduced Girard--Hutchinson estimator.
For problems with decaying eigenvalues, \HutchPP and \XTrace can converge must faster.
As an advantage, once the stochastic Lanczos quadrature rule has been computed, one can estimate $\tr(f(\mat{A}))$ for as many different functions $f$ as desired.
This feature contrasts with \HutchPP and \XTrace, which must be applied separately for each function $f$.

\subsection{Krylov-aware stochastic trace estimation}

The virtues of \HutchPP and stochastic Lanczos quadrature are combined in the \emph{Krylov-aware trace estimator} of Chen and Hallman \cite{CH23a}.
As with stochastic Lanczos quadrature, their algorithm applies to functions $f(\mat{A})$ of a \warn{Hermitian} matrix $\mat{A}$.
Fix a block size $b$, number of Lanczos steps $t$, rank $k$, and number of Girard--Hutchinson steps $m$.
The Krylov-aware trace estimator begins by applying (roughly) \warn{$\lceil k/b \rceil+t$ steps} of the \warn{block} Lanczos method \cite{GU77} to a matrix $\mat{\Omega}$ with iid standard Gaussian columns to build an approximation
\begin{equation*}
    f(\mat{A}) \approx \mat{Q}f(\mat{T})\mat{Q}^*.
\end{equation*}
This approximation is then truncated to rank $k$ by symmetrically projecting onto the matrix $\mat{Q}_k \coloneqq \mat{Q}(:,1:k)$ containing the first $k$ columns of $\mat{Q}$:
\begin{equation*}
    f(\mat{A}) \approx \mat{M} \coloneqq \outprod{\mat{Q}_k}\cdot \mat{Q}f(\mat{T})\mat{Q}^*\cdot \outprod{\mat{Q}_k} = \mat{Q}_k^{\vphantom{*}} [f(\mat{T})](1:k,1:k)\mat{Q}_k^*.
\end{equation*}
Then, they estimate the residual trace $\tr(f(\mat{A}) - \mat{M})$ using the resphered Girard--Hutchinson estimator, resulting in the estimate
\begin{equation*}
    \hat{\tr}_{\mathrm{KA}} = \tr(\mat{M}) + \frac{n-k}{m} \sum_{i=1}^m \frac{\vec{\nu}_i^*(f(\mat{A})\vec{\nu}_i^{\vphantom{*}})}{\vec{\nu}_i^*\vec{\nu}_i^{\vphantom{*}}} \:\: \text{for } \vec{\nu}_i = (\Id - \outprod{\mat{Q}_k})\vec{\gamma}_i, \: \vec{\gamma}_i \simiid \Normal_\real(\vec{0},\Id).
\end{equation*}
Matvecs $f(\mat{A})\vec{\nu}_i$ in this display are performed using $t$ steps of the Lanczos algorithm.
See \cite[\S3]{CH23a} for details and \cite{PCM25} for analysis of the approximation $\mat{M}$.

The Krylov-aware trace estimator and \XTrace accelerate the \HutchPP method in different ways.
The Krylov-aware estimator consolidates the computation of matrix--vector products $f(\mat{A})\mat{\Omega}$ and $f(\mat{A})\mat{Q}$ into a single invokation of the block Lanczos algorithm, but retains the need to perform matrix--vector products $f(\mat{A})\vec{\nu}_i$ for residual trace estimation.
\XTrace still needs two batches $f(\mat{A})\mat{\Omega}$ and $f(\mat{A})\mat{Q}$ of matvecs to form a low-rank approximation of $f(\mat{A})$ but removes the need for additional matvecs for residual trace estimation.

The main advantage of the Krylov-aware trace estimator over \XTrace and \XNysTrace is that it can be used to evaluate $\tr(f(\mat{A}))$ for a family of functions $f$ at minimal additional effort.
For this reason, the Krylov-aware trace estimator is often the best algorithm for the $\tr(f(\mat{A}))$ problem.

\subsection{Probing methods}

Another family of methods for trace estimation are based on \emph{probing}.
Here, is one version of the idea. 
We know that the Girard--Hutchinson estimator with random sign vectors has a variance depending only on the sum of the off-diagonal entries:
\begin{equation*}
    \Var(\vec{\omega}^*\mat{A}\vec{\omega}) = 2 \sum_{i\ne j} a_{ij}^2 \quad \text{for } \mat{A} \in \real^{n\times n} \text{ symmetric and } \vec{\omega} \sim \Unif \{\pm 1\}^n.
\end{equation*}
Let us use this fact to our advantage.
Suppose that we are able to partition the index set $\{1,\ldots,n\}$ as a disjoint union of a small number of sets $\set{S}_1,\ldots,\set{S}_m$ with the property that the off-diagonal entries of each submatrix $\mat{A}(\set{S}_i,\set{S}_i)$ are small.
Using such an approximation, we can obtain an accurate approximation to the trace by applying the Girard--Hutchinson estimator separately to each submatrix $\mat{A}(\set{S}_i,\set{S}_i)$ where the random vectors $\vec{\omega}$ are supported on each partition $\set{S}_i$.
In particular, we can define the \emph{stochastic probing estimator}:
\begin{equation*}
    \hat{\tr}_{\mathrm{SP}} \coloneqq \frac{1}{t} \sum_{i=1}^m \sum_{j=1}^t \vec{\omega}_{ij}^*\mat{A} \vec{\omega}_{ij}^{\vphantom{*}} \quad \text{with } \begin{cases}
        \vec{\omega}_{i1}(\set{S}_i),\ldots,\vec{\omega}_{it}(\set{S}_i) \simiid \Unif \{\pm 1\}^{\set{S}_i}, \\
        \vec{\omega}_{i1}(\overline{\set{S}_i})=\cdots=\vec{\omega}_{it}(\overline{\set{S}_i}) = \vec{0}.
    \end{cases} 
\end{equation*}
Here, $t$ denotes the number of matvecs per set $\set{S}_i$, and $\overline{\set{S}_i} \coloneqq \{1,\ldots,n\} \setminus \set{S}_i$ denotes the complement of each $\set{S}_i$.

Stochastic probing methods were initially developed in quantum chromodynamics \cite{MBF+11,BBC+12}, and they reached computational mathematics shortly thereafter \cite{ASE14,GSO17}.
There are also deterministic versions of probing \cite{TS12}.
See the paper \cite{FRS25} for analysis of stochastic probing methods and a literature review on probing.

Probing methods require a good partition $\set{S}_1,\ldots,\set{S}_m$, which demands prior information about the matrix $\mat{A}$.
In a typical application of probing, one is interested in evaluating the trace (or diagonal) of $\mat{A} = f(\mat{H})$ for a sparse matrix $\mat{H}$.
In this setting, one can construct a partitioning $\set{S}_1,\ldots,\set{S}_m$ from a (distance-$d$) graph coloring of the adjacency graph of $\mat{H}$, with the hope that $f(\mat{H})$ possesses decaying entries away from the sparsity pattern of $\mat{H}$.

Access to a good partitioning is a nontrivial requirement for probing methods.
However, when such an a partitioning is known, probing methods can be very effective \cite[Fig.~3]{FRS25}.

\chapter{Diagonal estimation} \label{ch:diagonal}

\epigraph{Extracting [the diagonal of the inverse of $\mat{A}$] is considered a challenging task, in part because the problem cannot be easily expressed in the form of a system of equations. The problem is usually harder to solve than a linear system with the
same matrix $\mat{A}$.
So far, not much literature has been devoted to this topic, whereas the related problem of estimating the trace of $\mat{B}$ [$=\mat{A}^{-1}$] received much more attention.}{Jok M. Tang and Yousef Saad, \textit{A probing method for computing the diagonal of a matrix inverse} \cite{TS12}}

\Cref{ch:loo} presented the basic leave-one-out\index{leave-one-out randomized algorithm} methodology and applied it to the problem of trace estimation.
This section will use the same approach to develop diagonal estimators, \XDiag and \XNysDiag, for general and psd matrices.
We will also encounter another approach to matrix attribute estimation based on unbiased low-rank approximation.

\myparagraph{Sources}
The \XDiag estimator was developed in the \XTrace paper:

\fullcite{ETW24}.

The \XNysDiag algorithm and the unbiased randomized SVD algorithms are newly developed in this thesis.

\myparagraph{Outline}
\Cref{sec:variance-reduced-diag} introduces existing approaches to diagonal estimation not using the leave-one-out\index{leave-one-out randomized algorithm} approach, \cref{sec:xdiag} presents the \XDiag estimator for the diagonal of a general matrix, \cref{sec:xnysdiag} develops the \XNysDiag estimator for the diagonal of a psd matrix, and \cref{sec:diag-experiments,sec:subgraph-1} provides an experimental comparison.
We conclude with an alternative approach to matrix attribute estimation based on \emph{unbiased low-rank approximation} in \cref{sec:unbiased-lra}.
The \XDiag estimator can also be derived using this approach.

\section{The BKS and \DiagPP diagonal estimators} \label{sec:variance-reduced-diag}

To begin this chapter, we briefly discuss existing approaches to the diagonal estimation problem and propose a few modifications.
Let $\mat{B} \in \field^{n\times n}$ be a general square matrix.

In \cref{sec:bks}, we reviewed the BKS diagonal estimator \cref{eq:bks} for $\diag(\mat{B})$.
Given iid isotropic vectors $\vec{\omega}_1,\ldots,\vec{\omega}_s$, we defined 
\begin{equation} \label{eq:updated-bks}
    \hat{\diag}_{\mathrm{BKS}} = \frac{1}{s} \sum_{i=1}^s (\mat{B} \vec{\omega}_i) \odot \overline{\vec{\omega}_i}.
\end{equation}
The BKS estimator is an unbiased Monte Carlo estimator for $\diag(\mat{B})$.
Given its similarity to the Girard--Hutchinson trace estimator, this estimator is also sometimes called \emph{Hutchinson's diagonal estimator}.
Analysis appears in \cite{BN22,DM23a}.

\begin{remark}[Original BKS estimator]
    The original work \cite{BKS07} of Bekas, Kokiopoulou, and Saad treated the real case $\vec{\omega}_i \in \real^n$ and considered the estimator
    \begin{equation} \label{eq:original-bks}
        \hat{\diag}_{\mathrm{BKS},\mathrm{orig}} = \frac{\sum_{i=1}^s (\mat{B} \vec{\omega}_i) \odot \vec{\omega}_i}{\sum_{i=1}^s \vec{\omega}_i\odot \vec{\omega}_i}.
    \end{equation}
    Division is performed entrywise.
    The original BKS estimator \cref{eq:original-bks} and simplified BKS estimator \cref{eq:updated-bks} coincide when the vectors $\vec{\omega}_i$ are drawn iid from the random sign distribution $\Unif \{\pm 1\}^n$.
    This thesis will work exclusively with the simplified BKS estimator $\hat{\diag}_{\mathrm{BKS}}$ given in \cref{eq:updated-bks}.
    %
\end{remark}

The \HutchPP variance reduction framework extends naturally to the diagonal estimation problem.
Given a budget of $s$ matvecs, begin by computing the randomized SVD approximation of rank $s/3$.
In detail, generate a test matrix $\mat{\Omega} \in \field^{n\times (s/3)}$ with isotropic random columns and form $\mat{Q} = \orth(\mat{B}\mat{\Omega})$.
As usual, $\mat{Q}$ defines a low-rank approximation $\Bhat = \mat{Q}(\mat{Q}^*\mat{B})$.
Now, we decompose $\diag(\mat{B})$ as
\begin{equation*}
    \diag(\mat{B}) = \diag(\Bhat) + \diag(\mat{B} - \Bhat),
\end{equation*}
and apply the BKS estimator to estimate the residual diagonal.
The resulting \emph{unbiased \DiagPP estimator} is
\begin{equation} \label{eq:udiagpp}
    \udiagppest \coloneqq \diag(\mat{Q}(\mat{Q}^*\mat{B})) + \frac{1}{s/3} \sum_{i=1}^{s/3} [(\Id - \mat{Q}\mat{Q}^*)\mat{B}\vec{\gamma}_i] \odot \overline{\vec{\gamma}_i}.
\end{equation}
The vectors $\vec{\gamma}_1,\ldots,\vec{\gamma}_{s/3}$ are chosen to be isotropic and iid.
The computational cost of the unbiased \DiagPP estimator is $2s/3$ matvecs with $\mat{B}$, $s/3$ matvecs with $\mat{B}^*$, and $\order(s^2n)$ additional arithmetic operations.
See \cref{prog:udiagpp} for code.

\myprogram{Unbiased \DiagPP estimator for the diagonal of a matrix.}{Subroutine \texttt{random\_signs} is provided in \cref{prog:random_signs}.}{udiagpp}

In the original paper on variance-reduced diagonal estimation \cite{BN22}, Baston and Nakatsukasa use the following alternative definition of the \DiagPP estimator:
\begin{equation} \label{eq:diagpp}
    \diagppest \coloneqq \diag(\mat{Q}\mat{Q}^*(\mat{B}\mat{Q})\mat{Q}^*) + \frac{1}{s/3} \sum_{i=1}^{s/3} [(\Id - \mat{Q}\mat{Q}^*)\mat{B}(\Id - \mat{Q}\mat{Q}^*)\vec{\omega}_i] \odot \overline{\vec{\omega}_i}.
\end{equation}
This estimator requires only matvecs with $\mat{B}$ and not its adjoint $\mat{B}^*$, an advantage over the unbiased \DiagPP estimator.
In general, the original \DiagPP estimator is biased, although this bias may be small relative to the diagonal if the matrix $\mat{B}$ has rapid singular value decay.
To distinguish between the two estimators \cref{eq:udiagpp,eq:diagpp}, ``the \DiagPP estimator'' will always refer to the original \DiagPP estimator \cref{eq:diagpp}.

\section{\XDiag: A leave-one-out diagonal estimator for general matrices} \label{sec:xdiag}

The idea of \XDiag is similar to \XTrace.
We proceed with a derivation.

\myparagraph{Step 1: Compute a low-rank approximation to the input matrix by multiplying it against a collection of test vectors}
We begin by computing a low-rank approximation to $\mat{B}$  using the randomized SVD.
Let $\mat{\Omega} \in \field^{n\times s/2}$ be a random matrix with isotropic columns, form the product $\mat{B}\mat{\Omega}$, and orthonormalize $\mat{Q} = \orth(\mat{B}\mat{\Omega})$.
We obtain the approximation
\begin{equation*}
    \Bhat \coloneqq \mat{Q}\mat{Q}^*\mat{B}.
\end{equation*}

\myparagraph{Step 2: Decompose the quantity of interest into known piece depending on the low-rank approximation plus a residual}
The diagonal operation is linear.
As such,
\begin{equation*}
    \diag(\mat{B}) = \diag(\Bhat) + \diag(\mat{B} - \Bhat) = \diag(\mat{Q}\mat{Q}^*\mat{B}) + \diag((\Id - \mat{Q}\mat{Q}^*)\mat{B}).
\end{equation*}

\myparagraph{Step 3: Construct a Monte Carlo estimate of the residual using a single random vector}
To approximate the residual diagonal $\diag((\Id - \mat{Q}\mat{Q}^*)\mat{B})$, we employ the single-vector BKS estimator:
\begin{equation*}
    \diag((\Id - \mat{Q}\mat{Q}^*)\mat{B}) \approx (\Id - \mat{Q}\mat{Q}^*)\mat{B}\vec{\omega} \odot \overline{\vec{\omega}}.
\end{equation*}
This step leads to the unbiased diagonal estimate 
\begin{equation*} 
    \hat{\diag} \coloneqq \diag(\mat{Q}\mat{Q}^*\mat{B}) + (\Id - \mat{Q}\mat{Q}^*)\mat{B}\vec{\omega} \odot \overline{\vec{\omega}}.
\end{equation*}

\myparagraph{Step 4: Downdate the low-rank approximation by recomputing it with a test vector removed, and use the left-out test vector as the random vector in step 3}
Now, rather than using a fresh vector $\vec{\omega}$, we downdate the randomized SVD, obtaining a new matrix $\mat{Q}_{(i)} = \orth(\mat{B}\mat{\Omega}_{-i})$ that is independent from the $i$th test vector $\vec{\omega}_i$.
This process defines a family of basic \XDiag estimators
\begin{equation*}
    \hat{\diag}_i \coloneqq \diag(\outprod{\mat{Q}_{(i)}}\mat{B}) + (\Id - \outprod{\mat{Q}_{(i)}})\mat{B}\vec{\omega}_i \odot \overline{\vec{\omega}_i} \quad \text{for } i=1,\ldots,s/2.
\end{equation*}

\myparagraph{Step 5: Average over all choices of vectors to leave out}
Finally, we average over the index $i$ to obtain the full \XDiag estimate
\begin{equation} \label{eq:xdiag}
    \hat{\diag}_{\mathrm{X}} \coloneqq \frac{1}{s/2} \sum_{i=1}^{s/2} \left[ \diag(\outprod{\mat{Q}_{(i)}}\mat{B}) + (\Id - \outprod{\mat{Q}_{(i)}})\mat{B}\vec{\omega}_i \odot \overline{\vec{\omega}_i}\right].
\end{equation}

\myparagraph{Efficient formula and implementation}
To implement \XDiag efficiently, we employ the randomized SVD downdating formula \cref{eq:rsvd-downdate}.
In particular, we form $\mat{Y} = \mat{B}\mat{\Omega}$, compute the \QR decomposition $\mat{Y} = \mat{Q}\mat{R}$, and obtain the downdating matrix $\mat{S}$ by normalizing the columns of $\mat{R}^{-*}$.
The basic \XDiag estimator is then
\begin{equation*}
    \hat{\diag}_i = \diag(\mat{Q}(\Id - \vec{s}_i^{\vphantom{*}}\vec{s}_i^*) \mat{Q}^*\mat{B}) + (\Id - \mat{Q}\mat{Q}^* + \mat{Q}\outprod{\vec{s}_i}\mat{Q}^*) \mat{B} \vec{\omega}_i^{\vphantom{*}} \odot \overline{\vec{\omega}}_i^{\vphantom{*}}.
\end{equation*}
We can simplify this equation is several ways.
First, define and compute 
\begin{equation*}
    \mat{Z} \coloneqq \mat{B}^*\mat{Q}.
\end{equation*}
Second, observe that $\Id - \mat{Q}\mat{Q}^*$ annihilates $\mat{B}\vec{\omega}_i = \vec{y}_i$ and $\mat{Q}^*\mat{B}\vec{\omega}_i = \mat{Q}^*\vec{y}_i = \vec{r}_i$.
Introducing $\mat{Z}$ and these observations, the basic \XDiag estimators simplify as
\begin{equation*}
    \hat{\diag}_i = \diagprod(\mat{Q}^*,\mat{Z}^*) + \mat{Q}\vec{s}_i^{\vphantom{*}} \odot (-\overline{\mat{Z}\vec{s}_i^{\vphantom{*}}} + \vec{s}_i^*\vec{r}_i^{\vphantom{*}} \cdot \overline{\vec{\omega}}_i^{\vphantom{*}}).
\end{equation*}
Averaging, we obtain a formula for the full \XDiag estimator.
\begin{equation*}
    \hat{\diag}_{\mathrm{X}} = \diagprod(\mat{Q}^*,\mat{Z}^*) + \frac{1}{s/2} [\mat{Q}\mat{S} \odot (-\overline{\mat{Z}\mat{S}} + \overline{\mat{\Omega}} \cdot \Diag(\diagprod(\mat{S},\mat{R})))]\onevec.
\end{equation*}
Code is provided in \cref{prog:xdiag}.

\myprogram{Efficient implementation of \XDiag.}{Subroutines \texttt{cnormc} and \texttt{diagprod} are provided in \cref{prog:cnormc,prog:diagprod}.}{xdiag}

\begin{remark}[Original \XDiag estimator]
    In our original paper \cite{ETW24}, my coauthors and I used the original BKS estimator \cref{eq:original-bks} for residual trace estimation in \XDiag.
    However, all of the experiments used the random sign distribution and were thus insensitive to the choice between the original and unbiased BKS estimators.
    I now prefer the version of \XDiag presented in this thesis, because it is unbiased for any choice of isotropic test vectors.
    The publicly released code for the paper \cite{ETW24}---available at \url{https://github.com/eepperly/XTrace}---implements the same version of \XDiag that appears in this thesis.
\end{remark}

\section{\XNysDiag: A leave-one-out diagonal estimator for psd matrices} \label{sec:xnysdiag}

We can develop an improved estimator \XNysDiag for the diagonal of a psd matrix by using Nystr\"om approximation in place of the randomized SVD.
The relationship between \XNysDiag and \XDiag is analogous to the relationship between \XNysTrace and \XTrace.

Having presented many detailed algorithm derivations in the leave-one-out framework, we now proceed more quickly.
Let $\mat{A} \in \field^{n\times n}$ be a psd matrix.
Begin by drawing a matrix $\mat{\Omega} \in \field^{n\times s}$ with iid isotropic columns, and define the Nystr\"om approximation
\begin{equation*}
	\Ahat \coloneqq \mat{A}\langle \mat{\Omega}\rangle.
\end{equation*}
Leaving out each column of $\mat{\Omega}$ in turn gives downdated approximations
\begin{equation*}
	\Ahat_{(i)} \coloneqq  \mat{A}\langle \mat{\Omega}_{-i}\rangle \quad \text{for } i=1,2,\ldots,s.
\end{equation*}
For each $i$, the diagonal may be decomposed as
\begin{equation*}
	\diag(\mat{A}) = \diag(\Ahat_{(i)}) + \diag(\mat{A} - \Ahat_{(i)}).
\end{equation*}
Estimating the second term using the single-vector BKS diagonal estimator yields the basic \XNysDiag estimators
\begin{equation*}
	\hat{\diag}_i \coloneqq  \diag(\Ahat_{(i)}) + (\mat{A} - \Ahat_{(i)})\vec{\omega}_i\odot\overline{ \vec{\omega}_i} \quad \text{for } i=1,2,\ldots,s.
\end{equation*}
Averaging over the index $i$ yields the full \XNysDiag estimate
\begin{equation*}
	\hat{\diag}_{\mathrm{XN}} = \frac{1}{s} \sum_{i=1}^s \hat{\diag}_i .
\end{equation*}

\myparagraph{Efficient formula and implementation}
To implement \XNysDiag efficiently, we compute the Nystr\"om approximation $\Ahat = \mat{F}\mat{F}^*$ using \cref{eq:stable-nystrom} and form the matrix $\mat{Z}$ in \cref{eq:nystrom-downdate-practice-2}.
The basic \XNysDiag estimators are
\begin{equation*}
    \hat{\diag}_i = \diag(\mat{F}\mat{F}^* - \vec{z}_i^{\vphantom{*}} \vec{z}_i^*) + (\mat{A} - \mat{F}\mat{F}^* + \vec{z}_i^{\vphantom{*}}\vec{z}_i^*)\vec{\omega}_i^{\vphantom{*}} \odot \overline{\vec{\omega}}_i^{\vphantom{*}}
\end{equation*}
In view of the interpolatory property \cref{prop:nystrom-properties}\ref{item:nystrom-interpolatory}, we have the identity $\mat{A}\vec{\omega}_i = (\mat{F}\mat{F}^*)\vec{\omega}_i$, so this formula simplifies to
\begin{equation*}
    \hat{\diag}_i = \srn(\mat{F}) - |\vec{z}_i^{\vphantom{*}}|^2 + (\vec{z}_i^{\vphantom{*}} \odot \overline{\vec{\omega}}_i^{\vphantom{*}}) \cdot \vec{z}_i^*\vec{\omega}_i^{\vphantom{*}}.
\end{equation*}
Averaging over the index $i$ yields a formula for the full \XNysDiag estimator:
\begin{equation*}
    \hat{\diag}_{\mathrm{XN}} = \srn(\mat{F}) - \frac{1}{s} \srn(\mat{Z}) + \frac{1}{s} (\mat{Z} \odot \overline{\mat{\Omega}}) \cdot \Diag(\diagprod(\mat{Z},\mat{\Omega})) \cdot \onevec.
\end{equation*}

\myprogram{Efficient and stable implementation of \XNysDiag estimator.}{Subroutine \texttt{nystrom} is provided in \cref{prog:nystrom}.}{xnysdiag}

For stability, we apply \XNysDiag to $\mat{A} + \mu \Id$ with the shift from \cref{eq:nys-mu}. 
We correct the shift by subtracting $\mu \onevec$ from the diagonal estimate.
Code is provided in \cref{prog:xnysdiag}.

\begin{remark}[\XNysDiag and \XNysTrace]
	Observe that the \XNysTrace estimator is the sum of the entries of the \XNysDiag estimator.
	In contrast, the \XTrace estimate does not equal to the sum of the entries of the \XDiag estimate.
\end{remark}

\section{Synthetic Experiments} \label{sec:diag-experiments}

We compare diagonal estimators on two sets of synthetic examples.
First, we test with the four synthetic test matrices from \cref{eq:test-matrices} with varying spectral properties.
Second, we test on two matrices with similar spectra but very disparate distributions of diagonal entries.

\subsection{Tests on matrices with different spectra}

\begin{figure}
    \centering
    \includegraphics[width=0.99\linewidth]{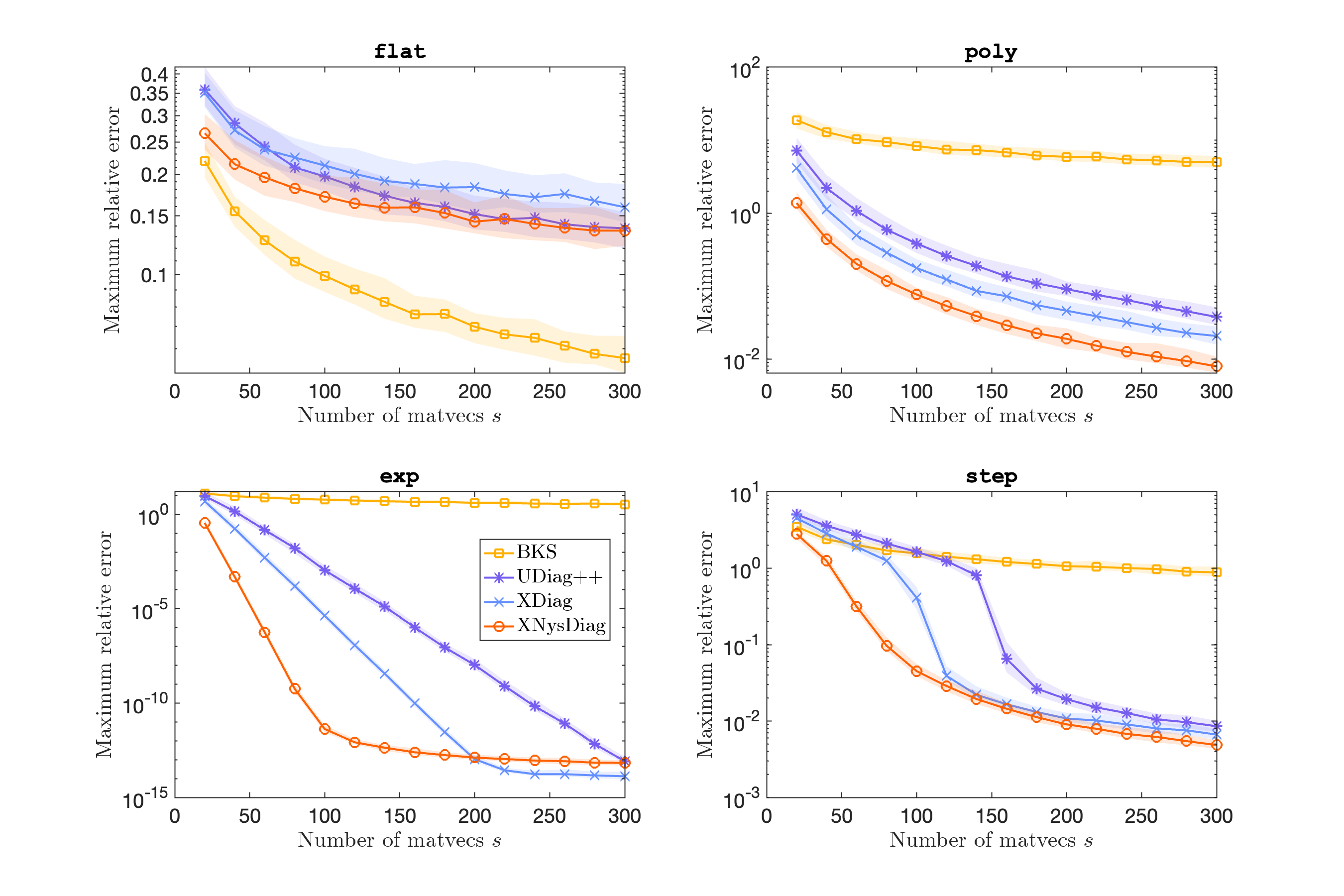}
    \caption[Comparison of BKS, unbiased \DiagPP, \XDiag, and \XNysDiag for estimating the diagonal of matrices with different spectra]{Maximum relative error of diagonal estimates by the BKS estimator (yellow squares), unbiased \DiagPP (purple asterisks), \XDiag (blue crosses), and \XNysDiag (orange circles) for the four test matrices \cref{eq:test-matrices} as a function of the numbers of matvecs $s$.
    The error metric is the maximum relative error, defined in \cref{eq:max-rel-err-diag}.
    Lines show the median of 100 trials, and shaded regions range from the 10\% to the 90\% quantile.}
    \label{fig:diagonal-comparison}
\end{figure}

\Cref{fig:diagonal-comparison} provides a comparison of the BKS, unbiased \DiagPP, \XDiag, and \XNysDiag algorithms on the synthetic matrices from \cref{eq:test-matrices}.
To quantify the discrepancy between a diagonal approximation $\hat{\diag}$ and the true diagonal $\diag$, we use the maximum relative error:
\begin{equation} \label{eq:max-rel-err-diag}
    \mathrm{maxrelerr} = \max_{1\le i \le n} \frac{|\mathrm{\hat{diag}}_i - \mathrm{diag}_i|}{|\mathrm{diag}_i|}.
\end{equation}
The comparison of diagonal estimators in \cref{fig:diagonal-comparison} is similar with the comparison of trace estimators in \cref{fig:trace-comparison}: The \XNysDiag and \XDiag estimators outperform unbiased \XDiag and the BKS estimator on problems with spectral decay (i.e., the \texttt{poly} and \texttt{exp} matrices), and the BKS estimator slightly outperforms the other estimators when the spectrum is flat (i.e., the \texttt{flat} example).
The \texttt{step} example demonstrates a dissimilarity between the diagonal estimation and trace estimation problems: For diagonal estimation, the \XNysDiag estimator outperforms \XDiag on this example for any number of matvecs $s$, whereas \XNysTrace does worse than \XTrace on this example for sufficiently large values of $s$.

\subsection{Tests on matrices with similar spectra but different diagonals}

A key distinction between diagonal estimation and trace estimation is that the diagonal entries may vary across many orders of magnitude.
In such cases, it is harder to estimate all of the diagonal entries up to a small relative error.
To demonstrate this point, we compare the diagonal estimators on two matrices with similar spectra but very different diagonals: the \texttt{poly} matrix \cref{eq:poly} defined above and the following \texttt{scaledWishart} matrix
\begin{equation} \label{eq:scaled-Wishart}
    \texttt{scaledWishart} \coloneqq \mat{D}^{1/2} \mat{G}^*\mat{G} \mat{D}^{-1/2}
\end{equation}
where $\mat{G} \in \real^{2n\times n}$ is a matrix with iid (real) standard Gaussian entries and $\mat{D} \coloneqq \Diag(i^{-2} : i=1,\ldots,n)$ is a diagonal matrix with entries decaying at an algebraic rate.
We again set $n\coloneqq 10^3$.
These two matrices have the following properties:
\begin{itemize}
    \item \textbf{\textit{Similar eigenvalue decay.}} Both \texttt{poly} and \texttt{scaledWishart} have polynomially decreasing eigenvalues spanning a similar range: $\lambda_{\mathrm{max}}/\lambda_{\mathrm{min}}$ is 1.0e6 for \texttt{poly} and 5.0e6 for \texttt{scaledWishart}.
    \item \textbf{\textit{Dustubct diagonal values.}} Despite their spectral similarities, \texttt{poly} and \texttt{scaledWishart} have very different diagonals.
    The diagonal entries of \texttt{poly} are comparable, with $\max a_{ii} / \min a_{ii}=\text{6.9e1}$, while the diagonal entries of \texttt{scaledWishart} span many orders of magnitude, with $\max a_{ii} / \min a_{ii}=\text{1.0e6}$.
\end{itemize}

\begin{figure}
    \centering
    \includegraphics[width=0.99\linewidth]{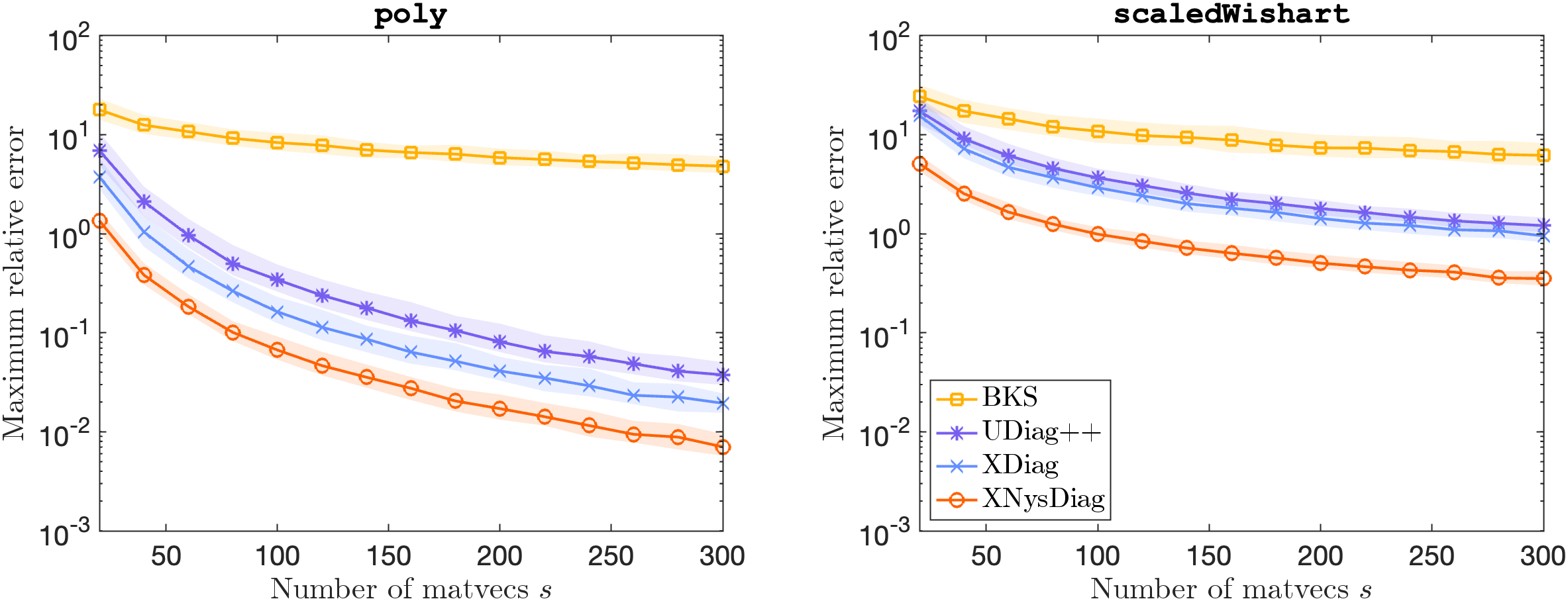}
    \caption[Comparison of BKS, unbiased \DiagPP, \XDiag, and \XNysDiag for matrices with different diagonal distributions]{Maximum relative error of diagonal estimates by the BKS estimator (yellow squares), unbiased \DiagPP (purple asterisks), \XDiag (blue crosses), and \XNysDiag (orange circles) on the \texttt{poly} matrix (\cref{eq:poly}, \emph{left}) and \texttt{scaledWishart} matrix (\cref{eq:scaled-Wishart}, \emph{right}).
    The error metric is the maximum relative error, defined in \cref{eq:max-rel-err-diag}.
    Lines show the median of 100 trials, and shaded regions show the 10\% and 90\% quantiles.}
    \label{fig:diagonal-multiple-scale}
\end{figure}

Results appear in \cref{fig:diagonal-multiple-scale}.
We see that while the \emph{ordinal ranking} of the algorithms remains unchanged (\XNysDiag > \XDiag > unbiased \DiagPP > BKS), the maximum relative error is significantly higher for the \texttt{scaledWishart} example than the \texttt{poly} examples.
These experiments demonstrate that accurately approximating the small diagonal entries of a matrix whose diagonal entries span multiple orders of magnitude can be challenging for variance-reduced diagonal estimators.
This behavior makes sense: Unbiased \DiagPP, \XDiag, and \XNysDiag all work by computing a low-rank approximation of the input matrix, which naturally capture large matrix entries better than smaller ones.

\section{Application: Subgraph centralities} \label{sec:subgraph-1}

In \cref{sec:estrada}, we used stochastic trace estimators to compute the Estrada index of a graph with adjacency matrix $\mat{M}$, which is defined as
\begin{equation*}
    \mathrm{estr} \coloneqq \tr(\exp(\mat{M})).
\end{equation*}
The Estrada index serves as a \emph{global} measure of how ``centralized'' or ``clustered'' the nodes of a graph are.
The \emph{subgraph centralities} \cite{Est22} provide an analogous notion of how ``central'' each node is to a graph.
The vector $\mathbf{sc}$ of subgraph centralities is the diagonal of $\mat{A} = \exp(\mat{M})$:
\begin{equation*}
    \mathbf{sc} \coloneqq \diag(\exp(\mat{M})).
\end{equation*}
Consequently, the sum of the subgraph centralities is the Estrada index.

\begin{figure}
    \centering
    \includegraphics[width=0.99\linewidth]{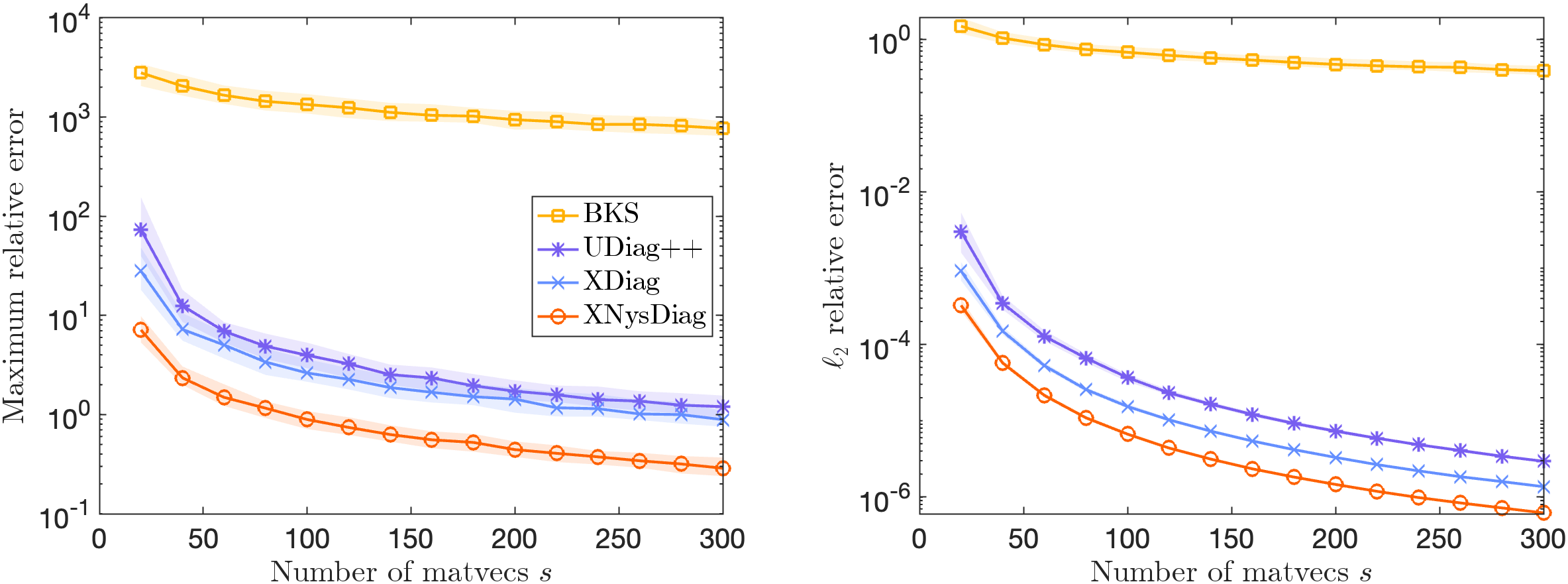}
    \caption[Comparison of BKS, unbiased \DiagPP, \XDiag, and \XNysDiag for estimating subgraph centralities]{Maximum relative error of diagonal estimates by the BKS estimator (yellow squares), unbiased \DiagPP (purple asterisks), \XDiag (blue crosses), and \XNysDiag (orange circles) for estimating the subgraph centralities.
    We show both the the maximum relative error (\emph{left}, \cref{eq:max-rel-err-diag}) and the relative $\ell_2$ error (\emph{right}).
    Lines show median of 100 trials, and error bars show 10\% and 90\% quantiles.}
    \label{fig:subgraph-1}
\end{figure}

\Cref{fig:subgraph-1} shows estimates of the the subgraph centralities calculated with different diagonal estimation algorithms.
As in \cref{fig:estrada}, we use the \texttt{yeast} dataset \cite{BZC+03} and perform matvecs with $\mat{A} = \exp(\mat{M})$ by forty steps of the Lanczos algorithm \cite[Ch.~6]{Che24}.
The left panel shows the maximum relative error \cref{eq:max-rel-err-diag} and the right panel shows the $\ell_2$ error $\norm{\smash{\hat{\diag} - \diag}}/\norm{\diag}$.
We see that \XNysDiag clearly outperforms \XDiag and unbiased \DiagPP, which in turn outperform the BKS estimator.
In absolute terms, the maximum relative error for all methods is fairly large, with the \XNysDiag estimator achieving a maximum relative error of 0.3 with $s = 300$ matvecs.
The $\ell_2$ errors are much smaller for all methods.
The reason for the large discrepancy between the maximum relative error and $\ell_2$ error is that the diagonal entries of $\mat{A}$ span many orders of magnitude, similar to the scaled Wishart example \cref{eq:scaled-Wishart}. 
Thus, the variance-reduced methods compute the largest diagonal entries to high-accuracy but struggle to accurately approximate the smaller entries.
We will explore a partial solution to this ``small diagonal entries'' problem in \cref{sec:square-root-trick,sec:subgraph-2}.

\section{Another view: Unbiased low-rank approximation} \label{sec:unbiased-lra} 

In this section, we will develop another way of deriving the \XDiag estimator.
This strategy extends naturally to construct estimators for other matrix entries or linear functionals.
So far, we have approached the problem of matrix attribute estimation by designing purpose-built unbiased estimators that target each particular matrix attribute.
But one could instead form an unbiased estimate $\Bhat$ of the input matrix $\mat{B} \in \field^{m\times n}$ itself.
Given such a $\Bhat$, one can obtain unbiased estimates for any \warn{linear} matrix attribute $L(\mat{B})$ by computing $L(\Bhat)$.
Here, $L : \field^{m\times n} \to \set{V}$ denotes an arbitrary linear map.
Trace estimation and diagonal estimation are special cases of this framework with $L = \tr$ and $L = \diag$.
Other special cases include $L(\mat{B}) = b_{ij}$ (computing a specified entry) or $L(\mat{B}) = \tr(\mat{C}^*\mat{B})$ (computing the inner product with another matrix $\mat{C}$).

We can use the leave-one-out\index{leave-one-out randomized algorithm} approach to develop unbiased low-rank approximation primitives.
We will illustrate by developing an unbiased version of the randomized SVD.
Consider the task of computing an unbiased rank-$k$ approximation to a general matrix $\mat{B} \in \field^{m\times n}$.
We follow the leave-one-out\index{leave-one-out randomized algorithm} approach:

\myparagraph{Step 1: Compute a low-rank approximation to the input matrix by multiplying it against a collection of test vectors}
Begin by generating a matrix $\mat{\Omega} \in \field^{n\times k}$ with iid isotropic columns, and form the randomized SVD approximation: $\mat{Q} = \orth(\mat{B}\mat{\Omega})$ and $\mat{Q}\mat{Q}^*\mat{B}\approx \mat{B}$.

\myparagraph{Step 2: Decompose the quantity of interest into known piece depending on the low-rank approximation plus a residual}
The matrix $\mat{B}$ can decomposed as
\begin{equation*}
    \mat{B} = \mat{Q}\mat{Q}^*\mat{B} + (\Id - \mat{Q}\mat{Q}^*)\mat{B}.
\end{equation*}

\myparagraph{Step 3: Construct a Monte Carlo estimate of the residual using a single random vector}
To approximate the residual $(\Id - \mat{Q}\mat{Q}^*)\mat{B}$, introduce a copy of the identity matrix and introduce the unbiased stochastic approximation $\Id \approx \vec{\omega}\vec{\omega}^*$ furnished by an isotropic vector $\vec{\omega} \in \field^n$:
\begin{equation*}
    (\Id - \mat{Q}\mat{Q}^*)\mat{B} = (\Id - \mat{Q}\mat{Q}^*)\mat{B} \cdot \Id \approx (\Id - \mat{Q}\mat{Q}^*)\mat{B} \cdot \vec{\omega}\vec{\omega}^*.
\end{equation*}
This results in a stochastic approximation to $\mat{B}$,
\begin{equation*}
    \mat{B} \approx \mat{Q}\mat{Q}^*\mat{B} + (\Id - \mat{Q}\mat{Q}^*)\mat{B}\vec{\omega}\vec{\omega}^*.
\end{equation*}

\myparagraph{Step 4: Downdate the low-rank approximation by recomputing it with a test vector removed, and use the left-out test vector as the random vector in step 3}
For each $i$, downdate the randomized SVD approximation to produce a downdated low-rank approximation $\mat{Q}_{(i)}^{\vphantom{*}}\mat{Q}_{(i)}^*\mat{B}$ that is independent of $\vec{\omega}_i$. 
For each $i$, this procedure yields an unbiased low-rank approximation
\begin{equation} \label{eq:usvd-i}
    \Bhat_i \coloneqq \mat{Q}_{(i)}^{\vphantom{*}}\mat{Q}_{(i)}^*\mat{B} + (\Id - \mat{Q}_{(i)}^{\vphantom{*}}\mat{Q}_{(i)}^*)\mat{B}\vec{\omega}_i^{\vphantom{*}}\vec{\omega}_i^*.
\end{equation}

\myparagraph{Step 5: Average over all choices of vectors to leave out}
Finally, average \cref{eq:usvd-i} over the index $i$ to obtain the unbiased randomized SVD low-rank approximation
\begin{equation} \label{eq:usvd}
    \Bhat \coloneqq \frac{1}{k} \sum_{i=1}^k \left[ \mat{Q}_{(i)}^{\vphantom{*}}\mat{Q}_{(i)}^*\mat{B} + (\Id - \mat{Q}_{(i)}^{\vphantom{*}}\mat{Q}_{(i)}^*)\mat{B}\vec{\omega}_i^{\vphantom{*}}\vec{\omega}_i^*\right].
\end{equation}

\myparagraph{Implementation}
To develop an efficient unbiased randomized SVD implementation, we substitute the downdating formula \cref{eq:rsvd-downdate} in \cref{eq:usvd} and simplify algebraically, \`a la \cref{sec:xtrace-implementation}.
As usual, the derivation is somewhat complicated; we omit the details.
Here is the procedure:
First, compute $\mat{Y} \coloneqq \mat{B}\mat{\Omega}$, the (economy-size) \QR decomposition $\mat{Y} = \mat{Q}\mat{R}$, the downdating matrix $\mat{S}$ in \cref{eq:rsvd-downdate}, and the product $\mat{C} \coloneqq \mat{B}^*\mat{Q}$.
The basic unbiased randomized SVD estimators \cref{eq:usvd-i} are
\begin{equation} \label{eq:usvd-i-formula}
    \Bhat_i = \mat{Q}\mat{C}^* + \mat{Q}\vec{s}_i^{\vphantom{*}} (-\vec{s}_i^*\mat{C}^* + \vec{s}_i^*\vec{r}_i^{\vphantom{*}} \cdot \vec{\omega}_i^*),
\end{equation}
and their mean is
\begin{equation*}
    \Bhat = \mat{Q}\left[ \left( \Id - \frac{\mat{S}\mat{S}^*}{k}\right)\mat{C}^* + \frac{1}{k} \cdot \mat{S} \cdot \Diag(\diagprod(\mat{S},\mat{R})) \cdot \mat{\Omega}^* \right].
\end{equation*}
An implementation appears in \cref{prog:usvd}.

\myprogram{Unbiased randomized SVD algorithm for unbiased low-rank approximation and matrix attribute estimation.}{Subroutines \texttt{diagprod}, \texttt{random\_signs}, and \texttt{cnormc} are provided in \cref{prog:diagprod,prog:random_signs,prog:cnormc}.}{usvd}

\myparagraph{Deriving the \XDiag estimator from the unbiased randomized SVD}
Comparing the \XDiag estimator \cref{eq:xdiag} to the unbiased randomized SVD approximation \cref{eq:usvd}, we see the \XDiag estimator is the diagonal of the unbiased randomized SVD approximation.
In this way, the unbiased randomized SVD algorithm (\cref{prog:usvd}) gives another way of implementing \XDiag: Simply run the unbiased randomized SVD and form \XDiag estimator by the formula $\hat{\diag} = \diagprod(\mat{U}^*,\mat{\Sigma}\mat{V}^*)$.
This approach has the advantage that one can reuse the same unbiased randomized SVD approximation to estimate other entries of the matrix on an as-needed basis via the formula
\begin{equation*}
    b_{ij} \approx \hat{b}_{ij} = \mat{U}(i,:)\mat{\Sigma}\mat{V}(j,:)^*.
\end{equation*}

\chapter{Row-norm estimation} \label{ch:row-norm}

\epigraph{I often say that when you can measure what you are speaking about, and express it in numbers, you know something about it; but when you cannot measure it, when you cannot express it in numbers, your knowledge is of a meagre and unsatisfactory kind; it may be the beginning of knowledge, but you have scarcely, in your thoughts, advanced to the stage of science, whatever the matter may be.}{Lord Kelvin \cite{Kel89}}

Among matrix attribute estimation problems, the trace estimation and diagonal estimation problems have received the lion's share of the attention.
But there is another problem of fundamental interest that also merits study: the problem of estimating the \emph{row norms} of a matrix $\mat{B}$.

This chapter uses the leave-one-out\index{leave-one-out randomized algorithm} approach to develop algorithms for estimating the row norms of an implicit matrix.
The resulting \XRowNorm and \XSymRowNorm estimators are applicable to general and Hermitian matrices, respectively.

\myparagraph{Sources}
This chapter is based on original research that has not yet been published.
It uses the leave-one-out\index{leave-one-out randomized algorithm} approach developed in the \XTrace paper \cite{ETW24}.

\myparagraph{Outline}
The squared row norms of $\mat{B}$ comprise the diagonal of $\mat{B}\mat{B}^*$, so it is natural to ask: Should we just use diagonal estimators for row-norm estimation?
We answer this question with a definitive ``no'' in \cref{sec:row-norm-diagonal}.
In fact, as we discuss in \cref{sec:square-root-trick}, sometimes one can (and should) use row-norm estimators to solve diagonal estimation problems!
With the importance of purpose-built row-norm estimators established, \cref{sec:variance-reduced-row-norm} describes existing variance-reduced but non-leave-on-out row-norm estimators, \cref{sec:xrownorm} develops the \XRowNorm algorithm for general matrices, and \cref{sec:xsymrownorm} introduces the \XSymRowNorm algorithm for Hermitian matrices.
\Cref{sec:row-norm-experiments,sec:subgraph-2} provide experimental comparisons.

\section{Is row-norm estimation just diagonal estimation?} \label{sec:row-norm-diagonal}

Viewed one way, the row-norm estimation problem is a special case of the diagonal estimation problem, owing to the identity $\srn(\mat{B}) = \diag(\mat{B}\mat{B}^*)$. 
Therefore, the row-norm estimation problem can be solved by applying a diagonal estimator, such as \XNysDiag, to $\mat{B}\mat{B}^*$.

However, the most accurate algorithms for the row-norm estimation problem do not proceed via diagonal estimation.
To see why, let us compare two estimates of $\srn(\mat{B}) = \diag(\mat{B}\mat{B}^*)$:
The BKS estimator 
\begin{equation} \label{eq:bks-row-norm}
    \hat{\srn}_{\mathrm{BKS}} = \frac{1}{k} \sum_{i=1}^k \mat{B}(\mat{B}^*\vec{\omega}_i) \odot \overline{\vec{\omega}_i}
\end{equation}
versus the Johnson--Lindenstrauss row-norm estimator
\begin{equation} \label{eq:mc-row-norm-compare}
    \hat{\srn}_{\mathrm{JL}} = \frac{1}{k} \sum_{i=1}^k |\mat{B}\vec{\omega}_i|^2 = \frac{1}{k} \sum_{i=1}^k (\mat{B}\vec{\omega}_i) \odot \overline{(\mat{B}\vec{\omega}_i)}.
\end{equation}
The first difference between these estimators is their computational cost; to form a $k$-term estimate, the BKS estimator $\hat{\srn}_{\mathrm{BKS}}$ requires $k$ matvecs with $\mat{B}$ and $k$ matvecs with $\mat{B}^*$.
The Johnson--Lindenstrauss estimate $\hat{\srn}_{\mathrm{JL}}$ requires just $k$ matvecs with $\mat{B}$.
Perhaps a more significant difference emerges from the error analysis of these two estimators.

\begin{theorem}[Comparison of Monte Carlo row norm estimates] \label{thm:row-norm-comparison}
    Let $\mat{B} \in\real^{m\times n}$ be a \warn{real} matrix with rows $\vec{\beta}_i^*$, and consider the squared row norm estimates \cref{eq:bks-row-norm,eq:mc-row-norm-compare} using iid vectors on the sphere $\sqrt{n}\sphere(\field^n)$.
    Then, for each $1\le i\le m$, we have the variance bounds
    \begin{align}
        \expect\left|\hat{\srn}_{\mathrm{JL}}(i) - \norm{\vec{\beta}_i}^2 \right|^2 &= \frac{1}{k} \cdot \frac{n-1}{n} K_\field \cdot \norm{\vec{\beta}_i}^4, \label{eq:mc-row-norm-error} \\
        \expect\left|\hat{\srn}_{\mathrm{BKS}}(i) - \norm{\vec{\beta}_i}^2 \right|^2 &\le \frac{1}{k} \cdot K_\field \cdot\sum_{j=1}^m |\vec{\beta}_i^*\vec{\beta}_j^{\vphantom{*}}|^2. \label{eq:bks-row-norm-error} 
    \end{align}
    The constants are $K_\real = 2n/(n+2)$ and $K_\complex = n/(n+1)$.
\end{theorem}

This result suggests that the Johnson--Lindenstrauss row-norm estimator \cref{eq:mc-row-norm-compare} is superior to the BKS diagonal estimator applied to $\mat{B}\mat{B}^*$ \cref{eq:bks-row-norm}.
The error of the Johnson--Lindenstrauss estimator of the $i$th row norm depends only on the $i$th row norm of $\mat{B}$; by contrast, the bound on the error of the BKS estimator depends on the magnitude of the inner products of the $i$th row of $\mat{B}$ against all of the other rows.
For this example $\mat{B} = \onevec\vec{\beta}^*$, mean-squared error for the Monte Carlo row-norm estimator \cref{eq:mc-row-norm-error} is roughly a factor of $m$ smaller than the bound \cref{eq:bks-row-norm-error} on the mean-squared error for the BKS estimator.
On this extreme case, the BKS estimator \cref{eq:bks-row-norm} is essentially worthless as it expends $\Omega(m)$ matvecs to achieve a nontrivial result. (Recall that any matrix attribute estimation problem can be solved at a trivial cost of $\min(m,n)$ matvecs.)
The Johnson--Lindenstrauss estimator \cref{eq:mc-row-norm-compare}, by contrast, computes any row norm of any matrix $\mat{B}$ to root-mean-squared relative error $\varepsilon$ in just $\order(1/\varepsilon^2)$ matvecs, independent of the dimensions of $\mat{B}$.

\Cref{thm:row-norm-comparison} demonstrates that dedicated algorithms for the row-norm estimation problem can significantly outperform diagonal estimators applied to the outer product matrix $\mat{B}\mat{B}^*$.
This result justifies our study of the row-norm estimation problem as distinct from diagonal estimation, and motivates the development of purpose-build row-norm estimators.

\begin{proof}[Proof of \cref{thm:row-norm-comparison}]
    Both the BKS and Johnson--Lindenstrauss row-norm estimators \cref{eq:bks-row-norm,eq:mc-row-norm-compare} are averages of $k$ iid copies of a single-vector unbiased estimator; ergo, the mean-squared error for each is $1/k$ times the variance of the single-vector estimator.

    Begin with the \RowNormName estimator
    \begin{equation*}
        \hat{\srn}_{\mathrm{JL}} = |\mat{B}\vec{\omega}|^2,
    \end{equation*}
    and restrict attention to the $i$th coordinate $\hat{\srn}_{\mathrm{MC}}(i) = |\vec{\beta}_i^*\vec{\omega}|^2$.
    We may rewrite this coordinate in the form of a Girard--Hutchinson trace estimator $\hat{\srn}_{\mathrm{MC}}(i) = \vec{\omega}^*(\vec{\beta}_i^{\vphantom{*}}\vec{\beta}_i^*)\vec{\omega}$.
    The matrix has simple eigenvalue $\norm{\vec{\beta}_i}^2$ and a multiple eigenvalue $0$ with multiplicity $n-1$.
    The average eigenvalue is $\overline{\lambda} = \norm{\vec{\beta}_i}^2/n$.
    Therefore, by \cref{fact:gh-variance}, we obtain
    \begin{equation*}
        \Var(\hat{\srn}_{\mathrm{MC}}(i)) = K_\field \left[\left(\norm{\vec{\beta}_i}^2 - \frac{\norm{\vec{\beta}_i}^2}{n}\right)^2 + (n-1)\left(\frac{\norm{\vec{\beta}_i}^2}{n}\right)^2\right] = \frac{n-1}{n} K_\field \cdot \norm{\vec{\beta}_i}^4.
    \end{equation*}

    Now, we treat the single-vector BKS estimator
    \begin{equation*}
        \hat{\srn}_{\mathrm{BKS}} = \mat{B} \mat{B}^*\vec{\omega} \odot \overline{\vec{\omega}}
    \end{equation*}
    and restrict attention to the $i$th coordinate $\hat{\srn}_{\mathrm{BKS}}(i) = \vec{\beta}_i^* \mat{B}^*\vec{\omega}\cdot \overline{\omega}_i$.
    Observe that we can write $\overline{\omega_i} = \vec{\omega}^*\evec_1$ and $\vec{\beta}_i^* = \evec_i^*\mat{B}$.
    Therefore,
    \begin{equation*}
        \hat{\srn}_{\mathrm{BKS}}(i) = \overline{\omega}_1\mat{B}(i,:) \mat{B}^*\vec{\omega} = \vec{\omega}^*(\evec_i^{\vphantom{*}}\evec_i^*\mat{B}\mat{B}^*)\vec{\omega}.
    \end{equation*}
    Invoking \cref{fact:gh-variance}, we obtain
    \begin{equation*}
        \Var(\hat{\srn}_{\mathrm{BKS}}(i)) \le K_\field \norm{\evec_i^{\vphantom{*}}\evec_i^*\mat{B}\mat{B}^*}_{\mathrm{F}}^2 = K_\field \norm{\vec{\beta}_i^*\mat{B}^*}^2 = K_\field \sum_{j=1}^m |\vec{\beta}_i^*\vec{\beta}_j^{\vphantom{*}}|^2.
    \end{equation*}
    This completes the proof of the stated results.
\end{proof}

\begin{remark}[History]
    Standard probabilistic proofs of the Johnson--Lindenstrauss lemma \cite{JL84} show that $\hat{\srn}_{\mathrm{JL}}$ approximates each squared row norm of $\mat{B}$ up to a constant factor, provided the number of matvecs is $k = \order(\log n)$.
    As such, it is natural to attribute the estimator $\hat{\srn}_{\mathrm{JL}}$ to Johnson and Lindenstrauss.
    However, despite being an immediate corollary of the Johnson–Lindenstrauss lemma (1984), I was unable to find a reference where $\hat{\srn}_{\mathrm{JL}}$ was used to estimate the row norms of an implicit matrix earlier than Spielman and Srivastava's work in 2008 \cite{SS08}.
    The papers \cite{DMMW12,LMP13} also use this technique.
    The value of purpose-built row-norm estimators over diagonal estimators is emphasized by Mathur, Moka, and Botev \cite{MMB21} and in recent works of Michael Lindsey \cite{Lin23a,FL24}.
\end{remark}

\subsection{Numerical demonstration and application: Leverage-score estimation}

As an application and a demonstration of \cref{thm:row-norm-comparison}, we compare the Johnson--Lindenstrauss row-norm estimator and the BKS and \XNysDiag diagonal estimators for the problem of \emph{leverage-score estimation}.
Recall that the leverage scores of a matrix $\mat{B}$ are the squared row norms of $\mat{Q} \coloneqq \orth(\mat{B})$.
Estimating leverage scores is major application for row-norm estimators.
Following \cite{CP19,MMM+23}, we focus on leverage-score estimation for polynomial regression.

\begin{figure}
    \centering
    \includegraphics[width=0.49\linewidth]{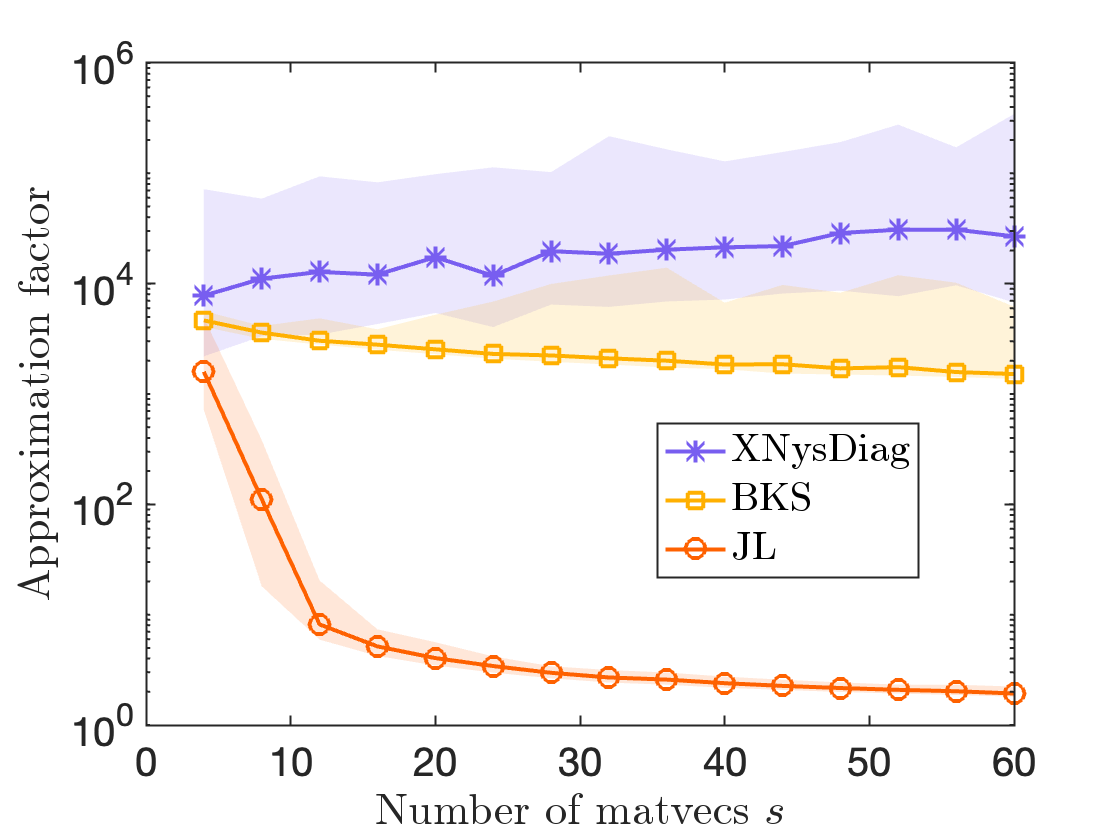}
    \includegraphics[width=0.49\linewidth]{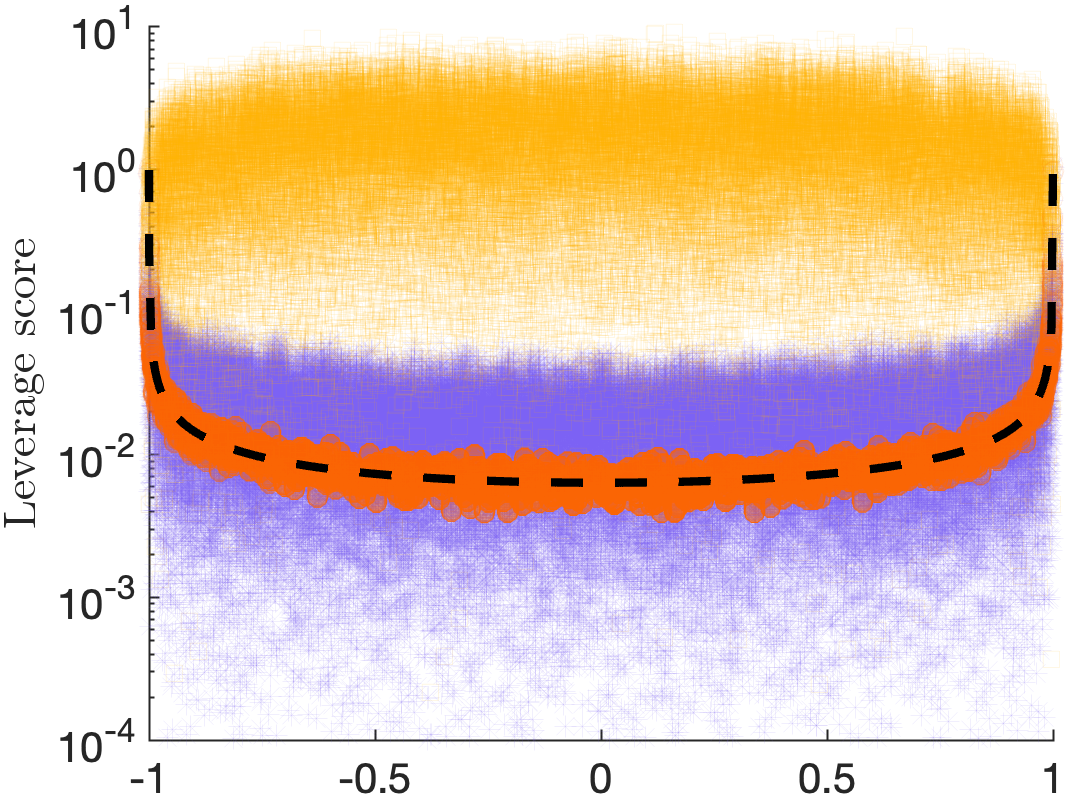}
    \caption[Comparison of Johnson--Lindenstrauss row-norm estimator to BKS and \XNysDiag diagonal estimators for leverage-score estimation]{\emph{Left:} Approximation factor for Johnson--Lindenstrauss row-norm estimator (orange circles) and BKS (yellow squares) and \XNysDiag (purple asterisks) diagonal estimators for estimating the leverage scores of the polynomial regression matrix.
    Lines trace the median of 100 trials, and shaded regions show 10\% and 90\% quantiles.
    \emph{Right:} Scatter plot of \warn{positive} leverage-score estimates by these methods (single trial with $s=60$). 
    The black dashed line indicates the true leverage scores.}
    \label{fig:lev-estimate}
\end{figure}

Let $x_1,\ldots,x_m$ denote equally spaced points on $[-1,1]$, and define the polynomial regression matrix $\mat{B} \in \real^{m\times n}$ with entries $b_{ij} = T_{j-1}(x_i)$; here, $T_k$ denotes the $k$th Chebyshev polynomial (of the first kind).
We set $m\coloneqq 10^5$ and $n\coloneqq 10^3$.
\Cref{fig:lev-estimate} shows the results of applying the Johnson--Lindenstrauss row-norm estimator to $\mat{Q}$ and the BKS and \XNysDiag diagonal estimators to $\mat{Q}\mat{Q}^*$, where $\mat{Q} = \Orth(\mat{B}).$
For each method, we report the approximation ratio between the true leverage scores $\vec{\ell} \coloneqq \srn(\mat{Q})$ and the approximate leverage scores $\vec{\hat{\ell}}$, calculated as
\begin{equation*}
    \text{approximation ratio} \coloneqq \max \left\{ \max \left\{ \frac{\hat{\ell}_i}{\ell_i}, \frac{\ell_i}{\hat{\ell}_i}\right\} : i=1,\ldots,n \right\}
\end{equation*}
As seen in the left panel, the approximation ratio for the Johnson--Lindenstrauss row-norm estimator is dramatically lower than for the diagonal estimators.
The diagonal estimators also produce many negative leverage-score estimates, which are clearly vacuous.
The right panel shows the (positive) leverage-score estimates as a function of the position $x$.
The Johnson--Lindenstrauss estimates closely track the true leverage scores, whereas the other estimates carry essentially no information.

\begin{remark}[Estimating leverage scores: Better algorithms]
    Given a matrix $\mat{B}\in\field^{m\times n}$, calculating $\mat{Q} = \Orth(\mat{B})$ and applying a row-norm estimate is a poor algorithm.
    Forming $\mat{Q}$ expends an expensive $\order(mn^2)$ operations and, if one expends the effort to compute $\mat{Q}$, the leverage scores are cheap to compute exactly.
    To obtain leverage-score estimates quickly, one should instead apply stochastic row-norm estimates to a preconditioned version of the matrix $\mat{B}$.
    This fast approach requires roughly $\order(mn + n^3)$ operations; see \cite{DMMW12} for details.
    Additionally, the case of polynomial regression, continuous analogs of the leverage scores are known in closed form.
    For theser easons, \cref{fig:lev-estimate} serves as comparison point for different algorithms, not as a serious application for row-norm estimation.
\end{remark}

\section{The square-root trick: Diagonal estimation via row-norm estimation} \label{sec:square-root-trick}

The last section provided strong evidence that one should not, as a matter of practice, use diagonal estimators to solve row-norm estimation problems.
But what about the other way around?
Can it ever be beneficial to use row-norm estimators to solve diagonal estimation problems?

The following model describes one scenario where the answer to this question is yes.
\actionbox{\textbf{Square-root matvec model.} The matrix $\mat{A}$ is psd, and we have \warn{efficient} access to matrix--vector products $\vec{z} \mapsto \mat{A}^{1/2}\vec{z}$ with the square root $\mat{A}^{1/2}$.}
I know of two scenarios where the square-root matvec model applies.
First, if the matrix $\mat{A}$ is \warn{sufficiently well-conditioned}, then matvecs with $\mat{A}^{1/2}$ can be computed efficiently using the Lanczos algorithm \cite[Ch.~6]{Che24}.
Second, if $\mat{A} = f(\mat{M})$ is a \warn{nonnegative} function $f : \real \to \real_+$, then matvecs with $\mat{A}$ are often computed with the Lanczos algorithm.
In this setting, matvecs with $\mat{A}^{1/2} = f^{1/2}(\mat{M})$ can  often be obtained at the \warn{same} computational cost by applying the Lanczos method with $f^{1/2}$ rather than $f$.
The first setting was observed by Mathur, Moka, and Botev \cite{MMB21}, and the second is used in \cite{Lin23a}.

In the square-root matvec model, we can estimate $\diag(\mat{A}) = \srn(\mat{A}^{1/2})$ by estimating the squared row norms of $\mat{A}^{1/2}$.
This approach to diagonal estimation is called the \emph{square-root trick}.
It was developed by Mathur, Moka, and Botev \cite{MMB21} and has been used to great effect in the works of Michael Lindsey \cite{Lin23a,FL24}.
\Cref{sec:subgraph-2} will illustrate the power of this technique, using the subgraph centrality task from \cref{sec:subgraph-1} as an example.


\section{Variance-reduced row-norm estimators} \label{sec:variance-reduced-row-norm}

Row-norm estimators using a \HutchPP-style variance reduction technique were developed by Sobczyk and Luisier \cite{SL22,Sob24a}.
The original version of their estimator is constructed as follows.
Let $s$ be a number of matvecs, evenly divisible by four.
Generate a random matrix $\mat{\Omega} \in \field^{n\times (s/4)}$ and iid isotropic columns $\vec{\gamma}_1,\ldots,\vec{\gamma}_{s/4}$.
Begin by computing the product
\begin{equation*}
    \mat{Y} \coloneqq \mat{B}^*(\mat{B}\mat{\Omega});
\end{equation*}
%
Then, orthogonalize $\mat{Q} \coloneqq \orth(\mat{Y})$, obtaining a low-rank approximation $\Bhat \coloneqq \mat{B}\mat{Q}\mat{Q}^*\approx \mat{B}$.
We may decompose the squared row norms by incoking the identity
\begin{equation*}
    \srn(\mat{B}) = \srn(\mat{B}\mat{Q}) + \srn(\mat{B}(\Id - \mat{Q}\mat{Q}^*)),
\end{equation*}
which is immediate from the relation $\srn(\mat{B}) = \diag(\mat{B}\mat{B}^*)$.
The first term can be computed exactly by forming $\mat{B}\mat{Q}$, and the second term can be estimated with the Johnson--Lindenstrauss row-norm estimator \cref{eq:mc-row-norm-compare}.
This process yields the variance-reduced estimator
\begin{equation*}
    \srn_{\mathrm{SL}4} \coloneqq \srn(\mat{B}\mat{Q}) + \frac{1}{s/4}\sum_{i=1}^{s/4} |\mat{B}(\Id - \mat{Q}\mat{Q}^*)\vec{\gamma}_i|^2.
\end{equation*}
We shall call this the Sobczyk--Luisier 4 (SL4) row-norm estimator, as it apportions its matvecs in four batches of size $s/4$.
The SL4 method produces an unbiased estimate of the squared row norms, and it requires $3s/4$ matvecs with $\mat{B}$, $s/4$ matvecs with $\mat{B}^*$, and $\order(ns^2)$ additional arithmetic operations.
Code is provided in \cref{prog:sl4}.

\myprogram{Sobczyk--Luisier 4 algorithm for row-norm estimation.}{The \texttt{sqrownorms} subroutine is defined in \cref{prog:sqrownorms}.}{sl4}

Before moving on to leave-one-out\index{leave-one-out randomized algorithm} estimates of the row norms, let us mention a refinement to the Sobczyk--Luisier estimator.
If we generate a random matrix $\mat{\Omega} \in \field^{m\times (s/3)}$ with $s/3$ columns, we can compute a low-rank approximation $\Bhat = \mat{B}\mat{Q}\mat{Q}^*$ to $\mat{B}$ using $2s/3$ matvecs:
\begin{equation*}
    \mat{Y} \coloneqq \mat{B}^*\mat{\Omega}\quad \text{and} \quad \mat{Q} \coloneqq \orth(\mat{Y}),
\end{equation*}
leaving $s/3$ matvecs to estimate the row norms of the residual.
This procedure results in a modified version of the Sobczyk--Luisier estimator, which we call the Sobczyk--Luisier 3 (SL3) estimator:
\begin{equation*}
    \srn_{\mathrm{SL}3} \coloneqq \srn(\mat{B}\mat{Q}) + \frac{1}{s/3}\sum_{i=1}^{s/3} |\mat{B}(\Id - \mat{Q}\mat{Q}^*)\vec{\gamma}_i|^2.
\end{equation*}
The computational cost of the SL3 estimator is $2s/3$ matvecs with $\mat{B}$, $s/3$ matvecs with $\mat{B}^*$, and $\order(ns^2)$ additional arithmetic operations.
Code is provided in \cref{prog:sl3}.

\myprogram{Sobczyk--Luisier 3 algorithm for row-norm estimation.}{The \texttt{sqrownorms} subroutine is defined in \cref{prog:sqrownorms}.}{sl3}

How do the SL3 and SL4 estimates compare?
The SL3 estimator uses a cruder low-rank approximation than the SL4 estimator,  requiring one less step of subspace iteration.
For a fixed budget of $s$ matvecs, however, the rank of the SL3 estimator is higher than the SL4 estimator, $s/3$ versus $s/4$.
Its higher rank often makes the SL3 estimator preferable to the SL4 estimator.
For instance, when applied to a matrix $\mat{B}$ with exponentially decaying singular values $\sigma_i(\mat{B}) \le \alpha^i$, the SL4 estimator produces estimates of the squared row norms with root-mean-squared error of size $\approx \alpha^{s/2}$, whereas the SL3 estimator achieves a faster error decay $\approx \alpha^{2s/3}$.
Thus, for problems with rapid singular value decay, the SL3 estimator is preferable.

\section{\XRowNorm: A leave-one-out row-norm estimator} \label{sec:xrownorm}

The SL4 estimator is based on a rank-$(s/4)$ randomized SVD approximation to $\mat{B}^*$ with one step of subspace iteration, and the SL3 estimator is based on a rank-$(s/3)$ approximation with no subspace iteration.
By using a leave-one-out\index{leave-one-out randomized algorithm} design, we can enjoy the best of both worlds, employing a rank-$(s/3)$ approximation with one step of subspace iteration.
We call the resulting estimator \XRowNorm.
As we will see in \cref{sec:row-norm-experiments}, \XRowNorm and SL3 perform similarly, with \XRowNorm offering a small but noticeable benefit for some instances; see \cref{fig:row-norm-comparison}.

Begin by drawing a matrix $\mat{\Omega} \in \field^{n\times (s/3)}$ with iid isotropic columns, and define a rank-$(s/3)$ approximation to $\mat{B}^*$:
\begin{equation*}
    \mat{Y} \coloneqq \mat{B}^*(\mat{B}\mat{\Omega}), \quad \mat{Q} \coloneqq \orth(\mat{Y}), \quad \Bhat \coloneqq \mat{B}\mat{Q}\mat{Q}^*.
\end{equation*}
Now, employ the leave-one-out\index{leave-one-out randomized algorithm} method.
We obtain a family of downdated low-rank approximations
\begin{equation*}
    \mat{Q}_{(i)} \coloneqq \orth(\mat{Y}_{-i}), \quad \Bhat_{(i)} \coloneqq \mat{B}\QiQi.
\end{equation*}
The squared row norms can be decomposed as 
\begin{equation*}
    \srn(\mat{B}) = \srn\big(\mat{B}\mat{Q}_{(i)}^{\vphantom{*}}\big) + \srn\big(\mat{B}\big(\Id - \QiQi\big)\big).
\end{equation*}
We estimate the second term by the single-vector \RowNormName estimator using the left-out vector $\vec{\omega}_i$,
\begin{equation}  \label{eq:xrownorm-i}
    \hat{\srn}_i \coloneqq \srn\big(\mat{B}\mat{Q}_{(i)}^{\vphantom{*}}\big) + \big| \mat{B}\big(\Id - \QiQi\big) \vec{\omega}_i^{\vphantom{*}} \big|^2,
\end{equation}
and average to form the \XRowNorm estimator
\begin{equation} \label{eq:xrownorm}
    \hat{\srn}_{\mathrm{X}} \coloneqq \frac{1}{s/3} \sum_{i=1}^{s/3} \hat{\srn}_i = \frac{1}{s/3} \sum_{i=1}^{s/3} \left[\srn\big(\mat{B}\mat{Q}_{(i)}^{\vphantom{*}}\big) + \big| \mat{B}\big(\Id - \QiQi\big) \vec{\omega}_i^{\vphantom{*}} \big|^2 \right].
\end{equation}

\myprogram{Efficient implementation of the \XRowNorm estimator.}{Subroutines \texttt{diagprod}, \texttt{random\_signs}, and \texttt{cnormc} appear in \cref{prog:diagprod,prog:random_signs,prog:cnormc}.}{xrownorm}

\myparagraph{Efficient formula and implementation}
To evaluate the \XRowNorm estimator efficiently, we compute the products $\mat{G} \coloneqq \mat{B}\mat{\Omega}$ and $\mat{Y} \coloneqq \mat{B}^*\mat{G}$, then form the \QR decomposition $\mat{Y} = \mat{Q}\mat{R}$.
By \cref{thm:rsvd-downdate}, the randomized SVD downdating is described implicitly by the matrix $\mat{S}$, obtained by normalizing the columns of $\mat{R}^{-*}$.
We substitute the downdating formula \cref{eq:rsvd-downdate} into \cref{eq:xrownorm-i,eq:xrownorm}, and we simplify algebraically.
We omit the details.
The resulting formulas are
\begin{equation*}
    \hat{\srn}_i = \srn(\mat{Z}) - \srn(\mat{Z}\vec{s}_i) + | \vec{g}_i - \mat{Z}\vec{x}_i |^2
\end{equation*}
and 
\begin{equation*}
    \hat{\srn}_{\mathrm{X}} = \srn(\mat{Z}) + \frac{1}{s/3} \left[ -\srn(\mat{Z}\mat{S}) + \srn(\mat{G} - \mat{Z}\mat{X}) \right].
\end{equation*}
where $\mat{Z} \coloneqq \mat{B}\mat{Q}$ and 
\begin{equation} \label{eq:xrownorm-X}
    \mat{X} \coloneqq \mat{W} - \mat{S} \cdot \Diag(\diagprod(\mat{S},\mat{W})) \quad \text{with } \mat{W} \coloneqq \mat{Q}^*\mat{\Omega}.
\end{equation}
Observe that the formula for $\mat{X}$ is the same as the formula \cref{eq:xtrace-X} for \XTrace.
The computational cost of \XRowNorm is $2s/3$ matvecs with $\mat{B}$, $s/3$ matvecs with $\mat{B}^*$, and $\order((m+n)s^2)$ additional arithmetic operations.
An implementation of \XRowNorm is provided in \cref{prog:xrownorm}.

\section{\XSymRowNorm: Improved Hermitian row-norm estimation} \label{sec:xsymrownorm}

For Hermitian matrices, we can develop an improved version of \XRowNorm that uses an approximation with larger rank $s/2$.
We call the resulting estimator \XSymRowNorm.

Let $\mat{A} \in \field^{n\times n}$ be a Hermitian matrix.
Begin by drawing a matrix $\mat{\Omega} \in\field^{n\times (s/2)}$ and calculating
\begin{equation*}
    \mat{Y} \coloneqq \mat{A}\mat{\Omega}, \quad \mat{Q} \coloneqq \orth(\mat{Y}).
\end{equation*}
Since $\mat{A}$ is Hermitian, the left-sided randomized SVD approximation $\Ahat_{\mathrm{left}} = \mat{Q}\mat{Q}^*\mat{A}$ has the same quality as the right-sided approximation $\Ahat = \mat{A} \mat{Q}\mat{Q}^*$.
Using the right-sided approximation, we may decompose the squared row norms
\begin{equation*}
    \srn(\mat{A}) = \srn(\mat{A}\mat{Q}) + \srn(\mat{A} (\Id - \mat{Q}\mat{Q}^*)).
\end{equation*}
Now, invoke the leave-one-out\index{leave-one-out randomized algorithm} method and introduce downdated $\mat{Q}$ matrices 
\begin{equation*}
    \mat{Q}_{(i)} \coloneqq \orth(\mat{Y}_{-i}).
\end{equation*}
Using an orthogonal decomposition of the squared row norms and the a single-vector \RowNormName estimator for the residual squared norms yields the basic \XSymRowNorm estimators
\begin{equation*}
    \hat{\srn}_i \coloneqq \srn\big(\mat{A}\mat{Q}_{(i)}\big) + \big| \mat{A} \big(\Id - \QiQi\big)\vec{\omega}_i^{\vphantom{*}} \big|^2.
\end{equation*}
Finally, average to produce the \XSymRowNorm estimator
\begin{equation*}
    \hat{\srn}_{\mathrm{XS}} \coloneqq \frac{1}{s/2} \sum_{i=1}^{s/2} \hat{\srn}_i = \frac{1}{s/2} \sum_{i=1}^{s/2} \left[ \srn\big(\mat{A}\mat{Q}_{(i)}\big) + \big| \mat{A} \big(\Id - \QiQi\big)\vec{\omega}_i^{\vphantom{*}} \big|^2 \right].
\end{equation*}

\myprogram{Efficient implementation of the \XSymRowNorm estimator.}{Subroutines \texttt{diagprod}, \texttt{random\_signs}, and \texttt{cnormc} appear in \cref{prog:diagprod,prog:random_signs,prog:cnormc}.}{xsymrownorm}

\myparagraph{Efficient formula and implementation}
The efficient formula for \XSymRowNorm is similar to \XRowNorm.
Compute the matrix product $\mat{Y} \coloneqq \mat{A}\mat{\Omega}$, form the \QR decomposition $\mat{Y} = \mat{Q}\mat{R}$, and produce the downdating matrix $\mat{S}$ by normalizing the columns of $\mat{R}^{-*}$.
Set $\mat{Z} \coloneqq \mat{A}\mat{Q}$ and produce the matrix $\mat{X}$ from \cref{eq:xrownorm-X}.
The basic \XSymRowNorm estimates are
\begin{equation*}
    \hat{\srn}_i = \srn(\mat{Z}) - \srn(\mat{Z}\vec{s}_i) + |\vec{y}_i - \mat{Z}\vec{x}_i|^2
\end{equation*}
and the \XSymRowNorm estimator is
\begin{equation*}
    \hat{\srn}_{\mathrm{XS}} = \srn(\mat{Z}) + \frac{1}{s/2} \left[ -\srn(\mat{Z}\mat{S}) + \srn(\mat{Y} - \mat{Z}\mat{X}) \right].
\end{equation*}

The computational cost of \XSymRowNorm is $s$ matvecs with $\mat{A}$ and $\order(ns^2)$ additional arithmetic operations.

\begin{remark}[\XSymRowNorm for non-Hermitian matrices]
    While we have developed the \XSymRowNorm estimator for the purpose of estimating the row norms of a Hermitian matrix, the estimator also works for normal matrices or square matrices that are nearly Hermitian.
    The \XSymRowNorm estimator is always unbiased, is effective when low-rank approximations $\mat{Q}\mat{Q}^*\mat{B}$ and $\mat{B}\mat{Q}\mat{Q}^*$ obtained by applying a projector to the left or right side are similarly good low-rank approximations.
\end{remark}






\section{Synthetic Experiments} \label{sec:row-norm-experiments}

\begin{figure}
    \centering
    \includegraphics[width=0.99\linewidth]{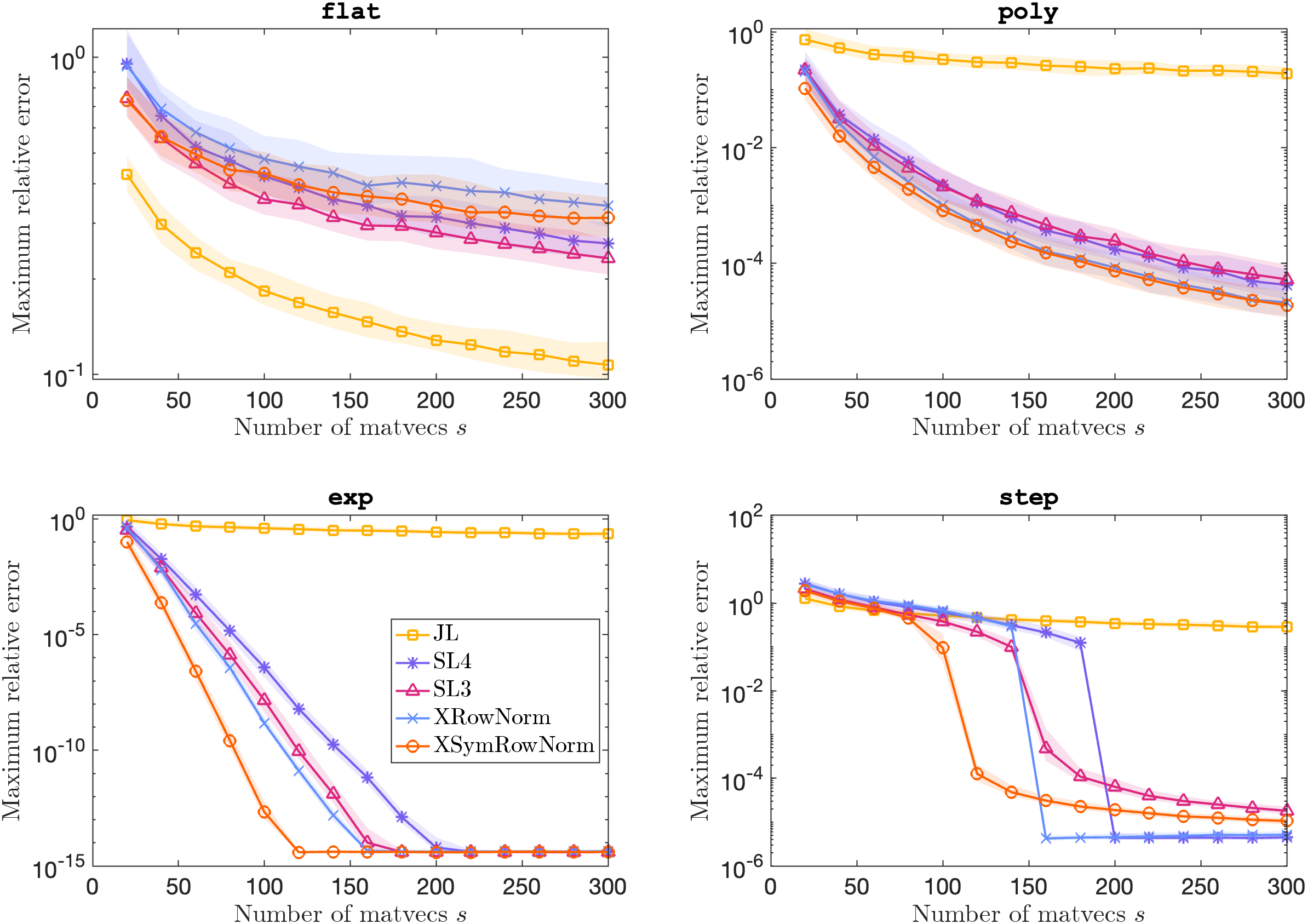}
    \caption[Comparison of JL, SL4, SL3, \XRowNorm, and \XSymRowNorm estimators for row-norm estimation of matrices with different spectra]{Relative error of squared row-norm estimates by the JL estimator (yellow squares), SL4 estimator (purple asterisks), SL3 (pink triangles), \XRowNorm (blue crosses), and \XSymRowNorm (orange circles) for the four test matrices \cref{eq:test-matrices} for different numbers of matvecs $s$.
    The error metric is maximum relative error, defined in \cref{eq:max-rel-err-srn}.
    Lines show median of 100 trials, and error bars show 10\% and 90\% quantiles.}
    \label{fig:row-norm-comparison}
\end{figure}

\Cref{fig:row-norm-comparison} shows a comparison of row-norm estimators for the four test matrices from \cref{eq:test-matrices}.
For each method, we plot the maximum relative error in the \warn{squared} row norms:
\begin{equation} \label{eq:max-rel-err-srn}
    \mathrm{maxrelerr} = \max_{1\le i \le n} \frac{|\mathrm{\hat{srn}}_i - \mathrm{srn}_i|}{|\mathrm{srn}_i|}.
\end{equation}
The results are broadly similar to the experimental results for trace and diagonal estimation in \cref{fig:trace-comparison,fig:diagonal-comparison}.
For problems with rapid spectral decay, most notably the \texttt{exp} matrix, the errors of the methods are sorted \XSymRowNorm < \XRowNorm < SL3 < SL4 $\ll$ JL.
The \XRowNorm and SL3 estimators perform similarly on this example, but the \XRowNorm estimator edges it out by incorporating an additional step of subspace iteration.
The results for the \texttt{step} matrix also demonstrate the benefit of this extra subspace iteration step, with SL4 and \XRowNorm outperforming SL3 and \XSymRowNorm for sufficiently large values of $s$.

\section{Application: Subgraph centralities, again} \label{sec:subgraph-2}

Now equipped with the square-root trick and varianced-reduced diagonal estimators, we return to the task of estimating the subgraph centralities from \cref{sec:subgraph-1}.
Recall that the subgraph centralities of a graph with adjacency matrix $\mat{M}$ are
\begin{equation*}
    \mathbf{sc} \coloneqq \diag(\exp(\mat{M})).
\end{equation*}
In \cref{sec:subgraph-1}, we performed matvecs with $\mat{A} = \exp(\mat{M})$ using the Lanczos method \cite[Ch.~6]{Che24} associated with the function $f(t) = \exp(t)$.
To apply matvecs with $\mat{A}^{1/2}$, we can simply change the function to $f(t) = \exp(t/2)$.
Thus, for this problem, matvecs with $\mat{A}^{1/2}$ are just as cheap as matvecs with $\mat{A}$, making this an ideal setting to apply the square-root trick.

\begin{figure}
    \centering
    \includegraphics[width=0.99\linewidth]{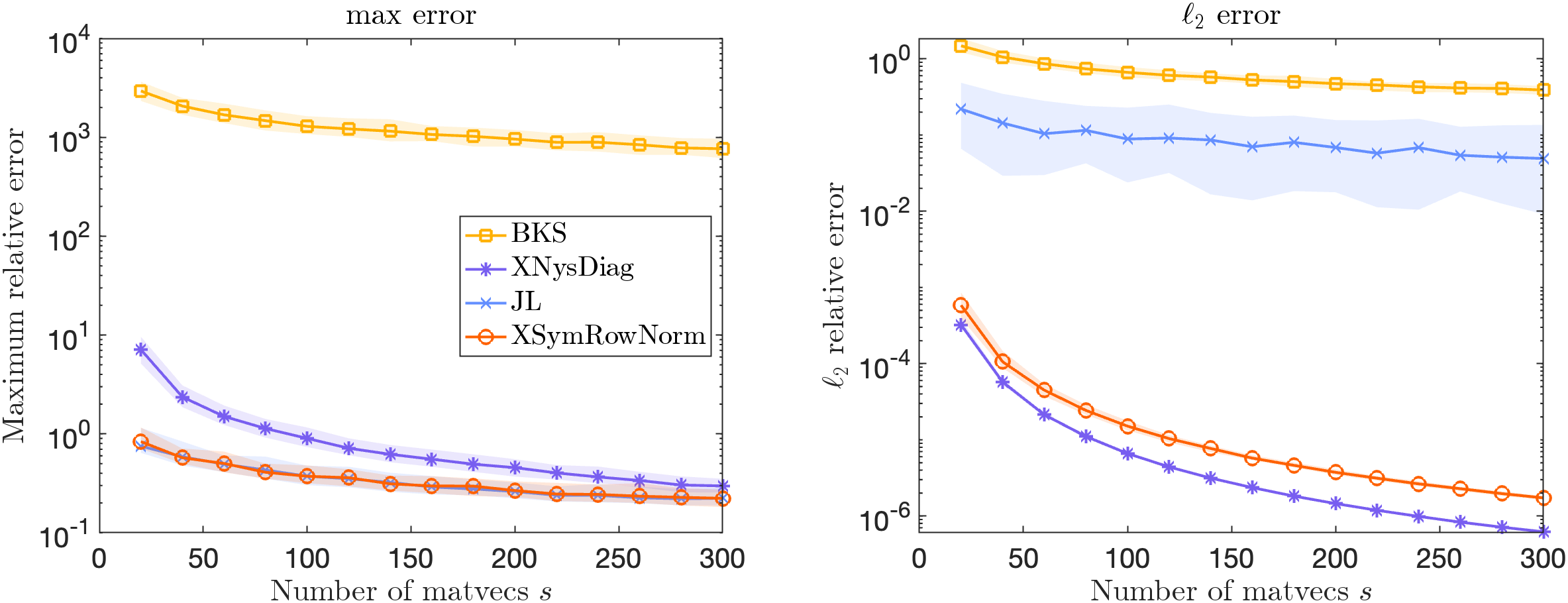}
    \caption[Comparison of BKS and \XNysDiag diagonal estimators to JL and \XSymRowNorm row-norm estimators for estimating subgraph centralities]{Maximum relative error of subgraph centrality estimates by the BKS diagonal estimator (yellow squares), \XNysDiag (purple asterisks), JL row-norm estimator (blue crosses), and \XSymRowNorm (orange circles) as a function of the number of matvecs $s$.
    We show both the the maximum relative error (\emph{left}, \cref{eq:max-rel-err-diag}) and the relative $\ell_2$ error (\emph{right}).
    Lines show the median of 100 trials, and shaded regions show 10\% and 90\% quantiles.}
    \label{fig:subgraph-2}
\end{figure}

\Cref{fig:subgraph-2} measures the accuracy of the subgraph centrality estimated by the BKS and \XNysDiag diagonal estimators and by the JL and \XSymRowNorm row-norm estimators.
When evaluated using the maximum relative error (\emph{left}), we see that the JL and \XSymRowNorm estimators dramatically outperform both diagonal estimators, achieving nearly $10\times$ smaller error using $s = 20$ matvecs.
These results demonstrate how the square-root trick can yield diagonal estimates that are much more accurate when the diagonal entries span many orders of magnitude.
When using the $\ell_2$ error (\emph{right}), the variance-reduced \XNysDiag and \XSymRowNorm algorithms significantly outperform the Monte Carlo-based JL and BKS estimates.
In particular, the \XSymRowNorm estimator achieves the best of both worlds, attaining both small maximum relative error and small $\ell_2$ error.

\chapter{Posterior error estimates for low-rank approximation}

\epigraph{\ldots{we} advocate using \emph{a posteriori} error estimators to assess the quality of the output of a randomized linear algebra computation. These error estimators are often quite cheap, yet they can give (statistical) evidence that the computation was performed correctly.
We also recommend adaptive algorithms that can detect when the accuracy is insufficient and make refinements.
With this approach, it is not pressing to produce theory that justifies all of the internal choices (e.g., the specific type of random embedding) in the NLA algorithm.}{Per-Gunnar Martinsson and Joel A.\ Tropp, \emph{Randomized numerical linear algebra: Foundations and algorithms} \cite{MT20a}}

In previous chapters, we employed the leave-one-out\index{leave-one-out randomized algorithm} approach to estimate attributes of a matrix $\mat{B}$.
In this chapter and the next one, we develop a conceptually distinct application of the leave-one-out\index{leave-one-out randomized algorithm} technique.
Rather than using leave-one-out\index{leave-one-out randomized algorithm} to \emph{estimate matrix attributes} as previously, we will use the leave-one-out\index{leave-one-out randomized algorithm} method to \emph{diagnose the quality} of randomized matrix approximations.
This chapter will deploy the leave-one-out\index{leave-one-out randomized algorithm} approach to compute error estimates for low-rank approximation, and next chapter will apply this machinery to estimate the variance of a general class of randomized matrix algorithms.
The downdating formulas (\cref{thm:rsvd-downdate,thm:nystrom-downdate}) will be essential tools in devising efficient algorithms.

\myparagraph{Sources}
This chapter is based on the paper

\fullcite{ET24}.

\myparagraph{Outline}
\Cref{sec:loo-error-estimation} describes the leave-one-out\index{leave-one-out randomized algorithm} approach to error estimation in generality, and \cref{sec:loo-rsvd} applies this technique to the randomized SVD.
\Cref{sec:loo-error-experiments} contains numerical experiments.

\section{Leave-one-out error estimation: General approach} \label{sec:loo-error-estimation}

We begin by developing the leave-one-out\index{leave-one-out randomized algorithm} error estimation technique in a general setting; \iffull\else the subsequent section \cref{eq:loo-rsvd-s} will apply this error estimation technique to the randomized SVD.\fi

\myparagraph{Setting}
Let $\mat{B} \in \field^{m\times n}$ be an input matrix, and let $\mat{\Omega} \in \field^{n\times k}$ be a random test matrix with iid isotropic columns. 
Consider a randomized approximation to $\mat{B}$ that depends on $\mat{\Omega}$:
\begin{equation} \label{eq:some-matrix-approximation}
    \Bhat = \Bhat_k = \Bhat(\mat{\Omega}).
\end{equation}
We assume that forming the approximation $\Bhat$ requires evaluating $\mat{B}\mat{\Omega}$, and we assume that $\Bhat(\mat{\Omega})$ is defined for input matrices $\mat{\Omega}$ with any number of columns.
Examples include the randomized SVD and the randomized Nystr\"om approximations, introduced in \cref{sec:lra}.
We will abuse notation and write $\Bhat_k$, $\Bhat(\mat{\Omega})$, and $\Bhat$ interchangeably to refer to the approximation \cref{eq:some-matrix-approximation}.
For simplicity, we shall assume that the approximation $\Bhat$ only depends on the deterministic input matrix $\mat{B}$ and the random test matrix $\mat{\Omega}$, although the methods in this section apply equally well to approximations $\Bhat(\mat{\Omega},\mat{\Psi},\ldots)$ that depend on additional random quantities $\mat{\Psi},\ldots$ that are independent of $\mat{\Omega}$.

We are interested in efficiently estimating the mean-squared Frobenius error
\begin{equation*}
\MSE_k \coloneqq \expect \, \norm{\mat{B} - \Bhat_k}_{\mathrm{F}}^2.
\end{equation*}
This type of estimate has several applications.
For example, it can be used to adaptively select the number of matvecs $k$ that suffice to meet a target accuracy $\tau$:
Simply increment $k$ until the estimate $\hat{\MSE}$ satisfies $\hat{\MSE} \le \tau$.

\myparagraph{Existing approach: Girard--Hutchinson norm estimator}
The Girard--Hutchinson norm estimator \cite[\S4.8]{MT20a} provides one way of estimating the error $\norm{\mat{B} - \Bhat}_{\mathrm{F}}$.
Given a matrix $\mat{C} \in \field^{m\times n}$, the Girard--Hutchinson estimate of the squared Frobenius norm is
\begin{equation}
    \hat{\mathrm{SN}} \coloneqq \frac{1}{s} \norm{\mat{C}\vec{\gamma}_i}^2, \label{eq:gh-norm}
\end{equation}
where $\vec{\gamma}_1,\ldots,\vec{\gamma}_s$ are iid isotropic vectors.
Observe that the Girard--Hutchinson squared norm estimate $\hat{\mathrm{SN}}$ coincides with the Girard--Hutchinson estimator for $\tr(\mat{C}^*\mat{C}) = \norm{\mat{C}}_{\mathrm{F}}^2$.
Applying the Girard--Hutchinson norm estimator to $\mat{C} = \mat{B} - \Bhat$ gives an unbiased estimate to the mean-squared error $\MSE_k$.
As disadvantages of this approach, computing the estimate $\hat{\mathrm{SN}}$ requires $s$ fresh matvecs with the matrix $\mat{B}$, and the only way to improve the estimate is to increase the number of matvecs $s$.
The leave-one-out\index{leave-one-out randomized algorithm} approach will give us an alternative error estimation technique that require no additional matvecs with $\mat{B}$ and that automatically improves with the approximation quality parameter $k$.


\myparagraph{Leave-one-out errorr estimation}
By now, our approach to designing a leave-one-out\index{leave-one-out randomized algorithm} error estimate should be routine.
Begin by downating the approximation $\Bhat$ by leaving out each column $\vec{\omega}_i$ out in turn, producing a family of approximations 
\begin{equation*}
    \Bhat_{(i)} \coloneqq \Bhat(\mat{\Omega}_{-i}) \quad \text{for } i =1,\ldots,k.
\end{equation*}
Then, use the left-out vector $\vec{\omega}_i$ to form a Girard--Hutchinson norm estimate of each downdated approximation:
\begin{equation*}
    \hat{\MSE}_i \coloneqq \norm{(\mat{B} - \Bhat)\vec{\omega}_i}^2 \approx \norm{\mat{B} - \Bhat_{(i)}}^2_{\mathrm{F}}.
\end{equation*}
Finally, average these basic estimates to obtain the leave-one-out\index{leave-one-out randomized algorithm} error estimator
\begin{equation*}
    \hat{\MSE} \coloneqq \frac{1}{k} \sum_{i=1}^k \hat{\MSE}_i = \frac{1}{k} \sum_{i=1}^k \norm{(\mat{B} - \Bhat_{(i)})\vec{\omega}_i}^2.
\end{equation*}

This leave-one-out\index{leave-one-out randomized algorithm} error estimation $\hat{\MSE}$ has several attractive features.
First, the error estimate requires with no additional matvecs with $\mat{B}$ beyond those necessary to form $\mat{B}\mat{\Omega}$.
Second, the quality of this error estimate automatically improves as the accuracy parameter $k$ is increased.
The leave-one-out error estimator does have one significant flaw:
The estimator \warn{$\hat{\MSE}$ is an unbiased estimate of the mean-squared error of the mean-squared error of the approximation $\Bhat_{k-1}$ with parameter $k-1$!}

\begin{proposition}[Leave-one-out error estimator]
    Instate the prevailing notation.
    The leave-one-out\index{leave-one-out randomized algorithm} estimator $\hat{\MSE}$ is an unbiased estimator for $\MSE_{k-1}$.
    That is,
    \begin{equation*}
        \expect[\hat{\MSE}] = \MSE_{k-1}.
    \end{equation*}
\end{proposition}

The mean-squared error for many types of approximation $\Bhat_k$, such as the randomized SVD and the randomized Nystr\"om approximations, are monotone $\MSE_k \le \MSE_{k-1}$.
For such algorithms, the leave-one-out\index{leave-one-out randomized algorithm} estimate $\hat{\MSE}$ \warn{overestimates} of the error, on average.

Since it is often more convenient to work with the error rather than the squared error, we define the leave-one-out\index{leave-one-out randomized algorithm} estimate of the \emph{error} to be
\begin{equation*}
    \hat{\Err} \coloneqq \smash{\hat{\MSE}}^{1/2}.
\end{equation*}


\section{Randomized SVD error estimation} \label{sec:loo-rsvd}

Efficient implementations of the leave-one-out\index{leave-one-out randomized algorithm} error estimates for the randomized SVD can be derived using the randomized SVD downdating formula \cref{eq:rsvd-downdate}.
The implementation takes different forms with, and without, additional subspace iteration.
We will present only the version without subspace iteration; see \cite{ET24} for the version with subspace iteration.

Generate a random matrix $\mat{\Omega}$ with iid isotropic columns, and form the randomized SVD approximation $\Bhat \coloneqq \mat{Q}\mat{Q}^*\mat{B}$ with a matrix product $\mat{Y} \coloneqq \mat{B}\mat{\Omega}$ and \QR decomposition $\mat{Y} = \mat{Q}\mat{R}$.
For the derivation, introduce the downdating matrix
\begin{equation} \label{eq:loo-rsvd-s}
    \mat{S} \coloneqq \mat{R}^{-*} \cdot \Diag\bigl( \norm{\mat{R}^{-*}(:,i)}^{-1} : 1\le i \le k\bigr),
\end{equation}
which contains the normalized columns of $\mat{R}^{-*}$.
The leave-one-out\index{leave-one-out randomized algorithm} error estimate is 
\begin{equation*}
    \MSE = \frac{1}{k} \sum_{i=1}^k \norm{ (\Id - \mat{Q} (\Id - \vec{s}_i^{\vphantom{*}}\vec{s}_i^*)\mat{Q}^* \mat{B}\vec{\omega}_i }^2.
\end{equation*}
By construction, $\Id - \mat{Q}\mat{Q}^*$ is an orthoprojector annihilating $\range(\mat{Y}) = \range(\mat{B}\mat{\Omega})$, from which $(\Id - \mat{Q}\mat{Q}^*)\mat{B}\vec{\omega}_i$.
In addition, $\mat{Q}^*\mat{B}\mat{\Omega} = \mat{Q}^*\mat{Y} = \mat{R}$, so $\mat{Q}^*\mat{B}\vec{\omega}_i = \vec{r}_i$.
Thus, 
\begin{equation*}
    \MSE = \frac{1}{k} \sum_{i=1}^k \norm{\mat{Q}\vec{s}_i^{\vphantom{*}}\vec{s}_i^*\vec{r}_i^{\vphantom{*}}}^2 = \frac{1}{k} \sum_{i=1}^k |\vec{s}_i^*\vec{r}_i^{\vphantom{*}}|^2.
\end{equation*}
In the second equality, we use the fact that $\mat{Q}\vec{s}_i$ is unit vector.
Finally, \cref{eq:loo-rsvd-s} implies that
\begin{equation*}
    \vec{s}_i^*\vec{r}_i^{\vphantom{*}} \!=\! (\mat{S}^*\mat{R})_{ii} \!=\! \Diag\bigl( \norm{\mat{R}^{-*}(:,i)}^{-1} : 1\le i \le k\bigr)_{ii} \!=\! \norm{\mat{R}^{-*}(:,i)}^{-1} = \norm{\smash{\mat{R}^{-1}(i,:)}}^{-1}.
\end{equation*}
Ergo, the leave-one-out\index{leave-one-out randomized algorithm} estimate of the mean-squared error is
\begin{equation*}
    \MSE = \frac{1}{k} \sum_{i=1}^k \norm{\mat{R}^{-1}(i,:)}^{-2},
\end{equation*}
and the leave-one-out\index{leave-one-out randomized algorithm} error estimate is
\begin{equation} \label{eq:loo-error-estimate-rsvd}
    \Err \coloneqq \left( \frac{1}{k} \sum_{i=1}^k \norm{\mat{R}^{-1}(i,:)}^{-2} \right)^{1/2}.
\end{equation}
The $\mat{R}$ matrix, often an unused byproduct of the standard randomized SVD implementation, contains enough information by itself to provide an estimate of the error.
See \cref{prog:rsvd_errest} for code.

\myprogram{Randomized SVD for matrix low-rank approximation together with leave-one-out error estimate.}{Subroutine \texttt{random\_signs} is provided in \cref{prog:random_signs}.}{rsvd_errest}

The computational cost of the leave-one-out\index{leave-one-out randomized algorithm} error estimator is just $\order(k^3)$ operations to invert the matrix $\mat{R}$.
In particular, the cost of the leave-one-out\index{leave-one-out randomized algorithm} error estimator is independent of \emph{both} the dimensions $m$ and $n$ of the input matrix $\mat{B}$, and it is always faster than the $\order((m+n)k^2)$ \warn{post-processing cost} of a standard randomized SVD implementation.

\begin{remark}[Even faster leave-one-out error estimation] \label{rem:faster-loo}
    One can accelerate the leave-one-out error estimator even further by estimating the row norms of $\mat{R}^{-1}$ using a stochastic estimator from \cref{ch:row-norm}.
    Using the \RowNormName row-norm estimator (\cref{eq:jl-row-norm}), $\order(\log k)$ matvecs suffice to estimate all row norms of $\mat{R}^{-1}$ up to a constant relative error \cite[Thm.~4]{Woo14a}.
    Further, each matvec $\mat{R}^{-1}\vec{\gamma}$ can be computed in $\order(k^2)$ operations.
    Therefore, the accelerated leave-one-out\index{leave-one-out randomized algorithm} error estimator requires just $\order(k^2 \log k)$ operations, which is much less than the randomized SVD.
\end{remark}

\section{Experiments} \label{sec:loo-error-experiments}

\begin{figure}
    \centering
    \includegraphics[width=0.99\linewidth]{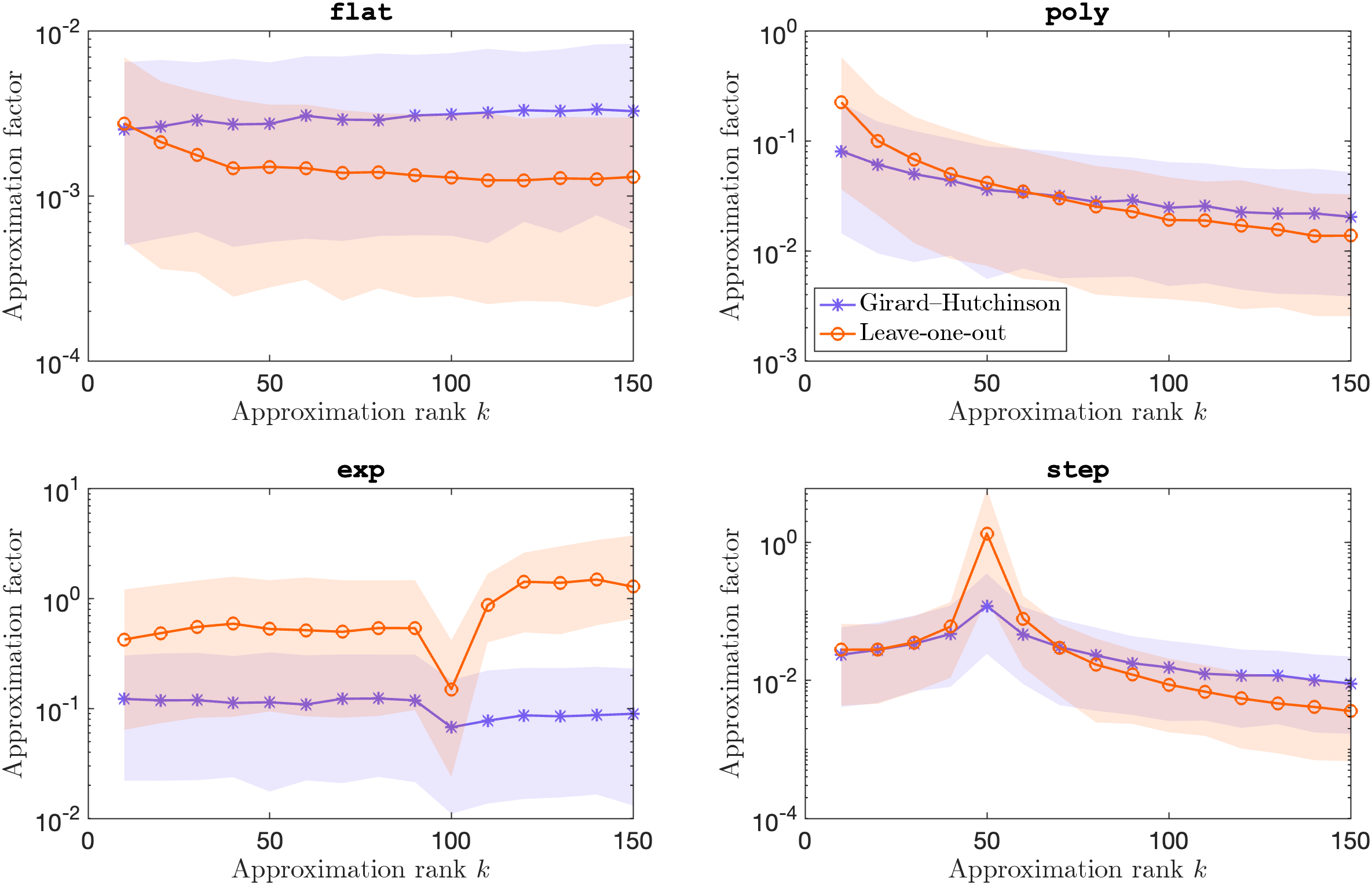}
    \caption[Comparison of leave-one-out and Girard--Hutchinson randomized SVD error estimates for matrices with different singular value profiles]{Approximation factor for the leave-one-out (orange circles) and Girard--Hutchinson (purple asterisks, $s=10$ matvecs) estimates for the error of the randomized SVD on the test matrices from \cref{eq:test-matrices}
    The error metric is the approximation factor, defined in \cref{eq:approximation-factor-errest}.
    Lines trace the median of 1000 trials, and error bars show the 10\% and 90\% quantiles.}
    \label{fig:err-est-comparison}
\end{figure}

\Cref{fig:err-est-comparison} presents a comparison between the leave-one-out\index{leave-one-out randomized algorithm} error estimator \cref{eq:loo-error-estimate-rsvd} and the Girard--Hutchinson norm estimator \cref{eq:gh-norm} with $s=10$ matvecs applied to the randomized SVD on the four test matrices in \cref{eq:test-matrices}.
We measure the quality of the error estimates using the approximation factor
\begin{equation} \label{eq:approximation-factor-errest}
    \alpha \coloneqq \max \left\{ \frac{\mathrm{est}}{\mathrm{err}}, \frac{\mathrm{err}}{\mathrm{est}} \right\} - 1.
\end{equation}
This is a sensible measure of the approximation quality, because an error estimate with an approximation factor of $\alpha$ is guaranteed to satisfy
\begin{equation*}
    \frac{1}{1+\alpha} \cdot \mathrm{est} \le \mathrm{err} \le (1+\alpha) \mathrm{est}.
\end{equation*}
Here are my conclusions:

\myparagraph{Leave-one-out\index{leave-one-out randomized algorithm} error-estimator: pretty good for being free}
Before commenting on the finer distinctions between the leave-one-out\index{leave-one-out randomized algorithm} error estimator and the Girard--Hutchinson norm estimator, it is worth emphasizing from the outset that the leave-one-out\index{leave-one-out randomized algorithm} error estimator has a basic advantage over the Girard--Hutchinson norm estimator: The leave-one-out\index{leave-one-out randomized algorithm} error estimator requires no additional matrix--vector products with the matrix $\mat{B}$, and it requires a small number of arithmetic operations to calculate.
The Girard--Hutchinson estimator requires fresh matrix--vector products with the matrix $\mat{B}$.
Thus, even on examples where the Girard--Hutchinson norm estimator is better, the leave-one-out\index{leave-one-out randomized algorithm} error estimator may still be preferable in many applications, since it is basically free to compute.

\myparagraph{Self-improving}
Another advantage of the leave-one-out\index{leave-one-out randomized algorithm} error estimator is that it is \emph{self-improving}.
When one increases the approximation rank $k$, the estimator automatically improves.
This differs from the Girard--Hutchinson estimator, which requires additional matvecs to meaningfully improve quality of the estimator.

\myparagraph{Overestimating the error when the singular values drop sharply}
The major weakness of the leave-one-out\index{leave-one-out randomized algorithm} error estimator is that it provides an unbiased estimate of the \warn{rank-$(k-1)$} randomized SVD error, which can overestimate the rank-$k$ randomized SVD error when the singular values decrease rapidly.
This phenomenon is visible at values $k\le 100$ in the \texttt{exp} example and at $k=50$ in the \texttt{step} example.

\myparagraph{When the approximation becomes accurate to machine precision}
Finally, for $k>100$ on the \texttt{exp} example, the randomized SVD approximation becomes accurate to machine precision $\norm{\mat{B} - \Bhat}_{\mathrm{F}} \approx 10^{-16}$.
When this happens, the leave-one-out\index{leave-one-out randomized algorithm} error estimator overestimates the true error by a factor of about 3 due to numerical issues.
Fortunately, for most applications, small factors do not matter when an approximation is accurate to machine-precision.

\chapter{Matrix jackknife variance estimation}

\epigraph{\emph{Good} simple ideas, of which the jackknife is a prime example, are
our most precious intellectual commodity, so there is no need to apologize for
the easy mathematical level.}{Bradley Efron, \textit{The Jackknife, the Bootstrap and Other Resampling Plans} \cite[p.~1]{Efr82}}

In the previous chapter, we discussed the leave-one-out\index{leave-one-out randomized algorithm} error estimator as a way to assess the error of an approximation $\Bhat\approx\mat{B}$ computed by a randomized algorithm.
In this chapter, we shall turn our attention to assessing the quality of randomized approximations to quantities $\mat{Q}(\mat{B})$ obtained from \emph{transforming} a matrix $\mat{B} \in \field^{m\times n}$ by a \emph{nonlinear} function $\mat{Q} : \field^{m\times n}\to\field^{m'\times n'}$.
There are many examples of such transformations:
\begin{enumerate}
    \item \textbf{\textit{Matrix functions.}} $\mat{Q}(\mat{B}) = f(\mat{B})$, where $f : \real \to \real$ is a function.
    \item \textbf{\textit{Best rank-$r$ approximation.}} $\mat{Q}(\mat{B}) = \lowrank{\mat{B}}_r$.
    \item \textbf{\textit{Top-$r$ singular values.}} $\mat{Q}(\mat{B}) = (\sigma_1(\mat{B}),\ldots,\sigma_r(\mat{B}))$ for given parameter $r$.
    \item \textbf{\textit{Dominant singular subspace projectors.}} $\mat{Q}(\mat{B}) = \mat{U}(:,1:r)\mat{U}(:,1:r)^*$ or $\mat{Q}(\mat{B}) = \mat{V}(:,1:r)\mat{V}(:,1:r)^*$, where $\mat{B}=\mat{U\Sigma V^*}$ is an SVD and $r$ is a given parameter.
\end{enumerate}
Each of these quantities can be estimated by applying a randomized algorithm to produce an initial matrix approximation $\Bhat\approx\mat{B}$ and then using $\mat{Q}(\Bhat)$ as a proxy for $\mat{Q}(\mat{B})$.
This chapter addresses the following question:
\actionbox{How can we assess the quality of the approximation $\mat{Q}(\Bhat)\approx\mat{Q}(\mat{B})$ at runtime?}

This question is nontrivial. 
In contrast to the leave-one-out\index{leave-one-out randomized algorithm} error estimation setting, we typically do not have access to the true quantity of interest $\mat{Q}(\mat{B})$, even via matvecs. 
Subject to this limitation, it may be infeasible or computationally intractable to directly estimate the \emph{error} $\norm{\mat{Q}(\mat{B})-\mat{Q}(\Bhat)}_{\mathrm{F}}$, and we must settle for other diagnostics of the quality of the computed solution.
One useful diagnostic is the \emph{variance} of the estimator $\mat{Q}(\Bhat)$.
The significance of the variance can be justified by appealing to the following intuition \cite[p.~A511]{ET24}:

\actionbox{In order to trust the answer provided by a randomized algorithm, the output should be insensitive to the randomness used by the algorithm}

The variance provides a quantitative measure of the random fluctuations in the output of a randomized algorithm, allowing us to put this principle into action.

This chapter develops a general approach for estimating the variance of a randomized matrix algorithm based on \emph{jackknife resampling}, a standard technique in statistics.
Jackknife resampling is a leave-one-out\index{leave-one-out randomized algorithm} approach, allowing us to use our matrix downdating results to obtain highly efficient algorithms.

\myparagraph{Sources}
This chapter is based on the paper

\fullcite{ET24}.

\myparagraph{Outline}
\Cref{sec:bias-variance} begins by describing the bias--variance decomposition of a matrix estimator.
We then introduce the matrix jackknife variance technique in generality in \cref{sec:matrix-jackknife}, and we discuss its use in \cref{sec:using-matrix-jackknife}.
\Cref{sec:spectral-downdating} describes how to efficiently compute matrix jackknife estimates for a broad class of ``spectral transformations'' of a Nystr\"om approxiamtion to a psd matrix.

\section{Bias--variance decomposition} \label{sec:bias-variance}

Throughout this chapter, we shall be interested in estimating the \emph{variance} of a random matrix $\mat{X} \in \field^{m\times n}$, defined as
\begin{equation*}
    \Var(\mat{X}) \coloneqq \expect \, \norm{\mat{X} - \expect \mat{X}}_{\mathrm{F}}^2.
\end{equation*}
The matrix variance is the sum of the scalar variances of each entry of $\mat{X}$:
\begin{equation} \label{eq:matrix-variance-decomposition}
    \Var(\mat{X}) = \sum_{i=1}^m \sum_{j=1}^n \Var(x_{ij}) = \sum_{i=1}^m \sum_{j=1}^n \expect |x_{ij} - \expect x_{ij}|^2.
\end{equation}
The variance measures the variability of $\mat{X}$ with respect to the \warn{Frobenius norm}; for discussion of other Schatten norms, see \cite[\S6]{ET24}.
The standard deviation is defined as $\SD(\mat{X}) \coloneqq \Var(\mat{X})^{1/2}$.

Suppose $\Bhat$ is a random estimator for the matrix $\mat{B}$.
One convenient measure for the approximation quality $\Bhat\approx \mat{B}$ is the \emph{mean-squared error}
\begin{equation*}
    \MSE \coloneqq \expect \norm{\mat{B} - \Bhat}_{\mathrm{F}}^2.
\end{equation*}
The mean-squared error admits a \emph{bias--variance decomposition}:
\begin{equation} \label{eq:bias-variance}
    \MSE = \Bias^2 + \Var(\Bhat) \quad \text{where } \Bias \coloneqq \norm{\mat{B} - \expect \Bhat}_{\mathrm{F}}.
\end{equation}
Consequently, the variance is a \warn{lower bound} on the mean-squared error.

\section{Matrix jackknife variance estimation} \label{sec:matrix-jackknife}

Jackknife resampling is an established approach for approximating the variance of a statistical estimator.
Typically, the method is applied to scalar statistics computed from observational data.
Here, we will hijack the basic machinery for the matrix-valued outputs of randomized algorithms.

We begin with the general formalism.
Let $\omega_1,\ldots,\omega_k$ be iid random variables taking values in a (measurable) space $\set{S}$, and use them to form a matrix-valued estimator $\mat{X}_k = \mat{X}_k(\omega_1,\ldots,\omega_k) \in \field^{m\times n}$.
We will abuse notation and use $\mat{X}_k$ to denote both the \warn{deterministic} function
\begin{equation*}
    \mat{X}_k : (\omega_1,\ldots,\omega_k) \in \set{S}^k \longmapsto \mat{X}_k (\omega_1,\ldots,\omega_k) \in \field^{m\times n}.
\end{equation*}
and the output of that function for the random inputs $\omega_1,\ldots,\omega_k$.
We assume that the estimator $\mat{X}_k$ is a \warn{permutation-invariant function} of its inputs $\omega_1,\ldots,\omega_k$.

In statistics and randomized matrix computations, we are often interested in a \emph{family} of estimators $(\mat{X}_k : k=1,2,\ldots)$ that can be instantiated for any natural number $k$.
Examples in statistics include the sample mean and variance:
\begin{align*}
    m_k(\omega_1,\ldots,\omega_k) &= \frac{1}{k} \sum_{i=1}^k \omega_i;\\
    v_k(\omega_1,\ldots,\omega_k) &= \frac{1}{k-1} \sum_{i=1}^k |\omega_i - m_k(\omega_1,\ldots,\omega_k)|^2. 
\end{align*}
Examples in randomized matrix computations include the randomized SVD and the randomized Nystr\"om approximation
\begin{align*}
    \Bhat(\vec{\omega}_1,\ldots,\vec{\omega}_k) &= \mat{\Pi}_{\mat{B}[\vec{\omega}_1 \:\cdots\: \vec{\omega}_k]} \mat{B} \quad \text{and} \quad \Ahat(\vec{\omega}_1,\ldots,\vec{\omega}_k) &= \mat{A}\left\langle [\vec{\omega}_1 \:\cdots\: \vec{\omega}_k] \right\rangle.
\end{align*}
We often assume the variance is monotone decreasing $\Var(\mat{X}_k) \le \Var(\mat{X}_{k-1})$ in the sample size $k$.

The jackknife variance estimate was introduced by Tukey in 1958 \cite{Tuk58}.
Here is the idea.
Define a family of \emph{jackknife replicates} by recomputing the estimator with each one of the samples $\omega_1,\ldots,\omega_k$ left out in turn:
\begin{equation*}
    \mat{X}^{(j)} \coloneqq \mat{X}_{k-1}(\omega_1,\ldots,\omega_{j-1},\omega_{j+1},\ldots,\omega_k).
\end{equation*}
The jackknife variance estimate is defined as
\begin{equation} \label{eq:jack}
    \Jack^2(\mat{X}_{k-1}) \coloneqq \sum_{j=1}^k \norm{\mat{X}^{(j)} - \mat{X}^{(\cdot)}}_{\mathrm{F}}^2\quad \text{where } \mat{X}^{(\cdot)} \coloneqq \frac{1}{k} \sum_{j=1}^k \mat{X}^{(j)}.
\end{equation}
The jackknife estimate serves as an approximation of the variance of the \warn{$(k-1)$-sample estimator $\mat{X}_{k-1}$}, and it typically provides a modest overestimate of the variance of the $k$-sample estimator $\mat{X}_k$ as well.

Observe that the formula \cref{eq:jack} for the jackknife variance estimator resembles the traditional sample covariance estimate
\begin{equation} \label{eq:sample-variance}
    \Var(\mat{X}_{k-1}) \approx \frac{1}{k-1}\sum_{j=1}^k \norm{\mat{Z}^{(j)} - \mat{Z}^{(\cdot)}}_{\mathrm{F}}^2 \quad \text{where } \mat{Z}^{(\cdot)} \coloneqq \frac{1}{k} \sum_{j=1}^k \mat{Z}^{(j)}
\end{equation}
for iid copies $\mat{Z}^{(1)},\ldots,\mat{Z}^{(k)}\sim \mat{X}_{k-1}$ of the estimator.
An important difference between the jackknife estimate \cref{eq:jack} and the sample covariance \cref{eq:sample-variance} is that the jackknife estimate \cref{eq:jack} is \emph{not} divided by $k-1$.
This difference may be justified intuitively.
The replicates $\mat{X}^{(j)}$ are \warn{not independent}, each one differs in a single input coordinate.
Since only one out of $k-1$ inputs differs between each pair of replicates, the replicates $\mat{X}^{(1)},\ldots,\mat{X}^{(k)}$ should have a ``$(k-1)^{-1}$ fraction'' of the variance of true iid copies $\mat{Z}^{(1)},\ldots,\mat{Z}^{(k)}\sim \mat{X}_{k-1}$.

The quality of jackknife variance estimates can be analyzed using the Efron--Stein--Steele inequality \cite{ES81,Ste86}:

\begin{fact}[Efron--Stein--Steele inequality]
    Let $\omega_1,\ldots,\omega_k$ be independent random variables in a measurable space $\set{S}$, let $f : \set{S}^k \to \field$ be a measurable function, and introduce iid copies $\omega_1',\ldots,\omega_k'$.
    Then
    \begin{equation*}
        \Var(Z) \le \frac{1}{2} \sum_{j=1}^k \expect \left|Z - Z^{(j)}\right|^2
    \end{equation*}
    where $Z \coloneqq f(\omega_1,\ldots,\omega_k)$ and $Z^{(j)} \coloneqq f(\omega_1,\ldots,\omega_{j-1},\omega_j',\omega_{j+1},\ldots,\omega_k)$.
\end{fact}

The Efron--Stein--Steele inequality is a powerful result that can be used to great effect in high-dimensional probability.
See \cite[\S2.1]{van14} and \cite[\S2.7]{Tro21} for accessible introductions and \cite{BLM13} for many uses of this inequality in probability theory.
It yields the following consequence for the jackknife estimator \cref{eq:jack}.

\begin{theorem}[Matrix jackknife overestimates variance] \label{thm:jackknife-overestimate}
    With the prevailing notation and assumptions, 
    \begin{equation*}
        \Var(\mat{X}_{k-1}) \le \expect \Jack^2(\mat{X}_{k-1}).
    \end{equation*}
\end{theorem}

This result states the matrix jackknife variance estimator \cref{eq:jack} overestimates the variance of the $(k-1)$-sample estimator \warn{on average}.
The fact that the jackknife estimate overestimates the true variance is a bit disappointing, but it is remarkable that the jackknife variance estimator possesses guarantees at this level of generality.
In practice, on the examples we consider in this thesis, the jackknife variance estimate tends to be accurate to within an order of magnitude, which is good enough to provide actionable information in applications.

\begin{proof}[Proof of \cref{thm:jackknife-overestimate}]
    First consider a scalar statistic $x_{k-1} = x_{k-1}(\omega_1,\ldots,\omega_{k-1})$ and introduce replicates $x^{(j)} = x_{k-1}(\omega_1,\ldots,\omega_{j-1},\omega_{j+1},\ldots,\omega_k)$ with mean $x^{(\cdot)}$.
    Observe that the replicates $x^{(j)}$ share the same distribution.
    Additionally,
    \begin{equation*}
        x_{k-1} = x_{k-1}(\omega_1,\ldots,\omega_{k-1}) = x^{(k)}.
    \end{equation*}
    By the Efron--Stein--Steele inequality,
    \begin{equation*}
        \Var(x_{k-1}) \le \frac{1}{2} \sum_{j=1}^{k-1} \expect \left| x_{k-1} - x_{k-1}(\omega_1,\ldots,\omega_{j-1},\omega_j',\omega_{j+1},\ldots,\omega_{k-1})\right|^2,
    \end{equation*}
    where $\omega_j'$ denotes an independent copy of $\omega_j$.
    The $\omega_j$ are iid random variables, so wecan replace $\omega_j'$ by $\omega_k$.
    Thus, we obtain
    \begin{equation*}
        \Var(x_{k-1}) \le \frac{1}{2} \sum_{j=1}^{k-1} \expect \left| x^{(k)} - x^{(j)}\right|^2.
    \end{equation*}
    Since each element of $\{ x^{(\ell)} - x^{(j)} : \ell \ne j \}$ has the same distribution, we can replace $x^{(k)}$ by $x^{(\ell)}$ and average over the index $\ell$, obtaining
    \begin{equation*}
        \Var(x_{k-1}) \le \frac{1}{2k} \sum_{j,\ell=1}^k \expect \left| x^{(\ell)} - x^{(j)}\right|^2.
    \end{equation*}
    Next, we compute
    \begin{align*}
        \expect \sum_{j=1}^k \left| x^{(j)} - x^{(\cdot)} \right|^2 &= \expect\sum_{j=1}^k \left| x^{(j)} - \frac{1}{k}\sum_{\ell=1}^k x^{(\ell)} \right|^2 \\
        &= \expect \sum_{j=1}^k \left[\left| x^{(j)}\right|^2 -\frac{2}{k}\sum_{\ell=1}^k \Re\left( \overline{x^{(j)}} x^{(\ell)} \right) - \frac{1}{k^2 }\sum_{\ell,p=1}^k \Re\left(\overline{x^{(\ell)}}x^{(k)}\right)\right] \\
        &= \expect \left[\frac{1}{2}\sum_{j=1}^k \left| x^{(j)}\right|^2 + \frac{1}{2}\sum_{\ell=1}^k \left| x^{(j)}\right|^2 - \frac{1}{k} \sum_{j,\ell=1}^k \Re\left( \overline{x^{(j)}} x^{(\ell)} \right)\right] \\
        &= \frac{1}{2k} \sum_{j,\ell=1}^k \expect \left| x^{(\ell)} - x^{(j)}\right|^2 = \Jack^2(x_{k-1}).
    \end{align*}
    In the second line, we expand the square.
    In the third line we consolidate the second and third terms, break the sum $\sum_{j=1}^k \left| \smash{x^{(j)}}\right|^2$ into two equal pieces, and reindex by replacing $j$ with $\ell$.
    In the last line, we recombine.
    Joining the two previous displays leads to the the desired result
    \begin{equation} \label{eq:jacknife-overestimate-scalar}
        \Var(x_{k-1}) \le \Jack^2(x_{k-1}).
    \end{equation}
    The matrix case follows by decomposing the matrix variance as a sum of the variance of its entries \cref{eq:matrix-variance-decomposition} and invoking \cref{eq:jacknife-overestimate-scalar} entry-by-entry.
\end{proof}

\begin{remark}[What about the bias?]
    Jackknife estimates of the bias were developed by Quenouille in 1949 \cite{Que49}, predating jackknife variance estimates by nearly a decade.
    Quenouille's bias estimate can be formed for a randomized matrix approximation, but its usefulness in this setting is unclear.
    Quenouille estimate is typically analyzed under the assumption that the bias can be expanded in reciprocal powers of $k$.
    Randomized low-rank approximation algorithms seemingly do not satisfy this property.
    Thus, it is unclear what, if any, insight Quenouille's bias estimate provides for randomized matrix algorithms. 
\end{remark}

\section{Using matrix jackknife variance estimation} \label{sec:using-matrix-jackknife}

In principle, the use cases for matrix jackknife variance estimation could be very broad. 
In this thesis, we will apply jackknife variance estimation to quantities $\mat{X}_k = \mat{Q}(\Ahat)$ computed from a randomized Nystr\"om approximation $\Ahat \approx \mat{A}$; see \cite{ET24} for further examples.
To instantiate the jackknife methodology, we treat 
\begin{equation*}
    \mat{X}_k = \mat{X}_k(\vec{\omega}_1,\ldots,\vec{\omega}_k)
\end{equation*}
as a function of the iid columns $\vec{\omega}_i$ of the random test matrix $\mat{\Omega}$.
We will focus on this setting for the remainder of the chapter.

A limitation of the matrix jackknife variance estimator is that it captures just one term in the bias--variance decomposition \cref{eq:bias-variance}.
To approximate the mean-squared error, a separate estimate of the bias is necessary.
However, the information provided by the jackknife variance estimate is still actionable.
\actionbox{If the variance is high, the mean-squared error is high, and the output of the algorithm should not be trusted.
If the variance is low, then the jackknife variance estimate is equivocal, as the mean-squared error could be low or high.}

For a user, what does it mean if the variance is high?
There are at least two possibilities:
\begin{enumerate}
    \item \textbf{Too few samples.} The approximation $\Ahat$ is under-resolved.
    More samples (i.e., more columns in the test matrix $\mat{\Omega}$) are needed to produce an approximation $\Ahat$ of sufficiently high quality to estimate the quantity of interest.
    \item \textbf{Bad conditioning.} The quantity of interest $\mat{Q}(\mat{A})$ could be \emph{poorly conditioned} in the sense that small changes to the input $\mat{A}$ result in large changes to $\mat{Q}(\mat{A})$. 
    In this case, even if the approximation $\mat{A} \approx \Ahat$ is highly accurate, $\mat{Q}(\Ahat)$ can be far from $\mat{Q}(\mat{A})$ owing to the inherent sensitivity of the quantity of interest.
    When dealing with a badly conditioned problems, adding additional samples may not help improve the quality of the approximation $\mat{Q}(\Ahat)\approx \mat{Q}(\mat{A})$.
\end{enumerate}
In either case, high variance---diagnosed by a jackknife estimate---flashes a signal not to trust the computed output.

Jackknife variance estimation can be deployed at runtime to detect these bad behaviors.
One possible use case is to run the jackknife variance estimate and provide a warning to the user if high variance is detected.
As an alternative, jackknife variance estimation can be used to adaptively determine the approximation rank $k$ until the variance drops below a specified threshold.

The main advantage of jackknife variance estimation is it provides a way of estimating the quality of a very general class of quantities produced by randomized matrix approximation algorithms.
Like the leave-one-out\index{leave-one-out randomized algorithm} error estimate, it requires \emph{no additional matrix--vector products beyond those used to run the original algorithm}.
To achieve the maximum efficiency for the jackknife variance estimator, one can develop purpose-build fast implementations using the randomized SVD or the randomized Nystr\"om downdating results (\cref{thm:rsvd-downdate,thm:nystrom-downdate}).
The remainder of this chapter develops an example of application of the matrix jackknife for assessing the quality of ``spectral transformations'' of Nystr\"om approximations; see \cite{ET24} for several more examples.

\section{Example: Spectral transformations of Nystr\"om approximations} \label{sec:spectral-downdating}

Many of the ways we \emph{use} Nystr\"om approximations to psd matrices can be encompassed under the umbrella of applying a \emph{spectral transformation} to the matrix.
We make the following definition:

\begin{definition}[Spectral transformation]
    Let $\vec{f}$ denote a vector-valued function $\vec{f} : \real_+^n \to \complex^n$, and let $\mat{A} \in \field^{n\times n}$ denote a psd matrix with eigendecomposition $\mat{A} = \mat{U}\Diag(\vec{\lambda}(\mat{A}))\mat{U}^*$.
    The \emph{spectral transformation} of $\mat{A}$ by $\vec{f}$ is
    \begin{equation*}
        \vec{f}[\mat{A}] \coloneqq \mat{U}\Diag(\vec{f}(\vec{\lambda}(\mat{A})))\mat{U}^*.
    \end{equation*}
    We remind the reader that the eigenvalues $\vec{\lambda}(\mat{A})$ are sorted in \warn{nonincreasing order}.
    The spectral transformation may fail to be uniquely defined if $\mat{A}$ has repeated eigenvalues.
    
    The spectral transformation $\vec{f}$ is \emph{rank-$k$} if 
    \begin{equation*}
        \vec{d}(k+1:n) = \vec{0} \implies [\vec{f}(\vec{d})](k+1:n) \quad \text{for every } \vec{d} \in \real_+^n.
    \end{equation*}
    For a rank-$k$ spectral transformation, we can extend its domain of definition to matrices of size $k$ by
    \begin{equation} \label{eq:spectral-transformation-rank-k-extension}
        \vec{f}(d_1,\ldots,d_k) \coloneqq [\vec{f}(d_1,\ldots,d_k,0\ldots,0)](1:k).
    \end{equation}
\end{definition}

Many of the application of Nystr\"om approximations can be reinterpreted as spectral transformations.
Here are several examples:

\begin{enumerate}
    \item \textbf{\textit{Spectral projectors:}} One of the main uses for Nystr\"om approximation is to approximate eigenvectors.
    However, even when the eigenvalues are distinct, the eigenvectors are only defined up to a sign (or phase, if $\field = \complex$); when there are repeated eigenvalues, the non-uniqueness of eigenvectors become more severe.
    To resolve many of these uniqueness issues, we can use \emph{spectral projectors}. 
    Given a subset $\set{S} \subseteq \{1,\ldots,n\}$ and a matrix $\mat{A} = \mat{U}\Diag(\vec{\lambda}(\mat{A}))\mat{U}^*$, the spectral projector associated with the index set $\set{S}$ is
    \begin{equation*}
        \mat{\Pi}_{\set{S}}(\mat{A}) = \mat{U}(:,\set{S})\mat{U}(:,\set{S})^*.
    \end{equation*}
    Natural examples include $\set{S} = \{1,\ldots,k\}$ (top-$k$ subspace) or $\set{S} = \{j\}$ (single eigenvector).
    The spectral projector is the orthprojector on the \emph{invariant subspace} spanned by the eigenvectors $\operatorname{span} \{ \vec{u}_s : s \in \set{S} \}$ indexed by $\set{S}$.
    Spectral projectors are spectral transformation associated with the constant function $\vec{f}(\vec{d}) = \sum_{s \in \set{S}} \evec_s$.
    As a spectral transformation, the spectral projector function $\vec{f}$ is rank-$k$ if and only if $\set{S} \subseteq \{1,\ldots,k\}$.

    \item \textbf{\textit{Rank-$r$ truncation.}}
    We can truncate rank-$k$ Nystr\"om approximations to smaller rank $r < k$ to reduce storage.
    The rank-$r$ approximation operation $\mat{A} \mapsto \lowrank{\mat{A}}_r$ is also a spectral transformation, associated with the function
    \begin{equation*}
        \vec{f}(d_1,\ldots,d_n) = (d_1,\ldots,d_r,0,\ldots,0) \quad \text{for every } \vec{d} \in \real_+^n.
    \end{equation*}
    This spectral transformation is rank-$k$ for every $1\le k\le n$.

    \item \textbf{\textit{Matrix functions.}}
    Another use of Nystr\"om approximations, developed by David Persson and collaborators \cite{PK23,PMM25}, is to approximate a (standard) matrix functions $f(\mat{A})$.
    The idea is straightforward: Given a Nystr\"om approximation $\Ahat \approx \mat{A}$, call $f(\Ahat)$ as the \emph{funNystr\"om} approximation to $f(\mat{A})$.
    The funNystr\"om approximation is the spectral transformation associated with the function
    \begin{equation*}
        \vec{f}(d_1,\ldots,d_n) = (f(d_1),\ldots,f(d_n)) \quad \text{for all } \vec{d} \in \real^n_+.
    \end{equation*}
    For every $1\le k\le n$, this spectral transformation is rank-$k$ if and only if $f(0) = 0$.
    (The requirement that $f(0) = 0$ can always be satisfied by shifting $\tilde{f} \coloneqq f - f(0)$.)
\end{enumerate}

Given these various settings, it is natural to seek ways of estimating the variance of spectral transformations of a Nystr\"om approximation.

For the rest of this section, we consider the task of estimating a general spectral transformation $\vec{f}[\Ahat]$ using a single-pass Nystr\"om approximation $\Ahat = \mat{A}\langle \mat{\Omega}\rangle$ to a psd matrix $\mat{A}$.
We assume the test matrix $\mat{\Omega}\in \field^{n\times k}$ is composed of \emph{independent} columns, and we assume that $\vec{f}$ has rank $k$.

\subsection{Variance estimation for spectral transformations: $\order(k^4)$ algorithm}

The simplest approach to variance estimation for spectral transformations passes through a variant of the downdating formula \cref{eq:nystrom-downdate-practice} for Nystr\"om approximations.

Begin by computing the Nystr\"om approximation as in \cref{prog:nystrom_shiftcor}:
We form $\mat{Y} \coloneqq \mat{A}\mat{\Omega}$, compute $\mat{H} \coloneqq \mat{\Omega}^*\mat{Y}$, factorize $\mat{H} = \mat{R}^*\mat{R}$, define $\mat{F} \coloneqq \mat{Y}\mat{R}^{-1}$, and compute an SVD $\mat{F} = \mat{U}\mat{\Sigma}\mat{V}^*$.
The Nystr\"om approximation is now given as 
\begin{equation*}
    \Ahat = \mat{U}\mat{D}\mat{U}^* \quad \text{where } \mat{D} = \mat{\Sigma}^2.
\end{equation*}
(For numerically stable implementation in practice, one should add shift correction to this procedure; see \cref{eq:stable-nystrom}.)
Following \cref{eq:nystrom-downdate}, the family of downdated Nystr\"om approximations $\Ahat_{(i)} = \mat{A}\langle \mat{\Omega}_{-i}\rangle$ is 
\begin{equation*}
    \Ahat_{(i)} = \mat{U}\mat{D}\mat{U}^* - \outprod{\vec{z}_i} \quad \text{where } \mat{Z} \coloneqq \mat{F}\mat{R}^{-*} \cdot \Diag\big(\srn(\mat{R}^{-1})\big)^{-1/2}.
\end{equation*}
Using the SVD $\mat{F} = \mat{U}\mat{\Sigma}\mat{V}^*$, we can rewrite the downdating rule as
\begin{equation*}
    \Ahat_{(i)} = \mat{U}(\mat{D} - \outprod{\vec{w}_i})\mat{U}^* \quad \text{where } \mat{W} = \mat{\Sigma}\mat{V}^*\mat{R}^{-*}\cdot \Diag\big(\srn(\mat{R}^{-1})\big)^{-1/2}.
\end{equation*}
Consequently, the entire family of Nystr\"om approximations is given by the $\mat{W}$ matrix, which can be computed in $\order(k^3)$ operations.

Using the downdating matrix $\mat{W}$ in a straightforward way, we can compute the jackknife variance estimate in $\order(k^4)$ operations.
The spectral transformation of each $\Ahat_{(i)}$ is given as
\begin{equation*}
    \vec{f}[\Ahat_{(i)}] = \mat{U}\vec{f}[\mat{D} - \outprod{\vec{w}_i}]\mat{U}^*,
\end{equation*}
where we have used the extension \cref{eq:spectral-transformation-rank-k-extension} of spectral functions to $k\times k$ matrices.
Thus, we have
\begin{equation*}
    \Jack^2(\vec{f}[\Ahat]) = \sum_{j=1}^k \norm{ \vec{f}[\Ahat_{(j)}] - \frac{1}{k} \sum_{i=1}^k \vec{f}[\Ahat_{(i)}] }_{\mathrm{F}}^2 = \sum_{j=1}^k \norm{ \mat{X}^{(j)} - \mat{X}^{(\cdot)} }_{\mathrm{F}}^2,
\end{equation*}
where
\begin{equation*}
    \mat{X}^{(j)} \coloneqq \vec{f}[\mat{D} - \outprod{\vec{w}_j}] \quad \text{and} \quad \mat{X}^{(\cdot)} \coloneqq \frac{1}{k} \sum_{j=1}^k \mat{X}^{(j)}.
\end{equation*}
Computing each $\mat{X}^{(j)}$ directly requires a spectral decomposition of $\mat{X}^{(j)}$ at $\order(k^3)$ cost.
Since there are $k$ replicates $\mat{X}^{(j)}$, the total cost is $\order(k^4)$ operations.
Once the replicates have been formed, the jackknife variance estimate requires an additional $\order(k^3)$ arithmetic operations.
Code is provided in \cref{prog:nystrom_jack}.

\myprogram{Single-pass Nystr\"om approximation with jackknife variance estimation for spectral transformation.}{Subroutine \texttt{nystrom} is provided in \cref{prog:nystrom}.}{nystrom_jack}

\subsection{Acceleration to $\order(k^3)$ operations}

In some cases, variance estimation for spectral transformations can be accelerated to $\order(k^3)$ operations.
Specifically, consider a spectral transformation for which $\vec{f}(\vec{d})$ has only $r$ nonzero entries for every $\vec{d} \in \real_+^n$. 
To compute $\vec{f}[\mat{D} - \outprod{\vec{w}_j}]$ efficiently, observe that $\mat{D} - \outprod{\vec{w}_j}$ is a rank-one modification to a diagonal matrix $\mat{D}$.
The spectral decomposition of each $\mat{D} - \outprod{\vec{w}_j} = \mat{Q}_{(j)}^{\vphantom{*}}\Diag(\vec{d}_{(j)})\mat{Q}_{(j)}^*$ can be computed in $\order(k^2)$ operations; see \cite[\S5.3.3]{Dem97}.
Then, each $\mat{X}^{(j)} = \vec{f}[\mat{D} - \outprod{\vec{w}_j}]$ can be computed in $\order(k^2r)$ operations.
Thus, the total runtime using this approach is $\order(k^3r)$ operations.
If $r = \order(1)$, the total cost is $\order(k^3)$ operations, as promised.
See the code of \cite{ET24} for an implementation of this faster algorithm.

\subsection{Experiment}

To demonstrate the effectiveness of the jackknife variance estimation technique, we use it to assess the variance of spectral transformations of Nystr\"om approximations to a psd matrix $\mat{A}$.
We set $\mat{A}$ to be the \texttt{exp} matrix \cref{eq:exp}, and we consider two spectral transformations: the dominant rank-$10$ projector $\vec{f}[\mat{A}] \coloneqq \mat{\Pi}_{1:10}(\mat{A})$ and the matrix square root $\vec{f}[\mat{A}] \coloneqq \mat{A}^{1/2}$.
We test approximation ranks $25\le k \le 150$.

\begin{figure}
    \centering
    \includegraphics[width=0.99\linewidth]{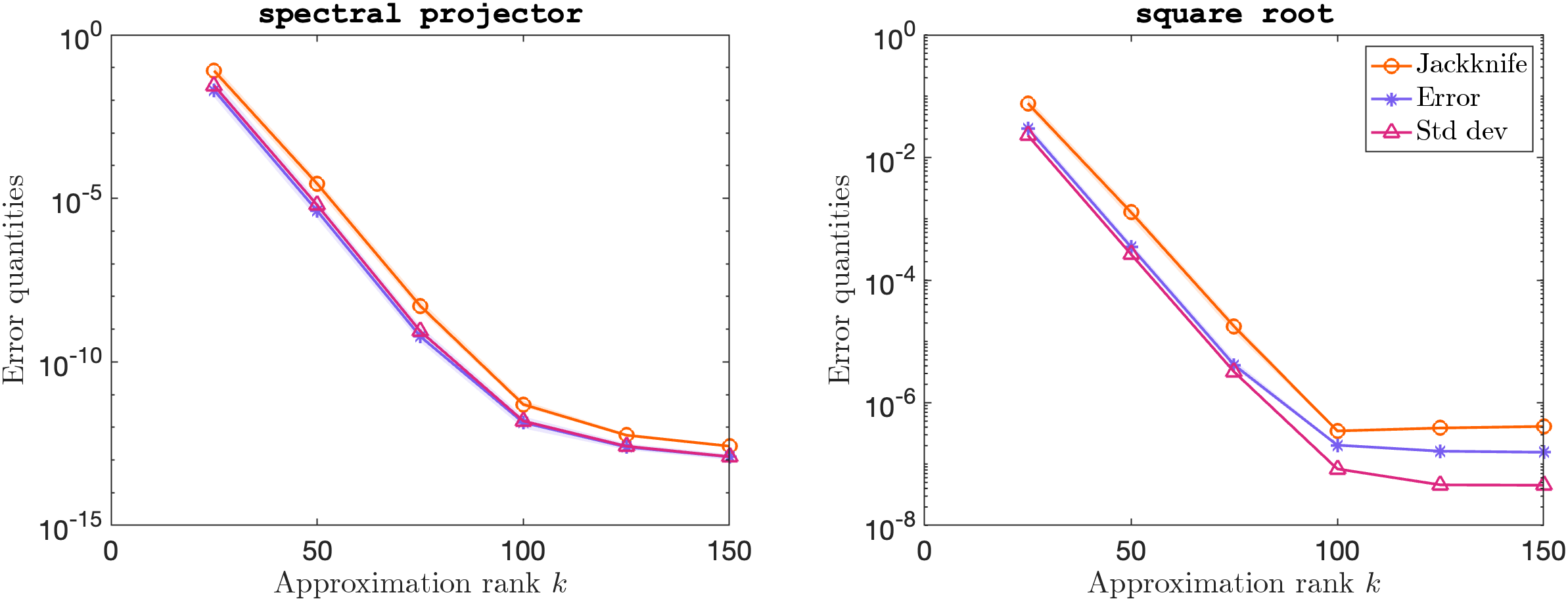}
    \caption[Jackknife variance estimate for spectral transformations of Nystr\"om approximations]{Jackknife standard deviation estimate (orange circles), error (purple asterisks), and standard deviation (pink triangles) for spectral transformations of Nystr\"om approximations of ranks $25\le k\le 150$.
    Spectral transformations are the dominant rank-$10$ projector $\vec{f}[\mat{A}] \coloneqq \mat{\Pi}_{1:10}(\mat{A})$ (\emph{left}) and the matrix square root $\vec{f}[\mat{A}] \coloneqq \mat{A}^{1/2}$ (\emph{right}).}
    \label{fig:jackknife}
\end{figure}

Results appear in \cref{fig:jackknife}.
For each spectral transformation and each Nystr\"om approximation, we show the error $\norm{\vec{f}[\mat{A}] - \vec{f}[\Ahat]}_{\mathrm{F}}$, the standard deviation $\Var(\vec{f}[\Ahat])^{1/2}$, and the jackknife standard deviation estimate $\Jack(\vec{f}([\Ahat])$.
For both examples, we see that the jackknife variance estimate serves as a modest overestimate of both the error and standard deviation, providing a useful diagnostic of the quality of the computed approximation.
Many more examples of the jackknife methodology are provided in \cite[\S\S1.2 \& 5]{ET24}.

\chapter{Leave-one-out randomized matrix algorithms: Open problems}

This chapter discusses open problems related to leave-one-out\index{leave-one-out randomized algorithm} randomized matrix algorithms including error analysis (\cref{sec:loo-error-analysis}), leave-one-out\index{leave-one-out randomized algorithm} algorithms for Hermitian indefinite matrices (\cref{sec:hermitian-indefinite}), and numerically stable downdating formulas for subspace iteration (\cref{sec:open-rsi-downdating}).

\section{Open problem: Error analysis} \label{sec:loo-error-analysis}

This part of the thesis has discussed matrix attribute estimation algorithms.
The basic versions of \XTrace and \XNysTrace have \emph{a priori} bounds on the mean-squared error, and most of the remaining algorithms have no error analysis at all.
In my view, theoretical analysis of these algorithms is not a major limitation; they are based on sound principles---unbiased Monte Carlo approximation, variance reduction by low-rank approximation, and exchangeable leave-one-out\index{leave-one-out randomized algorithm} design---and they achieve \emph{spectral accuracy} in practice.
Still, it is natural to desire for a more comprehensive mathematical analysis of these algorithms.

There are several natural topics for future research:
\begin{itemize}
    \item Develop sharp error bounds for \XTrace and \XNysTrace that hold with high-probability.
    \item Obtain \emph{a priori} error bounds, either high-probability or on the mean-squared error, for the randomized diagonal estimators and row-norm estimators from \cref{ch:diagonal,ch:row-norm}.
    \item Obtain bounds on the variance of the leave-one-error estimate from \cref{ch:loo}.
\end{itemize}

\section{Open problem: Hermitian indefinite matrices} \label{sec:hermitian-indefinite}

This thesis has developed leave-one-out\index{leave-one-out randomized algorithm} matrix attribute estimation algorithms for general matrices using the randomized SVD (e.g., \XTrace) and algorithms for psd matrices using single-pass Nystr\"om approximation (e.g., \XNysTrace).
Can optimized leave-one-out\index{leave-one-out randomized algorithm} algorithms be devised for matrices that have additional structure, but are not psd?
In particular, is there a way to improve leave-one-out\index{leave-one-out randomized algorithm} matrix algorithms for Hermitian \emph{indefinite} matrices?

\myprogram{Implementation of the \XNysTrace algorithm designed to work with Hermitian indefinite matrices.}{Subroutines \texttt{diagprod} and \texttt{random\_signs} appear in \cref{prog:diagprod} and \cref{prog:random_signs}.}{xsymtrace}

\begin{figure}
    \centering
    \includegraphics[width=0.95\linewidth]{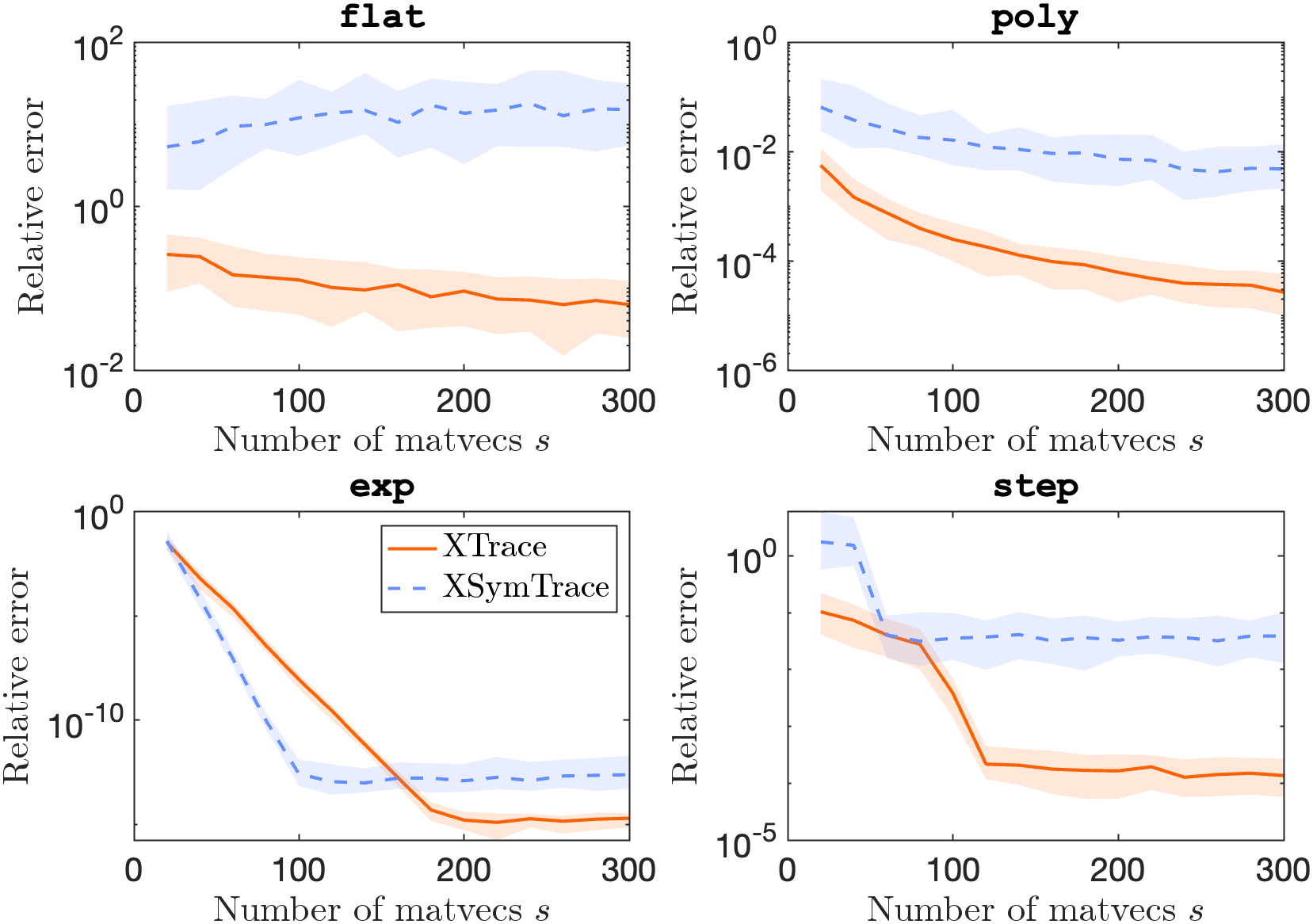}
    \caption[Comparison of \XTrace and bad \XSymTrace algorithm for estimating the trace of matrices with different spectra]{Comparison of \XTrace (orange solid) and \XSymTrace (blue dashed) on four test matrices \cref{eq:test-matrices} with different spectral characteristics.
    Lines show median of 100 trials, and shaded regions show 10\% and 90\% quantiles.
    }
    \label{fig:xsymtrace}
\end{figure}

As we discussed in \cref{rem:indefinite}, one \emph{can} apply single-pass Nystr\"om approximation to Hermitian indefinite matrices without modification, although it can be inaccurate on certain instances.
\Cref{prog:xsymtrace} provides an implementation of \XSymTrace, a version of \XNysTrace designed to work with Hermitian indefinite matrices.
Results for this algorithm on our testbed of synthetic matrices with different spectral characteristics appear in \cref{fig:xsymtrace}.
We use the same test matrices from previous sections, except that we randomize the sign of each eigenvalue so the matrices are indefinite.
We see that the \XSymTrace estimator is substantially less accurate than \XTrace on three of four examples, demonstrating the failure of (single-pass) Nystr\"om approximation for indefinite matrices.
Thus, we do not regard \XSymTrace as an effective general-purpose algorithm for trace estimation of Hermitian indefinite matrices.
(Nevertheless, its performance may be acceptable for matrices with only a \emph{small} number of negative eigenvalues or with very rapid spectral decay.)

For this reason, determining the ``right'' approach to developing leave-one-out\index{leave-one-out randomized algorithm} algorithms for Hermitian indefinite matrices remains an open problem.
The oversampling technique discussed in \cref{rem:indefinite} may be a promising strategy, although downdating the regularized pseudoinverses  $\lowrank{\mat{\Omega}^*\mat{A}\mat{\Omega}}_k^\dagger$ may be a nontrivial task.
(The methods from \cref{sec:spectral-downdating} could prove helpful.)
For the oversampled Nystr\"om method to lead to substantial speedups over \XTrace, it is important to keep the oversampling factor as small as possible, potentially using smaller oversampling factors than considered by Park and Nakatsukasa \cite{PN24}.
I also believe there could be new approaches to low-rank approximation of Hermitian indefinite matrices that could be more natural and powerful than the truncated Nystr\"om method.

\section{Open problem: Numerically stable downdating for subspace iteration} \label{sec:open-rsi-downdating}

In \cref{sec:loo-rsvd}, we developed an efficient implementation of the leave-one-out\index{leave-one-out randomized algorithm} error estimator for the randomized SVD (without subspace iteration).
This method immediately yields error estimates for randomized subspace iteration (\cref{sec:rsi}); see \cref{prog:rsi_errest} for an implementation.

\myprogram{Randomized subspace iteration for producing a low-rank approximation to a general matrix with leave-one-out error estimation.}{Subroutines \texttt{cnormc} and \texttt{diagprod} are provided in \cref{prog:cnormc} and \cref{prog:diagprod}.}{rsi_errest}

\begin{figure}
    \centering
    \includegraphics[width=0.7\linewidth]{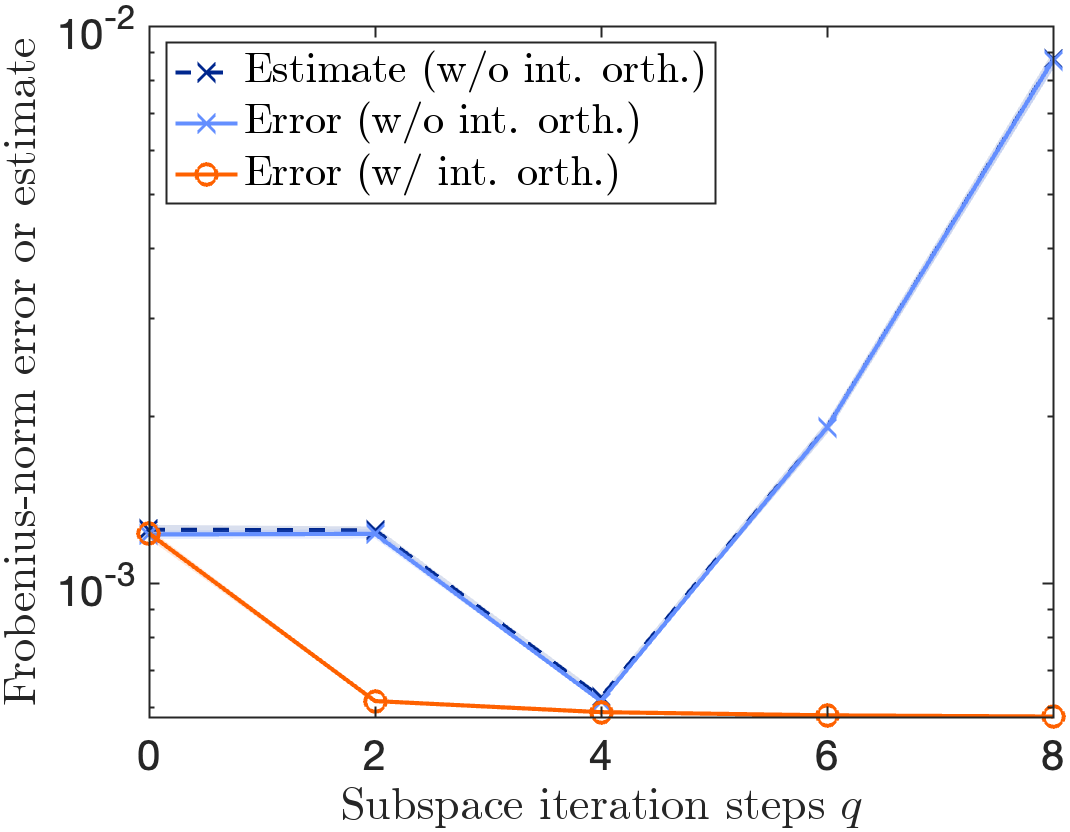}
    \caption[Error and leave-one-out error estimate for randomized subspace iteration with and without intermediate reorthogonalization]{Error for randomized SVD with subspace iteration with (solid orange circles) and without (solid light blue crosses) intermediate reorthogonalization for rank-100 approximation of the \texttt{poly} matrix \cref{eq:poly} as a function of the subspace iteration count $q$.
    leave-one-out\index{leave-one-out randomized algorithm} error estimate (without intermediate reorthogonalization) is shown as dashed dark blue crosses.
    Solid lines trace the median of 100 trials, and shaded regions show 10\% and 90\% quantiles.}
    \label{fig:rsi-loo}
\end{figure}

\Cref{fig:rsi-loo} charts performance of randomized subspace iteration and its error estimate as function of the subspace iteration count $q$.
We see that the leave-one-out\index{leave-one-out randomized algorithm} error estimator closely tracks the true value of the error across all values of $q$.
However, for $q \ge 6$, the error (and error estimates) for randomized subspace iteration actually \emph{increase} with the number of subspace iteration steps $q$.
This behavior arise from numerical issues.
As shown in \cref{fig:rsi-loo}, these numerical issues are easily cured by using intermediate orthogonalization, discussed in \cref{sec:rsi}.

This brings us to the open problem.
At present, it is not known whether there exists a \warn{numerically stable} way of downdating the randomized approximation computed via orthogonalized subspace iteration.
As such, the leave-one-out\index{leave-one-out randomized algorithm} error estimation technique is currently limited to randomized subspace iteration without intermediate reorthogonalization.
The open problem is to remove this limitation by extending leave-one-out\index{leave-one-out randomized algorithm} randomized matrix algorithms to be compatible with subspace iteration.
There could be other benefits to stable downdating approaches for randomized subspace iteration, including numerically stable versions of \XNysTrace incorporating subspace iteration.

\partepigraph{Dedicated to my grandparents Bob and Sarah Epperly, Kathie and Bill Quarles, and John and Christine Weidner.}
\part{Sketching, solvers, and stability}
\label{part:sketching}

\chapter{Algorithms for least squares, a brief history} \label{ch:ls-algs-history}

\epigraph{But recent research led to new uses of randomization: random mixing and random sampling, which can be combined to form random projections. These ideas have been explored theoretically and have found use in some specialized applications (e.g., data mining), but they have had little influence so far on mainstream numerical linear algebra. Our paper answers a simple question, Can these new techniques beat state-of-the-art numerical linear algebra libraries \emph{in practice}? Through careful engineering of a new least-squares solver, which we call Blendenpik, and through extensive analysis and experimentation, we have been able to answer this question, yes.}{Haim Avron, Petar Maymounkov, and Sivan Toldeo, \textit{Blendenpik: Supercharging LAPACK's least-squares solver} \cite[\S1]{AMT10}}

The third part of this thesis will investigate randomized algorithms for linear least-squares fitting and related problems.
Our core algorithmic tool will be (linear) \emph{randomized dimensionality reduction}, also known as \emph{sketching}.
This introductory chapter will set the stage, introducing the linear least-squares problem and algorithms for solving it from a historical perspective.
Our historical journey will be brief emphasizing topics that will inform the technical development in subsequent chapters.

\myparagraph{Sources}
This is an introductory chapter, not based on any particular research article.
The reference \cite{BMR22} provides a definitive history of numerical linear algebra from antiquity to the first years of the twenty-first century.
Our tale will hew closely to their account for the history of deterministic methods for least squares.
In the past twenty-five years, randomized methods have greatly expanded the class of algorithms for least squares.
My understanding of the history of these methods is informed by my reading of the original literature, surveys \cite{Woo14a,MT20a,MDM+23a}, and discussions with the researchers involved.

\myparagraph{Outline}
We begin our tour in \cref{sec:least-squares-intro} with the origins of linear least squares and its solution by the normal equations.
We encounter the numerical instabilities of this approach in \cref{sec:normal-equations-unstable} and Golub's \QR factorization as a stable alternative in \cref{sec:qr}.
But how accurate is the \QR factorization method?
This question is answered in \cref{sec:perturbation-theory}, which reviews perturbation theory, and \cref{sec:forward-backward-stability}, which describes forward and backward stability.
For the largest problems, direct methods like \QR factorization are prohibitively slow, and we must turn to iterative methods, discussed in \cref{sec:ls-iterative-methods}.
Randomization enters in \cref{sec:sketch-and-solve}, which describes the sketch-and-solve method as a quick way to obtain an approximate least-squares solution.
To achieve higher accuracy, one can use randomization to precondition an iterative method, resulting in the sketch-and-precondition method (\cref{sec:spre}).
This chapter concludes in \cref{sec:spre-unstable?} with questions about the numerical stability of sketch-and-precondition raised in the recent paper \cite{MNTW24}.
These concerns will be resolved in \cref{ch:stable-ls}, which describes fast, numerically stable randomized least-squares solvers.

\section{The overdetermined linear least-squares problem} \label{sec:least-squares-intro}

The overdetermined linear least-squares problem
\begin{equation} \label{eq:least-squares}
    \vec{x} = \argmin_{\vec{z} \in \field^n} \norm{\vec{c} - \mat{B}\vec{z}}
\end{equation}
is ubiquitous in modern science, engineering, mathematics, and computer science.
Throughout this part of the thesis, $\mat{B} \in \field^{m\times n}$ denotes a matrix of coefficients, $\vec{c} \in \field^m$ is the right-hand side, and $\vec{x} = \mat{B}^\dagger \vec{c} \in\field^n$ denotes \warn{optimal} least-squares solution (i.e., the solution to \cref{eq:least-squares}).
Approximate minimizers will be called $\vec{z}$ or $\hatvector{x}$.
We will typically assume that $\mat{B}$ is full-rank, so that \cref{eq:least-squares} has a unique solution.

Legendre published the first description of the method of least squares in 1805 \cite[p.~57]{BMR22}.
Legendre was motivated by problems in geodesy, and he used the method of least squares to estimate the length of the Paris meridian between Dunkirk and Barcelona.
In this early work, Legendre took a ``bare hands'' approach, writing down an expression for the total squared deviation and finding a linear system for the least-squares system by computing the partial derivatives and setting them to zero.

We can repeat Legendre's approach, aided by modern matrix notation (for now, focusing on the \warn{real} case $\field = \real$.)
The gradient of the least-squares objective function \cref{eq:least-squares} is
\begin{equation*}
    \vec{\nabla}_{\vec{z}} \norm{\vec{c} - \mat{B}\vec{z}}^2 = 2\mat{B}^*(\vec{c} - \mat{B}\vec{z}).
\end{equation*}
Setting this derivative to zero, and rearranging, we obtain the \emph{normal equations} for the least-squares problem:
\begin{equation} \label{eq:normal-equations}
    (\mat{B}^*\mat{B}) \vec{x} = \mat{B}^*\vec{c}.
\end{equation}
Provided that $\mat{B}$ has full rank, the normal equations comprise a square, consistent system of linear equations, and the unique least-squares solution is $\vec{x} = (\mat{B}^*\mat{B})^{-1}(\mat{B}^*\vec{c})$.
A similar computation with Wirtinger derivatives establishes the normal equations in the complex case ($\field = \complex$).

Legendre's claim to the discovery of the method of least squares was disputed by Gauss, who first published on least squares in 1811.
(He claimed to have discovered the method in unpublished work in 1795 \cite[pp.~58--60]{BMR22}..) 
Gauss was motivated by problems in astronomy, and he successfully used the method of least squares to localize the orbit of Ceres.
Gauss's work expanded on Legendre's by establishing a connection between least squares, statistics, and the ``normal'' distribution that now bears Gauss's name.
Unfortunately, the disagreement between Legendre and Gauss was not cordial.
To quote from \cite[p.~60]{BMR22},

\begin{quote}
    [Gauss] claimed that he already used his method [of least squares] (Unser Prinzip) in 1795 but that its publication was delayed by the Napoleonic wars.
    Although his name was mentioned by Gauss, Legendre was badly offended.
    In a letter to Gauss dated May 31, 1809, he rightfully stated that priority is only established by publication.
    Gauss did not answer him.
\end{quote}

\section{Numerical instabilities of the normal equations} \label{sec:normal-equations-unstable}

As the nineteenth century became the twentieth, the generally accepted practice for solving a linear least-squares problems were to solve ``Gauss' normal equations'' \cite[p.~17]{van74} by the ``Choleski square-root method'' \cite[p.~339]{Hou58}.
In the 1960s, issues of \emph{numerical stability} associated with this approach became wildly known; Golub's paper \cite[\S1]{Gol65} provides a crisp example (see also \cite{Laeu61}).

\myparagraphnp{How accurately can we solve systems of linear equations?}
In order to understand the numerical instabilities associated with solving the normal equations, we must first understand how accurately we can solve linear systems of equations.
We provide a brief review of numerical stability for linear systems, a subject pioneered by von Neumann and Goldstine \cite{vG47} and Turing \cite{Tur48}.

To solve a system of linear equations
\begin{equation} \label{eq:example-linsys}
    \mat{M} \vec{x} = \vec{f}
\end{equation}
on a digital computer requires the real or complex entries of the matrix $\mat{M}$ and the vector $\vec{f}$ to be stored using finite-precision representations.
Most modern computing platforms use \emph{floating-point representations} of numbers; see \cite[\S2]{Hig02} or \cite{Ove01} for an introduction to floating-point numbers.

The accuracy of a system of floating-point numbers is characterized by the \emph{unit roundoff} $u$, which captures the magnitude of errors in representing real numbers using floating-point representations and the scale of \emph{rounding errors} incurred when arithmetic operations are performed.
The unit roundoff is $u\approx 10^{-16}$ in double-precision arithmetic and $u\approx 10^{-8}$ in single precision.
The precision of a floating-point number system is also sometimes characterized using the \emph{machine epsilon}, which is twice the unit roundoff $u$.

Even \emph{storing} the matrix $\mat{M}$ incurs errors on the order of the unit roundoff:
\begin{equation*}
    \mat{M}_{\mathrm{stored}} = \mat{M} + \mat{E} \quad \text{for } \mat{E} \lessapprox \norm{\mat{M}} u.
\end{equation*}
Solving the linear system \cref{eq:example-linsys} with the stored matrix $\mat{M}_{\mathrm{stored}}$ results in an approximate solution
\begin{equation*}
    \hatvector{x} = \mat{M}_{\mathrm{stored}}^{-1}\vec{f} \approx \mat{M}^{-1}\vec{f} - \mat{M}^{-1} \mat{E} \mat{M}^{-1}\vec{f}.
\end{equation*}
Here, we expanded the inverse $(\mat{M} + \mat{E})^{-1} \approx \mat{M}^{-1} - \mat{M}^{-1}\mat{E}\mat{M}^{-1}$ to first order in the perturbation $\mat{E}$.
Therefore, the (relative) \emph{forward error} is roughly of size
\begin{equation*}
    \frac{\norm{\hatvector{x} - \mat{M}^{-1}\vec{f}}}{\norm{\mat{M}^{-1}\vec{f}}} \approx \frac{\norm{\mat{M}^{-1} \mat{E} \mat{M}^{-1}\vec{f}}}{\norm{\mat{M}^{-1}\vec{f}}} \le \norm{\smash{\mat{M}^{-1}}} \cdot \norm{\mat{E}} \approx \norm{\smash{\mat{M}^{-1}}}\norm{\mat{M}}u \eqqcolon \cond(\mat{M})u.
\end{equation*}
We see that the \warn{relative} forward error in solving the linear system of equations \cref{eq:example-linsys} using the stored floating-point representation of $\mat{M}$ is on the order of the \emph{condition number} $\cond(\mat{M}) \coloneqq \norm{\smash{\mat{M}^{-1}}}\norm{\mat{M}}$ times the unit roundoff.

This conclusion is striking.
Even the act of \emph{storing the matrix $\mat{M}$} has already introduced errors in the computed solution $\hatvector{x}$ of size $\cond(\mat{M})u$ \cite[\S2.7.11]{GV13}.
Of course, solving a linear system of equations requires many arithmetic operations, each of which introduces additional \emph{rounding errors}.
Fortunately, it can be shown that---even in the presence of such rounding errors---our heuristic calculation describes the correct magnitude of the numerical errors produced in solving $\mat{M}\vec{x} = \vec{f}$ with stable methods such as Cholesky or \QR factorization:
\actionbox{Numerically stable algorithms for solving a linear system $\mat{M}\vec{x}=\vec{f}$ produce a solution that is accurate up to a relative forward error of roughly $\cond(\mat{M})u$.}
See \cite[Chs.~7--10]{Hig02} for a rigorous treatment of perturbation theory and stability analysis for the solution of linear systems of equations.

\myparagraph{Numerical instabilities arising from the normal equations}
From our discussion of linear systems, the numerical issues with solving the normal equations \cref{eq:normal-equations} are now apparent. 
By the boxed maxim above, we would expect that solving the normal equations produces a solution with a relative error of size roughly $\cond(\mat{B}^*\mat{B})u$.
Defining the condition number of a rectangular matrix as 
\begin{equation*}
    \cond(\mat{B}) \coloneqq \norm{\mat{B}} \norm{\smash{\mat{B}^\dagger}} = \frac{\sigma_{\mathrm{max}}(\mat{B})}{\sigma_{\mathrm{min}}(\mat{B})},
\end{equation*}
we see that the condition number of $\mat{B}^*\mat{B}$ is the \emph{square} of the condition number of $\mat{B}$:
\begin{equation*}
    \cond(\mat{B}^*\mat{B}) = \cond(\mat{B})^2.
\end{equation*}
Thus, we obtain the following rule-of-thumb for solution of least-squares problems by the normal equations:
\actionbox{Solutions to the least-squares problem \cref{eq:least-squares} that \warn{explicitly} form the normal equations \cref{eq:normal-equations} yield relative forward errors of roughly $\cond(\mat{B})^2u$.}
This problem of \emph{squaring the condition number} is the main deficiency with solving least-squares problems by forming the normal equations.
Golub's paper \cite[\S1]{Gol65} identifies a related potential problem, called \emph{numerical rank deficiency}, where the \warn{numerically computed $\mat{B}^*\mat{B}$ matrix} becomes rank-deficient. 

\section{Numerically stable algorithms by \QR factorization} \label{sec:qr}

In addition to elucidating the numerical instabilities of the normal equations, Golub's paper \cite{Gol65} also identifies a solution: \QR factorization.

Recall that a \QR factorization of a matrix $\mat{B}$ refers to one of two factorizations.
The \emph{full \QR factorization} takes the form
\begin{equation*}
    \mat{B} = \mat{Q}_{\mathrm{full}} \twobyone{\mat{R}}{\mat{0}}
\end{equation*}
where $\mat{Q}_{\mathrm{full}} \in \field^{m\times m}$ is unitary and $\mat{R} \in \field^{n\times n}$ is upper triangular.
The \emph{thin} or \emph{economy-size \QR factorization} is
\begin{equation*}
    \mat{B} = \mat{Q}\mat{R},
\end{equation*}
where $\mat{Q} \in \field^{m\times n}$ has orthonormal columns and $\mat{R} \in \field^{n\times n}$ is upper triangular.
We emphasize that each of these representations constitutes an \emph{exact} decomposition of the matrix $\mat{B}$, in contrast to the partial \QR decompositions discussed in \cref{ch:low-rank-general}.
A full \QR decomposition can be converted to an economy-size \QR decomposition by extracting the first $n$ columns of $\mat{Q}_{\mathrm{full}}$; that is, $\mat{Q} \coloneqq \mat{Q}_{\mathrm{full}}(:,1:n)$.
Unless otherwise stated, \QR decompositions in this thesis are economy-sized.

The \QR factorization for a \warn{square} matrix and stable algorithms for it them was developed by Givens in the 1950s for use in computations at Oak Ridge National Laboratory.
(The technical report \cite{Giv54} describes the use of ``Givens rotations'' in eigenvalue computations.)
In 1958, Householder developed a faster and equally stable approach using relection matrices now known as Householder reflectors \cite{Hou58}.
Golub's insight was to extend Householder's method to compute \QR factorizations of a tall rectangular matrix $\mat{B}$.
From there, the unique solution $\vec{x}$ to a full-rank least-squares problem \cref{eq:least-squares} is readily obtained as 
\begin{equation*}
    \vec{x} = \mat{R}^{-1}(\mat{Q}^*\vec{b}).
\end{equation*}
Golub's procedure has dramatically better numerical stability properties than direct solution of the normal equations.
With appropriate implementation, it runs in $\order(mn^2)$ operations.

\section{Perturbation theory for least-squares} \label{sec:perturbation-theory}

Had \QR factorization truly defeated the normal equations fundamental flaw, \emph{squaring the condition number}?
Just a year after Golub's paper \cite{Gol65} was published, Golub and Wilkinson provided a somewhat disappointing answer \cite{GW66} by proving the following perturbation theorem:

\begin{fact}[First-order perturbation theory for least squares] \label{fact:ls-perturb-first-order}
    Consider the least-squares problem \cref{eq:least-squares}.
    Let $\varepsilon > 0$, and consider perturbations $\mat{\Delta B}$ and $\vec{\Delta c}$ of magnitude
    \begin{equation*}
        \norm{\mat{\Delta B}} \le \varepsilon \norm{\mat{B}}, \quad \norm{\vec{\Delta c}} \le \varepsilon \norm{\vec{c}}.
    \end{equation*}
    Then, the solution of the perturbed least-squares problem
    \begin{equation*}
        \hatvector{x} = \argmin_{\vec{z} \in \field^n} \norm{(\vec{c}+\vec{\Delta c}) - (\mat{B} +\mat{\Delta B})\vec{z} }^2,
    \end{equation*}
    satisfies
    \begin{equation*}
        \norm{\hatvector{x} - \vec{x}} \le \cond(\mat{B})\varepsilon \cdot \left(\frac{\norm{\vec{c}}}{\norm{\mat{B}}} + \norm{\vec{x}}\right) + \cond(\mat{B})^2\varepsilon \cdot \frac{\norm{\vec{c} - \mat{B}\vec{x}}}{\norm{\mat{B}}} + \order(\varepsilon^2).
    \end{equation*}
    Moreover, this bound is approximately attainable.
\end{fact}

Golub and Wilkinson summarize the message of their theorem, writing \cite[p.~144]{GW66}
\begin{quote}
    We conclude that although the use of the orthogonal transformation [i.e, \QR factorization] avoids some of the ill effects inherent in the use of the normal equations the value of $\cond(\mat{B})^2$ is still relevant to some extent.
\end{quote}
(We have amended the quote to use the present notation.)
In particular, Golub and Wilkinson's analysis implies that, for a good least-squares solver, the (forward) error should scale like $\cond(\mat{B})u$ when the residual is small and like $\cond(\mat{B})^2u$ when the residual is large.

Golub and Wilkinson's first-order perturbation theorem (\cref{fact:ls-perturb-first-order}) describes the essential nature of the sensitivity of the least-squares problems to small perturbations in the inputs.
Genuine perturbation \emph{bounds} (rather than first-order \emph{estimates}) were developed later by Wedin in his 1973 paper \cite{Wed73}; see \cite[Thm.~20.1]{Hig02} for a clean modern statement and a proof.
Here is a simplication of the (Golub--Wilkinson--)Wedin result that I have used in my work \cite{Epp24a,EMN24}: 

\begin{fact}[Perturbation bounds for least squares] \label{fact:ls-perturb}
    Consider a perturbed least-squares problem \cref{eq:least-squares}:
    \begin{equation*}
        \hatvector{x} = \argmin_{\vec{z} \in \field^n} \norm{(\vec{c} + \vec{\Delta c}) - (\mat{B}+\mat{\Delta B}) \hatvector{x} } \quad \text{with } \norm{\mat{\Delta B}} \le \varepsilon \norm{\mat{B}} ,\:\norm{\vec{\Delta c}} \le \varepsilon \norm{\vec{c}}.
    \end{equation*}
    Then, provided that $\cond(\mat{B})\varepsilon \le 0.1$, the following bounds hold.
    \begin{align}
        \norm{\vec{x} - \hatvector{x}} &\le 2.23 \, 
        \cond(\mat{B}) 
        \left(  \norm{\vec{x}} + \cond(\mat{B}) \frac{\norm{\vec{c} - \mat{B}\vec{x}}}{\norm{\mat{B}}}  \right) \varepsilon, \label{eq:forward-error-wedin} \\
        \norm{\mat{B}(\vec{x} - \hatvector{x})} &\le 2.23 \left(  \norm{\mat{B}} \norm{\vec{x}} + \cond(\mat{B}) \norm{\vec{c} - \mat{B}\vec{x}} \right) \varepsilon.\label{eq:residual-error-wedin}
    \end{align}
\end{fact}

A useful feature of this result is the bound on the error $\norm{\mat{B}(\vec{x} - \hatvector{x})}$.
I call this quantity the \emph{residual error}.
It is the norm distance between the computed and the true residual
\begin{equation*}
    \norm{\mat{B}(\vec{x} - \hatvector{x})} = \norm{(\vec{c} - \mat{B}\vec{x}) - (\vec{c} - \mat{B}\hatvector{x})}.
\end{equation*}
In addition, the residual norm of the solution $\hatvector{x}$ admits the \emph{Pythagorean decomposition}
\begin{equation} \label{eq:pythagorean}
    \norm{\vec{c} - \mat{B}\hatvector{x}}^2 = \norm{\vec{c} - \mat{B}\vec{x}}^2 + \norm{\mat{B}(\vec{x} - \hatvector{x})}^2.
\end{equation}
Observe that the bound \cref{eq:forward-error-wedin} on the \emph{forward error} $\norm{\vec{x} - \hatvector{x}}$ follows immediately from the bound \cref{eq:residual-error-wedin} on the residual error, together with the identity
\begin{equation} \label{eq:y-to-By}
    \norm{\vec{y}} \le \frac{\norm{\mat{B}\vec{y}}}{\sigma_{\mathrm{min}}(\mat{B})} = \cond(\mat{B}) \cdot \frac{\norm{\mat{B}\vec{y}}}{\norm{\mat{B}}} \quad \text{for every } \vec{y} \in \field^n.
\end{equation}

\section{Forward and backward stability} \label{sec:forward-backward-stability}

Let us take a break from our tour through history and stop to ask:
What does it mean to solve a least-squares problem accurately?
These questions have been lurking in the background of our discussion.
Now, we address them head-on.

In \cref{sec:normal-equations-unstable}, we saw that the mere act of \emph{storing} a matrix already introduces perturbations on the order of the unit roundoff $u$.
Since \emph{backward perturbations} (i.e., perturbations to the inputs $\mat{B}$ and $\vec{c}$) are inevitable, it is natural to suggest that a least-squares problem has been solved accurately if the computed solution $\hatvector{x}$ is the \warn{exact} solution to a corrupted version of the least-squares problem
\begin{equation} \label{eq:least-squares-perturbed}
    \hatvector{x} = \argmin_{\vec{z} \in \field^n} \norm{(\vec{c}+\vec{\Delta c}) - (\mat{B} + \mat{\Delta B})\vec{z}}^2.
\end{equation}
Here, the backward perturbations $\vec{\Delta c}$ and $\mat{\Delta B}$ should be of norm ``roughly'' $\order(u)$. 
Specifically, we make the following definition:

\begin{definition}[Backward stable] \label{def:backward-stability}
    A vector $\hatvector{x} \in \field^n$ is a \emph{backward stable} solution to the least-squares problem \cref{eq:least-squares} if it is the \warn{exact} solution to a modified least-squares problem \cref{eq:least-squares-perturbed}, where the perturbations admit bounds $\norm{\vec{\Delta c}} \lesssim \norm{\vec{c}}u$ and $\norm{\mat{\Delta B}} \lesssim \norm{\mat{B}}u$.
    A least-squares algorithm is \emph{backward stable} if it produces backward stable solutions on every input matrix $\mat{B}$ that is \emph{numerically full-rank}; that is, $\cond(\mat{B})u \ll 1$.
\end{definition}

Informally, we use the notation $\alpha\lesssim \beta$ to indicate that $\alpha$ is at most $\gamma\beta$ for some ``modest'' prefactor $\gamma$.
Since we are considering numerical algorithms on $m\times n$ matrices, we will allow the prefactor to depend \warn{polynomially} on the parameters $m$ and $n$.
Many authors also stipulate that the prefactor $\gamma$ be ``modest'' in the sense that it is a low-degree polynomial in $m$ and $n$ with small constants.
We will not put such a requirement our formal definition of stability, though a method with a prefactor $\gamma = 10^{100}$ is certainly not stable in any practically useful sense.
We defer a formal description of the notations $\lesssim$ and $\ll$ to \cref{sec:stability-notation}.

Backward stability is generally considered to be the gold standard stability property for a numerical algorithm (with one small caveat, discussed in \cref{rem:componentwise-stability}).
The (Householder) \QR factorization method of Golub is backward stable \cite[Thm.~20.3]{Hig02}.

The distance $\norm{\hatvector{x} - \vec{x}}$ from the computed solution $\hatvector{x}$ to the true solution $\vec{x}$ is called the \emph{forward error}, and the relative forward error is $\norm{\hatvector{x} - \vec{x}}/\norm{\vec{x}}$.
If $\hatvector{x}$ is computed by a backward stable method, the size of the forward error (and the residual error $\norm{\mat{B}(\vec{x} - \hatvector{x})}$) is controlled by the Golub--Wilkinson--Wedin theorem (\cref{fact:ls-perturb}):
\begin{align}
\norm{\vec{x} - \hatvector{x}} &\lesssim
    \cond(\mat{B}) 
    \left(  \norm{\vec{x}} + \cond(\mat{B}) \frac{\norm{\vec{c} - \mat{B}\vec{x}}}{\norm{\mat{B}}}  \right) u; \label{eq:forward-error-stable}\\
    \norm{\mat{B}(\vec{x} - \hatvector{x})} &\lesssim \left(  \norm{\mat{B}} \norm{\vec{x}} + \cond(\mat{B}) \norm{\vec{c} - \mat{B}\vec{x}} \right) u.\label{eq:residual-error-stable} 
\end{align}

As we will see later in this chapter, there are many interesting least-squares algorithms that satisfy weaker notions of numerical stability than backward stability.
In particular, we have the concept of \emph{forward stability}:

\begin{definition}[(Strong) forward stability]
    An approximate least-squares solution is $\hatvector{x} \in \field^n$ is \emph{forward stable} if it satisfies \cref{eq:forward-error-stable} and \emph{strongly forward stable} if it satisfies \cref{eq:residual-error-stable}.
    Similarly, a least-squares algorithm is (strongly) forward stable if it always produces a (strongly) forward stable solution whenever $\mat{B}$ is numerically full-rank; that is, $\cond(\mat{B})u \ll 1$.
\end{definition}

Backward stability implies strong forward stability implies forward stability, but the reverse implications do not hold.
We will encounter methods that are strongly forward stable but not backward stable in \cref{sec:sketch-and-descend,sec:spre-with-sketch-solve}.

\begin{remark}[Columnwise stability] \label{rem:componentwise-stability}
    The notion of backward stability (\cref{def:backward-stability}) is more precisely called \emph{normwise backward stability}.
    It states that the computed solution to a least-squares problem is the exact solution to a perturbed least-squares problem \cref{eq:least-squares-perturbed} that is close to the original least-squares problem \cref{eq:least-squares} \warn{in norm}.
    The Householder \QR algorithm also satisfies a stronger stability property called \emph{columnwise backward stable} \cite[Thm.~20.3]{Hig02}, where the perturbation $\mat{\Delta B}$ satisfies
    \begin{equation*}
        \norm{\vec{\Delta b}_i} \lesssim \norm{\vec{b}_i} u \quad \text{for each } i=1,\ldots,n.
    \end{equation*}
    The distinction between normwise and columnwise backward stability can be important for least-squares problems in which the column norms $\vec{b}_i$ span multiple orders of magnitude.
    
    Fortunately, any normwise backward stable least-squares solver can be converted into a columnwise backward stable solver.
    First, preprocess the matrix $\mat{B}$ by scaling its columns to have unit norm.
    Then solve the least-squares problem using a normwise stable solver.
    Finally, post-processing the computed least-squares solution $\vec{x}$ by counter-scaling its entries.
    The columnwise backward stable of this approach follows from the theorems of van der Sluis \cite{van69}.
    The conclusion is that, by using column scaling, the distinction between normwise and columnwise stability can be elided.
    For the rest of this thesis, we will only refer to normwise concepts of stability.
\end{remark}

\section{Krylov iterative methods: CGNE, CGLS, and LSQR} \label{sec:ls-iterative-methods}

In the decades after Golub introduced the Householder \QR factorization method for solving least-squares problems, the fundamental paradigm for the solution of least-squares problem remained largely unchanged.
If one requested the solution a least-squares problem in a programming environment such as MATLAB, it would be computed by Householder way of \QR factorization.
Certainly, there were improvements---researchers developed fast blocked algorithms for Householder \QR \cite{Lan98}, applied pivoting to solve sparse \cite{GN86a} and (nearly) rank-deficient \cite{BG65} and packaged high-quality \QR implementations were into optimized software packages like \LAPACK \cite{ABB+99}---but the fundamental structure ``\QR factorize then solve'' remained in place.

While the period 1966--2000 saw relatively little development in direct, high-accuracy least-squares algorithms, researchers began exploring iterative methods for least-squares problems, most notably \emph{Krylov subspace methods}.
These methods produce approximate solutions to the least-squares problem \cref{eq:least-squares} via a sequence of matrix--vector products $\vec{z} \mapsto \mat{B}\vec{z}$ and $\vec{v} \mapsto \mat{B}^*\vec{v}$.
The first methods, variously known as ``conjugate gradient least squares'' (CGLS) or ``conjugate gradient normal equations'' (CGNE), worked by applying the conjugate gradient method \cite{HS52} to the normal equations \cref{eq:normal-equations}.
Paige and Saunder provided improvement with their LSQR method \cite{PS82}.
LSQR is based on the Golub--Kahan process for bidiagonalizing a matrix \cite[\S10.4]{GV13}; it produces the same output as CGLS and CGNE in exact arithmetic, but it has improved stability properties in finite precision.
Code for LSQR is provided in \cref{prog:mylsqr}.
Research on Krylov solvers for least-squares problems continued into the 2010s, producing the LSMR \cite{FS11} and LSMB \cite{HG18} algorithms.

\myprogram{Implementation of LSQR for solving least-squares problems, iteratively.}{}{mylsqr}

We have already seen that the \emph{accuracy} of solutions to least-squares problems depends on the condition number $\cond(\mat{B})$ of the matrix $\mat{B}$.
The conditioning also controls the rate of convergence for iterative methods.
In particular, we have the following convergence result for CGNE and LSQR, derived from standard results for conjugate gradient \cite[eq.~(6.128)]{Saa03}:
\begin{fact}[CGNE and LSQR convergence] \label{fact:lsqr}
    Let $\vec{x}_1,\vec{x}_2,\ldots$ denote the iterates produced by the CGNE or LSQR algorithm (in \warn{exact arithmetic}).
    Then
    \begin{equation*}
        \norm{\mat{B}(\vec{x} - \vec{x}_k)} \le 2 \left( \frac{\cond(\mat{B}) - 1}{\cond(\mat{B}) + 1}\right)^k \norm{\mat{B}(\vec{x}-\vec{x}_0)} \le 2 \e^{-k/\cond(\mat{B})} \norm{\mat{B}(\vec{x}-\vec{x}_0)}.
    \end{equation*}
\end{fact}
This bound shows that CGNE and LSQR take at most $\order(\cond(\mat{B}) \log(1/\varepsilon))$ iterations to solve a least-sqaures problem to ``$\varepsilon$ accuracy'' in exact arithmetic.
While the bound \cref{fact:lsqr} is not always quantitatively sharp, it has some predictive power: Except in rare cases, iterative methods converge slowly for ill-conditioned problems.

Iterative methods for least-squares problems are important tools.
In particular, they provide a path to solve large sparse least-squares problems that would be impossible to solve using direct methods.
However, by the turn of the century, direct methods based on \QR factorization---proposed by Golub to solve least-squares problems in 1965---still remained the most effective algorithms for solving most least-squares problems.

\section{Randomization enters: The sketch-and-solve method} \label{sec:sketch-and-solve}

The late 1990s and early 2000s saw the origins of the modern field of randomized numerical linear \cite{FKV98,PTRV98}.
The first researchers in this effort were computer scientists, who saw randomization as a way of speeding up computations for large-scale data analysis.

Drineas, Mahoney, and Muthukrishnan developed the first randomized methods for least-squares in 2006  \cite{DMM06}.
They proposed computing an \emph{approximate} least-squares solution by subsampling a collection of rows of $\mat{B}$ (and entries of $\vec{c}$), reweighting, then solving the subsampled system.
Later that year, Sarl\'os \cite{Sar06} proposed algorithms that multiply the matrix $\mat{B}$ and $\vec{c}$ by a \emph{randomized dimensionality reduction matrix} (also called a \emph{sketching matrix}) and solve the ``sketched'' least-squares problem.
These randomized algorithms offered new strategies for quickly obtaining an approximate solution to a least-squares problem.

In modern language, the algorithms of Drineas et al.\ and Sarl\'os are examples of the \emph{sketch-and-solve paradigm}.
Sketch-and-solve algorithms for the least-squares problem \cref{eq:least-squares} proceed as follows:
\begin{enumerate}
    \item \textbf{\emph{Generate}} a random \emph{sketching matrix} $\mat{S} \in \field^{m\times d}$ with $n\le d \ll m$ columns.
    \item \textbf{\emph{Sketch}} the data $\mat{B}$ and $\vec{c}$ by computing the matrix products $\mat{S}^*\mat{B}$ and $\mat{S}^*\vec{c}$.
    \item \textbf{\emph{Solve}} the reduced least-squares problem
    \begin{equation} \label{eq:sketch-solve-def}
        \hatvector{x} = \argmin_{\vec{z} \in \field^n} \norm{(\mat{S}^*\mat{B})\vec{z} - (\mat{S}^*\vec{c})}^2.
    \end{equation}
\end{enumerate}
An implementation appears in \cref{prog:sketch_solve}.

\myprogram{Sketch-and-solve method for solving overdetermined linear least-squares problems.}{}{sketch_solve}

\section{Which sketch should I use? The subspace embedding property}

Our discussion of sketch-and-solve raises two questions: ``Which random matrix $\mat{S}$ should we pick as the sketching matrix?'' and
``How accurate is the accurate is the sketch-and-solve solution?''
To answer both of these questions, we appeal to the concept of a \emph{subspace embedding}, which is essentially due to Sarl\'os \cite{Sar06}.

\begin{definition}[Subspace embedding] \label{def:subspace-embed}
    A matrix $\mat{S} \in \field^{m\times d}$ is called a \emph{subspace embedding} for a matrix $\mat{F} \in \field^{m\times n}$ with \emph{distortion} $\eta \in (0,1)$ when
    \begin{equation} \label{eq:subspace-embedding}
        (1-\eta) \norm{\mat{F}\vec{z}} \le \norm{(\mat{S}^*\mat{F})\vec{z}} \le (1+\eta) \norm{\mat{F}\vec{z}} \quad \text{for every } \vec{z} \in \field^n.
    \end{equation}
    A sketching matrix $\mat{S}$ is said to be an \emph{oblivious subspace embedding} with \emph{dimension} $n$, distortion $\eta$, and \emph{failure probability} $\delta \in [0,1)$ if, for every matrix $\mat{F} \in \field^{m\times n}$, the condition \cref{eq:subspace-embedding} holds with probability at least $1-\delta$.
\end{definition}

A subspace embedding maps a high-dimensional vector $\vec{v} \in \field^m$ to a low-dimensional vector $\mat{S}^*\vec{v} \in \field^d$ while \emph{simultaneously} preserves the norm of each vector $\vec{v} \in \range(\mat{F})$.
The range of the $\mat{F}$ is the titular subspace, for which $\mat{S}$ is an embedding.
The following theorem characterizes the accuracy of the sketch-and-solve method when implemented with a subspace embedding:
\begin{theorem}[Sketch-and-solve] \label{thm:sketch-and-solve}
    Let $\mat{S}\in\field^{d\times m}$ be a subspace embedding \warn{for $\flatonebytwo{\mat{B}}{\vec{c}}$} of distortion $\eta$, and let $\hatvector{x}$ be the sketch-and-solve solution \cref{eq:sketch-solve-def}.
    Then
    \begin{align*}
        \norm{\vec{c} - \mat{B}\hatvector{x}} &\le \min \left\{ \left(1 + \frac{9\eta^2}{(1-\eta)^4}\right)^{1/2},\frac{1+\eta}{1-\eta} \right\} \cdot \norm{\vec{c} - \mat{B}\vec{x}} \\ &= (1+\order(\eta^2))\norm{\vec{c} - \mat{B}\vec{x}}.
    \end{align*}
\end{theorem}
The theory for sketch-and-solve is surprisingly nuanced; see \cref{app:sketch-and-solve-analysis} for an introduction this theory, a proof of this result, and a discussion of history.

\begin{figure}
    \centering
    \includegraphics[width=0.7\linewidth]{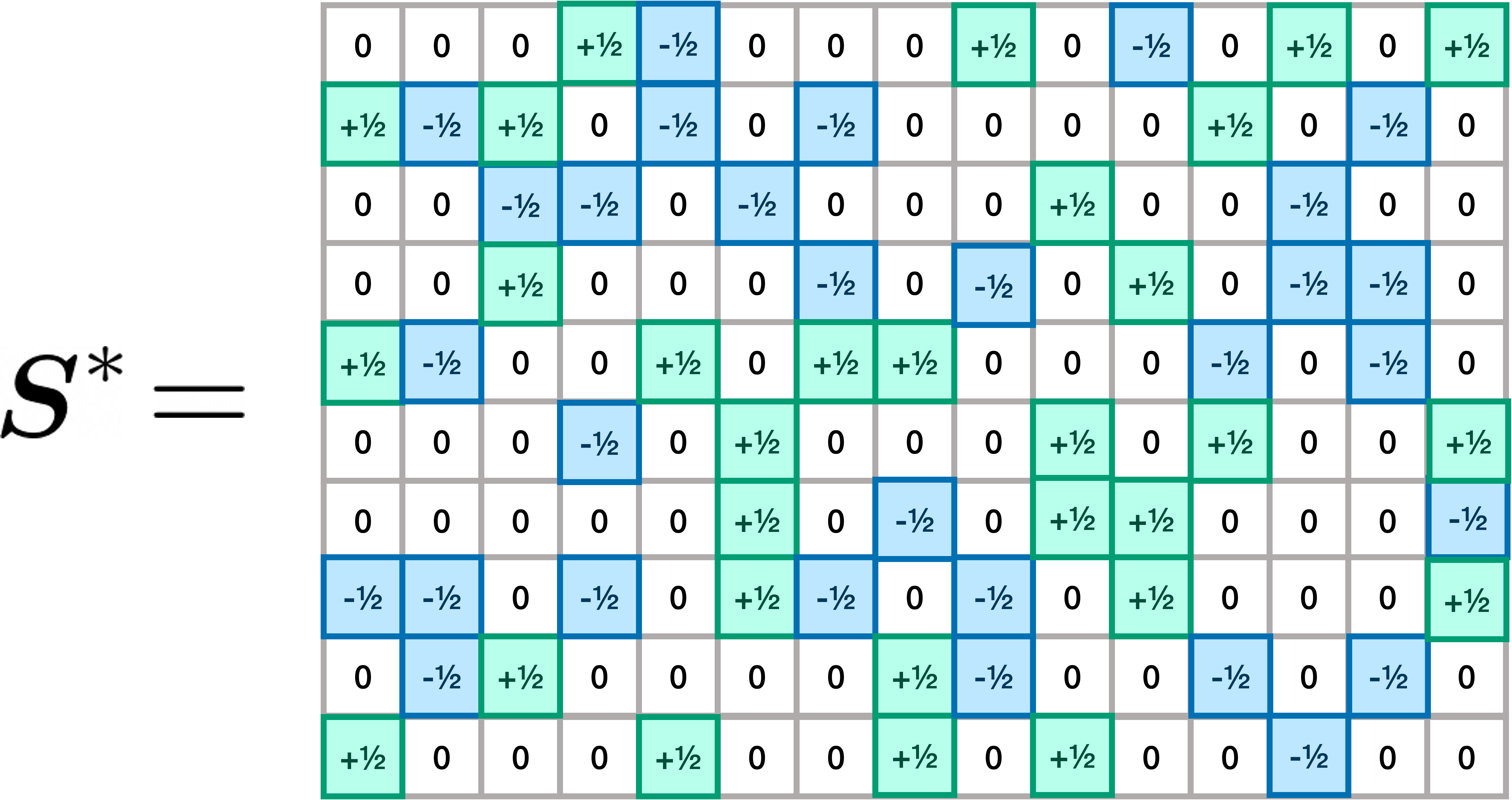}
    \caption[Illustration of sparse sign embedding]{Example of a sparse sign embedding, illustrated with $\zeta = 4$, $m = 15$, and $d = 10$.
    Observe that there are precisely $\zeta$ nonzero entries per \warn{row} of $\mat{S}$ or \warn{column} of $\mat{S}^*$.}
    \label{fig:sparse-sign-cartoon}
\end{figure}

\Cref{thm:sketch-and-solve} provides guarantees for the that the sketch-and-solve method when implemented with any subspace embedding.
The literature contains many constructions of subspace embeddings (see \cite[\S\S8--9]{MT20a} and \cite[Chs.~2, 6, \& 7]{MDM+23a}), and debates about which construction to use continue to this day.
\Cref{ch:which-sketch} weighs in on this debate by providing detailed discussion of the most popular subspace embedding constructions and presenting numerical comparisons. 
Based on the results documented in this appendix, I recommend \emph{sparse sign embeddings} as the subspace embedding of choice for most applications:
\begin{definition}[Sparse sign embedding] \label{def:sparse-sign-body}
    A \emph{sparse sign embedding} $\mat{S} \in \real^{m\times d}$ with \emph{sparsity} $\zeta$ is a random matrix constructed as follows.
    Each row is independent and possesses exactly $\zeta$ nonzero entries.
    The nonzero entries are placed in uniformly random positions (selected \warn{without replacement}) and have uniform $\pm \zeta^{-1/2}$ values.
\end{definition}
As we document in \cref{ch:which-sketch}, sparse-sign embeddings have a nearly ideal distortion parameter $\eta$ and support a fast product operation $\mat{B} \mapsto \mat{S}^*\mat{B}$, owing to sparsity.
See \cref{fig:sparse-sign-cartoon} for a cartoon illustration of the sparse sign embedding construction.

The best-known analysis fo \emph{sparse} oblivious subspace embeddings due to Cohen \cite[Thm.~4.2]{Coh16}:

\begin{fact}[Sparse sign embeddings] \label{thm:ose-existence}
    Fix parameters $n\ge 1$, $\eta \in(0,1)$, and $\delta \in (0,1)$.
    With embedding dimension $d = \order(n\log (n/\delta)/\eta^2)$ and sparsity $\zeta = \order(\log(n/\delta)/\eta)$, a sparse sign embedding $\mat{S}\in\real^{m\times d}$ is an oblivious subspace embedding with dimension $n$, distortion $\eta$, and failure probability $\delta$.
    Since $\mat{S}$ is sparse, matrix--vector products $\mat{S}^*\vec{v}$ can be computed in $\order(n\log(n/\delta)/\eta)$ operations.
\end{fact}
Evidence provided in \cref{ch:which-sketch} suggests that sparse sign embeddings are reliable with more aggressive parameter choices (e.g., $d = 2n$ and $\zeta = 8$).
Providing a theoretical explanation for the reliable performance of sparse sketching matrices with smaller values of $d$ and $\zeta$ has been a major research effort in recent years; see \cite{CDDR24,CDD24,Tro25}.

As a corollary of \cref{thm:ose-existence}, we get the following runtime guarantees for sketch-and-solve with a sparse sign embedding:

\begin{corollary}[Sketch-and-solve] \label{cor:sketch-and-solve-with-sparse-sign}
    Using a sparse sign embedding with Cohen's parameter settings (\cref{thm:ose-existence}), one can achieve a $(1+\varepsilon)$-approximate least-squares solution
    \begin{equation*}
        \norm{\vec{c} - \mat{B}\hatvector{x}} \le (1+\varepsilon)\norm{\vec{c} - \mat{B}\vec{x}}
    \end{equation*}
    with 99\% probability with a sketching matrix of embedding dimension
    \begin{equation*}
        d = \order\left(\frac{n\log n}{\varepsilon}\right).
    \end{equation*}
    The total runtime is 
    \begin{equation} \label{eq:sketch-and-solve-runtime}
        \order\left( \frac{mn \log n}{\varepsilon^{1/2}} + \frac{n^3 \log n}{\varepsilon} \right).
    \end{equation}
\end{corollary}

\Cref{cor:sketch-and-solve-with-sparse-sign} shows that, using a sketching matrix like a sparse sign embedding with a fast multiply operation, sketch-and-solve runs in roughly $\order(mn + n^3)$ operations.
For large, highly overdetermined problems ($m\gg n\gg 1$), this is a large speedup over Householder \QR factorization.
However, the runtime sketch-and-solve depends inversely on the least-squares suboptimality $\varepsilon$.
As such, improving the accuracy of sketch-and-solve is costly, and sketch-and-solve is mainly useful as a way of obtaining a cheap, approximate solution to a least-squares problem for a relatively large value of $\varepsilon \in [0.25,1]$.
Indeed, we have already seen sketch-and-solve type methods used to construct approximate interpolative decompositions\index{interpolative decomposition} in \iffull
\else\cref{sec:sketchy-pivoting} and \CUR decompositions in \cref{sec:cur-algs}.
(Note that weighted row sampling is a type of sketching.)\fi 

\section{Making randomized least-squares accurate: Sketch-and-precondition} \label{sec:spre}

The papers \cite{DMM06,Sar06} demonstrate that randomized algorithms offer a fast way of obtaining an approximate least-squares solutions.
Can randomized methods also solve least-squares problems to high accuracy?
Rokhlin and Tygert answered this question by developing the sketch-and-precondition method \cite{RT08}.
This method was popularized by Avron, Maymounkov, and Toledo \cite{AMT10}, who developed a variant of the algorithm they called \emph{Blendenpik}.
The papers \cite{MSM14,CFS21} propose further variants of sketch-and-precondition.

Rokhlin and Tygert's core observation is that the sketched matrix $\mat{S}^*\mat{B}$ provides an excellent preconditioner for the least-squares problem \cref{eq:least-squares}.
This preconditioner can be combined with a Krylov iterative method (\cref{sec:ls-iterative-methods}) such as CGNE or LSQR to solve the least-squares problem rapidly to high precision.

Here is the main technical result for sketch-and-precondition, in modern language.

\begin{proposition}[Randomized preconditioning] \label{prop:randomized-preconditioning}
    Let $\mat{S} \in \field^{d\times m}$ be a subspace embedding \warn{for $\mat{B}$} with distortion $\eta < 1$.
    Construct the sketched matrix $\mat{S}^*\mat{B}$, and let
    \begin{equation*}
        \mat{S}^*\mat{B} = \mat{U}\mat{M}
    \end{equation*}
    be \emph{any} factorization $\mat{S}^*\mat{B}$ into a matrix $\mat{U} \in \field^{d\times n}$ with orthonormal columns and a nonsingular matrix $\mat{M}$.
    Then
    \begin{equation*}
        \sigma_{\mathrm{max}}(\mat{B}\mat{M}^{-1}) \le \frac{1}{1-\eta}, \quad \sigma_{\mathrm{min}}(\mat{B}\mat{M}^{-1}) \ge \frac{1}{1+\eta}.
    \end{equation*}
    Consequently,
    \begin{equation*}
        \cond(\mat{B}\mat{M}^{-1}) \le \frac{1+\eta}{1-\eta} = 1 + \order(\eta).
    \end{equation*}
\end{proposition}

\begin{proof}
    We establish the upper bound on $\sigma_{\mathrm{max}}(\mat{B}\mat{M}^{-1})$; the lower bound on $\sigma_{\mathrm{min}}(\mat{B}\mat{M}^{-1})$ follows on similar lines.
    By the variational characterization of the largest singular value, we have
    \begin{equation*}
        \sigma_{\mathrm{max}}(\mat{B}\mat{M}^{-1}) = \max_{\norm{\vec{v}} = 1} \, \norm{\mat{B}\mat{M}^{-1}\vec{v}}.
    \end{equation*}
    Now, apply the subspace embedding property to bound
    \begin{equation*}
        \sigma_{\mathrm{max}}(\mat{B}\mat{M}^{-1}) = \max_{\norm{\vec{v}} = 1} \, \norm{\mat{B}\mat{M}^{-1}\vec{v}} \le \frac{1}{1-\eta}\max_{\norm{\vec{v}} = 1} \, \norm{\mat{S}^*\mat{B}\mat{M}^{-1}\vec{v}}.
    \end{equation*}
    But $\mat{S}^*\mat{B}\mat{M}^{-1} = \mat{U}$, so
    \begin{multline*}
        \sigma_{\mathrm{max}}(\mat{B}\mat{M}^{-1}) = \max_{\norm{\vec{v}} = 1}\, \norm{\mat{B}\mat{M}^{-1}\vec{v}} \le \frac{1}{1-\eta}\cdot\max_{\norm{\vec{v}} = 1}\,  \norm{\mat{S}^*\mat{B}\mat{M}^{-1}\vec{v}} \\ = \frac{1}{1-\eta}\cdot\max_{\norm{\vec{v}} = 1} \, \norm{\mat{U}\vec{v}} = \frac{1}{1-\eta}\cdot\max_{\norm{\vec{v}} = 1} \, \norm{\vec{v}} = \frac{1}{1-\eta}.
    \end{multline*}
    In the penultimate equality, we used the fact that $\mat{U}$ has orthonormal columns, and thus $\norm{\mat{U}\vec{v}} = \norm{\vec{v}}$.
\end{proof}

Provided that $\mat{S}$ is a subspace embedding with small distortion $\eta$, randomized preconditioning brings the condition number of $\mat{B}$ below an \emph{absolute constant}.
For example, $\cond(\mat{B}\mat{M}^{-1})\le 3$ when $\eta\le 1/2$.
This results shows that randomization produces a preconditioner of \emph{extraordinary} quality.
Using this preconditioner with a Krylov method like LSQR or CGNE produces a solution that converges very rapidly (see \cref{cor:sketch_and_pre} below).
This observation motivates the following (prototype) sketch-and-precondition algorithm:
\begin{enumerate}
    \item \textbf{\textit{Generate}} a sketching matrix $\mat{S} \in \real^{n\times d}$.
    \item \textbf{\textit{Sketch}} the matrix $\mat{B}$, computing $\mat{S}^*\mat{B}$.
    \item \textbf{\textit{Form}} an orthonormal decomposition $\mat{S}^*\mat{B} = \mat{U}\mat{M}$.
    \item \textbf{\textit{Solve}} the least-squares problem \cref{eq:least-squares} using LSQR with preconditioner $\mat{M}$.
\end{enumerate}
Implementations details will follow in next section and code is provided in \cref{prog:sketch_precondition}.

As we will see below in \cref{cor:sketch_and_pre}, with an appropriate implementation, the sketch-and-precondition algorithm runs in at most operations $\order(mn \log (m/\varepsilon) + n^3 \log n)$ to produce a $(1+\varepsilon)$-approximate least-squares solution.
Since the accuracy parameter $\varepsilon$ appears inside the logarithm, we see that the sketch-and-precondition method can achieve high accuracy (machine-precision accuracy, even) in a small number of LSQR steps.
Contrast this with the sketch-and-solve method, whose runtime scales with the inverse-accuracy parameter $1/ \varepsilon$ and thus cannot achieve high accuracy without an astronomical slowdown.
For large, highly over-determined problems ($1\ll n \ll m$) sketch-and-precondition is asymptotically faster than the \QR factorization method, which requires $\order(mn^2)$ operations.

I consider the sketch-and-precondition method to be one of the most striking examples of a randomized algorithm for matrix computations.
Here, randomness can be used to solve a least-squares problems to any desired accuracy level (in exact arithmetic at least) at a cost that is (asymptotically) faster than the \QR factorization method.
Forty-five years after its ascendance to the throne, the \QR method faces its first serious challenger as the method of choice for solving general dense least-squares problems, at least ones satisfying $1\ll n \ll m$.

\section{Implementation and analysis of sketch-and-precondition}

The sketch-and-precondition algorithm will be important in the rest of this thesis.
In this section, we will discuss implementation choices and present analysis.

\myparagraph{Choice of orthonormal decomposition}
To put the randomized preconditioning idea into practice, the two natural orthonormal decompositions are the \QR decomposition 
\begin{equation*}
    \mat{S}^*\mat{B} = \mat{Q} \mat{R}
\end{equation*}
and the (economy-size) SVD
\begin{equation*}
    \mat{S}^*\mat{B} = \mat{U}\mat{M} \quad \text{for } \mat{M} = \mat{\Sigma}\mat{V}^*.
\end{equation*}
In the former case, the inverse-preconditioning operation $\mat{R}^{-1}\vec{z}$ can efficiently performed by triangular substitution.
In the latter case, $\mat{M}^{-1} = \mat{V}\mat{\Sigma}^{-1}$ has a simple closed form.
For reasons that will emerge in \cref{sec:backerr-estimate}, we will use the SVD-based implementation for this thesis.

\myparagraph{Initialization}
In their original paper, Rokhlin and Tygert suggested initializing via the sketch-and-solve and solve solution
\begin{equation*}
    \vec{x}_0 \coloneqq (\mat{S}^*\mat{B})(\mat{S}^*\vec{c}) = \mat{M}^{-1}(\mat{U}^*(\mat{S}^*\vec{c})).
\end{equation*}
to initialize the sketc
Their motivation for using this initialization was to reduce the number of iteration steps required to converge.
As we will see, the choice of initialization also has \emph{numerical stability} implications for the sketch-and-precondition algorithm (\cref{sec:spre-with-sketch-solve}).

Interestingly, the matter of initialization is not mentioned in the influential paper \cite{AMT10}.
By default, their code uses the trivial, zero initialization $\vec{x}_0 \coloneqq \vec{0}$, although they do provide the sketch-and-solve initialization as an optional flag.

\myparagraph{Choice of iterative method}
The sketch-and-precondition algorithm can be implemented using any preconditioned Krylov iterative method.
In this thesis, we will use preconditioned LSQR in all our implementations.

In this thesis, we will precondition LSQR ``explicitly'' by performing the change of variables
\begin{equation} \label{eq:change-of-vars-x-y}
    \vec{x} = \mat{M}^{-1}\vec{y}
\end{equation}
and applying standard LSQR to the preconditioned least-squares problem 
\begin{equation} \label{eq:least-squares-pre}
    \vec{y} = \argmin_{\vec{z} \in \field^n} \norm{\vec{c} - \mat{B}(\mat{M}^{-1}\vec{z})}.
\end{equation}
We emphasize that \warn{the product $\mat{B}\mat{M}^{-1}$ should not be formed explicitly}.
Instead, its action is achieved by applying the matrices in sequence $\vec{z} \mapsto \mat{B}(\mat{M}^{-1}\vec{z})$.
Likewise, the adjoint is applied as $\vec{v} \mapsto \mat{M}^{-*}(\mat{B}^*\vec{v})$.
Once \cref{eq:least-squares-pre} has been solved, the solution $\vec{x}$ to the least-squares problem \cref{eq:least-squares} is obtained from the change-of-variables formula \cref{eq:change-of-vars-x-y}.
See also \cite[Alg.~8]{Mei24} for an alternate approach to preconditioning LSQR.

\myprogram{Sketch-and-precondition method for solving overdetermined linear least-squares problems.}{Subroutine \texttt{mylsqr} is provided in \cref{prog:mylsqr}.}{sketch_precondition}

\myparagraph{Implementation}
An implementation of sketch-and-precondition using the above recommendation appears in \cref{prog:sketch_precondition}.

\myparagraph{Analysis}
By combining \cref{fact:lsqr,thm:sketch-and-solve,prop:randomized-preconditioning}, we can analyze the sketch-and-precondition method.
\begin{corollary}[Sketch-and-precondition: Exact arithmetic] \label{cor:sketch_and_pre}
    Let $\mat{S} \in \field^{d\times m}$ be a subspace embedding \warn{for $\flatonebytwo{\mat{B}}{\vec{c}}$} with distortion $\eta < 1$.
    The sketch-and-precondition algorithm (\cref{prog:sketch_precondition}) \warn{in exact arithmetic} satisfies the bound
    \begin{equation} \label{eq:sketch-and-pre-guarantee}
        \norm{\vec{c} - \mat{B}\vec{x}_k}^2 \le \left( 1 + \frac{36}{(1-\eta)^4} \cdot \eta^{2+k}\right)\norm{\vec{c} - \mat{B}\vec{x}}^2.
    \end{equation}
    In particular, using a sparse sign embedding with $\eta = 1/2$ (see \cref{thm:ose-existence}), sketch-and-precondition produces a $(1+\varepsilon)$-approximate least-squares solution
    \begin{equation*}
        \norm{\vec{c} - \mat{B}\vec{x}_k} \le (1+\varepsilon)\norm{\vec{x} - \mat{B}\vec{x}}.
    \end{equation*}
    after $k = \order(\log(1/\varepsilon))$ iterations and $\order(mn \log(n/\varepsilon) + n^3\log n)$ operations.
\end{corollary}

\begin{proof}
    By \cref{thm:sketch-and-solve},
    \begin{equation*}
        \cond(\mat{B}\mat{M}^{-1}) \le \frac{1+\eta}{1-\eta}.
    \end{equation*}
    Therefore,
    \begin{equation*}
        \frac{\cond(\mat{B}\mat{M}^{-1}) - 1}{\cond(\mat{B}\mat{M}^{-1}) + 1} \le \eta.
    \end{equation*}
    Therefore, by the Pythagorean identity \cref{eq:pythagorean} and the LSQR guarantee (\cref{fact:lsqr}), we obtain the bound 
    \begin{equation*}
        \norm{\vec{c} - \mat{B}\vec{x}_k}^2 \le 4 \eta^{2k} \norm{\mat{B}(\vec{x} - \vec{x}_0)}^2 + \norm{\vec{c} - \mat{B}\vec{x}}^2.
    \end{equation*}
    Finally, invoking the Pythagorean identity \cref{eq:pythagorean} and the sketch-and-solve guarantee (\cref{thm:sketch-and-solve}),
    \begin{equation*}
        \norm{\mat{B}(\vec{x} - \vec{x}_0)}^2 = \norm{\vec{c} - \mat{B}\vec{x}_0}^2 - \norm{\vec{c} - \mat{B}\vec{x}}^2 \le \frac{9\eta^2}{(1-\eta)^4}\norm{\vec{c} - \mat{B}\vec{x}}^2.
    \end{equation*}
    Combining the two previous displays, we obtain the sketch-and-precondition guarantee \cref{eq:sketch-and-pre-guarantee}.
    The runtime guarantee for producing a $(1+\varepsilon)$-approximate least-squares solution follows from \cref{thm:ose-existence}.
\end{proof}

\section{Is sketch-and-precondition numerically unstable?} \label{sec:spre-unstable?}

In 2023, Meier, Nakasukasa, Townsend, and Webb released a preprint demonstrating that the version of sketch-and-precondition from the popular paper \cite{AMT10}---that is, the version \warn{with the zero initialization $\vec{x}_0 \coloneqq \vec{0}$}---is numerically unstable \cite{MNTW24}.


\begin{figure}
    \centering
    \includegraphics[width=0.98\linewidth]{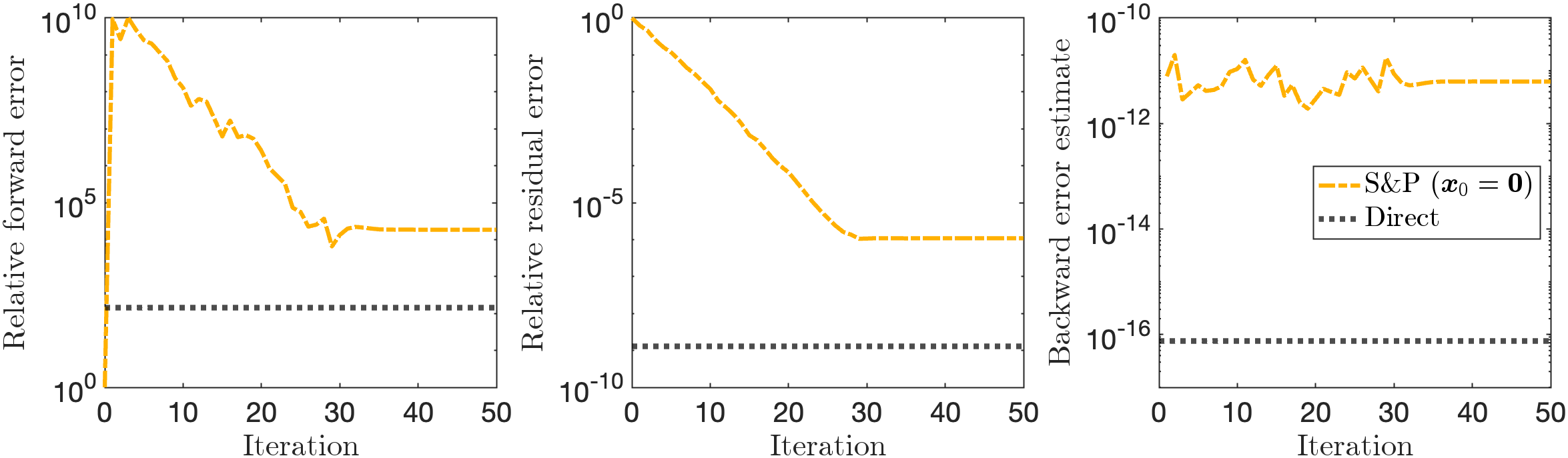}
    \caption[Forward, residual, and backward error for sketch-and-precondition with the zero initialization]{Forward (\emph{left}), residual (\emph{middle}), and backward (\emph{right}) error of sketch-and-precondition \warn{with the zero initialization} (single execution).
    Errors for direct method (MATLAB's \texttt{mldivide}) shown for reference.}
    \label{fig:sketch-and-precondition-unstable}
\end{figure}

\myprogram{Generate a random least-squares problem with prescribed dimensions, condition number, and residual norm.}{Subroutine \texttt{haarorth} is provided in \cref{prog:haarorth}.}{random_ls_problem}

A variant of their experiment appears in \cref{fig:sketch-and-precondition-unstable}.
Here, I generated a random least-squares problem of dimensions $4000\times 50$ using \cref{prog:random_ls_problem} with condition number $\texttt{condB} = \cond(\mat{B}) = 10^{12}$ and residual norm $\texttt{rnorm} = \norm{\vec{c} - \mat{B}\vec{x}} = 10^{-4}$.
The sketching matrix $\mat{S}$ is a sparse sign embedding (\cref{def:sparse-sign-body}) with dimensions $4000\times 100$ and sparsity level $\zeta = 4$.

The left panel of \cref{fig:sketch-and-precondition-unstable} plots the relative forward error $\norm{\vec{x} - \hatvector{x}} / \norm{\vec{x}}$; the middle panel plots the relative residual error $\norm{\mat{B}(\vec{x} - \hatvector{x})} / \norm{\mat{B}\vec{x}}$; and the right panel plots the \emph{backward error} $\norm{\mat{\Delta B}}_{\mathrm{F}} / \norm{\mat{B}}_{\mathrm{F}}$, the magnitude of the \warn{minimal} backward perturbation $\mat{\Delta B}$ in \cref{eq:least-squares-perturbed}.
In fact, we report an \emph{estimate} of this quantity, computed using \cref{prog:backerr_est} described below.
As we see, all of these error quantities are well above the error levels for MATLAB's \QR-based \texttt{mldivide}, shown as a black dotted line.
Therefore, we recognize that sketch-and-precondition \warn{with the zero initialization $\vec{x}_0 \coloneqq \vec{0}$} is numerically unstable, in both the backward and forward senses (see \cref{sec:forward-backward-stability}).

The work of Meier et al.\ raised serious questions about the viability of sketch-and-precondition as a replacement for least-squares solvers based on \QR factorization, at least in applications that demand high numerical accuracy.
Perhaps \QR factorization's place on the throne as the default least-squares solver for general-purpose software remained safe after all.

Fortunately, the instabilities of sketch-and-precondition are curable.
In next chapter, we will develop fast, numerically stable randomized least-squares solvers, including a simple modification of the sketch-and-precondition that is backwards.
Using these stable algorithms, we can solve least-squares problems to machine accuracy up to an order of magnitude faster than traditional least-squares solvers based on \QR factorization.

\chapter{Fast, stable randomized least-squares solvers} \label{ch:stable-ls}

\epigraph{For me, \texttt{eig(A)} epitomizes the successful contribution of numerical analysis to our technological world. Physicists, chemists, engineers, and mathematicians know that computing eigenvalues of matrices is a solved problem. Simply invoke \texttt{eig(A)}\ldots and you tap into the work of generations of numerical analysts. The algorithm involved, the \QR algorithm, is completely reliable, utterly nonobvious, and amazingly fast.}{Lloyd N.\ Trefethen, \textit{An Applied Mathematician's Apology} \cite[p.~52]{Tre22}}

Last chapter, we saw the results of Meier, Nakatsuksa, Townsend, and Webb \cite{MNTW24}, who showed that a common implementation of the sketch-and-precondition algorithm is numerically unstable, achieving errors in floating-point arithmetic that are are orders of magnitude higher than standard direct methods.
Meier et al.'s work left open the question of whether \emph{any} randomized least-squares solver is numerically stable, even in the forward sense.

This chapter answers this question in the affirmative.
We discuss three algorithms: the \emph{sketch-and-descend algorithm} and sketch-and-precondition algorithm \warn{with sketch-and-solve initialization}, which are (strongly) forward stable but not backward stable, and \emph{sketch-and-precondition with iterative refinement} (SPIR), which is stable in both the (strong) forward and backward senses.

\myparagraph{Sources}
This chapter is based on the papers

\fullcite{Epp24a}

and

\fullcite{EMN24}.

\myparagraph{Outline}
After some notation in \cref{sec:stability-notation}, we discuss the sketch-and-descend, sketch-and-precondition, and SPIR algorithms in \cref{sec:sketch-and-descend,sec:spre-with-sketch-solve,sec:spir}.
For our stability proofs, we will need to analyze modified versions of sketch-and-precondition using the Lanczos algorithm in place of LSQR.
For this reason, our treatment will be briefly interrupted in \cref{sec:lanczos} with an introduction to the Lanczos algorithm and its relationship to other Krylov methods.
\Cref{sec:stable-ls-experiments} contains experiments, and \cref{sec:backerr-estimate} investiates estimation of the backward error.

\section{Notation} \label{sec:stability-notation}

Before we go further, let us establish some notation.
We are interested in understanding the behavior of iterative least-squares algorithms that operate on $m\times n$ matrices and may require $\order(\log(1/u))$ iterations to converge to the maximum numerically achievable accuracy.
As such, the total number of arithmetic operations for algorithms in this chapter will be a polynomial in these three parameters: the dimensions $m$ and $n$ and the iteration count $\log(1/u)$.
Since each arithmetic operation incurs some rounding error, it is natural to expect that iterative least-squares algorithms should be only accurate up to error $\poly(m,n,\log(1/u))u$.
(Here, $\poly$ denotes an unspecified polynomial function.)

With this in mind, we write $\alpha \lesssim \beta$ or $\beta \gtrsim \alpha$ if $\alpha \le \gamma \beta$ for a prefactor $\gamma = \poly(m,n,\log(1/u))$.
We write $\alpha \asymp \beta$ if $\alpha \lesssim \beta \lesssim \alpha$.
Similarly, we write $\alpha \ll \beta$ if $\alpha \le \beta/\gamma$ for a \emph{sufficiently large} polynomial factor $\gamma = \poly(m,n,\log(1/u))$.
(That is, the factor $\gamma$ is as large a polynomial function as needed to complete the analysis.) 
For an algorithm to be stable in a practical sense, we hope the polynomial factors $\gamma$ are small.
As we will see in \cref{ch:stability-proofs}, the proofs are already complicated when we make no effort to track the prefactors.

We work in the standard model of floating-point arithmetic \cite[\S2.2]{Hig02}.
We write $\fl(\vec{z})$ for the value of a quantity $\vec{z}$ evaluated in floating-point, and we introduce $\err(\vec{z}) \coloneqq \fl(\vec{z}) - \vec{z}$.
When using the $\fl$ and $\err$ notations, we will always assume that the quantity is computed most natural way, i.e., $\mat{A}\mat{R}^{-1}\vec{z}$ is evaluated as \texttt{A*(R\textbackslash{}z)}, not \texttt{A*(inv(R)*z)} and certainly not \texttt{(A*inv(R))*z}.

Throughout, we focus on the least-squares problem \cref{eq:least-squares} with matrix $\mat{B}$ and right-hand side $\vec{c}$.
The condition number of $\mat{B}$ will be denoted $\kappa \coloneqq \cond(\mat{B})$.

\section{Sketch-and-descend: A strongly forward stable least-squares solver} \label{sec:sketch-and-descend}

After Meier, Nakatsukasa, Townsend, and Webb demonstrated instabilities in sketch-and-precondition \cite{MNTW24}, it was natural to ask whether other least-squares algorithms were numerically stable.
(The stabilizing effects of the sketch-and-solve initialization, discussed below in \cref{sec:spre-with-sketch-solve}, were not known at the time.)
Joel Tropp suggested to me the possibility that \emph{iterative sketching} methods would be numerically stable.
As it happens, these methods are strongly forward stable (but not backward stable) when implemented appropriately.

\subsection{Iterative sketching methods and the sketch-and-descend method}

We now describe the family of iterative sketching methods, culminating in a method we call the \emph{sketch-and-descend} method.

\myparagraph{Iterative Hessian sketch}
In 2016, Pilanci and Wainwright developed an alternative to the sketch-and-precondition paradigm for rapidly solving least-squares problems to high-accuracy.
They called their approach \emph{iterative Hessian sketch} \cite{PW16}.
Their basic idea is to reformulate the least-squares problem \cref{eq:least-squares} as a sum of a quadratic term and a linear term:
\begin{equation} \label{eq:ihs-1}
    \vec{x} = \argmin_{\vec{z} \in \field^n} \left[ \norm{\mat{B}\vec{z}}^2 - 2\Re((\mat{B}^*\vec{c})^*\vec{z}) \right].
\end{equation}
Then, reformulate this equation by introducing an offset $\vec{x}_i$:
\begin{equation} \label{eq:ihs-2}
    \vec{x} =  \argmin_{\vec{z} \in \field^n} \left[ \norm{\mat{B}(\vec{z} - \vec{x}_i)}^2 - 2\Re((\mat{B}^*[\vec{c} - \mat{B}\vec{x}_i])^*(\vec{z}-\vec{x}_i)) \right].
\end{equation}
Both \cref{eq:ihs-1,eq:ihs-2} are exact reformulations of the least-squares problem \cref{eq:least-squares}, with all three optimization problems sharing the same solution $\vec{x}$ (which is unique when $\mat{B}$ has full rank).

Now, we can introduce sketching.
Beginning from an initialization $\vec{x}_0 \in \field^n$, we iterate for $i = 0,1,\ldots$.
In their original algorithm, Pilanci and Wainwright draw a \warn{new sketch $\mat{S}_i \in \field^{m\times d}$} at each iteration, and they sketch only the quadratic term, resulting in the iteration
\begin{equation*}
    \vec{x}_{i+1} = \argmin_{\vec{z} \in \field^n} \left[ \norm{\mat{S}_i^*\mat{B}(\vec{z} - \vec{x}_i)}^2 - 2\Re((\mat{B}^*[\vec{c} - \mat{B}\vec{x}_i])^*(\vec{z}-\vec{x}_i)) \right].
\end{equation*}
Solving this optimization problem, we see that the resulting iteration is 
\begin{equation*}
    \vec{x}_{i+1} \coloneqq \vec{x}_i + [(\mat{S}_i^*\mat{B})^*(\mat{S}_i^*\mat{B})]^{-1}\mat{B}^*(\vec{c} - \mat{B}\vec{x}_i).
\end{equation*}
This iteration constitutes Pilanci and Wainwright's iterative Hessian sketch algorithm for the least-squares problem \cref{eq:least-squares}.

An advantage of the iterative Hessian sketch approach is that it generalizes nicely to more general optimization problems, including \emph{constrained} least-squares problems; see \cite{PW16} for details.
For this thesis, we will focus only on the unconstrained least-squares problem \cref{eq:least-squares}.

\myparagraph{From iterative Hessian sketch to sketch-and-descend}
In the years following Pilanci and Wainwright's original paper, researchers developed optimizations to the iterative Hessian sketch method for the unconstrained linear least-squares problems \cref{eq:least-squares}.
First, the iteration-dependent sketching matrix $\mat{S}_i$ was replaced by a single sketching matrix $\mat{S} \in\field^{m\times d}$ used across all iterations.
Second, the procedure was accelerated by incorporating dampling and momentum
\begin{equation} \label{eq:ihs-optimized}
    \vec{x}_{i+1} \coloneqq \vec{x}_i + \alpha [(\mat{S}^*\mat{B})^*(\mat{S}^*\mat{B})]^{-1}\mat{B}^*(\vec{c} - \mat{B}\vec{x}_i) + \beta (\vec{x}_i - \vec{x}_{i-1}).
\end{equation}
When implemented with a subspace embedding of distortion $\eta$, the optimal values of the coefficients $\alpha$ and $\beta$ are 
\begin{equation} \label{eq:optimal-coeffs}
    \alpha = (1-\eta^2)^2, \quad \beta = \eta^2.
\end{equation}
These refinements were developed in concurrent papers of Ozaslan, Pilanci, and Arikan \cite{OPA19} and Lacotte and Pilanci \cite{LP21}.

In my view, the name iterative Hessian sketch is somewhat of a misnomer for the modern version of the procedure \cref{eq:ihs-optimized} for least squares, as only a single sketch is computed.
(When I have publicly presented on this work, the names \emph{iterative Hessian sketch} and \emph{iterative sketching} for this procedure caused confusion.)
With the optimizations of \cite{OPA19,LP21}, the iteration \cref{eq:ihs-optimized} can be interpreted as preconditioned gradient descent with heavy-ball momentum for the least-squares objective \cref{eq:least-squares}.
Thus, the procedure \cref{eq:ihs-optimized} can be interpreted as a version of sketch-and-precondition where a (momentum-accelerated) gradient method is used in place of a Krylov solver.
To avoid potential confusion, we will refer to the iteration \cref{eq:ihs-optimized} as the \emph{sketch-and-descend} method for this thesis.

\subsection{Numerically stable implementation of sketch-and-descend}

The original works \cite{PW16,OPA19,LP21} do not offer precise guidance on how to implement the iterative Hessian sketch/sketch-and-descend method. 
Directly forming and manipulating the matrix $(\mat{S}^*\mat{B})^*(\mat{S}^*\mat{B})$ can lead to numerical issues, analogous to those involved in solving the normal equations (\cref{sec:normal-equations-unstable}).
To address this issues, my paper \cite{Epp24a} proposes implementing the iteration \cref{eq:ihs-optimized} by first computing a \QR decomposition (or other orthonormal factorization) of the sketched matrix
\begin{equation*}
    \mat{S}^*\mat{B} = \mat{Q}\mat{R},
\end{equation*}
and using the update rule
\begin{equation} \label{eq:ihs-3}
    \vec{x}_{i+1} \coloneqq \vec{x}_i + \alpha \mat{R}^{-1}(\mat{R}^{-*}(\mat{B}^*(\vec{c} - \mat{B}\vec{x}_i))) + \beta (\vec{x}_i - \vec{x}_{i-1}).
\end{equation}
The matrices $\mat{R}^{-1}$ and $\mat{R}^{-*}$ are applied by triangular substitution, not explicit inversion.

For numerical stability, it is important to evaluate \cref{eq:ihs-3} in the parenthesized order.
The following implementation
\begin{equation}
    \vec{x}_{i+1} \coloneqq \vec{x}_i + \alpha \mat{R}^{-1}(\mat{R}^{-*}(\mat{B}^*\vec{c} - \mat{B}^*(\mat{B}\vec{x}_i)))) + \beta (\vec{x}_i - \vec{x}_{i-1}), \tag{Bad!}
\end{equation}
equivalent to \cref{eq:ihs-3} in exact arithmetic, can be numerically disastrous.

As a refinement, we can choose the initial iterate $\vec{x}_0$ to be the sketch-and-solve solution (\cref{sec:sketch-and-solve}).
Using this good initialization can reduce the iteration count needed to convergence.
Code for the sketch-and-descend method using these improvements is provided in \cref{prog:sketch_descend}.
Evidence that these modifications improve stability appears in \cite[Fig.~3]{Epp24a}.

\begin{remark}[Sketch-and-descend versus sketch-and-precondition] \label{rem:descend-vs-precondition}
    The sketch-and-descend method with optimal parameters \cref{eq:optimal-coeffs} and the sketch-and-precondition method \warn{with sketch-and-solve initialization} (\cref{sec:spre-with-sketch-solve}) typically perform quite similarly in practice.
    As such, I usually recommend the sketch-and-precondition method in practice, as the algorithm has no free parameters.
    An exception where sketch-and-descend may be preferable is in parallel computing environments \cite{MSM14}, because sketch-and-descend requires fewer parallel synchronization steps.
\end{remark}

\myprogram{Sketch-and-descend method for solving overdetermined linear least-squares problems.}{}{sketch_descend}

\subsection{Analysis in exact arithmetic}
With an appropriate embedding and suitable parameter choices $\alpha,\beta$, the sketch-and-descend method converges geometrically, and the total runtime is $\order(mn \log (n/\varepsilon) + n^3\log n)$, the same as sketch-and-precondition.
Indeed, with the optimal parameters $\alpha,\beta$, the \emph{rate of convergence} is even the same for sketch-and-descent as for sketch-and-precondition.
We have the following result:

\begin{theorem}[Sketch-and-descend: Exact arithmetic]
    Let $\mat{S} \in \field^{d\times m}$ be a subspace embedding \warn{for $\flatonebytwo{\mat{B}}{\vec{c}}$} with distortion $\eta < 1$, and let $\vec{x}_i$ denote the iterates produced by the sketch-and-descend method (\cref{prog:sketch_descend}) \warn{in exact arithmetic}.
    For the trivial parameter choices $\alpha = 1$ and $\beta = 0$, we have
    \begin{equation*}
        \norm{\vec{c} - \mat{B}\vec{x}_k}^2 \le \left( 1 + 27\eta \cdot [(2+\sqrt{2})\eta]^{2k} \right) \norm{\vec{c} - \mat{B}\vec{x}}^2,
    \end{equation*}
    provided the convergence rate is less than one: $(2+\sqrt{2})\eta < 1$.
    With the optimal parameter choices \cref{eq:optimal-coeffs}, we have the improved bound
    \begin{equation*}
        \norm{\vec{c} - \mat{B}\vec{x}_k}^2 \le \left( 1 + \frac{256}{(1-\eta)^2} \cdot k\eta^{2k-1} \right) \norm{\vec{c} - \mat{B}\vec{x}}^2.
    \end{equation*}
    In particular, when implemented with a sparse sign embedding with Cohen's parameter choices (\cref{thm:ose-existence}), sketch-and-descend produces a $(1+\varepsilon)$-approximate least-squares solution
    \begin{equation*}
        \norm{\vec{c} - \mat{B}\vec{x}_k} \le (1+\varepsilon)\norm{\vec{x} - \mat{B}\vec{x}}.
    \end{equation*}
    with 99\% probability after $k = \order(\log(1/\varepsilon))$ iterations and $\order(mn \log(n/\varepsilon) + n^3\log n)$ operations.
\end{theorem}

The essence of this result appears in the works of Pilanci and coauthors \cite{PW16,OPA19,LP21}.
The precise result stated here is an immediate consequence of \cite[Thms.~3.1 \& B.1]{Epp24a}.

\subsection{Analysis in finite precision arithmetic}
In my paper \cite{Epp24a}, I proved that sketch-and-descend \warn{with trivial parameter choices $\alpha = 1$ and $\beta = 0$} is strongly forward stable.
We have the following result \cite[Thm.~4.1]{Epp24a}:

\begin{theorem}[Sketch-and-descend is strongly forward stable]
    Let $\mat{S} \in \field^{m\times d}$ be a subspace embedding \warn{for $\flatonebytwo{\mat{B}}{\vec{c}}$} of distortion \warn{$\eta \le 0.29$} with $d\le m$, suppose $\mat{B}$ is numerically full-rank $\kappa u \ll 1$, and assume the multiply operation is stable $\norm{\err(\mat{S}^*\vec{v})} \lesssim \norm{\vec{v}}u$.
    The sketch-and-descend method \warn{with trivial parameter choices $\alpha = 1$ and $\beta = 0$} is geometrically convergent until reaching strong forward stability
    \begin{equation} \label{eq:strong-forward-sketch-descend}
        \norm{\mat{B}(\vec{x} - \fl(\vec{x}_k))} \le 20\sqrt{\eta} [(2+\sqrt{2})\eta + \gamma\kappa u]^k + \gamma \left( \norm{\mat{B}}\norm{\vec{x}} + \kappa \norm{\vec{c} - \mat{B}\vec{x}} \right)u.
    \end{equation}
    Here, the prefactor $\gamma = \poly(m,n)$.
    In particular, for real data and implemented with a sparse sign embedding with Cohen's parameter settings (\cref{thm:ose-existence}), sketch-and-descend produces a strongly forward stable solution in $\order(mn \log(n/u) + n^3 \log n)$ operations.
\end{theorem}

The statement in \cite{Epp24a} is for real arithmetic, but the proof transfers to complex arithmetic without issue.
The bound \cref{eq:strong-forward-sketch-descend} on the residual error leads to a bound on the forward error $\norm{\vec{x} - \fl(\vec{x}_k)}$ by using the identity \cref{eq:y-to-By}.
We will omit proof of this result because our upcoming results for sketch-and-precondition are of greater practical relevance.

\section{Intermezzo: Lanczos, conjugate gradient, and LSQR} \label{sec:lanczos}

Error analysis of Krylov iterative methods in finite precision arithmetic remains a notoriously difficult and incomplete area of study.
Moreover, our theoretical understanding is uneven across different methods: The Lanczos algorithm has been extensively analyzed in finite precision, but virtually nothing is known about LSQR. 
Consequently, for analytical purposes, it is useful to study a version of sketch-and-precondition based on the Lanczos algorithm.

An introduction the Lanczos algorithm is beyond the scope of this thesis; see Tyler Chen's monograph \cite{Che24}.
The Lanczos algorithm is a Krylov subspace method that takes as input a \warn{Hermitian} matrix $\mat{M} \in \field^{n\times n}$ (or, more precisely, a subroutine implementing matrix--vector products $\vec{v} \mapsto \mat{M}\vec{z}$) and a vector $\vec{f} \in \field^n$.
It produces a \warn{Hermitian} tridiagonal matrix $\mat{T} \in \field^{k\times k}$ and a matrix $\mat{Q} \in \field^{n\times k}$ with orthonormal columns such that
\begin{itemize}
    \item The first column of $\mat{Q}$ is $\vec{q}_1 = \vec{f} / \norm{\vec{f}}$.
    \item For each $1 \le t \le k$, the leading $t$ columns of $\mat{Q}$ form a basis for the Krylov subspace $\operatorname{span} \{ \vec{f},\mat{M}\vec{f}, \ldots,\mat{M}^{t-1}\vec{f} \}$.
    \item The matrix $\mat{T} = \mat{Q}^*\mat{M}\mat{Q}$.
\end{itemize}
The computational cost of the Lanczos method is $k$ matvecs with $\mat{M}$ plus $\order(nk)$ additional storage and arithmetic.
See \cref{prog:lanczos}.

\myprogram{Lanczos algorithm for partially tridiagonalizing a matrix.}{}{lanczos}

The Lanczos method is frequently used to approximate the action $g(\mat{M})\vec{f}$ of a function of $\mat{M}$ on the vector $\vec{f}$.
Specifically, the \emph{Lanczos function approximation} is defined as
\begin{equation*}
    \hat{g(\mat{M})\vec{f}} \coloneqq \mat{Q}\, g(\mat{Q}^*\mat{M}\mat{Q}) \, \mat{Q}^*\vec{f} = \norm{\vec{f}} \cdot \mat{Q} g(\mat{T})\evec_1.
\end{equation*}
In particular, the choice $g(x) = x^{-1}$ yields the \emph{Lanczos linear solver}:
\begin{equation*}
    \hat{\mat{M}^{-1}\vec{f}} = \norm{\vec{f}} \cdot \mat{Q} \mat{T}^{-1}\evec_1.
\end{equation*}

We have the following equivalence:
\actionbox{\warn{In exact arithmetic}, the LSQR algorithm applied to the least-squares problem
\begin{equation*}
    \vec{x} = \argmin_{\vec{z} \in \field^n} \norm{\vec{c} - \mat{B}\vec{x}}
\end{equation*}
produces the same sequence of iterates $\vec{x}_0,\vec{x}_1,\ldots$ as either the conjugate gradient method or the Lanczos method applied to the normal equations
\begin{equation*}
    (\mat{B}^*\mat{B})\vec{x} = \mat{B}^*\vec{c}
\end{equation*}}
Given both the strong theoretical results for Lanczos and the equivalence between LSQR, CG, and Lanczos, it is natural to analyze versions of sketch-and-precondition that use Lanczos on the normal equations in place of LSQR.

\begin{remark}[Reorthogonalization]
    When the Lanczos algorithm is executed in finite-precision arithmetic, rounding errors leads the matrix $\mat{Q}$ to deviate from orthonormality.
    This loss of orthogonality can lead to significant problems for some applications of the Lanczos method, which may be mitigated using different types of \emph{reorthogonalization} strategies; see \cite[\S7.4--7.6]{Dem97} for a discussion of several options.
    Remarkably, for Lanczos function approximation and Lanczos linear solves, the Lanczos algorithm produces satisfactory results in finite precision arithmetic even without reorthogonalization \cite{DK91,DGK98,MMS18}; see \cite[Ch.~4]{Che24} for a discussion of this surprising phenomenon.
    All results in this thesis apply to the plain version of the Lanczos algorithm without any form of reorthogonalization.
\end{remark}



\section{Sketch-and-precondition with sketch-and-solve initialization} \label{sec:spre-with-sketch-solve}

After the initial release of Meier et al.'s manuscript, a curious thing was observed: Using the sketch-and-solve initialization greatly improves the numerical stability of sketch-and-precondition.
The published version of \cite{MNTW24} provides evidence that the \emph{forward error} of sketch-and-precondition with sketch-and-solve initialization is close to that of \QR-based direct methods (i.e., evidence of \emph{weak} forward stability).
My paper \cite{Epp24a} shows that the residual error is also comparable with \QR-based methods, providing evidence of \emph{strong} forward stability.
For difficult problems, however, sketch-and-precondition with sketch-and-solve initialization is still not to backward stable; see \cref{fig:rand-ls-stable} ahead.

In \cite{MNTW24,Epp24a}, the stabilizing effects of the sketch-and-solve initialization are presented as empirical findings, and no rigorously explanation was available.
Maike Meier, Yuji Nakatsukasa, and myself provided such an explanation by proving the following result \cite[Thm.~8.2]{EMN24}:

\begin{theorem}[Sketch-and-precondition: Strong forward stability]
    Assume we work over the \warn{real} field $\field = \real$.
    Let $\mat{S} \in \real^{m\times d}$ be a subspace embedding \warn{for $\flatonebytwo{\mat{B}}{\vec{c}}$} of distortion \warn{$\eta \le 0.9$} and with $d \le m$, suppose $\mat{B}$ is numerically full-rank $\kappa u \ll 1$, and assume the multiply operation is stable $\norm{\err(\mat{S}^*\vec{v})} \lesssim \norm{\vec{v}}u$.
    Then the \warn{Lanczos-based} implementation of sketch-and-precondition with sketch-and-solve initialization produces a strongly forward stable solution to the least-squares problem \cref{eq:least-squares} after $\order(\log(1/u))$ iterations.
    To produce such a solutionwith a sparse sign embedding with Cohen's parameter settings (\cref{thm:ose-existence}), sketch-and-precondition requires at most $\order(mn \log(n/u) + n^3 \log n)$ operations.
\end{theorem}

A proof of a simplified version of this result is provided in \cref{ch:stability-proofs}.

\section{Sketch-and-precondition with iterative refinement} \label{sec:spir}

In matrix computations, \emph{iterative refinement} is a standard approach to improving the quality of a computed solution \cite[Ch.~12 \& \S20.5]{Hig02}.
In its basic form, iterative refinement executes the following sequence of steps:
\begin{equation} \label{eq:iterative-refinement}
    \begin{split}
    \vec{x}_1 &\coloneqq \Call{NotFullyStableSolver}{\mat{B},\vec{c}}; \\
    \vec{x}_2 &\coloneqq \vec{x}_1 + \Call{NotFullyStableSolver}{\mat{B},\vec{c} - \mat{B}\vec{x}_1}.
    \end{split}
\end{equation}
Here, \Call{NotFullyStableSolver}{} is any subroutine producing an approximate least-squares solution: $\Call{NotFullyStableSolver}{\mat{B},\vec{c}} \approx \mat{B}^\dagger \vec{c}$.
Equation \cref{eq:iterative-refinement} describes one step of iterative refinement.
In general, multiple steps of iterative refinement may be necessary to achieve full accuracy.
(That is, we continue on as $\vec{x}_3 \coloneqq \vec{x}_2 + \Call{NotFullyStableSolver}{\mat{B},\vec{c} - \mat{B}\vec{x}_2}$, etc.).

The most conventional use of iteration refinement is \emph{mixed-precision iterative refinement}.
Here, the \Call{NotFullyStableSolver}{} routine is furnished by a factorization of $\mat{B}$ computed using a lower numerical precision (i.e., single precision):
\begin{equation*}
    (\mat{Q},\mat{R}) \gets \Call{LowPrecisionQR}{\mat{B}}; \:\: \Call{NotFullyStableSolver}{\mat{B},\vec{c}} = \mat{R}^{-1}(\mat{Q}^*\vec{c}).
\end{equation*}
In this thesis, we will perform iterative refinement in a \emph{uniform numerical precision} $u$ (e.g., double precision) and use sketch-and-precondition as the \Call{NotFullyStableSolver}{} subroutine.

This motivates us to consider \emph{sketch-and-precondition with iterative refinement} (SPIR, \cite{EMN24}).
We begin with the sketch-and-solve initialization 
\begin{equation*}
    \vec{x}_0 \coloneqq (\mat{S}^*\mat{B})^\dagger (\mat{S}^*\vec{c}).   
\end{equation*}
Then, we produce a first approximate solution by running sketch-and-precondition with initialization $\vec{x}_0$:
\begin{equation*}
    \vec{x}_1 \coloneqq \Call{SketchAndPrecondition}{\mat{B},\vec{c},\texttt{initialization} = \vec{x}_0},
\end{equation*}
Then, we apply a single step of iterative refinement
\begin{subequations}
\begin{equation} \label{eq:spir-traditional}
    \vec{x}_2 \coloneqq \vec{x}_1 + \Call{SketchAndPrecondition}{\mat{B},\vec{c} - \mat{B}\vec{x}_1,\texttt{initialization} = \vec{0}}.
\end{equation}
Perhaps confusingly, it is imperative that we use the \warn{zero initialization} for the iterative refinement step.
Alternatively, rather than the traditional refinement step \cref{eq:spir-traditional}, we can use the following \emph{restarting refinement step}:
\begin{equation} \label{eq:spir-restarting}
    \vec{x}_2 \coloneqq \vec{x}_1 + \Call{SketchAndPrecondition}{\mat{B},\vec{c},\texttt{initialization} = \vec{x}_1}.
\end{equation}
\end{subequations}
In our experience, the traditional \cref{eq:spir-traditional} and restarting \cref{eq:spir-restarting} refinement steps have similar stability properties.
An implementation of SPIR appears in \cref{prog:spir}.

\myprogram{Sketch-and-precondition with iterative refinement for solving overdetermined linear least-squares problems.}{}{spir}

The backward stability of SPIR may be surprising in view of the classic work of Golub and Wilkinson \cite{GW66}, who showed that standard \emph{mixed-precision} iterative refinement \cref{eq:iterative-refinement} does not greatly substantially improve the numerical accuracy for least-squares problems except for nearly consistent systems.
The setting of SPIR is differs from Golub and Wilkinson's setting. 
SPIR performs iterative refinement in \emph{uniform precision}, while using preconditioned Krylov methods (appropriately initialized) as the solvers.
This change in setting explains the difference in stability behavior.

Unfortunately, we do not have a \emph{completely} unconditional proof of backward stability for SPIR.
We need one small assumption:

\begin{definition}[Non-pathological rounding error assumption] \label{def:non-pathological}
    A vector $\hatvector{x} \in \field^n$ satisfies the \emph{non-pathological rounding error assumption} if
    \begin{equation*}
        \norm{\hatvector{x}} \gtrsim \norm{\vec{x}} + \kappa u \left(\norm{\vec{x}} + \kappa \frac{\norm{\vec{c} - \mat{B}\vec{x}}}{\norm{\mat{B}}} \right).
    \end{equation*}
\end{definition}

The Golub--Wilkinson--Wedin theorem (\cref{fact:ls-perturb}) shows that the errors for a strongly forward or backward approximate least-squares solution have magnitude
\begin{equation*}
    \norm{\hatvector{x} - \vec{x}} \lesssim \kappa u \left((\norm{\vec{x}} + \kappa \frac{\norm{\vec{c} - \mat{B}\vec{x}}}{\norm{\mat{B}}} \right).
\end{equation*}
Moreover, this bound is approximately attained for ``most'' backward perturbations $\mat{\Delta B}$ and $\vec{\Delta c}$ \cite[\S5]{GW66}.
So the non-pathological rounding error assumption just states that (1) the error $\hatvector{x} - \vec{x}$ is of the typical size expected for a backward stable method (that is, $\norm{\smash{\hatvector{x} - \vec{x}}}$ is not abnormally small) and (2) the errors $\hatvector{x} - \vec{x}$ do not catastrophically cancel with the true solution $\vec{x}$ in such a way that $\hatvector{x}$ has an unusually small norm.

With this assumption, we state the following backward stability theorem:

\begin{theorem}[SPIR: Backward stability] \label{thm:spir-backward}
    Assume we work over the \warn{real} field $\field = \real$.
    Let $\mat{S} \in \real^{m\times d}$ be a subspace embedding \warn{for $\flatonebytwo{\mat{B}}{\vec{c}}$} of distortion \warn{$\eta \le 0.9$} and with $d \le m$, suppose $\mat{B}$ is numerically full-rank $\kappa u \ll 1$, and assume the multiply operation is stable $\norm{\err(\mat{S}^*\vec{v})} \lesssim \norm{\vec{v}}u$.
    Finally, assume that the computed solution satisfies the non-pathological rounding error assumption.
    Then the \warn{Lanczos-based} implementation of SPIR with traditional refinement \cref{eq:spir-traditional} produces a backward stable solution to the least-squares problem \cref{eq:least-squares} after $\order(\log(1/u))$ total iterations.
    In particular, implemented with a sparse sign embedding with Cohen's parameter settings (\cref{thm:ose-existence}), the runtime is $\order(mn \log(n/u) + n^3 \log n)$ operations.
\end{theorem}

\Cref{ch:stability-proofs} provides a proof of a simplified version of this result.

\begin{remark}[FOSSILS]
    Our paper \cite{EMN24} also introduces FOSSILS, a backward stable version of iterative Hessian sketch/sketch-and-descend.
    While there are a few scenarios when one might prefer FOSSILS over SPIR (cf.\ \cref{rem:descend-vs-precondition}), SPIR is generally the preferable method in practice.
    See \cite[\S4.4]{EMN24} for discussion.
\end{remark}

\begin{remark}[(Numerically) rank-deficient problems]
    For problems that are numerically rank-deficient ($\kappa u \gtrapprox 1$), the sketch-and-precondition algorithm does not work properly.
    In this case, there are three options:
    \begin{enumerate}
        \item \textbf{\textit{``Smoothed analysis.''}} 
        One way of addressing numerical rank-deficiency is to add a small independent random perturbation to each matrix entry (say, a centered Gaussian with standard deviation $10u\norm{\mat{B}}_{\mathrm{F}}$). 
        This operation has a regularizing effect on a matrix, causing it to become numerically full rank with high probability.
        This idea was proposed in the context of randomized least-squares solvers by Meier et al.\ \cite{MNTW24}.
        The idea of perturbing the input to an algorithm to avoid bad instances is known as smoothed analysis \cite{ST01}, and it has established itself as a core tool in the design of randomized matrix algorithms \cite[\S6.3]{KT24}.
        A disadvantage of this approach is that it destroys any benevolent structure that may be present in the matrix $\mat{B}$, such as sparsity.
        \item \textbf{\textit{Ridge regularization.}}
        Another way of regularizing a least-squares problem is to add a small ridge-regularization penalty to the least-squares objective:
        \begin{equation*}
            \vec{x}_{\mathrm{reg}} = \argmin_{\vec{z} \in \field^n} \norm{\vec{c} - \mat{B}\vec{z}}^2 + \mu^2 \norm{\vec{z}}^2.
        \end{equation*}
        The ridge-regularized least-squares problem is equivalent to an ordinary least-squares problem involving an augmented matrix:
        \begin{equation*}
            \vec{x}_{\mathrm{reg}} = \argmin_{\vec{z} \in \field^n} \norm{\twobyone{\vec{c}}{\vec{0}} - \twobyone{\mat{B}}{\mu \Id} \vec{z}}.
        \end{equation*}
        To solve rank-deficient problems, we can add a small amount of ridge regularization (say, $\mu = 10u\norm{\mat{B}}_{\mathrm{F}}$).
        See \warn{arXiv version 2} of \cite{EMN24} for details.
        \item \textbf{\textit{Truncation.}}
        A third approach is to truncate the least-squares problem by computing an SVD $\mat{S}^*\mat{B} = \mat{U}\mat{\Sigma}\mat{V}^*$ of the sketched matrix.
        Then, instead of using the standard-inverse preconditioner $\mat{M}^{-1} = \mat{V}\mat{\Sigma}^{-1}$, one executes preconditioned LSQR with a \emph{rectangular} inverse-preconditioner $\mat{P} \coloneqq \mat{V}(:,1:r)\mat{\Sigma}(1:r,1:r)^{-1}$, where $1\le r \le n$ is the number of singular values of $\mat{S}^*\mat{B}$ above a level, say, $10 \sigma_{1}u$.
        Truncating in this way has the effect of restricting the solution to lie in the range of $\mat{V}(:,1:r)$, a subspace on which the linear transformation $\mat{B}$ is numerically full-rank.
        See \cite{EMN24} for details.
    \end{enumerate}
    %
\end{remark}

\section{Experiments} \label{sec:stable-ls-experiments}

\begin{figure}
    \centering
    \includegraphics[width=0.49\textwidth]{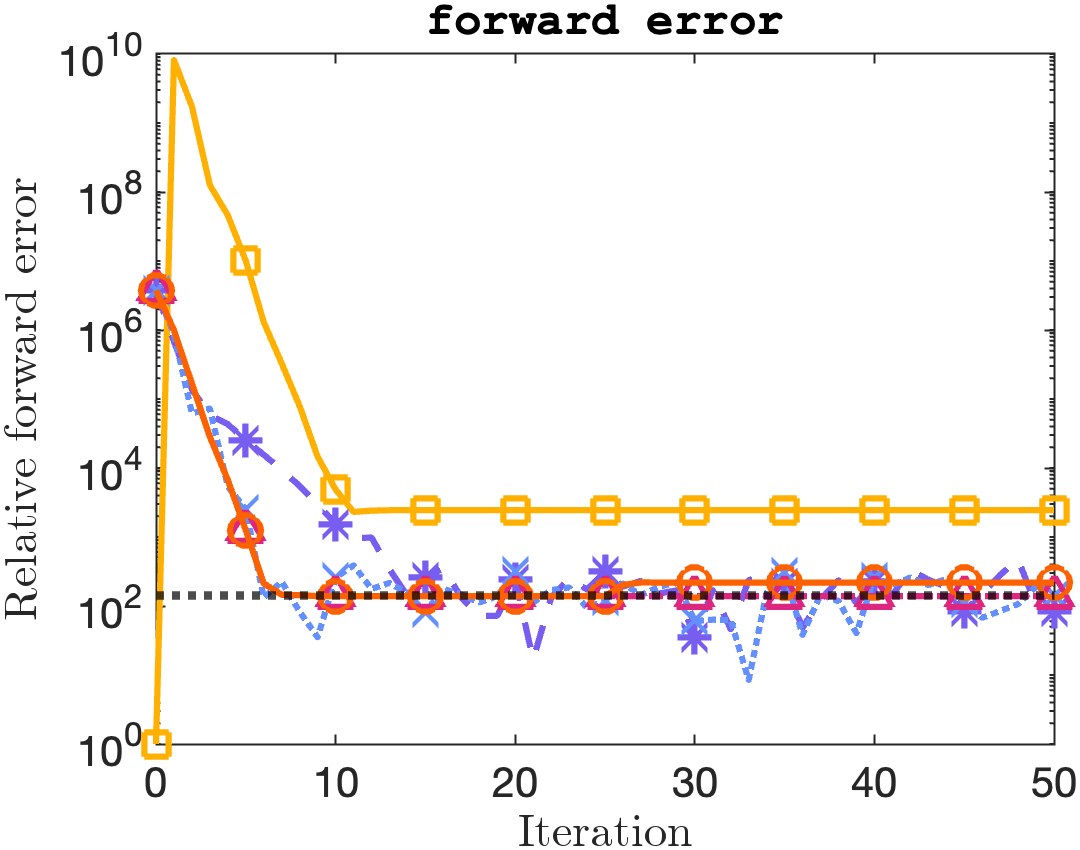}
    \includegraphics[width=0.49\textwidth]{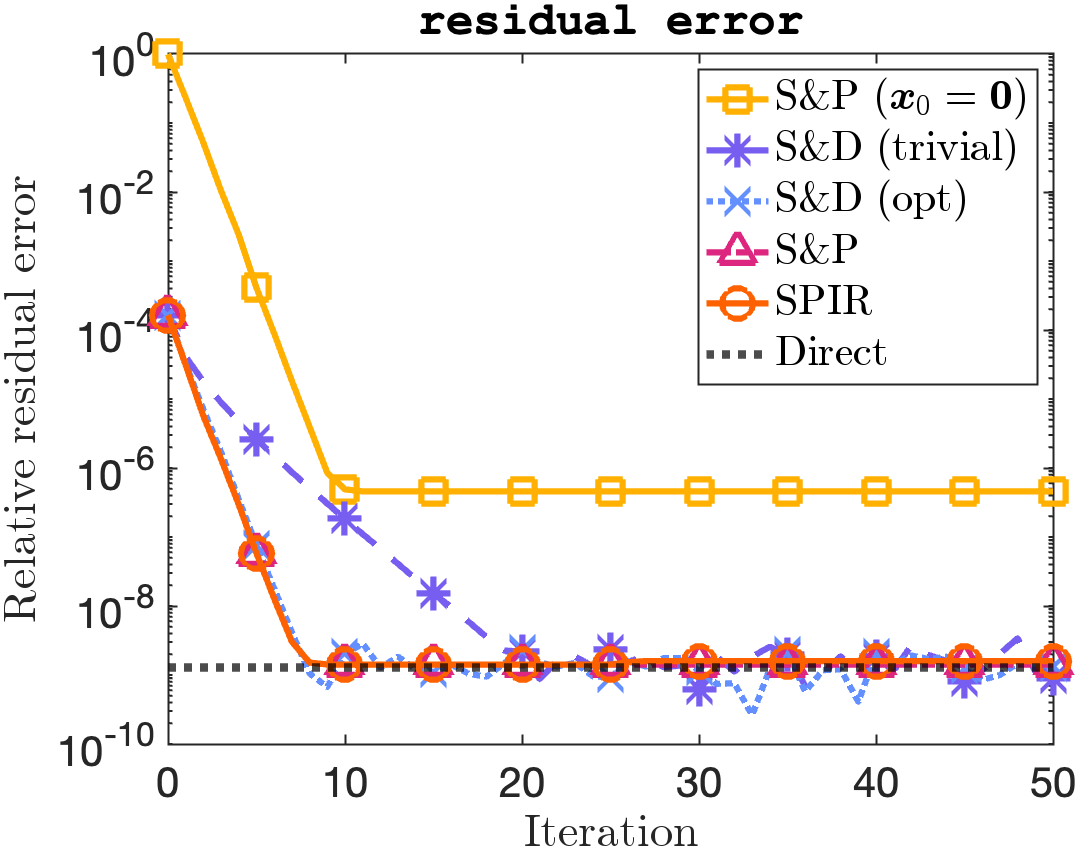}

    \includegraphics[width=0.49\textwidth]{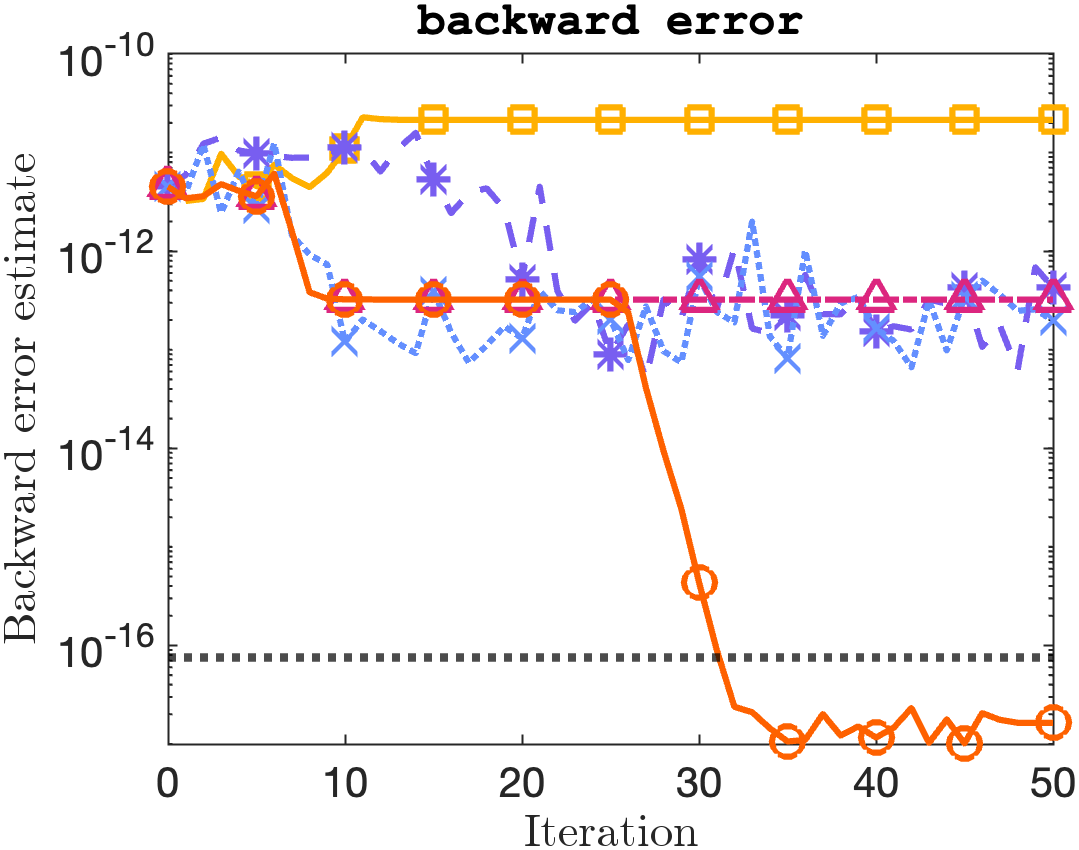}

    \caption[Forward, residual, and backward error for different randomized least-squares solvers]{Forward (\emph{top left}), residual (\emph{top right}), and backward (\emph{bottom}) errors for five randomized least-squares solvers: sketch-and-precondition with zero initialization (yellow solid squares), sketch-and-descend with trivial parameters $\alpha = 1$, $\beta = 0$ (purple dashed asterisks), sketched and descend with optimized parameters \cref{eq:optimal-coeffs} (blue dotted crosses), sketch-and-precondition with sketch-and-solve initialization (pink dash-dotted triangles), and SPIR (restarted after 25 iterations, orange solid circles).
    We present a single execution.
    Errors for direct method (MATLAB's \texttt{mldivide}) shown for reference.}
    \label{fig:rand-ls-stable}
\end{figure}

\Cref{fig:rand-ls-stable} compares the accuracy, stability, and convergence rates of different randomized least-squares problems.
We use the same problem instance from \cref{fig:sketch-and-precondition-unstable}, and we use a sparse sign embedding with embedding dimension $d\coloneqq 20n$ and sparsity level $\zeta = 4$ for all methods.
(A larger embedding dimension is necessary to ensure convergence of sketch-and-descend with trivial parameter choices.)
We see that all methods, except sketch-and-precondition with the zero initialization, are strongly forward stable.
The convergence rate of sketch-and-descend with the trivial parameter settings is slower than the other methods; all other methods have the same geometric rate of convergence.
We see that all methods except SPIR fail to be backward stable, with the backward error stagnating far above the level of the direct method.
SPIR, with refinement performed starting at iteration 25, converges to full backward stability.
In fact, SPIR achieves a slightly smaller backward error than even the direct method.

\section{The backward error and its estimation} \label{sec:backerr-estimate}

As stated in \cref{def:backward-stability}, an approximate least-squares solution $\hatvector{x} \in \field^n$ is backward stable if there are small perturbations $\vec{\Delta c}$ and $\mat{\Delta B}$, of unit roundoff size, such that $\hatvector{x}$ is the exact solution of the perturbed system
\begin{equation*}
    \hatvector{x} = \argmin_{\vec{z} \in \field^n} \norm{(\vec{c} + \vec{\Delta c}) - (\mat{B}+\mat{\Delta B})\vec{z}}.
\end{equation*}
The minimal size of these perturbations is called the backward error.

\begin{definition}[Backward error]
    Introduce a parameter $\theta \in [0,+\infty]$.
    The \emph{backward error} is 
    \begin{equation*}
        \BE_\theta(\hatvector{x}) \coloneqq \min \big\{ \|\flatonebytwo{\mat{\Delta B}}{\theta \cdot \vec{\Delta c}}\|_{\mathrm{F}} : \hatvector{x} = \argmin_{\vec{z} \in \field^n} \norm{(\vec{c} + \vec{\Delta c}) - (\mat{B}+\mat{\Delta B})\vec{z}} \big\}.
    \end{equation*}
\end{definition}

The most natural values for $\theta$ are $\theta = +\infty$ (i.e., $\vec{\Delta c}$ is required to be set to zero) or, provided $\mat{B}$ and $\vec{c}$ are commensurately scaled, $\theta = 1$ (equal weight placed on both $\mat{\Delta B}$ and $\vec{\Delta c}$).
Assuming the normalization $\norm{\mat{B}} = \norm{\vec{c}} = 1$, an approximate least-squares solution $\hatvector{x}$ is backward stable if and only if $\BE_1(\hatvector{x}) \lesssim u$.

To deploy randomized least-squares solvers, it is desirable to have \emph{estimates} of the backward error that are computable at runtime.
These error estimates can be used as \emph{stopping criteria} for the algorithm.

The question of how to compute or estimate the backward error of an approximate least-squares solution has a surprisingly deep history.
As Higham writes \cite[\S20.7]{Hig02}, ``Although it has been known since the 1960s that a particular method for solving the LS problem, namely the Householder QR factorization method, yields a small normwise backward error\ldots it was for a long time an open problem to obtain a formula for the backward error of an arbitrary approximate solution.''
This problem was resolved in 1995 by Wald\'en, Karlson, and Sun \cite{WKS95}, who discovered a formula for it.
Higham developed a numerically stable version of Wald\'en, Karlson, and Sun's formula.
Unfortunately, the Wald\'en--Karlson--Sun--Higham formula requires a steep $\order(m^3)$ operations to evaluate, making it impractical for many use cases.

\myparagraph{The Karlson--Wald\'en estimate}
The great cost of the Wald\'en--Karlson--Sun--Higham formula inspired a flurry of research effort over two decades to develop and analyze computationally efficient \emph{estimates} for the backward error.
This line of work was initiated by Karlson and Wald\'en \cite{KW97} and refined over a series of follow-up works \cite{Gu98,Grc03,GSS07,GJT12}.
These papers converged on estimates of the following form \cite[\S2]{GSS07}:
\begin{equation*}
    \hat{\BE}_\theta(\hatvector{x}) \coloneqq \frac{\theta}{(1 + \theta^2 \norm{\hatvector{x}}^2)^{1/2}} \norm{ \left( \mat{B}^*\mat{B} + \frac{\theta^2 \norm{\vec{c} - \mat{B}\hatvector{x}}^2}{1 + \theta^2 \norm{\hatvector{x}}^2} \right)^{-1/2} \mat{B}^*(\vec{c} - \mat{B}\hatvector{x}) }.
\end{equation*}
Following previous papers \cite{GJT12}, we shall this the \emph{Karlson--Wald\'en estimate} for the backward error.
This estimate can be evaluated in $\order(mn^2)$ operations.
Given a \QR factorization or SVD of $\mat{B}$, then $\order(mn)$ operations suffice \cite{GSS07}.
An implementation of the Karlson--Wald\'en estimate appears in \cref{prog:backerr_est}.

\myprogram{Karlson--Wald\'en estimate of the backward error for an least-squares solution.}{}{backerr_est}

The series of papers \cite{Gu98,Grc03,GSS07} gave increasingly sharp characterizations of the quality of the Karlson--Wald\'en estimate.
Gratton, Jir\'anek, and Titley-Peloquin conclusively resolved this line of research, providing a sharp analysis of the estimate \cite{GJT12}.
Here is a simplified version of their result:
\begin{fact}[Karlson--Wald\'en estimate] \label{fact:kw-estimate}
    Over the \warn{real} field $\field = \real$, the Karlson--Wald\'en estimate is sharp within a factor of $\sqrt{2}$:
    \begin{equation} \label{eq:GJT}
        1\cdot \hat{\BE}_\theta(\hatvector{x}) \le \BE_\theta(\hatvector{x}) \le \sqrt{2} \cdot \hat{\BE}_\theta(\hatvector{x})
    \end{equation}
    for every $\theta \in [0,+\infty]$ and $\hatvector{x} \in \real^n$.
    The constants $1$ and $\sqrt{2}$ in \cref{eq:GJT} cannot be improved.
\end{fact}

\myparagraph{Improving the Karlson--Wald\'en estimate with sketching}
The cost of $\order(mn^2)$ operations for the Karlson--Wald\'en estimate is still prohibitive for deployment with randomized least-squares solvers, since our goal is to achieve runtimes of roughly $\order(mn + n^3)$ operations.
Fortunately, we can accelerate the Karlson--Wald\'en estimate with sketching.
My coauthors and I in \cite{EMN24} proposed the \emph{sketched Karlson--Wald\'en estimate}:
\begin{equation*}
    \hat{\BE}_{\mathrm{sk},\theta}(\hatvector{x}) \coloneqq \frac{\theta}{(1 + \theta^2 \norm{\hatvector{x}}^2)^{1/2}} \norm{ \left( (\mat{S}^*\mat{B})^*(\mat{S}^*\mat{B}) + \frac{\theta^2 \norm{\vec{c} - \mat{B}\hatvector{x}}^2}{1 + \theta^2 \norm{\hatvector{x}}^2} \right)^{-1/2} \mat{B}^*(\vec{c} - \mat{B}\hatvector{x}) }.
\end{equation*}
Here, $\mat{S}$ is a subspace embedding.
We have the following analysis for the sketched Karlson--Wald\'en estimate, adapted from \cite[Prop.~4.2]{EMN24}:

\begin{proposition}[Sketched Karlson--Wald\'en estimate]
    Consider the \warn{real} field $\field = \real$, let $\theta \in [0,+\infty]$, $\hatvector{x} \in \real^n$, and let $\mat{S} \in \real^{m\times d}$ be a subspace embedding with distortion $\eta$.
    The sketched Karlson--Wald\'en estimate is sharp up to the following bounds:
    \begin{equation*}
        (1-\eta)\cdot \hat{\BE}_\theta(\hatvector{x}) \le \BE_\theta(\hatvector{x}) \le \sqrt{2}(1+\eta) \cdot \hat{\BE}_\theta(\hatvector{x}).
    \end{equation*}
\end{proposition}

To compute the sketched Karlson--Wald\'en estimate, we can modify \cref{prog:backerr_est} to use the SVD of the sketched $\mat{B}$ matrix, $\mat{S}^*\mat{B} = \mat{U}\mat{\Sigma}\mat{V}^*$, in place of the original $\mat{B}$ matrix.
With a sparse sign embedding (\cref{def:sparse-sign-body}), the estimator requires a one-time set-up cost of roughly $\order(mn + n^3)$ operations, after which the cost of computing the estimator is $\order(mn)$ operations.

\myparagraph{Adaptively stopping sketch-and-precondition and SPIR}
We can use the sketched Karlson--Wald\'en estimate to adaptively stop SPIR when backward stability has been achieved.
We can also use automatic methods to switch between the two stages of SPIR.
See \cite[\S4]{EMN24} for details.

\chapter{Proofs of stability} \label{ch:stability-proofs}

\epigraph{Once in a lifetime a user of computer arithmetic should examine the details of a backward error analysis.}{Beresford N.\ Parlett, \textit{The Symmetric Eigenvalue Problem}, \cite[\S2.6.1]{Par98}}

\epigraph{The two main classes of rounding error analysis are not, as my audience might imagine, `backwards' and `forwards', but rather `one's own' and `other people's'. One's own is, of course, a model of lucidity; that of others serves only to obscure the essential simplicity of the matter in hand.}{James Wilkinson, NAG 1984 Annual General Meeting \cite{Wil85}}

In this chapter, we will provide stability analysis for sketch-and-precondition and sketch-and-precondition with iterative refinement (SPIR).
As Wilkinson's quote at the beginning of the section highlights, every numerical analyst has their own approach to rounding error analysis.
My approach more qualitative, and I prefer to avoid carefully tracking the size of constants and prefactors.
Even at this level of granularity, the numerical stabilility analysis of randomized iterative methods involves lengthy calculations and estimates.

For this thesis, I will provide the simplest argument for strong forward stability of sketch-and-precondition and backward stability of SPIR that I am aware of.
To make the analysis as simple as possible, this thesis will develop first-order perturbation estimates.
The original paper \cite{EMN24} contains detailed analysis which accounts for the higher-order terms.
Throughout this section, important intermediate results are boxed:
\importantresult{\begin{equation*}
    \text{I am an important result.}
\end{equation*}}

\myparagraph{Sources}
This chapter is adapted from the following paper:

\fullcite{EMN24}.

\myparagraph{Outline}
\Cref{sec:stability-proofs-assumptions} introduces standing assumptions and notation for this section.
\Cref{sec:sketching-QR-factorizing,sec:sketch-and-solve-numerical} analyze the floating-point errors incurred during sketching.
Next, in \cref{sec:stability-multiplication-solves,sec:stability-interleaved}, we analyze the stability of multiplications by the matrices $\mat{B}$ and $\mat{R}^{-1}$, as well as their adjoints.
We take a brief pause from analysis in \cref{sec:lanczos-stability}, which describes existing stability results for Lanczos linear solves.
\Cref{sec:error-formula} contains the core of the stability analysis, establishing a general formula the floating point errors introduced in a single iterative refinement step with sketch-and-precondition.
As a consequence of this general formula, we derive strong forward stability for sketch-and-precondition (with sketch-and-solve initialization) in \cref{sec:forward-stability-spre-proof} and backward stability of SPIR in \cref{sec:backward-stability-spir-proof}.

\section{Standing assumptions and more notation} \label{sec:stability-proofs-assumptions}

Throughout this section, we will assume the normalization 
\begin{equation*}
    \norm{\mat{B}} = \norm{\vec{c}} = 1.
\end{equation*}
As in \cref{ch:stable-ls}, we will denote $\kappa \coloneqq \cond(\mat{B})$, and we will use the standing assumption that $\kappa u \ll 1$.

A vector $\vec{z} \in \field^n$ is \emph{exactly represented} if $\fl(\vec{z}) = \vec{z}$.
Throughout, $\vec{e}$, $\vec{e}'$, etc.\ denote arbitrary vectors of norm $\lesssim u$ (likewise, for matrices $\mat{E}$, $\mat{E}'$, etc.).
We will freely interchange the value of $\vec{e}$ from line-to-line, even in the same line (i.e., we could write $\vec{e} + \vec{e}' = \vec{e}$).
We will use $\mathrm{HOT}$ (``higher order terms'') to denote any term proportional to $u^2$, suppressing all prefactors depending on $m$, $n$, the norms of various objects, and $\kappa$.
So, for our purposes $\norm{\vec{z}}\kappa^{100} u^2 = \mathrm{HOT}$.
The ``$\mathrm{HOT}$'' notation helps to streamline stability arguments to their essence, at the cost of losing precise information about the size of higher-order terms.

\section{Sketching and \QR factorizing} \label{sec:sketching-QR-factorizing}

Our first order of business will be the assess the quality of the randomized preconditioner $\mat{R}$ in the presence of rounding errors.

First, we investigate the rounding errors incurred in sketching and computing the \QR decomposition.
We have assumed the sketching operation is numerically stable in the sense that
\begin{equation*}
    \fl(\mat{S}^*\mat{B}) = \mat{S}^*\mat{B} + \mat{E}.
\end{equation*}
\QR factorization is backward stable \cite[Thm.~19.4]{Hig02}, so there exists a matrix $\mat{U}$ for which
\begin{equation} \label{eq:sketching-qr-errors}
    \mat{U} \fl(\mat{R}) = \fl(\mat{S}^*\mat{B}) + \mat{E} = \mat{S}^*\mat{B} + \mat{E}.
\end{equation}

Now, we use \cref{eq:sketching-qr-errors} to assess the quality of the \emph{numerically computed} $\mat{R}$ matrix as a preconditioner for $\mat{B}$.
By the subspace embedding property, we have
\begin{multline*}
    \norm{\mat{B} \fl(\mat{R})^{-1}} \le \frac{1}{1-\eta} \cdot \norm{\mat{S}^*\mat{B} \fl(\mat{R})^{-1}} \\ = \frac{1}{1-\eta} \cdot \norm{\mat{U} + \mat{E}\fl(\mat{R})^{-1}} \le \frac{1}{1-\eta} + \cdot \frac{1}{1-\eta} \frac{\norm{\mat{E}}}{\sigma_{\mathrm{min}}(\fl(\mat{R}))}.
\end{multline*}
To bound the minimum singular value of $\fl(\mat{R})$, we compute
\begin{multline} \label{eq:R-minsing-numerical}
    \sigma_{\mathrm{min}}(\fl(\mat{R})) = \sigma_{\mathrm{min}}(\mat{S}^*\mat{B} + \mat{E}) \ge \sigma_{\mathrm{min}}(\mat{S}^*\mat{B}) - \norm{\mat{E}} \\ \ge \frac{1}{1+\eta} \cdot \sigma_{\mathrm{min}}(\mat{B}) - \norm{\mat{E}} \ge \frac{1}{(1+\eta)\kappa} - \order(u) \ge \frac{\mathrm{const}}{\kappa}.
\end{multline}
In the last inequality, we used the hypothesis $\kappa u \ll 1$ (so that $u \ll 1/\kappa$).
Therefore,
\begin{equation} \label{eq:BR-norm-numerical}
    \norm{\mat{B} \fl(\mat{R})^{-1}} \le \frac{1}{1-\eta} + \order(\kappa u) \le \mathrm{const}.
\end{equation}
A similar argument shows
\begin{equation} \label{eq:R-maxsing-numerical}
    \norm{\mat{R}} \le \mathrm{const}.
\end{equation}
and
\begin{equation} \label{eq:BR-minsing-numerical}
    \sigma_{\mathrm{min}}(\mat{B} \fl(\mat{R})^{-1}) \ge \frac{1}{1+\eta} - \order(\kappa u) \ge \mathrm{const} > 0.
\end{equation}
Thus, the numerically computed matrix $\fl(\mat{R})$ is a good preconditioner for $\mat{B}$, satisfying the following bound:
\importantresult{\begin{equation*} 
    \cond(\mat{B}\fl(\mat{R})^{-1}) = \cond(\mat{B}\mat{R}^{-1}) + \order(\kappa u) \le \mathrm{const} < 1.
\end{equation*}}
To make our lives easier, we make the following notational affordance: \warn{Going forward, the symbol $\mat{R}$ will denote the \emph{numerically computed} preconditioner $\mat{R}$.}

\section{The sketch-and-solve solution} \label{sec:sketch-and-solve-numerical}

Now, we analyze the numerically computed sketch-and-solve solution.
Householder \QR factorization is a backward stable for solving a least-squares problem \cite[Thm.~20.3]{Hig02}.
Therefore, the numerically computed sketch-and-solve solution $\vec{x}_0 = (\mat{S}^*\mat{B})^\dagger (\mat{S}^*\vec{c})$ is
\begin{equation*}
    \fl(\vec{x}_0) = (\mat{S}^*\mat{B} + \mat{E})^\dagger (\mat{S}^*\vec{c})
\end{equation*}
To first order, the pseudoinverse of a perturbation of a matrix $\mat{F}$ is
\begin{equation*}
    (\mat{F} + \mat{E})^\dagger = \mat{F}^\dagger + (\mat{F}^*\mat{F})^{-1} \mat{E}^*(\Id - \mat{F}\mat{F}^\dagger) - \mat{F}^\dagger \mat{E} \mat{F}^\dagger + \mathrm{HOT}.
\end{equation*}
Instantiating this result with $\mat{F} \coloneqq \mat{S}^*\mat{B}$ and invoking the identity $\vec{x}_0 = (\mat{S}^*\mat{B})^\dagger (\mat{S}^*\vec{c})$, we obtain
\begin{multline*}
    \norm{\mat{B}(\vec{x} - \fl(\vec{x}_0))} = \norm{\mat{B}(\vec{x} - \vec{x}_0)} + \norm{\mat{B}(\mat{B}^*\mat{S}\mat{S}^*\mat{B})^{-1}\mat{E}^*\mat{S}^*(\vec{c} - \mat{B}\vec{x}_0)} \\ + \norm{\mat{B}(\mat{S}^*\mat{B})^\dagger \mat{E} \vec{x}_0} + \mathrm{HOT}.
\end{multline*}
By the analysis of sketch-and-solve (\cref{thm:sketch-and-solve}), the residual error is bounded: $\norm{\mat{B}(\vec{x} - \vec{x}_0)} \lesssim \norm{\vec{c} - \mat{B}\vec{x}}$.
Now, invoking the subspace embedding property and the submultiplicative property of the spectral norm, we have
\begin{multline*}
    \norm{\mat{B}(\vec{x} - \fl(\vec{x}_0))} \lesssim \norm{\vec{c} - \mat{B}\vec{x}} + \norm{
    \smash{\mat{S}^*\mat{B}(\mat{B}^*\mat{S}\mat{S}^*\mat{B})^{-1}}} \norm{\mat{E}} \norm{\vec{c} - \mat{B}\vec{x}} \\ + \norm{\smash{\mat{S}^*\mat{B}(\mat{S}^*\mat{B})^\dagger}} \norm{\mat{E}} \norm{\vec{x}_0} + \mathrm{HOT}.
\end{multline*}
We have $\mat{S}^*\mat{B}(\mat{B}^*\mat{S}\mat{S}^*\mat{B})^{-1} = (\mat{S}\mat{B})^{\dagger *}$, which has norm $\lesssim \kappa$.
Thus,
\begin{equation*}
    \norm{\mat{B}(\vec{x} - \fl(\vec{x}_0))} \lesssim \norm{\vec{c} - \mat{B}\vec{x}} + \norm{\vec{x}_0} u + \mathrm{HOT}.
\end{equation*}
Finally, we note that $\norm{\vec{x}_0} u \le \norm{\vec{x}} u + \norm{\vec{x} - \vec{x}_0}u = \norm{\vec{x}} u + \mathrm{HOT}$.
We conclude
\importantresult{\begin{equation} \label{eq:sketch-and-solve-finite-precision}
    \norm{\mat{B}(\vec{x} - \fl(\vec{x}_0))} \lesssim \norm{\vec{c} - \mat{B}\vec{x}} + \norm{\vec{x}} u + \mathrm{HOT}.
\end{equation}}
The residual error of the numerically computed sketch-and-solve solution is bounded by a multiple of the residual $\norm{\vec{c} - \mat{B}\vec{x}}$ plus $\norm{\vec{x}} u$.

\section{Stability of multiplication and triangular solves} \label{sec:stability-multiplication-solves}

The basic primitive in all of these algorithms are taking linear combinations of vectors and performing matrix multiplications.
The stability of these basic operations is analyzed in \cite[\S\S2--3]{Hig02}.
We will use the crudest versions of these results.
For exactly represented vectors $\vec{z},\vec{w} \in \field^n$ and $\vec{v} \in \field^m$ and scalar $\xi \in \field$, we have
\begin{subequations}
\begin{align}
    \norm{\err(\xi \cdot \vec{z})} &\lesssim |\xi| \cdot \norm{\vec{z}} + \mathrm{HOT}, \label{eq:err-scalar-mul}\\
    \norm{\err(\vec{z} \pm \vec{w})} &\lesssim \norm{\vec{z} \pm \vec{w}} u + \mathrm{HOT},  \label{eq:err-zw}\\
    \norm{\err(\mat{B}\vec{z})} &\lesssim \norm{\vec{z}}u + \mathrm{HOT}, \label{eq:err-Bz} \\
    \norm{\err(\mat{B}^*\vec{v})} &\lesssim \norm{\vec{v}}u + \mathrm{HOT}. \label{eq:err-B*v}
\end{align}
\end{subequations}

Solves by triangular matrices will be another important primitive for us.
The key fact about triangular solves is that they are \emph{backward stable}.
In particular, for an exactly represented vector $\vec{z}$,
\begin{equation*}
    \fl(\mat{R}^{-1}\vec{z}) = (\mat{R} + \mat{\Delta R})^{-1} \vec{z} \quad \text{for } \norm{\mat{\Delta R}} \lesssim u.
\end{equation*}
Here, we used the bound $\norm{\mat{R}} \lesssim 1$ (equation \cref{eq:R-maxsing-numerical}).
Using the first-order expansion $(\mat{R} + \mat{\Delta R})^{-1} = \mat{R}^{-1} - \mat{R}^{-1} \cdot \mat{\Delta R} \cdot \mat{R}^{-1} + \mathrm{HOT}$, it follows that
\begin{equation*}
    \fl(\mat{R}^{-1}\vec{z}) = [\mat{R}^{-1} - \mat{R}^{-1} \cdot \mat{\Delta R} \cdot \mat{R}^{-1}] \vec{z} + \mathrm{HOT} = \mat{R}^{-1}\vec{z} + \norm{\mat{R}^{-1}\vec{z}} \cdot \mat{R}^{-1}\vec{e} + \mathrm{HOT}.
\end{equation*}
We conclude that 
\begin{subequations}
\importantresult{\begin{equation} \label{eq:Rinv-err}
    \err(\mat{R}^{-1}\vec{z}) = \norm{\mat{R}^{-1}\vec{z}} \cdot \mat{R}^{-1}\vec{e} + \mathrm{HOT}.
\end{equation}}
A similar argument shows that
\begin{equation} \label{eq:R*inv-err}
    \err(\mat{R}^{-*}\vec{z}) = \norm{\mat{R}^{-*}\vec{z}} \cdot \mat{R}^{-*}\vec{e} + \mathrm{HOT}.
\end{equation}
\end{subequations}

\section{Stability of interleaved multiplications} \label{sec:stability-interleaved}

The sketch-and-precondition algorithm consists sequences of multiplies with the matrix $\mat{B}$, the preconditioner $\mat{R}$, and their adjoints.

Let us first analyze the floating-point errors incurred in evaluating the product $\mat{R}^{-*}\mat{B}^*\vec{z}$.
Combining \cref{eq:err-B*v} and \cref{eq:R*inv-err}, we obtain
\begin{equation*}
    \err(\mat{R}^{-*}\mat{B}^*\vec{z}) = \underbrace{\norm{\mat{R}^{-*}\mat{B}^*\vec{z}} \cdot \mat{R}^{-*}\vec{e}}_{\text{errors from applying $\mat{R}^{-*}$}} + \underbrace{\norm{\vec{z}} \cdot \mat{R}^{-*}\vec{e}'}_{\mat{R}^{-*}\cdot \err(\mat{B}^*\vec{z})} + \mathrm{HOT}.
\end{equation*}
Applying the bound $\norm{\smash{\mat{B}\mat{R}^{-1}}} \lesssim 1$ (equation \cref{eq:BR-norm-numerical}), we conclude
\importantresult{\begin{equation} \label{eq:R-*B*z}
    \err(\mat{R}^{-*}\mat{B}^*\vec{z}) = \norm{\vec{z}} \cdot\mat{R}^{-*}\vec{e} + \mathrm{HOT}.
\end{equation}}

Next, let us analyze $\mat{R}^{-*}\mat{B}^*\mat{B}\mat{R}^{-1}\vec{z}$.
First, applying \cref{eq:err-Bz} and \cref{eq:Rinv-err}, observe that
\begin{equation*}
    \err(\mat{B}\mat{R}^{-1}\vec{z}) = \underbrace{\norm{\mat{R}^{-1}\vec{z}} \cdot \mat{B}\mat{R}^{-1}\vec{e}}_{\mat{B}\cdot \err(\mat{R}^{-1}\vec{z})} + \underbrace{\norm{\mat{R}^{-1}\vec{z}} \cdot \vec{e}'}_{\text{errors from applying $\mat{B}$}} + \mathrm{HOT}.
\end{equation*}
Deploying the bound $\norm{\smash{\mat{B}\mat{R}^{-1}}} \lesssim 1$ (equation \cref{eq:BR-norm-numerical}) and the inequality $\norm{\smash{\mat{R}^{-1}\vec{z}}} \le \norm{\vec{z}} / \sigma_{\mathrm{min}}(\mat{R}) \lesssim \kappa \norm{\vec{z}}$ (equation \cref{eq:R-minsing-numerical}), we obtain
\begin{equation*}
    \norm{\err(\mat{B}\mat{R}^{-1}\vec{z})} \lesssim \kappa u \norm{\vec{z}} + \mathrm{HOT}.
\end{equation*}
Combining this result with \cref{eq:R-*B*z} and using the bound  $1/\sigma_{\mathrm{min}}(\mat{R}) \lesssim \kappa$, we conclude
\importantresult{\begin{equation} \label{eq:R-*B*BR-1z}
    \norm{\err(\mat{R}^{-*}\mat{B}^*\mat{B}\mat{R}^{-1}\vec{z})} \lesssim \kappa u \norm{\vec{z}} + \mathrm{HOT}.
\end{equation}}
The results \cref{eq:R-*B*z,eq:R-*B*BR-1z} will be used later in our analysis.

\section{Accuracy of Lanczos linear solves} \label{sec:lanczos-stability}

The finite-precision analysis of Lanczos function approximation and Lanczos linear solves is the subject of papers \cite{DK91,DGK98,MMS18}.
We draw on the results of these works, though only in a qualitative manner.
Specifically, we will use the following consequence of the analysis of these papers:

\begin{fact}[Lanczos linear solves: Well-conditioned matrix] \label{fact:lanczos-numerical}
    Let $\mat{M}$ be a positive definite matrix with condition number bounded by an absolute constant
    \begin{equation*}
        \cond(\mat{M}) \le \mathrm{const},
    \end{equation*}
    and assume matrix--vector products are computed to \emph{effective precision} $\tilde{u}\ge u$:
    \begin{equation*}
        \norm{\err(\mat{M}\vec{z})} \lesssim \norm{\mat{M}}\norm{\vec{z}}\tilde{u} \quad \text{for every exactly represented } \vec{z} \in \real^n.
    \end{equation*}
    Then $\order(\log(1/\tilde{u}))$ steps of the Lanczos linear solver produces an approximation to $\mat{M}^{-1}\vec{z}$ that is forward stable in the effective precision $\tilde{u}$:
    \begin{equation} \label{eq:lan-lin-solve-finite-precision}
        \norm{\smash{\err(\mat{M}^{-1}\vec{z})}} \lesssim \norm{\mat{M}}^{-1}\norm{\vec{z}}\tilde{u}.
    \end{equation}
    The $\lesssim$ notation in \cref{eq:lan-lin-solve-finite-precision} suppresses polylogarithic factors in the inverse-effective precision $1/\tilde{u}$ and the dimension of $\mat{M}$.
\end{fact}

See \cite[App.~C]{EMN24} for a detailed explanation of how this result can be derived from \cite{MMS18}.

In the Lanczos-based sketch-and-precondition method, we use Lanczos to solve the normal equations
\begin{equation*}
    \mat{M}\vec{y} = \vec{f} \quad \text{for } \mat{M} \coloneqq \mat{R}^{-*}\mat{B}^*\mat{B}\mat{R}^{-1}, \: \vec{f} \coloneqq \mat{R}^{-*}\mat{B}^*(\vec{c} - \mat{B}\vec{x}_i).
\end{equation*}
As shown in \cref{eq:R-*B*BR-1z}, multiplies $\mat{M}\vec{z} = \mat{R}^{-*}(\mat{B}^*(\mat{B}(\mat{R}^{-1}\vec{z})))$ is forward stable in effect precision $\tilde{u} = \kappa u$, up to higher order terms.
Further, by \cref{eq:BR-norm-numerical,eq:BR-minsing-numerical}, the matrix $\mat{M}$ is well-conditioned, $\cond(\mat{M}) \le \mathrm{const}$.
Thus, \cref{fact:lanczos-numerical} implies that Lanczos produces a solution to the normal equations satisfying
\importantresult{\begin{equation} \label{eq:lanczos-normal-num-error}
    \norm{\fl(\vec{y}) - \mat{M}^{-1}\fl(\vec{f})} \lesssim \kappa u \norm{\fl(\vec{f})} + \mathrm{HOT}.
\end{equation}}

\section{The error formula} \label{sec:error-formula}

We now analyze an iterative refinement step using sketch-and-precondition:
\begin{equation} \label{eq:iterative-refinement-step}
    \vec{x}_{i+1} \gets \vec{x}_i + \mat{R}^{-1} \, \underbrace{(\mat{R}^{-*}\mat{B}^*\mat{B}\mat{R}^{-1})^{-1}}_{\text{applied via Lanczos}} \, \mat{R}^{-*}\mat{B}^*(\vec{c} - \mat{B}\vec{x}_i).   
\end{equation}
Here, $\vec{x}_i$ denotes an arbitrary $i$th iterate.
Throughout this section, we let $\fl$ and $\err$ denote only the floating-point errors incurred during the refinement step \cref{eq:iterative-refinement-step}.

Begin by using the error bounds for vector subtraction \cref{eq:err-zw} and matrix multiplication \cref{eq:err-Bz} to obtain
\begin{equation*}
    \norm{\err(\vec{c} - \mat{B}\vec{x}_i)}\lesssim (1 + \norm{\vec{x}_i}) u + \mathrm{HOT}.
\end{equation*}
Next, invoking \cref{eq:R-*B*z} yields
\begin{equation*}
    \err(\mat{R}^{-*}\mat{B}^*(\vec{c} - \mat{B}\vec{x}_i)) = \underbrace{(1 + \norm{\vec{x}_i}) \mat{R}^{-*}\mat{B}^*\vec{e}}_{\mat{R}^{-*}\mat{B}^*\cdot \err(\vec{c} - \mat{B}\vec{x}_i)} + \underbrace{\norm{\vec{c} - \mat{B}\vec{x}_i} \cdot \mat{R}^{-*}\vec{e}'}_{\text{errors from applying } \mat{R}^{-*}\mat{B}^*} + \mathrm{HOT}.
\end{equation*}
Employing the bound $\norm{\mat{R}^{-*}\mat{B}^*} \lesssim 1$ (\cref{eq:BR-norm-numerical}), this equation simplifies 
\begin{equation*}
    \err(\mat{R}^{-*}\mat{B}^*(\vec{c} - \mat{B}\vec{x}_i)) = (1 + \norm{\vec{x}_i}) \vec{e} + \norm{\vec{c} - \mat{B}\vec{x}_i} \cdot \mat{R}^{-*}\vec{e}' + \mathrm{HOT}.
\end{equation*}
Now, introduce the shorthand $\mat{M} \coloneqq \mat{R}^{-*}\mat{B}^*\mat{B}\mat{R}^{-1}$ and apply the Lanczos linear solve bound \cref{eq:lanczos-normal-num-error} to obtain
\begin{multline} \label{eq:error-intermediate}
    \err(\mat{M}^{-1}\mat{R}^{-*}\mat{B}^*(\vec{c} - \mat{B}\vec{x}_i)) = \underbrace{(1 + \norm{\vec{x}_i}) \vec{e} + \norm{\vec{c} - \mat{B}\vec{x}_i} \cdot \mat{M}^{-1}\mat{R}^{-*}\vec{e}'}_{\mat{M}^{-1}\cdot \err(\mat{R}^{-*}\mat{B}^*(\vec{c} - \mat{B}\vec{x}_i))} \\ + \underbrace{\kappa \norm{\mat{R}^{-*}\mat{B}^*(\vec{c} - \mat{B}\vec{x}_i)} \vec{e}''}_{\text{errors from applying } \mat{M}^{-1}} + \mathrm{HOT}.
\end{multline}
We have simplified by noting that $\norm{\smash{\mat{M}^{-1}}} \le \mathrm{const}$.
To reduce further, we note that 
\begin{multline*}
    \norm{\mat{R}^{-*}\mat{B}^*(\vec{c} - \mat{B}\vec{x}_i)} = \norm{\mat{R}^{-*}\mat{B}^*(\vec{c} - \mat{B}\vec{x} + \mat{B}(\vec{x} - \vec{x}_i))} \\ = \norm{\mat{R}^{-*}\mat{B}^*(\mat{B}\vec{x} - \mat{B}\vec{x}_i)} \lesssim \norm{\mat{B}(\vec{x} - \vec{x}_i)}.
\end{multline*}
The second equality is the orthogonality of the residual $\vec{c} - \mat{B}\vec{x}$ to the range of $\mat{B}$ and the inequality is \cref{eq:BR-norm-numerical}.
We also note that $\mat{M}^{-1}\mat{R}^{-*} = \mat{R}(\mat{B}^*\mat{B})^{-1}$.
Applying all of these observations and combining with \cref{eq:error-intermediate}, we obtain
\begin{multline*}
    \err(\mat{M}^{-1}\mat{R}^{-*}\mat{B}^*(\vec{c} - \mat{B}\vec{x}_i)) = (1 + \norm{\vec{x}_i} + \kappa \norm{\mat{B}(\vec{x} - \vec{x}_i)}) \vec{e} \\ + \norm{\vec{c} - \mat{B}\vec{x}_i} \cdot \mat{R}(\mat{B}^*\mat{B})^{-1}\vec{e}' + \mathrm{HOT}.
\end{multline*}
Now, employ the inversion bound \cref{eq:Rinv-err} and addition bound \cref{eq:err-zw} to conclude
\begin{multline*}
    \err(\vec{x}_{i+1}) = \underbrace{(1 + \norm{\vec{x}_i} + \kappa \norm{\mat{B}(\vec{x} - \vec{x}_i)}) \mat{R}^{-1}\vec{e} + \norm{\vec{c} - \mat{B}\vec{x}_i} \cdot (\mat{B}^*\mat{B})^{-1}\vec{e}'}_{\mat{R}^{-1}\cdot \err(\mat{M}^{-1}\mat{R}^{-*}\mat{B}^*(\vec{c} - \mat{B}\vec{x}_i))} \\ + \underbrace{\norm{\mat{R}^{-1}\mat{M}^{-1}\mat{R}^{-*}\mat{B}^*(\vec{c} - \mat{B}\vec{x}_i)}\cdot \mat{R}^{-1}\vec{e}''}_{\text{errors from applying} \mat{R}^{-1}} + \underbrace{\norm{\vec{x}} \vec{e}'''}_{\text{errors from final addition}} + \mathrm{HOT}.
\end{multline*}
This display constitutes an error bound for a single sketch-and-precondition step.
To simplify, note that
\begin{equation*}
    \norm{\mat{R}^{-1}\mat{M}^{-1}\mat{R}^{-*}\mat{B}^*(\vec{c} - \mat{B}\vec{x}_i)} = \norm{\mat{B}^\dagger(\vec{c} - \mat{B}\vec{x}_i)} = \norm{\vec{x} - \vec{x}_i} \le \norm{\vec{x}} + \norm{\vec{x}_i}.
\end{equation*}
Thus, we have shown
\importantresult{\begin{multline} \label{eq:error-formula-spre}
    \err(\vec{x}_{i+1}) = (1 + \norm{\vec{x}} + \norm{\vec{x}_i} + \kappa \norm{\mat{B}(\vec{x} - \vec{x}_i)}) \mat{R}^{-1}\vec{e} \\ + \norm{\vec{c} - \mat{B}\vec{x}_i} \cdot (\mat{B}^*\mat{B})^{-1}\vec{e}' + \mathrm{HOT}.
\end{multline}}

\section{Forward stability of sketch-and-precondition} \label{sec:forward-stability-spre-proof}

Now we show strong forward stability of sketch-and-precondition with the sketch-and-solve initialization.
By \cref{eq:sketch-and-solve-finite-precision}, the numerically computed sketch-and-solve solution satisfies
\begin{equation*}
    \norm{\fl(\vec{x}_0)} \le \norm{\vec{x}} + \norm{\vec{x} - \fl(\vec{x}_0)} \le \norm{\vec{x}} + \kappa \norm{\mat{B}(\vec{x} - \fl(\vec{x}_0))} \lesssim  \norm{\vec{x}} + \kappa \norm{\vec{c} - \mat{B}\vec{x}}.
\end{equation*}
Here, we used the hypothesis $\kappa u\ll 1$.
Similarly,
\begin{equation*}
    1 = \norm{\vec{c}} \le \norm{\vec{c} - \mat{B}\vec{x}} + \norm{\vec{x}}.
\end{equation*}
Substituting the previous two displays in the error formula \cref{eq:error-formula-spre} gives
\begin{equation*}
    \err(\vec{x}_1) = (\norm{\vec{x}} + \kappa \norm{\vec{c} - \mat{B}\vec{x}}) \mat{R}^{-1}\vec{e} + \norm{\vec{c} - \mat{B}\vec{x}} \cdot (\mat{B}^*\mat{B})^{-1}\vec{e}' + \mathrm{HOT}.
\end{equation*}
Finish up by multiplying by $\mat{B}$, taking norms, and using the bound s$\norm{\mat{B}\mat{R}^{-1}} \lesssim 1$ (equation \cref{eq:BR-norm-numerical}) and $\norm{\mat{B}(\mat{B}^*\mat{B})^{-1}} = \norm{\mat{B}^\dagger} = \kappa$.
We obtain
\importantresult{\begin{equation} \label{eq:strong-forward-spre-result}
    \norm{\mat{B}(\vec{x} - \fl(\vec{x}_1))} \lesssim \norm{\vec{x}} u + \kappa \norm{\vec{c} - \mat{B}\vec{x}} u + \mathrm{HOT}.
\end{equation}}
This bound is strong forward stability, up to the higher-order terms.

\section{Backward stability of SPIR} \label{sec:backward-stability-spir-proof}

Now, we show backward stability of SPIR.
Our main tool will be the following characterization \cite[Thm.~2.8 \& Cor.~2.9]{EMN24}: 
\begin{theorem}[Backward stability: Componentwise errors] \label{thm:backward-stable-characterization}
    Consider a \warn{real} least-squares problem \cref{eq:least-squares} with the normalization $\norm{\mat{B}} = \norm{\vec{c}} = 1$, and introduce the SVD $\mat{B} = \sum_{i=1}^n \sigma_i^{\vphantom{*}}\vec{u}_i^{\vphantom{*}} \vec{v}_i^*$.
    A vector $\hatvector{x}$ is backward stable solution to the least-squares problem \cref{eq:least-squares} if and only if 
    \begin{equation*}
        \left| \vec{v}_i^*(\vec{x} - \hatvector{x}) \right| \lesssim \sigma_i(\mat{B})^{-1} \cdot (1+\norm{\hatvector{x}})u + \sigma_i(\mat{B})^{-2} \cdot \norm{\vec{c} - \mat{B}\hatvector{x}} u \quad \text{for each } i=1,\ldots,n.
    \end{equation*}
    Consequently, if the error takes the form
    \begin{equation*}
        \vec{x} - \hatvector{x} = (1 + \norm{\hatvector{x}}) \cdot \fl(\mat{R})^{-1}\vec{e} + \norm{\vec{c} - \mat{B}\hatvector{x}} \cdot (\mat{B}^*\mat{B})^{-1}\vec{e}',
    \end{equation*}
    then $\hatvector{x}$ is backward stable.
\end{theorem}

The proof follows by comparing the backward error to the Karlson--Wald\'en estimate using \cref{fact:kw-estimate}, then decomposing the Karlson--Wald\'en estimate using the singular value decomposition of $\mat{B}$.
The computation is straightforward, and we omit the details here.

To prove backward stability of SPIR, we use the following ingredients: the characterization theorem (\cref{thm:backward-stable-characterization}), the error formula \cref{eq:error-formula-spre}, and the strong forward stability of sketch-and-precondition \cref{eq:strong-forward-spre-result}.
Begin by using strong forward stability \cref{eq:strong-forward-spre-result} of the numerically computed solution $\fl(\vec{x}_1)$ together with the error formula \cref{eq:error-formula-spre} to obtain
\begin{equation*}
    \err(\vec{x}_2) = (1 + \norm{\vec{x}}) \mat{R}^{-1}\vec{e} + \norm{\vec{c} - \mat{B}\vec{x}} \cdot (\mat{B}^*\mat{B})^{-1}\vec{e}' + \mathrm{HOT}.
\end{equation*}
We see that $\fl(\vec{x}_2)$ and $\norm{\vec{x}}$ differ by an amount of size $\order(u)$, so we may replace $\norm{\vec{x}}u$ by $\norm{\fl(\vec{x}_2)}u$ up to higher order terms.
Similarly, the $\vec{x}$ minimizes the least-squares residual norm, so $\norm{\vec{c} - \mat{B}\vec{x}} \le \norm{\vec{c} - \mat{B}\hatvector{x}}$.
Thus, we have shown
\importantresult{\begin{equation*}
    \err(\vec{x}_2) = (1 + \norm{\fl(\vec{x}_2)}) \mat{R}^{-1}\vec{e} + \norm{\vec{c} - \mat{B}\fl(\vec{x}_2)} \cdot (\mat{B}^*\mat{B})^{-1}\vec{e}' + \mathrm{HOT}.
\end{equation*}}
This conclusion is precisely backward stability, up to higher-order terms.

\begin{remark}[What about the benign rounding error assumption?]
    The careful reader may notice that our derivation did not use the non-pathological rounding error assumption (\cref{def:non-pathological}).
    Thus, we conclude that backward stability \warn{to first order} holds without this assumption.
    As far as I am aware, this assumption is needed to control higher-order terms and prove the backward stability of SPIR (\cref{thm:spir-backward}).
\end{remark}

\chapter{Sketching, solvers, and stability: Open problems}

This chapter presents open problems related to this part of the thesis.
First, in \cref{sec:lanczos-to-cg-to-lsqr}, we discuss improving our stability analysis to treat CG or LSQR in place of Lanczos.
Second, \cref{sec:krylov-stability} investigates open problems in the numerical stability of Krylov methods.
We conclude in \cref{sec:pcg-stable} by pondering the stability properties of preconditioned iterative methods for square systems of linear equations.

\section{From Lanczos to conjugate gradient to LSQR} \label{sec:lanczos-to-cg-to-lsqr}

The major limitation of our analysis sketch-and-precondition and SPIR in the last chapter was that we analyzed a Lanczos-based version of the procedure.
In practice, we instead use LSQR or perhaps conjugate gradient on the normal equations.
To close this gap between theory and practice, we ask: Could we analyze sketch-and-precondition-type methods with conjugate gradient or LSQR?

The path to analyzing sketch-and-precondition with conjugate gradient seems more clear.
The conjugate gradient algorithm has already been analyzed many times in finite-precision arithmetic \cite{Gre89,GS92,Gre97,MS06}.
Further, the pseudocodes for conjugate gradient and the Lanczos algorithm are fairly similar.
Given this existing body of work and these connections, it is natural to hope that a version of \cref{fact:lanczos-numerical} should be attainable for conjugate gradient.
I am actively working on this question in collaboration with Deeksha Adil, Anne Greenbaum, and Christopher Musco.

Analyzing sketch-and-precondition with LSQR seems harder as a practical matter, as there is no existing rigorous analysis of the LSQR method in finite precision arithmetic (to the best of my knowledge).
However, in principle, any tools used to analyze conjugate gradient or Lanczos should be extensible to LSQR, as well.

\section{Numerical stability of Krylov methods} \label{sec:krylov-stability}

In the Lanczos stability result (\cref{fact:lanczos-numerical}), the $\lesssim$ notation suppresses prefactors that, apon closer analysis, are quantitatively large.
The analysis of \cite{MMS18}, which underrides the proof of \cref{fact:lanczos-numerical}, shows that the Lanczos algorithm solves an $n\times n$ positive definite linear system $\mat{M}\vec{f}$ up accuracy roughly
\begin{equation} \label{eq:mms}
    \norm{\mat{M}} \cdot \norm{\fl(\mat{M}^{-1}\vec{f}) - \mat{M}^{-1}\vec{f}} \lessapprox \poly(n) \cdot \cond(\mat{M})^3\tilde{u} \log^4(1/\tilde{u})\cdot \norm{\vec{f}}.
\end{equation}
Here, as in \cref{fact:lanczos-numerical}, $\tilde{u}\ge u$ denotes the effective precision at which matvecs with $\mat{M}$ are computed.
The analysis of \cite[Thm.~2.2]{DGK98} is stronger, roughly implying that
\begin{equation} \label{eq:dgk}
    \norm{\mat{M}} \cdot \norm{\fl(\mat{M}^{-1}\vec{f}) - \mat{M}^{-1}\vec{f}} \lessapprox \poly(n) \cdot \cond(\mat{M})^{2.25}\tilde{u}\log^{0.5}(1/\tilde{u})\cdot \norm{\vec{f}}.
\end{equation}
(The full version of \cite[Thm.~2.2]{DGK98} is stronger \cref{eq:dgk}, as it bounds the residual $\norm{\vec{f} - \mat{M}\fl(\mat{M}^{-1}\vec{f})}$.)
As a point of comparison, recall that a (forward) stable algorithm for the linear system $\mat{M}\vec{y} = \vec{f}$ achieves a much higher accuracy
\begin{equation} \label{eq:lanczos-if-forward-stable}
    \norm{\mat{M}} \cdot \norm{\fl(\mat{M}^{-1}\vec{f}) - \mat{M}^{-1}\vec{f}} \lessapprox \poly(n) \cdot \cond(\mat{M})u\cdot \norm{\vec{f}}.
\end{equation}
In my numerical experience, I have always found the Lanczos algorithm to be forward stable or nearly so.
As such, I believe the the polylogarithmic factors $\poly(\log(1/\tilde{u}))$ and superlinear powers of $\cond(\mat{M})$ in \cref{eq:mms,eq:dgk} are parasitic.

In our analysis, the unappealing aspects of \cref{eq:mms,eq:dgk} are hidden by the $\lesssim$ notation.
The condition number of $\mat{M} = \mat{R}^{-*}\mat{B}^*\mat{B}\mat{R}^{-1}$ is $\order(1)$ and polylogarithmic factors in $\log(1/u) \ge \log(1/\tilde{u})$ are suppressed under the $\lesssim$ relation.
Still, if one were to develop versions of stability theorems in this thesis with explicit prefactors, those prefactors would pay a steep price for the weak bounds \cref{eq:mms,eq:dgk}.
As such, developing improved versions of \cref{fact:lanczos-numerical} with better dependence on $\cond(\mat{M})$ and $\log(1/\tilde{u})$ is a natural question for future work.
In particular, a precise open question is to determine whether the Lanczos algorithm is forward stable \cref{eq:lanczos-if-forward-stable}.

\section{From least squares to linear systems} \label{sec:pcg-stable}

On inspection, the analysis in \cref{ch:stability-proofs} uses randomness only to \warn{construct} the preconditioner $\mat{R}$ and the initialization $\vec{x}_0$.
As such, if one had another means of producing a high-enough quality preconditioner $\mat{R}$ and initialization $\vec{x}_0$, our analysis would also establish backward stability in this case (with a single step of iterative refinement).
Remarks to these effects are made in \cite[\S6.7]{EMN24}.

Unfortunately, finding a good initialization $\vec{x}_0$ satisfying $\norm{\vec{c} - \mat{B}\vec{x}_0} \lesssim \norm{\vec{c} - \mat{B}\vec{x}}$ appears challenging beyond the setting of highly overdetermined least squares.
Indeed, for a least-squares problem which is consistent ($\vec{c} = \mat{B}\vec{x}$), finding such a good initialization amounts to solving the problem!
As such, randomized preconditioning represents the only setting of which I am aware for which one has access to both high-quality initialization and high-quality preconditioner.

Still, the techniques we studied in this thesis could be useful beyond the setting of randomized least squares.
In particular, it is natural to ask:
\actionbox{Can we say anything about the solution to \emph{linear systems of equations} by preconditioned iterative methods from trivial initialization $\vec{x}_0 = \vec{0}$?}
I investigated this question with Anne Greenbaum and Yuji Nakatsukasa in \cite{EGN25}, and we showed that iterative methods on the \warn{preconditioned normal equations} converge to backward stability after sufficiently many iterative refinement.
As with this thesis, the analysis only holds for Lanczos, but the numerical results suggest it holds for LSQR as well.
I also make the caveat that our analysis requires that our analysis requires that applying the inverse-preconditioner $\vec{z} \mapsto \mat{R}^{-1}\vec{z}$ is done in a backward stable manner.

This leaves open a fundamental question: What about other iterative methods?
Is preconditioned conjugate gradient (with iterative refinement) backward stable, provided with a high-enough quality preconditioner?
What about preconditioned GMRES?
These are natural and basic questions about the attainable accuracy of some of the most fundamental preconditioned iterative methods of computational mathematics, and they remain open.
Some preliminary numerical evidence for preconditioned conjugate gradient is provided in \cite[Fig.~4]{EGN25}, showing that preconditioned conjugate gradient is backward stable \warn{with, and only with, iterative refinement}.

\appendix

\chapter{\texorpdfstring{Incremental \QR decomposition}{Incremental QR decomposition}} \label{app:incremental-qr}

In both \cref{part:random-pivoting,part:loo} of this thesis, we are interested in maintaining a \QR decomposition of a tall matrix $\mat{Y} \in \field^{m\times k}$ that evolves by having new columns appended to its end.
It is straightforward using appropriate calls to a linear algebra library like \LAPACK \cite{ABB+99}, though I am unaware of a reference.
For completeness, I provide a description of this ``incremental \QR decomposition'' primitive here.
I also provide a MATLAB implementation using \warn{internal undocumented functions}, which provide an interface to the corresponding \LAPACK calls.
(Warning: MATLAB internal functions are subject to be removed or changed between releases.)
In the following discussion, we use names of \LAPACK routines for real, double-precision arithmetic.

\begin{figure}[t]
    \centering
    \includegraphics[width=0.98\textwidth]{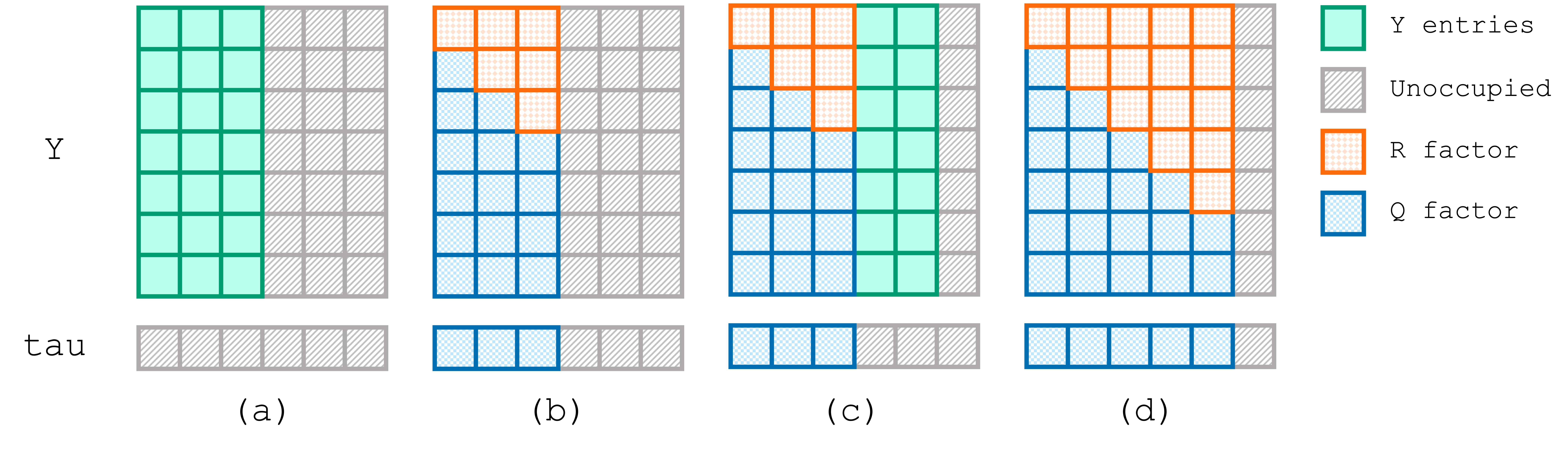}
    \caption[Data layout for incremental \QR decomposition]{Data layout for incremental \QR decomposition, illustrated with a buffer of size 42 for \texttt{Y} and size 6 for \texttt{tau}.
    (a) Initial matrix $\mat{Y} \in \field^{7\times 3}$ stored in first 21 entries of \texttt{Y} buffer; remaining buffer positions are unoccupied.
    (b) \QR decomposition of $\mat{Y}$ overwritten on $\mat{Y}$ in place; upper triangle of $\mat{Y}$ stores the R factor, and lower triangle (and first three entries of \texttt{tau} buffer) store the Q factor implicitly using Householder reflectors.
    (c) Matrix $\mat{Y}$ is expanded to size $7\times 5$; new columns are appended to the end of the \texttt{Y} buffer.
    (d) New columns are brought into \QR decomposition format.}
    \label{fig:incremental-qr}
\end{figure}

Suppose the entries of $\mat{Y} \in \field^{m\times k}$ are stored in column-major order in a buffer \texttt{Y}.
New (blocks of) columns are added by appending to the end of the buffer.
If the allocated space for the buffer is ever exhausted, a new buffer of twice the size is allocated and the contents of the existing buffer are copied over.
To store a \QR decomposition using \LAPACK, we allocate a second small buffer \texttt{tau}, which also can be dynamically resized as necessary.
See \cref{fig:incremental-qr}(a) for an illustration.

A (full) \QR decomposition
\begin{equation*}
    \mat{Y} = \mat{Q}\twobyone{\mat{R}}{\mat{0}} \quad \text{for } \mat{Q} \in \field^{m\times m} \text{ and } \mat{R} \in \field^{k\times k}
\end{equation*}
can now be computed in place using the \LAPACK command \texttt{dgeqrf}.
The matrix $\mat{R}$ is overwritten over the upper triangular portion of the buffer \texttt{Y}, and the matrix $\mat{Q}$ is represented implicitly as a product of Householder reflectors, which are stored in the strictly lower triangular parts of the \texttt{Y} buffer and the first $k$ entries of \texttt{tau}.
See \cref{fig:incremental-qr}(b) for illustration.

New columns or blocks of columns can be adjoined to the matrix $\mat{Y}$ by appending to the buffer \texttt{Y}.
The newly added columns are stored as-is, whereas the existing columns have been overwritten by the \QR decomposition.
One can continue to add columns in this way for as long as desired.
See \cref{fig:incremental-qr}(c).

Suppose at some future time one wishes to compute a \QR decomposition of the matrix $\mat{Y}$ with its new columns.
Denote $\mat{Y}_1 \in \field^{m\times k}$ the previous matrix with its full \QR decomposition 
\begin{equation*}
    \mat{Y}_1 = \mat{Q}_1\twobyone{\mat{R}_{11}}{\mat{0}}
\end{equation*}
already computed and stored in place, and denote by $\mat{Y}_2 \in \field^{m\times \ell}$ the newly added columns, so that $\mat{Y} = \onebytwo{\mat{Y}_1}{\mat{Y}_2}$.
As a first step, apply the matrix $\mat{Q}_1^*$ to $\mat{Y}_2$ in place, which can be accomplished using the \LAPACK routine \texttt{dormqr}.
The result of this operation can be written as
\begin{equation*}
    \mat{Q}_1^*\mat{Y} = \mat{Q}_1^*\onebytwo{\mat{Y}_1}{\mat{Y}_2} \eqqcolon \twobytwo{\mat{R}_{11}}{\mat{R}_{12}}{\mat{0}}{\mat{\tilde{Y}}_2} \quad \text{for } \mat{R}_{12} \in \field^{k\times \ell} \text{ and } \mat{\tilde{Y}}_2 \in \field^{(m-k)\times \ell}.
\end{equation*}
The matrix $\mat{Y}_2$ in the buffer \texttt{Y} has been overwritten by two matrices, $\mat{R}_{12}$ and $\mat{\tilde{Y}}_2$ stacked on top of each other.
The top part of this matrix, $\mat{R}_{12}$, shall serve as the upper right block of the R factor of the full matrix $\mat{Y}$.
The bottom part $\mat{\tilde{Y}}_2$ is \QR factorized in place, again using the \LAPACK routine \texttt{dgeqrf}.
(The Q factor requires $\ell$ additional elements of storage, which are placed into the next $\ell$ entries of the buffer \texttt{tau}.)
Symbolically, this \QR decomposition may be written
\begin{equation*}
    \mat{\tilde{Y}_2} = \mat{Q}_2 \twobyone{\mat{R}_{22}}{\mat{0}} \quad \text{for } \mat{Q}_2 \in \field^{(m-k)\times (m-k)} \text{ and } \mat{R}_{22} \in \field^{\ell\times \ell}.
\end{equation*}
Combining, the two previous displays, we have obtained a \QR decomposition of the new $\mat{Y}$ matrix
\begin{equation*}
    \mat{Y} = \onebytwo{\mat{Y}_1}{\mat{Y}_2} = \mat{Q}_1 \twobytwo{\Id}{\mat{0}}{\mat{0}}{\mat{Q}_2} \begin{bmatrix}
        \mat{R}_{11} & \mat{R}_{12} \\
        \mat{0} & \mat{R}_{22} \\
        \mat{0} & \mat{0}
    \end{bmatrix} \eqqcolon \mat{Q} \twobyone{\mat{R}}{0}.
\end{equation*}
The dimensions are $\mat{Q} \in \field^{m\times m}$ and $\mat{R}\in\field^{(k+\ell)\times (k+\ell)}$.
We have updated the \QR decomposition in place; see \cref{fig:incremental-qr}(d) for the data layout at the end of this procedure.

\myparagraph{Computational cost}
The computational cost of the update procedure consists of $\order(mk\ell )$ operations to multiply $\mat{Y}_2$ by $\mat{Q}_1^*$, and $\order(m\ell^2)$ operations to compute the \QR decomposition of $\mat{\tilde{Y}}_2$ in place.
The total cost of updating the \QR decomposition of an $m\times k$ matrix with $\ell$ column appends is $\order(m\ell (k+\ell))$ operations.

This incremental \QR procedure is very efficient.
To see this, suppose we aggregate columns in increments of sizes $\ell_1,\ldots,\ell_t$ resulting in a total of $p = \sum_{i=1}^t \ell_i$ columns.
Then the total cost of computing the \QR decomposition is 
\begin{equation*}
    \order \left( \sum_{i=1}^t m\ell_i(\ell_i + p) \right) = \order(mp^2).
\end{equation*}
The total cost of the incremental \QR decomposition, $\order(mp^2)$, is the same asympotically as computing the \QR decomposition of an $m\times p$ matrix $\mat{Y}$ all at once.
However, the practical speed of this procedure very much depends on the block sizes $\ell_i$; many small column updates are much slower to process than one larger one.

\myparagraph{Inverting the R factor}
For the applications of the incremental \QR decomposition in \cref{part:loo} of this thesis, we need to have access to the \emph{inverse} of the $\mat{R}$ matrix.
While computing the inverse of a matrix explicitly is usually heresy in numerical linear algebra \cite[\S14.1]{Hig02}, this instance may be an exception where matrix inversion is well-justified.
(The trace estimates produced by \XTrace do reach machine accuracy---see \cref{fig:trace-comparison}---so using explicit inversion in that application does not appear to resulting in numerical issues!)

The inverse of the R factor can be computed in place using the \LAPACK routine \texttt{dtrtri}.
To perform the update steps, observe that we have the block matrix inverse identity
\begin{equation*}
    \twobytwo{\mat{R}_{11}}{\mat{R}_{12}}{\mat{0}}{\mat{R}_{22}}^{-1} = \twobytwo{\mat{R}_{11}^{-1}}{-\mat{R}_{11}^{-1}\mat{R}_{12}\mat{R}_{22}^{-1}}{\mat{0}}{\mat{R}_{22}^{-1}}.
\end{equation*}
The inverses $\mat{R}_{11}^{-1}$ and $\mat{R}_{22}^{-1}$, after which the triple product $-\mat{R}_{11}^{-1}\mat{R}_{12}\mat{R}_{22}^{-1}$ can be computed in place using the BLAS routine \texttt{dtrmm}. 

\myparagraph{Working with the Q factor}
This implementation presents the Q factor implicitly, represented in the lower triangular portion of the buffer \texttt{Y} and the buffer \texttt{tau}.
For most purposes, this representation is sufficient, as one can compute matrix products with $\mat{Q}$ and $\mat{Q}^*$ using the \LAPACK routine \texttt{dormqr}.
If desired, an explicit representation the economy Q factor can be extracted using the \LAPACK routine \texttt{dorgqr}; the cost is $\order(mp^2)$ operations, where $p$ denotes the current number of columns.
For computational efficiency, one should make sure to only extract the Q factor only once, when the matrix $\mat{Y}$ has reached its maximal size.

\myparagraph{The design of \LAPACK}
Many of \LAPACK's design decisions, such as the infamous ``leading dimension'' input for matrix algorithms and the representation of \QR decomposition via packed, in-place formats are confusing and unintuitive for new users.
The effortlessness by which the incremental \QR decomposition primitive can be implemented using \LAPACK routines speaks to the wisdom of these aspects of \LAPACK. 
Indeed, the efficient implementation described here would not be possible, or would be much less efficient, without \LAPACK's in-place \QR decomposition and leading dimension argument, allowing for seamless execution of matrix routines in place on submatrices.
As a technical community, we owe a great debt to the designers of well-designed pieces of computational mathematics software like \LAPACK.

\myparagraph{Implementation}
A MATLAB implementation of the incremental QR primitive is provided in Program~\ref{prog:incremental_qr}.
Instead of \LAPACK routines, this implementation uses undocumented MATLAB functions, which provide an interface to the low-level \LAPACK routines.
This implementation works in MATLAB 2023b, but their is no guarantee the relevant private functions are not removed or changed in upcoming MATLAB releases.
This program defines an \texttt{incremental\_qr} class, whose usage is demonstrated in the following code segment:

\begin{lstlisting}
iqr = incremental_qr(Y); % Initialize incremental_qr object
iqr.addcols(Ynew);       % Append Ynew to matrix, update QR
Qx = iqr.applyQ(x);      % Compute Q*x
Qtx = iqr.applyQt(x);    % Compute Q'*x
y = iqr.projectOut(x);   % Compute x - Q*Q'*x
Q = iqr.getQ();          % Extract economy Q factor
R = iqr.getR();          % Extrace square R factor
\end{lstlisting}

Some of the functionality of  the \texttt{incremental\_qr} class is broken out into subroutines \texttt{hhqr}, \texttt{apply\_Qt}, and \texttt{get\_Q}, which are provided in \cref{prog:hhqr,prog:apply_Qt,prog:get_Q}.
These subroutines are used in the bespoke incremental \QR implementation in the Householder reflector-based accelerated randomly pivoted \QR implementation (\cref{prog:acc_rpqr}).

\myprogram{Compute a Householder \QR decomposition of input matrix, represented in compact format.}{Warning: Uses \warn{undocumented internal MATLAB functions}.}{hhqr}

\myprogram{Apply the adjoint of (full) Q matrix for a compactly represented Householder \QR decomposition.}{Warning: Uses \warn{undocumented internal MATLAB functions}.}{apply_Qt}

\myprogram{Get the (thin) Q matrix for a compactly represented Householder \QR decomposition.}{Warning: Uses \warn{undocumented internal MATLAB functions}.}{get_Q}

\newpage

\begin{lstlisting}[
    frame=single, 
    caption={An implementation of the incremental \QR primitive described in \cref{app:incremental-qr}. Warning: This implementation is based off of \warn{undocumented internal MATLAB functions} that are subject to be removed or change between releases. Subroutines \texttt{hhqr}, \texttt{apply\_Qt}, and \texttt{get\_Q} are provided in \cref{prog:hhqr,prog:apply_Qt,prog:get_Q}.}, 
    label={prog:incremental_qr}
  ]
classdef incremental_qr < handle
% Implementation of the incremental QR decomposition of a tall matrix
% with an increasing number of columns

    properties
        H   % Matrix to store evolving QR decomposition
        tau % Vector to store scalars for Householder reflections
        k   % Number of columns in current matrix
        n   % Height of matrix
    end
    
    methods
        function obj = incremental_qr(Y,varargin)
            % Initialize the QR decomposition with matrix Y
            [obj.H, obj.tau] = hhqr(Y);
            [obj.n,obj.k] = size(Y);
            if ~isempty(varargin) % Optionally: Set size
                obj.H = [obj.H zeros(obj.n,varargin{1}-obj.k)];
            end
        end
        
        function obj = addcols(obj, Ynew)
            % Add new columns to the existing QR decomposition
            l = size(Ynew,2);

            % Double and copy if necessary
            if obj.k + l > size(obj.H,2)
                obj.H = [obj.H zeros(size(obj.H))];
                obj.tau = [obj.tau;zeros(size(obj.tau))];
                obj.addcols(Ynew);
                return
            end

            Ynew = obj.applyQtfull(Ynew);   % Apply Qt to Ynew
            Ybot = Ynew(obj.k+1:end,:);     % Extract bottom of Ynew

            % Form QR decomposition of bottom
            [Ynew(obj.k+1:end,:),obj.tau(obj.k+1:obj.k+l)]...
                = hhqr(Ybot);

            obj.H(:,obj.k+1:obj.k+l) = Ynew; % Write to buffer
            obj.k = obj.k + l; % Increase size of k
        end
        
        function Qx = applyQ(obj, x)
            % Apply the Q matrix to a vector x
            Qx = matlab.internal.decomposition.applyHouseholder(...
                obj.H, obj.tau, [x;zeros(obj.n-obj.k,...
                size(x,2))], false, obj.k);
        end
        
        function Qtx = applyQt(obj, x)
            % Apply the transpose of the Q matrix to a vector x
            Qtx = apply_Qt(obj.H.obj.tau,x);
            Qtx = Qtx(1:obj.k,:);
        end
        
        function y = projectOut(obj, x)
            % Compute y = (I-Q*Q')*x
            y = apply_Qt(obj.H.obj.tau,x);
            y(1:obj.k,:) = 0;
            y = matlab.internal.decomposition.applyHouseholder(...
                obj.H, obj.tau, y, false, obj.k);
        end
        
        function Q = getQ(obj)
            % Return the Q matrix from the current QR decomposition
            Q = get_Q(obj.H,obj.tau);
        end
        
        function R = getR(obj)
            % Return the R matrix from the current QR decomposition
            R = triu(obj.H);  % Extract the upper triangular part
            R = R(1:obj.k,1:obj.k);
        end
    end
end
\end{lstlisting}

\chapter{Which sketch should I use?} \label{ch:which-sketch}

\epigraph{Our implementation of \textsf{CQRRPT} uses structured sparse sketching operators since these can provide exceptional speed without sacrificing reliability of the algorithm.}{Maksim Melnichenko, Oleg Balabanov, Riley Murray, James Demmel, Michael W.\ Mahoney, and Piotr Luszczek, \emph{CholeskyQR with randomization and pivoting for tall matrices (CQRRPT)} \cite{MBM+23}}


This appendix surveys the various constructions sketching matrices and compares them.
Echoing \cite{DM23}, my conclusion is that, among the available options, \emph{sparse sign embeddings} (\cref{sec:sse}) are the fastest and most reliable for general-purpose use.

\myparagraph{Sources}
This appendix is a significantly expanded version of a blog post I wrote on this topic \cite{Epp23c}, and it is informed by my numerical experience writing the papers \cite{Epp24a,EMN24}.
The original inspiration for this line of research was the numerical experiments of \cite{DM23}, which first demonstrated to me the significant speed advantage of sparse random embeddings over alternatives.

\myparagraph{Outline}
\Cref{sec:sketching-properties} begins by describing, in more detail, the kind of properties we desire for sketching matrices to be useful for applications in matrix computations.
The next sections provide a tour through standard options for sketching matrices: Gaussian embeddings (\cref{sec:gaussian-embed}), iid embeddings (\cref{sec:iid-embeddings}), subsampled randomized trigonometric transforms (\cref{sec:srtt}), and iid sparse embeddings (\cref{sec:iid-sparse}).
Each of these constructions is seen to have both merits and demerits.
We conclude our tour with sparse sign embeddings in \cref{sec:sse}, which we observe to have a balance of speed and accuracy (i.e., low distortion) not matched by other types of dimensionality maps.
We end with some concluding thoughts in \cref{sec:sketching_conclusions} and a postscript on recent developments in \cref{sec:postscript}.

\section{What properties do we want sketching matrices to have?} \label{sec:sketching-properties}

Before studying different types of sketching matrices, let us revisit the question: What properties does a sketching matrix need for use in matrix computations?
This section answers this question in three subsections.
First, we redefine the subspace embedding to allow different ``upper'' and ``lower'' distortion parameters, and we investigate the asymmetric importance of these parameters.
Next, we discuss the distribution of singular values for a sketched matrix.
Finally, we discuss the differences between subspace embeddings and Johnson--Lindenstrauss embeddings.

\subsection{Redefining subspace embeddings and the sketching asymmetry principle}

In \cref{ch:ls-algs-history}, we defined a subspace embedding $\mat{S} \in \field^{m\times d}$ for a matrix $\mat{F} \in \field^{m\times n}$ as a matrix for which 
\begin{equation} \label{eq:single-eta-subspace}
    (1-\eta) \norm{\mat{F}\vec{z}} \le \norm{\mat{S}^*(\mat{F}\vec{z})} \le (1+\eta) \norm{\mat{F}\vec{z}} \quad \text{for all } \vec{z} \in \field^n.
\end{equation}
The parameter $\eta \in [0,1]$, called the \emph{distortion}, measures the quality of the embedding.
This definition is somewhat lacking in that it uses the same parameter $\eta$ both upper and lower bound $\norm{\mat{S}^*(\mat{F}\vec{z})}$.
In this appendix, we will use the following two-parameter version of the subspace embedding property:

\begin{definition}[Subspace embeddings, again] \label{def:subspace-embed-2}
    A matrix $\mat{S} \in \field^{m\times d}$ is said to be a \emph{subspace embedding} for a matrix $\mat{F} \in \field^{m\times n}$ with \emph{lower distortion} $\eta_- \in [0,1)$ and \emph{upper distortion} $\eta_+ \ge 0$ provided that 
    \begin{equation*}
        (1-\eta_-) \norm{\mat{F}\vec{z}} \le \norm{(\mat{S}^*\mat{F})\vec{z}} \le (1+\eta_+) \norm{\mat{F}\vec{z}} \quad \text{for every } \vec{z} \in \field^n.
    \end{equation*}
    If $\eta \coloneqq \max \{\eta_-, \eta_+\} < 1$, we say that $\mat{S}$ is a subspace embedding for $\mat{F}$ with \emph{distortion} $\eta$.
\end{definition}

This simple change in definition reframes our perspective on what properties we should desire for a sketching matrix, and we have new versions of \cref{thm:sketch-and-solve,prop:randomized-preconditioning} consistent with this new definition:

\begin{theorem}[Sketch-and-solve, again]
    Let $\mat{S}$ be a subspace embedding \warn{for $\flatonebytwo{\mat{B}}{\vec{c}}$} with distortions $\eta_-,\eta_+$, and let $\hatvector{x}$ be the sketch-and-solve solution.
    Then 
    \begin{equation*}
        \norm{\vec{c} - \mat{B}\hatvector{x}} \le \max \left\{ \frac{1+\eta_+}{1-\eta_-}, \frac{4.5 \max\{\eta_+,\eta_-\}^2}{(1-\eta_-)^4}\right\} \cdot \norm{\vec{c} - \mat{B}\vec{x}}.
    \end{equation*}
\end{theorem}

\begin{proposition}[Randomized preconditioning, again] \label{prop:randomized-preconditioning-again}
    Let $\mat{S}$ be a subspace embedding \warn{for $\mat{B}$} with distortions $\eta_-,\eta_+$.
    Construct the sketched matrix $\mat{S}^*\mat{B}$, and let
    \begin{equation*}
        \mat{S}^*\mat{B} = \mat{U}\mat{M}
    \end{equation*}
    be \emph{any} orthonormal decomposition of $\mat{S}^*\mat{B}$, as in \cref{prop:randomized-preconditioning}.
    Then
    \begin{equation*}
        \sigma_{\mathrm{max}}(\mat{B}\mat{M}^{-1}) \le \frac{1}{1-\eta_-}, \quad \sigma_{\mathrm{min}}(\mat{B}\mat{M}^{-1}) \ge \frac{1}{1+\eta_+},\quad
        \cond(\mat{B}\mat{M}^{-1}) \le \frac{1+\eta_+}{1-\eta_-}.
    \end{equation*}
\end{proposition}

We see that the bounds on the sketch-and-solve residual norm and the condition number diverge as $\eta_- \uparrow 1$, but remain bounded when $\eta_+ \uparrow 1$.
This demonstrates a fundamental asymmetry in the importance of the lower and upper distortions:

\actionbox{\textbf{The sketching asymmetry principle.} For most applications of sketching in linear algebra, the most important requirement for a sketching matrix is that the \warn{lower distortion is bounded away from $1$}, e.g., $\eta_- \le 0.9$. Sketching algorithms typically produce meaningful results even with weak control on the upper distortion $\eta_+$, say, $\eta_+ \le \order(\log n)$.}

The sketching asymmetry principle is a recent perspective, as most foundational works on sketching and subspace embeddings are concerned with establishing the subspace embedding property with a single distortion parameter $\eta$, which is sometimes defined differently than we have (e.g., \cite[Def.~1.1]{CDDR24}).
The sketching asymmetry principle was developed and advocated for in recent work of Joel Tropp \cite{Tro25}, who suggests the name \emph{subspace injection}\index{subspace injection} for a sketching matrix satisfying $\eta_- < 1$.

\subsection{The distribution of singular values}
As the name suggests, a subspace embedding for $\mat{F}$ preserves the lengths of all vectors in the \emph{subspace} $\range(\mat{F})$.
The distortions $\eta_-,\eta_+$ describe the minimum ($1-\eta_-$) and maximum ($1+\eta_+$) factors by which $\mat{S}^*$ can rescale vectors in this subspace.
As the following result shows, these factors must always bound an interval containing there singular values of the matrix $\mat{S}^*\mat{Q}$, where $\mat{Q} = \orth(\mat{F})$:

\begin{proposition}[Singular values to subspace embedding] \label{prop:svals-subspace}
    Instate the notation of \cref{prop:randomized-preconditioning-again}.
    Let $\mat{F} \in \field^{m\times n}$ be a matrix, and let $\mat{Q}$ be \emph{any} matrix whose columns for an orthonormal basis for $\range(\mat{F})$.
    Then $\mat{S}\in\field^{m\times d}$ is a subspace embedding for $\mat{F}$ with distortions $\eta_-,\eta_+$ if and only if 
    \begin{equation*}
        1-\eta_- \le \sigma_{\mathrm{min}}(\mat{S}^*\mat{Q}) \le \sigma_{\mathrm{max}}(\mat{S}^*\mat{Q}) \le 1+\eta_+.
    \end{equation*}
\end{proposition}

This result is standard, and its proof is straightforward; see, e.g., \cite[Prop.~5.2]{KT24}.

Can the \emph{distribution} of the singular values of $\mat{S}^*\mat{Q}$ tell us anything about sketching in matrix computations?
Perhaps.
We have the following result:

\begin{proposition}[Singular values and randomized preconditioning]
    Instate the notation of \cref{prop:randomized-preconditioning-again}.
    The singular values of the preconditioned matrix $\mat{B}\mat{M}^{-1}$ are the reciprocals of the singular values of $\mat{S}^*\mat{Q}$.
\end{proposition}

\begin{proof}
    Let $\mat{B} = \mat{Q}(\mat{B}^*\mat{B})^{1/2}$ be the polar decomposition, and let $\mat{S}^*\mat{B} = \mat{U}\mat{M}$ be any orthonormal decomposition.
    Now, taking a polar decomposition $\mat{M} = \mat{V} (\mat{B}^*\mat{S}^*\mat{S}\mat{B})^{1/2}$, we have
    \begin{equation*}
        \mat{B}\mat{M}^{-1} = \mat{Q}(\mat{B}^*\mat{B})^{1/2}(\mat{B}^*\mat{S}\mat{S}^*\mat{B})^{-1/2}\mat{V}^*.
    \end{equation*}
    Therefore, the squared singular values of $\mat{B}\mat{M}^{-1}$ are 
    \begin{equation*}
        \vec{\sigma}\bigl(\mat{B}\mat{M}^{-1}\bigr)^2 = \vec{\lambda}\left((\mat{B}^*\mat{B})^{1/2}(\mat{B}^*\mat{S}\mat{S}^*\mat{B})^{-1}(\mat{B}^*\mat{B})^{1/2}\right).
    \end{equation*}
    (Recall that we extend nonlinear functions to vectors entrywise.)
    Since $\mat{S}^*\mat{B} = \mat{S}^*\mat{Q}(\mat{B}^*\mat{B})^{1/2}$, we have the simplification:
    \begin{equation*}
        \vec{\sigma}\bigl(\mat{B}\mat{M}^{-1}\bigr)^2 = \vec{\lambda}\left((\mat{Q}^*\mat{S}\mat{S}^*\mat{Q})^{-1}\right) = \vec{\sigma}(\mat{S}^*\mat{Q})^{-2}.
    \end{equation*}
    Conclude by observing that the eigenvalues of $(\mat{Q}^*\mat{S}\mat{S}^*\mat{Q})^{-1}$ are the inverses of the eigenvalues of $\mat{Q}^*\mat{S}\mat{S}^*\mat{Q}$, which are the inverse-squares of the singular values of $\mat{S}^*\mat{Q}$.
\end{proof}

This result shows that the singular values of $\mat{S}^*\mat{Q}$ determine the singular values of the preconditioned matrix $\mat{B}\mat{M}^{-1}$.
In principle, knowing the singular values of $\mat{B}\mat{M}^{-1}$ can lead to bounds on the convergence of sketch-and-precondition-type algorithms that are more precise than the simple bounds we derived in \cref{cor:sketch_and_pre}.
See \cite[\S3.1]{Gre97a} for a discussion of such \emph{spectrum-adaptive} bounds for Krylov iterative methods.
(Be aware that these spectrum-adaptive are valid in exact arithmetic but often fail in finite-precision arithmetic.)

In my experience, I understanding the precise distribution of the singular values $\mat{S}^*\mat{Q}$ is usually not \emph{that} informative for studying the behavior of linear algebraic algorithms.
Still, plotting the distribution of singular values of $\mat{S}^*\mat{Q}$ can be a valuable way of understanding the behavior of sketching matrices, particularly those that \emph{lack} the subspace embedding property for given parameters $\eta_+$, $\eta_-$.
Are there just a few stray singular values that escape the interval $[1-\eta_-,1+\eta_+]$?
A small clump?
Or does the bulk of singular values extend past the endpoints of the interval?
These are all different ways an embedding can fail to have the subspace embedding property, and visualizing the singular values can help diagnose \emph{how} a sketching matrix fails to have the subspace embedding property.

\subsection{Subspace embeddings vs.\ Johnson--Lindenstrauss embeddings}

In matrix computations, the term ``sketching'' is generally accepted to have a very broad definition.
To some, \emph{any} form of \warn{linear} dimensionality reduction is an instance of sketching.
Typically, we expect that the sketching process should ``preserve lengths'' of the matrix or vectors being sketched.
The subspace embedding property (\cref{def:subspace-embed-2}) is one way of formalizing this length-preservation property, but there are others.
In particular, we have the following standard definition:

\begin{definition}[Johnson--Lindenstrauss embedding]
    A matrix $\mat{S} \in \field^{m\times d}$ is said to be a \emph{Johnson--Lindenstrauss} embedding\index{Johnson--Lindenstrauss embedding} with distortion $\eta \in [0,1)$ for a \warn{finite} subset $\set{E} \subseteq \field^m$ if 
    \begin{equation*}
        (1-\eta)\norm{\vec{v}} \le \norm{\mat{S}^*\vec{v}} \le (1+\eta) \norm{\vec{v}} \quad \text{for every } \vec{v} \in \set{E}.
    \end{equation*}
\end{definition}

The Johnson--Lindenstrauss lemma\index{Johnson--Lindenstrauss lemma} \cite{JL84} states that, e.g., Gaussian embeddings (\cref{sec:gaussian-embed}) are \warn{oblivious} Johnson--Lindenstrauss embeddings\index{Johnson--Lindenstrauss embedding} with embedding dimension $d = \order(\eta^{-2} \log |\set{E}|)$ (with failure probability, say, $1/n$) \cite[Lem.~18]{Woo14a}.
We see the embedding dimension need only be \emph{logarithmic} in the \warn{cardinality} of the set of points $\set{E}$.

Both the Johnson--Lindenstrauss embedding\index{Johnson--Lindenstrauss embedding} property and subspace embedding property are useful conditions for a sketching matrix to satisfy.
But it is important not to confuse them.
Johnson--Lindenstrauss embeddings\index{Johnson--Lindenstrauss embedding} preserve lengths of vectors in a \warn{finite set} $\set{E}$ and the embedding dimension $d\sim (\log |\set{E}|)/\eta^2$ must be \warn{logarithmic} in the \warn{cardinality $|\set{E}|$}.
Subspace embeddings preserve lengths of vectors in a \warn{finite-dimensional subspace} $\range(\mat{F})$ and the embedding dimension $d \sim \operatorname{dim} \range(\mat{F})/\eta^2$ must be \warn{linear} in the \warn{dimension $\operatorname{dim} \range(\mat{F})$}.
In \cref{part:sketching} of this thesis, we are exclusively concerned with the subspace embedding property.

\index{Gaussian embedding|(}
\section{Gaussian embeddings} \label{sec:gaussian-embed}

To understand sketching, the best place to start is with Gaussian embeddings, which admit a beautiful and mathematically precise theory.
While Gaussian embeddings are often less computationally efficient than other types of embeddings, they set the standard to which other types of dimensionality reduction maps will aspire.
\begin{definition}[Gaussian embedding]
    A \emph{Gaussian embedding} (over $\field$) is a matrix $\mat{S} \in \field^{m\times d}$ with iid $\Normal_\field(0,1/d)$ entries.
\end{definition}
The Gaussian distribution has many beautiful properties.
These properties can be applied to give very sharp quantitative analysis of Gaussian embeddings.
In particular, we will employ the following fundamental fact:

\actionbox{For any matrix $\mat{Q} \in \field^{m\times n}$ with orthonormal columns and a Gaussian matrix $\mat{S} \in \field^{m\times d}$ over $\field$, the sketch $\mat{S}^*\mat{Q} \in \field^{d\times n}$ has iid $\Normal_\field(0,1/d)$ entries.}

On the basis of this result, analysis of Gaussian embeddings reduces to questions about the singular values of rectangular matrices with iid Gaussian entries, which are among the most well-studied in random matrix theory.

\myparagraph{Asymptotic theory}
In the large-data limit $n,d \to \infty$, the behavior of Gaussian embeddings is captured by standard asymptotic results from random matrix theory.
We have the following result:

\begin{fact}[Gaussian embeddings, distribution of singular values: large-data limit] \label{fact:gaussian-bulk-asymptotic}
    Consider a family of matrices $\mat{Q}_i \in \field^{m_i\times n_i}$ and Gaussian embeddings $\mat{S}_i \in \field^{m_i\times d_i}$ over $\field$ with dimensions $m_i,n_i,d_i \to \infty$ tending to infinity and with limiting aspect ratio $n_i / d_i \to \varrho \in (0,1)$.
    Then the distribution of singular values of $\mat{S}_i^*\mat{Q}_i^{\vphantom{*}}$ converges (weakly in probability) to the \emph{Mar\v{c}enko--Pastur singular value distribution},\index{Mar\v{c}enko--Pastur singular value distribution} a continuous probability distribution on $\real_+$ with density
    \begin{equation*}
        f_{\mathrm{mp}}(\sigma) = \frac{\sqrt{(\sigma_+^2-\sigma^2)_+(\sigma^2 - \sigma_-^2)_+}}{\pi \varrho \sigma}.
    \end{equation*}
    Here, 
    \begin{equation*} 
        \sigma_\pm = 1 \pm \sqrt{\varrho}
    \end{equation*}
    and $a_+ \coloneqq \max \{a,0\}$ denotes the positive part of a real number.
\end{fact}

\Cref{fact:gaussian-bulk-asymptotic} is a deep fact of random matrix theory, proven (in a somewhat different form) in Mar\v{c}enko and Pastur's seminal paper \cite{MP67}.
Its proof requires a significant mathematical apparatus.

This result demonstrates that the \emph{distribution} of singular values converges to a continuous distribution supported on the interval $[\sigma_-,\sigma_+]$, but it does not rule out the possibility of a few extraneous singular values exiting this interval.
The following result \cite{BY93}\index{Bai--Yin law} forbids it.

\begin{fact}[Gaussian embeddings, extreme singular values: large-data limit]
    Assume the setting of \cref{fact:gaussian-bulk-asymptotic} and assume in addition that $\mat{S}_i$ are \emph{independent}.
    Then $\sigma_{\mathrm{max}}(\mat{S}_i^*\mat{Q}_i) \to \sigma_+$ and $\sigma_{\mathrm{min}}(\mat{S}_i^*\mat{Q}_i) \to \sigma_-$, almost surely as $i\to\infty$.
    In particular, the upper and lower distortions $\eta_+,\eta_-$ both converge almost surely to $\sqrt{\varrho} = \sqrt{n/d}$.
\end{fact}

Rearranging the relation $\eta \approx \sqrt{n/d}$, we derive the relation
\begin{equation*}
    d \approx n / \eta^2.   
\end{equation*}
This gives us guidance for how to set the embedding dimension $d$ for a Gaussian embedding to achieve a desired distortion $\eta$.

\begin{figure}
    \centering
    \includegraphics[width=0.98\linewidth]{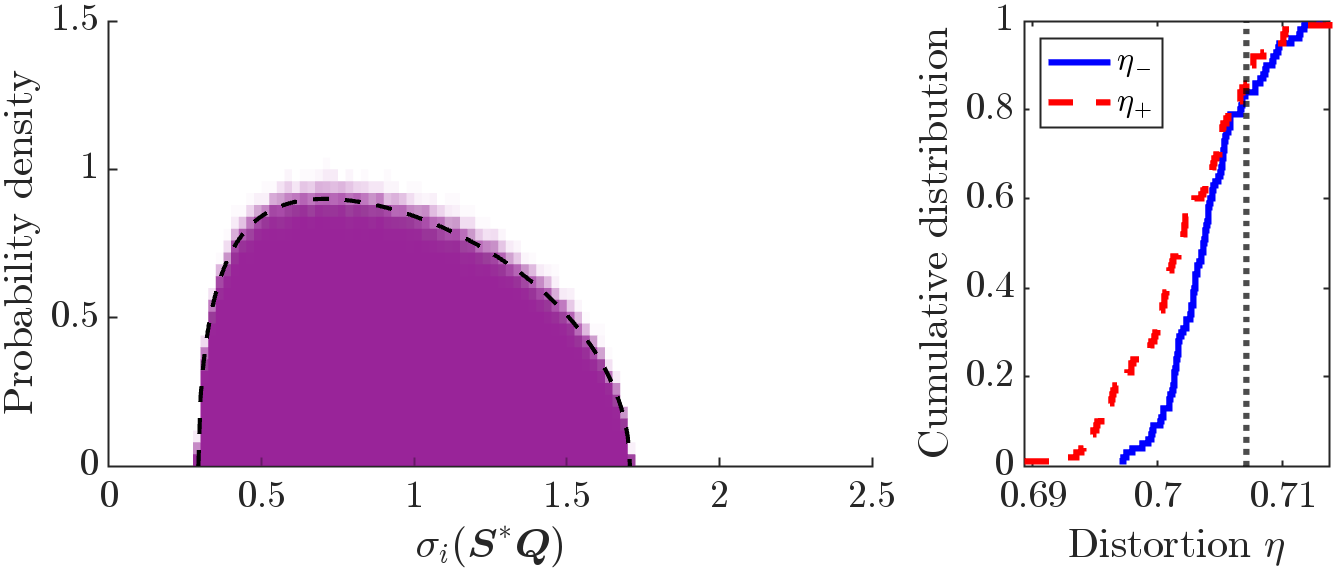}
    \caption[Singular values $\vec{\sigma}(\mat{S}^*\mat{Q})$ and distortions for Gaussian embeddings]{\emph{Left:} Superimposed singular value distributions $\vec{\sigma}(\mat{S}^*\mat{Q})$ for 100 independent realizations of a Gaussian embedding applied to a matrix $\mat{Q} \in \real^{m\times n}$ with dimension $n = 1000$.
    The embedding dimension is $d=2n$.
    The Mar\v{c}enko--Pastur singular value distribution is shown as a black dashed line for reference.
    \emph{Right:} Empirical cumulative distribution function of lower and upper distortions.
    The asymptotic limiting value $\eta = \sqrt{n/d}$ is shown in black.}
    \label{fig:gaussian-sketch}
\end{figure}

\myparagraph{Empirical results}
The asymptotic theory for Gaussian embeddings accurately predicts the empirical results.
The left panel of \cref{fig:gaussian-sketch} shows the distribution of singular values $\vec{\sigma}(\mat{S}^*\mat{Q})$ for a Gaussian embedding over 100 independent random trials.
We set $n \coloneqq 1000$ and $d \coloneqq 2n$.
(The value of $m$ and the matrix $\mat{Q}$ are immaterial, owing to the boxed fact above.)
Each trial is plotted as a single transluscent histogram, and the histograms for all 100 trials are superimposed.
The empirically computed singular value distributions show remarkable fidelty to the Mar\v{c}enko--Pastur singular value distribution,\index{Mar\v{c}enko--Pastur singular value distribution} marked as a dotted black line.
Not one of the 100 random realizations shows a significant excursion from the the bounds of the Mar\v{c}enko--Pastur density curve.

The right panel plots the cumulative distribution of the distortion parameters $\eta_+$ and $\eta_-$ for these 100 trials
The distortion parameters are consistently close to the almost-sure limit value $\eta = \sqrt{n/d} = 1/\sqrt{2}$, marked as a dotted vertical line.

\myparagraph{Nonasymptotic results}
In practice, we do not apply sketching matrices to a series of problems asymptotically growing to infinity; we apply sketching to a given problem of a fixed size.
As such, we value theoretical results that predict the performance of Gaussian embeddings for a specific problem size.

\begin{fact}[Gaussian embeddings, extreme singular values: nonasymptotic results] \label{fact:gaussian-nonasymptotic}
    Let $\mat{S} \in \real^{m\times d}$ be a \warn{real} Gaussian embedding, let $\mat{B} \in \real^{m\times n}$ be a \warn{real} matrix, and let $\eta_+,\eta_-$ be the distortions.
    Then 
    \begin{equation*}
        \expect [\eta_+],\expect[\eta_-]\le \sqrt{n/d}.
    \end{equation*}
    Furthermore, we have the tail bound
    \begin{equation*}
        \prob \big\{ \eta_+ \ge \sqrt{n/d} + t \big\}\le \exp \left( - \frac{dt^2}{2}  \right) \quad \text{and} \quad\prob \big\{ \eta_- \ge \sqrt{n/d} + t \big\} \le \exp \left( - \frac{dt^2}{2}  \right).
    \end{equation*}
\end{fact}

This result is due to \cite{DS01}.
See \cite[Cor.~6.38]{AS17} for an exposition and \cite[Thm.~8.4]{MT20a} for a generalization.
It's proof requires many techniques.
The bound on $\expect[\eta_+]$ is a direct consequence of Chevet's theorem \cite{Che77}, and the bound on $\expect[\eta_-]$ follows from Gordon's inequalities \cite{Gor85}.
The tail bounds follow from concentration for Lipschitz continuous functions of Gaussian random variables \cite[Thm.~3.25]{van14}.
I recommend the references \cite{van14,AS17,Ver18,Tro21} for learning about this family of ideas.
A similar (but slightly worse) result is also known for complex Gaussian embeddings \cite[Prop.~6.3.3]{AS17}.\index{Gaussian embedding|)}

\myparagraph{Universality}
As we will see in the rest of this chapter, the most effective random embeddings tend to exhibit performance similar to Gaussian embeddings, a manifestation of the probabilistic phenomenon of universality.
Providing mathematical support to this type of universality phenomenon has been an active subject for research (see, e.g., \cite{OT18,DL19a,Bv22,CDD24,CDDR24,Tro25}), but even the best-known mathematical tools are not precise enough to fully explain empirically observed universality phenomena.
Notwithstanding these theoretical gaps, Martinsson and Tropp recommend using the Gaussian theory to inform parameter decisions \cite[\S9.7]{MT20a}:

\begin{quote}
    Even so, we would like to have a priori predictions about how our algorithms will behave.
    Beyond that, we need reliable methods for selecting algorithm parameters, especially in the streaming setting where we cannot review the data and repeat the computation.
    
    Here is one answer to these concerns.
    As a practical matter, we can simply invoke the lessons from the Gaussian theory, even when we are using a different type of random embedding. The universality result [of \cite{OT18}] gives a rationale for this approach in one special case.
    We also recommend undertaking computational experiments to verify that the Gaussian theory gives an adequate description of the observed behavior of an algorithm.
\end{quote}

I fully agree with this viewpoint, and the experiments provided in this chapter should provide insights to practitioners interested in the differences between types of random embeddings.

\section{IID embeddings} \label{sec:iid-embeddings}

\begin{figure}
    \centering
    \includegraphics[width=0.98\linewidth]{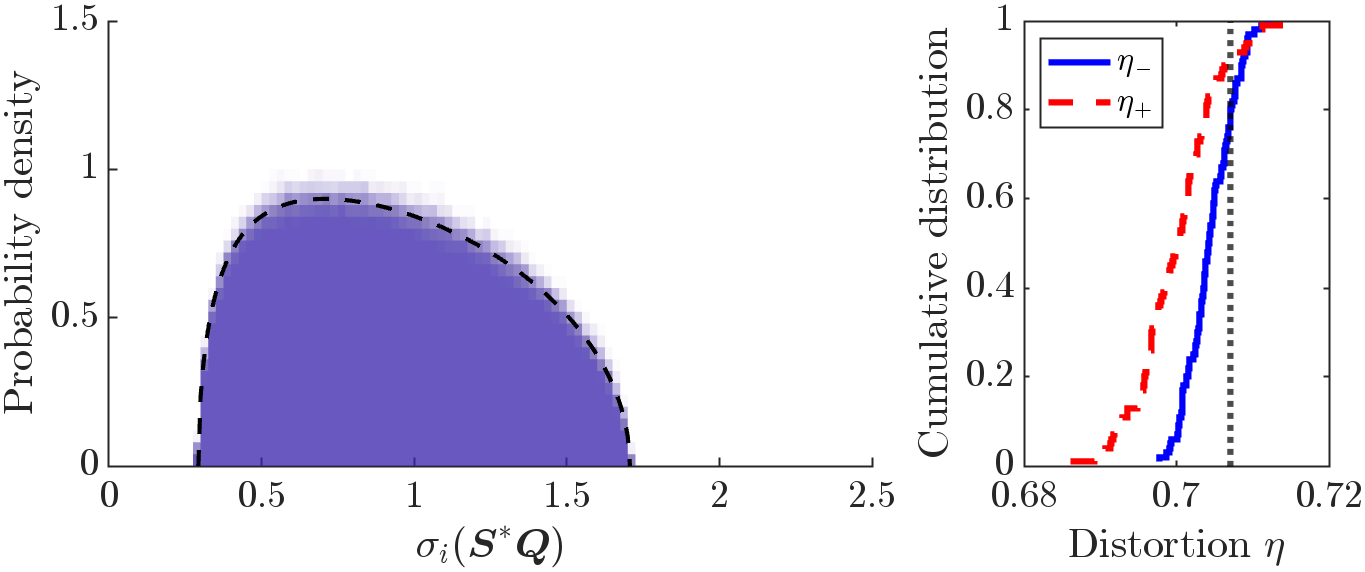}
    \caption[Singular values $\vec{\sigma}(\mat{S}^*\mat{Q})$ and distortions for iid scaled sign embeddings]{\emph{Left:} Superimposed singular value distributions $\vec{\sigma}(\mat{S}^*\mat{Q})$ for 100 independent realizations of a iid scaled sign embedding applied to the adversarial matrix \cref{eq:difficult-example-Q} with embedding dimension $d=2n$.
    The Mar\v{c}enko--Pastur singular value distribution is shown as a black dashed line for reference.
    \emph{Right:} Empirical cumulative distribution function of lower and upper distortions.
    The asymptotic limiting value $\eta = \sqrt{n/d}$ is shown in black.}
    \label{fig:rademacher-sketch}
\end{figure}

One simple way of constructing a non-Gaussian sketching matrix is to replace the Gaussian distribution with another random distribution, such as the (scaled) random sign distribution $\Unif \{\pm d^{-1/2}\}$.
As long as the new distribution has mean zero, variance $1/d$, and is sufficiently well-behaved (e.g., exhibits rapidly decaying subgaussian tails), the performance is generally comparable to a Gaussian embedding.
\Cref{fig:rademacher-sketch} shows the behavior of a sketching matrix with iid $\Unif \{\pm d^{-1/2}\}$ entries applied to the matrix
\begin{equation} \label{eq:difficult-example-Q}
    \mat{Q} \coloneqq \twobyone{\Id_n}{\mat{0}} \in \real^{m\times n} \quad \text{for } m \coloneqq 10^5, \: n \coloneqq 10^3;
\end{equation}
we use embedding dimension $d=2n$.
(This matrix presents a hard case for certain other sketching matrices, as we will see below.)
The results are essentially identical to those for Gaussian embeddings in \cref{fig:gaussian-sketch}, providing a clear demonstration of the universality phenomenon.

\section{Subsampled trigonometric transforms} \label{sec:srtt}

Subsampled randomized trigonometric transforms (SRTTs) are widely viewed as among the most effective types of sketching matrices for use in matrix computations.
Speaking to other researchers, one hears claims such as ``SRTTs behave exactly like Gaussian embeddings''.
Do these claims hold up to scrutiny?
We shall see.

Let us begin with the definition:

\begin{definition}[Subsampled randomized trigonometric transform] \label{def:srtt}
    A (standard) \emph{subsampled randomized trigonometric transform} is a sketching matrix $\mat{S} \in \field^{m\times d}$ defined as a scaled product of three matrices.
    \begin{equation*}
        \mat{S}^* = \sqrt{\frac{m}{d}}\mat{R}\mat{F}\mat{D}.
    \end{equation*}
    These matrices have the following definitions:
    \begin{itemize}
        \item $\mat{D} = \Diag(\vec{\varepsilon})$ is a diagonal matrix with entries $\varepsilon_1,\ldots,\varepsilon_m$ drawn iid from either the random sign $\Unif \{\pm 1\}$ or random phase $\Unif \unitcircle(\complex)$ distributions.
        \item $\mat{F}$ is a unitary fast trigonometric transform for which $\mat{F}\vec{z}$ can be computed in $\order(m\log m)$ operations.
        Examples include the discrete Fourier transform, the (type-II) discrete cosine transform, the (Walsh--)Hadamard transform, or the Hartley transform.
        \item $\mat{R}$ is a restriction to $d$ coordinates, sampled uniformly \warn{without replacement}.
    \end{itemize}
\end{definition}

This construction is due to Ailon and Chazelle \cite{AC06} (see also \cite{WLRT08}).
SRTTs are often called SRFTs or SRHTs when $\mat{F}$ is the discrete Fourier transform or (Walsh--)Hadamard transform.
In fact, the name SRFT is often used even when with other types of trigonometric transforms.
In this work, we exclusively use the discrete cosine transform as our trigonometric transform $\mat{F}$ and the random signs distribution for $\mat{D}$.

\myparagraph{Runtime}
The operation count for sketching using SRTTs is smaller than for Gaussian embeddings.
Since we use a \emph{fast} trigonmetric transform $\mat{F}$, the sketch $\mat{S}^*\vec{v} = \sqrt{n/d} \cdot \mat{R}(\mat{F}(\mat{D}\vec{v}))$ can be computed in $\order(m\log m)$ operations.
Compare with $\order(dm)$ operations for Gaussian embeddings.
In principle, the runtime for SRTTs can be reduced to $\order(m \log d)$ operations \cite{WLRT08}.

Fast trigonometric transforms are notoriously difficult to implement \cite{Van92}.
Even with a world-class implementation, the computational throughput of fast trigonometric transforms (as measured by floating point operations per second) is much less than for matrix multiplication.\index{matrix--matrix operations}
For sufficiently large problems, SRTTs are substantially faster than Gaussian embeddings, but the speedups can be less impressive than one might hope.
In my experience, sparse sketching matrices (\cref{sec:iid-sparse,sec:sse}) are typically significantly faster than SRTTs.

\myparagraph{Analysis}
Guarantees for the \warn{SRHT} were proven by Tropp \cite{Tro11}.
Here is a simplified version of his result:

\begin{fact}[Standard SRHTs: Subspace embedding property] \label{fact:srhts}
    Standard SRHTs of size $m\times d$ are oblivious subspace embeddings for dimension $n$ with distortion $\eta = 0.6$ and failure probability $\delta = \order(1/n)$ provided that
    \begin{equation*}
        d \ge \mathrm{const} \cdot (n + \log m) \log n.
    \end{equation*}
\end{fact}

The full versions of Tropp's results have explicit constants.
The analysis of \cite{Tro11} can be extended to prove subspace embedding bounds for general distortions $\eta > 0$ (such a resulted is stated in \cite[Thm.~7]{Woo14a}) and adapted to the random phase distribution for $\mat{D}$ or other types of trigonometric transforms.

Tropp's result suggests (correctly) that an embedding dimension of {$d\sim$}\,\:\warn{$n\log n$} is \warn{necessary} for standard SRTTs to achieve the subspace embedding property (for constant $\eta < 1$).
Compare with Gaussian embeddings, which allow {$d\sim$}\,\:\warn{$n$}, with no need for a logarithmic oversampling factor.

\myparagraph{Bad examples}
If we want to use sketching matrices as primitives in general-purpose software, then we must demand that they work reliably when applied to \emph{any} input matrix.
As such, even a single bad example can dissuade us from using a particular type of dimensionality reduction map.

\begin{figure}
    \centering
    \includegraphics[width=0.98\linewidth]{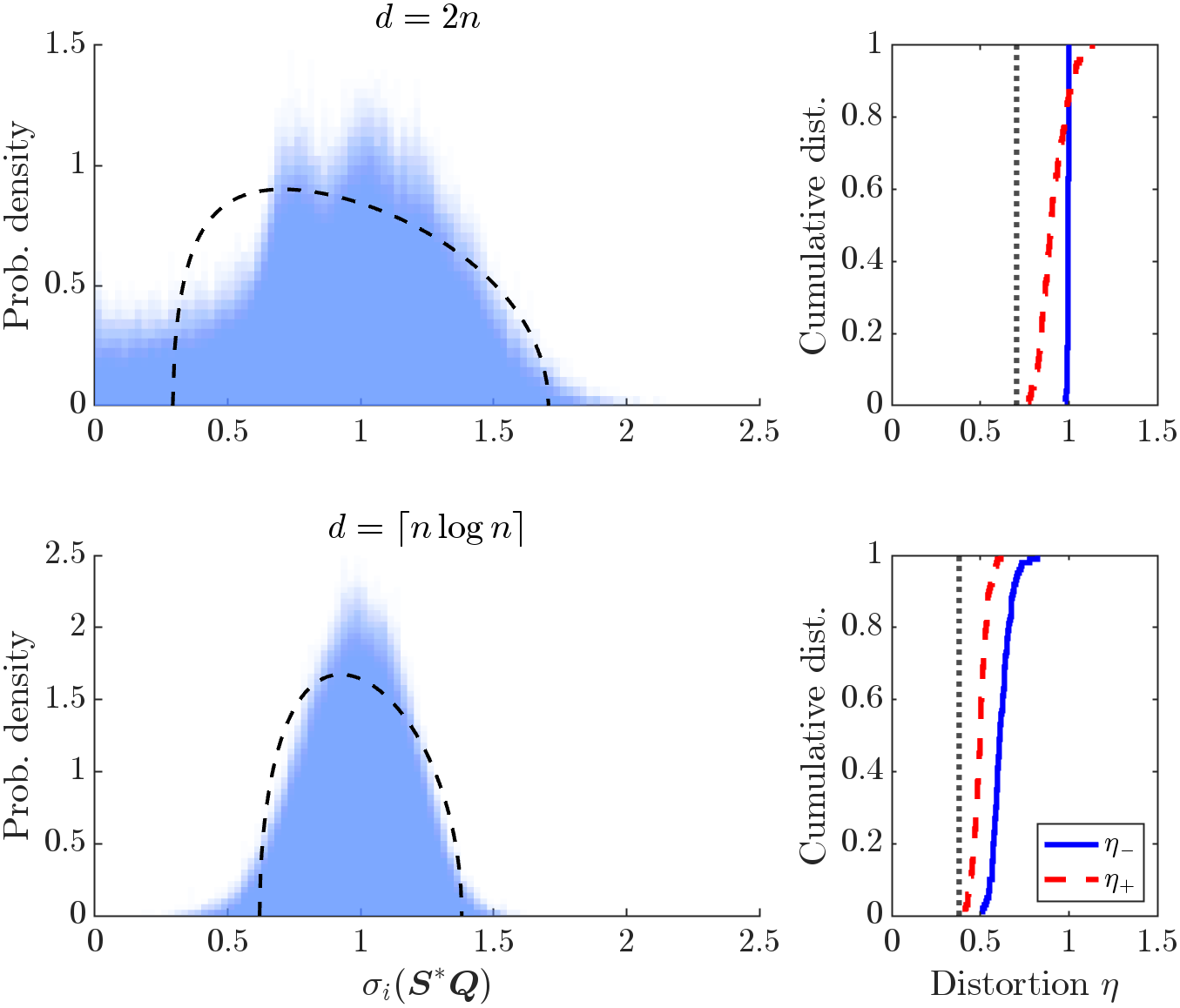}
    \caption[Singular values $\vec{\sigma}(\mat{S}^*\mat{Q})$ and distortions for standard SRTT embeddings applied to a tough matrix $\mat{Q}$]{\emph{Left:} Superimposed singular value distributions $\vec{\sigma}(\mat{S}^*\mat{Q})$ for 100 independent realizations of standard SRTT applied to the adversarial matrix \cref{eq:difficult-example-Q} with embedding dimension $d=2n$ (\emph{top}) and $d= \lceil n \log n \rceil$ (\emph{bottom}).
    The Mar\v{c}enko--Pastur singular value distribution is shown as a black dashed line for reference.
    \emph{Right:} Empirical cumulative distribution function of lower and upper distortions.
    The asymptotic limiting value $\eta = \sqrt{n/d}$ is shown in black.}
    \label{fig:srtt-1-sketch}
\end{figure}

With that context, we evaluate the adversarial example matrix $\mat{Q}$ defined in \cref{eq:difficult-example-Q}.
Results with embedding dimension $d=2n$ are shown in the top panels of \cref{fig:srtt-1-sketch}.
We see that, on this difficult example, the embedding has lower distortion $\eta_- \approx 1$ with high probability.
Thus, the standard SRTT \warn{with $d \sim n$} is not effective for this type of matrix.

Consistent with \cref{fact:srhts}, we may remedy this issue by setting $d\sim n \log n$.
The bottom panels of \cref{fig:srtt-1-sketch} show the same experiments with $d=\lceil n \log n \rceil$.
We now see that we reliably obtain an embedding with $\eta < 0.9$, although the singular value statistics still deviate significantly from the appropriate Mar\v{c}enko--Pastur singular value distribution.\index{Mar\v{c}enko--Pastur singular value distribution}

One can also confirm the failures of standard SRTTs with $d \ll n \log n$ theoretically.
The following result appears in \cite[\S3.3]{Tro11}:

\begin{fact}[Standard SRHTs: Subspace embedding property]
    For $n\ge 1$, construct the orthnormal matrix $\mat{Q} \coloneqq \Id_n \otimes \evec_1 \in \real^{n^2\times n}$, and consider the action of a standard SRHT with embedding dimension $d$ on $\mat{Q}$.
    Then $\eta_- = 1$ with high probability unless $d \ge \Omega(n\log n)$.
\end{fact}

This failure mode is a consequence of the coupon collector phenomenon\index{coupon collector problem} \cite[\S3.6]{MR95}; see \cite[\S3.3]{Tro11} for details.

The bad examples for the SRTT are unfortunate, and their existence appears to not be widely known.
SRTTs are routinely used with a small embedding dimension (say, $d=2n$ and $d=3n$), accompanied implicitly or explicitly by claims that they are ``just as good a Gaussian embedding''.

\myparagraph{Rerandomization}
The performance of subsampled randomized trigonometric transforms can be greatly improved by adding an additional layer of randomness.
We make the following definition:

\begin{definition}[Rerandomized subsampled randomized trigonometric transform] \label{def:rerandomized-srtt}
    A \emph{rerandomized subsampled randomized trigonometric transform} is a sketching matrix $\mat{S} \in \field^{m\times d}$ defined so that
    \begin{equation*}
        \mat{S}^* = \sqrt{\frac{m}{d}}\mat{R}\mat{F}\mat{D}_1 \mat{F}\mat{D}_2,
    \end{equation*}
    where $\mat{F}$ and $\mat{R}$ are as in \cref{def:srtt} and $\mat{D}_1$ and $\mat{D}_2$ are \warn{independent} random diagonal matrices populated with random signs or random phases.
\end{definition}

A version of the rerandomization trick using random Givens rotations appears in \cite{RT08}.
The \textsf{Blendenpik}\index{Blendenpik@\textsf{Blendenpik}} software from the paper \cite{AMT10} contains a version of rerandomization more similar to \cref{def:rerandomized-srtt}.
Another version of rerandomization is advocated in the paper \cite{TYUC17b}.

\begin{figure}
    \centering
    \includegraphics[width=0.98\linewidth]{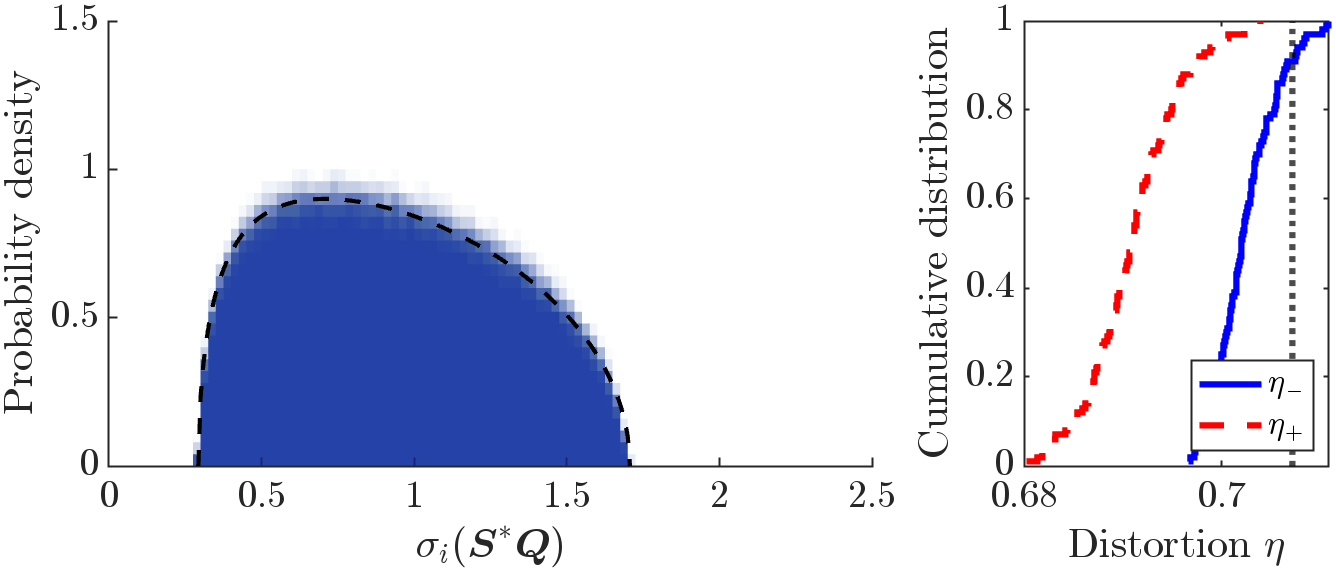}
    \caption[Singular values $\vec{\sigma}(\mat{S}^*\mat{Q})$ and distortions for \emph{rerandomized} SRTT embeddings applied to a tough matrix $\mat{Q}$]{\emph{Left:} Superimposed singular value distributions $\vec{\sigma}(\mat{S}^*\mat{Q})$ for 100 independent realizations of a rerandomized SRTT applied to the adversarial matrix \cref{eq:difficult-example-Q} with embedding dimension $d=2n$.
    The Mar\v{c}enko--Pastur singular value distribution is shown as a black dashed line for reference.
    \emph{Right:} Empirical cumulative distribution function of lower and upper distortions.
    The asymptotic limiting value $\eta = \sqrt{n/d}$ is shown in black.}
    \label{fig:srtt-2-sketch}
\end{figure}

Empirically, the rerandomized SRTT appears to achieve the touted ``just as good as Gaussian'' performance of SRTTs.
\Cref{fig:srtt-2-sketch} shows the results of rerandomized SRTTs applied to the adversarial matrix \cref{eq:difficult-example-Q}.
The singular values are observed to follow Mar\v{c}enko--Pastur statistics,\index{Mar\v{c}enko--Pastur singular value distribution} and the distortion is quite similar to Gaussian embeddings---or even slightly better.

Unfortunately, theoretical guarantees for the near-Gaussian behavior of rerandomized SRTTs are not known.
This is a natural subject for future work.

\begin{remark}[Adding even more randomness]
    Another variant of SRTTs, adding further randomness, is to use random \emph{signed permutation matrices} instead of diagonal matrices.
    This idea can mollify, but not fix, the issues with the standard SRTT (\cref{def:srtt}).
    It is also compatible with the rerandomized SRTT (\cref{def:rerandomized-srtt}), but I am unaware of an example where this trick substantially improves on rerandomized SRTTs.
\end{remark}

\myprogram{Class implementing subsampled randomized trigonmetric transforms.}{The argument \texttt{rounds} sets the number of rounds of randomization (\texttt{1} for standard SRTTs, \texttt{2} for rerandomized SRTTs). The subroutine \texttt{random\_signs} is provided in \cref{prog:random_signs}.}{srtt}

\myparagraph{Implementation}
A MATLAB class implementing SRTTs is shown in \cref{prog:srtt}.
We use random signs for the diagonal matrices and discrete cosine transforms as the fast trigonometric transform.
The first arguments to the constructor are the dimensions $d$ and $m$, and the third argument is the number of rounds of randomization, set to \texttt{1} for standard SRTTs (\cref{def:srtt}) or \texttt{2} for rerandomized SRTTs (\cref{def:rerandomized-srtt}).
We have overloaded the multiplication routine (\texttt{mtimes}), so $\mat{S}^*\vec{v}$ can be computed as \texttt{St * v}.

\section{IID sparse embeddings} \label{sec:iid-sparse}

With rerandomization, SRTTs appear to have distortions and singular value distributions that closely match Gaussian embeddings, even when applied to challenging instances.
Despite their $\order(m\log m)$ runtimes, these embeddings are often not as fast as one would hope; see \cref{fig:sketching-runtime-compare}.

Another way to obtain a sketching matrix with an efficient apply operation $\mat{B}\mapsto \mat{S}^*\mat{B}$ is through \emph{sparse} constructions.
In this section, we will explore one class of sparse random embeddings with \emph{iid} entries.
These embeddings already have some appealing properties, and they are (comparatively) easy to analyze.
In next section, we will consider another class of sparse random embeddings that have even better properties.
We make the following definition:

\begin{definition}[IID sparse embedding]
    An \emph{iid sparse embedding} $\mat{S} \in \real^{m\times d}$ with \emph{expected sparsity} $\zeta \in [0,d]$ is a random matrix with iid entries drawn from the \emph{scaled sparse sign distribution}
    \begin{equation*}
        s = \begin{cases}
            +\zeta^{-1/2} & \text{with probability } \zeta/2d, \\
            0,            & \text{with probability } 1-\zeta/d, \\
            -\zeta^{-1/2} & \text{with probability } \zeta/2d.
        \end{cases}
    \end{equation*}
\end{definition}

As the name suggests, the parameter $\zeta$ is equal to the expected number of nonzero entries in each row of $\mat{S}$.
As such, the expected number of nonzero entries in the entire matrix is $\nnz(\mat{S}) = \zeta m$, and the expected fraction of nonzero entries is $\nnz(\mat{S}) / md = \zeta / d$.


\myparagraph{Runtime and implementation}
We will not focus on how to implement and analyze the runtime of iid sparse embeddings, because sparse sign embeddings (\cref{sec:sse}) are more effective in practice.

IID sparse maps can be generated efficiently by (1) randomly generating the number $\nnz(\mat{S})$ (which has a binomial distribution with parameters $nd$ and $\zeta/d$), (2) generating the positions of those entries, and (3) then generating the values $\pm \zeta^{-1/2}$.
There are several ways to implement step (2); one way is to randomly generate random positions $(i,j) \sim \Unif (\{1,\ldots,d\} \times \{1,\ldots,m\})$ and add them to a hashtable until $\nnz(\mat{S})$ distinct positions are found.
(In the unusual case where $\nnz(\mat{S}) \ge md/2$, one should instead select random positions to \emph{disinclude}).
The expected runtime is $\order(\zeta m)$.
Once generated, iid sparse maps can be applied to a vector in $\order(\nnz(\mat{S}) + d)$ operations, so the expected runtime is $\order(\zeta m + d)$ operations.

\myparagraph{Analysis}
The subspace \warn{injection} properties of iid sparse embeddings were studied in a recent paper of Tropp \cite[Thm.~6.3]{Tro25}.
Here is a simplified version of his result:
\begin{fact}[IID sparse embeddings: Subspace injection] \label{fact:iid-sparse}
    Let $\mat{B} \in \real^{m\times n}$ be a \warn{real} matrix, and fix $\eta_- > 0$.
    An iid sparse embeddings $\mat{S} \in \real^{m\times d}$ with expected sparsity $\zeta$ satisfies the (oblivious) subspace injection property
    \begin{equation*}
        \norm{\mat{S}^*(\mat{B}\vec{z})} \ge (1-\eta_-) \norm{\mat{B}\vec{z}} \quad \text{for every } \vec{z} \in \real^n
    \end{equation*}
    with probability at least $1-\delta$ provided that
    \begin{equation*}
        d \ge \frac{16 \max\{n,6\log(2n/\delta)\}}{\eta_-^2} \quad \text{and} \quad \zeta \ge \frac{32 \log(n/\delta)}{\eta_-^2}.
    \end{equation*}
\end{fact} 

This result shows that an embedding dimension of $d \sim n$ and sparsity $\zeta \sim \log n$ to obtain a nontrivial (lower) distortion $\eta_- < 1$.
Tropp's result has explicit constants; however, the bounds are pessemistic.

\begin{remark}[Coherence--sparsity tradeoff]
    The full version of Tropp's bound establishes a \emph{coherence--sparsity tradeoff}\index{coherence--sparsity tradeoff} for iid sparse embeddings.
    The smaller the \emph{coherence}\index{coherence, of a matrix} of $\mat{B}$ (defined to be its largest leverage score, \cref{def:ridge-general}), the smaller the sparsity parameter $\zeta$ can be.
    Indeed, for the most incoherent matrices,\index{coherence, of a matrix} Tropp's result shows that an expected sparsity of $\zeta = \order(\log (n) / m)$ is sufficient to obtain the subspace injection property.\index{subspace injection}
    In this extreme case, the embedding has just $\expect[\nnz(\mat{S})] = \order(\log n)$ nonzero entries!
    These results allow the practitioner to employ very aggressive choices for the parameter $\zeta$ if they have strong assurances that the matrix is incoherent.
    On the other hand, for general use, it is advisable to use the ``safe'' parameter settings $d \sim n$ and $\zeta \sim \log n$.
\end{remark}

\myparagraph{Limitations}
A first limitation of iid sparse embeddings is that, for a worst-case subspace, they provably require an expected sparsity parameter $\zeta \sim \log n$ to achieve the subspace injection condition $\eta_- < 1$.

\begin{proposition}[IID sparse embeddings: Logarithmic sparsity is necessary]
    Consider a iid sparse embedding $\mat{S} \in \real^{m\times d}$ with expected sparsity $\zeta \le 0.5d$.
    There exists a matrix $\mat{Q}$ for which $\mat{S}$ is not a subspace injection for $\mat{Q}$ (i.e., the lower distortion has the trivial value $\eta_- = 1$) with 50\% probability unless $\zeta \ge 0.5 \log n$.
\end{proposition}

\begin{proof}
    Consider the adversarial matrix $\mat{Q} \in \real^{m\times n}$ from \cref{eq:difficult-example-Q}.
    The matrix $\mat{Q}$ has the effect of selecting the first $n$ columns of $\mat{S}^*$, so $\mat{S}^*\mat{Q} \in \real^{d\times n}$ is an iid sparse embedding with expected sparsity $\zeta$.

    To show that $\mat{S}$ has lower distortion $\eta_- = 1$ with at least 50\% probability, it is sufficient to show that $\mat{S}^*\mat{Q}$ has a zero column with at least 50\% probability. 
    Each column of $\mat{S}^*\mat{Q}$ has the following probability of being fully zero:
    \begin{equation*}
        p \coloneqq \left( 1 - \frac{\zeta}{d} \right)^d \ge \exp(-2\zeta).
    \end{equation*}
    Here, we used the numerical inequality $1-\alpha \ge \exp(-2\alpha)$, which is valid for all $\alpha \in [0,0.5]$.
    Thus, the probability that some column of $\mat{S}^*\mat{Q}$ is zero is $1 - (1-p)^n$.
    For this probability to be at most $1/2$, we must have $p \le 1 - 2^{-1/n}$.
    Observe that
    \begin{equation*}
        1 - 2^{-1/n} = 1 - \exp \left( - \frac{\log 2}{n}\right) \le \frac{\log 2}{n}.
    \end{equation*}
    Thus, by the chain of inequalities $\exp(-2\zeta) \le p \le 1 - 2^{-1/n} \le (\log 2)/n$ and some algebra, we deduce
    %
    \begin{equation*}
        \zeta \ge \frac{1}{2} \log \left( \frac{n}{\log 2} \right) > \frac{1}{2} \log n.
    \end{equation*}
    The stated result is proven.
\end{proof}

A second limitation of iid sparse embeddings, less significant in my eyes, is that it appears to require $\zeta \sim 1/\eta_-^2$ to achieve lower distortion $\eta_-$ \cite[Rem.~3.6]{CDD24}.
In contrast, sign embeddings (\cref{sec:sse}) allow $\zeta \sim 1/\eta$ to obtain distortion $\eta$.

\begin{figure}
    \centering
    \includegraphics[width=0.98\linewidth]{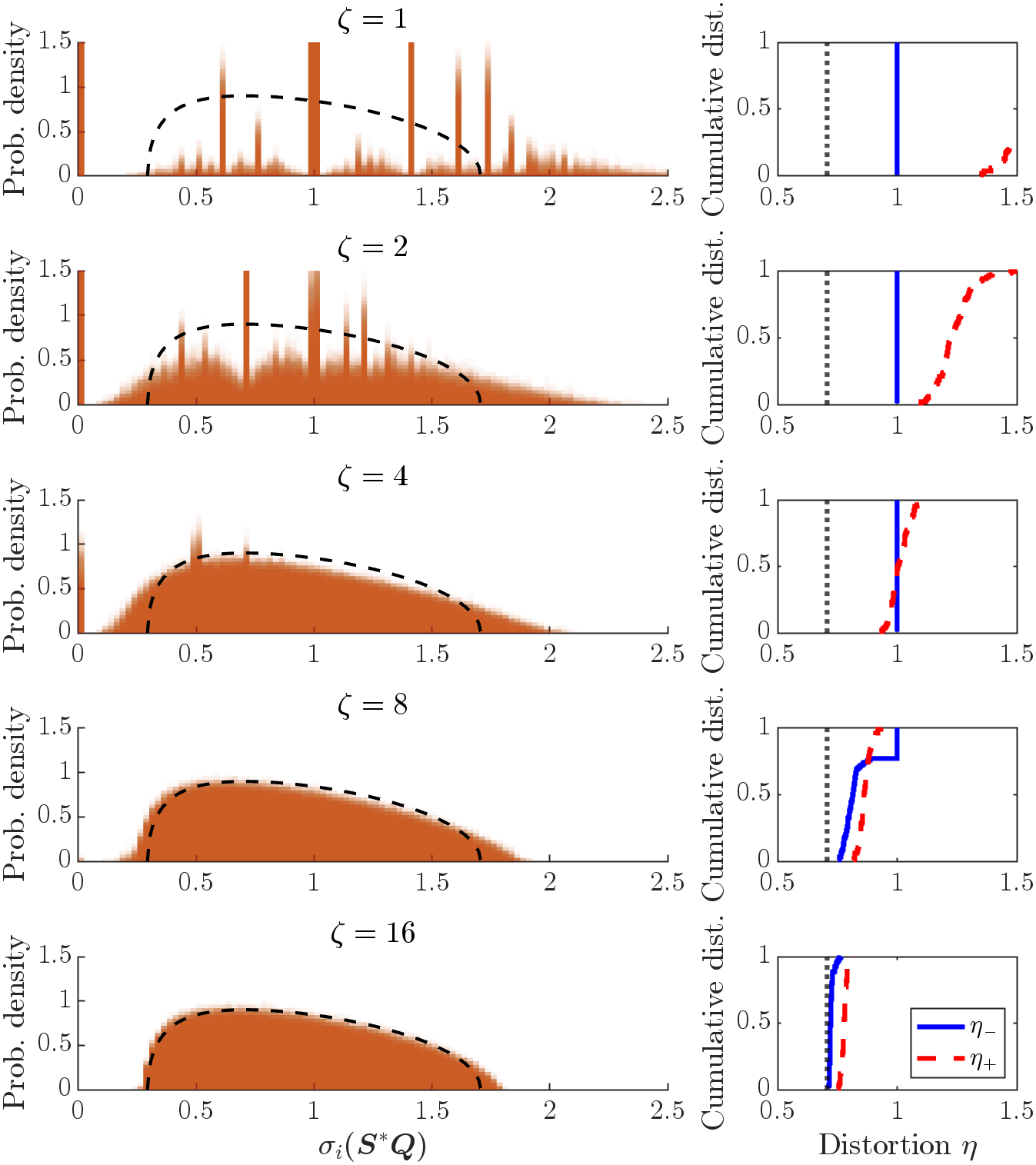}
    \caption[Singular values $\vec{\sigma}(\mat{S}^*\mat{Q})$ and distortions for iid sparse embeddings applied to a tough matrix $\mat{Q}$]{\emph{Left:} Superimposed singular value distributions $\vec{\sigma}(\mat{S}^*\mat{Q})$ for 100 independent realizations of a iid sparse embedding applied to the adversarial matrix \cref{eq:difficult-example-Q} with embedding dimension $d=2n$ for expected sparsity levels $\zeta \in \{1,2,4,8,16\}$ (\emph{top} to \emph{bottom}).
    The Mar\v{c}enko--Pastur singular value distribution is shown as the black dashed line.
    \emph{Right:} Empirical cumulative distribution function of lower and upper distortions.
    The asymptotic limiting value $\eta = \sqrt{n/d}$ is shown in black.}
    \label{fig:iid-sparse-sketch}
\end{figure}

\myparagraph{Empirical results}
\Cref{fig:iid-sparse-sketch} shows the results for iid sparse embeddings on the difficult instance \cref{eq:difficult-example-Q}.
We set $d \coloneqq 2n$ and we test sparsity parameters $\zeta = 2^i$ for $0\le i\le 4$.
We observe gradual convergence towards the Mar\v{c}enko--Pastur singular value distribution\index{Mar\v{c}enko--Pastur singular value distribution} as $\zeta$ increases.
To achieve a high probability of success, we require expected sparsity $\zeta > 8$ (for this problem size).

\section{Sparse sign embeddings} \label{sec:sse}

The weakness of iid sparse embeddings is in the name---the entries are generated iid.
Because of this, it is possible---likely, even---that full rows of $\mat{S}$ are zero unless $\zeta$ is chosen sufficiently large.
We can remedy this issue by using embeddings that have a fixed number of nonzeros per row.
This motivates the following definition:

\begin{definition}[Sparse sign embedding] \label{def:sparse-sign-appendix}
    A \emph{sparse sign embedding} $\mat{S} \in \real^{m\times d}$ with \emph{sparsity} $\zeta$ is a random matrix constructed as follows.
    Each row is independent and possesses exactly $\zeta$ nonzero entries.
    The nonzero entries are placed in uniformly random positions (selected \warn{without replacement}) and have uniform $\pm \zeta^{-1/2}$ values.
\end{definition}

See \cref{fig:sparse-sign-cartoon} in the main body for an illustration.

Sparse sign embeddings go by many names: sparse Johnson--Lindenstrauss embeddings (SJLTs, \cite{KN12}), oblivious sparse norm-approximating projections (OSNAPs, \cite{NN13a}), sparse sign embeddings \cite{MT20a}, hashing matrices \cite{CFS21}, short-axis sketching operators (SASOs, \cite{MDM+23a}), and sparse maps \cite{NT24}.
(We note, however, that some of these terms refer to a broader class of embeddings, with sparse sign embeddings---as defined above---serving as just one instance.)
As best I can tell, this construction is due to Kane and Nelson \cite{KN12}.
The case $\zeta = 1$ of the sparse sign embedding is known as a CountSketch,\index{CountSketch} which was introduced in the paper \cite{CCF04}.

\begin{figure}
    \centering
    \includegraphics[width=0.98\linewidth]{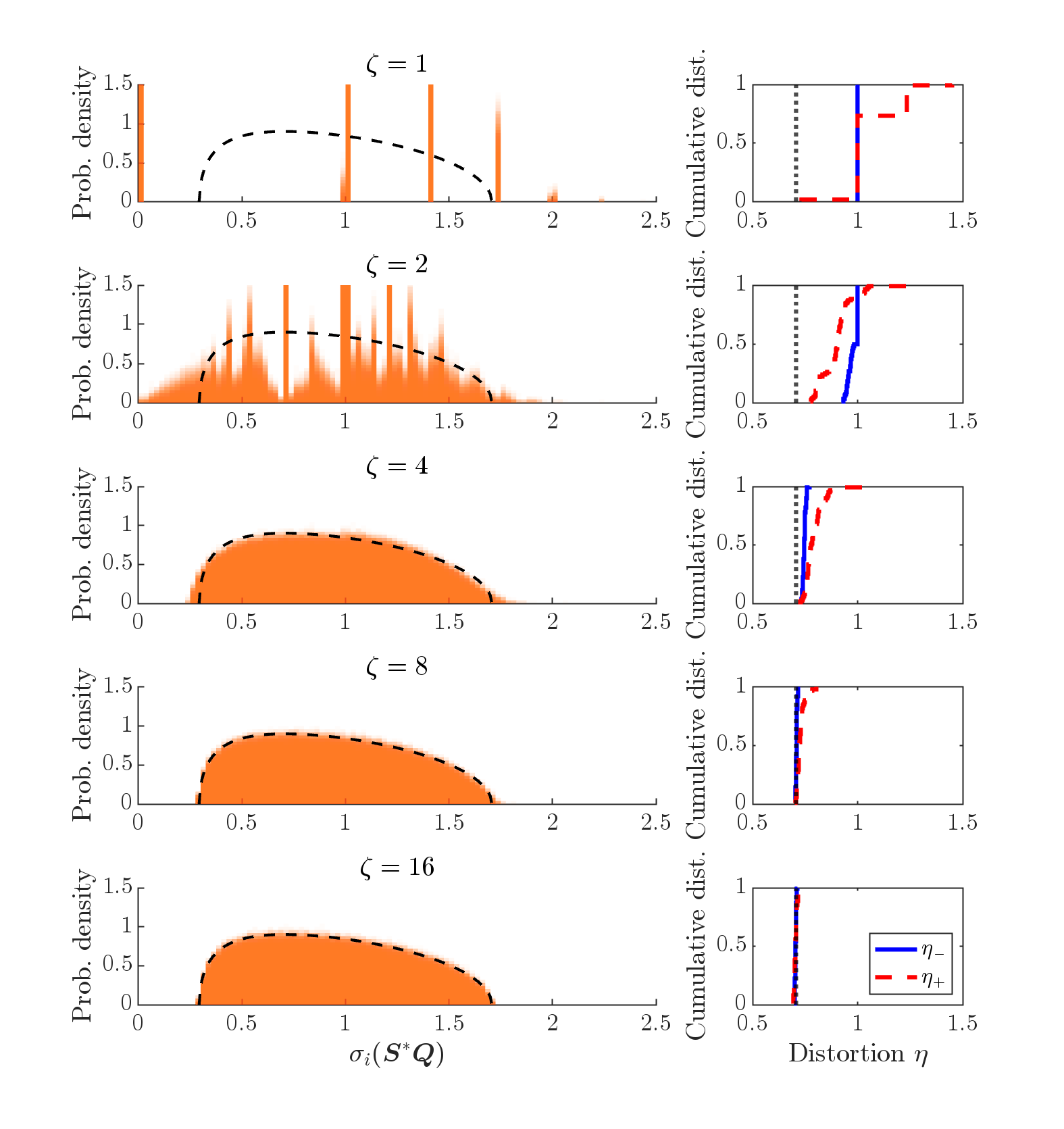}
    \caption[Singular values $\vec{\sigma}(\mat{S}^*\mat{Q})$ and distortions for sparse sign embeddings applied to a tough matrix $\mat{Q}$]{\emph{Left:} Superimposed singular value distributions $\vec{\sigma}(\mat{S}^*\mat{Q})$ for 100 independent realizations of a sparse sign embedding applied to the adversarial matrix \cref{eq:difficult-example-Q} with embedding dimension $d=2n$ for expected sparsity levels $\zeta \in \{1,2,4,8,16\}$ (\emph{top} to \emph{bottom}).
    The Mar\v{c}enko--Pastur singular value distribution is shown as a black dashed line for reference.
    \emph{Right:} Empirical cumulative distribution function of lower and upper distortions.
    The asymptotic limiting value $\eta = \sqrt{n/d}$ is shown in black.}
    \label{fig:sse-sketch}
\end{figure}

Let us begin this time with empirical results.
\Cref{fig:sse-sketch} shows the results for sparse sign embeddings on the difficult example \cref{eq:difficult-example-Q}.
We set $d \coloneqq 2n$ and we test sparsity parameters $\zeta = 2^i$ for $0\le i \le 4$.
We see that, for $\zeta \ge 4$, the sparse sign embedding is completely reliable, showing only slight excursions from the Mar\v{c}enko--Pastur singular value distribution\index{Mar\v{c}enko--Pastur singular value distribution} and producing a lower distortion $\eta_- \le 0.8$.
Compare with the $\zeta = 4$ and $\zeta = 8$ panels of \cref{fig:iid-sparse-sketch}, which show nonzero probbabilities of complete failure, $\eta_- = 1$.
This suggests that sparse sign embeddings ``work'' for smaller values of the sparsity parameter $\zeta$ than for iid sparse embeddings.

For smaller values of $\zeta$, the performance of sparse sign embeddings deteriorates.
The CountSketch\index{CountSketch} matrix (i.e., $\zeta=1$) fails catastrophically (i.e., $\eta_- = 1$) in every trial in this example.
The $\zeta = 2$ embedding also fails catastrophically about 50\% of the time.
Based on these examples, $\zeta = 3$ is the smallest value I would accept for general-purpose use (with embedding dimension $d\sim n$).
The CountSketch\index{CountSketch} embedding ($\zeta = 1$) achieves the subspace embedding property only with a \warn{quadratic} embedding dimension $d \sim n^2$ \cite[Thm.~23]{Woo14a}!

The experiments in thesis add to a growing body of evidence \cite{TYUC19,DM23,Epp23c,MBM+23,EMN24,CNR+25,CEMT25} that sparse sign embeddings are fast and completely reliable for linear algebraic computations.
The insights of this literature may be crystallized into a hypothesis:

\actionbox{Of the existing constructions, sparse sign embeddings are the most effective oblivious subspace embedding for general-purpose use in matrix computations.
They are reliable (i.e., achieve $\eta_-$ bounded away from 1) when implemented with an embedding dimension as low as $d = 2n$ and \warn{constant} sparsity parameter $\zeta=4$ or $\zeta = 8$.}

In the rest of this section, I will support this claim with evidence.

\subsection{Implementation and runtime comparison}

Making sparse sign embeddings fast in practice requires careful implementation and a good sparse linear algebra library.
Often, to obtain the greatest speed in a high-level programming language, such as MATLAB or Python, the construction of the sparse sign embedding should be written in a low-level programming language, like C or C\texttt{++}.

The challenge of implementing sparse sign embeddings is selecting the locations of the nonzero entries in each row \emph{without replacement}.
One simple approach, appropriate when $1 \approx \zeta \ll d$ is to use \emph{rejection sampling}.
Generate a stream indices $i \sim \Unif \{1,\ldots,d\}$, and check whether each index has already been picked, stopping when $\zeta$ indices have been selected.
This approach is illustrated in the following code segment.

\begin{lstlisting}
positions = zeros(zeta,1); i = 0;
while i <= zeta                               % Pick zeta positions
    idx = randsample(d,1);                    % Random index in [1,d]
    for j = 1:(i-1)                           % Check for duplicates
        if positions(j) == idx; continue; end % Reject duplicate
    positions(i) = idx; i = i + 1;            % Accept new position
end
\end{lstlisting}

Assuming $\zeta \le d/2$, the expected cost of generating the nonzero locationsis $\order(\zeta^2)$ operations per row.
For small values of $\zeta$ and with a low-level implementation in C, the overhead of this rejection sampling approach is minimal, and it produces sparse sign embeddings rapidly in practice.
C code for generating sparse sign embeddings with this construction appears in \cref{ch:sparse-sign-implementation}.

More efficient implementations of sparse sign embeddings have been developed as part of the \RandBLAS/\RandLAPACK software project.
In particular, using a careful implementation of Fisher--Yates sampling, the \RandBLAS software builds sparse sign embeddings $\order(\zeta m + d)$ operations, which is nearly optimal.
See \cite[App.~A]{MDM+23a} for details.

\begin{figure}
    \centering
    \includegraphics[width=0.99\linewidth]{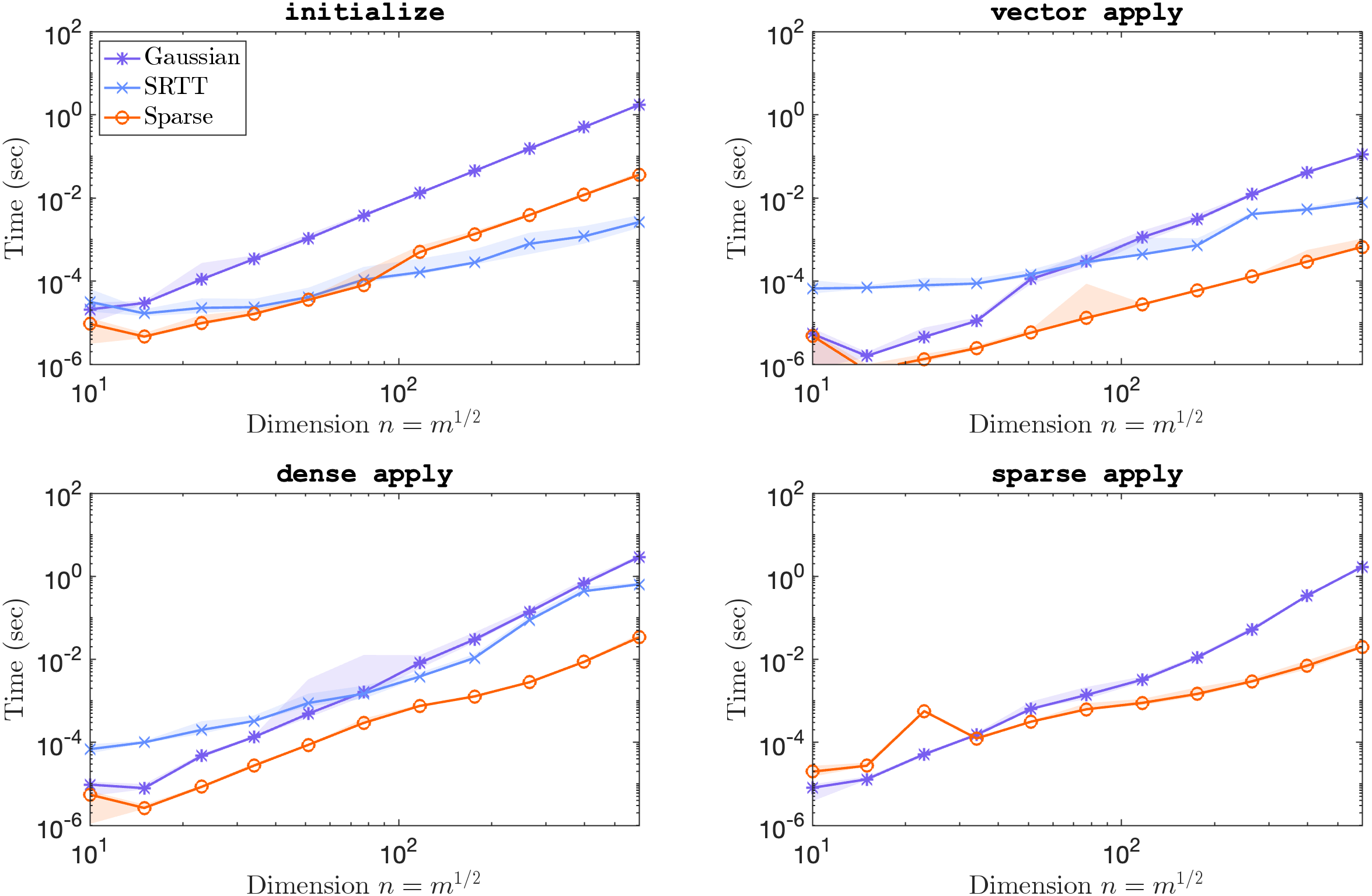}
    \caption[Runtimes for different sketching matrices]{Runtimes for initialization (\emph{top left}), matrix--vector apply (\emph{top right}), dense matrix--matrix apply (\emph{bottom left}), and sparse matrix--matrix apply (\emph{bottom right}) for Gaussian embeddings (purple asterisks), standard SRTTs (blue crosses), and sparse sign embeddings (orange circles).
    Lines show median of 100 trials, and shaded regions show 10\% and 90\% quantiles.}
    \label{fig:sketching-runtime-compare}
\end{figure}

\Cref{fig:sketching-runtime-compare} compares the runtimes of Gaussian embeddings, (standard) subsampled trigonometric transforms, and sparse sign embeddings (with sparsity $\zeta = 4$).
We consider a family of problems with dimensions $n = m^{1/2} \in [10,600]$, and we set $d\coloneqq 2n$ for each value of $n$.
We evaluate four components of the runtime: time to construct the embedding, i.e., generate the necessary random variables (\emph{top left}), time to compute a matrix--vector product $\mat{S}^*\vec{c}$ (\emph{top right}), time to compute the product $\mat{S}^*\mat{B}$ with a dense matrix $\mat{B} \in \real^{m\times n}$ (\emph{bottom left}), and time to multiply with a sparse matrix $\mat{B} \in \real^{m\times n}$ with 10 nonzeros per column (\emph{bottom right}).
(For the last of these, we omit the SRTT because it does not possess a fast multiply operation with sparse matrices.)
Sparse sign embeddings are definitively the fastest on all metrics except initialization.
In particular, for $n=600$, the dense matrix--matrix product operation for sparse sign embeddings is 19$\times$ faster than SRTTs and 86$\times$ faster than Gaussian embeddings.

\subsection{Analysis}

The best analysis of sparse sign embeddings remains the work of Cohen \cite{Coh16} (\cref{thm:ose-existence}).
I restate his theorem below.

\begin{fact}[Sparse sign embeddings: Subspace embedding property] \label{fact:cohen}
    Sparse sign embeddings are \warn{oblivious} subspace embeddings with dimension $n$, distortion $\eta \in (0,1)$, and failure probability $\delta \in (0,1)$ provided
    \begin{equation*}
        d \ge \Omega \left( \frac{n \log (n/\delta)}{\eta^2} \right) \quad \text{and} \quad \zeta \ge \Omega \left( \frac{\log(n/\delta)}{\eta} \right).
    \end{equation*}
\end{fact}

In particular, Cohen's results show that \emph{logarithmic sparsity} $\zeta = \order(\log n)$ and a \emph{log-linear embedding dimension} $d = \order(n\log n)$ suffice to obtain the subspace embedding guarantee.

Cohen's analysis helps explain the success of sparse embeddings.
However, it falls short of explaining the near-Gaussian behavior of sparse sign embeddings when implemented with a \emph{linear embedding dimension} $d = \order(n)$ and constant sparsity $\zeta = 4$.
Closing this gap remained an open problem in the years after Cohen's analysis.

While no further progress has been made on analyzing sparse sign embeddings, recent work has improved our understanding of related constructions.
We saw above the work of Tropp (\cref{fact:iid-sparse}, \cite{Tro25}), which considered iid sparse embeddings.
Another construction is given as follows:

\begin{definition}[SparseStack] \label{def:sparsestack}
    Let $d,m,\zeta \ge 1$ be integer parameters such that $\zeta$ divides evenly into $d$.
    An \emph{SparseStack} is a sparse matrix with exactly one nonzero entry per \emph{subcolumn} $\mat{S}^*((i-1)\zeta + 1 : i\zeta, j)$ of $\mat{S}^*$ for $1\le i \le d/\zeta$ and $1\le j \le m$.
    This entry is placed in a uniformly random position and has a uniformly random value $\pm 1/\sqrt{\zeta}$.
\end{definition}

SparseStacks were introduced alongside sparse sign embeddings in the original paper of Kane and Nelson \cite{KN12}.
The SparseStack model is studied in \cite{CDDR24,CDD24}.
The paper \cite{CDD24} contains the following result.

\begin{fact}[SparseStack] \label{fact:ind-subcols}
    SparseStacks are \warn{oblivious} subspace embeddings with dimension $n$, distortion $\eta \in (0,1)$, and failure probability $\delta \in (0,1)$ provided
    \begin{equation*}
        \Omega \left( \frac{n + \log(1/\delta)}{\eta^2} \right) \le d \le \order(\e^n) \quad \text{and} \quad \zeta \ge \Omega \left( \frac{\log^2(n/(\eta\delta))}{\eta} + \log^3(n/(\eta\delta) \right).
    \end{equation*}
\end{fact}

This bound shows that an embedding dimension of $d \sim n$ suffices for SparseStacks to achieve the subspace embedding property provided that the sparsity level is $\zeta \sim \log^3 n$.
Compare with Tropp's bound (\cref{fact:iid-sparse}), which states that $d \sim n$ and $\zeta \sim \log n$ imply the \warn{subspace injection} property for iid sparse maps.

To summarize, no theoretical results yet explain the success of sparse sign embeddings with the aggressive parameters $d=2n$ and a $\zeta = 4$.
The bound \cref{fact:cohen} of Cohen \cite{Coh16}, showing $d \sim n \log n$ and $\zeta \sim \log n$, almost reaches the parameter settings we use in practice, and recent works \cite{CDDR24,CDD24,Tro25} (analyzing slightly types of sparse embeddings) show that a linear embedding dimension $d \sim n$ suffices, provided one has a large enough sparsity parameter (e.g., $\zeta \sim \log n$ in \cite{Tro25}).

\subsection{Impossibility results: What they say and what they don't say}

Particularly sophisticated readers may be aware that there are impossibility results showing that sparse embeddings with $\zeta = 4$, as we recommend.
Cohen \cite[\S1]{Coh16}, summarizing \cite{NN14}, writes:
\begin{quote}
    Lower bounds (\cite{NN14}), on the other hand, \emph{suggest} that the true tradeoff allows $d = \order(B(n + \log(1/\delta))/ \eta^2)$ with $\zeta = \order(\log_B(n/\delta)/\eta)$.
\end{quote}
I have rephrased Cohen's notation to be consistent with ours, and the emphasis is added.
Cohen suggests that \cite{NN14} establishes a \emph{sparsity--embedding dimension tradeoff}, controlled by a continuous parameter $B > 0$.
In particular, to obtain a linear embedding dimension $d \sim  n$, Cohen suggests that a logarithmic embedding dimension $\zeta \sim \log n$ is \emph{necessary} to achieve the subspace embedding property.
If true, this provide given strong theoretical evidence that constant-sparsity embeddings could fail for some instances.

A careful reading of the paper \cite{NN14} yields a more nuanced conclusion.
Specifically, upon close inspections, the arguments of \cite[\S3]{NN14} yield the following result.

\begin{fact}[Sparse embeddings: Lower bound] \label{fact:sparse-lower-bound}
    Assume $m\ge 100n^2$.
    Any sparse oblivious subspace embedding $\mat{S} \in \real^{m\times 2n}$ with embedding dimension $d = 2n$, distortions $\eta_+, \eta_-$, and failure probability $\delta = 0.02$ must have
    \begin{equation*}
        \zeta \ge \mathrm{const} \cdot \left( \frac{1-\eta_-}{1+\eta_+} \right)^2 \frac{\log n}{\log \log n}
    \end{equation*}
\end{fact}

We make a few observations.
First, to achieve a subspace embedding with $\eta < 1$ and embedding dimension $d = 2n$, this bound shows we need a sparsity level of $\zeta \sim (\log n)/ \log \log n$.
This is slightly weaker than the suggested $\zeta = \Omega(\log n)$ lower bound in Cohen's quote above.
Second, and more importantly, \cref{fact:sparse-lower-bound} leaves open the possibility of constant sparsity embeddings $\zeta = \mathrm{const}$ satisfying the subspace embedding property with
\begin{equation} \label{eq:conj-distortions}
    \eta_- = \mathrm{const} < 1 \quad \text{and} \quad \eta_+ = \order \left( \sqrt{\frac{\log n}{\log \log n}} \right).
\end{equation}
Under the sketching asymmetry principle, an embedding with these parameters will produce meaningful results in most applications.
Recall: For a sketching matrices to be useful, the lower distortion $\eta_-$ must be away from $1$, but we can tolerate the upper distortion $\eta_+$ to be logarithmically large

Let us clarify that \cref{fact:sparse-lower-bound} does not \emph{prove} that sparse sign embeddings with parameters $d=2n$ and $\zeta = \mathrm{const}$ achieve the distortion scaling in \cref{eq:conj-distortions}; it just states that these distortions are consistent with the best-known information-theoretic lower bounds from \cite[\S3]{NN14}.
I conjecture that sparse sign embeddings attain these bounds.

\begin{conjecture}[Sparse sign embeddings: Linear embedding dimension and constant sparsity]
    Sparse sign embeddings with parameters $d = 2n$ and $\zeta = 4$ are oblivious subspace embeddings with dimension $n$, failure probability $\delta = 0.01$, and distortions given by \cref{eq:conj-distortions}.
\end{conjecture}

As some evidence for this c,kaim, we can use the techniques of \cite{NN14} to show that the upper distortion satisfies $\eta_+ \ge \Omega((\log n)/\log \log n)$ with high probability.

\begin{proposition}[Sparse sign embeddings: Upper distortion]
    Consider a sparse sign embedding $\mat{S} \in \real^{m\times d}$ with parameters $d = 2n$ and $\zeta = \mathrm{const}$.
    Then, there exists a matrix $\mat{Q}$ for which the upper distortion of $\mat{S}$ is 
    \begin{equation*}
        \eta_+ \ge \mathrm{const} \cdot \sqrt{\frac{\log n}{\log \log n}}
    \end{equation*}
    with 99\% probability.
\end{proposition}

\begin{proof}[Proof sketch]
    The idea is similar to the analysis of \cite{NN14}.
    Set $\mat{Q}$ to be the ``the adversarial matrix'' \cref{eq:difficult-example-Q}.
    Then $\mat{S}^*\mat{Q}$ is a sparse random matrix with $\zeta$ entries in each column in uniformly random positions, with random values $\pm \zeta^{-1/2}$.
    Think of each row of $\mat{S}^*\mat{Q}$ as a bin and each nonzero entry as a ball dropped into a uniformly random bin (with the constraint that, within each column, no two balls are dropped in the same bin).
    Up to the dependencies within each column (which become insignificant in the limit $n\to\infty$, since $\zeta = \mathrm{const}$), this is an instance of the classical balls-and-bins problem, for which it is known the heaviest bin has at least $\mathrm{const} \cdot (\log n)/ \log \log n)$ balls with high probability (see, e.g., \cite[Lem.~14]{NN14}).
    Thus, with high probability,
    \begin{multline*}
        \eta_+ = \sigma_{\mathrm{max}}(\mat{S}^*\mat{Q}) \ge \max_{1\le i \le d} \norm{(\mat{S}^*\mat{Q})(i,:)} \\ \ge \sqrt{\mathrm{const} \cdot \frac{\log n}{\log \log n} \cdot \left( \pm \zeta^{-1/2}\right)^2} = \mathrm{const}\cdot  \sqrt{\frac{\log n}{\log \log n}}.
    \end{multline*}
    The argument is complete.
\end{proof}







\section{Conclusions} \label{sec:sketching_conclusions}

In my experience working, the question ``Which sketch should you use?'' is shrouded in folklore and accounts of personal experience.
Despite the importance of sketching to modern randomized matrix computations, there seems to be little consensus on which types of sketching matrices to use in practice.
In this chapter, I have evaluated different sketching matrices on a challenging test matrix and summarized the relevant theoretical literature.
On the basis of this evidence, I and others believe sparse sign embeddings are a strong candidate for the sketching matrix of choice for most applications.
Still, there remain open questions, and it is worth further investigations into the theoretical properties and empirical behavior of sparse sign embeddings and other sketching matrices.

\section{Postscript: Recent developments} \label{sec:postscript}

Subsequent to my defense and initial drafting of this thesis, Chris Cama\~no, Raphael Meyer, Joel Tropp, and myself released a paper \cite{CEMT25} revisiting the theory of sketching.
In particular, it shows that subspace injection matrices satisfying the isotropy condition
\begin{equation*}
    \expect [\norm{\smash{\mat{S}^*\vec{v}}}^2] = \norm{\vec{v}}^2 \quad \text{for all } \vec{v} \in \field^n
\end{equation*}
produce sketch-and-solve solutions $\hatvector{x} \in \field^n$ that are accurate up to a constant factor
\begin{equation*}
    \norm{\vec{c} - \mat{B}\hatvector{x}} \le \mathrm{const} \cdot \norm{\vec{c} - \mat{B}\vec{x}}.
\end{equation*}
It also analyzes the use of such random matrices in other randomized linear algebra primitives, including the randomized SVD (\cref{sec:rsvd}), Nystr\"om approximation (\cref{sec:nystrom}), and generalized Nystr\"om approximation (\cref{sec:gen-nys}).
This paper also provides more evidence for the reliability of sparse sketching matrices, and it discusses how SparseStack embeddings (\cref{def:sparsestack}) may be easier to construct than sparse sign embeddings (\cref{def:sparse-sign-appendix}).
Perhaps it is the SparseStack, not the sparse sign embedding, that should be the sketching matrix of choice for most applications!

As another late-breaking update, Chenakkod, Derezi\'nski, and Dong released a new preprint \cite{CDD25a} containing an even sharper analysis of SparseStack embeddings.
They show that SparseStacks are oblivious subspace embeddings with parameters
\begin{equation*}
    d = \order\left(\frac{n\log^{\order(1/\log\log\log\log n)}(n)}{\eta^2}\right) \quad \text{and} \quad \zeta = \order\left( \frac{\log^{1+\order(1/\log\log\log\log n)}(n)}{\eta} \right).
\end{equation*}
This result improves on \cref{fact:ind-subcols}, and it achieves the lower bound of \cite{NN14} \warn{for oblivious subspace embeddings} up to \emph{subpolylogarithmic} factors.

\chapter{Analysis of sketch-and-solve} \label{app:sketch-and-solve-analysis}

The classical result for sketch-and-solve, cited in many popular surveys and expository works \cite{MT20a,MDM+23a,KT24} and many of my own works, is as follows:

\begin{proposition}[Sketch-and-solve: Simple bound] \label{prop:sketch-and-solve-simple}
    Let $\mat{S}$ be a subspace embedding \warn{for $\flatonebytwo{\mat{B}}{\vec{c}}$} with lower and upper distortions $\eta_-,\eta_+$ (\cref{def:subspace-embed-2}), and let $\hatvector{x}$ be the corresponding sketch-and-solve solution.
    Then
    \begin{equation*}
        \norm{\vec{c} - \mat{B}\hatvector{x}} \le \frac{1+\eta_+}{1-\eta_-} \cdot \norm{\vec{c} - \mat{B}\vec{x}}.
    \end{equation*}
    In particular, using a sparse sign embedding with Cohen's parameter settings (\cref{fact:cohen}), sketch-and-solve produces $(1+\varepsilon)$-approximate least-squares solution with 99\% probability provided 
    \begin{equation*}
        d = \order\left( \frac{n \log n}{\varepsilon^2} \right).
    \end{equation*}
    The runtime is 
    \begin{equation*}
        \order \left( mn \log n + \frac{n^3 \log n}{\varepsilon^2} \right) \text{ operations}.
    \end{equation*}
\end{proposition}

\begin{proof}
    We make the following chain of inequalities
    \begin{equation*}
        \norm{\vec{c} - \mat{B}\hatvector{x}} \le \frac{1}{1-\eta_-} \cdot \norm{\mat{S}^*(\vec{c} - \mat{B}\hatvector{x})} \le \frac{1}{1-\eta_-} \cdot \norm{\mat{S}^*(\vec{c} - \mat{B}\vec{x})} \le \frac{1+\eta_+}{1-\eta_-} \cdot \norm{\vec{c} - \mat{B}\vec{x}}.
    \end{equation*}
    The first and third inequality are subspace embedding property, and the second inequality is the optimality of $\hatvector{x}$ as a solution to the sketched least-squares problem $\argmin_{\vec{z} \in \field^n} \norm{\mat{S}^*(\vec{c} - \mat{B}\hatvector{x})}$.
\end{proof}

This bound suggests that the residual norm of the sketch-and-solve solution is within a factor $1+\order(\eta)$ of the optimal residual norm and that an embedding dimension $d \sim 1/\varepsilon^2$ is necessary to obtain a $(1+\varepsilon)$-approximate least-squares solution.
In fact, this scaling can be improved using a sharper analysis:

\begin{theorem}[Sketch-and-solve: Sharper bound] \label{thm:sketch-and-solve-sharp}
    Let $\mat{S}$ be a subspace embedding \warn{for $\flatonebytwo{\mat{B}}{\vec{c}}$} with lower and upper distortion $\eta_-,\eta_+$ (\cref{def:subspace-embed-2}) and let $\hatvector{x} = (\mat{S}^*\mat{B})^\dagger(\mat{S}^*\vec{c})$ be the sketch-and-solve solution.
    Assume $\eta \coloneqq \max \{ \eta_+,\eta_-\} \le 1$.
    Then
    \begin{equation*}
        \norm{\vec{c} - \mat{B}\hatvector{x}} \le \left(1 + \frac{4.5\eta^2}{(1-\eta_-)^4}\right) \norm{\vec{c} - \mat{B}\vec{x}}.
    \end{equation*}
    In particular, using a sparse sign embedding with Cohen's parameter settings (\cref{fact:cohen}), sketch-and-solve produces $(1+\varepsilon)$-approximate least-squares solution with 99\% probability provided 
    \begin{equation*}
        d = \order\left( \frac{n \log n}{\varepsilon} \right).
    \end{equation*}
    The runtime is 
    \begin{equation*}
        \order \left( mn \log n + \frac{n^3 \log n}{\varepsilon} \right) \text{ operations}.
    \end{equation*}
\end{theorem}

This result establishes represents a significant improvement over \cref{prop:sketch-and-solve-simple} in the limit where $\eta \to 0$ (equivalently $d \to \infty$), and the resulting $1+\order(\eta^2)$ scaling is sharp, as demonstrated by exact computations for Gaussian embeddings \cite{BP20,Epp24d}.
The sharp $d \sim 1/\varepsilon$ scaling for sketch-and-solve has been known since the work of Drineas, Mahoney, Muthukrishnan, and Sarl\'os \cite{DMMS11}, originally released in 2007; see also the expository notes \cite{LK22,Mey23}.
However, these earlier bounds pass through a different technical argument and require hypotheses on the sketching matrix $\mat{S}$ beyond the subspace embedding property.
This bound, a version of which was published in the blog post \cite{Epp25}, is the only result I know which established the correct $1 + \order(\eta^2)$ scaling with the subspace embedding property as the sole hypothesis.

We proceed now with its proof.

\begin{proof}[Proof of \cref{thm:sketch-and-solve-sharp}]
    We shall bound $\norm{\mat{B}(\hatvector{x} - \vec{x})}$.
    To do so, we use the representation
    \begin{equation*}
        \mat{B}(\hatvector{x} - \vec{x}) = \mat{B}(\mat{S}^*\mat{B})^\dagger \mat{S}^*(\vec{c} - \mat{B}\vec{x}).
    \end{equation*}
    Letting $\mat{B} = \mat{Q}\mat{R}$ be an (economy-size) \QR decomposition, we have
    \begin{equation*}
        \mat{B}(\hatvector{x} - \vec{x}) = \mat{Q} (\mat{S}^*\mat{Q})^\dagger \mat{S}^*(\vec{c} - \mat{B}\vec{x}) = \mat{Q} (\mat{Q}^*\mat{S}\mat{S}^*\mat{Q})^{-1} \mat{Q}^*\mat{S}\mat{S}^*(\vec{c} - \mat{B}\vec{x}).
    \end{equation*}
    Thus, introducing $\vec{\overline{r}} \coloneqq (\vec{c} - \mat{B}\vec{x}) / \norm{\vec{c} - \mat{B}\vec{x}}$, we have
    \begin{equation} \label{eq:resid-error-upper-bound}
        \norm{\mat{B}(\hatvector{x} - \vec{x})} \le \norm{(\mat{Q}^*\mat{S}\mat{S}^*\mat{Q})^{-1}} \cdot \norm{\mat{Q}^*\mat{S}\mat{S}^*\vec{\overline{r}}} \cdot \norm{\vec{c} - \mat{B}\vec{x}}.
    \end{equation}

    We now bound the right-hand side of \cref{eq:resid-error-upper-bound}.
    By \cref{prop:svals-subspace}, we have
    \begin{equation} \label{eq:resid-error-upper-bound-1}
        \norm{(\mat{Q}^*\mat{S}\mat{S}^*\mat{Q})^{-1}} = \frac{1}{\sigma_{\mathrm{min}}(\mat{S}^*\mat{Q})^2} \le \frac{1}{(1-\eta_-)^2}.
    \end{equation}
    Bounding $\norm{\mat{Q}^*\mat{S}\mat{S}^*\vec{\overline{r}}}$ requires a more sophisticated argument.
    Observe that
    \begin{equation*}
        \norm{\mat{Q}^*\mat{S}\mat{S}^*\vec{\overline{r}}} \text{ is a submatrix of }\flatonebytwo{\mat{Q}}{\vec{\overline{r}}}^*\mat{S}\mat{S}^*\flatonebytwo{\mat{Q}}{\vec{\overline{r}}} - \Id.
    \end{equation*}
    The columns of the matrix $\flatonebytwo{\mat{Q}}{\vec{\overline{r}}}$ form an orthonormal basis for $\range(\flatonebytwo{\mat{B}}{\vec{c}})$.
    Therefore, by \cref{prop:svals-subspace} again, we have
    \begin{equation} \label{eq:resid-error-upper-bound-2}
        \norm{\mat{Q}^*\mat{S}\mat{S}^*\vec{\overline{r}}} \le \norm{\flatonebytwo{\mat{Q}}{\vec{\overline{r}}}^*\mat{S}\mat{S}^*\flatonebytwo{\mat{Q}}{\vec{\overline{r}}} - \Id} \le \max \{ (1+\eta_+)^2 - 1, 1 - (1-\eta_-)^2 \}.
    \end{equation}
    
    Combining \cref{eq:resid-error-upper-bound-1,eq:resid-error-upper-bound-2} and substituting in \cref{eq:resid-error-upper-bound} yields
    \begin{equation*}
        \norm{\mat{B}(\hatvector{x} - \vec{x})} \le \frac{2\eta + \eta^2}{(1-\eta_-)^2} \norm{\vec{c} - \mat{B}\vec{x}}.
    \end{equation*}
    Combining with the Pythagorean identity,
    \begin{equation*}
        \norm{\vec{c} - \mat{B}\hatvector{x}}^2 = \norm{\vec{c} - \mat{B}\vec{x}}^2 + \norm{\mat{B}( \hatvector{x} - \vec{x})}^2 \le \left(1 + \frac{(2\eta+\eta^2)^2}{(1-\eta_-)^4} \right)\norm{\vec{c} - \mat{B}\vec{x}}^2.
    \end{equation*}
    Finally, with the hypothesis, $\eta \le 1$, we have $(2\eta+\eta^2)^2 \le 9\eta^2$.
    Take square roots and use the identity $(1+x)^{1/2} \le 1 + 0.5x$ for $x\ge -1$ to complete the argument.
\end{proof}
  
\chapter{Deferred proofs}

This section contains deferred proofs of results from the main text.

\section{\texorpdfstring{Proof of \cref{thm:optimal-nys-approx-2}}{Proof of Theorem 3.12}} \label{sec:optimal-nys-approx-2}

We use a more algebraic version of the proof of \cite[Prop.~3.3]{DRVW06}.
Fix a sufficiently small parameter $\delta > 0$ and generate the matrix
\begin{equation*}
    \mat{A} = (k+1)\Id_{k+1} - (1-\delta) \outprod{\onevec_{k+1}} \in \real^{(k+1)\times (k+1)}.
\end{equation*}
This matrix has eigenvalue $(k+1)$ with multiplicity $k$ and eigenvalue $(k+1)\delta$ with multiplicity one.
Therefore, 
\begin{equation*}
    \tr(\mat{A} - \lowrank{\mat{A}}_k) = (k+1)\delta.
\end{equation*}
Since the matrix $\mat{A}$ is the same under any conjugation by a permutation matrix, we can assume without loss of generality that $\set{S} = \{1,\ldots,k\}$.
We compute
\begin{equation*}
    \tr(\mat{A} - \mat{A}\langle \set{S}\rangle) = a_{(k+1)(k+1)} - \vec{a}_{k+1}(1:k)^* \mat{A}(1:k,1:k)^{-1}\vec{a}_{k+1}^{\vphantom{*}}(1:k).
\end{equation*}
Applying the Sherman--Morrison formula gives
\begin{equation*}
    \mat{A}(1:k,1:k)^{-1} = \frac{1}{k+1} \cdot \Id_k + \frac{1-\delta}{(k+1)(1 + k\delta)} \cdot \outprod{\onevec_k}.
\end{equation*}
Thus,
\begin{align*}
    \tr(\mat{A} - \mat{A}\langle \set{S}\rangle) &= (k + \delta) - \frac{(1-\delta)^2k}{k+1} - \frac{(1-\delta)^3k^2}{(k+1)(1+k\delta)} = (k+1)^2 \delta + \order(\delta^2).
\end{align*}
Taking $\delta$ sufficiently small proves the claim.\hfill $\qedsymbol$

\section{\texorpdfstring{Proof of \cref{thm:weighted-cur}}{Proof of Theorem 10.8}} \label{sec:weighted-cur-proof}

    The existence result for weighted \CUR decompositions follows from \cref{fact:lev-ls}.
    We focus on establishing the nonexistence result for unweighted \CUR approximations.
    Our argument is based on a construction of Derezi\'nski, Warmuth, and Hsu \cite[Thm.~1]{DWH18}.
    First consider the case $k=1$.
    
    Consider the matrix
    \begin{equation*}
        \mat{B} = \begin{bmatrix}
            1 & 1 \\
            p^{-1/2} & 0 \\
            p^{-1/2} & 0 \\
            \vdots & \vdots \\
            p^{-1/2} & 0
        \end{bmatrix} \in \real^{(p+1)\times 2},
    \end{equation*}
    and choose column subset $\set{S} \coloneqq \{1\}$.
    The column projection approximation is
    \begin{equation*}
        \Bhat_{\mathrm{CP}} = \mat{B}(:,\set{S})\mat{B}(:,\set{S})^\dagger \mat{B} = \begin{bmatrix}
            \vec{b}_1 & \frac{1}{2}\vec{b}_1
        \end{bmatrix},
    \end{equation*}
    with error
    \begin{equation*}
        \norm{\mat{B} - \Bhat_{\mathrm{CP}}}_{\mathrm{F}}^2 = \frac{1}{2}.
    \end{equation*}
    
    Let $\set{T} \subseteq \{1,\ldots,p+1\}$ denote any row subset of $\ell$ elements of $\mat{B}$.
    If $1\notin \set{T}$, then the unweighted \CUR cross approximation is zero and thus $\norm{\mat{B} - \Bhat}_{\mathrm{F}}^2 = 3$.
    Thus, we may assume $1 \in \set{T}$.
    Then the unweighted \CUR cross approximation is
    \begin{equation*}
        \Bhat = \mat{B}(:,\set{S})\mat{B}(\set{T},\set{S})^\dagger \mat{B}(\set{T},:) = \onebytwo{\vec{b}_1}{\frac{\vec{b}_1(\set{T})^*\vec{b}_2(\set{T})}{\vec{b}_1(\set{T})^*\vec{b}_1(\set{T})}\vec{b}_1} = \onebytwo{\vec{b}_1}{\frac{\vec{b}_1}{1+(\ell-1)/p}}
    \end{equation*}
    with error
    \begin{equation} \label{eq:cur-k-1}
        \norm{\mat{B} - \Bhat}_{\mathrm{F}}^2 = 1 - \frac{2}{1+(\ell-1)/p} + \frac{2}{(1+(\ell-1)/p)^2} \eqqcolon f(\ell).
    \end{equation}

    The function $f$ is strictly decreasing on the interval $[1,p+1]$ and it achieves the value $3/4$ at $\ell = 1 + (3-2\sqrt{2})p$.
    The result is proven for $k=1$.

    To prove the result for $k > 1$, generate a block diagonal matrix with $k$ copies of $\mat{B}/\sqrt{k}$, and set $\set{S} = \{1,p+1,\ldots,(k-1)p+1\}$.
    The optimal squared Frobenius column projection error is again $1/2$, and the \CUR error is bounded by 
    \begin{equation} \label{eq:cur-lower-bound}
        \norm{\mat{B} - \Bhat}_{\mathrm{F}}^2 \ge \frac{1}{k} \left(f(\theta_1 \ell) + \cdots + f(\theta_2\ell) + \cdots f(\theta_k\ell)\right),
    \end{equation}
    where $\theta_j\ell$ denotes the number of indices assigned to block $j$.
    Here, to handle boundary cases, the values $f(t)$ are defined as \cref{eq:cur-k-1} for $0<t\le p$, $f(0) = 3$, and $f(t) = 1/2$ for $t > p$.
    With this definition, $f$ is convex on $\real_+$.
    Therefore, by Jensen's inequality, the right-hand side of \cref{eq:cur-lower-bound} is minimized by the equal apportionment $\theta_1 = \cdots = \theta_k = 1/k$:
    \begin{equation*}
        \norm{\mat{B} - \Bhat}_{\mathrm{F}}^2 \ge \frac{1}{k} \left(f(\theta_1 \ell) + \cdots + f(\theta_2\ell) + \cdots f(\theta_k\ell)\right) \ge f(\ell / k).
    \end{equation*}
    Ergo, the optimal unweighted \CUR decomposition has a squared Frobenius norm error at least $1.5\times$ the squared Frobenius column projection error provided $\ell \le k+(3-2\sqrt{2})pk$.\hfill $\qedsymbol$

\chapter{Implementation of sparse random embeddings} \label{ch:sparse-sign-implementation}

\begin{lstlisting}[
    frame=single, 
    language=C,
    caption={MEX (MATLAB-callable C function) implementation of sparse sign embeddings.}, 
    label={prog:sparse_sign}
  ]
#include <math.h> 
#include "mex.h"

int log_base_d(long a, int d) {
  int output = 0;
  while (a > 0) {
    a /= d;
    output += 1;
  }
  return output-1;
}

void mexFunction(int nlhs, mxArray *plhs[],
                 int nrhs, const mxArray *prhs[]) {

  mwSize d = (mwSize) mxGetScalar(prhs[0]);
  mwSize m = (mwSize) mxGetScalar(prhs[1]);
  mwSize zeta = (mwSize) mxGetScalar(prhs[2]);
  mwSize nnz = m*zeta;

  int idx_per_rand = log_base_d((long) RAND_MAX + 1, d);
  int bit_per_rand = log_base_d((long) RAND_MAX + 1, 2);

  if (zeta > d) zeta = d;

  double lowval = -1/sqrt((double) zeta);
  double increment = -2*lowval;

  plhs[0] = mxCreateSparse(d,m,nnz,false);
  double *vals  = mxGetPr(plhs[0]);
  mwIndex *rows = mxGetIr(plhs[0]);
  mwIndex *colstarts = mxGetJc(plhs[0]);

  // Set values
  int myrand = rand();
  for (int i = 0; i+bit_per_rand < nnz; i += bit_per_rand) {
    for (int j = i; j < i+bit_per_rand; ++j) {
      vals[j] = (myrand % 2) * increment + lowval;
      myrand = myrand >> 1;
    }
    myrand = rand();
  }
  for (int i = bit_per_rand*(nnz/bit_per_rand); i < nnz; ++i) {
    vals[i] = (myrand % 2) * increment + lowval;
    myrand = myrand >> 1;
  }

  // Set column starts
  for (int i = 1; i < m+1; ++i){
    colstarts[i] = i*zeta;
  }

  // Set row indices
  myrand = rand();
  int ir = 0;
  for (int i = 0; i < m*zeta; i += zeta) {
    int idx = 0;
    while (idx < zeta) {
      rows[i+idx] = myrand % d;
      ir++;
      if (ir == idx_per_rand) {
	ir = 0;
	myrand = rand();
      } else {
	myrand /= d;
      }
      int j = 0;
      for (; j < idx; ++j) {
	if (rows[i+idx] == rows[i+j]) break;
      }
      idx += (int) (j == idx);
    }
  }

}
\end{lstlisting}

\chapter{Helpful MATLAB subroutines}

\myprogram[h]{Code to normalize the columns of a matrix.}{}{cnormc}

\myprogram[h]{Generate matrix of random signs.}{}{random_signs}

\myprogram[h]{Compute squared row norms of a matrix.}{}{sqrownorms}

\myprogram[h]{Compute squared column norms of a matrix.}{}{sqcolnorms}

\myprogram[h]{Preconditioned conjugate gradient for solving positive definite linear systems.}{}{mypcg}

\myprogram[h]{Generate a Hermitian matrix with Haar-random eigenvectors and the specified eigenvalues.}{}{rand_with_evals}

\myprogram[h]{Generate a Haar-random matrix with orthonormal columns.}{}{haarorth}

\printbibliography[heading=bibintoc]

\end{document}